\pdfoutput=1

\documentclass[onecolumn, 12 pt, letterpaper]{report}

\usepackage{amssymb}
\usepackage{titletoc,tocloft}
\usepackage{textcomp}
\usepackage{graphicx}
\usepackage{graphics}
\usepackage{epsfig}
\usepackage{epstopdf}
\usepackage{float}
\usepackage{color}
\usepackage[cmex10]{amsmath}
\usepackage{latexsym,amsfonts}
\usepackage{amsthm}
\usepackage{url}
\usepackage{longtable}
\usepackage[figuresright]{rotating}
\usepackage{listings}
\usepackage{etoolbox}
\usepackage{sectsty}
\usepackage{array}
\usepackage{bibentry}
\nobibliography*

\makeatletter\let\saved@bibitem\@bibitem\makeatother
\usepackage[colorlinks,linkcolor=.,citecolor=.,pagebackref=true,backref=true,bookmarks=false]{hyperref}
\renewcommand*{\backrefalt}[4]{%
    \ifcase #1 \footnotesize{(Not cited.)}%
    \or        \footnotesize{(Cited on page~#2)}%
    \else      \footnotesize{(Cited on pages~#2)}%
    \fi}

\usepackage{cleveref}
\makeatletter\let\@bibitem\saved@bibitem\makeatother

\usepackage{threeparttable, adjustbox}
\usepackage{booktabs}       

\chapterfont{\fontsize{14}{15}\selectfont}
\sectionfont{\fontsize{14}{15}\selectfont}
\subsectionfont{\fontsize{14}{15}\selectfont}
\subsubsectionfont{\fontsize{14}{15}\selectfont}

\usepackage{fancyhdr}
\fancyheadoffset[L]{0.25in}

\makeatletter
\patchcmd{\@makechapterhead}{50\p@}{20pt}{}{}
\patchcmd{\@makeschapterhead}{50\p@}{20pt}{}{}
\makeatother

\fancypagestyle{plain}{
\fancyhf{} 
\fancyhead[C]{\thepage} 

}

\usepackage{setspace}

   \usepackage{ifthen,xkeyval,xfor,amsgen}
   \usepackage[acronym,toc, nogroupskip]{glossaries}
   \newglossary[slg]{symbols}{syi}{sbl}{List of Symbols}
 
   \makeglossaries
   
 \usepackage{multirow}

\usepackage{amsfonts}
\usepackage{amsmath}
\usepackage{amsthm}
\usepackage{amssymb} 
\usepackage{mathrsfs} 

\usepackage{mdframed} 
\usepackage{thmtools}

\usepackage{enumerate}
\usepackage[shortlabels]{enumitem} 
\usepackage{pifont}
\usepackage{enumitem}
\usepackage{import}

\usepackage{mathtools} 

\newcolumntype{C}[1]{>{\centering\let\newline\\\arraybackslash\hspace{0pt}}m{#1}}

 \usepackage{caption}
 \usepackage{subcaption}

 \usepackage{xcolor}
\usepackage{color}
\usepackage{graphicx}

\usepackage{algpseudocode} 
\usepackage{algorithmicx}
\usepackage{algorithm}
\usepackage{xspace} 

\newcommand{\R}{\mathbb{R}} 
\newcommand{\N}{\mathbb{N}} 

\newcommand{\cB}{{\cal B}}
\newcommand{\cC}{{\cal C}}
\newcommand{\cD}{{\cal D}}

\newcommand{\cL}{{\cal L}}
\newcommand{\cM}{{\cal M}}
\newcommand{\cN}{{\cal N}}
\newcommand{\cO}{{\cal O}}

\newcommand{\cV}{{\cal V}}
\newcommand{\cX}{{\cal X}}
\newcommand{\cY}{{\cal Y}}

\newcommand{\cZ}{{\cal Z}}

\newcommand{\mA}{{\bf A}}
\newcommand{\mB}{{\bf B}}

\newcommand{\mD}{{\bf D}}

\newcommand{\mI}{{\bf I}}

\newcommand{\mK}{{\bf K}}
\newcommand{\mL}{{\bf L}}

\newcommand{\mP}{{\bf P}}
\newcommand{\mQ}{{\bf Q}}

\newcommand{\mU}{{\bf U}}
\newcommand{\mV}{{\bf V}}
\newcommand{\mW}{{\bf W}}
\newcommand{\mX}{{\bf X}}

 \newcommand{\argmin}{\arg \min}

\DeclareMathOperator{\dom}{dom}         

\newcommand{\E}[1]{\mathbb{E}\left[#1\right] }

\newcommand{\EE}[2]{{\bf E}_{#1}\left[#2\right] } 
\providecommand{\Range}[1]{\mathbf{Range}\left( #1\right)}

\declaretheorem[within=section]{definition}
\declaretheorem[sibling=definition]{theorem}
\declaretheorem[sibling=definition]{proposition}
\declaretheorem[sibling=definition]{assumption}
\declaretheorem[sibling=definition]{corollary}
\declaretheorem[sibling=definition]{lemma}

\theoremstyle{remark}
\newtheorem{example}{Example} 
\newtheorem{remark}{Remark} 

\providecommand{\Range}[1]{\mathbf{Range}\left( #1\right)}

\providecommand{\E}[1]{\mathbf{E}\left[ #1\right]}

\newcommand{\sumin}{\sum_{i=1}^n}
\newcommand{\sumiM}{\sum_{i=1}^M}
\newcommand{\sumjn}{\sum_{j=1}^n}

\newcommand{\avein}{\frac{1}{n}\sum_{i=1}^n}

\newcommand{\RR}{\mathbb{R}}

\def\<#1,#2>{\left\langle #1,#2\right\rangle}

\newcolumntype{M}[1]{>{\centering\arraybackslash}m{#1}}
\newcolumntype{N}{@{}m{0pt}@{}}

\newcommand{\algname}[1]{#1} 

\def\dom{{\rm dom \,}}

\newcommand\ev[1]{\left \langle #1 \right \rangle}
\newcommand\ec[2][]{\ensuremath{\mathbb{E}_{#1} \left[#2\right]}}
\newcommand\ecn[2][]{\ec[#1]{\norm{#2}^2}}
\newcommand\br[1]{\left ( #1 \right )}
\newcommand\pbr[1]{\left \{ #1 \right \} }

\newcommand{\e}{\varepsilon}

\newcommand{\Om}{\Omega}

\newcommand{\mc}{\mathcal}

\newcommand{\norm}[1]{\left\lVert#1\right\rVert}
\newcommand{\sqn}[1]{{\left\lVert#1\right\rVert}^2}
\newcommand{\abs}[1]{\left\lvert#1\right\rvert}

\newcommand{\prox}{\mathrm{prox}}
\newcommand{\D}{\mathcal{D}}
\newcommand{\eqdef}{\overset{\text{def}}{=}}
\newcommand{\avemm}{\frac{1}{M}\sum_{m=1}^M}

\newcommand{\sigmaopt}{\sigma_{\mathrm{opt}}}
\newcommand{\sigmaf}{\sigma_{\mathrm{dif}}}

\usepackage{tcolorbox}
\usepackage{pifont}

\definecolor{mydarkgreen}{RGB}{39,130,67}
\definecolor{mygreen}{rgb}{0.0, 0.65, 0.31}
\newcommand{\green}{\color{mydarkgreen}}
\definecolor{mydarkred}{RGB}{192,47,25}
\newcommand{\red}{\color{mydarkred}}

\newcommand{\cmark}{\green\ding{51}}%
\newcommand{\xmark}{\red\ding{55}}%
\newcommand\floor[1]{\left \lfloor #1 \right \rfloor}

\newcommand{\prm}[1]{\pi_{#1}}
\newcommand{\emone}{\bf}
\newcommand{\emtwo}{\em}

\newcommand{\sigmaesc}{\sigma_{\ast}}
\newcommand{\sigmass}{\sigma_{\mathrm{Shuffle}}^2}
\newcommand{\ctO}{\tilde{\mathcal{O}}}
\newcommand\et[1]{{\mathbb{E}_{k} \left[#1\right]}}

\newcommand{\rrbox}[1]{\colorbox{myblue!30}{#1}}
\newcommand{\sobox}[1]{\colorbox{mygreen!30}{#1}}
\definecolor{myblue}{RGB}{55,126,184}
\newcommand{\sigmarr}{\sigma_{\mathrm{rad}}}
\DeclarePairedDelimiter{\ceil}{\lceil}{\rceil}

\newcommand{\Lave}{\overline{L}}
\def\xx{{\boldsymbol x}}
\def\yy{{\boldsymbol y}}

\newcommand{\ind}{\chi} 
\newcommand{\avejm}{\frac{1}{m}\sum_{j=1}^m}

\newcommand{\sumlm}{\sum_{l=1}^m}

\newcommand{\sumjm}{\sum_{j=1}^m}

\newcommand{\proxR}{\prox_{\gamma \psi}}
\newcommand{\proxj}{\prox_{\eta_j g_j}}

\newcommand{\ps}[1]{\langle #1 \rangle}
\newcommand{\cF}{{\mathscr F}}
\newcommand{\bP}{{\mathbb P}}
\newcommand{\bR}{{\mathbb R}}
\newcommand{\bE}{{\mathbb E}}
\newcommand{\Span}{\mathop{\mathrm{span}}\nolimits}

\newcommand{\n}[1]{\|#1 \|}
\newcommand{\lr}[1]{\langle #1\rangle}
\renewcommand{\th}{\theta}
\DeclareMathOperator{\clconv}{\overline{conv}}
\newcommand{\x}{\bar x}
\renewcommand{\a}{\alpha}
\renewcommand{\b}{\beta}
\renewcommand{\ker}{\mathbf{Ker}}

\newcommand{\sign}{\mathrm{sign}}

\renewcommand{\ast}{\star}

\newacronym{AdGD}{AdGD}{Adaptive Gradient Descent}
\newacronym{AdGD-accel}{AdGD-accel}{Adaptive Accelerated Gradient Descent}
\newacronym{ADMM}{ADMM}{Alternating Direction Method of Multipliers}
\newacronym{AdSGD}{AdSGD}{Adaptive Stochastic Gradient Descent}
\newacronym{BB}{BB}{Barzilai--Borwein}
\newacronym{DR}{DR}{Douglas--Rachford}
\newacronym{DYS}{DYS}{Davis--Yin Splitting}
\newacronym{ERM}{ERM}{Empirical Risk Minimization}
\newacronym{FedRR}{FedRR}{Federated Random Reshuffling}
\newacronym{FL}{FL}{Federated Learning}
\newacronym{GD}{GD}{Gradient Descent}
\newacronym{IG}{IG}{Incremental Gradient}
\newacronym{IAG}{IAG}{Incremental Average Gradient}
\newacronym{Lasso}{Lasso}{Least absolute shrinkage and selection operator}
\newacronym{LiCoSGD}{LiCoSGD}{Linearly Constrained Stochastic Gradient Descent}
\newacronym{MPI}{MPI}{Message Passing Interface}
\newacronym{PDDY}{PDDY}{Primal--Dual Davis--Yin}
\newacronym{PDHG}{PDHG}{Primal--Dual Hybrid Gradient}
\newacronym{PriLiCoSGD}{PriLiCoSGD}{Primal Linearly Constrained Stochastic Gradient Descent}
\newacronym{ProxRR}{ProxRR}{Proximal Random Reshuffling}
\newacronym{QSGD}{QSGD}{Quantized Stochastic Gradient Descent}
\newacronym{RR}{RR}{Random Reshuffling}
\newacronym{SAG}{SAG}{Stochastic Average Gradient}
\newacronym{SDCA}{SDCA}{Stochastic Dual Coordinate Ascent}
\newacronym{SDM}{SDM}{Stochastic Decoupling Method}
\newacronym{SGD}{SGD}{Stochastic Gradient Descent}
\newacronym{SVM}{SVM}{Support-Vector Machine}
\newacronym{SVRG}{SVRG}{Stochastic Variance Reduced Gradient}
\newacronym{SO}{SO}{Shuffle-Once}
\newacronym{VR}{VR}{Variance Reduction}

\usepackage[lmargin=1.5in, rmargin=1in, vmargin=1in, headsep=0.083334in]{geometry}
\newcommand{\mathsym}[1]{{}}
\newcommand{\unicode}[1]{{}}

\renewcommand\bibname{\centering BIBLIOGRAPHY}

\begin{document}


\vspace{2pt}
\thispagestyle{empty}
\addvspace{10mm}

\begin{center}

{\textbf{{\large On Seven Fundamental Optimization Challenges  \\ in Machine Learning }}}\vfill 
{Dissertation by}\\
{Konstantin Mishchenko}\vfill

{ In Partial Fulfillment of the Requirements}\\[12pt]
{ For the Degree of}\\[12pt]
{Doctor of Philosophy} \vfill
{King Abdullah University of Science and Technology }\\
{Thuwal, Kingdom of Saudi Arabia}
\vfill
{October 2021}

\end{center}

\newpage

\begin{center}

\end{center}

\begin{center}

{ \textbf{{\large EXAMINATION COMMITTEE PAGE}}}\\\vspace{1cm}

\end{center}
\noindent{The dissertation of Konstantin Mishchenko is approved by the examination committee}
\addcontentsline{toc}{chapter}{Examination Committee Page}

\vspace{4\baselineskip}

\noindent{Committee Chairperson: Peter Richt\'arik}\\
Committee Members: Lawrence Carin, Bernard Ghanem, Suvrit Sra, Wotao Yin \vfill

\newpage
\addcontentsline{toc}{chapter}{Copyright}
\vspace*{\fill}
\begin{center}
{ \copyright \; October 2021}\\
{Konstantin Mishchenko}\\
{All Rights Reserved}
\end{center}

\singlespacing

\begin{center}

\end{center}

\begin{center}
{{\bf\fontsize{14pt}{14.5pt}\selectfont \uppercase{ABSTRACT}}}
\end{center}

\addcontentsline{toc}{chapter}{Abstract}

\begin{center}
{{\fontsize{14pt}{14.5pt}\selectfont {On Seven Fundamental Optimization Challenges in Machine Learning\\
\phantom{On Some Fundamental Challenges in Optimization for Machine Learning}
 Konstantin Mishchenko}}}
\end{center}

Many recent successes of machine learning went hand in hand with advances in optimization. The exchange of ideas between these fields has worked both ways, with machine learning building on standard optimization procedures such as gradient descent, as well as with new directions in the optimization theory stemming from machine learning applications. In this thesis, we discuss new developments in optimization inspired by the needs and practice of machine learning, federated learning, and data science. In particular, we consider seven key challenges of mathematical optimization that are relevant to modern machine learning applications, and develop a solution to each.

Our first contribution is the resolution of a key open problem in Federated Learning: we establish the first theoretical guarantees for the famous Local SGD algorithm in the crucially important heterogeneous data regime. As the second challenge, we close the gap between the upper and lower bounds for the theory of two incremental algorithms known as Random Reshuffling (RR) and Shuffle-Once that are widely used in practice, and in fact set as the default data selection strategies for SGD in modern machine learning software. Our third contribution can be seen as a combination of our new theory for proximal RR and Local SGD yielding a new algorithm, which we call FedRR. Unlike Local SGD, FedRR is the first local first-order method that can provably beat gradient descent in communication complexity in the heterogeneous data regime. The fourth challenge is related to the class of adaptive methods. In particular, we present the first parameter-free stepsize rule for gradient descent that provably works for any locally smooth convex objective. The fifth challenge we resolve in the affirmative is the development of an algorithm for distributed optimization with quantized updates that preserves global linear convergence of gradient descent. Finally, in our sixth and seventh challenges, we develop new VR mechanisms applicable to the non-smooth setting based on proximal operators and matrix splitting.

In all cases, our theory is simpler, tighter and uses fewer assumptions than the prior literature. We accompany each chapter with numerical experiments to show the tightness of the proposed theoretical results.



\begin{center}

\end{center}

\begin{center}

{\bf\fontsize{14pt}{14.5pt}\selectfont \uppercase{Acknowledgements}}\\\vspace{1cm}
\end{center}

\addcontentsline{toc}{chapter}{Acknowledgements}
First and foremost, I would like to thank my supervisor Peter Richt\'{a}rik for his continous guidance throughout my PhD. The discussions and projects that I had with him helped me shape my research interests and navigate my efforts towards a combination of rigor, simplicity and mathematical beauty, which resulted in many papers that I can be proud of. Above all, with Peter Richt\'{a}rik's help, I quickly became able to work independently and collaborate with people from different countries and backgrounds, which lays a solid foundation for my future work.

I am grateful to my coauthors for making this thesis happen. I cannot imagine myself achieving the same research progress without the help of the people that I worked with.

I was lucky to have had many short and long research visits, and many eye-opening and inspiring discussions. I would like to particularly thank:\\
	 J{\'e}r{\^o}me Malick and Franck Iutzeler for hosting me at Universit\'{e} Grenoble Alpes when I was writing my MSc thesis in France. \\
	 Carola-Bibiane Sch\"{o}nlieb for hosting me at the University of Cambridge two times.\\
	 Matthias Ehrhardt who hosted me at the University of Bath. \\
	 Mert G\"{u}rb\"{u}zbalaban for hosting me at Rutgers University. \\
	 Aryan Mokhtari for hosting me at MIT. \\
	 Federico Vaggi for hosting me as an intern at Amazon. \\
	 Dmitriy Drusvyatskiy for letting me visit the University of Washington. \\
	 Lin Xiao for my short visit to Microsoft Research at Redmond. \\
	 Martin Jaggi for hosting me at EPFL. \\
	 Stephen Boyd for hosting me at Stanford University. \\
	 Alexander Gasnikov for hosting me at the Moscow Institute of Physics and Technology. \\
	 Courtney Paquette and Nicolas Le Roux for hosting my internship at Google Brain. I am also thankful to Fabian Pedregosa for his support before, during and after the internship.

I am very happy that my PhD allowed me to meet many brilliant people. I would not have succeeded in my projects without the help of Laurent Condat, Eduard Gorbunov, Filip Hanzely, Samuel Horv\'{a}th, Ahmed Khaled, Dmitry Kovalev, Yura Malitsky, Xun Qian, Adil Salim, Alibek Sailanbayev, Egor Shulgin, Martin Tak{\'a}{\v{c}}, and Bokun Wang. I was lucky to have had discussions with many more great minds, including Ahmet Alacaoglu, L{\'e}on Bottou, Aaron Defazio, Sai Praneeth Karimireddy, James Martens, Dmitrii Ostrovskii, Boris Polyak, Sebastian U.\ Stich, Junzi Zhang, and Nikita Zhivotovskiy. I also wish to say a big ``Thank you!''\ to Donald Goldfarb for his kind words about my work.

It was also a great pleasure to be part of our research group and I thank its past and current members for many interesting conversations: El Houcine Bergou, Konstantin Burlachenko, Aritra Dutta, Elnur Gasanov, Robert M.\ Gower, Slavom\'{i}r Hanzely, Jakub Kone\v{c}n\'{y}, Zhize Li, Nicolas Loizou, Grigory Malinovsky, Mher Safaryan, and Igor Sokolov.

Finally, I would like to thank Olivier Bousquet who recommended me to ``start with a PhD no matter what'' when I was in doubt about my future.


\addcontentsline{toc}{chapter}{Table of Contents}
\renewcommand{\contentsname}{\centerline{\textbf{{\large TABLE OF CONTENTS}}}}
\tableofcontents

\glsaddall
\printglossary[type=\acronymtype,style=long3col, title=\centerline{LIST OF ABBREVIATIONS}, toctitle=List of Abbreviations, nonumberlist=true]

\cleardoublepage
\addcontentsline{toc}{chapter}{\listfigurename} 
\renewcommand*\listfigurename{\centerline{LIST OF FIGURES}} 
\listoffigures

\cleardoublepage
\addcontentsline{toc}{chapter}{\listtablename}
\renewcommand*\listtablename{\centerline{LIST OF TABLES}} 
\listoftables
\hypersetup{linkcolor=.}


%

\singlespacing 
\chapter{Introduction}\label{chapter:intro}

One of the first questions asked by machine learning practitioners when trying to build a model for finding patterns in and gaining insights from real-world data pertains to the choice of an appropriate model class. This is typically an exercise in applying prior expert knowledge and experience. Subsequently, they need to choose the right model inside the selected model class, which is performed by a suitable training algorithm. 

One of the most universal approaches to building such algorithms is to pick an error function, and minimize it over the space of parameters characterizing the model class:
\begin{equation}\label{eq:biug9dgufd-98yf9ydf}\min_{x\in \R^d} f(x).\end{equation}
 That is, training is performed via optimization. As it turns out, there is an abundance of optimization algorithms applicable or even specifically designed for minimizing such error functions appearing in machine learning.  As explained in the very first sentences of the seminal book of Nesterov~\cite{Nesterov2013}, the explosion of theoretical research in mathematical optimization happened after it was observed in the mid eighties that the theoretical complexity of many optimization algorithms often captures their practical performance very well. Since then, theoretical complexity has by many been considered as a key efficiency measure of optimization algorithms, one that is equal and in some aspects superior to their performance in numerical experiments, which is the main reason why optimization theory is so highly regarded. For example, complexity theory leads to deep understanding, which acts as an important catalyst in the development of new and more efficient algorithms, and a guide in the selection of the most appropriate method for a given application. 

Motivated by this observation, in this thesis we set out to study the convergence rates of several optimization algorithms of key importance to machine learning in general, and distributed and federated learning in particular, in order to shed light on their superior empirical performance. Equally importantly, we also propose several new efficient algorithms, providing improved, tight and insightful theoretical convergence rates. 

\begin{quote}This thesis is motivated by seven key theoretical challenges and open problems in optimization for machine learning. 
\end{quote}

However, before formulating the challenges that we plan to address, let us first briefly outline the relevant optimization history, and provide the necessary background for the algorithms that we shall analyze in the remainder of this thesis.

The field of mathematical optimization has been around for a long time and has gone through several markedly distinct stages. Starting with the early works of Newton~\cite{newton1687} and Cauchy~\cite{cauchy1847methode}, the first optimization algorithms were used on small-dimensional problems that could be solved by hand. While the introduction of computers and new algorithms, such as the simplex method and Karmarkar's algorithm for linear programming, allowed the practitioners to solve larger problem instances, their sizes were still quite limited.

The modern state of the optimization field, however, goes hand in hand with the development of high-dimensional machine learning models, and relies on the common availability of extremely fast computers. This has led to a change of paradigm, and many new and exciting optimization algorithms were proposed in the past two decades. For example, larger models (i.e., models described by a larger number $d$ of parameters) required finer regularization procedures, which served as the inspiration for the development of many new proximal algorithms (Parikh and Boyd \cite{parikh2014proximal}). Further, the emergence of new parallel processing hardware and the coupling of computers into clusters paved the way for new parallel and distributed optimization algorithms (Bertsekas and Tsitsiklis \cite{bertsekas1989parallel}), which brought new theoretical challenges, most notably the challenge of taming the high communication cost of distributed methods.

Perhaps the most striking characteristic of the modern era of optimization lies in the use of huge-scale datasets to formulate optimization problems. Such problems are characterized by the finite-sum structure of the objective,
\begin{equation}\label{eq:sum-09uf0d9u0d}f(x) = \frac{1}{n}\sum_{i=1}^n f_i(x),\end{equation}
where $n$ is typically the number of training data points.  As both the dimension $d$ and the number of data points $n$ grow beyond  millions or even billions, standard optimization algorithms, such as Gradient Descent (GD) and Newton's method, become either inapplicable or inefficient. Indeed, it is often impossible to make more than a few hundred ``passes''\footnote{One pass typically refers to work equivalent to the  evaluation of the $n$ gradients $\nabla f_1, \dots, \nabla f_n$.} over large datasets, while Gradient Descent  may require orders of magnitude more passes to converge. Newton's method, on the other hand, does not scale well with the dimension, and even quasi-Newton methods, such as BFGS (Broyden \cite{broyden1967quasi}; Fletcher \cite{fletcher1970new}; Goldfarb \cite{goldfarb1970family}; Shanno \cite{shanno1970conditioning}), do not break this limit due to memory issues.

To make iterations cheap, a popular approach to solving problem \eqref{eq:biug9dgufd-98yf9ydf} is to use incremental updates that do not require processing all data in each iteration (Nedi{\'{c}} and Bertsekas \cite{Nedic2001}). The most famous example of such an algorithm is Stochastic Gradient Descent (\algname{SGD}) (Robbins and Monro \cite{robbins1951stochastic}). This method uses a noisy estimate of the gradient and runs much faster than its deterministic counterpart, Gradient Descent (Bottou \cite{bottou2012stochastic}). Moreover, Nemirovski and Yudin~\cite{nemirovsky1983problem} established that under a certain noise condition, SGD is provably optimal for the minimization of strongly convex functions. To provide the reader with early intuition, let us write the update of a generic SGD method explicitly:
\begin{align*}
	x^{k+1} = x^k - \gamma_k \nabla f(x^k; \xi^k).
\end{align*}
Here $x^{k+1}, x^k \in\mathbb{R}^d$ are the new and the current iterates, $d$ is the problem dimension, $\gamma_k>0$ is a stepsize, and $\nabla f(x^k; \xi^k)$ is a stochastic gradient, whose expectation over random sample $\xi^k$ is equal to the gradient $\nabla f(x^k)$ of function $f$ that we aim to minimize. The popularity of SGD in practice inspires us to study theoretically its extensions as well as to look at algorithms that may outperform SGD in specific scenarios.

The optimality of SGD has its limits too. For instance, for objectives that are defined as sums of a finite number of terms, as in \eqref{eq:sum-09uf0d9u0d}, SGD can take more time to converge than Gradient Descent if a low-error solution is required. The cheap iterations of SGD eventually become less useful than the expensive gradient-descent iteration due to the variance inherent in stochastic approximation of the gradients. Achieving the best of both worlds was an elusive goal until the breakthrough work by Le Roux et al.~\cite{sag} proposed the SAG algorithm that attained the iteration convergence rate of Gradient Descent and the iteration cost of SGD. The paper of Le Roux et al.~\cite{sag} received the Lagrange Prize in Continuous Optimization for this achievement and its impact on the field, awarded in 2018, just six years after it was published. This discovery subsequently sparked a series of works on the topic and many other practical variants appeared thereupon, most cited of which is the SVRG paper of Johnson and Zhang~\cite{SVRG}. At the same time, it continued to be challenging to combine variance reduction with many other techniques that matter in structured and statistics-based optimization, most importantly with stochastic proximal operators and matrix splitting. 

The goal of this thesis is to address some of the key challenges that are relevant to the practice of machine learning. Some of these challenges became apparent only recently while others have been around, unsolved, for a long time. Of special interest to us will be algorithms based on stochastic updates, including SGD, and variance-reduced methods. Having provided a quick and broad overview of the field, we shall now elaborate in more detail on the challenges considered in this thesis. For the reader interested in a quick summary, we provide \Cref{tab:intro_summary}.

\section{Challenge 1: Convergence of Local SGD for Federated Learning in the Heterogeneous Data Regime}

\paragraph{Motivation.}
As machine learning models consume a lot of data during training, preservation of privacy becomes an important goal. How can we train machine learning models on private data without directly sharing it? The current way to make this possible is to equip data owners with the task of training a model directly where the data are stored, and  ask them to periodically communicate model updates with an orchestrating server. This approach to training models is called federated learning (Kone\v{c}n\'{y} et al., \cite{FEDLEARN}; McMahan et al., \cite{McMahan17}).

Privacy is not the only motivation for federated learning. Soon after edge devices, such as mobile phones, became capable of computation, it became apparent that the communication costs of training on edge devices far exceeds that of computation. It is much more efficient to use the devices to perform training instead of transferring all the data for centralized processing. Unfortunately, this makes SGD inapplicable in its standard form because it requires a synchronization after every gradient step. A simple remedy that works surprisingly well is to run SGD locally on each device and only communicate the final iterate of the local run. This idea dates back to the nineties (Mangasarian, \cite{Mangasarian95}), yet only recently we started to understand why it works so well. Moreover, it became really useful mostly because federated learning posed unprecedented constraints on the algorithms that are supposed to perform it, while the popularity of the area exploded in the recent years due to quickly arising applications.

In October 2020, Forbes\footnote{\href{https://www.forbes.com/sites/robtoews/2020/10/12/the-next-generation-of-artificial-intelligence/}{https://www.forbes.com/sites/robtoews/2020/10/12/the-next-generation-of-artificial-intelligence/}} placed federated learning as the second of three emerging areas of AI. This placement followed the extremely rapid growth of the research on federated learning and how it can be used for various applications. As sometimes argued (Kairouz et al., \cite{kairouz2019advances}; McMahan et al., \cite{McMahan17}), federated learning has the potential to be the primary way every machine learning model for mobile devices is trained. These promises have attracted a lot of attention to the area: the pioneering works of McMahan et al.~\cite{McMahan17} and Kone\v{c}n\'{y} et al.~\cite{FEDLEARN} have been cited 3,400+ times and 1,600+ times in just 5 years, respectively,  and a recent survey of Kairouz et al.~\cite{kairouz2019advances} has been cited 900+ times in less than 2 years. As promising as it is, federated learning is also extremely hard as a research problem. This area overlaps with many mathematical and engineering fields, but its optimization side is particularly nontrivial. Indeed, in contrast to the standard centralized learning scenario, federated learning was formulated as a problem with data stored in private on unreliable devices that have limited ability to communicate (McMahan et al., \cite{McMahan17}).

At the heart of federated learning lies the Federated Averaging algorithm---a variant of Local SGD that can find a solution with little communication. While the theory of local methods started in the nineties with the work of Mangasarian~\cite{Mangasarian95}, it did not get full recognition until federated learning turned out to be the right application. Once this became clear, the tables turned and the theory got quickly behind the practice. Without understanding the convergence rates of Local SGD, it is particularly hard for any new method to be shown superior to it. And the demand for a solid theory is particularly high due to the interest in deploying federated learning in applications such as speech recognition, medical research and mobile applications (Kairouz et al., \cite{kairouz2019advances}). To at least address the present gap in theory versus practice, we are motivated to study the precise rates of Local SGD, which is the first challenge that we will address in this thesis.

\paragraph{Contributions.}
In \Cref{chapter:local_sgd}, we provide a new analysis of \algname{Local SGD}, removing unnecessary assumptions and elaborating on the difference between two data regimes: identical and heterogeneous. In the heterogeneous case, which is of key importance in federated learning, our analysis is the first one to show that Local SGD works even when the gradients do not satisfy any variant of bounded dissimilarity. This was considered an important open problem in federated learning. However, even in the case of identically sampled data, we improve the existing theory and provide values for the optimal stepsize and optimal number of local iterations. Our bounds are based on a new notion of variance that is specific to \algname{Local SGD} methods with heterogeneous data. The tightness of our results is guaranteed by recovering known statements when we specialize them to $H=1$, where $H$ is the number of local steps. Our empirical evidence further validates the severe impact of data  heterogeneity on the performance of \algname{Local SGD}.

Before our work, the main focus of research in optimization for federated learning was on the data sampled from identical (Stich, \cite{Stich2018}) or almost identical (Jiang and Agrawal, \cite{Jiang18}) distributions. Such results can guarantee an improvement when the goal is to train a model faster by parallelizing SGD over multiple devices. In federated learning, however, the data come from various sources; for instance, from mobile users living in different regions or even countries. In this case, assuming similarity between gradients is rather limiting and even unrealistic. Our results, in contrast, do not require any global bound on gradient dissimilarity, and depend  on the  gradient norms at the solution only.

\paragraph{Paper.}
The chapter is based on the paper:
\begin{quote}
	\cite{khaled2020tighter} \bibentry{khaled2020tighter}.
\end{quote}

\section{Challenge 2: Convergence of Random Reshuffling}

\paragraph{Motivation.}
Theoretical understanding of {\em incremental gradient methods} constitutes a rather long-standing challenge. These methods are tailored to optimization problems \eqref{eq:biug9dgufd-98yf9ydf} with $f$ being of the ``finite-sum'' form \eqref{eq:sum-09uf0d9u0d},  all relying on the  iteration
\[x^{k+1} = x^k - \gamma_k \nabla f_{i^k}(x^k).\]
While in SGD the index $i^k$ is picked independently in each iteration from some fixed probability law, i.e., one performs  data {\em sampling with replacement}, the world of incremental methods is much richer as it allows for virtually arbitrary rules for the selection of the next datapoint $i^k$ to learn from.

Of particular interest to us in this thesis are incremental gradient methods of the {\em random permutation/shuffling/reshuffling} variety. This is because despite the common recognition of SGD as the workhorse behind many successes in machine learning and deep learning, in practice, \algname{SGD} is virtually always superseded by the incremental algorithm known as {\em Random Reshuffling} (\algname{RR}), which is based on a {\em sampling without replacement} approach to the selection of the data indices $\{i^k\}_k$.  That is, in RR, a random permutation/reshuffling of the data points $\{1,2,\dots,n\}$ is selected at the beginning of each epoch, and this order subsequently dictates the selection of the indices $\{i^k\}_k$. It is worthwhile to remark that RR data selection strategies are the default in modern deep learning software.

Despite the practical superiority of RR over SGD, virtually all  theoretical effort is directed towards the understanding of SGD type methods, and a proper understanding of methods based on random shuffling/reshuffling of data in general, and RR in particular, remains elusive. One of the key reasons for the disproportionate focus on SGD is that it is much easier to explain theoretically.  For example, one can view it as an instance of the long-studied and well-understood {\em Stochastic Approximation} of Robbins and Monro~\cite{robbins1951stochastic}, and immediately obtain a convergence rate. SGD has also been combined with iterate averaging by Polyak and Juditsky~\cite{polyak1992acceleration}, which further refines the convergence rates and keeps the theory still very simple. Nevertheless, the theory for Random Reshuffling has been developing at a much slower pace and its theoretical superiority to \algname{SGD} was established only recently.

The history of incremental and shuffling-based methods is quite extensive. The algorithms gained their popularity in the eighties under the name of online backpropagation algorithm (G\"{u}rb\"{u}zbalaban et al., \cite{Gurbuzbalaban2019IG}), and the theoretical development of these methods started more than thirty years ago, see, for instance, the survey by Bertsekas~\cite{Bertsekas2011} for more details. Nevertheless, until now, even for the earliest versions of the most basic incremental algorithms, the convergence has not been fully understood. As we shall see later, it was long unknown even where exactly the intermediate iterates of those algorithms converge to. In particular, these algorithms exhibit strong oscillations, which turned out to be quite hard to explain. Explaining them, on the other hand, leads to a significantly improved analysis that we present in this thesis.

We also highlight that the generality of shuffling-based methods allows for them to be applied in federated learning. Just as Random Reshuffling outperforms SGD, a local variant of RR should at least be on par with Local SGD. However, no such algorithm is known in the literature and tackling it without understanding RR itself would not give us a clear comparison. 
Therefore, we first approach RR as a challenge itself, and only later proceed to obtain a local variant of RR, which we cal FedRR. Equipped with a tighter theory for Local SGD that we obtain in this thesis too, we will be able to give a comprehensive comparison of FedRR and Local SGD.

\paragraph{Contributions.} 
In \Cref{chapter:rr}, we improve upon existing theory of Random Reshuffling in several ways and provide guarantees that match existing lower bounds. In prior literature for strongly convex and smooth functions, \algname{RR} was shown to converge faster than \algname{SGD} if 1) the stepsize is small, 2) the gradients are bounded, and 3) the number of epochs is large. However, large stepsizes are crucial for fast convergence at initialization, the gradients cannot be bounded for strongly convex functions, and small-epoch convergence is of high value if the time budget is limited. Thus, we provide a theory without these 3 assumptions, and, in addition, improve the dependence on the condition number from $\kappa^2$ to $\kappa$ (respectively from $\kappa$ to $\sqrt{\kappa}$) . Furthermore, we show that the power of \algname{RR} comes from a fundamentally different type of variance, which is based on the notion of Bregman divergence. We argue through theory and experiments that the new variance type gives an additional justification of the superior performance of \algname{RR}. To go beyond strong convexity, we present several results for non-strongly convex and non-convex objectives. We show that in all cases, our theory improves upon existing literature. 

Finally, we prove fast convergence of the \algname{Shuffle-Once} (\algname{SO}) algorithm, which shuffles the data only once, at the beginning of the optimization process. Our theory for strongly convex objectives matches the known lower bounds for both \algname{RR} and \algname{SO} and substantiates the common practical heuristic of shuffling once or only a few times. As a byproduct of our analysis, we also get new results for the \algname{Incremental Gradient} algorithm (\algname{IG}), which does not shuffle the data at all.

\paragraph{Paper.} 
The chapter is based on the paper:
\begin{quote}
	\cite{MKR2020rr} \bibentry{MKR2020rr}.
\end{quote}

\section{Challenge 3: Going Beyond Local SGD in Federated Learning}

\paragraph{Motivation.} 
As outlined above, incremental algorithms have been widely used in practice, and we managed to close some of the gaps in the theory of  these methods in  \Cref{chapter:rr}. The practical and theoretical success of RR naturally raises the question whether it is possible to successfully employ sampling without replacement in federated learning as well. 

Since in \Cref{chapter:local_sgd} we studied Local SGD---a method of key importance to federated learning---and provided the first guarantees for it in the heterogeneous case, it is natural to ask whether the theoretical tools developed therein can be combined with the new tools developed in our study of RR in order to  improve upon the state-of-the-art Local SGD rate. This is highly desired, as even the best known communication complexity results for Local SGD  do not, in general, improve upon the communication complexity of simple Gradient Descent, which casts a deep shadow onto the current state of theory in federated learning.

\paragraph{Contributions.} 

We answer the above challenge in the affirmative. Our new method, FedRR, is the first local-based gradient method that beats Gradient Descent (and hence also Local SGD) in communication complexity.

In \Cref{chapter:proxrr}, we propose two new algorithms: \algname{Proximal Random Reshuffing} (\algname{ProxRR}) and \algname{Federated Random Reshuffing} (\algname{FedRR}). The first algorithm, \algname{ProxRR}, solves {\em composite} convex finite-sum minimization problems. These are problems in which the objective is the sum of the average of $n$ smooth objectives as in \eqref{eq:sum-09uf0d9u0d}, and of a (potentially non-smooth) convex regularizer. This problem is of independent interest, as ProxRR is the first RR-based method that can provably solve proximal problems. However, the development of ProxRR should also be seen as an intermediary step towards obtaining our second algorithm, FedRR. Indeed, we obtain the second algorithm, \algname{FedRR}, as a special case of \algname{ProxRR} applied to a carefully designed reformulation of  the distributed problem appearing in federated learning, allowing for both homogeneous and  heterogeneous data. 

We study the convergence properties of both methods with constant and decreasing stepsizes, and show that they have considerable advantages over Proximal and Local \algname{SGD}. In particular, our methods have superior complexities, and \algname{ProxRR} evaluates the proximal operator once per epoch only. When the proximal operator is expensive to compute, this small difference makes \algname{ProxRR} up to $n$ times faster than algorithms that evaluate the proximal operator in every iteration. We give examples of practical optimization tasks where the proximal operator is difficult to compute and \algname{ProxRR} has a clear advantage. When specializing to the federated learning setting, our FedRR algorithm needs to communicate only after every local pass over local data is done. In contrast to the theory of Local SGD, which requires dividing the stepsize by $n$ when $n$ local steps are performed, our theory allows for stepsizes that do not depend on $n$ at all.

We note that our results considerably improve upon the complexity of Local SGD. Since incremental algorithms use the finite-sum structure of the objective, they are not subject to the lower bounds established for Local SGD by Woodworth et al.~\cite{woodworth2020minibatch}. This allows our algorithm FedRR to beat Local SGD after a certain number of iterations, regardless of how heterogeneous the data are.

\paragraph{Paper.} 
The chapter is based on the paper:
\begin{quote}
	\cite{mishchenko2021proximal} \bibentry{mishchenko2021proximal}.
\end{quote}

\section{Challenge 4: The First Adaptive Stepsize Rule for Gradient Descent that Provably Works}

\paragraph{Motivation.} 
To run Gradient Descent, one needs to use a stepsize that depends on problem-specific constants such as the Lipschitz parameter of the objective's gradient (Nesterov, \cite{Nesterov2013}) or the objective value at the optimum (Polyak, \cite{polyak1963gradient}). How can we circumvent this requirement if we do not know the required constants in advance? This fundamental question has been important to the field of optimization for many decades despite the significant changes in applications. 

Loosely speaking, algorithms that are able to run without explicit knowledge of some problem-specific constants and that can perform comparably to a similar algorithm that has explicit access to these constants are called {\em adaptive}~\cite{lei2020adaptivity}. As we mentioned earlier, tight theoretical complexity of an algorithm is quite often closely reflected in its practical performance. Of all counterexamples to this rule, one that stands out the most is the class of adaptive methods. In fact, the methods that are provably {\em not convergent} are often among the most popular ones. For instance, Adam (Kingma and Ba, \cite{adam}), for which Reddi et al., \cite{adam2} gave convex counterexamples, is much more frequently used in practice than its theory-based-counterpart Adagrad (Duchi et al., \cite{duchi2011adaptive}). Indeed, as of now, the work that proposed Adam has 86,000+ citations, while the paper on Adagrad has only 9,100+ citations. Similarly, in a recent work, Burdakov et al.~\cite{burdakov2019stabilized} identified a simple convex counterexample for the widely adopted Barzilai--Borwein (BB) method (Barzilai and Borwein, \cite{barzilai1988two}), which has been praised for its numerical performance (Wright, \cite{wright2010optimization}) and has accumulated 2,600+ citations.

The history of adaptive methods started long ago before Adagrad and Adam were introduced. The early approaches to adaptive parameter estimation include the seminal works of Armijo~\cite{armijo1966} and Polyak~\cite{polyak1963gradient} that proposed what later became to be known as Armijo line-search and Polyak's stepsize. The Barzilai--Borwein method was proposed two decades later (Barzilai and Borwein, \cite{barzilai1988two}) and soon after gained popularity for its efficiency in practice. At first, it even appeared that BB might be provably convergent as a few years after its introduction it was shown to work for quadratic problems by Raydan~\cite{raydan1993barzilai}. Despite the fact that line search and Polyak's rule appeared more than half a century ago and BB has never had any theory for non-quadratic functions, their ability to estimate objective parameters still attracts considerable attention (Tan et al., \cite{tan2016barzilai}; Hazan and Kakade, \cite{hazan2019revisiting}; Vaswani et al., \cite{vaswani2019painless}; Loizou et al., \cite{loizou2021stochastic}). This ability is, in fact, even more important for machine learning applications, because the corresponding problems rarely admit closed-form expressions for problem constants such as gradient Lipschitzness. In classical problems, this constant can be often computed exactly. For instance, it is the maximum singular value of the data matrix in least-squares regression. In contrast, objectives such as those appearing in neural network training may have unknown or infinite global smoothness constants (Zhang et al., \cite{zhang2019gradient}), and hence only a local estimation can potentially work in practice.

The empirical success of adaptive methods gives ample motivation to start developing theory for adaptive methods. Surprisingly, however, there exists no known closed-form stepsize for Gradient Descent that is completely parameter-free and would provably converge. Methods such as normalized gradient (Shor, \cite{shor1962application}) and Adagrad (Duchi et al., \cite{duchi2011adaptive}) need constants related to the distance from the solution, while Polyak's stepsize rule requires the knowledge of the optimal function value. Further, while line search procedures are parameter-free, they are not given in a closed-form, require subroutines to be run, and work only for globally smooth objectives. 

\paragraph{Contributions.} 
Therefore, one of our goals in this thesis is to provide the first stepsize for Gradient Descent that requires access to no information beyond the gradients themselves.  In particular, in Chapter~\ref{chapter:adaptive}, we shall present a proof that Gradient Descent with our stepsize rule can provably minimize any locally smooth convex function.

More specifically, we prove that two simple rules are sufficient to automate \algname{Gradient Descent}: 1)~do not increase the stepsize too fast, and 2)~do not overstep the local curvature. Namely, our theory guarantees than for a stepsize sequence $\gamma_1,\dotsc, \gamma_k,\dotsc$ to work, the only two requirements are
\begin{eqnarray*}
  \gamma_k^2 & \leq & \left(1+\frac{\gamma_{k-1}}{\gamma_{k-2}} \right)\gamma_{k-1}^2,\\ \gamma_k & \leq & \frac{\n{x^{k}-x^{k-1}}}{2\n{\nabla
        f(x^{k})-\nabla f(x^{k-1})}},
\end{eqnarray*}
Our method does not need any line search, and works without knowing the functional values or any other information about the objective except for the gradients. By choosing $\gamma_k$ per the rules above, we obtain a method adaptive to the local geometry, with convergence guarantees depending  on the smoothness in a neighborhood of a solution only. Given that the problem is convex, our method converges even if the global smoothness constant is infinity. As an illustration, it can minimize an arbitrary twice continuously differentiable convex function. We examine its performance on a range of convex and non-convex problems, including logistic regression, matrix factorization, and neural network training.

To the best of our knowledge, our stepsize rule for Gradient Descent is the only one that provably gives convergence for non-quadratic functions. Many other attempts to obtain adaptive stepsize rule work either  only for quadratics  (Raydan, \cite{raydan1993barzilai}), self-concordant functions with known self-concordance parameter (Gao and Goldfarb, \cite{gao2019quasi}), or require knowledge of the problem conditioning (Tan et al., \cite{tan2016barzilai}). The immense interest in adaptive methods, which can be seen from the number of citations to the aforementioned papers, points to only one explanation for the apparent scarcity of the results on the topic: the technical difficulty of obtaining such results.

\paragraph{Paper.} 
The chapter is based on the paper:
\begin{quote}
	\cite{malitsky2019adaptive} \bibentry{malitsky2019adaptive}.
\end{quote}

\section{Challenge 5: Achieving Fast Rates in Distributed Optimization with Quantization}

\paragraph{Motivation.}
While Local SGD and similar algorithms have been very successful in addressing the communication challenge in federated learning, they also have limitations (Woodworth et al., \cite{woodworth2020minibatch}). An alternative way to make communication cheaper is to apply lossy compression to the communicated vectors (Alistarh et al., \cite{alistarh2017qsgd}; Wen et al., \cite{wen2017terngrad}). The choice of compression plays an important role in how the communication is performed, and the idea of sending only the sign of the update vector has been particularly popular (Bernstein et al., \cite{bernstein2018signsgd}). Another popular technique is to send  the coordinates with the largest magnitudes only (Stich et al., \cite{stich2018sparsified}), which requires sorting the uncompressed vector and may be thereby a bit slower.

To quantify the overall benefit of relying on compressed communication, we need to consider two factors: the per-iteration savings coming from compression, and the increase in the iteration complexity resulting from compression. Since we send less information, it is only reasonable to expect that the iteration complexity would get worse. Therefore, we ask the following natural question: if the compressed update requires $\omega\ge 1$ times fewer bytes to send, how does the iteration complexity depend on $\omega$? As it turns out, the main drawback of using lossy compression is the requirement to use $\omega$-times smaller stepsizes, and the resulting slowdown in the convergence rates of existing methods is proportional to $\omega$, too. Thus, in general, there might be little or no benefit from applying compression. The error-feedback technique proposed by Stich et al.~\cite{stich2018sparsified} and refined by Stich and Karimireddy~\cite{stich19errorfeedback} allows to partially improve the rate when the complexity is driven by the noise of stochastic gradients but, unfortunately, it is still incapable of fixing the stepsize requirement. Our next goal in this thesis, therefore, is to find a new algorithm that does not suffer from requiring small stepsizes when the noise is mild.

Several other methods based on the compression (e.g., sparsification and/or quantization) of the updates were recently proposed, including QSGD (Alistarh et al., \cite{alistarh2017qsgd}), TernGrad (Wen et al., \cite{wen2017terngrad}), SignSGD (Bernstein et al., \cite{bernstein2018signsgd}), and DQGD (Khirirat et al., \cite{khirirat2018distributed}). However, all of these methods suffer from severe issues, such as the inability to converge to the true optimum in the batch mode, inability to work with a non-smooth regularizer, and slow convergence rates. 

\paragraph{Contributions.} 
We propose a new distributed learning method---DIANA---which resolves these issues via a new algorithmic tool: {\em compression of gradient differences}. DIANA is the first variance-reduction mechanism for distributed training which can progressively reduce the variance introduced by gradient compression. In other words, DIANA is to DQGD what SVRG is to SGD.

We perform a theoretical analysis in the strongly convex and non-convex settings and show that our rates are superior to the existing rates. Moreover, our analysis of block-quantization and differences between $\ell_2$ and $\ell_{\infty}$ quantization closes one of the gaps in theory and practice. Finally, by applying our analysis technique to TernGrad, we establish the first convergence rate for this method.

The idea of difference quantization proposed in our work has proved to be very helpful and the results of Chapter~\ref{chapter:diana} were extended in a number of works. In a follow-up paper (Horv{\'a}th et al., \cite{horvath2019stochastic}), we generalized it to arbitrary unbiased compressors, and combined it with a secondary variance reduction mechanism which allows to compress stochastic gradients without suffering rate deterioration. 

Our idea  was extended by many in various other ways, in particular, to server-side compression by Liu et al.~\cite{liu2020double}, device sampling by Philippenko and Dieuleveut~\cite{philippenko2020bidirectional}, acceleration by Li et al.~\cite{li2020acceleration}, second-order methods for generalized linear models by Islamov et al.~\cite{NL2021}, second-order methods for federated learning by Safaryan et al.~\cite{FedNL}, acceleration via matrix smoothness by Safaryan et al.~\cite{Safaryan+2021}, integer compression for SwitchML by Mishchenko et al.~\cite{mishchenko2021intsgd}, and to biased compression by Gorbunov et al.~\cite{gorbunov2020linearly}. Even more variants as well as a unification for their analysis were obtained by Gorbunov et al.~\cite{gorbunov2020unified}.

\paragraph{Paper.} 
The chapter is based on the paper:
\begin{quote}
	\cite{mishchenko2019distributed} \bibentry{mishchenko2019distributed}.
\end{quote}

\section{Challenge 6: Developing Variance Reduction for Proximal Operators}

\paragraph{Motivation.} The continued interest in algorithms that process the data efficiently resulted in a steady development of incremental and stochastic algorithms. Nevertheless, these algorithms are not always optimal. A well-known limitation of \algname{SGD} and \algname{RR} is that they might escape from a solution if the stochastic gradients do not converge to zero. In 2005, the Incremental Average Gradient (\algname{IAG}) method, proposed by Blatt et al.~\cite{blatt2007convergent},  managed to bypass this limitation in practice. However, the corresponding convergence rates were not promising. A true revolution started when the stochastic counterpart of IAG, known as Stochastic Average Gradient (SAG), was shown by Le Roux et al.~\cite{sag} to converge at a much faster rate than IAG, SGD and other similar incremental algorithms. SAG and other algorithms that have cheap iterations but can achieve the rate of Gradient Descent became known as {\em variance-reduced methods}, and include SDCA (Shalev-Shwartz and Zhang, \cite{SDCA}), SAGA  (Defazio et al., \cite{defazio2014saga}), SVRG (Johnson and Zhang, \cite{SVRG}), S2GD (Kone\v{c}n\'{y} and Richt\'{a}rik, \cite{S2GD}), MISO (Mairal et al., \cite{MISO}; Qian et al., \cite{qian2019miso}), QUARTZ (Qu et al., \cite{QUARTZ}), and JacSketch (Gower et al., \cite{JacSketch}).  For a recent review of variance reduced methods for machine learning, we refer the reader to Gower et al.~\cite{VR-Review2020}. 

Many further extensions of variance reduction have been proposed. For instance, proximal (Xiao and Zhang, \cite{xiao2014proximal}), mini-batch  (Kone\v{c}n\'{y} et al., \cite{mS2GD}),  accelerated (Allen-Zhu, \cite{allen2017katyusha}; Snang et al., \cite{Shang2018}), and loopless SVRG (Kovalev et al., \cite{kovalev2019don}; Qian et al., \cite{L-SVRG-AS}). The focus of these works, however, was directed towards improving the use of stochastic gradients. 

At the same time, many convex problems come with complicated constraints and regularizers, such as Group Lasso, that require expensive computation of proximal operators. For such problems, the speed advantage coming from cheap gradient computation fades away once we take into account the time required to tackle the regularizers. To make variance reduction useful in these settings, one has to find a way to make the computation of the proximal operator inexpensive, too. Notwithstanding the efforts to tackle this challenge, for instance by Ryu and Yin~\cite{ryu2017proximal}, Defazio~\cite{defazio2016simple}, and Pedregosa et al.~\cite{pedregosa2019proximal}, it has been unresolved in its general form.

\paragraph{Contributions.}  Motivated by these needs, we consider the problem of minimizing the sum of three convex functions: i) a smooth function $f$ in the form of an expectation or a finite average, ii) a non-smooth function $g$ in the form of a finite average of proximable functions $g_j$, and iii) a proximable regularizer $\psi$. We design a variance reduced method which is able to progressively learn the proximal operator of $g$ via the computation of the proximal operator of a single randomly selected function $g_j$ in each iteration only. Our method can provably and efficiently accommodate many strategies for the estimation of the gradient of $f$, including via standard and variance-reduced stochastic estimation, effectively decoupling the smooth part of the problem from the non-smooth part. We prove a number of iteration complexity results, including a general $\cO(\frac{1}{K})$ rate, $\cO(\frac{1}{K^2})$ rate in the case of strongly convex smooth $f$, and several linear rates in special cases, including accelerated linear rate. For example, our method achieves a linear rate for the problem of minimizing a strongly convex function $f$ subject to linear constraints under no assumption on the constraints beyond consistency. When combined with \algname{SGD} or \algname{SAGA} estimators for the gradient of $f$, this leads to a very efficient method for empirical risk minimization. 

Our method generalizes several existing algorithms, including forward-backward splitting, Douglas--Rachford splitting, \algname{Proximal SGD}, \algname{Proximal SAGA}, \algname{SDCA}, \algname{Randomized Kaczmarz} and \algname{Point--SAGA}. However, our method leads to new methods in special cases; for instance, we obtain the first randomized variant of the \algname{Dykstra's method} for projection onto the intersection of closed convex sets. The unified analysis proposed in our work might be of interest on its own as is gives a simple way to derive multiple methods at once.

\paragraph{Paper.}
The chapter is based on the paper:
\begin{quote}
	\cite{mishchenko2019stochastic} \bibentry{mishchenko2019stochastic}.
\end{quote}

\section{Challenge 7: Designing Variance-Reduced Algorithms with Splitting}

\paragraph{Motivation.}
Regularized objectives, such as PC-Lasso (Tay et al., \cite{tay2018principal}), sometimes regularize a linear transformation of the parameter vector instead of regularizing the parameters themselves. If the regularizer is a non-smooth function, it is recommended to use algorithms based on matrix splitting and proximal operators (Davis and Yin, \cite{davis2017three}). Such algorithms provide efficient update rules for the regularizer, but, unfortunately, they ignore the ample cost of computing gradients of the other part of the objective function. Can we find a way to design an algorithm that combines the best of both worlds: use variance-reduced gradient updates and split matrix multiplication from proximal operators?

\paragraph{Contributions.}
To answer this question, in Chapter~\ref{chapter:pddy}, we consider the task of minimizing the sum of three convex functions, where the first one $f$ is smooth, the second one is non-smooth and proximable, and the third one is the composition of a non-smooth proximable function with a linear operator $\mL$. This template problem has many applications in machine learning and signal processing. 

We propose a new primal--dual algorithm called \algname{PDDY} to solve such problem. \algname{PDDY} can be seen as an instance of \algname{Davis--Yin Splitting} involving operators which are monotone under a new metric depending on $\mL$. This representation of \algname{PDDY} eases the non-asymptotic analysis of \algname{PDDY}: it allows us to prove its sublinear convergence (respectively linear convergence if strong convexity is involved). 

Moreover, our proof technique easily extends to the case where a variance reduced stochastic gradient of $f$ is used instead of the full gradient. 
Besides, we obtain as a special case a linearly converging algorithm for the minimization of a strongly convex function $f$ under linear constraints $\mL x = b$. This algorithm can be applied to decentralized optimization problems and competes with other approaches specifically designed for decentralized optimization. 

Finally, we show that three other primal--dual algorithms (the two \algname{Condat--V\~u} algorithms and the \algname{PD3O} algorithm) can be seen as \algname{Davis--Yin Splitting} under a metric depending on $\mL$. Such representation was not known for the \algname{Condat--V\~u} algorithms. We show again that this representation eases the non-asymptotic analysis of \algname{PD3O} in the case where a variance reduced stochastic gradient is used. Our theory covers several settings that are not tackled by any existing algorithm; we illustrate their importance with real-world applications and we show the efficiency of our algorithms by numerical experiments.

\paragraph{Paper.}
The chapter is based on the paper:
\begin{quote}
	\cite{salim2020dualize} \bibentry{salim2020dualize}.
\end{quote}

\begin{table}[t]
    \caption{A summary of the results obtained in this thesis.}
    \label{tab:intro_summary}
   \centering 
    \begin{tabular}{ccc}
        \toprule
        \textbf{Challenge} & \textbf{Summary} & \textbf{Chapter} \\
        \midrule
        \multirow[t]{3}{10em}{Convergence of Local SGD in heterogeneous federated learning} & \multirow[t]{3}{21em}{We provide the first convergence result for Local SGD when the data are heterogeneous. } & \multirow{3}{4em}{\Cref{chapter:local_sgd}} \\ \\ \\[0.3cm]
        \multirow[t]{5}{10em}{Tight convergence guarantees for Random Reshuffling} & \multirow[t]{5}{21em}{We provide upper bounds for strongly convex, convex and non-convex objectives that are simpler and tighter than prior work. Moreover, our upper bound is the first one to match the lower bound.} & \multirow{5}{4em}{\Cref{chapter:rr}} \\ \\ \\ \\ \\[0.3cm]
        \multirow[t]{5}{10em}{Going beyond Local SGD in federated learning} & \multirow[t]{5}{21em}{We design a new algorithm called FedRR that uses the update of random reshuffling in composition with periodic communication of Local SGD, but has a much better asymptotic communication complexity than Local SGD.} & \multirow{5}{4em}{\Cref{chapter:proxrr}} \\ \\ \\ \\ \\[0.3cm]
        \multirow[t]{4}{10em}{The first adaptive stepsize rule for Gradient Descent that provably works} & \multirow[t]{4}{21em}{We develop the first stepsize rule for Gradient Descent that can provably adapt to the local geometry of any convex objective by estimating the local Lipschitz constant of the gradients.} & \multirow{4}{4em}{\Cref{chapter:adaptive}} \\ \\  \\ \\[0.3cm]
        \multirow[t]{3}{10em}{Achieving fast rates in distributed optimization with quantization} & \multirow[t]{3}{21em}{We present the first distributed algorithm with compressed communication that provably preserves fast linear rates of Gradient Descent. 
        } & \multirow{3}{4em}{\Cref{chapter:diana}} \\ \\ \\[0.3cm]
        \multirow[t]{4}{10em}{Developing variance reduction for proximal operators} & \multirow[t]{4}{21em}{We develop a unified theory for a family of methods that converge linearly while having access only to stochastic gradients and proximal operators.} & \multirow{4}{4em}{\Cref{chapter:sdm}} \\ \\ \\ \\[0.3cm]
        \multirow[t]{4}{10em}{Designing variance-reduced algorithms with splitting} & \multirow[t]{4}{21em}{We combine gradients estimation with matrix splitting to obtain a number of linearly-convergent algorithms.} & \multirow{4}{4em}{\Cref{chapter:pddy}} \\ \\ \\ \\
        \bottomrule 
    \end{tabular}
\end{table}

\section{Excluded papers}
To make the thesis more self-consistent, some papers were excluded from it. In particular, we do not discuss here the following works that were finished during my PhD studies:
\begin{itemize}
	\item Stochastic algorithms for constrained minimization~\cite{mishchenko2018stochastic}.
	\item Variance reduction for coordinate descent with non-separable regularizer~\cite{hanzely2018sega}.
	\item A work on delay-tolerant asynchronous gradient method~\cite{mishchenko2018delay} and its extended version with more general analysis~\cite{mishchenko2020distributed}.
	\item Tackling communication bottleneck by update sparsification~\cite{mishchenko201999}.
	\item Variance-reduced extension of the algorithm presented in \Cref{chapter:diana}~\cite{horvath2019stochastic}.
	\item Improved analysis and extension of variance reduction based on the MISO algorithm\cite{qian2019miso}.
	\item Equivalence between the celebrated Sinkhorn algorithm and mirror descent applied to a certain objective~\cite{mishchenko2019sinkhorn}.
	\item A variance-reduced Newton method and its cubic-regularized version \cite{kovalev2019stochastic}.
	\item A quasi-Newton asynchronous method~\cite{soori2020dave}.
	\item Integer compression for communication-efficient distributed training~\cite{mishchenko2021intsgd}.
	\item Hierarchical time series regression~\cite{mishchenko2019self}.
\end{itemize}
We note that some of these papers have been presented as part of the PhD thesis by Filip Hanzely~\cite{hanzely2020optimization} and are excluded for this reason.

\section{Theoretical Background for Most Chapters}
To avoid introducing the same assumptions and standard results in each chapter, we shall now discuss some of them here, excluding those facts that are specific to a single chapter. We also briefly mention some key distinctions of the notation used in specific chapters in this introduction section.

The notions of convex analysis and operator theory that we introduce are very standard and can be found in textbooks. For instance, we recommend the reader to consult Boyd and Vandenberghe~\cite{boyd2004convex}, Nesterov~\cite{Nesterov2013}, and Bauschke and Combettes~\cite{bauschke2017convex} for details and proofs.

\section{Basic Facts and Notation}\label{sec:basic_notions}
In all chapters, we are going to be solving in one way or another minimization of a given function over $\R^d$. The main differentiable part of the objective is always denoted as $f$, while the overall objective, if different from $f$, is denoted by $P$. If there is a single regularization term in the objective, it is denoted by $\psi$, so most of the time we will be solving the problem
\[
	\min_{x\in\R^d} \left[P(x) \eqdef f(x)+\psi(x)\right].
\]
Function $f$ might have different forms depending on the context. In particular, we will sometimes consider it to be of an expectation form with respect to some random variable $\xi$:
\[
	f(x)\eqdef \mathbb{E}_{\xi} \left[ f(x;\xi) \right].
\]
Alternatively, $f$ may take a finite-sum form. For instance, we may use the notation $f(x)=\frac{1}{n}\sum_{i=1}^n f_i(x)$ or $f(x)=\frac{1}{M}\sum_{m=1}^M f_m(x)$, depending on context.

We denote the optimal value of $P$ (or $f$) as $P^*$ (or $f^*$). In case we need to work with a solution of the problem, we denote it by $x^*$, and the set of all solutions by  $\cX^\star$. We denote scalars and vectors with standard letters, for instance, $\alpha\in\R$ or $x\in\R^d$. The iteration index is usually denoted with $k$, so that $x^k$ is the main iterate of the considered algorithm. Matrices and linear operators are denoted with bold capital letters, for instance, $\mL$. For any positive integer $n\ge 1$ we define $[n]\eqdef \{1,2,\dots,n\}$. We denote by $(x^k)_k$ the infinite sequence of elements with values $x^0, x^1, x^2, \dotsc$.

We use $\<\cdot, \cdot>$ to denote the standard Euclidean scalar product of two vectors, and $\|\cdot\|$ to denote the associated Euclidean norm. For any $p\ge 1$, we denote by $\|\cdot\|_p$ the $\ell_p$ norm of a vector and we drop the subscript when $p=2$. For a matrix $\mA\in\R^{d\times m}$, we denote by $\|\mA\|$ its operator norm, and by $\|\mA\|_{2,1}=\sum_{j=1}^m \|\mA_j\|$ the $\ell_{2,1}$ norm, where $\mA_j$ is the $j$-th column of matrix $\mA$. We denote by $\lambda_{\min}(\mA)$ and $\lambda_{\min}^+(\mA)$ the smallest and the smallest positive eigenvalues of $\mA$, respectively.

\subsection{Random vectors}
For any fixed vector $h\in \R^d$, the variance of a random vector $X$ with finite second moment can be decomposed as follows: 
\begin{align}
 \mathbb{E}\left[\|X - \ec{ X}\|^2\right]  = 	\mathbb{E} \left[\|X - h\|^2\right] -  \|\mathbb{E} [X] - h\|^2 . \label{eq:second_moment_decomposition}
\end{align}
In particular, if we plug in $h=0$ and rearrange the terms, we get the standard variance decomposition formula
\begin{align}
    \label{eq:variance_def}
    \ecn {X} = \ecn { X - \ec{X} } + \sqn{\ec{X}}.
\end{align}
If, in addition, $X$ takes only a finite number of values, we get
\begin{align}
    \avein \norm{X_n}^2
    = \avein \norm{X_i - \frac{1}{n}\sumjn X_j}^2 + \norm{\avein X_i}^2.\label{eq:variance_m}
\end{align}
As a consequence of \eqref{eq:variance_def} we also have that
\begin{align}
    \label{eq:variance_sqnorm_upperbound}
    \ecn{X - \ec{X}} \leq \ecn{X}.
\end{align}
In the case when $X$ takes a finite number of values, inequality~\eqref{eq:variance_sqnorm_upperbound} simplifies to
\begin{equation}
  \label{eq:sqnorm-jensen}
  \biggl\| \frac{1}{n} \sum_{i=1}^{n} X_i \biggr\|^2 \leq \frac{1}{n} \sum_{i=1}^{n} \sqn{X_i}.
\end{equation}
After multiplying  both sides of \eqref{eq:sqnorm-jensen} by $n^2$, we  get
\begin{equation}
  \label{eq:sqnorm-sum-bound}
  \biggl\|\sum_{i=1}^{n} X_i\biggr\|^2 \leq n \sum_{i=1}^{n} \sqn{X_i}.
\end{equation}

\subsection{Norms and products} 
We now state some straightforward linear algebra results. Firstly, for any two vectors $a, b \in \R^d$, we have
\begin{equation}
    \label{eq:square-decompos}
    2 \ev{x, y} = \sqn{x} + \sqn{y} - \sqn{x - y}.
\end{equation}

We will also use the following facts, which are sometimes referred to as Young's inequality:
\begin{align}
    \sqn{x + y} &\leq 2 \sqn{x} + 2 \sqn{y}, \label{eq:sum_sqnorm} \\
    2 \ev{x, y} &\leq \zeta \sqn{x} + \zeta^{-1} \sqn{y} \text { for all } x, y \in \R^d \text { and } \zeta > 0. \label{eq:youngs-inequality}
\end{align}

Finally, for any $0\leq \alpha \leq 1$ and $x,y\in \R^d$, we have
\begin{align}\label{eq:sqaured_norm_of_lin_combination}
    \|\alpha x + (1 - \alpha) y  \|^2 =  \alpha \|x\|^2 + (1 - \alpha)\|y\|^2  - \alpha(1 - \alpha) \|x - y\|^2.
\end{align}

\subsection{Function properties} 
We say that an extended real-valued function $f\colon\R^d\to \R\cup \{+\infty\}$ is proper if its domain, $${\rm dom} \; f \eqdef \{x \colon f(x)<+\infty\},$$ is nonempty.  We say that it is convex (respectively closed) if its epigraph, $${\rm epi}\; f \eqdef \{(x,t) \in \R^d\times \R  \colon \; f(x) \leq t\},$$ is a convex (respectively closed) set. Equivalently, $f$ is convex if ${\rm dom} \; f $ is a convex set and $$f(\alpha x + (1-\alpha)y) \leq \alpha f(x) + (1-\alpha) f(y)$$ for all $x,y\in {\rm dom} \; f$ and $\alpha\in(0, 1)$. Finally, $f$ is $\mu$-strongly convex if $f (x) -\frac{\mu}{2}\norm{x}^2$ is convex.

We define the subdifferential of $f$ as the set-valued operator $$\partial f\colon x\in\R^d\mapsto \{g\in\R^d\ \colon\ (\forall y\in\R^d)\ f(x)+\langle y-x, g\rangle\leq f(y)\}.$$ If $f$ is differentiable at $x\in\R^d$, then $\partial f(x)=\{\nabla f(x)\}$, where $\nabla f(x)$ denotes the gradient of $f$ at $x$.

We denote by $f^*$ the conjugate of $f$, defined by $$f^*\colon x\mapsto \sup_{y\in\R^d} \{\langle x,y\rangle -f(y)\},$$ which is always convex, proper and closed. Finally, given any convex set $\cC\subset\R^d$, we define the indicator function 
\[
	\ind_{\cC}(x)\eqdef \begin{cases}0, &\text{if } x\in \cC \\ +\infty, & \text{otherwise} \end{cases}.
\]
Note that this function is always convex, proper and closed. For brevity, if $\cC=\{b\}$ with some $b\in\R^d$, we denote $\ind_{b}\eqdef \ind_{\{b\}}$.

Finally, let us introduce the following standard result for convex functions.
\begin{proposition}[Jensen's inequality]\label{pr:jensen}
    For any convex function $f$ and any vectors $x_1,\dotsc, x_M$ we have
    \begin{align}
        f\br{\avemm x_m}
        \le \avemm f(x_m). \label{eq:jensen}
    \end{align}
\end{proposition}

\subsection{Bregman divergence}
To simplify the notation and proofs, it is convenient to work with Bregman divergences.
It is important to note that the Bregman divergence of a convex function is always nonnegative and is a (non-symmetric) notion of ``distance'' between $x$ and $y$. For $x^\star\in \cX^\star$, the quantity $D_f(x, x^\star)$ serves as a generalization of the functional gap $f(x) - f(x^\star)$ in cases when  $\nabla f(x^\star)\neq 0$.

We denote the Bregman divergence associated with function $f$ and arbitrary $x, y$ as
\begin{align*}
D_f(x, y)
\eqdef f(x) - f(y) - \ev{\nabla f(y), x - y}.
\end{align*}
Moreover, a continuously differentiable function $f$ is called $L$-smooth if its gradient is $L$-Lipschitz, i.e., if
\begin{equation}
	\|\nabla f(x) - \nabla f(y)\|
	\le L\|x-y\|,  \qquad \forall x, y \in \R^d. \label{eq:nabla-Lip}
\end{equation}
A very important consequence of $L$-smoothness is that it implies an upper bound on the Bregman divergence of $f$:
\begin{equation} 
	D_f(x,y) \leq \frac{L}{2} \norm{x - y}^2, \qquad \forall x,y\in \R^d.  \label{eq:L-smooth-intro}
\end{equation}

If $f$ is $\mu$-strongly convex, then we also have
\begin{equation}
     \label{eq:asm-strong-convexity}
     \frac{\mu}{2} \sqn{y - x} \leq D_f(x,y), \qquad \forall x, y \in \R^d.
\end{equation}

The most basic consequence of function smoothness is that the squared gradient norm can be upper bounded with the functional gap. This is formalized in the following proposition.
\begin{proposition}
	Let $f$ be $L$-smooth and lower bounded by $f^\star\in\R$, then 
	\begin{align}
		\|\nabla f(x)\|^2
		\le 2L (f(x)-f^\star),  \qquad \forall x \in \R^d.\label{eq:grad-bound}
	\end{align}
\end{proposition}

Since many of our proofs are easier to write when one uses Bregman divergences, we will formulate the next two lemmas in terms of $D_f(\cdot, \cdot)$.

\begin{proposition}
	Let $f$ be convex and $L$-smooth, then 	\begin{eqnarray}
		\|\nabla f(x) - \nabla f(y)\|^2
		& \le & 2L D_f(x, y),  \qquad \forall x, y \in \R^d,\label{eq:grad_dif_bregman}\\
		\|\nabla f(x) - \nabla f(y)\|^2
		& \le & L\<\nabla f(x) - \nabla f(y), x - y>, \qquad \forall x, y \in \R^d.\label{eq:grad_dif_scalar_prod}
	\end{eqnarray}
\end{proposition}
Sometimes, to make the analysis tighter, we require the following statement.
\begin{proposition}
	Let $f$ be differentiable and $\mu$-strongly convex. Then 
	\begin{eqnarray}
		\frac{\mu}{2}\|x - y\|^2 + D_f(x, y)
		& \le & \<\nabla f(x) - \nabla f(y), x - y>,  \qquad \forall x, y \in \R^d,\label{eq:scal_prod_cvx}\\
		\mu\| x - y\|^2
		& \le & \<\nabla f(x) - \nabla f(y), x - y>,  \qquad \forall x, y \in \R^d. \label{eq:scal_prod_str_cvx}
	\end{eqnarray}
	Moreover, if $f$ is also $L$-smooth, then
	\begin{align}
		\frac{\mu L}{L + \mu}\|x - y\|^2 + \frac{1}{L + \mu}\|\nabla f(x) - \nabla f(y)\|^2
		\le \<\nabla f(x) - \nabla f(y), x - y> \label{eq:scal_prod_tight_str_cvx}
	\end{align}
	holds for all $x, y \in \R^d$.
\end{proposition}
The last inequality is the tightest inequality one can get and, in particular, it implies \eqref{eq:grad_dif_scalar_prod} when $\mu=0$.

\subsection{Proximal operator}

To solve problems with non-smooth regularizer, one of the best approaches is to use proximal operator. Given $\gamma>0$, the proximal operator for function $\psi$ is defined as
\[
	\proxR(u) \eqdef \arg\min_v \left\{ \gamma \psi(v) + \frac{1}{2}\|v - u\|^2 \right\}.
\]
Let us state some basic and well-known properties of the regularized objectives. Firstly, the following lemma explains why the solution of~\eqref{eq:finite-sum-min} is a fixed point of the proximal-gradient step for any stepsize.
\begin{proposition}
	\label{prop:fixed-point} Let $\psi$ be proper, closed and convex. Then point $x^\ast$ is a minimizer of $P(x) =f(x) + \psi (x)$ if and only if for any $\gamma > 0$ we have
	\[ x^\ast = \prox_{\gamma \psi} (x^\ast - \gamma \nabla f(x^\ast)). \]
\end{proposition}
\begin{proof}
	This follows by writing the first-order optimality conditions for problem~\eqref{eq:finite-sum-min}, see \cite[p.32]{parikh2014proximal} for a full proof.
\end{proof}

The proposition above only shows that proximal-gradient step does not hurt if we are at the solution. In addition, we will rely on the following a bit stronger result which postulates that the proximal operator is a contraction (respectively strong contraction) if the regularizer $\psi$ is convex (respectively strongly convex). 
\begin{lemma}
	\label{prop:prox-contraction} Let $\psi$ by proper and closed. If $\psi$ is $\mu$-strongly convex with $\mu\ge 0$, then for any $\gamma>0$ we have
	\begin{equation}
		\|\prox_{\gamma \psi}(x)-\prox_{\gamma \psi}(y)\|^2
		\le \frac{1}{1+2\gamma\mu }\|x - y\|^2, \label{eq:prox_non_exp}
    \end{equation}
    for all $x, y \in \R^d$.
\end{lemma}
\begin{proof}
	Let $u\eqdef \prox_{\gamma n\psi}(x)$ and $v\eqdef \prox_{\gamma n\psi}(y)$. By definition, $u=\argmin_w \{\psi(w) + \frac{1}{2\gamma n}\|w-x\|^2 \}$. By first-order optimality, we have $0\in \partial \psi(u)+ \frac{1}{\gamma n}(u - x)$ or simply $x-u\in  \gamma n\partial\psi(u)$. Using a similar argument for $v$, we get $x-u-(y-v)\in \gamma n(\partial \psi(u)-\partial\psi(v))$. Thus, by strong convexity of $\psi$, we get
	\[
		\<x-u-(y-v), u-v>\ge \gamma \mu n\|u-v\|^2.
	\]
	Hence,
	\begin{align*}
		\|x-y\|^2
		&= \|u - v + (x-u - (y-v)) \|^2  \\
		&=\|u - v\|^2 + 2\<x-u - (y-v), u-v> + \|x-u-(y-v)\|^2 \\
		&\ge \|u - v\|^2 + 2\<x-u - (y-v), u-v>  \\
		&\ge (1+2\gamma \mu n)\|u-v\|^2. 
	\end{align*}
\end{proof}

\begin{figure}[t]
	\center
	\includegraphics[scale=0.09]{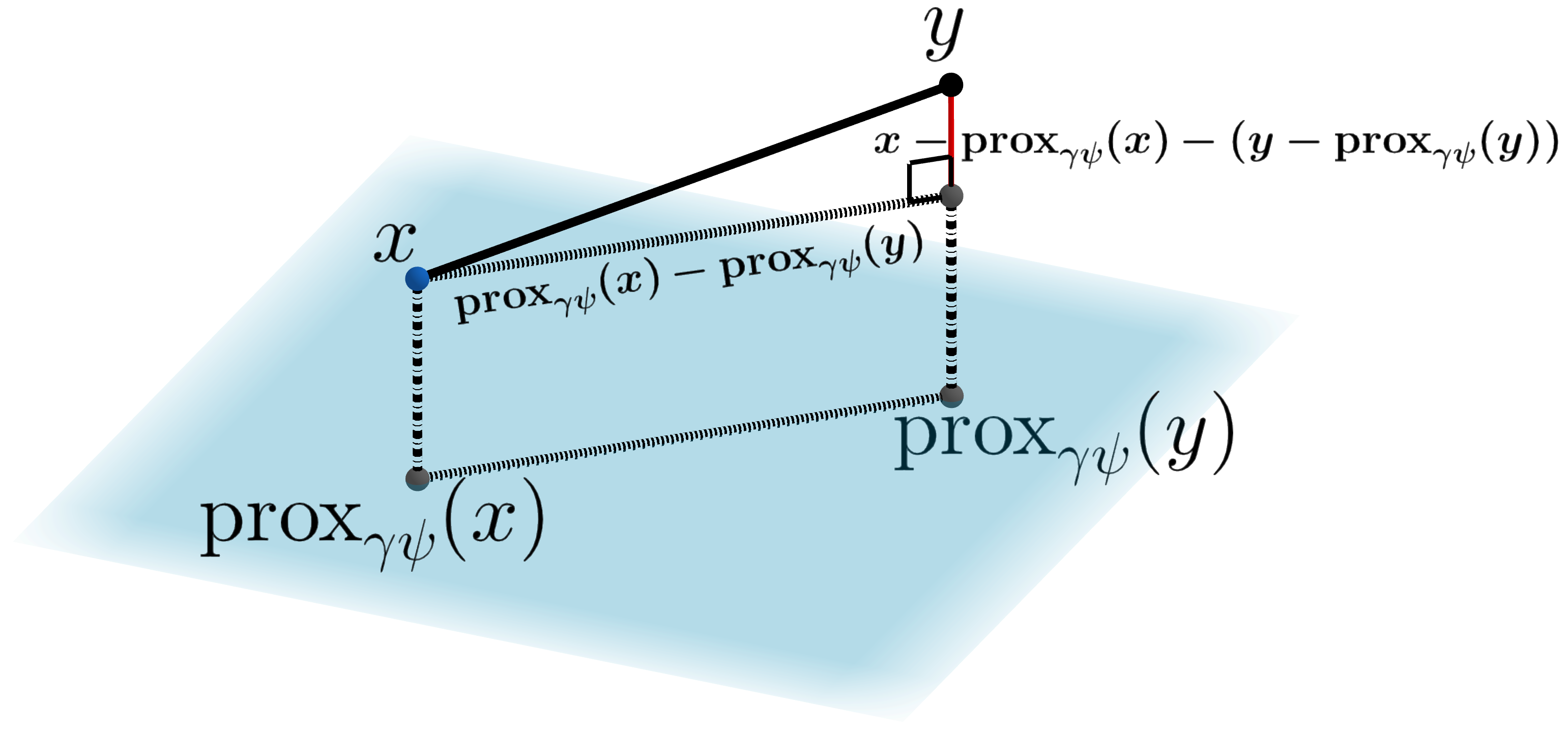}
	\caption{Illustration of property~\eqref{eq:nonexp} with characteristic function of a linear subspace, $\psi(x)=\ind_{\{x \;:\; a^\top x = b\}}$. In this case the proximal operator returns the projection of a point onto the subspace, and Inequality~\eqref{eq:nonexp} becomes identity and follows from Pythagorean theorem.}
	\label{fig:prox}
\end{figure}

Another important property of the proximal operator is firm nonexpansiveness given below.
\begin{proposition} Let $\psi\colon\RR^d \to \RR \cup \{+\infty\}$ be a proper closed  convex function. Then its proximal operator is firmly nonexpansive. That is, for all $\gamma\in \RR$, 
\begin{align}
    &\|\prox_{\gamma \psi}(x) - \prox_{\gamma \psi}(y)\|^2 \\
    &\qquad\le \|x - y\|^2 - \left(1 + \frac{1}{\gamma L_\psi}\right)\|x - \prox_{\gamma \psi}(x) - (y - \prox_{\gamma \psi}(y))\|^2, \label{eq:nonexp}
\end{align}
where $L_\psi\in\RR\cup \{+\infty\}$ is the smoothness constant of function $\psi$ (for non-smooth functions, $L_\psi=+\infty$). In particular, regardless of the smoothness properties of $\psi$, we have
\begin{align}
    \|\prox_{\gamma \psi}(x) - \prox_{\gamma \psi}(y)\|^2
    \le \|x - y\|^2, \label{eq:nonexpansive}
\end{align}
\end{proposition}
Moreover, we would like to note that Equation~\ref{eq:nonexp} is tight if $\psi(x) = \ind_{\{x\;:\;a^\top x = b\}}$ for some $a\in\R^d$ and $b\in\R$, as is shown in Figure~\ref{fig:prox}.

\subsection{Monotone operators}
Consider a set-valued operator $M\colon \cZ \rightrightarrows \cZ$. The inverse $M^{-1}$ of $M$ is defined by the relation $$z' \in M(z) \Leftrightarrow z \in M^{-1}(z').$$ The set of zeros of $M$ is $$\mathrm{zer}(M) = M^{-1}(0) = \{z \in \cZ, 0 \in M(z)\}.$$ The operator $M$ is monotone if $$\ps{w-w',z-z'} \geq 0$$ whenever $u \in M(z)$ and $u' \in M(z')$, and strongly monotone if there exists $\mu >0$, such that $$\ps{w-w',z-z'} \geq \mu\|z-z'\|^2.$$ The resolvent operator of $M$ is defined by $$J_{M} = (I + M)^{-1},$$ where $I$ denotes the identity. If $M$ is monotone, then $J_M(z)$ is either empty or single-valued. $M$ is maximal monotone  if $J_M(z)$ is single-valued, for every $z \in \cZ$. We identify single-valued operators as operators from $\cZ$ to $\cZ$. If $f$ is convex, proper and closed, then $\partial f$ is maximal monotone, $$J_{\partial f} = \prox_{f}, \quad \mathrm{zer}(\partial f) = \argmin f, \quad \text{and} \quad (\partial f)^{-1} = \partial f^*.$$

A single-valued operator $M$ on $\cZ$ is $\alpha$-cocoercive if $$\alpha\|M(z) - M(z')\|^2 \leq \ps{M(z) - M(z'),z-z'}.$$ 
 The resolvent of a maximal monotone operator is $1$-cocoercive. In addition, $\nabla f$ is $1/\nu$-cocoercive for any $\nu$-smooth function $f$.

Let $\cX, \cY$ be real Hilbert spaces and let $\mL \colon \cX \to \cY$ be a linear operator. The adjoint of  $\mL$ is denoted by $\mL^*\colon\mathcal{Y}\rightarrow \mathcal{X}$, and the operator norm of $\mL$ is defined as $$\|\mL\| = \sup\{\|\mL x\|, x \in \cX, \|x\|\leq 1\}.$$ The largest eigenvalue of $\mL\mL^*$ is $$\|\mL\mL^*\| = \|\mL\|^2 = \|\mL^*\|^2.$$ Let $\mP\colon \cZ \to \cZ$ be a linear and symmetric operator ($\mP^*=\mP$). $\mP$ is positive semidefinite if $$\ps{\mP z,z} \geq 0$$ for every $z \in \cZ$,  and positive definite if, additionally, $\ps{\mP z,z} = 0$ implies $z=0$. In this latter case, the inner product induced by $\mP$ is defined by $\ps{z,z'}_\mP = \ps{\mP z,z'}$ and the norm induced by $\mP$ is defined by $$\|z\|_\mP^2 = \ps{z,z}_\mP.$$ We denote by $\mathcal{Z}_\mP$ the space $\mathcal{Z}$ endowed with $\ps{\cdot,\cdot}_\mP$.

\subsection{Differences in section-specific notation}
Many of the chapters deal with finite-sum problems but use slightly different notation for the number of functions used. We use $M$ to denote the number of machines in Chapters~\ref{chapter:local_sgd} and \ref{chapter:proxrr}, so the finite sum in the objective has $M$ summands. At the same time, the number of terms that correspond to data points is denoted by $n$ in Chapters~\ref{chapter:rr}, \ref{chapter:proxrr}, \ref{chapter:sdm}, \ref{chapter:pddy}. The linear operators are denoted by $\mA_1,\dotsc, \mA_m$ in Chapter~\ref{chapter:sdm} and by $\mL$ in Chapter~\ref{chapter:pddy}.

In \Cref{chapter:local_sgd}, we pay particular attention to the nature of stochastic gradients. Therefore, we will be denoting different data/noise distributions as $\D_1, \dotsc, \D_M$. Correspondingly, we will use $\mathbb{E}_{\xi\sim \D_m}[\cdot]$ to denote the expectation with respect to distribution $\D_m$, where $m\in [M]$.

For the reader's convenience, we provide notation tables in \Cref{sec:table}. They summarize the notation we introduce above as well as provide additional details for specific chapters.

\chapter{Convergence of Local SGD for Federated Learning in the Heterogeneous Data Regime}
\label{chapter:local_sgd}

\graphicspath{{local_sgd/}}


\section{Introduction}
Modern hardware increasingly relies on the power of uniting many parallel units into one system. This approach requires optimization methods that target specific issues arising in distributed environments such as decentralized data storage. Not having data in one place implies that computing nodes have to communicate back and forth to keep moving toward the solution of the overall problem. A number of efficient first-, second-order and dual methods that are capable of reducing the communication overhead existed in the literature for a long time, some of which are in certain sense optimal.

Yet, when federated learning (FL) showed up, it turned out that the problem of balancing the communication and computation had not been solved. On the one hand, \algname{Mini-Batch Stochastic Gradient Descent} (\algname{SGD}), which averages the result of stochastic gradient steps computed in parallel on several machines, again demonstrated its computation efficiency. Seeking communication efficiency, Kone\v{c}n\'{y} et al.\ \cite{FEDLEARN} and McMahan et al.\ \cite{McMahan17}) proposed to use a natural variant of \algname{Mini-Batch SGD}---\algname{Local SGD} (Algorithm~\ref{alg:local_sgd}), which does a few \algname{SGD} iterations locally on each involved node and only then computes the average. This approach saves a lot of time on communication, but, unfortunately, in terms of theory things were not as great as in terms of practice and there are still gaps in our understanding of \algname{Local SGD}. 

The idea of \algname{Local SGD} in fact is not recent, it traces back to the work of Mangasarian~\cite{Mangasarian95} and has since been popular among practitioners from different communities. An asymptotic analysis can be found in~\cite{Mangasarian95} and quite a few recent papers proved new convergence results, making the bounds tighter with every work. The theory has been developing in two important regimes: identical and heterogeneous data.

The identical data regime is more of interest if the data are actually stored in one place. In that case, we can access it on each computing device at no extra cost and get a fast, scalable method. Although not very general, this framework is already of interest to a wide audience due to its efficiency in training large-scale machine learning models~\cite{TaoLin19}. The first contribution of this chapter is to provide the fastest known rate of convergence for this regime under weaker assumptions than in prior work.

Federated learning, however, is done on a very large number of mobile devices, and is operating in a highly non-i.i.d.\ regime. To address this, we present the first analysis of \algname{Local SGD} that applies to arbitrarily heterogeneous data, while all previous works assumed a certain type of similarity between the data or local gradients.

\begin{algorithm*}[t]
   \caption{\algname{Local Stochastic Gradient Descent} (\rrbox{\algname{Local SGD}})}
   \label{alg:local_sgd}
\begin{algorithmic}[1]
  \Require Stepsize $\gamma > 0$, initial vector $x^0 = x^0_m$ for all $m \in [M]$, number of steps $K$, synchronization timesteps $k_1, k_2, \ldots$
   \For{$k=0,1,\dotsc, K-1$}
      \For{$m=1,\dotsc, M$ in parallel}
         \State Sample $\xi_m \overset{\text{i.i.d.}}{\sim} \D_m$
         \If{data are identical}
            \State Compute $g^k_m = \nabla f(x^k_m; \xi_m)$ such that $\ec{g^k_m \mid x^k_m} = \nabla f(x^k_m)$
         \Else
            \State Compute $g^k_m = \nabla f_m(x^k_m; \xi_m)$ such that $\ec{g^k_m\mid x^k_m}=\nabla f_m (x^k_m)$
         \EndIf
         \State $x^{k+1}_m=
         \begin{cases}
         	\frac{1}{m}\sum_{j=1}^M (x^k_j - \gamma g^k_j), & \text{ if } k = k_p \text { for some } p \in \N \\
         x^k_m - \gamma g^k_m, & \text{ otherwise. }
         \end{cases}$
      \EndFor
   \EndFor
\end{algorithmic}
\end{algorithm*}

To explain the challenge of heterogeneity better, let us introduce the problem we are trying to solve. Given that there are $M$ devices and corresponding local losses $f_m\colon\R^d\to \R$, we want to find
\begin{equation}
     \label{eq:optimization-problem}
     \min_{x \in \R^d} \left[ f(x) = \frac{1}{M} \sum_{m=1}^{M} f_m (x)\right].
\end{equation}
In the case of identical data, we are able to obtain on each node an unbiased estimate of the gradient $\nabla f$. In the case of heterogeneous data, $m$-th node can only obtain an unbiased estimate of the gradient $\nabla f_m$. Data similarity can then be formulated in terms of the differences between functions $f_1,\dotsc, f_M$. If the underlying data giving rise to the loss functions are i.i.d., the function share optima and one could even minimize them separately, averaging the results at the end. We will demonstrate this rigorously later in this chapter.

If the data are dissimilar, however, we need to be much more careful since running \algname{SGD} locally will yield solutions of local problems. Clearly, their average might not minimize the true objective~\eqref{eq:optimization-problem}, and this poses significant issues for the convergence of \algname{Local SGD}. 

To properly discuss the efficiency of \algname{Local SGD}, we also need a practical way of quantifying it. Normally, a method's efficiency is measured by the total number of times each function $f_m$ is touched and the cost of the touches. On the other hand, in distributed learning we also care about how much information each computing node needs to communicate. In fact, when communication is as expensive as is the case in FL, we predominantly care about communication. The question we address in this chapter, thus, can be posed as follows: how many times does each node need to communicate if we want to solve~\eqref{eq:optimization-problem} up to accuracy $\e$? Equivalently, we can ask for the optimal \textit{synchronization interval length} between communications, $H$, i.e., how many computation steps per one communication we can allow for. We next review related work and then present our contributions.
\begin{table*}[h]
    \caption{Existing theoretical bounds for \algname{Local SGD} for identical data with convex objectives.}
    \begin{threeparttable}[b]
     \centering
     \def\arraystretch{1.2}
     \resizebox{\textwidth}{!}{
     \begin{tabular}{ccccccccc}
     \toprule
          \begin{tabular}{c}Unbounded\\ gradient \end{tabular} & \begin{tabular}{c} $H=K$ \\ convergent \end{tabular} & \begin{tabular}{c} \# communications \\ $f$ strongly convex \end{tabular} & \begin{tabular}{c} \# communications \\ $f$ convex \end{tabular} & Reference \\
     \midrule    
          \xmark & \xmark &  $\Om(\sqrt{ M K })$ & \xmark & \cite{Stich2018} \\
          \xmark & \xmark & $\Om\br{\sqrt{M K}}$ & \xmark & \cite{Basu2019} \\
          \cmark & \xmark & $\tilde{\Om}\br{M}$ & $\Om(M^{3/2} K^{1/2})$ & \cite{stich19errorfeedback}\\
          \cmark & \xmark & $\tilde{\Om}\br{M^{1/3} K^{1/3} }$\tnote{b} & --- & \cite{Haddadpour2019local} \\
     \hline
          \cmark & \cmark & $\tilde{\Om}(M)$ & $\Om(M^{3/2} K^{1/2})$ & THIS THESIS \\
     \bottomrule
     \end{tabular}
     }
     \begin{tablenotes}
         \item [a] $C(K)$ denotes the minimum number of communication steps required at each of $K$\\ iterations to achieve a linear speedup in the number of nodes $M$.
         \item [b] The PL inequality, a generalization of strong convexity, is assumed in \cite{Haddadpour2019local}, but for\\ comparison we specialize to strong convexity.
     \end{tablenotes}
    \end{threeparttable}
    \label{tab:related_work_convexity}
\end{table*}

\section{Related Work}
While \algname{Local SGD} has been used among practitioners for a long time, see, e.g., \cite{Coppola15, McDonald10}, its theoretical analysis has been limited until recently. Early theoretical work on the convergence of local methods exists as in~\cite{Mangasarian95}, but no convergence rate was given there. The previous work can mainly be divided into two groups: those assuming identical data (that all nodes have access to the same dataset) and those that allow each node to hold its own dataset. As might be expected, the analysis in the latter case is more challenging, more limited, and usually shows worse rates. We note that in recent work more sophisticated local stochastic gradient methods have been considered, for example with momentum~\cite{yu2019linear, wang2019slowmo}, with quantization~\cite{Reisizadeh19,Basu2019}, with adaptive stepsizes~\cite{xie2019local} and with various variance-reduction methods \cite{liang2019variance, sharma2019parallel, karimireddy2019scaffold}. Our work is complimentary to these approaches, and provides improved rates and analysis for the vanilla method.

\renewcommand{\thefootnote}{\arabic{footnote}}
\subsection{\algname{Local SGD} with identical data}
The analysis of \algname{Local SGD} in this setting shows that a reduction in communication is possible without affecting the asymptotic convergence rate of \algname{Mini-Batch SGD} with $M$ nodes (albeit with usually worse dependence on constants). An overview of related work on \algname{Local SGD} for convex objectives is given in Table~\ref{tab:related_work_convexity}. We note that analysis for non-convex objectives has been carried out in a few recent works~\cite{Zhou18,Wang18,Jiang18}, but our focus in this thesis is on convex objectives and hence they were not included in Table~\ref{tab:related_work_convexity}. The comparison shows that we attain superior rates in the strongly convex setting to previous work with the exception of the concurrent work of Stich and Karimireddy~\cite{stich19errorfeedback} and we attain these rates under less restrictive assumptions on the optimization process compared to them. We further provide a novel analysis in the convex case, which has not been previously explored in the literature, with the exception of \cite{stich19errorfeedback}. Their analysis attains the same communication complexity but is much more pessimistic about possible values of $H$. In particular, it does not recover the convergence of one-shot averaging, i.e., substituting $H = K$ or even $H = K/M$ gives noninformative bounds, unlike our Theorem~\ref{thm:sc-convergence-theorem}.

In addition to the works listed in the table, Dieuleveut and Patel~\cite{Patel19}  also analyze \algname{Local SGD} for identical data under a Hessian smoothness assumption in addition to gradient smoothness, strong convexity, and uniformly bounded variance. However, we believe that there are issues in their proof that we explain in Section~\ref{sec:patel-discuss} in the supplementary material. As a result, the work is excluded from the table. 

\subsection{\algname{Local SGD} with heterogeneous data}
\begin{table*}[t]
    \caption{Existing theoretical bounds for \algname{Local SGD} with heterogeneous data.}
    \centering
    \def\arraystretch{1.4}
    \resizebox{\textwidth}{!}{
    \begin{tabular}{ccccccccccc}
    \toprule
         \begin{tabular}{c}Unbounded\\ gradient \end{tabular} & \begin{tabular}{c}Unbounded\\ dissimilarity/diversity \end{tabular} & \begin{tabular}{c} \# communications \\ $f$ strongly convex \end{tabular} & \begin{tabular}{c} \# communications \\ $f$ convex \end{tabular} & \begin{tabular}{c} \# communications \\ $f$ non-convex \end{tabular} & Reference \\
    \midrule
         \xmark & \xmark & --- & --- & $\Om\br{M^{3/4} K^{3/4}}$ & \cite{Yu18} \\
         \cmark & \xmark & --- & --- & $\Om(K)$ & \cite{Jiang18} \\
         \xmark & \xmark & $\Om\br{\sqrt{MK}}$ & --- & $\Om\br{M^{3/4} K^{3/4}}$ & \cite{Basu2019} \\
         \cmark & \xmark & $\Om\br{M^{1/3} K^{1/3}}$ & --- & $\Om\br{M^{3/2} K^{1/2}}$ & \cite{Haddadpour19FL}\\
    \hline
         \cmark & \cmark & --- & $\Om\br{M^{3/4} K^{3/4}}$ & --- & THIS THESIS \\
    \bottomrule
    \end{tabular}
    }
    \label{tab:related_work_non-iid}
\end{table*}

An overview of related work on \algname{Local SGD} in this setting is given in Table~\ref{tab:related_work_non-iid}. In addition to the works in Table~\ref{tab:related_work_non-iid}, Wang et al.~\cite{WangTuor18} analyze a local gradient descent method under convexity, bounded dissimilarity, and bounded gradients, but do not show convergence to arbitrary precisions. Li et al.~\cite{Li2019} analyze Federated Averaging (discussed below) in the strongly convex and non-convex cases under bounded gradient norms. However, their result is not included in Table~\ref{tab:related_work_non-iid} because in the more general setting of Federated Averaging, their analysis and experiments suggest that retaining a linear speedup is not possible. 

\algname{Local SGD} is at the core of the Federated Averaging algorithm which is popular in federated learning applications~\cite{FEDLEARN}. Essentially, Federated Averaging is a variant of \algname{Local SGD} with participating devices sampled randomly. This algorithm has been used in several machine learning applications such as mobile keyboard prediction~\cite{HardRao18}, and strategies for improving its communication efficiency were explored in~\cite{FEDLEARN}. Despite its empirical success, little is known about convergence properties of this method and it has been observed to diverge when too many local steps are performed~\cite{McMahan17}. This is not so surprising as the majority of common assumptions are not satisfied; in particular, the data are typically very non-i.i.d.~\cite{McMahan17}, so the local gradients can point in different directions. This property of the data can be written for any vector $x$ and indices $i,j$ as
\begin{align*}
	\norm{\nabla f_i(x) - \nabla f_j(x)} \gg 1.
\end{align*}
Unfortunately, it is very hard to analyze local methods without assuming a bound on the dissimilarity of $\nabla f_i(x)$ and $\nabla f_j(x)$. For this reason, almost all prior work assumed some regularity notion over the functions such as bounded dissimilarity \cite{yu2019linear, Li2019, Yu18, WangTuor18} or bounded gradient diversity \cite{Haddadpour19FL} and addressed other less challenging aspects of federated learning such as decentralized communication, non-convexity of the objective or unbalanced data partitioning. In fact, a common way to make the analysis simple is to assume Lipschitzness of local functions,
$
	\norm{\nabla f_i(x)} \le G
$
for any $x$ and $i$. We argue that this assumption is pathological and should be avoided when seeking a meaningful convergence bound. First of all, in unconstrained strongly convex minimization this assumption cannot be satisfied, and both are assumed by Stich~\cite{Stich2018}. Second, there exists at least one method, whose convergence is guaranteed under bounded variance~\cite{juditsky2011solving}, but in practice the method diverges~\cite{chavdarova2019reducing, mishchenko2019revisiting}. Finally, under the bounded gradients assumption we have
\begin{align*}
    \label{eq:related-work-2}
    \norm{\nabla f_i (x) - \nabla f_{j} (x)} \leq \norm{\nabla f_i (x)} + \norm{\nabla f_j (x)} \leq 2G.
\end{align*}
In other words, we lose control over the difference between the functions. Since $G$ bounds not just dissimilarity, but also the gradients themselves, it makes the statements less insightful or even vacuous. For instance, it is not going to be tight if the data are actually i.i.d.\ since $G$ in that case will remain a positive constant. In contrast, we will show that the rate should depend on a much more meaningful quantity, 
$$
	\sigmaf^2 \eqdef \frac{1}{M} \sum_{m=1}^{M} \ecn[\xi_m \sim \D_m]{\nabla f_m (x^\ast; \xi_m)},
$$
where $x^\ast$ is a fixed minimizer of $f$ and $f_m (\cdot; \xi_m)$ for $\xi_m \sim \D$ are stochastic realizations of $f_m$ (see the next section for the setting). Obviously, for all nondegenerate sampling distributions $\D_m$ the quantity $\sigmaf$ is finite and serves as a natural measure of variance in local methods. We note that an attempt to get more general convergence statement has been made by Li et al.~\cite{Sahu18}, but unfortunately their guarantee is strictly worse than that of \algname{Mini-Batch Stochastic Gradient Descent} (\algname{SGD}). In the overparameterized regime where $\sigmaf = 0$, Zhang and Li~\cite{zhang2019distributed} prove the convergence of \algname{Local SGD} with arbitrary $H$.

In an earlier workshop paper~\cite{khaled2019analysis}, we explicitly analyzed Local Gradient Descent (\algname{Local GD}) as opposed to \algname{Local SGD}, where there is no stochasticity in the gradients. An analysis of \algname{Local GD} for non-convex objectives with the PL inequality and under bounded gradient diversity was subsequently carried out by Haddadpour and Mahdavi~\cite{Haddadpour19FL}. 

\section{Settings and Contributions}
We use a notation similar to that of Stich~\cite{Stich2018} and denote the sequence of time stamps when synchronization happens as  $(k_{p})_{p=1}^{\infty}$. Given stochastic gradients $g^k_1, g^k_2, \ldots, g^k_M$ at time $k \geq 0$ we define
\begin{align*}
    g^k \eqdef \avemm g^k_m, && \bar{g}^k_m \eqdef \ec{g^k_m} = \begin{cases}
     \nabla f(x^k_m) & \text{for identical data }  \\
     \nabla f_m (x^k_m) & \text{otherwise. } 
\end{cases}, && \bar{g}^k \eqdef \ec{g^k}.
\end{align*}

We define an epoch to be a sequence of timesteps between two synchronizations: for $p \in \N$ an epoch is the sequence $k_{k_p}, k_{k_p + 1}, \ldots, k_{k_{p+1} - 1}$.  We summarize some of the notation used in Table~\ref{tab:notation-summary}.

The assumptions of this chapter are a bit different from that of others because we consider $M$ potentially different distribution $\D_1,\dotsc, \D_M$. As we shall show below, this distinction is quite important because when the data are identical, the convergence becomes significantly faster.
\begin{assumption}
     \label{asm:convexity-and-smoothness}
     Assume that the set of minimizers of \eqref{eq:optimization-problem} is nonempty. Each $f_m$ is $\mu$-strongly convex for $\mu \geq 0$ and $L$-smooth. That is, for all $x, y \in \R^d$
     \begin{align*}
          \frac{\mu}{2} \sqn{x - y} 
          &\leq f_m (x) - f_m (y) - \ev{\nabla f_m (y), x - y} 
          \leq \frac{L}{2} \sqn{x - y}.
     \end{align*}
     When $\mu = 0$, we say that each $f_m$ is just convex. When $\mu \neq 0$, we define $\kappa \eqdef \frac{L}{\mu}$, the condition number.
\end{assumption}

Assumption~\ref{asm:convexity-and-smoothness} formulates our requirements on the overall objective. Next, we have two different sets of assumptions on the stochastic gradients that model different scenarios, which also lead to different convergence rates.

\begin{assumption}
     \label{asm:uniformly-bounded-variance}
     Given a function $h$, a point $x \in \R^d$, and a sample $\xi \sim \D$ drawn i.i.d.\ according to a distribution $\D$, the stochastic gradients $\nabla h(x; \xi)$ satisfy
     $
          \ec[\xi \sim \D]{\nabla h(x; \xi)} = \nabla h (x),
          \ecn[\xi \sim \D]{\nabla h(x; \xi) - \nabla h(x)} \leq \sigma^2.
     $
\end{assumption}

Assumption~\ref{asm:uniformly-bounded-variance} holds for example when $\nabla h(x; \xi) = \nabla h(x) + \xi$ for a random variable $\xi$ of expected bounded squared norm: $\ecn[\xi \sim \D]{\xi} \leq \sigma^2$. Assumption~\ref{asm:uniformly-bounded-variance}, however, typically does not hold for finite-sum problems where $g(x; \xi)$ is a gradient of one of the functions in the finite-sum. To capture this setting, we consider the following assumption:

\begin{assumption}
    \label{asm:finite-sum-stochastic-gradients}
    Given an $L$-smooth and $\mu$-strongly convex (possibly with $\mu = 0$) function $h\colon \R^d \to \R$ written as an expectation $h = \ec[\xi \sim \D]{h(x; \xi)}$, we assume that a stochastic gradient $g = g(h, x; \xi)$ is computed by $g(h, x; \xi) = \nabla h(x; \xi).$
    We assume that $h(\cdot, \xi)\colon \R^d \to \R$ is almost-surely $L$-smooth and $\mu$-strongly convex (with the same $L$ and $\mu$ as $h$).
\end{assumption}

When Assumption~\ref{asm:finite-sum-stochastic-gradients} is assumed in the identical data setting, we assume it is satisfied on each node $m \in [M]$ with $h = f$ and distribution $\D_m$, and we define as a measure of variance at the optimum
\[ \sigmaopt^2 \eqdef \frac{1}{M} \sum_{m=1}^{M} \ecn[\xi_m \sim \D_m]{\nabla f (x^\ast; \xi_m)}. \] 
Whereas in the heterogeneous data setting we assume that it is satisfied on each node $m \in [M]$ with $h = f_m$ and distribution $\D_m$, and we analogously define
\[ \sigmaf^2 \eqdef \frac{1}{M} \sum_{m=1}^{M} \ecn[\xi_m \sim \D_m]{\nabla f_m (x^\ast; \xi_m)}. \]

Assumption~\ref{asm:finite-sum-stochastic-gradients} holds, for example, for finite-sum optimization problems with uniform sampling and permits direct extensions to more general settings such as expected smoothness \cite{Gower2019}.

\textbf{Our contributions.} In this chapter, we achieve the following:
\begin{enumerate}
     \item In the identical data setting under Assumptions~\ref{asm:convexity-and-smoothness} and~\ref{asm:uniformly-bounded-variance} with $\mu > 0$, we prove that the iteration complexity of \algname{Local SGD} to achieve $\e$-accuracy is 
     \[ \mathcal{\tilde{O}}\br{ \frac{\sigma^2}{\mu^2 M \e} } \]
     in squared distance from the optimum provided that $K = \Omega\br{ \kappa \br{H  - 1} }$. This improves the communication complexity in prior work (see Table~\ref{tab:related_work_convexity}) with a tighter results compared to concurrent work (recovering convergence for $H = 1$ and $H = K$). When $\mu = 0$ we have that the iteration complexity of \algname{Mini-Batch SGD} to attain an $\e$-accurate solution in functional suboptimality is 
     \[ \mathcal{O} \br{\ \frac{L^2 \norm{x^0 - x^\ast}^4}{M \e^2} + \frac{\sigma^4}{L^2 M \e^2}}, \] 
     provided that $K = \Omega\br{M^3 H^2}$. We further show that the same $\e$-dependence holds in both the $\mu > 0$ and $\mu = 0$ cases under Assumption~\ref{asm:finite-sum-stochastic-gradients}. This has not been explored in the literature on \algname{Local SGD} before, and hence we obtain the first results that apply to arbitrary convex and smooth finite-sum problems. 
     \item When the data on each node is different and Assumptions~\ref{asm:convexity-and-smoothness} and \ref{asm:finite-sum-stochastic-gradients} hold with $\mu = 0$, the iteration complexity needed by \algname{Local SGD} to achieve an $\e$-accurate solution in functional suboptimality is 
     \[ \mathcal{O} \br{ \frac{L^2 \norm{x^0 - x^\ast}^4}{M \e^2} + \frac{\sigmaf^4}{L^2 M \e^2}} \]
     provided that $K = \Omega(M^3 H^4)$. This improves upon previous work by not requiring any restrictive assumptions on the gradients and is the first analysis to capture true data heterogeneity between different nodes.
     \item We verify our results by experimenting with logistic regression on multiple datasets, and investigate the effect of heterogeneity on the convergence speed.
\end{enumerate}

\section{Convergence Theory}

The following quantity is crucial to the analysis of both variants of \algname{Local SGD}, and measures the deviation of the iterates from their average $\hat{x}^k$ over an epoch:
\[ V^k \eqdef \avemm \sqn{x^k_m - \hat{x}^k}\ \text { where }\ \hat{x}^k \eqdef \avemm x^k_m. \]
To prove our results, we follow the line of work started by Stich~\cite{Stich2018} and first show a recurrence similar to that of \algname{SGD} up to an error term proportional to $V^k$, then we bound each $V^k$ term individually or the sum of $V^k$'s over an epoch. All proofs are relegated to the supplementary material.

\subsection{Identical data}
Our first lemma presents a bound on the sequence of the $V^k$ in terms of the synchronization interval $H$.
\begin{lemma}
     \label{lemma:uniform-var-iterate-variance-bound}
     Choose a stepsize $\gamma > 0$ such that $\gamma \leq \frac{1}{2L}$. Under Assumptions~\ref{asm:convexity-and-smoothness}, and \ref{asm:uniformly-bounded-variance} we have that for Algorithm~\ref{alg:local_sgd} with $\max_{p} \abs{k_p - k_{p+1}} \leq H$ and with identical data, for all $k \ge 1$
     \begin{align*}
         \ec{V^{k}} \leq \br{H - 1} \gamma^2 \sigma^2.
     \end{align*}
 \end{lemma}

Combining Lemma~\ref{lemma:uniform-var-iterate-variance-bound} with perturbed iterate analysis as in \cite{Stich2018} we can recover the convergence of \algname{Local SGD} for strongly-convex functions:
 \begin{theorem}
     \label{thm:sc-convergence-theorem}
     Suppose that Assumptions~\ref{asm:convexity-and-smoothness}, and \ref{asm:uniformly-bounded-variance} hold with $\mu > 0$. Then for Algorithm~\ref{alg:local_sgd} run with identical data, a constant stepsize $\gamma > 0$ such that $\gamma \leq \frac{1}{4L}$, and $H \geq 1$ such that $\max_{p} \abs{k_p - k_{p+1}} \leq H$,
     \begin{align}
         \label{eq:thm-sc-convergence-rate}
         \begin{split}
          \ecn{\hat{x}^K - x^\ast} \leq (1 &- \gamma \mu)^{K} \sqn{x^0 - x^\ast} + \frac{\gamma \sigma^2}{\mu M} + \frac{2 L \gamma^2 \br{H - 1} \sigma^2}{\mu}.
         \end{split}
     \end{align}
 \end{theorem}

By \eqref{eq:thm-sc-convergence-rate} we see that the convergence of \algname{Local SGD} is the same as \algname{Mini-Batch SGD} plus an additive error term which can be controlled by controlling the size of $H$, as the next corollary and the successive discussion show.

\begin{corollary}
    \label{corrollary:sc-convergence-comm-complexity}
     Choosing $\gamma = \frac{1}{\mu a}$, with $a = {4 \kappa + k}$ for $k > 0$ and we take $K = 2 a \log a$ steps. Then substituting in \eqref{eq:thm-sc-convergence-rate} and using that $1 - x \leq \exp(-x)$ and some algebraic manipulation we can conclude that,
     \begin{equation*}
        \ecn{\hat{x}^K - x^\ast} = \tilde{\mathcal{O}} \br{ \frac{\sqn{x^0 - x^\ast}}{K^2} + \frac{\sigma^2}{\mu^2 M K} + \frac{\kappa \sigma^2 (H-1)}{\mu^2 K^2}}.
     \end{equation*}
    where $\tilde{\mathcal{O}} (\cdot)$ ignores polylogarithmic and constant numerical factors. 
\end{corollary}

\textbf{Recovering fully synchronized \algname{Mini-Batch SGD}.} When $H = 1$ the error term vanishes and we obtain directly the ordinary rate of \algname{Mini-Batch SGD}. 

\textbf{Linear speedup in the number of nodes $M$.} We see that choosing $H = \mathcal{O}(K/M)$ leads to an asymptotic convergence rate of $\mathcal{\tilde{O}}\br{ \frac{\sigma^2 \kappa}{\mu^2 M K} }$ which shows the same linear speedup of \algname{Mini-Batch SGD} but with worse dependence on $\kappa$. The number of communications in this case is then $C(K) = K/H = \tilde{\Omega}(M)$.

\textbf{\algname{Local SGD} vs.\ \algname{Mini-Batch SGD}.} We assume that the statistical $\sigma^2/K$ dependence dominates the dependence on the initial distance $\sqn{x^0 - x^\ast} / K^2$. From Corollary~\ref{corrollary:sc-convergence-comm-complexity}, we see that in order to achieve the same convergence guarantees as \algname{Mini-Batch SGD}, we must have $H = \mathcal{O} \br{\frac{K}{\kappa M}}$, achieving a communication complexity of $\mathcal{O} \br{\kappa M}$. This is only possible when $K > \kappa M$. It follows that given a number of steps $K$ the optimal $H$ is $H = 1 + \floor{K/(\kappa M)}$ achieving a communication complexity of $\tilde{\Omega} \br{ \min (K, \kappa M) }$. 

\textbf{One-shot averaging.} Putting $H = K + 1$ yields a convergence rate of $\mathcal{\tilde{O}}(\sigma^2 \kappa /(\mu^2 K))$, showing no linear speedup but showing convergence, which improves upon all previous work. However, we admit that simply using Jensen's inequality to bound the distance of the average iterate $\ecn{\hat{x}^K - x^\ast}$ would yield a better asymptotic convergence rate of $\mathcal{\tilde{O}}(\sigma^2 / (\mu^2 K))$. Under a Lipschitz Hessian assumption, Zhang et al.~\cite{Zhang13} show that one-shot averaging can attain a linear speedup in the number of nodes, so one may do analysis of \algname{Local SGD} under this additional assumption to try to remove this gap, but this is beyond the scope of our work.

Similar results can be obtained for weakly convex functions, as the next Theorem shows.

\begin{theorem}
     \label{thm:weakly-convex-thm}
     Suppose that Assumptions~\ref{asm:convexity-and-smoothness}, \ref{asm:uniformly-bounded-variance} hold with $\mu = 0$ and that a constant stepsize $\gamma$ such that $\gamma \geq 0$ and $\gamma \leq \frac{1}{4 L}$ is chosen and that Algorithm~\ref{alg:local_sgd} is run for identical data with $H \geq 1$ such that $\sup_{p} \abs{k_p - k_{p+1}} \leq H$, then for $\bar{x}^K = \frac{1}{K} \sum_{k=1}^{K} \hat{x}^k$,
     \begin{align}
         \label{eq:thm-wc-convergence-rate}
         \begin{split}
               \ec{f(\bar{x}^K) - f(x^\ast)} \leq \frac{2}{\gamma K} &\norm{x^0 - x^\ast}^2  + \frac{2 \gamma \sigma^2}{M} + 4 \gamma^2 L \sigma^2 \br{H - 1}.
         \end{split}
     \end{align}
\end{theorem}

Theorem~\ref{thm:weakly-convex-thm} essentially tells the same story as Theorem~\ref{thm:sc-convergence-theorem}: convergence of \algname{Local SGD} is the same as \algname{Mini-Batch SGD} up to an additive constant whose size can be controlled by controlling $H$.
 
\begin{corollary}
    \label{corollary:weakly-convex-convergence-iid}
     Assume that $K \geq M$. Choosing $\gamma = \frac{\sqrt{M}}{4 L \sqrt{K}}$, then substituting in \eqref{eq:thm-wc-convergence-rate} we have,
     \[ 
     	\ec{f(\bar{x}^K) - f(x^\ast)} \leq \frac{8 L \sqn{x^0 - x^\ast}}{\sqrt{M K}} + \frac{\sigma^2}{2 L \sqrt{M K}} + \frac{\sigma^2 M \br{H - 1}}{L K}.
      \]
\end{corollary}

\textbf{Linear speedup and optimal $H$.} From Corollary~\ref{corollary:weakly-convex-convergence-iid} we see that if we choose $H = \mc{O}(\sqrt{K} M^{-3/2})$ then we obtain a linear speedup, and the number of communication steps is then $C = K/H = \Omega\br{M^{3/2} K^{1/2}}$, and we get that the optimal $H$ is then $H = 1 + \floor{ K^{1/2} M^{-3/2} }$.

The previous results were obtained under Assumption~\ref{asm:uniformly-bounded-variance}. Unfortunately, this assumption does not easily capture the finite-sum minimization scenario where $f(x) = \frac{1}{n} \sum_{i=1}^{n} f_i (x)$ and each stochastic gradient $g^k_m$ is sampled uniformly at random from the sum. 

Using smaller stepsizes and more involved proof techniques, we can show that our results still hold in the finite-sum setting. For strongly-convex functions, the next theorem shows that the same convergence guarantee as Theorem~\ref{thm:sc-convergence-theorem} can be attained.
\begin{figure*}[t]
\centering
	\includegraphics[scale=0.33]{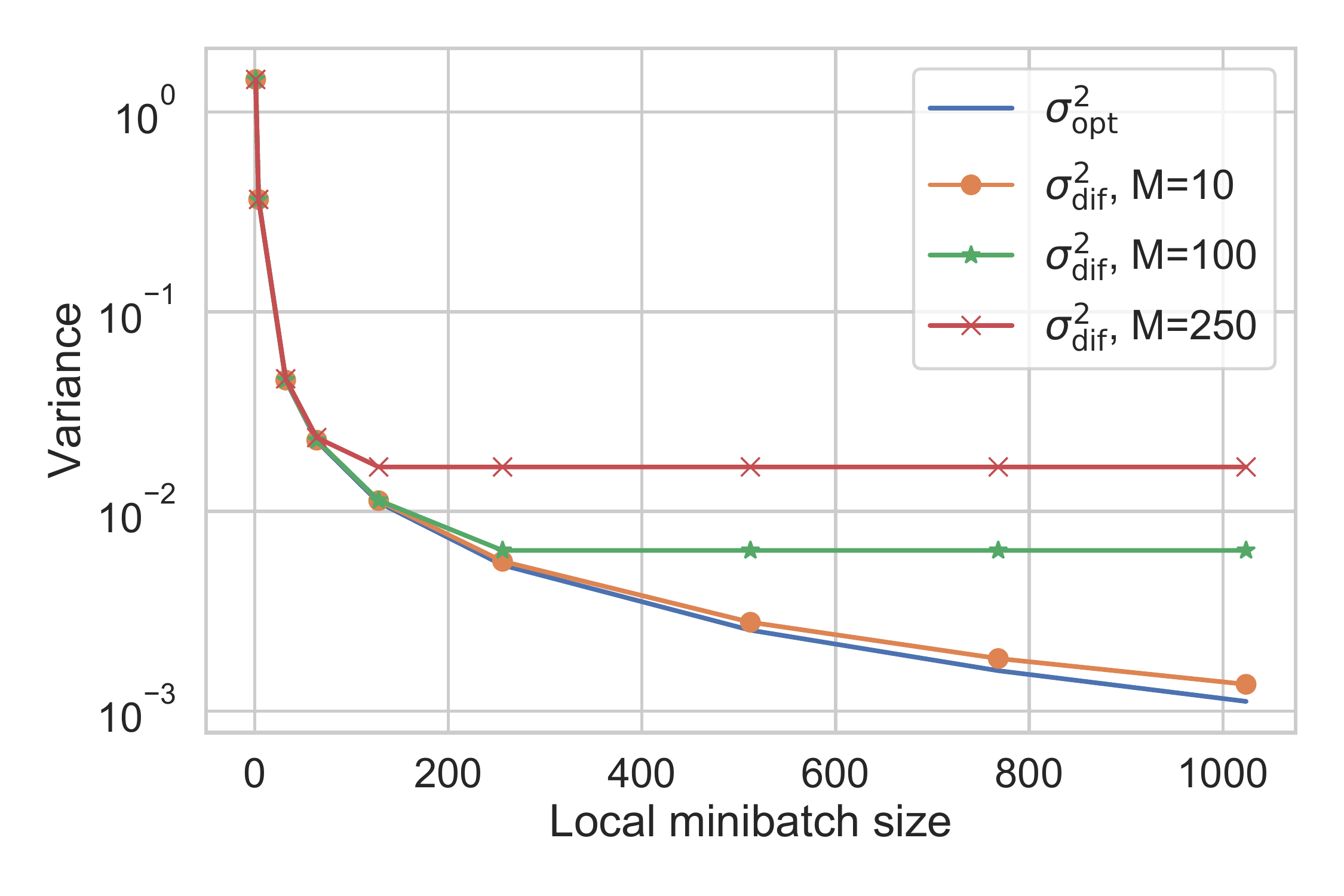}
	\includegraphics[scale=0.33]{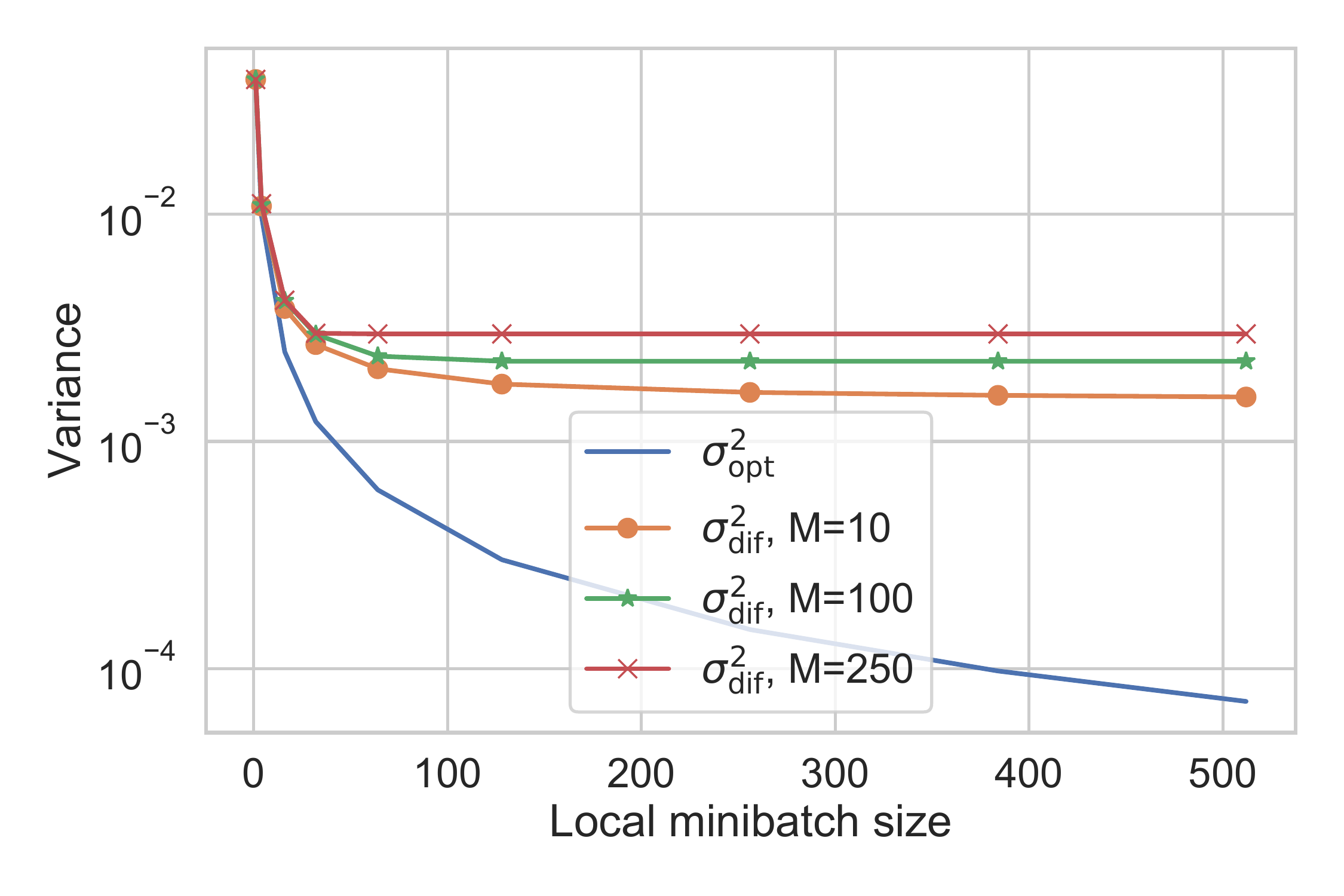}
	\includegraphics[scale=0.33]{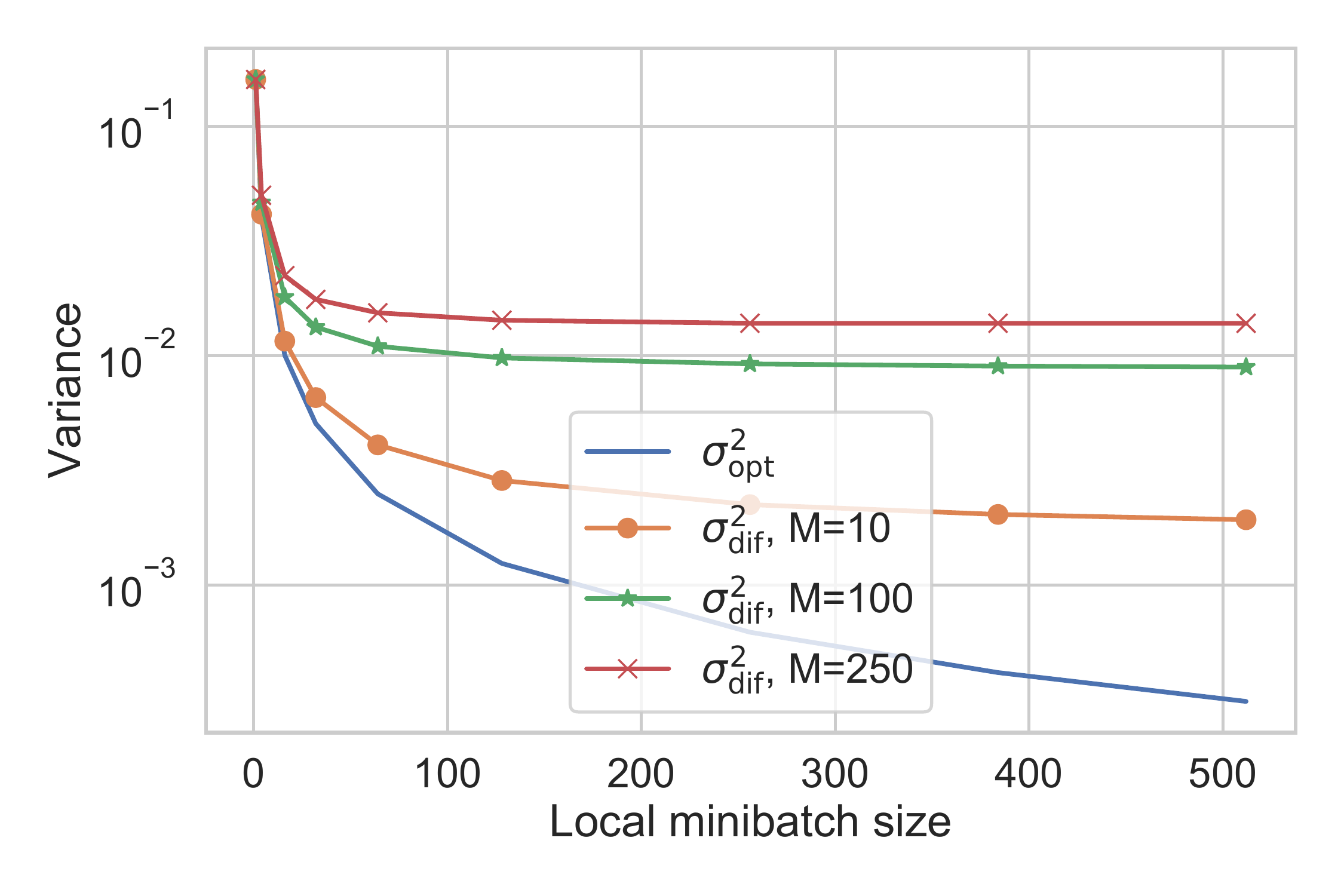}
	\caption{The effect of the dataset and number of workers $M$ on the variance parameters. Left: `a8a', middle: `mushrooms', right: `w8a' dataset. We use uniform sampling of data points, so $\sigmaopt^2$ is the same as $\sigmaf^2$ with $M=1$, while for higher values of $M$ the value of $\sigmaf^2$ might be drastically larger than $\sigmaopt^2$.}
	\label{fig:local_sgd_variance}
\end{figure*}
\begin{theorem}
    \label{theorem:sc-unbounded-variance-iid}
    Suppose that Assumptions~\ref{asm:convexity-and-smoothness} and \ref{asm:finite-sum-stochastic-gradients} hold with $\mu > 0$. Suppose that Algorithm~\ref{alg:local_sgd} is run for identical data with $\max_{p} \abs{k_{p} - k_{p+1}} \leq H$ for some $H \geq 1$ and with a stepsize $\gamma > 0$ chosen such that 
    $ \gamma \leq \min \pbr{\frac{1}{4L \br{1 + \frac{2}{M}}}, \frac{1}{\mu + 8 L \br{H-1}} }. $
    Then for any timestep $k$ such that synchronization occurs,
    \begin{align}
        \label{eq:thm-sc-const-gamma-unbounded-variance}
        \begin{split}
            \ecn{\hat{x}^k - x^\ast} &\leq \br{1 - \gamma \mu}^{k} \sqn{x^0 - x^\ast} + \frac{2 \gamma \sigmaopt^2}{\mu M} + \frac{4 \sigmaopt^2 \gamma^2 \br{H - 1} L}{\mu}.
        \end{split}
    \end{align}
\end{theorem}

As a corollary, we can obtain an asymptotic convergence rate by choosing specific stepsizes $\gamma$ and $H$.

\begin{corollary}
    \label{corollary:unbounded-var-convergence}
    Let $a = 18 \kappa k$ for some $k > 0$, let $H \leq k$ and choose $\gamma = \frac{1}{\mu a} \leq \frac{1}{9 L H}$. We substitute in \eqref{eq:thm-sc-const-gamma-unbounded-variance} and take $K = 18 a \log a$ steps, then
    \begin{align*}
           \ecn{\hat{x}^K - x^\ast} = \mathcal{\tilde{O}} \Bigg ( \frac{\sqn{x^0 - x^\ast}}{K^2} &+ \frac{\sigmaopt^2}{\mu^2 M K} + \frac{\sigmaopt^2 \kappa (H-1)}{\mu^2 K^2} \Bigg ).
    \end{align*}
\end{corollary}

Substituting $H = 1 + \floor{t/M} = 1 + \floor{K/(18 \kappa M)}$ in Corollary~\ref{corollary:unbounded-var-convergence} we get an asymptotic convergence rate of $\mc{\tilde{O}} \br{ \frac{\sigmaopt^2}{K M} }$. This preserves the rate of \algname{Mini-Batch SGD} up to problem-independent constants and polylogarithmic factors, but with possibly fewer communication steps.

\begin{theorem}
    \label{thm:wc-iid-unbounded-var}
    Suppose that Assumptions~\ref{asm:convexity-and-smoothness} and \ref{asm:finite-sum-stochastic-gradients} hold with $\mu = 0$, that a stepsize $\gamma \leq \frac{1}{10 L H}$ is chosen and that Algorithm~\ref{alg:local_sgd} is run on $M \geq 2$ nodes with identical data and with $\sup_{p} \abs{k_p - k_{p+1}} \leq H$, then for any timestep $K$ such that synchronization occurs we have for $\bar{x}^K = \frac{1}{K} \sum_{k=1}^{K} \hat{x}^k$ that
    \begin{equation}
        \label{eq:thm-weakly-convex-case-iid-unbounded-var}
            \ec{ f(\bar{x}^K) - f(x^\ast) } \leq \frac{10 \sqn{x^0 - x^\ast} }{\gamma K} + \frac{20 \gamma \sigmaopt^2}{M} + 40 \gamma^2 L \sigmaopt^2 \br{H-1}.
    \end{equation}
\end{theorem}

\begin{corollary}
    \label{corollary:wc-iid-conv-unbounded-var}
    Let $H \leq \frac{\sqrt{K}}{\sqrt{M}}$, then for $\gamma = \frac{\sqrt{M}}{10 L \sqrt{K}}$ we see that $\gamma \leq \frac{1}{10 L H}$, and plugging it into \eqref{eq:thm-weakly-convex-case-iid-unbounded-var} yields 
    \[
            \ec{f(\bar{x}^{K}) - f(x^\ast)} \leq \frac{100 L \sqn{x^0 - x^\ast}}{\sqrt{K M}} + \frac{2 \sigmaopt^2}{L \sqrt{K M}} \frac{2 \sigmaopt^2 M (H-1)}{5 L K}.
    \]
\end{corollary}

This is the same result as Corollary~\ref{corollary:weakly-convex-convergence-iid}, and hence we see that choosing $H = \mathcal{O}\br{K^{1/2} M^{-3/2}}$ (when $K > M^3$) yields a linear speedup in the number of nodes $M$.

\subsection{Heterogeneous data}
We next show that for arbitrarily heterogeneous convex objectives, the convergence of \algname{Local SGD} is the same as \algname{Mini-Batch SGD} plus an error that depends on $H$.

\begin{theorem}
    \label{thm:wc-noniid-unbounded-var}
    Suppose that Assumptions~\ref{asm:convexity-and-smoothness} and \ref{asm:finite-sum-stochastic-gradients} hold with $\mu = 0$ and for heterogeneous data. Then for Algorithm~\ref{alg:local_sgd} run for different data with $M \geq 2$, $\max_{p} \abs{k_p - k_{p+1}} \leq H$, and a stepsize $\gamma > 0$ such that $\gamma \leq \min \pbr{ \frac{1}{4L}, \frac{1}{8 L (H-1)}}$, then we have
    \begin{align*}
        \begin{split}
            \ec{f(\bar{x}^K) -f(x^\ast)} \leq &\frac{4 \sqn{x^0 - x^\ast}}{\gamma K} + \frac{20 \gamma \sigmaf^2}{M} + 16 \gamma^2 L (H-1)^2 \sigmaf^2.
        \end{split}
    \end{align*}
    where $\bar{x}^K \eqdef \frac{1}{K} \sum_{i=0}^{K-1} \hat{x}_{i}$.
\end{theorem}

\textbf{Dependence on $\sigmaf$.} We see that the convergence guarantee given by Theorem~\ref{thm:wc-noniid-unbounded-var} shows a dependence on $\sigmaf$, which measures the heterogeneity of the data distribution. In typical (non-federated) distributed learning settings where data are distributed before starting training, this term can vary quite significantly depending on how the data are distributed.

\textbf{Dependence on $H$.} We further note that the dependence on $H$ in Theorem~\ref{thm:wc-noniid-unbounded-var} is quadratic rather than linear. This translates to a worse upper bound on the synchronization interval $H$ that still allows convergence, as the next corollary shows.

\begin{figure*}[t]
	\centering
	\includegraphics[scale=0.33]{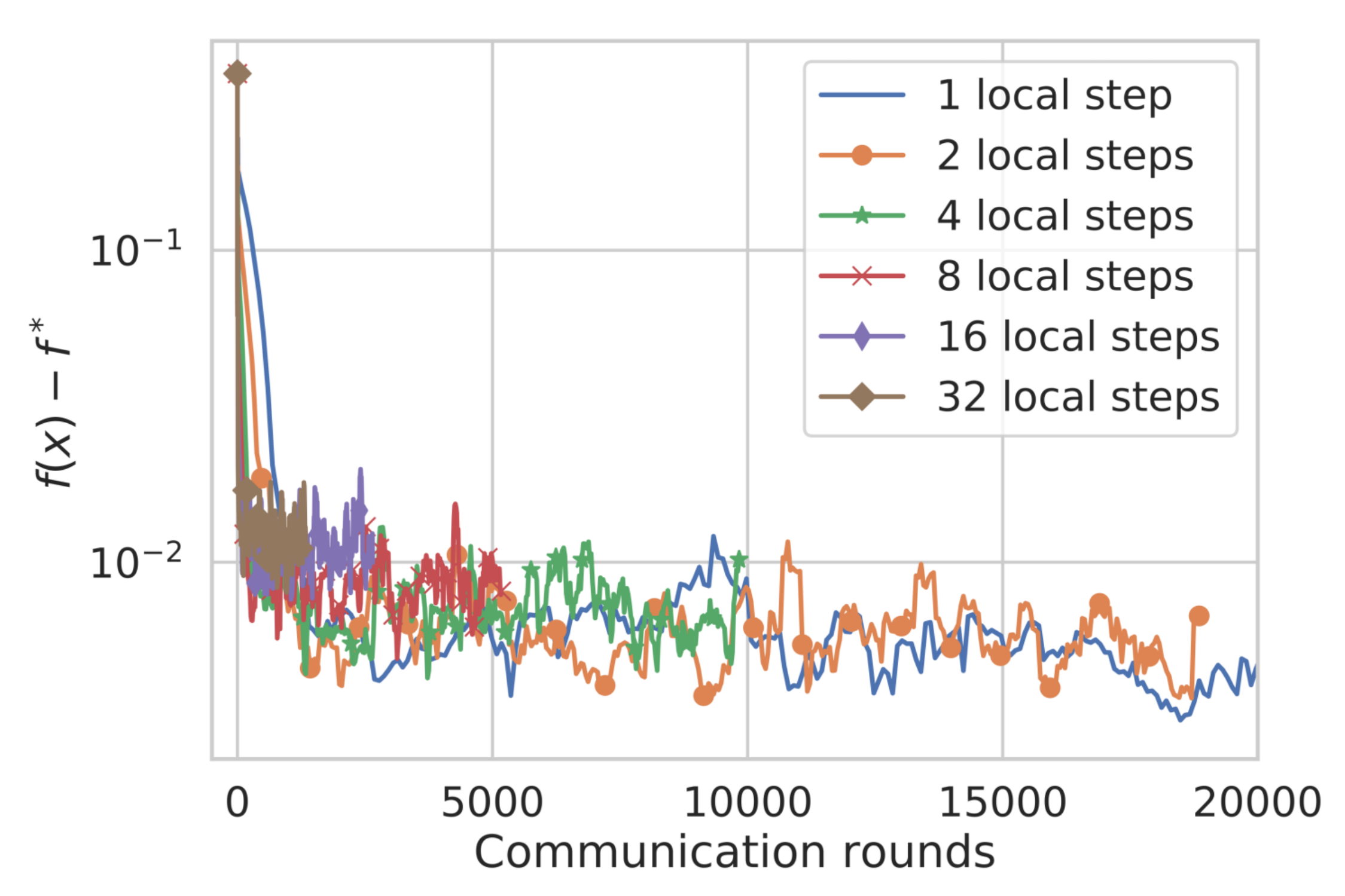}
	\includegraphics[scale=0.33]{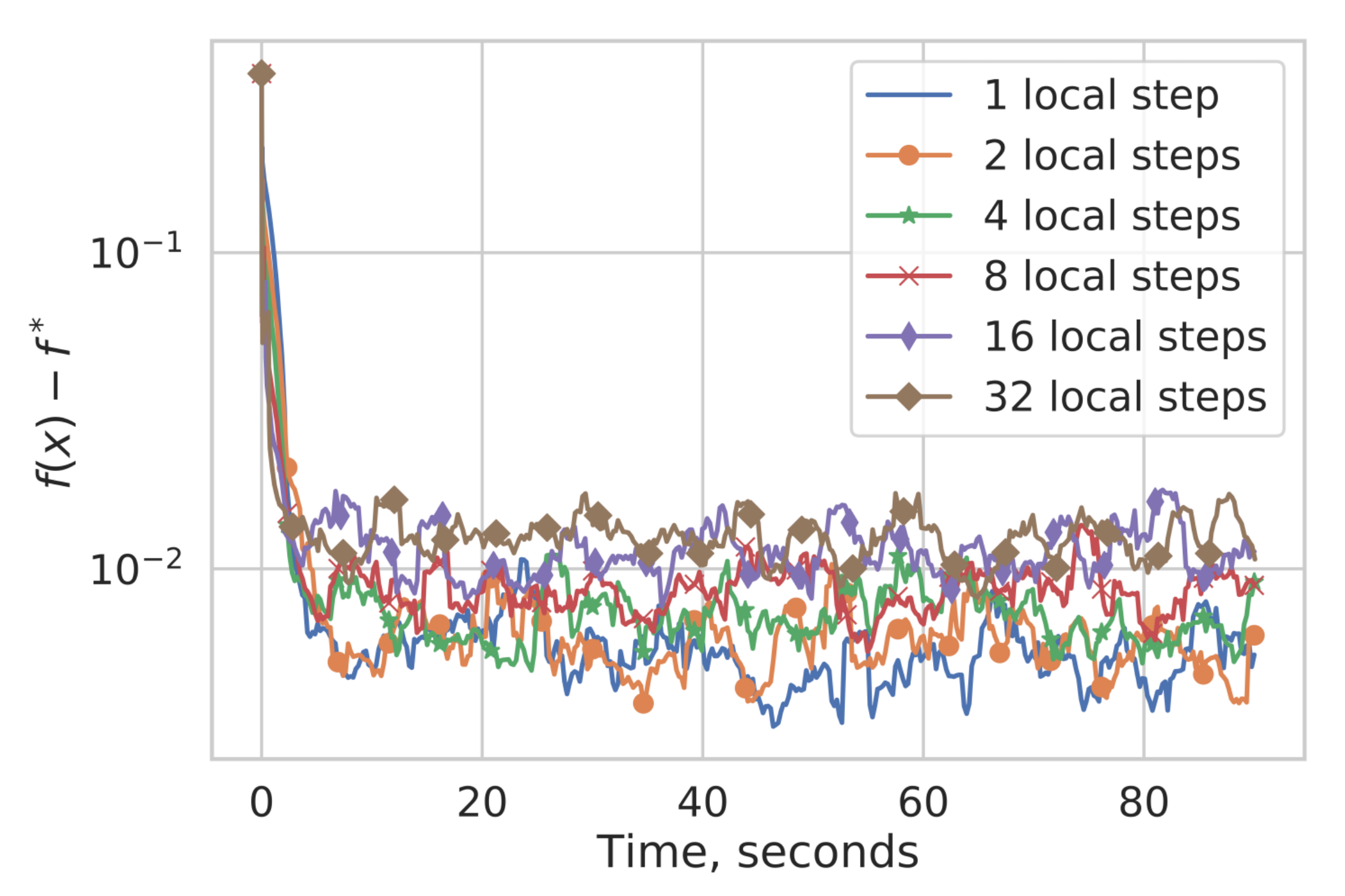}
	\caption{Results on `a9a' dataset, with stepsize $\frac{1}{L}$. For any value of local iterations $H$ the method converged to a neighborhood within a small number of communication rounds due to large stepsizes.}
	\label{fig:a5a_same_data_01}
\end{figure*}
\begin{figure*}[!t]
	\centering
	\includegraphics[scale=0.33]{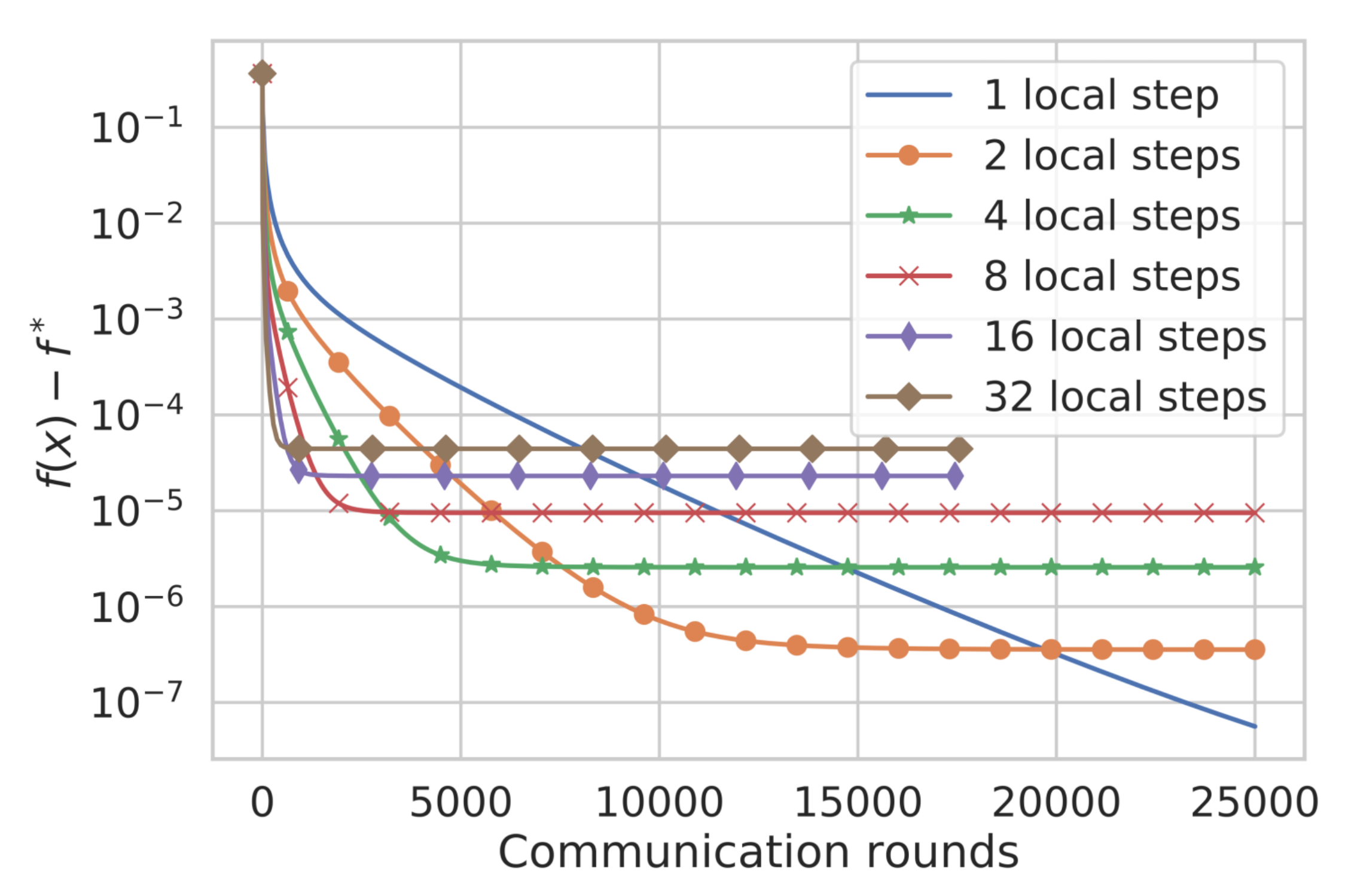}
	\includegraphics[scale=0.33]{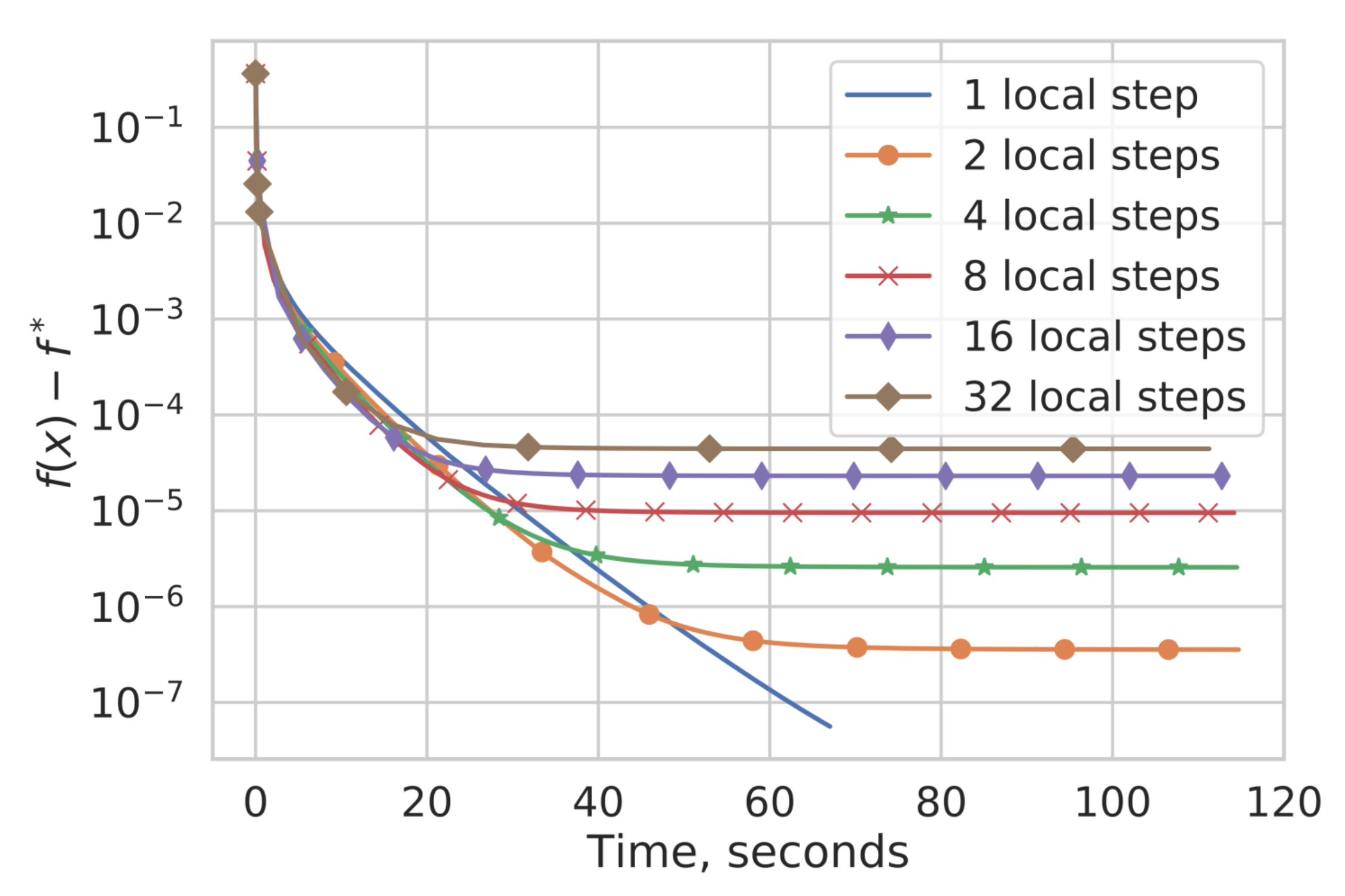}
	\includegraphics[scale=0.33]{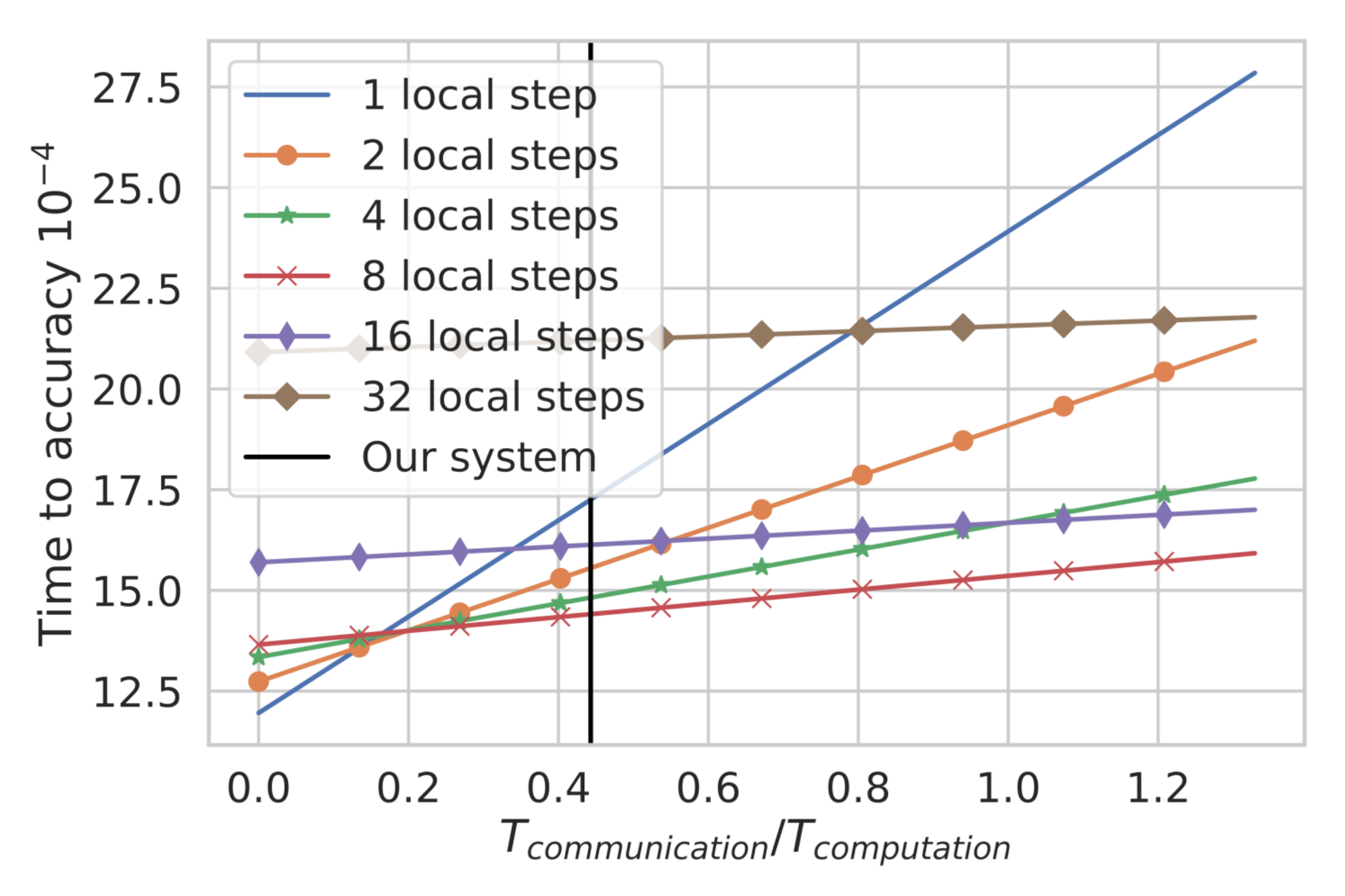}
	\caption{Convergence on heterogeneous data with different number of local steps on the `a5a' dataset. 1 local step corresponds to fully synchronized gradient descent. Left: convergence in terms of communication rounds, which shows a clear advantage of \algname{Local GD} when only limited accuracy is required. Mid plot: wall-clock time might improve only slightly if communication is cheap. Right: what changes with different communication cost.}
	\label{fig:a5a_different_H}
\end{figure*}

\begin{corollary}
    \label{corollary:wc-noniid-unbounded-var}
    Choose $H \leq \frac{\sqrt{K}}{\sqrt{M}}$, then $\gamma = \frac{\sqrt{M}}{8 L \sqrt{K}} \leq \frac{1}{8 H L}$, and hence applying the result of Theorem~\ref{thm:wc-noniid-unbounded-var},
    \[
           \ec{f (\bar{x}^K) - f (x^\ast)} \leq \frac{32 L \sqn{x^0 - x^\ast}}{\sqrt{M K}} + \frac{5 \sigmaf^2}{2 L \sqrt{MK}} + \frac{\sigmaf^2 M (H-1)^2}{4 L K}.
    \]
\end{corollary}

\textbf{Optimal $H$.} By Corollary~\ref{corollary:wc-noniid-unbounded-var} we see that the optimal value of $H$ is $H = 1 + \floor{K^{1/4} M^{-3/4}}$, which gives $\mathcal{O}\br{\frac{1}{\sqrt{MK}}}$ convergence rate. Thus, the same convergence rate is attained provided that communication is more frequent compared to the identical data regime. 

\section{Experiments}
All experiments described below were run on logistic regression problem with $\ell_2$ regularization of order $\frac{1}{n}$. The datasets were taken from the LIBSVM library~\cite{chang2011libsvm}. The code was written in Python using MPI~\cite{dalcin2011parallel} and run on Intel(R) Xeon(R) Gold 6146 CPU @3.20GHz cores in parallel.
\subsection{Variance measures}
We provide values of $\sigmaf^2$ and $\sigmaopt^2$ in Figure~\ref{fig:local_sgd_variance} for different datasets, mini-batch sizes and $M$. The datasets were split evenly without any data reshuffling and no overlaps. For any $M>1$, the value of $\sigmaf$ is lower bounded by $\avemm \|\nabla f_m(x^\ast)\|^2$ which explains the difference between identical and heterogeneous data. 

\subsection{Identical data}
For identical data we used $M=20$ nodes and 'a9a' dataset.
We estimated $L$ numerically and ran two experiments, with stepsizes $\frac{1}{L}$ and $\frac{0.05}{L}$ and minibatch size equal 1. In both cases we observe convergence to a neighborhood, although of a different radius. Since we run the experiments on a single machine, the communication is very cheap and there is little gain in time required for convergence. However, the advantage in terms of required communication rounds is self-evident and can lead to significant time improvement under slow communication networks. The results are provided here in Figure~\ref{fig:a5a_same_data_01} and in the supplementary material in Figure~\ref{fig:a5a_same_data_001}.

\subsection{Heterogeneous data}

\begin{figure}[t]
	\centering
	\includegraphics[scale=0.33]{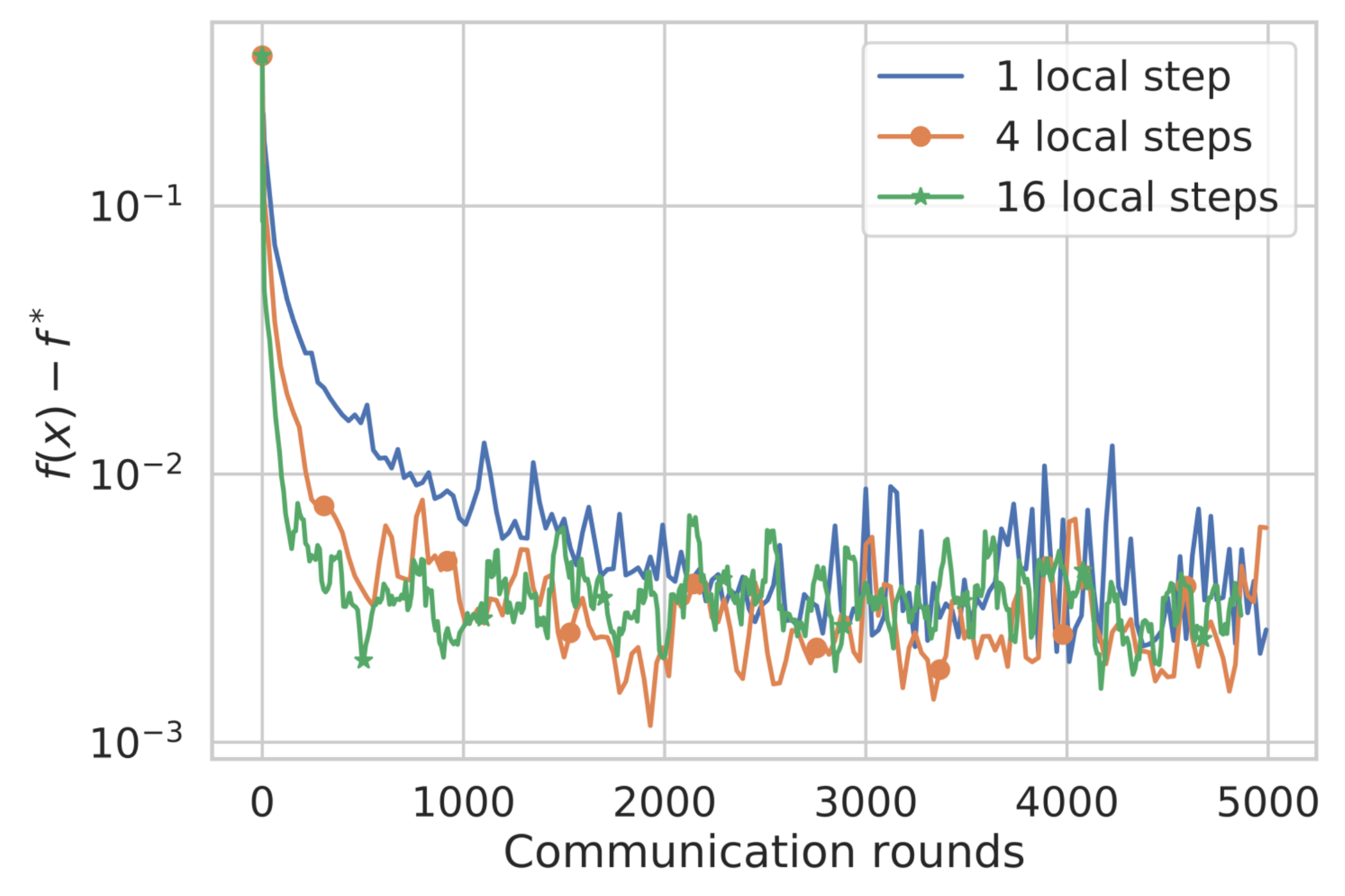}
	\label{fig:a5a_different_H2}
	\caption{Convergence of \algname{Local SGD} on heterogeneous data with different number of local steps on the `a5a' dataset.}
\end{figure}
Since our architecture leads to a very specific trade-off between computation and communication, we  provide plots for the case the communication time relative to gradient computation time is higher or lower. To see the impact of $\sigmaf$, in all experiments we use full gradients $\nabla f_m$ and constant stepsize $\frac{1}{L}$. The data partitioning is not i.i.d.\ and is done based on the index in the original dataset. The results are provided in Figure~\ref{fig:a5a_different_H} and in the supplementary material in Figure~\ref{fig:mushrooms_different_H}.
In cases where communication is significantly more expensive than gradient computation, local methods are much faster for imprecise convergence.

\chapter{Convergence of Random Reshuffling}
\label{chapter:rr}

\graphicspath{{rr/}}

\section{Introduction}

We study the finite-sum minimization problem
\begin{equation}
  \label{eq:finite-sum-min}
  \min_{x \in \R^d} \Bigl [ f(x) = \frac{1}{n} \sum \limits_{i=1}^{n} f_{i} (x) \Bigr ],
\end{equation}
where each $f_{i}\colon \R^d \to \R$ is differentiable and smooth, and we  are particularly interested in the big data machine learning setting where  the number of functions $n$ is large. Thanks to their scalability and low memory requirements, first-order methods are especially popular in this setting \cite{Bottou18}. Stochastic first-order algorithms in particular have attracted a lot of attention in the machine learning community and are often used in combination with various practical heuristics. Explaining these heuristics may lead to further development of stable and efficient training algorithms. In this chapter, we aim at better and sharper theoretical explanation of one intriguingly simple but notoriously elusive heuristic: data permutation/shuffling.


In particular, the goal of this chapter is to obtain deeper theoretical understanding of  methods for solving  \eqref{eq:finite-sum-min} which rely on random or deterministic permutation/shuffling of the data  $\{1,2,\dots,n\}$ and perform incremental  gradient updates following the permuted order.  We study three methods which belong to this class, described next.

An immensely popular but theoretically elusive  method  belonging to the class of data permutation methods  is the {\bf \algname{Random Reshuffling} (\algname{RR})} algorithm (see Algorithm~\ref{alg:rr}). This is the method we pay most attention to in this chapter. In each epoch $t$ of \algname{RR}, we sample indices $\pi_{0}, \pi_{1}, \ldots, \pi_{n-1}$ without replacement from $\{ 1, 2, \ldots, n \}$, i.e., $\{\pi_{0}, \pi_{1}, \ldots, \pi_{n-1}\}$ is a random permutation of the set $\{ 1, 2, \ldots, n \}$, and proceed with $n$ iterates of the form
\[ x^{k}_{i+1} = x^k_i - \gamma \nabla f_{\pi_{i}} (x^k_i), \]
where $\gamma > 0$ is a stepsize. We then set $x^{k+1} = x^k_n$,  and repeat the process for a total of $T$ epochs. Notice that in  \algname{RR}, a new permutation/shuffling is generated at the beginning of each epoch, which is why the term reshuffling is used.

Furthermore, we consider the {\bf \algname{Shuffle-Once} (\algname{SO})} algorithm, which is identical to \algname{RR} with the exception that it shuffles the dataset only once---at the very beginning---and then reuses this random permutation in all subsequent epochs (see Algorithm~\ref{alg:so}). Our results for \algname{SO} follow as corollaries of the tools we developed in order to conduct a sharp analysis of \algname{RR}.

Finally, we also consider the {\bf \algname{Incremental Gradient} (\algname{IG})} algorithm, which is identical to \algname{SO}, with the exception that the initial permutation is not random but deterministic. Hence, \algname{IG} performs incremental gradient steps through the data in a cycling fashion. The ordering could be arbitrary, e.g., it could be  selected implicitly by the ordering the data comes in, or chosen adversarially. Again, our results for \algname{IG} follow as a byproduct of our efforts to understand \algname{RR}.

\begin{minipage}[t]{0.485\textwidth}
  \begin{algorithm}[H]
    \caption{\algname{Random Reshuffling} \rrbox{(\algname{RR})}}
    \label{alg:rr}
 \begin{algorithmic}[1]
   \Require Stepsize $\gamma > 0$, initial vector $x^0 = x_0^0 \in \R^d$, number of epochs $T$
    \For{epochs $k=0,1,\dotsc,T-1$}
   \State {\color{mydarkred} Sample a permutation $\prm{0}, \prm{1}, \ldots, \prm{n-1}$ of $\{ 1, 2, \ldots, n\}$}
       \For{$i=0, 1, \ldots, n-1$}
          \State $x^{k}_{i+1} = x^k_{i} - \gamma \nabla f_{\prm{i}} (x^k_i)$
       \EndFor
       \State $x^{k+1} = x^{k}_{n}$
    \EndFor
 \end{algorithmic}
 \end{algorithm}
\end{minipage}
  \hfill
\begin{minipage}[t]{0.485\textwidth}
  \begin{algorithm}[H]
    \caption{\algname{Shuffle-Once} \sobox{(\algname{SO})}}
    \label{alg:so}
 \begin{algorithmic}[1]
   \Require Stepsize $\gamma > 0$, initial vector $x^0 = x_0^0 \in \R^d$, number of epochs $T$
   \State {\color{mydarkred} Sample a permutation $\prm{0}, \prm{1}, \ldots, \prm{n-1}$ of $\{ 1, 2, \ldots, n \}$}
    \For{epochs $k=0,1,\dotsc,T-1$}
       \For{$i=0, 1, \ldots, n-1$}
          \State $x^{k}_{i+1} = x^k_{i} - \gamma \nabla f_{\prm{i}} (x^k_i)$
       \EndFor
       \State $x^{k+1} = x^{k}_{n}$
    \EndFor
 \end{algorithmic}
 \end{algorithm}
\end{minipage}

\section{Related Work}

\algname{RR} is usually contrasted with its better-studied sibling \algname{Stochastic Gradient Descent} (\algname{SGD}), in which each $\pi_{i}$ is sampled uniformly with replacement from $\{ 1, 2, \ldots, n \}$. \algname{RR} often converges faster than \algname{SGD} on many practical problems~\cite{Bottou2009,Recht2013}, is more friendly to cache locality~\cite{Bengio2012}, and is in fact standard in deep learning~\cite{Sun2020}.

The convergence properties of \algname{SGD} are well-understood, with tightly matching lower and upper bounds in many settings~\cite{Rakhlin2012,Drori2019,HaNguyen2019}. Sampling without replacement allows \algname{RR} to leverage the finite-sum structure of \eqref{eq:finite-sum-min} by ensuring that each function contributes to the solution once per epoch. On the other hand, it also introduces a significant complication: the steps are now biased. Indeed, in any iteration $i>0$ within an epoch, we face the challenge of not having (conditionally) unbiased gradients since
\[ \ec{\nabla f_{\pi_{i}} (x^k_i) \mid x^k_i} \neq \nabla f(x^k_i). \]
This bias implies that individual iterations do not necessarily approximate a full gradient descent step. Hence, in order to obtain meaningful convergence rates for \algname{RR}, it is necessary to resort to more involved proof techniques. In recent work, various convergence rates have been established for \algname{RR}. However, a satisfactory, let alone complete, understanding of the algorithm's convergence remains elusive. For instance, the early line of attack pioneered by Recht and R{\'e}~\cite{RR-conjecture2012} seems to have hit the wall as their noncommutative arithmetic-geometric mean conjecture is not true~\cite{lai2020recht}. The situation is even more pronounced with the \algname{SO} method, as Safran and Shamir~\cite{Safran2020good} point out that there are no convergence results specific for the method, and the only convergence rates for \algname{SO} follow by applying the worst-case bounds of \algname{IG}. Rajput et al.~\cite{Rajput2020} state that a common practical heuristic is to use methods like \algname{SO} that do not reshuffle the data every epoch. Indeed, they add that 
``current theoretical bounds are insufficient to explain this phenomenon, and a new theoretical breakthrough may be required to tackle it''.

\algname{IG} has a long history owing to its success in training neural networks~\cite{Luo1991, Grippo1994}, and its asymptotic convergence has been established early~\cite{Mangasarian1994,Bertsekas2000}. 
Several rates for non-smooth and smooth cases were established by \cite{Nedic2001, Li2019, Gurbuzbalaban2019IG, Ying2019} and \cite{Nguyen2020}.
Using \algname{IG} poses the challenge of choosing a specific permutation for cycling through the iterates, which Nedi{\'{c}} and Bertsekas~\cite{Nedic2001} note to be difficult. Bertsekas~\cite{Bertsekas2011} gives an example that highlights the susceptibility of \algname{IG} to bad orderings compared to \algname{RR}. Yet, thanks to G\"{u}rb\"{u}zbalaban et al.~\cite{Gurbuzbalaban2019RR} and Haochen and Sra~\cite{haochen2018random}, \algname{RR} is known to improve upon both \algname{SGD} and \algname{IG} for twice-smooth objectives. Nagaraj et al.~\cite{Nagaraj2019} also study convergence of \algname{RR} for smooth objectives, and Safran and Shamir~cite{Safran2020good} as well as Rajput et al.~{Rajput2020} give lower bounds for \algname{RR} and related methods.

\section{Settings and Contributions} \label{sec:contributions}

In this chapter, we study the convergence behavior of the data-permutation methods \algname{RR}, \algname{SO} and \algname{IG}. While existing proof techniques succeed in obtaining insightful bounds for \algname{RR} and \algname{IG}, they fail to fully capitalize on the intrinsic power reshuffling and shuffling offers, and are not applicable to \algname{SO} at all\footnote{As we have mentioned before, the best known bounds for \algname{SO} are those which apply to \algname{IG} also, which means that the randomness inherent in \algname{SO} is wholly ignored.}. Our proof techniques are dramatically novel, simple, more insightful, and lead to improved convergence results, all under weaker  assumptions on the objectives than prior work.

We will derive results for strongly convex, convex as well as non-convex  objectives. 
To compare between the performance of first-order methods, we define an $\e$-accurate solution as a point $\tilde{x} \in \R^d$ that satisfies (in expectation if $\tilde{x}$ is random)
\begin{align*}
  \norm{\nabla f(\tilde{x})} \leq \e, && \text { or } && \norm{\tilde{x} - x^\ast}^2 \leq \e, && \text { or } && f(\tilde{x}) - f(x^\ast) \leq \e
\end{align*}
for non-convex, strongly convex, and non-strongly convex objectives, respectively, and where $x^\ast$ is assumed to be a minimizer of $f$ if $f$ is convex. We then measure the performance of first-order methods by the number of individual gradients $\nabla f_{i} (\cdot)$ they access to reach an $\e$-accurate solution. 

Our first assumption is that the objective is bounded from below and smooth. This assumption is used in all of our results and is widely used in the literature.
\begin{assumption}
  \label{asm:f-smoothness}
The objective $f$ and the individual losses $f_{1}, \ldots, f_{n}$ are all $L$-smooth, i.e., their gradients are $L$-Lipschitz. Further, $f$ is lower bounded by some $f^\ast \in \R$. If $f$ is convex, we also assume the existence of a minimizer $x^\ast \in \R^d$.
\end{assumption}

Assumption~\ref{asm:f-smoothness} is necessary in order to obtain better convergence rates for \algname{RR} compared to \algname{SGD}, since without smoothness the \algname{SGD} rate is optimal and cannot be improved \cite{Nagaraj2019}.

\subsection{New and improved convergence rates for \algname{RR}, \algname{SO} and \algname{IG}} In Section~\ref{sec:conv-theory}, we analyze the \algname{RR} and \algname{SO} methods and present novel convergence rates for strongly convex, convex, and non-convex  smooth objectives. Our results for \algname{RR} are summarized in Table~\ref{tab:conv-rates}. 
  
  \begin{itemize}[leftmargin=0.15in,itemsep=0.01in,topsep=0pt]
    \item {\emone Strongly convex case.} If each $f_i$ is strongly convex, we introduce a {\emtwo new proof technique} for studying the convergence of \algname{RR}/\algname{SO}   that allows us to obtain a {\emtwo better dependence on problem constants}, such as the number of functions $n$ and the condition number $\kappa$, compared to prior work (see Table~\ref{tab:conv-rates}). Key to our results is a {\emtwo new notion of variance specific to \algname{RR}/\algname{SO}} (see Definition~\ref{def:bregman-div-noise}), which we argue explains the superior convergence of \algname{RR}/\algname{SO} compared to \algname{SGD} in many practical scenarios. Our result for \algname{SO} tightly {\emtwo matches the lower bound} of \cite{Safran2020good}. We prove similar results in the more general setting when each $f_{i}$ is convex and $f$ is strongly convex (see Theorem~\ref{thm:only-f-sc-rr-conv}), but in this case we are forced to use smaller stepsizes. 

    \item {\emone Convex case.}  For convex but not necessarily strongly convex objectives $f_{i}$, we give the first result showing that {\emtwo \algname{RR}/\algname{SO} can provably achieve better convergence than \algname{SGD}} for a large enough number of iterations. This holds even when comparing against results that assume second-order smoothness, like the result of \cite{haochen2018random}.
    
    \item {\emone Non-convex case.}  For non-convex objectives $f_{i}$, we obtain for \algname{RR} a  {\emtwo much better dependence on the number of functions $n$} compared to the prior work of \cite{Nguyen2020}. 
  \end{itemize}

Furthermore, in the appendix we formulate and prove convergence results for \algname{IG} for strongly convex objectives, convex, and non-convex objectives as well. The bounds are worse than \algname{RR} by a factor of $n$ in the noise/variance term, as \algname{IG} does not benefit from randomization. Our result for strongly convex objectives tightly matches the lower bound of \cite{Safran2020good} up to an extra iteration and logarithmic factors, and is the first result to tightly match this lower bound.

\subsection{More general assumptions on the function class} Previous non-asymptotic convergence analyses of \algname{RR} either obtain worse bounds that apply to \algname{IG}, e.g., \cite{Ying2019,Nguyen2020}, or depend crucially on the assumption that each $f_i$ is Lipschitz \cite{Nagaraj2019,haochen2018random,Ahn2020}. Unfortunately, requiring each $f_i$ to be Lipschitz contradicts strong convexity~\cite{nguyen2018sgd} and is furthermore not satisfied in least square regression, matrix factorization, or for neural networks with smooth activations. In contrast, our work is the first to show how to leverage randomization to obtain better rates for \algname{RR} without assuming each $f_i$ to be Lipschitz. In concurrent work, Ahn et al.~\cite{AhnYun2020} also obtain a result for non-convex objectives satisfying the Polyak-\L{}ojasiewicz inequality, a generalization of strong convexity. Their result holds without assuming bounded gradients or bounded variance, but unfortunately with a worse dependence on $\kappa$ and $n$ when specialized to $\mu$-strongly convex functions.
  \begin{itemize}[leftmargin=0.15in,itemsep=0.01in,topsep=0pt]
    \item {\emone Strongly convex and convex case.} For strongly convex and convex objectives {\emtwo we do not require any assumptions on the functions used beyond smoothness and convexity.} 
    \item {\emone Non-convex case.} For non-convex objectives we obtain our results under a significantly more general assumption than the bounded gradients assumptions employed in prior work. Our assumption is also provably satisfied when each function $f_i$ is lower bounded, and hence is {\emtwo not only more general but also a more realistic assumption to use.}
\end{itemize}

\begin{table}[]
  \begin{threeparttable}[b]
    \centering
    \caption{Number of individual gradient evaluations needed by \algname{RR} to reach an $\e$-accurate solution (defined in Section~\ref{sec:conv-theory}). Logarithmic factors and constants that are not related to the assumptions are ignored. For non-convex objectives, $A$ and $B$ are the constants given by Assumption~\ref{asm:2nd-moment}.}
    \label{tab:conv-rates}
    \begin{tabular}{@{}C{1cm}@{}C{1cm}@{}cccc@{}}
    \toprule
      \multicolumn{2}{c}{Assumptions} & $\mu$-Strongly & Non-strongly & \multirow[c]{2}{*}{Non-convex} & \multirow[c]{2}{*}{Citation}\\ \cmidrule(lr){1-2}
      N.L.\tnote{\color{red}(1)} & U.V.\tnote{\color{red}(2)} & convex & convex & & \\ \midrule
        \cmark & \cmark & {\footnotesize $\kappa^2 n  + \frac{\kappa n  \sigmaesc}{\mu\sqrt{\e}}$} & --  & -- & {\footnotesize \cite{Ying2019}} \\[5pt]
        \xmark & \xmark & {\footnotesize $\kappa^2 n + \frac{\kappa \sqrt{n}  G}{\mu \sqrt{\e}}$} & {\footnotesize $\frac{L D^2 }{\e} + \frac{G^2 D^2}{\e^2}$}\tnote{\color{red}(3)} & -- & {\footnotesize \cite{Nagaraj2019}} \\[5pt]
        \xmark & \xmark & -- & -- & {\footnotesize $\frac{ L n}{\e^2} + \frac{ L n G}{\e^3}$} & {\footnotesize \cite{Nguyen2020}} \\[5pt]
        \cmark & \cmark & {\footnotesize $\frac{ \kappa^2 n }{\sqrt{\mu \e}}  + \frac{ \kappa^2 n \sigmaesc}{\mu \sqrt{\e}}$} \tnote{\color{red}(4)} & -- & -- & {\footnotesize \cite{Nguyen2020}} \\[5pt]
        \xmark & \xmark & {\footnotesize$\frac{\kappa \alpha }{\e^{1/\alpha}} + \frac{\kappa  \sqrt{n} G \alpha^{3/2}}{\mu\sqrt{\e}}$}\tnote{\color{red}(5)} & -- & -- & {\footnotesize \cite{Ahn2020}} \\[5pt]
        \cmark & \cmark & { \footnotesize$ \kappa n + \frac{\sqrt{n}}{\sqrt{\mu\e}} + \frac{\kappa \sqrt{n} G_0}{\mu \sqrt{\e}}$ }\tnote{\color{red}(6)} & -- & -- & {\footnotesize \cite{AhnYun2020}} \\[5pt]
        \cmidrule(lr){1-5}
        \multirow[c]{2}{*}[-4pt]{\cmark} & \multirow[c]{2}{*}[-4pt]{\cmark} & {\footnotesize $\kappa + \frac{\sqrt{\kappa n} \sigmaesc}{\mu \sqrt{\e}}$}\tnote{\color{red}(7)} &  \multirow[c]{2}{*}{$\frac{ L n}{\e} + \frac{ \sqrt{L n} \sigmaesc}{\e^{3/2}}$} & { \footnotesize \multirow[c]{2}{*}{$\frac{ L n}{\e^2} + \frac{ L \sqrt{n} (B+\sqrt{A})}{\e^3}$}} & \multirow[c]{2}{*}[-4pt]{This thesis} \\[5pt]
        & & { \footnotesize $\kappa n + \frac{\sqrt{\kappa n} \sigmaesc}{\mu \sqrt{\e}}$ } &  &  & \\
        \bottomrule
    \end{tabular}
    \begin{tablenotes}
      {\footnotesize
        \item [{\color{red}(1)}] Support for non-Lipschitz functions (N.L.): proofs without assuming that $\max_{i=1,\dotsc,n}\norm{\nabla f_{i} (x)} \leq G$ for all $x \in \R^d$ and some $G>0$. Note that $ \frac{1}{n} \sum_{i=1}^{n} \sqn{\nabla f_{i} (x^\ast)} \eqdef \sigmaesc^2 \leq G^2$ and $B^2 \leq G^2$.
        \item [{\color{red}(2)}] Unbounded variance (U.V.): there may be no constant $\sigma$ such that Assumption~\ref{asm:2nd-moment} holds with $A = 0$ and $B = \sigma$. Note that when the individual gradients are bounded, the variance is automatically bounded too.
        \item [{\color{red}(3)}] Nagaraj et al.~\cite{Nagaraj2019} require, for non-strongly convex functions, projecting at each iteration onto a bounded convex set of diameter $D$. We study the unconstrained problem.
        \item [{\color{red}(4)}] For strongly convex, Nguyen et al.~\cite{Nguyen2020} bound $f(x)-f(x^\ast)$ rather than squared distances, hence we use strong convexity to translate their bound into a bound on $\|x-x^\ast\|^2$.
        \item [{\color{red}(5)}] The constant $\alpha > 2$ is a parameter to be specified in the stepsize used by Ahn et al.~\cite{Ahn2020}. Their full bound has several extra terms but we include only the most relevant ones. 
        \item [{\color{red}(6)}] The result of \cite{AhnYun2020} holds when $f$ satisfies the Polyak-{\L}ojasiewicz inequality, a generalization of strong convexity. We nevertheless specialize it to strong convexity for our comparison. The constant $G_0$ is defined as $G_0 \eqdef \sup_{x: f(x) \leq f(x^0)} \max_{i \in [n]} \norm{\nabla f_{i} (x)}$. Note that $\sigmaesc \leq G_0$. We show a better complexity for PL functions under bounded variance in Theorem~\ref{thm:rr-nonconvex}.
        \item [{\color{red}(7)}] This result is the first to show that \algname{RR} and \algname{SO} work with any $\gamma\le \frac{1}{L}$, but it asks for each $f_i$ to be strongly convex. The second result assumes that only $f$ is strongly convex.
      }
    \end{tablenotes}
  \end{threeparttable}
\end{table}

\section{Convergence Theory}
\label{sec:conv-theory}

 The following quantity is key to our analysis and serves as an asymmetric distance between two points measured in terms of functions.
\begin{definition}
	For any $i$, the quantity $D_{f_i}(x, y)\eqdef f_i(x) - f_i(y) - \<\nabla f_i(y), x - y>$ is the Bregman divergence between $x$ and $y$ associated with $f_i$.
\end{definition}
It is well-known that if $f_i$ is $L$-smooth and $\mu$-strongly convex, then for all $x, y \in \R^d$
\begin{equation}
  \label{eq:bregman-sc-smooth-properties}
   \frac{\mu}{2}\|x-y\|^2\le D_{f_i}(x, y)\le \frac{L}{2}\|x-y\|^2,
\end{equation}
so each Bregman divergence is closely related to the Euclidian distance. Moreover, the difference between the gradients of a convex and $L$-smooth $f_i$ is related to its Bregman divergence by
\begin{equation}
  \label{eq:bregman-sc-grad-smooth-properties}
  \sqn{\nabla f_i (x) - \nabla f_i (y)} \leq 2 L \cdot D_{f_i} (x, y).
\end{equation}

\subsection{Main result: strongly convex objectives}
\label{sec:strongly-convex}

Before we proceed to the formal statement of our main result, we need to present the central finding of this chapter. The analysis of many stochastic methods, including \algname{SGD}, rely on the fact that the iterates converge to $x^\ast$ up to some noise. This is exactly where we part ways with the standard analysis techniques, since, it turns out, the intermediate iterates of shuffling algorithms converge to some other points. Given a permutation $\pi$, the real limit points are defined below,
\begin{align}
	x^\ast_i \eqdef x^\ast - \gamma \sum \limits_{j=0}^{i-1} \nabla f_{\pi_{j}} (x^\ast), \qquad i=1,\dotsc, n-1. \label{eq:x_ast_i}
\end{align}
In fact, it is predicted by our theory and later validated by our experiments that within an epoch the iterates go away from $x^\ast$, and closer to the end of the epoch they make a sudden comeback to $x^\ast$.

The second reason the vectors introduced in Equation~\eqref{eq:x_ast_i} are so pivotal is that they allow us to define a new notion of variance. Without it, there seems to be no explanation for why \algname{RR} sometimes overtakes \algname{SGD} from the very beginning of optimization process. We define it below.
\begin{definition}[Shuffling variance]
  \label{def:bregman-div-noise}
  Given a stepsize $\gamma>0$ and a random permutation $\pi$ of $\{ 1, 2, \ldots, n \}$, define $x^\ast_i$ as in \eqref{eq:x_ast_i}. Then, the shuffling variance is given by
  \begin{equation}\label{eq:bregman-div-noise} \sigmass \eqdef \max \limits_{i=1, \ldots, n-1} \left [ \frac{1}{\gamma}\ec{D_{f_{\pi_{i}}} (x^\ast_i, x^\ast)} \right ], \end{equation}
  where the expectation is taken with respect to the randomness in the permutation $\pi$.
\end{definition}

Naturally, $\sigmass$ depends on the functions $f_1, \ldots, f_n$, but, unlike \algname{SGD}, it also depends in a non-trivial manner on the stepsize $\gamma$. The easiest way to understand the new notation is to compare it to the standard definition of variance used in the analysis of \algname{SGD}. We argue that $\sigmass$ is the natural counter-part for the standard variance used in \algname{SGD}. We relate both of them by the following upper and lower bounds:

\begin{proposition}
  \label{prop:shuffling-variance-normal-variance-bound}
  Suppose that each of $f_1, f_2, \ldots, f_{n}$ is $\mu$-strongly convex and $L$-smooth. Then
  $\frac{\gamma\mu n}{8}\sigmaesc^2
  \le \sigmass \leq \frac{\gamma L n}{4} \sigmaesc^2, $
  where $\sigmaesc^2 \eqdef \frac{1}{n} \sum_{i=1}^{n} \sqn{\nabla f_{i} (x^\ast)}$.
\end{proposition}

In practice, $\sigmass$ may be much closer to the lower bound than the upper bound; see Section~\ref{sec:experiments}. This leads to a dramatic difference in performance and provides additional evidence of the superiority of \algname{RR} over \algname{SGD}. The next theorem states how exactly convergence of \algname{RR} depends on the introduced variance.

\begin{theorem}
  \label{thm:all-sc-rr-conv}
  Suppose that the functions $f_1,\dotsc, f_n$ are $\mu$-strongly convex and that Assumption~\ref{asm:f-smoothness} holds. Then for Algorithms~\ref{alg:rr} or~\ref{alg:so} run with a constant stepsize $\gamma \leq \frac{1}{L}$, the iterates generated by either of the algorithms satisfy
  \[  \ecn{x^{T} - x^\ast} \leq \br{1 - \gamma \mu}^{nT} \sqn{x^0 - x^\ast} + \frac{2\gamma \sigmass}{\mu}. \]
  \vspace{-2em}
\end{theorem}
\begin{proof}[Proof]
  The key insight of our proof is that the intermediate iterates $x^k_1, x^k_2, \ldots$ do not converge to $x^\ast$, but rather converge to the sequence $x^\ast_1, x^\ast_2, \ldots$ defined by \eqref{eq:x_ast_i}. Keeping this intuition in mind, it makes sense to study the following recursion:
  \begin{align}
    &\ec{\|x^k_{i+1}-x^\ast_{i+1}\|^2} \notag\\
    &=\ec{\|x^k_{i}-x^\ast_{i}\|^2-2\gamma\<\nabla f_{\pi_i}(x^k_i)-\nabla f_{\pi_i}(x^\ast), x^k_i - x^\ast_i>+\gamma^2\|\nabla f_{\pi_i}(x^k_i) - \nabla f_{\pi_i}(x^\ast)\|^2}.\label{eq:big89fg9h9d}
  \end{align}
  Once we have this recursion, it is useful to notice that the scalar product can be decomposed as
  \begin{align}
    \<\nabla f_{\pi_i}(x^k_i)-\nabla f_{\pi_i}(x^\ast), x^k_i - x^\ast_i>
    &= [f_{\pi_i}(x^\ast_i)-f_{\pi_i}(x^k_i)-\<\nabla f_{\pi_i}(x^k_i), x^\ast_i-x^k_i>] \notag \\
    & \quad + [f_{\pi_i}(x^k_i)-f_{\pi_i}(x^\ast)-\<\nabla f_{\pi_i}(x^\ast), x^k_i-x^\ast>] \notag  \\
    & \quad - [f_{\pi_i}(x^\ast_i)-f_{\pi_i}(x^\ast)-\<\nabla f_{\pi_i}(x^\ast), x^\ast_i-x^\ast>] \notag \\
    &= D_{f_{\pi_i}}(x^\ast_i, x^k_i)+D_{f_{\pi_i}}(x^k_i, x^\ast) - D_{f_{\pi_i}}(x^\ast_i, x^\ast).\label{eq:bi89gfdb09hff}
  \end{align}
  This decomposition is, in fact, very standard and is a special case of the so-called three-point identity \cite{Chen1993}. So, it should not be surprising that we use it. \\
  The rest of the proof relies on obtaining appropriate  bounds for the terms in the recursion. Firstly, we bound each of the three Bregman divergence terms appearing in \eqref{eq:bi89gfdb09hff}. By $\mu$-strong convexity of $f_i$, the first term in \eqref{eq:bi89gfdb09hff} satisfies
  \[
      \frac{\mu}{2}\|x^k_i - x^\ast_i\|^2 \overset{\eqref{eq:bregman-sc-smooth-properties}}{\leq} D_{f_{\pi_i}}(x^\ast_i, x^k_i)		,
  \]
  which we will use to obtain contraction.	The second term in \eqref{eq:bi89gfdb09hff} can be bounded via
  \[
    \frac{1}{2L}\|\nabla f_{\pi_i}(x^k_i) - \nabla f_{\pi_i}(x^\ast)\|^2
    \overset{\eqref{eq:bregman-sc-grad-smooth-properties}}{\le} D_{f_{\pi_i}}(x^k_i, x^\ast),
  \]
  which gets absorbed in the last term in the expansion of $\|x^k_{i+1}-x^\ast_{i+1}\|^2$.   The expectation of the third divergence term in 
  \eqref{eq:bi89gfdb09hff}   is trivially bounded as follows:
  \[ \ec{D_{f_{\pi_{i}}} (x^\ast_i, x^\ast) } \leq \max_{i=1, \ldots, n-1} \left [ \ec{D_{f_{\pi_{i}}} (x^\ast_i, x^\ast)} \right ] = \gamma\sigmass. \]
  Plugging these three bounds back into \eqref{eq:bi89gfdb09hff}, and the resulting inequality into \eqref{eq:big89fg9h9d}, we obtain
  \begin{align}
    \ec{\|x^k_{i+1}-x^\ast_{i+1}\|^2}
    &\le \ec{(1-\gamma\mu)\|x^k_i-x^\ast_i\|^2 - 2\gamma(1-\gamma L)D_{f_{\pi_i}}(x^k_i, x^\ast)} + 2 \gamma^2 \sigmass \nonumber \\
    \label{eq:thm_str_cvx_main-proof-1}
    &\le (1-\gamma\mu)\ec{\|x^k_i-x^\ast_i\|^2}+ 2 \gamma^2 \sigmass.
  \end{align}
  The rest of the proof is just solving this recursion, and is relegated to Section~\ref{sec:proof-of-thm-1} in the appendix.
\end{proof}

We show (\Cref{corr:all-sc-rr-conv} in the appendix) that by carefully controlling the stepsize, the final iterate of \algname{RR} after $T$ epochs satisfies
\begin{equation}
  \label{eq:all-sc-rr-conv-bound-T}
  \ecn{x^{T} - x^\ast} = \mathcal{\tilde{O}} \br{ \exp\left(-\frac{\mu n T}{L} \right) \sqn{x^0 - x^\ast} + \frac{\kappa \sigmaesc^2}{\mu^2 n T^2} },
\end{equation}
where the $\tilde{\cO}(\cdot)$ notation suppresses absolute constants and polylogarithmic factors. Note that Theorem~\ref{thm:all-sc-rr-conv} covers both \algname{RR} and \algname{SO}, and for \algname{SO}, \cite{Safran2020good} give almost the same lower bound. Stated in terms of the squared distance from the optimum, their lower bound is 
\[ \ecn{x^{T} - x^\ast} = \Omega\br{\min \pbr{ 1, \frac{\sigmaesc^2}{\mu^2 n T^2} }}, \]
where we note that in their problem $\kappa = 1$. This translates to sample complexity $\mathcal{O}\br{ \sqrt{n} \sigmaesc/(\mu \sqrt{\e}) }$ for $\e \leq 1$\footnote{In their problem, the initialization point $x^0$ satisfies $\sqn{x^0 - x^\ast} \leq 1$ and hence asking for accuracy $\e > 1$ does not make sense.}. Specializing $\kappa = 1$ in \Cref{eq:all-sc-rr-conv-bound-T} gives the sample complexity of $\mathcal{\tilde{O}}\br{1 + \sqrt{n} \sigmaesc/(\mu \sqrt{\e}) }$, matching the optimal rate up to an extra iteration. More recently, Rajput et al.~\cite{Rajput2020} also proved a similar lower bound for \algname{RR}. We emphasize that Theorem~\ref{thm:all-sc-rr-conv} is not only tight, but it is also the first convergence bound that applies to \algname{SO}. Moreover, it also immediately works if one permutes once every few epochs, which interpolates between \algname{RR} and \algname{SO} mentioned by Rajput et al.~\cite{Rajput2020}.

\textbf{Comparison with \algname{SGD}} To understand when \algname{RR} is better than \algname{SGD}, let us borrow a convergence bound for the latter. Several works have shown (e.g., see \cite{needell2014stochastic, stich2019unified}) that for any $\gamma\le \frac{1}{2L}$ the iterates of \algname{SGD} satisfy
\[
	 \ecn{x^{nT}_{\mathrm{SGD}} - x^\ast} \leq \br{1 - \gamma \mu}^{nT} \sqn{x^0 - x^\ast} + \frac{ 2\gamma \sigmaesc^2}{\mu}.
\]
Thus, the question as to which method will be faster boils down to which variance is smaller: $\sigmass$ or $\sigmaesc^2$. According to Proposition~\ref{prop:shuffling-variance-normal-variance-bound}, it depends on both $n$ and the stepsize. Once the stepsize is sufficiently small, $\sigmass$ becomes smaller than $\sigmaesc^2$, but this might not be true in general. Similarly, if we partition $n$ functions into $n/\tau$ groups, i.e., use mini-batches of size $\tau$, then $\sigmaesc^2$ decreases as $\cO\br{1/\tau}$ and $\sigmass$ as $\cO\br{1/\tau^2}$, so \algname{RR} can become faster even without decreasing the stepsize. We illustrate this later with numerical experiments.

While \Cref{thm:all-sc-rr-conv} requires each $f_i$ to  be strongly convex, we can also obtain results in the case where the individual strong convexity  assumption is replaced by convexity. However, in such a case, we need to use a smaller stepsize, as the next theorem shows.
\begin{theorem}
  \label{thm:only-f-sc-rr-conv}
  Suppose that each $f_i$ is convex, $f$ is $\mu$-strongly convex, and Assumption~\ref{asm:f-smoothness} holds. Then provided the stepsize satisfies $\gamma \leq \frac{1}{\sqrt{2} L n}$ the final iterate generated by Algorithms~\ref{alg:rr} or~\ref{alg:so} satisfies
  \[ \ecn{x^T - x^\ast} \leq \br{ 1 - \frac{\gamma \mu n}{2} }^T \sqn{x^0 - x^\ast} + \gamma^2 \kappa n \sigmaesc^2. \]
\end{theorem}
It is not difficult to show that by properly choosing the stepsize $\gamma$, the guarantee given by Theorem~\ref{thm:only-f-sc-rr-conv} translates to a sample complexity of $\ctO\br{ \kappa n + \frac{\sqrt{\kappa n} \sigmaesc}{\mu \sqrt{\e}} }$, which matches the dependence on the accuracy $\e$ in Theorem~\ref{thm:all-sc-rr-conv} but with $\kappa (n-1)$ additional iterations in the beginning. For $\kappa = 1$, this translates to a sample complexity of $\ctO\br{n + \frac{\sqrt{n} \sigmaesc}{\mu \sqrt{\e}}}$ which is worse than the lower bound of \cite{Safran2020good} when $\e$ is large. In concurrent work, Ahn et al.~\cite{AhnYun2020} obtain in the same setting a complexity of $\ctO\br{1/\e^{1/\alpha} + \frac{\sqrt{n} G}{\mu \sqrt{\e}}}$ (for a constant $\alpha > 2$), which requires that each $f_i$ is Lipschitz and matches the lower bound only when the accuracy $\e$ is large enough that $1/\e^{1/\alpha} \leq 1$. Obtaining an optimal convergence guarantee for all accuracies $\e$ in the setting of Theorem~\ref{thm:only-f-sc-rr-conv} remains open.

\subsection{Non-strongly convex objectives}
\label{sec:weakly-convex}
We also make a step towards better bounds for \algname{RR}/\algname{SO} without any strong convexity at all and provide the following convergence statement.
\begin{theorem}
  \label{thm:weakly-convex-f-rr-conv}
  Let functions $f_1, f_2, \ldots, f_n$ be convex. Suppose that Assumption~\ref{asm:f-smoothness} holds. Then for Algorithm~\ref{alg:rr} or Algorithm~\ref{alg:so} run with a stepsize $\gamma \leq \frac{1}{\sqrt{2} L n}$, the average iterate $\hat{x}_{T} \eqdef \frac{1}{T} \sum_{j=1}^{T} x_j$ satisfies
\[ \ec{f(\hat{x}^T) - f(x^\ast)} \le \frac{\sqn{x^0 - x^\ast}}{2 \gamma n T} + \frac{\gamma^2 L n \sigmaesc^2}{4}. \] 
\end{theorem}

Unfortunately, the theorem above relies on small stepsizes, but we still deem it as a valuable contribution, since it is based on a novel analysis. Indeed, the prior works showed that \algname{RR} approximates a full gradient step, but we show that it is even closer to the implicit gradient step, see the appendix.

To translate the recursion in Theorem~\ref{thm:weakly-convex-f-rr-conv} to a complexity, one can choose a small stepsize and obtain (Corollary~\ref{corr:weakly-convex-f} in the appendix) the following bound for \algname{RR}/\algname{SO}:
\[  \ec{f(\hat{x}^T) - f(x^\ast)} = \mathcal{O}\br{ \frac{L \sqn{x^0 - x^\ast}}{T} + \frac{L^{1/3} \norm{x^0 - x^\ast}^{4/3} \sigmaesc^{2/3}}{n^{1/3} T^{2/3}} }. \]
Stich~\cite{stich2019unified} gives a convergence upper bound of $\mathcal{O}\br{ \frac{L \sqn{x^0 - x^\ast}}{nT} + \frac{\sigmaesc \norm{x^0 - x^\ast}}{\sqrt{n T}}}$ for \algname{SGD}. Comparing upper bounds, we see that \algname{RR}/\algname{SO} beats \algname{SGD} when the number of epochs satisfies $T \geq \frac{L^2 \sqn{x^0 - x^\ast} n}{\sigmaesc^2}$. To the best of our knowledge, there are no strict lower bounds in this setting. Safran and Shamir~\cite{Safran2020good} suggest a lower bound of $\Omega\br{ \frac{\sigmaesc}{\sqrt{n T^3}} + \frac{\sigmaesc}{n T}} $ by setting $\mu$ to be small in their lower bound for $\mu$-strongly convex functions, however this bound may be too optimistic.

\subsection{Non-convex objectives}
\label{sec:non-convex}

For non-convex objectives, we formulate the following assumption on the gradients variance.
\begin{assumption}
  \label{asm:2nd-moment}
  There exist nonnegative constants $A, B \geq 0$ such that for any $x \in \R^d$ we have,
  \begin{equation}
    \label{eq:2nd-moment-bound}
    \frac{1}{n} \sum \limits_{i=1}^{n} \sqn{\nabla f_{i} (x) - \nabla f(x)} \leq 2 A \br{f(x) - f^\ast} + B^2.
  \end{equation}
\end{assumption}
Assumption~\ref{asm:2nd-moment} is quite general: if there exists some $G > 0$ such that $\norm{\nabla f_{i} (x)} \leq G$ for all $x \in \R^d$ and $i \in \{ 1, 2, \ldots, n \}$, then Assumption~\ref{asm:2nd-moment} is clearly satisfied by setting $A = 0$ and $B = G$. Assumption~\ref{asm:2nd-moment} also generalizes the uniformly bounded variance assumption commonly invoked in work on non-convex \algname{SGD}, which is equivalent to~\eqref{eq:2nd-moment-bound} with $A=0$. Assumption~\ref{asm:2nd-moment} is a special case of the Expected Smoothness assumption of \cite{Khaled2020}, and it holds whenever each $f_{i}$ is smooth and lower-bounded, as the next proposition shows.

\begin{proposition}
  \label{prop:2nd-moment-bound}
  \cite[special case of Proposition 3]{Khaled2020} Suppose that $f_{1}, f_2, \ldots, f_n$ are lower bounded by $f_1^\ast, f_2^\ast, \ldots, f_n^\ast$ respectively and that Assumption~\ref{asm:f-smoothness} holds. Then there exist constants $A, B \geq 0$ such that Assumption~\ref{asm:2nd-moment} holds.
\end{proposition}

We now give our main convergence theorem for \algname{RR} without assuming convexity.
\begin{theorem}
  \label{thm:rr-nonconvex}
  Suppose that Assumptions~\ref{asm:f-smoothness} and~\ref{asm:2nd-moment} hold. Then for Algorithm~\ref{alg:rr} run for $T$ epochs with a stepsize $\gamma \leq \min \br{ \frac{1}{2 L n}, \frac{1}{(A L^2 n^2 T)^{1/3}} }$ we have
  \[ 
    \min \limits_{t=0, \ldots, T-1} \ecn{\nabla f(x^k)} \leq \frac{12 \br{f(x^0) - f^\ast}}{\gamma n T} + 2 \gamma^2 L^2 n B^2. 
  \]
  If, in addition, $A=0$ and $f$ satisfies the Polyak-\L{}ojasiewicz inequality with $\mu>0$, i.e., $\|\nabla f(x)\|^2 \ge 2\mu(f(x)-f^\star)$ for any $x\in\R^d$, then
  \[
	   \ec{f(x^T)-f^\star}
	  \le \br{1-\frac{\gamma\mu n}{2}}^T(f(x^0)-f^\star) + \gamma^2\kappa L n B^2.
  \]
\end{theorem}
\textbf{Comparison with \algname{SGD}.} From Theorem~\ref{thm:rr-nonconvex}, one can recover the complexity that we provide in Table~\ref{tab:conv-rates}, see Corollary~\ref{corr:rr-nonconvex} in the appendix. Let's ignore some constants not related to our assumptions and specialize to uniformly bounded variance. Then, the sample complexity of \algname{RR}, 
$ K_{\mathrm{RR}} \geq \frac{ L \sqrt{n}}{\e^2} ( \sqrt{n}+ \frac{\sigma}{\e} ), $
becomes better than that of \algname{SGD},
$ K_{\mathrm{SGD}} \geq \frac{L}{\e^2} (1+ \frac{\sigma^2}{\e^2})$, whenever $\sqrt{n} \e \leq \sigma$.

\section{Experiments}
\begin{figure}[t]
	\noindent\makebox[\textwidth]{
	\includegraphics[width=0.33\linewidth]{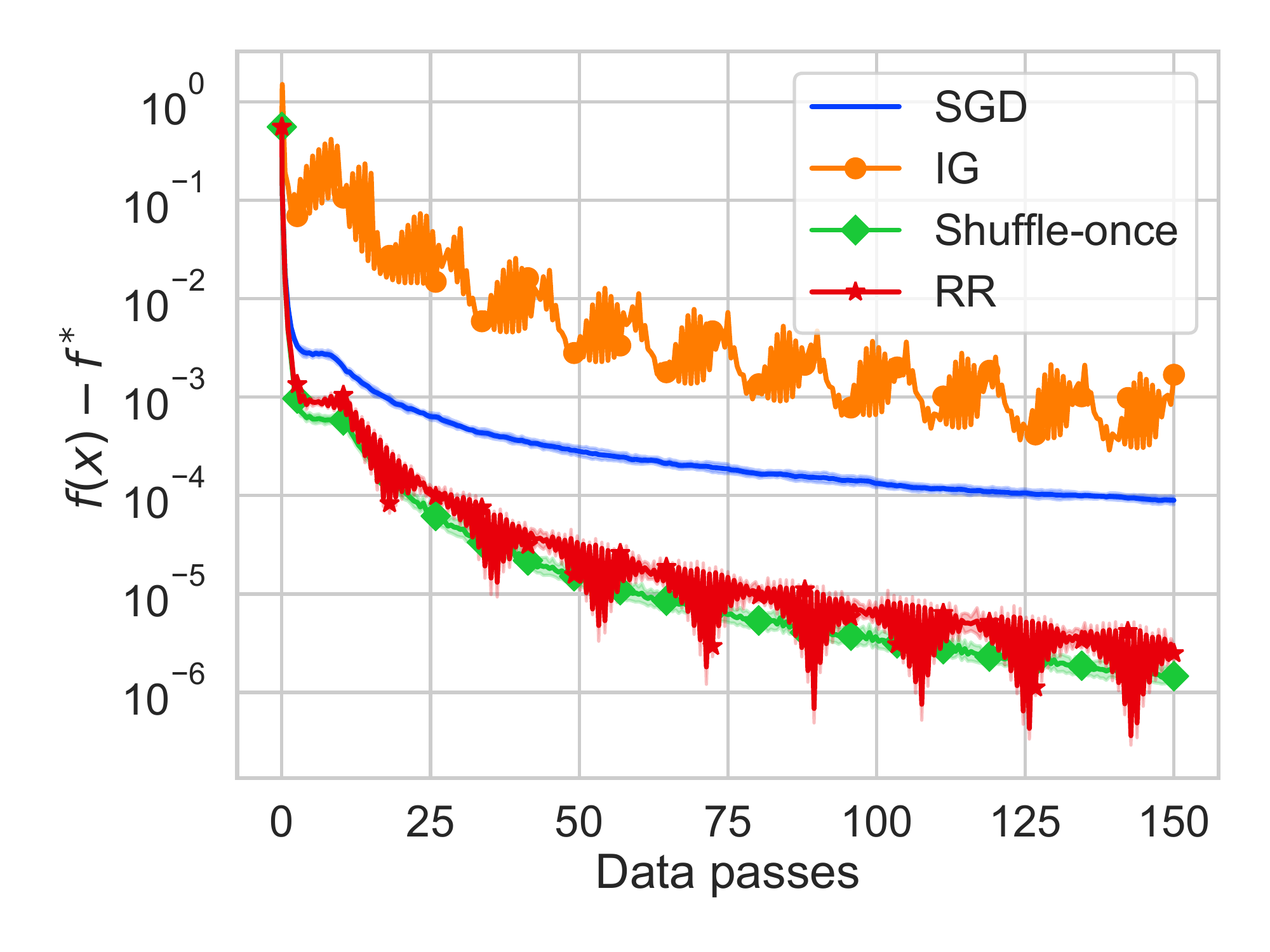}
	\includegraphics[width=0.33\linewidth]{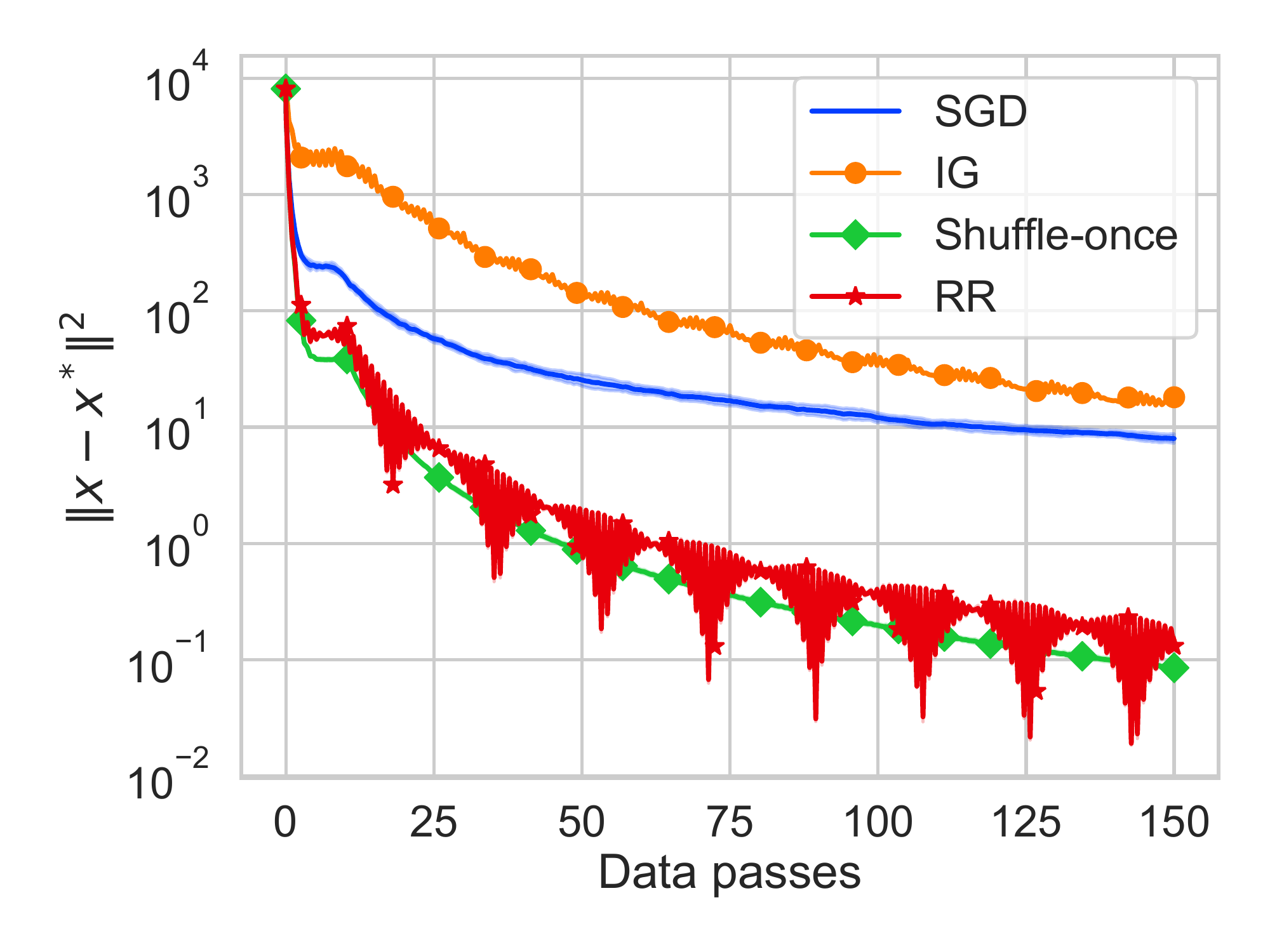}
	\includegraphics[width=0.33\linewidth]{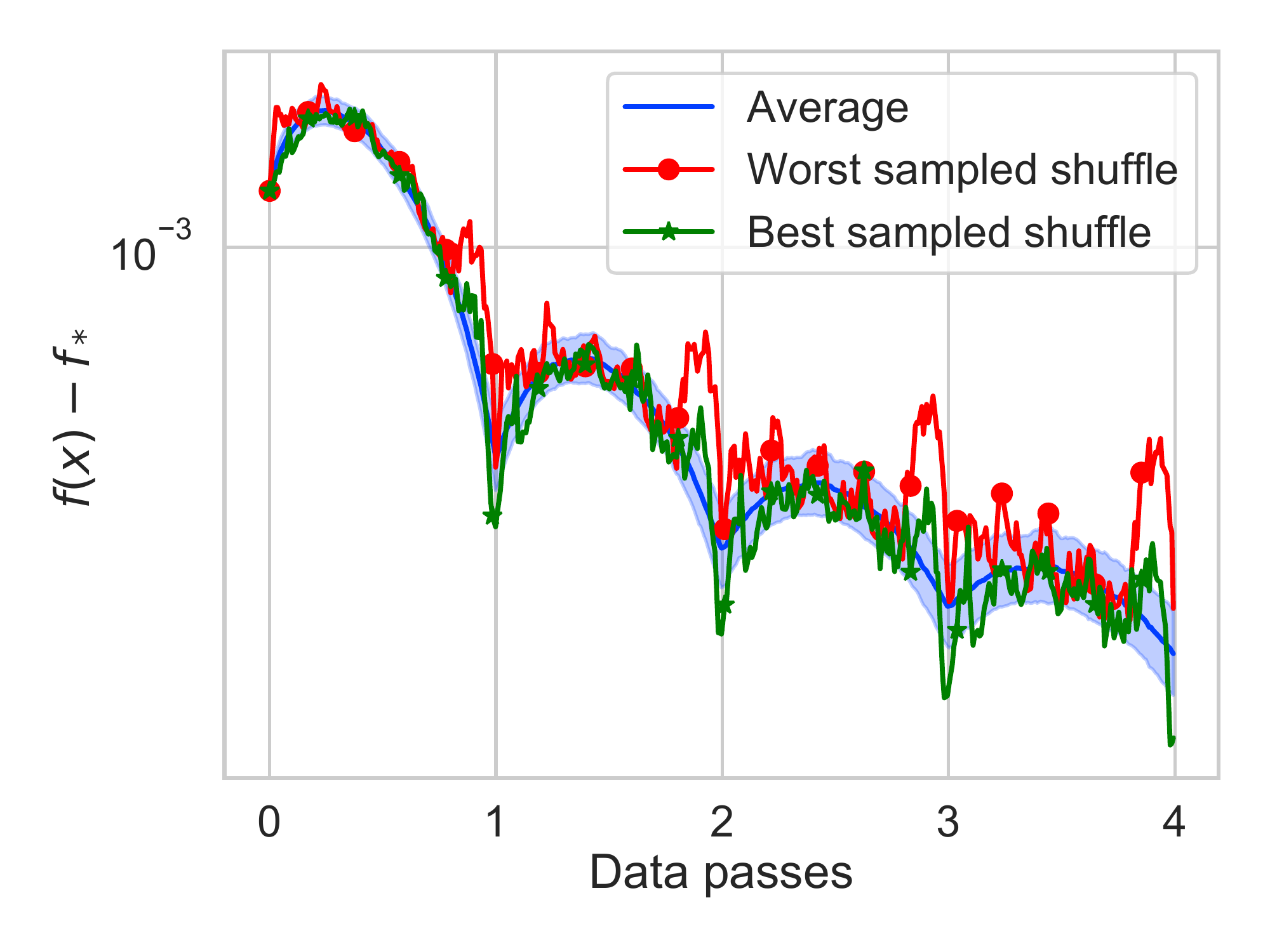}}
	\noindent\makebox[\textwidth]{
	\includegraphics[width=0.33\linewidth]{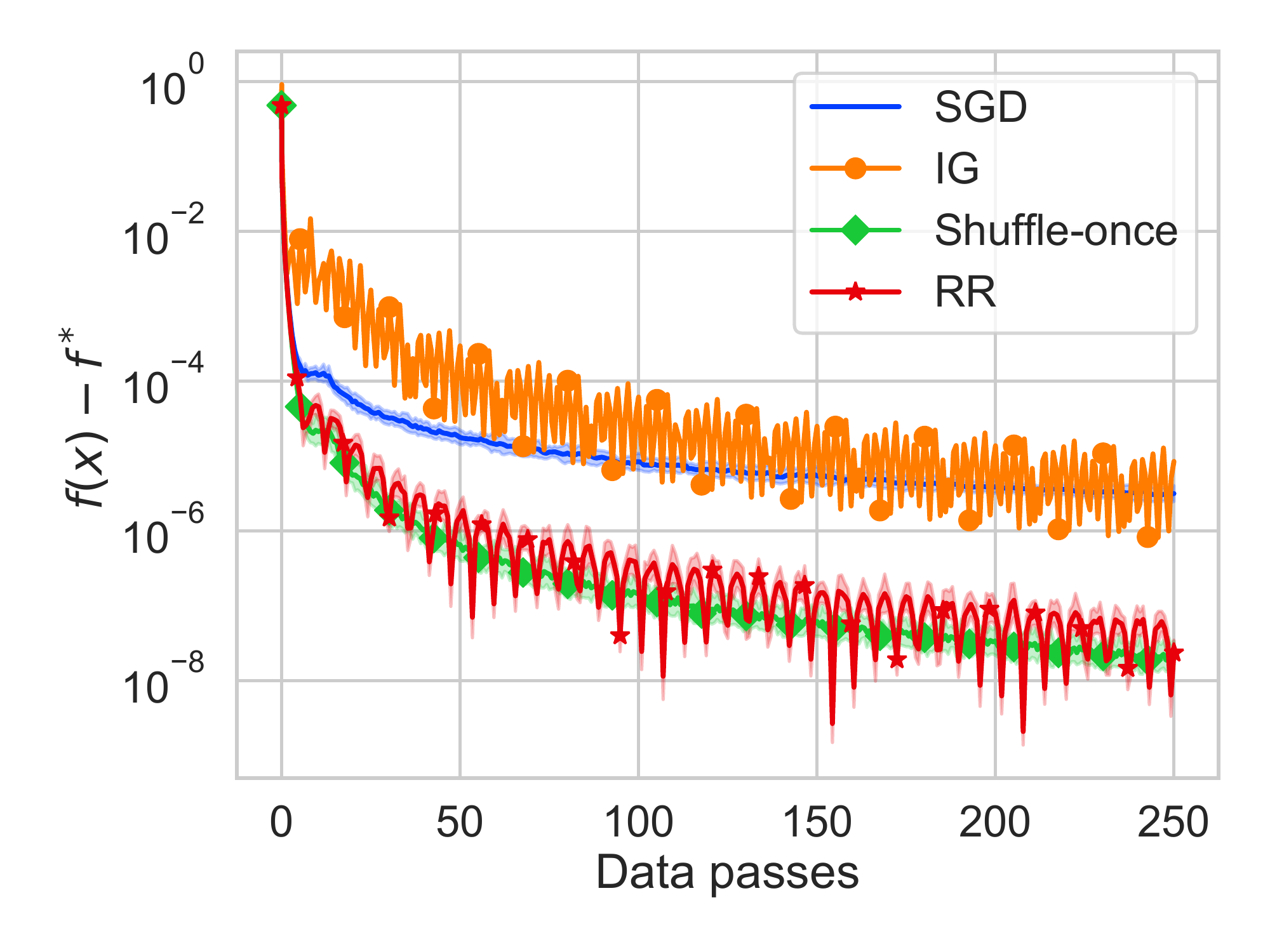}
	\includegraphics[width=0.33\linewidth]{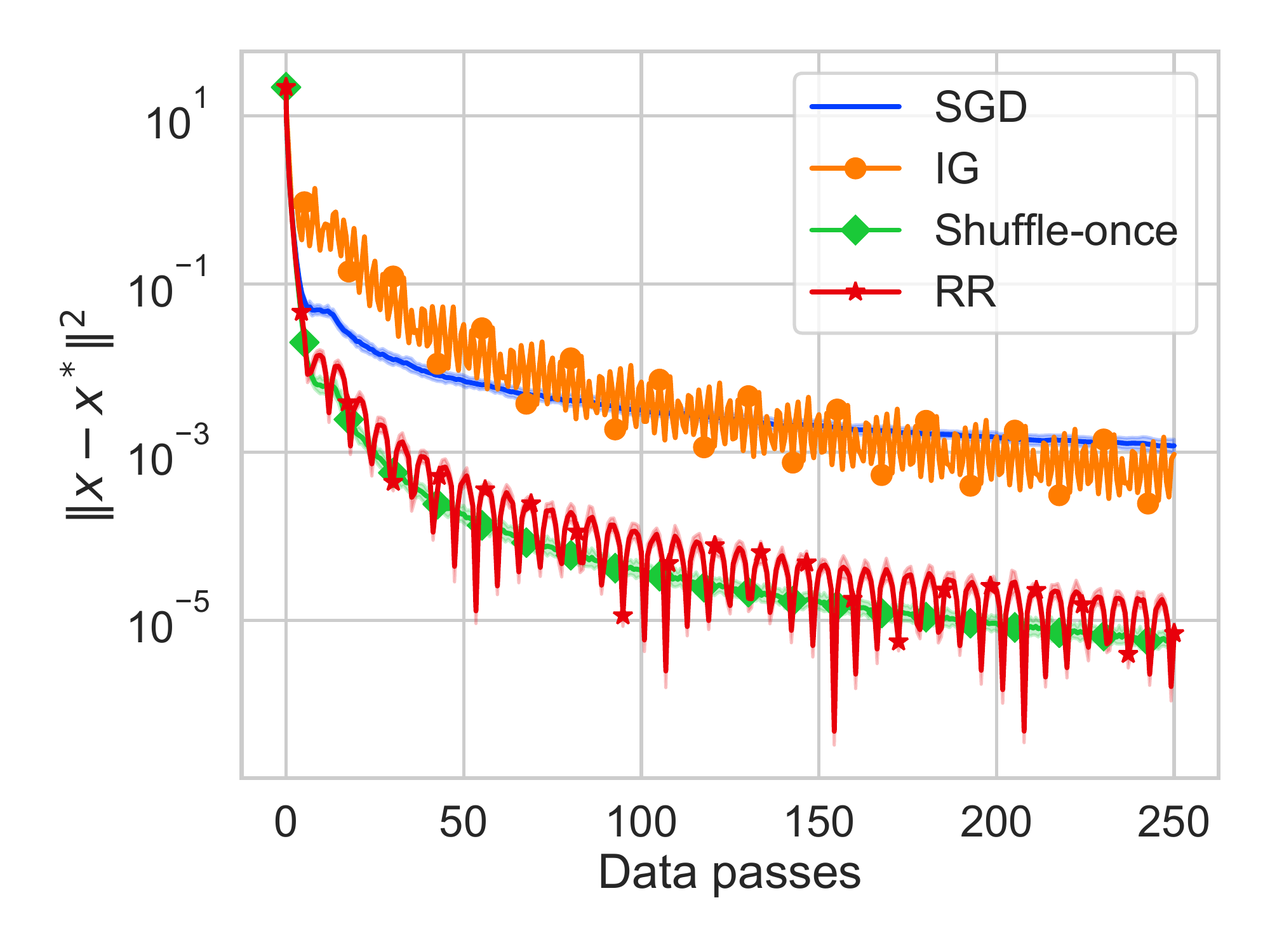}
	\includegraphics[width=0.33\linewidth]{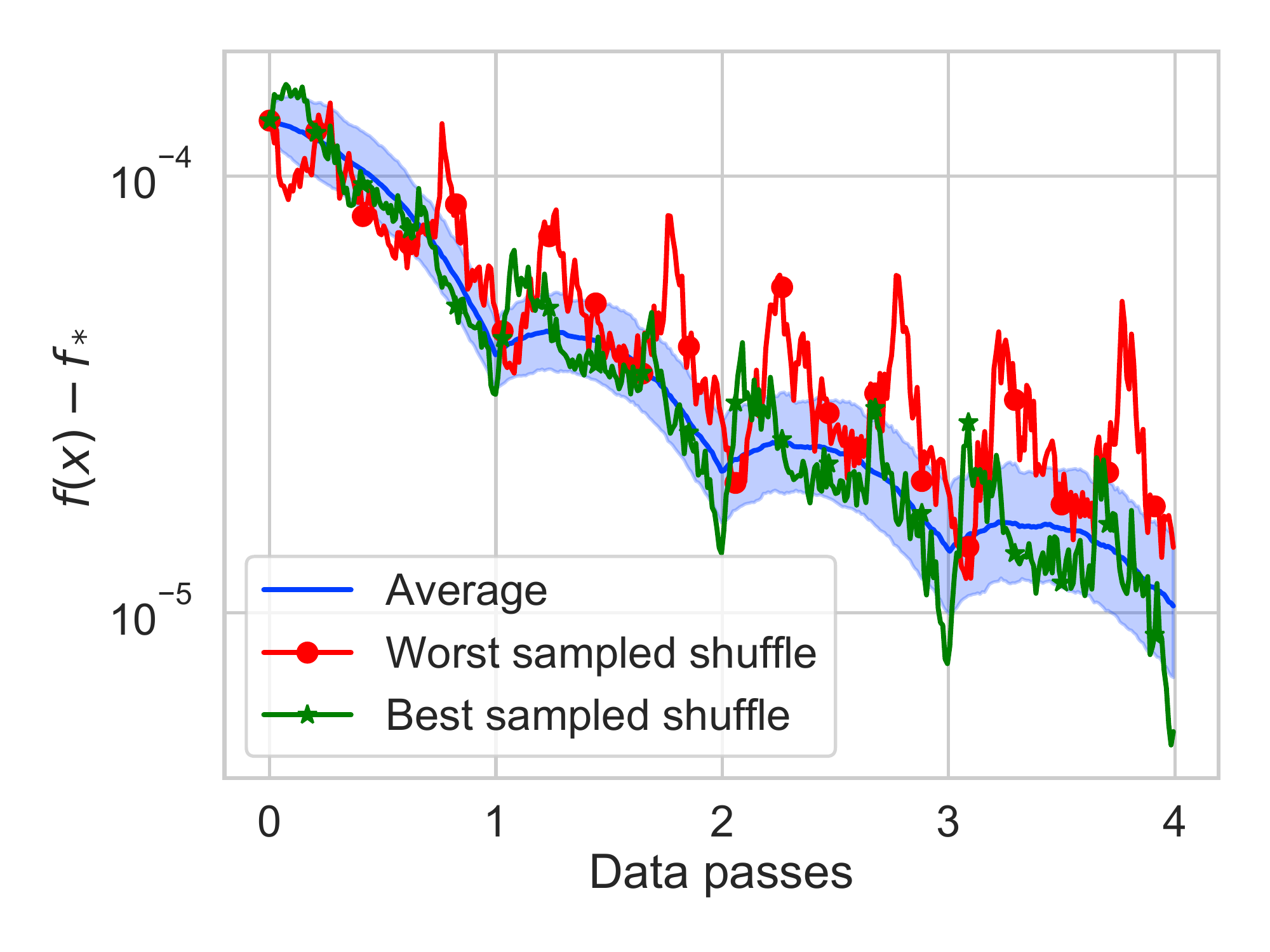}}
	\noindent\makebox[\textwidth]{
	\includegraphics[width=0.33\linewidth]{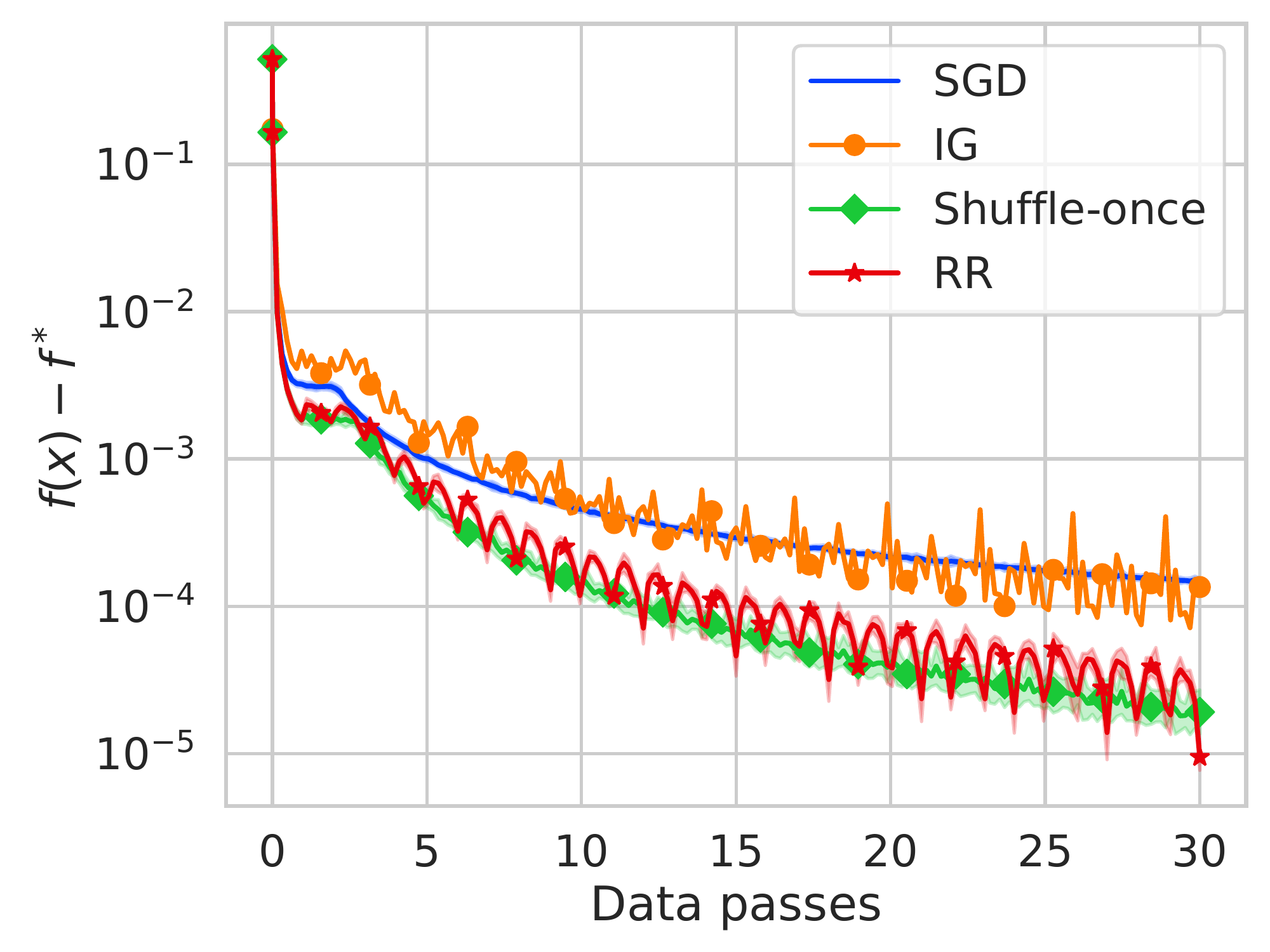}
	\includegraphics[width=0.33\linewidth]{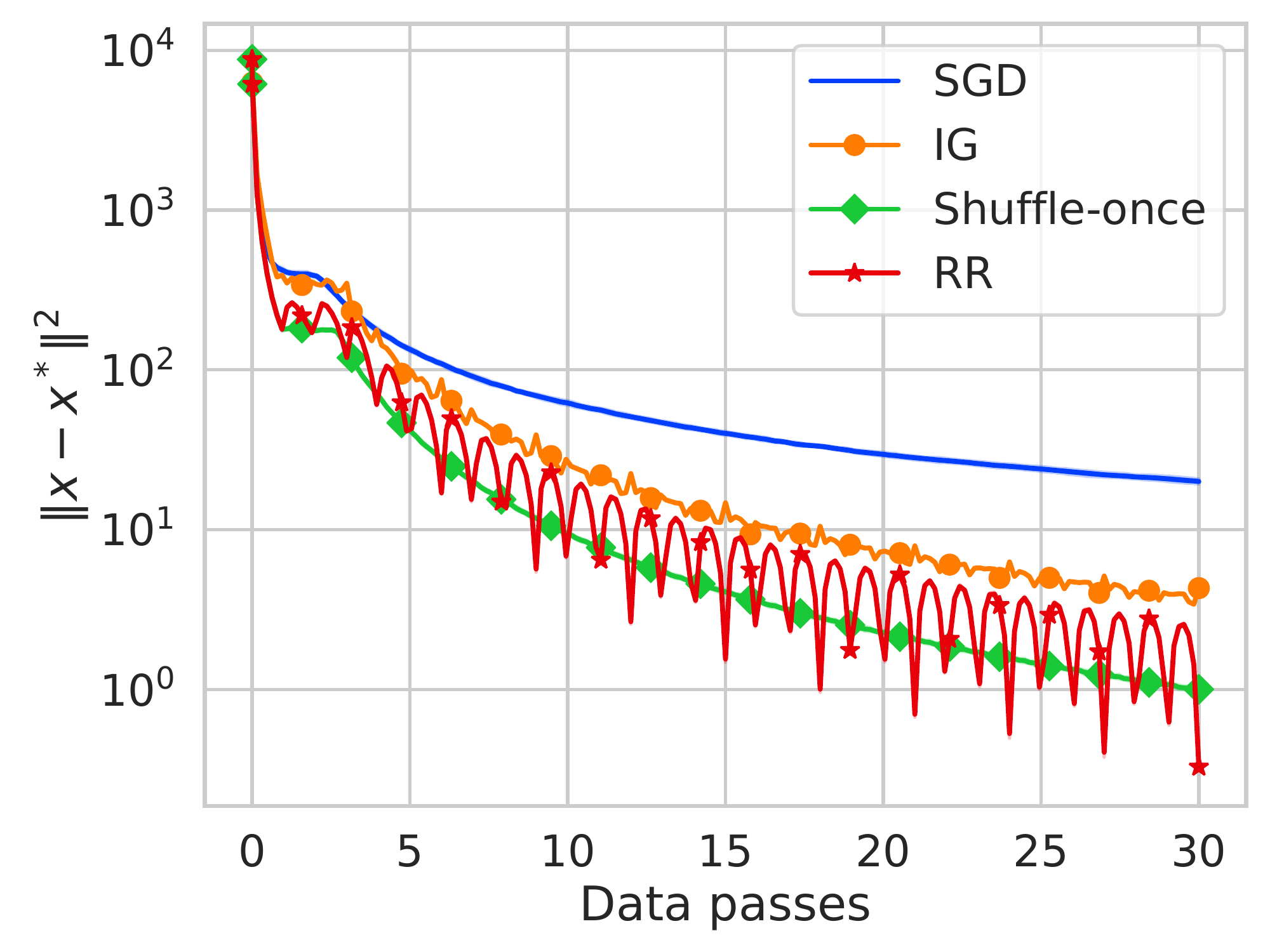}
	\includegraphics[width=0.33\linewidth]{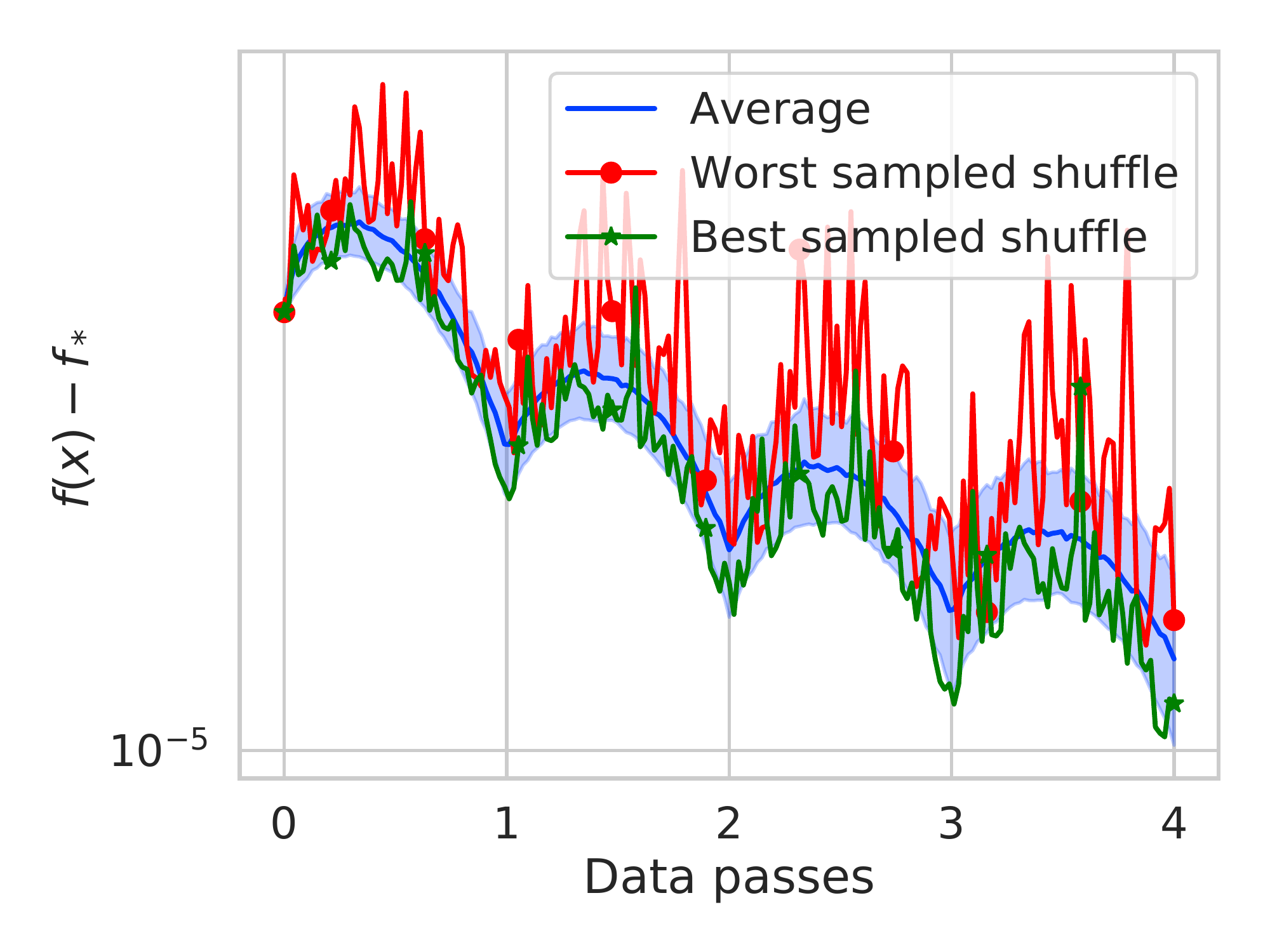}}
\caption{Top: `real-sim' dataset ($N=72,309$; $d=20,958$), middle row: `w8a' dataset ($N=49,749$; $d=300$), bottom: `RCV1' dataset ($N=804,414$; $d=47,236$). Left: convergence of $f(x^k_i)$, middle column: convergence of $\|x^k_i-x^\ast\|^2$, right: convergence of \algname{SO} with different permutations	.}
\label{fig:conv_plots}
\end{figure}

\label{sec:experiments}
We run our experiments on the $\ell_2$-regularized logistic regression problem given by
\[
	\frac{1}{N}\sum \limits_{i=1}^N \br{-\big(b_i \log \big(h(a_i^\top x)\big) + (1-b_i)\log\big(1-h(a_i^\top x)\big)\big)}+\frac{\lambda}{2}\|x\|^2,
\]
where $(a_i, b_i)\in \R^d\times \{0, 1\}$, $i=1,\dotsc, N$ are the data samples and $h\colon t\to1/(1+e^{-t})$ is the sigmoid function. For better parallelism, we use mini-batches of size 512 for all methods and datasets. 
We set $\lambda=L/\sqrt{N}$ and use stepsizes decreasing as $\cO(1/t)$. See the appendix for more details on the parameters used, implementation details, and reproducibility.

\textbf{Reproducibility.} Our code is provided at \href{https://github.com/konstmish/random_reshuffling}{https://github.com/konstmish/random\_reshuffling}. All used datasets are publicly available and all additional implementation details are provided in the appendix.

\textbf{Observations.} One notable property of all shuffling methods is that they converge with oscillations, as can be seen in Figure~\ref{fig:conv_plots}. There is nothing surprising about this as the proof of our Theorem~\ref{thm:all-sc-rr-conv} shows that the intermediate iterates converge to $x^\ast_i$ instead of $x^\ast$. It is, however, surprising how striking the difference between the intermediate iterates within one epoch can be.

\begin{figure}[h]
	\noindent\makebox[\textwidth]{
	\includegraphics[width=0.33\linewidth]{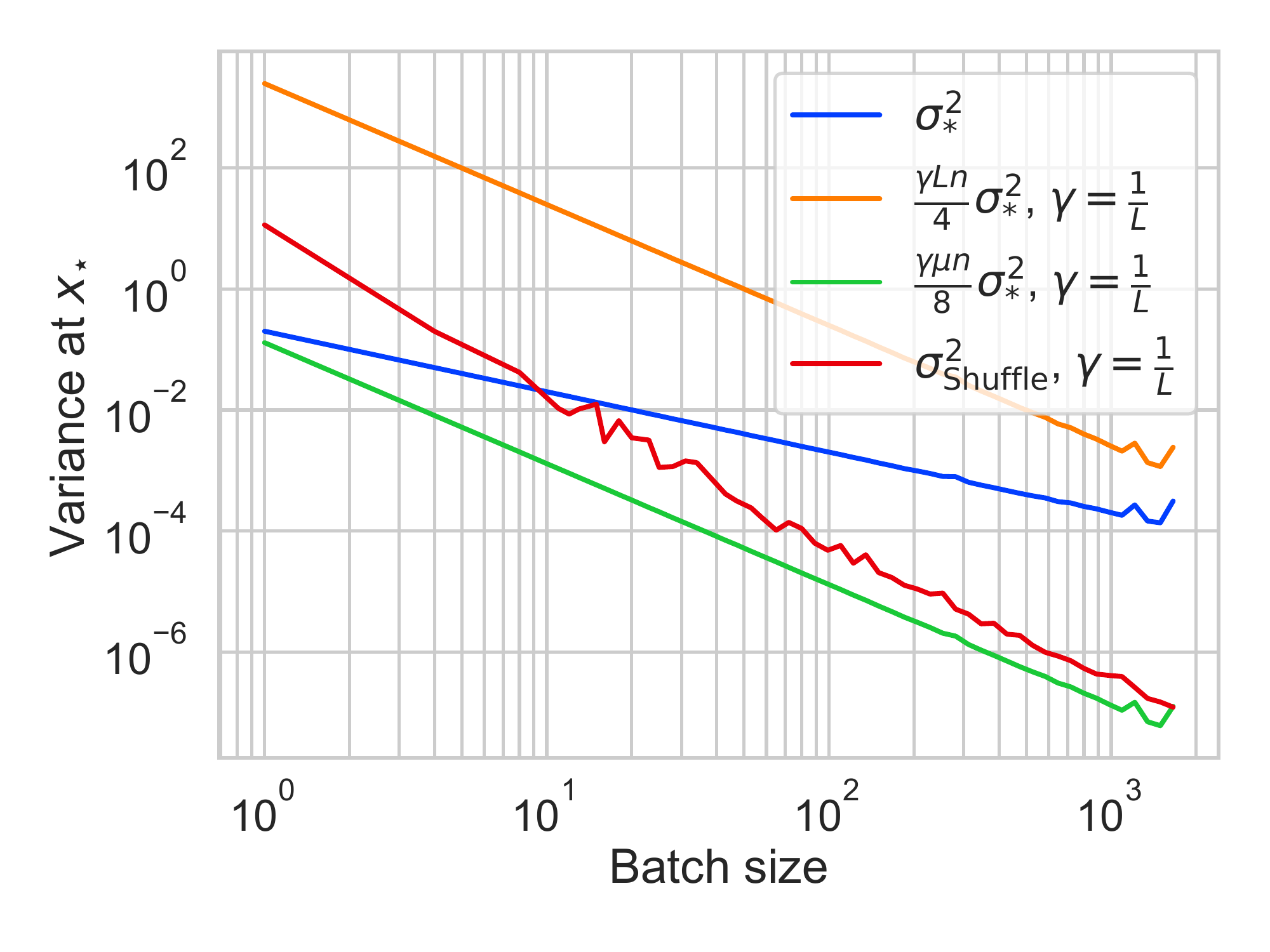}
	\includegraphics[width=0.33\linewidth]{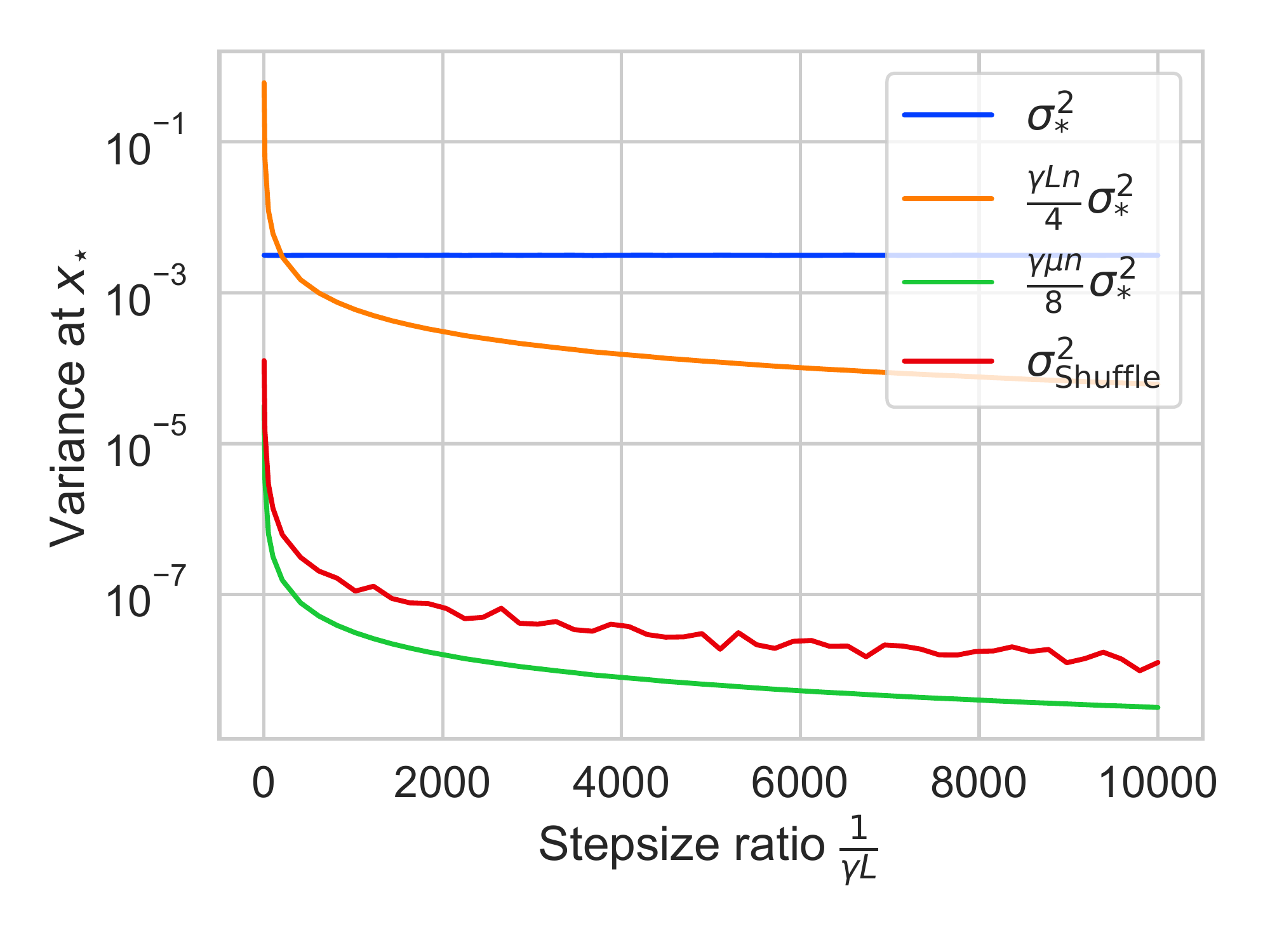}
	\includegraphics[width=0.33\linewidth]{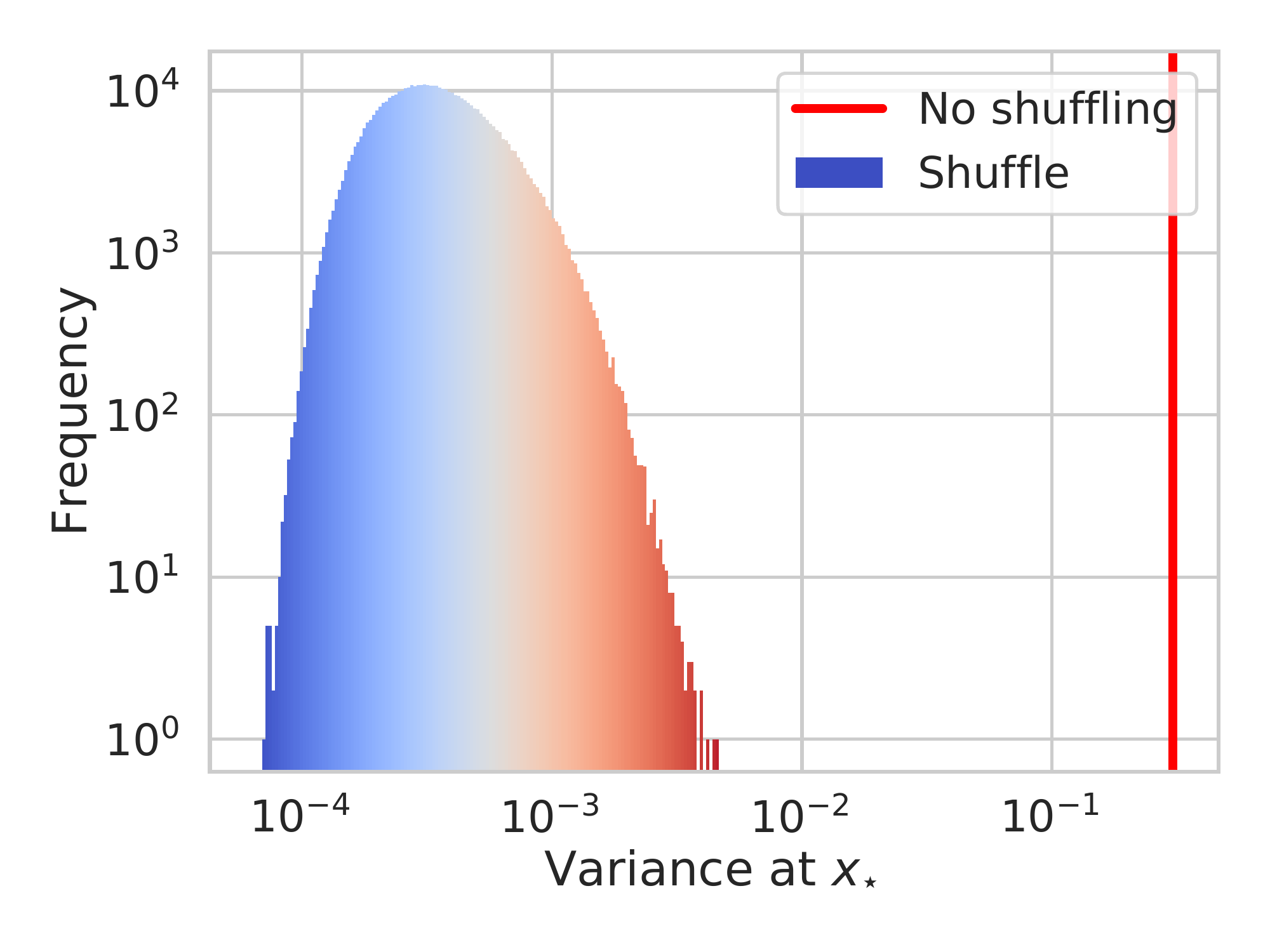}
	}
\caption{Estimated variance at the optimum, $\sigmass$ and $\sigmaesc^2$, for the `w8a' dataset. Left: the values of variance for different mini-batch sizes with $\gamma=\frac{1}{L}$. Middle: variance with fixed mini-batch size 64 for different $\gamma$, starting with $\gamma=\frac{1}{L}$ and ending with $\gamma=\frac{10^{-4}}{L}$. Right: the empirical distribution of $\sigmass$ for $500,000$ sampled permutations with $\gamma=\frac{1}{L}$ and mini-batch size 64.}
	\label{fig:variance}
\end{figure}

Next, one can see that \algname{SO} and \algname{RR} converge almost the same way, which is in line with Theorem~\ref{thm:all-sc-rr-conv}. On the other hand, the contrast with \algname{IG} is dramatic, suggesting existence of bad permutations. The probability of getting such a permutation seems negligible; see the right plot in Figure~\ref{fig:variance}.

Finally, we remark that the first two plots in Figure~\ref{fig:variance} demonstrate the importance of the new variance introduced in Definition~\ref{def:bregman-div-noise}. The upper and lower bounds from Proposition~\ref{prop:shuffling-variance-normal-variance-bound} are depicted in these two plots and one can observe that the lower bound is often closer to the actual value of $\sigmass$ than the upper bound. And the fact that $\sigmass$ very quickly becomes smaller than $\sigmaesc^2$ explains why \algname{RR} often outperforms \algname{SGD} starting from  early iterations.

\chapter{Going Beyond Local SGD in Federated Learning}\label{chapter:proxrr}
\graphicspath{{rr_prox/}}


\section{Introduction}

Modern theory and  practice of training supervised machine learning models is based on the paradigm of regularized empirical risk minimization (ERM)~\cite{shai_book}. While the ultimate goal of supervised learning is  to train models that generalize well to unseen data, in practice only a finite dataset is available during training. Settling for a model merely minimizing the average loss on this training set---the empirical risk---is insufficient, as this often leads to over-fitting and poor generalization performance in practice. Due to this reason, empirical risk is virtually always amended with a suitably chosen regularizer whose role is to encode prior knowledge about the learning task at hand, thus biasing the training algorithm towards better performing models.

The regularization framework is quite general and perhaps surprisingly it also allows us to consider methods for federated learning (FL)---a paradigm in which we aim at training model for a number of clients that do not want to reveal their data~\cite{FEDLEARN, McMahan17, kairouz2019advances}. The training in FL usually happens on devices with only a small number of model updates being shared with a global host. To this end, Federated Averaging algorithm has emerged that performs \algname{Local SGD} updates on the clients' devices and periodically aggregates their average. Its analysis usually requires special techniques and deliberately constructed sequences hindering the research in this direction. We shall see, however, that the convergence of our \algname{FedRR} follows from merely applying our algorithm for regularized problems to a carefully chosen reformulation.

Formally, regularized ERM problems are optimization problems of the form \begin{equation}
  \label{eq:prox-finite-sum-min}  
  \min \limits_{x \in \R^d} \Bigl [ P(x) \eqdef \frac{1}{n} \sum \limits_{i=1}^{n} f_{i} (x)  + \psi(x)\Bigr ],
\end{equation}
where $f_i\colon\R^d \to \R$ is the loss of model parameterized by vector $x\in \R^d$ on the $i$-th training data point, and $\psi:\R^d\to \R \cup \{+\infty\}$ is a regularizer.  Let $[n]\eqdef \{1,2,\dots,n\}$. We shall make the following assumption throughout the paper without explicitly mentioning it:

\begin{assumption}\label{as:smooth_fi_proper_psi}
The functions $f_i$ are $L_i$-smooth and convex, and the regularizer $\psi$ is proper, closed and convex. Let $L_{\max} \eqdef \max_{i\in [n]} L_i$.
\end{assumption}

In some results we will additionally assume that either the individual functions $f_i$, or their average $f\eqdef \frac{1}{n}\sum_i f_i$, or the regularizer $\psi$ are $\mu$-strongly convex.  Whenever we need such additional assumptions, we will make this explicitly clear. While all these concepts are standard, we review them briefly in Section~\ref{sec:basic_notions}.

{\bf \algname{Proximal SGD}.} 
When the number $n$ of training data points  is huge, as is increasingly common in practice, the most efficient algorithms for solving \eqref{eq:prox-finite-sum-min} are stochastic first-order methods, such as Stochastic Gradient Descent (\algname{SGD})~\cite{bordes2009sgd}, in one or another of its many variants proposed in the last decade~\cite{Shang2018,pham2020proxsarah}. These method almost invariably rely on alternating stochastic gradient  steps with the evaluation of the proximal operator
\[
	\prox_{\gamma\psi}(x)
	\eqdef \argmin_{z\in \R^d} \left\{\gamma\psi(z) + \frac{1}{2}\|z-x\|^2\right\}.
\]
The simplest of these has the form
\begin{equation}\label{eq:prox-SGD}
	x^{k+1}_{\mathrm{SGD}} = \prox_{\gamma_k \psi}(x^k_{\mathrm{SGD}} - \gamma_k \nabla f_{i_k}(x^k_{\mathrm{SGD}})),
\end{equation}
where $i_k$ is an index from $\{1,2,\dots,n\}$ chosen uniformly at random, and $\gamma_k>0$ is a properly chosen stepsize. Our understanding of \eqref{eq:prox-SGD} is quite mature; see \cite{gorbunov2020unified} for a general  treatment which considers  methods of this form in conjunction with more advanced stochastic gradient estimators in place of $\nabla f_{i_k}$.

Applications such as training sparse linear models~\cite{tibshirani1996regression}, nonnegative matrix factorization~\cite{lee1999learning}, image deblurring~\cite{rudin1992nonlinear, bredies2010total}, and training with group selection~\cite{yuan2006model} all rely on the use of hand-crafted  regularizes. For most of them, the proximal operator can be evaluated efficiently, and \algname{SGD} is near or at the top of the list of efficient training algorithms.

{\bf \algname{Random Reshuffling}.} A particularly successful variant of \algname{SGD} is based on the idea of random shuffling (permutation) of the training data followed by $n$ iterations of the form \eqref{eq:prox-SGD}, with the index $i_k$ following the pre-selected  permutation~\cite{bottou2012stochastic}. This process is repeated several times, each time using a new freshly sampled random permutation of the data, and the resulting  method is known under the name \algname{ Random Reshuffling (RR)}.\footnote{While we will comment on this in more detail later, \algname{RR} is not known to converge in the proximal setting, i.e., if $\psi \neq 0$. Moreover, it is not even clear if this is the right proximal extension of \algname{RR}.} When the same permutation is used throughout, the technique is known under the name \algname{Shuffle-Once (SO)}.

One of the main advantages of this approach is rooted in its intrinsic ability to avoid cache misses when reading the data from memory, which enables a significantly faster implementation.  Furthermore, \algname{RR} is often observed to converge in fewer iterations than \algname{SGD} in practice. This can intuitively be ascribed to the fact that while due to its sampling-with-replacement approach \algname{SGD} can miss to learn from some data points in any given epoch,  \algname{RR} will necessarily learn from each data point in each epoch. 

\section{Related work}
Understanding the random reshuffling trick, and why it works, has been a non-trivial open problem for a  long time~\cite{Bottou2009,RR-conjecture2012,Gurbuzbalaban2019RR,haochen2018random}. Until recent development which lead to a significant simplification of the convergence analysis technique and proofs (see \Cref{chapter:rr}), prior state of the art relied on long and  elaborate proofs requiring sophisticated arguments and tools, such as  analysis via the Wasserstein distance~\cite{Nagaraj2019}, and relied on a significant number of strong assumptions about the objective \cite{shamir2016without, haochen2018random}. In alternative recent development, Ahn et al.~\cite{Ahn2020} also develop new tools for analyzing the convergence of \algname{Random Reshuffling}, in particular using decreasing stepsizes and for objectives satisfying the Polyak-{\L}ojasiewicz condition, a generalization of strong convexity~\cite{polyak1963gradient, lojasiewicz1963topological}.

The difficulty of analyzing \algname{RR} has been the main obstacle in the development of even some of the most seemingly benign extensions of the method. Indeed, while all these are well understood in combination with its much simpler-to-analyze cousin \algname{SGD},  to the best of our knowledge, there exists no theoretical analysis of proximal, parallel, and importance sampling variants of \algname{RR} with both constant and decreasing stepsizes, and in most cases it is not even clear how should such methods be constructed.
 Empowered by and building on the advances Chapter~\ref{chapter:rr}, in this chapter we address all these challenges.

\section{Contributions}

In this section we outline the key contributions  of our work, and also offer a few intuitive explanations motivating some of the development.

{\bf From \algname{RR} to \algname{ProxRR}.}
Despite rich literature on \algname{Proximal SGD} \cite{gorbunov2020unified}, it is not obvious how one should extend \algname{RR} to solve problem~\eqref{eq:prox-finite-sum-min} when a  nonzero regularizer $\psi$ is present. Indeed, the standard practice for \algname{SGD} is to apply the proximal operator after each stochastic step~\cite{duchi2009efficient}, i.e., in analogy with \eqref{eq:prox-SGD}. On the other hand, \algname{RR} is motivated by the fact that a data pass approximates the full gradient step. If we apply the proximal operator after each iteration of \algname{RR}, we would no longer approximate the full gradient after an epoch, as illustrated by the next example.
\begin{example}
	Let $n=2$, $\psi(x)=\frac{1}{2}\|x\|^2$, $f_1(x)=\<c_1, x>$, $f_2(x)=\<c_2, x>$ with some $c_1, c_2\in\mathbb{R}^d$, $c_1\neq c_2$. Let $x^0\in\mathbb{R}^d$, $\gamma>0$ and define $x^1 = x^0 - \gamma \nabla f_1(x^0)$, $x^2=x^1-\gamma\nabla f_2(x^1)$. Then, we have $\prox_{2\gamma\psi}(x^2)=\prox_{2\gamma\psi}(x^0-2\gamma\nabla f(x^0))$. However, if $\tilde x^1=\prox_{\gamma\psi}(x^0 - \gamma \nabla f_1(x^0))$ and $\tilde x^2=\prox_{\gamma \psi}(x_1-\gamma \nabla f_2(\tilde x^1))$, then $\tilde x^2\neq \prox_{2\gamma\psi}(x^0-2\gamma\nabla f(x^0))$.
\end{example}

Motivated by this observation, we propose \algname{ProxRR} (\Cref{alg:proxrr}), in which the proximal operator is applied at the end of each epoch of \algname{RR}, i.e., after each pass through all randomly reshuffled data. 

\begin{algorithm}[t]
    \caption{\algname{Proximal Random Reshuffling} \rrbox{(\algname{ProxRR})} and \algname{Shuffle-Once} \sobox{(\algname{ProxSO})}}
    \label{alg:proxrr}
    \begin{algorithmic}[1]
    \Require Stepsizes $\gamma_k > 0$, initial vector $x^0 \in \R^d$, number of epochs $T$
    \State Sample a permutation $\pi =(\prm{0}, \prm{1}, \ldots, \prm{n-1})$ of $[n]$   \sobox{(Do step 1 only for \algname{ProxSO})} 
        \For{epochs $k=0,1,\dotsc,T-1$}
            \State  Sample a permut.\ $\pi =(\prm{0}, \prm{1}, \ldots, \prm{n-1})$ of $[n]$  \rrbox{(Do step 3 only for \algname{ProxRR})}
            \State $x^k_0 = x^k$
            \For{$i=0, 1, \ldots, n-1$}
                \State $x^{k}_{i+1} = x^k_{i} - \gamma_k \nabla f_{\prm{i}} (x^k_i)$
            \EndFor
        \State $x^{k+1} = \prox_{\gamma_k n \psi}(x^{k}_{n})$
        \EndFor
    \end{algorithmic}
\end{algorithm}

A notable property of \Cref{alg:proxrr} is that only a single proximal operator evaluation is needed during each data pass. This is in sharp contrast with the way \algname{Proximal SGD} works, and offers significant advantages in regimes where the evaluation of the proximal mapping is expensive (e.g., comparable to the evaluation of $n$ gradients $\nabla f_1, \dots, \nabla f_n$).

We establish several convergence results for \algname{ProxRR}, of which we highlight two here. Both offer a linear convergence rate with a fixed stepsize to a neighborhood of the solution. Firstly, in the case when each $f_i$ is $\mu$-strongly convex, we prove the rate (see Theorem~\ref{thm:f-strongly-convex-psi-convex})
\[
			\ecn{x^T - x^\ast} \leq \br{1 - \gamma \mu}^{n T} \sqn{x^0 - x^\ast} + \frac{2 \gamma^2 \sigmarr^2}{\mu},
\]
where $\gamma_k = \gamma \leq \frac{1}{L_{\max}}$ is the stepsize,   and $\sigmarr^2$ is a shuffling radius constant (for precise definition, see \eqref{eq:bregman-div-rad}). In Theorem~\ref{thm:shuffling-radius-bound} we bound the shuffling radius in terms of  $\norm{\nabla f(x^\ast)}^2$, $n$,  $L_{\max}$ and the more common quantity $\sigmaesc^2 \eqdef \frac{1}{n} \sum_{i=1}^{n} \sqn{\nabla f_{i} (x^\ast) - \nabla f(x^\ast)}$.

 Secondly, if  $\psi$ is $\mu$-strongly convex, we prove the rate \[
			\ecn{x^T - x^\ast} \leq \br{1 + 2\gamma \mu n}^{-T} \sqn{x^0 - x^\ast} + \frac{ \gamma^2 \sigmarr^2}{\mu},
\]
	where $\gamma_k = \gamma \leq \frac{1}{L_{\max}}$ is the stepsize (see Theorem~\ref{thm:psi-strongly-convex-f-convex}).

Both mentioned rates show exponential (linear in logarithmic scale) convergence to a neighborhood whose size is proportional to $\gamma^2 \sigmarr^2$. Since we can choose $\gamma$ to be arbitrarily small or periodically decrease it, this implies that the iterates converge to $x^\ast$ in the limit. Moreover, we show in \Cref{sec:strongly_convex} that when $\gamma = \cO(\frac{1}{T})$ the error is $\cO(\frac{1}{T^2})$, which is superior to the $\cO(\frac{1}{T})$ error of \algname{SGD}.

{\bf Decreasing stepsizes.}
The convergence of \algname{RR} is not always exact and depends on the parameters of the objective. Similarly, if the shuffling radius $ \sigmarr^2$ is positive, and we wish to find an $\e$-approximate solution, the optimal choice of a fixed stepsize for \algname{ProxRR} will depend on $\e$. This deficiency can be fixed by using decreasing stepsizes in both vanilla \algname{RR}~\cite{Ahn2020} and in \algname{SGD}~\cite{stich2019unified}. We adopt the same technique to our setting. However, we depart from \cite{Ahn2020} by only adjusting the stepsize once per epoch rather than at every iteration, similarly to the concurrent work of \cite{tran2020shuffling} on \algname{RR} with momentum.  For details, see Section~\ref{sec:extensions}.

{\bf Importance sampling for \algname{ProxRR}.}
While importance sampling is a well established technique for speeding up the convergence of \algname{SGD} \cite{IProx-SDCA, Khaled2020},  no importance sampling variant of \algname{RR} has been proposed nor analyzed. This is not surprising since the key property of importance sampling in \algname{SGD}---unbiasedness---does not hold for \algname{RR}. Our approach to equip \algname{ProxRR} with importance sampling is via a reformulation of problem \eqref{eq:prox-finite-sum-min} into a similar problem with a larger number of summands. In particular, for each $i\in [n]$ we include $n_i$ copies of the function $\frac{1}{n_i}f_i$, and then take average of all $N = \sum_i n_i$ functions constructed this way. The value of $n_i$ depends on the ``importance'' of $f_i$, described below.  We then apply \algname{ProxRR} to this reformulation.

If $f_i$ is $L_i$-smooth for all $i\in [n]$ and we let $\bar{L}\eqdef \frac{1}{n}\sum_i L_i$, then we choose $n_i= \lceil \frac{L_i}{\bar{L}} \rceil$. It is easy to show that $N\leq 2n$, and hence our reformulation leads to at most a doubling of the number of functions forming the finite sum. However, the overall complexity of \algname{ProxRR} applied to this reformulation will depend on $\bar{L}$ instead of $\max_i L_i$ (see Theorem~\ref{thm:IS}), which can lead to a significant improvement.  For details of the construction and our complexity results, see Section~\ref{sec:extensions}.

{\bf Application to federated learning.}
In Section~\ref{sec:FL} we describe an application of our results to federated learning \cite{FEDLEARN, McMahan17, kairouz2019advances}. 

{\bf Results for \algname{SO}.}
All of our results apply to the Shuffle-Once algorithm as well. For simplicity, we center the discussion around \algname{RR}, whose current theoretical guarantees in the non-convex case are better than that of \algname{SO}. Nevertheless, the other results are the same for both methods, and \algname{ProxRR} is identical to \algname{ProxSO} in terms of our theory too. A study of the empirical differences between \algname{RR} and \algname{SO} can be found in Chapter~\ref{chapter:rr}.

\section{Preliminaries}

In our analysis, we build upon the notions of limit points and shuffling variance introduced in \Cref{chapter:rr} for vanilla (i.e., non-proximal) \algname{RR}. Given a stepsize $\gamma > 0$ (held constant during each epoch) and a permutation $\pi$ of $\{1,2,\dots,n\}$, the inner loop iterates of \algname{RR}/\algname{SO} converge to a neighborhood of intermediate limit points $x^\ast_1, x^\ast_2, \ldots, x^\ast_{n}$ defined by
\begin{equation}
	x^\ast_i \eqdef x^\ast - \gamma \sum \limits_{j=0}^{i-1} \nabla f_{\pi_{j}} (x^\ast), \quad i=1,\dotsc, n. \label{eq:x_ast_i_prox}
\end{equation}
The intuition behind this definition is fairly simple: if we performed $i$ steps starting at $x^\star$, we would end up close to $x^\star_i$. To quantify the closeness, we define the shuffling radius.

\begin{definition}[Shuffling radius]
	\label{def:bregman-div-rad}
	Given a stepsize $\gamma>0$ and a random permutation $\pi$ of $\{ 1, 2, \ldots, n \}$ used in Algorithm~\ref{alg:proxrr}, define $x^\ast_i = x^\ast_i (\gamma, \pi)$ as in \eqref{eq:x_ast_i_prox}. Then, the shuffling radius is defined by
	\begin{equation}\label{eq:bregman-div-rad}  	
	\sigmarr^2 (\gamma) \eqdef \max \limits_{i=0, \ldots, n-1} \left [ \frac{1}{\gamma^2} \mathbb{E}_\pi\bigl[D_{f_{\pi_{i}}} (x^\ast_i, x^\ast)\bigr] \right ], 
	\end{equation}
	where the expectation is taken with respect to the randomness in the permutation $\pi$. If there are multiple stepsizes $\gamma_1, \gamma_2, \ldots$ used in Algorithm~\ref{alg:proxrr}, we take the maximum of all of them as the shuffling radius, i.e., 
	\[ \sigmarr^2 \eqdef \max_{k \geq 1}  \; \sigmarr^2 (\gamma_k). \]
\end{definition}

The shuffling radius is related by a multiplicative factor in the stepsize to the shuffling variance introduced in the previous Chapter. When the stepsize is held fixed, the difference between the two notions is minimal but when the stepsize is decreasing, the shuffling radius is easier to work with, since it can be upper bounded by problem constants independent of the stepsizes. To prove this upper bound, we rely on a lemma from \Cref{chapter:rr} that bounds the variance when sampling without replacement.

\begin{lemma}[Lemma 1 in Chapter~\ref{chapter:rr}]
	\label{lemma:sampling-without-replacement}
	Let $X_1, \ldots, X_n \in \R^d$ be fixed vectors, let $\bar{X} = \frac{1}{n} \sum_{i=1}^{n} X_i$ be their mean, and let $\sigma^2 = \frac{1}{n} \sum_{i=1}^{n} \sqn{X_i - \bar{X}}$ be their  variance. Fix any $i \in [n]$ and let $X_{\pi_0}, \ldots, X_{\pi_{i-1}}$ be sampled uniformly without replacement from $\{ X_1, \ldots, X_n \}$ and $\bar{X}_{\pi}=\frac{1}{i}\sum_{j=0}^{i-1} X_{\pi_j}$ be their average. Then, the sample average and variance are given by
\begin{equation}\label{eq:b97fg07gdf_08yf8d}
		\ec{\bar{X}_{\pi}} = \bar{X}, \qquad \ecn{\bar{X}_\pi - \bar{X}} = \frac{n-i}{i (n-1)} \sigma^2.
\end{equation}
\end{lemma}

Armed with \Cref{lemma:sampling-without-replacement}, we can upper bound the shuffling radius using the smoothness constant $L_{\max}$, size of the vector $\nabla f(x^\ast)$ and the variance $\sigmaesc^2$ of the gradient vectors $\nabla f_1(x^\ast)$, $\nabla f_2(x^\ast)$, \dots, $\nabla f_n(x^\ast)$. 

\begin{theorem}
	\label{thm:shuffling-radius-bound}
	For any stepsize $\gamma > 0$ and any random permutation $\pi$ of $\{1,2,\dots,n\}$ we have
	\[	
		\sigmarr^2 \le \frac{L_{\max}}{2}n\Bigl(n\|\nabla f(x^\ast)\|^2 + \frac{1}{2}\sigmaesc^2\Bigr),
	\]
	where $x^\ast$ is a solution of Problem~\eqref{eq:prox-finite-sum-min} and $\sigmaesc^2$ is the population variance at the optimum 
	\begin{equation}\label{eq:UG(*G(DG(*DGg87gf7ff}
		\sigmaesc^2 \eqdef \frac{1}{n} \sum \limits_{i=1}^{n} \sqn{\nabla f_{i} (x^\ast) - \nabla f(x^\ast)}.
		\end{equation}
\end{theorem}

All proofs are relegated to the supplementary material. In order to better understand the bound given by \Cref{thm:shuffling-radius-bound}, note that if there is no proximal operator (i.e., $\psi = 0$) then $\nabla f(x^\ast) = 0$ and we get that $\sigmarr^2 \leq \frac{L_{\max} n \sigmaesc^2}{4}$. This recovers the existing upper bound on the shuffling variance from \Cref{chapter:rr} for vanilla \algname{RR}. On the other hand, if $\nabla f(x^\ast) \neq 0$ then we get an additive term of size proportional to the squared norm of $\nabla f(x^\ast)$. 

\section{Convergence Theory for Strongly Convex Losses $f_1, \dots,f_n$}\label{sec:strongly_convex}

Our first theorem establishes a convergence rate for Algorithm~\ref{alg:proxrr} applied with a constant stepsize to Problem~\eqref{eq:prox-finite-sum-min} when each objective $f_i$ is strongly convex. This assumption is commonly satisfied in machine learning applications where each $f_i$ represents a regularized loss on some data points, as in $\ell_2$ regularized linear regression and $\ell_2$ regularized logistic regression.

\begin{theorem}
	\label{thm:f-strongly-convex-psi-convex}
Let \Cref{as:smooth_fi_proper_psi} be satisfied. 		Further, assume that each $f_i$  is $\mu$-strongly convex. If Algorithm~\ref{alg:proxrr} is run with constant stepsize $\gamma_k = \gamma \leq \frac{1}{L_{\max}}$, then the iterates generated by the algorithm satisfy
\[
			\ecn{x^T - x^\ast} \leq \br{1 - \gamma \mu}^{n T} \sqn{x^0 - x^\ast} + \frac{2 \gamma^2 \sigmarr^2}{\mu}.
\]
\end{theorem}

We can convert the guarantee of Theorem~\ref{thm:f-strongly-convex-psi-convex} to a convergence rate by properly tuning the stepsize and using the upper bound of \Cref{thm:shuffling-radius-bound} on the shuffling radius. In particular, if we choose the stepsize as
$
\gamma = \min \pbr{ \frac{1}{L_{\max}}, \frac{\sqrt{\e \mu}}{\sqrt{2} \sigmarr} } ,$
then we obtain $\ecn{x^T - x^\ast} = \mathcal{O}\br{\e}$ provided that the total number of iterations $K_{\mathrm{RR}} = n T$ is at least
\begin{equation}
	\label{eq:proxrr-complexity-f-sc}
        K_{\mathrm{RR}} \geq \biggl (\kappa +  \frac{\sqrt{\kappa n}}{\sqrt{\e} \mu} ( \sqrt{n} \norm{\nabla f(x^\ast)} + \sigmaesc )   \biggr) \log\br{\frac{2 \|x^0 - x^\ast\|^2}{\e}},
\end{equation}
where $\kappa\eqdef L_{\max}/\mu$.

{\bf Comparison with vanilla \algname{RR}.} If there is no proximal operator, then $\norm{\nabla f(x^\ast)} = 0$ and we recover the earlier result from \Cref{chapter:rr} on the convergence of \algname{RR} without proximal operator, which is optimal in $\e$ up to logarithmic factors. On the other hand, when the proximal operator is nonzero, we get an extra term in the complexity proportional to $\norm{\nabla f(x^\ast)}$: thus, even when all the functions are the same (i.e., $\sigmaesc = 0$), we do not recover the linear convergence of \algname{Proximal Gradient Descent} \cite{Karimi2016,Beck2017}. This can be easily explained by the fact that \Cref{alg:proxrr} performs $n$ gradient steps per one proximal step. Hence, even if $f_1=\dotsb=f_n$, \Cref{alg:proxrr} does not reduce to \algname{Proximal Gradient Descent}. We note that other algorithms for composite optimization which may not take a proximal step at every iteration (for example, using stochastic projection steps) also suffer from the same dependence~\cite{patrascu2020stochastic}.

{\bf Comparison with \algname{Proximal SGD}.} 
In order to compare \eqref{eq:proxrr-complexity-f-sc} against the complexity of \algname{Proximal SGD} (Algorithm~\ref{alg:proxsgd}), we recall the following simple result on the convergence of \algname{Proximal SGD}. The result is standard \cite{needell2014stochastic,Gower2019}, with the exception that we present it in a slightly generalized in that we also consider the case when $\psi$ is strongly convex. Our proof is a minor modification of that in \cite{Gower2019}, and we offer it in the appendix for completeness.
\begin{algorithm}[t]
    \caption{\rrbox{Proximal SGD}}
    \label{alg:proxsgd}
    \begin{algorithmic}[1]
    \Require Stepsizes $\gamma_k > 0$, initial vector $x^0 \in \R^d$, number of steps $K$
        \For{steps $k=0,1,\dotsc,K-1$}
            \State Sample $i_k$ uniformly at random from $[n]$
            \State $x^{k+1} = \prox_{\gamma_k \psi}(x^{k} - \gamma_k \nabla f_{i_k} (x^k))$
        \EndFor
    \end{algorithmic}
\end{algorithm}

\begin{theorem}[\algname{Proximal SGD}]
	\label{thm:conv-prox-sgd}
		Let \Cref{as:smooth_fi_proper_psi} hold. 
	Further, suppose that either $f \eqdef \frac{1}{n} \sum_{i=1}^{n} f_i$ is $\mu$-strongly convex or that $\psi$ is $\mu$-strongly convex. If Algorithm~\ref{alg:proxsgd} is run with a constant stepsize $\gamma_k = \gamma > 0$ satisfying $\gamma \leq \frac{1}{2 L_{\max}}$, then the final iterate returned by the algorithm after $K$ steps satisfies
	\[ 
	\ecn{x^K - x^\ast} \leq \br{1 - \gamma \mu}^{K} \sqn{x^0 - x^\ast} + \frac{2 \gamma \sigmaesc^2}{\mu}.  \] 
\end{theorem}

Furthermore, by choosing the stepsize $\gamma$ as
$\gamma = \min \pbr{ \frac{1}{2 L_{\max}}, \frac{\e \mu}{4 \sigmaesc} } $, 
we get that $\ecn{x^K - x^\ast} = \mathcal{O}\br{\e}$ provided that the number of iterations is at least
\begin{equation}
	\label{eq:proxsgd-complexity}
	K_{\mathrm{SGD}} \geq \left( \kappa + \frac{\sigmaesc^2}{\e \mu^2} \right) \log\br{\frac{2 \|x^0 - x^\ast\|^2}{\e}}.
\end{equation}

By comparing between the iteration complexities $K_{\mathrm{SGD}}$ (given by \eqref{eq:proxsgd-complexity}) and $K_{\mathrm{RR}}$ (given by~\eqref{eq:proxrr-complexity-f-sc}), we see that \algname{ProxRR} converges faster than \algname{Proximal SGD} whenever the target accuracy $\e$ is small enough to satisfy
\[  \e \leq \frac{1}{L_{\max} n \mu} \br{\frac{\sigmaesc^4}{n \sqn{\nabla f(x^\ast)} + \sigmaesc^2}}.  \]
Furthermore, the comparison is much better when we consider proximal iteration complexity (number of proximal operator access), in which case the complexity of \algname{ProxRR}~\eqref{eq:proxrr-complexity-f-sc} is reduced by a factor of $n$ (because we take one proximal step every $n$ iterations) while the proximal iteration complexity of \algname{Proximal SGD} remains the same as~\eqref{eq:proxsgd-complexity}. In this case, \algname{ProxRR} is better whenever the accuracy $\e$ satisfies
\begin{equation*}
	\begin{split} 
        \e & \geq \frac{n}{L_{\max} \mu} \left [ n \sqn{\nabla f(x^\ast)} + \sigmaesc^2 \right ] \\ 
        \qquad \text {or}, 
		\qquad \e & \leq \frac{n}{L_{\max} \mu} \left [ \frac{\sigmaesc^4}{n \sqn{\nabla f(x^\ast)} + \sigmaesc^2} \right ].
	\end{split}
\end{equation*}
Therefore we can see that if the target accuracy is large enough or small enough, and if the cost of proximal operators dominates the computation, \algname{ProxRR} is much quicker to converge than \algname{Proximal SGD}.

\section{Convergence Theory for Strongly Convex Regularizer $\psi$}
In Theorem~\ref{thm:f-strongly-convex-psi-convex}, we assume that each $f_i$ is $\mu$-strongly convex. This is motivated by the common practice of using $\ell_2$ regularization in machine learning. However, applying $\ell_2$ regularization in every step of Algorithm~\ref{alg:proxrr} can be expensive when the data are sparse and the iterates $x^k_i$ are dense, because it requires accessing each coordinate of $x^k_i$ which can be much more expensive than computing sparse gradients $\nabla f_{i} (x^k_i)$. Alternatively, we may instead choose to put the $\ell_2$ regularization inside $\psi$ and only ask that $\psi$ be strongly convex---this way, we can save a lot of time as we need to access each coordinate of the dense iterates $x^k_i$ only once per epoch rather than every iteration. Theorem~\ref{thm:psi-strongly-convex-f-convex} gives a convergence guarantee in this setting.

\begin{theorem}
	\label{thm:psi-strongly-convex-f-convex}
	Let \Cref{as:smooth_fi_proper_psi} be satisfied. Further, assume that $\psi$ is $\mu$-strongly convex. If Algorithm~\ref{alg:proxrr} is run with constant stepsize $\gamma_k = \gamma \leq \frac{1}{L_{\max}}$, where $L_{\max} = \max_i L_i$, then the iterates generated by the algorithm satisfy
\[
			\ecn{x^T - x^\ast} \leq \br{1 + 2\gamma \mu n}^{-T} \sqn{x^0 - x^\ast} + \frac{ \gamma^2 \sigmarr^2}{\mu}.
\]
\end{theorem}

By making a specific choice for the stepsize used by Algorithm~\ref{alg:proxrr}, we can obtain a convergence guarantee using Theorem~\ref{thm:psi-strongly-convex-f-convex}. Choosing the stepsize as
\begin{equation}
	\label{eq:stepsize-choice}
	\gamma = \min \pbr{ \frac{1}{L_{\max}}, \frac{\sqrt{\e \mu}}{\sigmarr} }.
\end{equation}

Then $\ecn{x^T - x^\ast} = \mathcal{O}\br{\e}$ provided that the total number of iterations satisfies
\begin{equation}
	\label{eq:convergence-psi-sc-f-c}
	K \geq \left ( \kappa + \frac{\sigmarr/\mu}{\sqrt{\e \mu}} + n \right ) \log\br{\frac{2 \|x^0 - x^\ast\|^2}{\e}}.
\end{equation}
This can be converted to a bound similar to \eqref{eq:proxrr-complexity-f-sc} by using Theorem~\ref{thm:shuffling-radius-bound}, in which case the only difference between the two cases is an extra $n \log\br{\frac{1}{\e}}$ term when only the regularizer $\psi$ is $\mu$-strongly convex. Since for small enough accuracies the $1/\sqrt{\e}$ term dominates, this difference is minimal.

\section{Extensions}\label{sec:extensions}

Before turning to applications, we discuss two extensions to the theory that significantly matter in practice: using decreasing stepsizes and applying importance resampling.  

{\bf Decreasing stepsizes.}
Using the theoretical stepsize \eqref{eq:stepsize-choice} requires knowing the desired accuracy $\e$ ahead of time as well as estimating $\sigmarr$. It also results in extra polylogarithmic factors in the iteration complexity~\eqref{eq:convergence-psi-sc-f-c}, a phenomenon observed and fixed by using decreasing stepsizes in both vanilla \algname{RR}~\cite{Ahn2020} and in \algname{SGD}~\cite{stich2019unified}. We show that we can adopt the same technique to our setting. However, we depart from the stepsize scheme of \cite{Ahn2020} by only varying the stepsize once per epoch rather than every iteration. This is closer to the common practical heuristic of decreasing the stepsize once every epoch or once every few epochs~\cite{Sun2020, tran2020shuffling}. The stepsize scheme we use is inspired by the schemes of \cite{stich2019unified,Khaled2020}: in particular, we fix $T > 0$, 
let $k_0 = \ceil{T/2}$, and choose the stepsizes $\gamma_k > 0$ by
\begin{equation}
	\label{def:dec-stepsizes}
	\gamma_{k} = \begin{cases}
		\frac{1}{L_{\max}} & \text { if } T \leq \frac{L_{\max}}{2 \mu n} \text { or } k \leq k_0, \\
		\frac{7}{\mu n \br{s + k - k_0}} & \text { if } T > \frac{L_{\max}}{2 \mu n} \text { and } k > k_0,
	\end{cases}
\end{equation}
where $s \eqdef 7 L_{\max}/(4 \mu n)$. Hence, we fix the stepsize used in the first $T/2$ iterations and then start decreasing it every epoch afterwards. Using this stepsize schedule, we can obtain the following convergence guarantee when each $f_i$ is smooth and convex and the regularizer $\psi$ is $\mu$-strongly convex.

\begin{theorem}
	\label{thm:psi-strongly-cvx-dec-stepsizes}
	Suppose that each $f_i$ is $L_{\max}$-smooth and convex, and that the regularizer $\psi$ is $\mu$-strongly convex. Fix $T > 0$. Then choosing stepsizes $\gamma_k$ according to \eqref{def:dec-stepsizes} we have that $\gamma_k \leq \frac{1}{L_{\max}}$ for all $t$ and the final iterate generated by Algorithm~\ref{alg:proxrr} satisfies
	\[
			\ecn{x^{T} - x^\ast} = \mathcal{O}\br{ \exp\br{-\frac{n T}{\kappa + 2n}} \|x^0 - x^\ast\|^2 + \frac{\sigmarr^2}{\mu^3 n^2 T^2}},
	\]
	where $\kappa\eqdef L_{\max}/\mu$and $\mathcal{O}(\cdot)$ hides absolute (non-problem-specific) constants.
\end{theorem}

This guarantee holds for any number of epochs $T > 0$. We believe a similar guarantee can be obtained in the case each $f_i$ is strongly-convex and the regularizer $\psi$ is just convex, but we did not include it as it adds little to the overall message.

{\bf Importance resampling.}
Suppose that each $f_i$ is $L_i$-smooth. Then the iteration complexities of both \algname{SGD} and \algname{RR} depend on $L_{\max}/\mu$, where $L_{\max}$ is the maximum smoothness constant among the smoothness constants $L_1, L_2, \ldots, L_n$. The maximum smoothness constant can be arbitrarily worse than the average smoothness constant $\Lave = \frac{1}{n} \sum_{i=1}^{n} L_i$. This situation is in contrast to the complexity of gradient descent which depends on the smoothness constant $L_{f}$ of $f = \frac{1}{n} \sum_{i=1}^{n} f_i$, for which we have $L_{f} \leq \Lave$. This is a problem commonly encountered with stochastic optimization methods and may cause significantly degraded performance in practical optimization tasks in comparison with deterministic methods~\cite{tang2019practicality}.

Importance sampling is a common technique to improve the convergence of \algname{SGD} (Algorithm~\ref{alg:proxsgd}): we sample function $\frac{\Lave}{L_i}f_i$ with probability $p_i$ proportional to $L_i$, where $\Lave\eqdef \frac{1}{n}\sum_{i=1}^n L_i$. In that case, the \algname{SGD} update is still unbiased since 
\[
 \mathbb{E}_i\left[\frac{\Lave}{L_i}f_i\right] = \sum \limits_{i=1}^n p_i \frac{\Lave}{L_i}f_i=f .
 \]
Moreover, the smoothness of function $\frac{\Lave}{L_i}f_i$ is $\Lave$ for any $i$, so the guarantees would depend on $\Lave$ instead of $\max_{i=1,\dotsc, n}L_i$. Importance sampling successfully improves the iteration complexity of \algname{SGD} to depend on $\Lave$~\cite{needell2014stochastic}, and has been investigated in a wide variety of settings~\cite{gower2018stochastic,gorbunov2020unified}.

Importance sampling is a neat technique but it relies heavily on the fact that we use unbiased sampling. How can we obtain a similar result if inside any permutation the sampling is biased? The answer requires us to think again as to what happens when we replace $f_i$ with $\frac{\Lave}{L_i}f_i$. To make sure the problem remains the same, it is sufficient to have $\frac{\Lave}{L_i}f_i$ inside a permutation exactly $\frac{L_i}{\Lave}$ times. And since $\frac{L_i}{\Lave}$ is not necessarily integer, we should use $n_i=\Bigl\lceil\frac{L_i}{\Lave}\Bigr \rceil$ and solve
\begin{equation}
	\min \limits_{x\in\mathbb{R}^d} \frac{1}{N}\sum \limits_{i=1}^n \Bigl(\underbrace{\frac{1}{n_i}f_i(x) +\dotsb + \frac{1}{n_i}f_i(x)}_{n_i \text{ times}}  \Bigr) + \psi(x), \label{eq:importance}
\end{equation}
where $N\eqdef n_1+\dotsb +n_n = \Bigl\lceil\frac{L_1}{\Lave}\Bigr \rceil + \dotsb + \Bigl\lceil\frac{L_n}{\Lave}\Bigr \rceil$. Clearly, this problem is equivalent to the original formulation in~\ref{eq:prox-finite-sum-min}. At the same time, we have improved all smoothness constants to $\Lave$. It might seem that that the new problem has more functions, but it turns out that the new number of functions satisfies $N\le 2n$, so any related costs, such as longer loops or storing duplicates of the data, are negligible, as the next theorem shows.

\begin{theorem}\label{thm:IS}
	For every $i$, assume  that each $f_i$ is convex and $L_i$-smooth, and let $\psi$ be $\mu$-strongly convex. Then, the number of functions $N$ in~\eqref{eq:importance} satisfies $N\le 2n$, and Algorithm~\ref{alg:proxrr} applied to problem~\eqref{eq:importance} has the same complexity as \eqref{eq:convergence-psi-sc-f-c} but proportional to $\Lave$ rather than $L_{\max}$.
\end{theorem}

\section{Federated Learning} \label{sec:FL}

Let us consider now the problem of minimizing the average of $N= \sum_{m=1}^M N_m$ functions that are stored on $M$ devices, which have $N_1,\dotsc, N_M$ samples correspondingly,
\[
	\min \limits_{x\in\R^d} \frac{1}{N}\sum \limits_{m=1}^M F_m(x) + R(x), \quad
	F_m(x) = \sum \limits_{j=1}^{N_m} f_{mj}(x).\label{eq:fed_learning}
\]
For example, $f_{mj}(x)$ can be the loss associated with a single sample $(X_{mj}, y_{mj})$, where pairs $(X_{mj}, y_{mj})$ follow a distribution $D_m$ that is specific to device $m$. An important instance of such formulation is federated learning, where $M$ devices train a shared model by communicating periodically with a server. We normalize the objective in~\eqref{eq:fed_learning} by $N$ as this is the total number of functions after we expand each $F_m$ into a sum. We denote the solution of~\eqref{eq:fed_learning} by $x^\ast$.

{\bf Extending the space.}
To rewrite the problem as an instance of~\eqref{eq:prox-finite-sum-min}, we are going to consider a bigger product space, which is sometimes used in distributed optimization~\cite{bianchi2015coordinate}. Let us define $n\eqdef \max\{N_1, \dotsc, N_m\}$ and introduce $\psi_C$, the consensus constraint,
\[
	\psi_C(x_1,\dotsc, x_M)
	= \begin{cases}
	0, & x_1=\dotsb= x_M\\
	+\infty, & \text{otherwise}
	\end{cases}.
\]
By introducing dummy variables $x_1,\dotsc, x_M$ and adding the constraint $x_1=\dotsb =x_M$, we arrive at the intermediate problem
\[
	\min \limits_{x_1,\dotsc, x_M\in \R^p} \frac{1}{N}\sum \limits_{m=1}^M F_m(x_m) + (R + \psi_C)(x_1, \dotsc, x_M),
\]
where $R+\psi_C$ is defined, with a slight abuse of notation, as
\[
	(R+\psi_C)(x_1,\dotsc, x_M)
	= \begin{cases}
	R(x_1), & x_1=\dotsb= x_M\\
	+\infty, & \text{otherwise}.
	\end{cases}
\]
Since we have replaced $R$ with a more complicated regularizer $R+\psi_C$, we need to understand how to compute the proximal operator of the latter. We show (\Cref{lem:ext-proximal-operator} in the supplementary) that the proximal operator of $(R+\psi_C)$ is merely the projection onto $\{(x_1,\dotsc, x_M) \mid x_1=\dotsb = x_M\}$ followed by the proximal operator of $R$ with a smaller stepsize.

{\bf Reformulation.}
To have $n$ functions in every $F_m$, we write $F_m$ as a sum with extra $n-N_m$ zero functions, $f_{mj}(x)\equiv 0$ for any $ j > N_m$, so that
\[
	F_m(x_m) = \sum \limits_{j=1}^n f_{mj}(x_m) = \sum \limits_{j=1}^{N_m} f_{mj}(x_m)+\sum \limits_{j=N_m+1}^n 0.
\]
We can now stick the vectors together into $\xx=(x_1,\dotsc, x_M)\in\R^{M\cdot d}$ and multiply the objective by $\frac{N}{n}$, which gives the following reformulation:
\begin{equation}
	\min \limits_{\xx\in\R^{M\cdot d}} \frac{1}{n} \sum \limits_{i=1}^n f_i(\xx) + \psi(\xx), \label{eq:fed_reformulation}
\end{equation}
where $ \psi(\xx)\eqdef \frac{N}{n}(R+\psi_C)$ and
\begin{align*}
	& f_i(\xx) = f_i(x_1,\dotsc, x_M) \eqdef \sum \limits_{m=1}^M f_{mi}(x_m).
\end{align*}
In other words, function $f_i(\xx)$ includes $i$-th data sample from each device and contains at most one loss from every device, while $F_m(x)$ combines all data losses on device $m$. Note that the solution of~\eqref{eq:fed_reformulation} is $\xx^\ast\eqdef ((x^\ast)^\top, \dotsc, (x^\ast)^\top)^\top$ and the gradient of the extended function $f_i(\xx)$ is given by
\[
	\nabla f_i (\xx)	= (
		\nabla f_{1i}(x_1)^\top,
		\cdots , 
		\nabla f_{Mi}(x_M)^\top )^\top
\]
Therefore, a stochastic gradient step that uses $\nabla f_i(\xx)$ corresponds to updating all local models with the gradient of $i$-th data sample, without any communication.

\Cref{alg:proxrr} for this specific problem can be written in terms of $x_1,\dotsc, x_M$, which results in \Cref{alg:fed_rr}. Note that since $f_{mi}(x_i)$ depends only on $x_i$, computing its gradient does not require communication. Only once the local epochs are finished, the vectors are averaged as the result of projecting onto the set $\{(x_1,\dotsc, x_M) \mid x_1=\dotsb = x_M\}$. The full description of our \algname{FedRR} is given in~\Cref{alg:fed_rr}.
\begin{algorithm}[t]
    \caption{\algname{Federated Random Reshuffling} \rrbox{(\algname{FedRR})} and \algname{Shuffle-Once} \sobox{(\algname{FedSO})}}
    \label{alg:fed_rr}
\begin{algorithmic}[1]
   \Require Stepsize $\gamma > 0$, initial vector $x^0 = x_0^0 \in \R^d$, number of epochs $T$
   \State For each $m$, sample permutation $\prm{0, m}, \prm{1, m}, \ldots, \prm{N_m-1, m}$ of $\{ 1, 2, \ldots, N_m \}$ \sobox{(Only \algname{FedSO})}
    \For{epochs $k=0,1,\dotsc,T-1$}
    	  \For{$m=1,\dotsc, M$ locally in parallel}
    	  \State $x^{k, m}_0=x^k$
    	  \State Sample permutation $\prm{0, m}, \prm{1, m}, \ldots, \prm{N_m-1, m}$ of $\{ 1, 2, \ldots, N_m \}$ \rrbox{(Only \algname{FedRR})}
       \For{$i=0, 1, \ldots, N_m-1$}
          \State $x^{k, m}_{i+1} = x^{k, m}_{i} - \gamma \nabla f_{\prm{i, m}} (x^{k, m}_i)$
       \EndFor
       \State $x^{k, m}_n = x^{k,m}_{N_m}$
       \EndFor
       \State $z^{k+1}=\frac{1}{M}\sum_{m=1}^M x^{k, m}_n$
       \State $x^{k+1}=\prox_{\gamma \frac{N}{M} R}(z^{k+1})$
    \EndFor
\end{algorithmic}
\end{algorithm}

{\bf Reformulation properties.}
To analyze \algname{FedRR}, the only thing that we need to do is understand the properties of the reformulation~\eqref{eq:fed_reformulation} and then apply~\Cref{thm:f-strongly-convex-psi-convex} or \Cref{thm:psi-strongly-convex-f-convex}. The following lemma gives us the smoothness and strong convexity properties of \eqref{eq:fed_reformulation}.
\begin{lemma}\label{lem:fed_reform_properties}
	Let function $f_{mi}$ be $L_i$-smooth and $\mu$-strongly convex for every $m$. Then, $f_i$ from reformulation~\eqref{eq:fed_reformulation} is $L_i$-smooth and $\mu$-strongly convex.
\end{lemma}

The previous lemma shows that the conditioning of the reformulation is $\kappa=\frac{L_{\max}}{\mu}$ just as we would expect. Moreover, it implies that the requirement on the stepsize remains exactly the same: $\gamma\le \frac{1}{L_{\max}}$. What remains unknown is the value of $\sigmarr^2$, which plays a key role in the convergence bounds for \algname{ProxRR} and \algname{ProxSO}. Our next goal, thus, is to obtain an upper bound on $\sigmarr^2$, which would allow us to have a complexity for \algname{FedRR} and \algname{FedSO}. To find it, let us define
\[
	\sigma_{m, \ast}^2
	\eqdef \frac{1}{N_m}\sum \limits_{j=1}^{n} \bigl\|\nabla f_{mj}(x^\ast) - \frac{1}{N_m}\nabla F_m(x^\ast)\bigr\|^2,
\]

which is the variance of local gradients on device $m$. This quantity characterizes the convergence rate of \algname{Local SGD}~\cite{yuan2020federated}, so we should expect it to appear in our bounds too. The next lemma explains how to use it to upper bound $\sigmarr^2$.
\begin{lemma}\label{lem:fed_sigma}
	The shuffling radius $\sigmarr^2$ of the reformulation~\eqref{eq:fed_reformulation} is upper bounded by
	\[
		\sigmarr^2
		\le L_{\max} \sum \limits_{m=1}^M\Bigl( \|\nabla F_m(x^\ast)\|^2 + \frac{n}{4}\sigma_{m, \ast}^2\Bigr).
	\]
\end{lemma}

The lemma shows that the upper bound on $\sigmarr^2$ depends on the sum of local variances $\sum_{m=1}^M \sigma_{m,\ast}^2$ as well as on the local gradient norms $\sum_{m=1}^M\|\nabla F_m(x^\ast)\|^2$. Both of these sums appear in the existing literature on convergence of \algname{Local GD}/\algname{SGD}~\cite{khaled2019analysis, woodworth2020minibatch, yuan2020federated}.

Equipped with the variance bound, we are ready to present formal convergence results. For simplicity, we will consider heterogeneous and homogeneous cases separately and assume that $N_1=\dotsb=N_M=n$. To further illustrate generality of our results, we will present the heterogeneous assuming strong convexity $R$ and the homogeneous under strong convexity of functions $f_{mi}$.

{\bf Heterogeneous data.}
In the case when the data are heterogeneous, we provide the first local \algname{RR} method. We can apply either ~\Cref{thm:f-strongly-convex-psi-convex} or \Cref{thm:psi-strongly-convex-f-convex}, but for brevity, we give only the corollary obtained from \Cref{thm:psi-strongly-convex-f-convex}.
\begin{theorem}\label{thm:fed_hetero}
    Assume that functions $f_{mi}$ are convex and $L_i$-smooth for each $m$ and $i$. If $R$ is $\mu$-strongly convex and $\gamma\le \frac{1}{L_{\max}}$, then we have for the iterates produced by \Cref{alg:fed_rr}
    \begin{align*}
         \ecn{x^T - x^\ast} \leq \br{1 + 2\gamma \mu n}^{-T} \sqn{x^0 - x^\ast} + \frac{ \gamma^2 L_{\max}}{M\mu} \sum \limits_{m=1}^M\Bigl(  \|\nabla F_m(x^\ast)\|^2 + \frac{N}{4M}\sigma_{m, \ast}^2\Bigr).
    \end{align*}
\end{theorem}

\begin{figure}[t]
 	\centering
 	\includegraphics[scale=0.23]{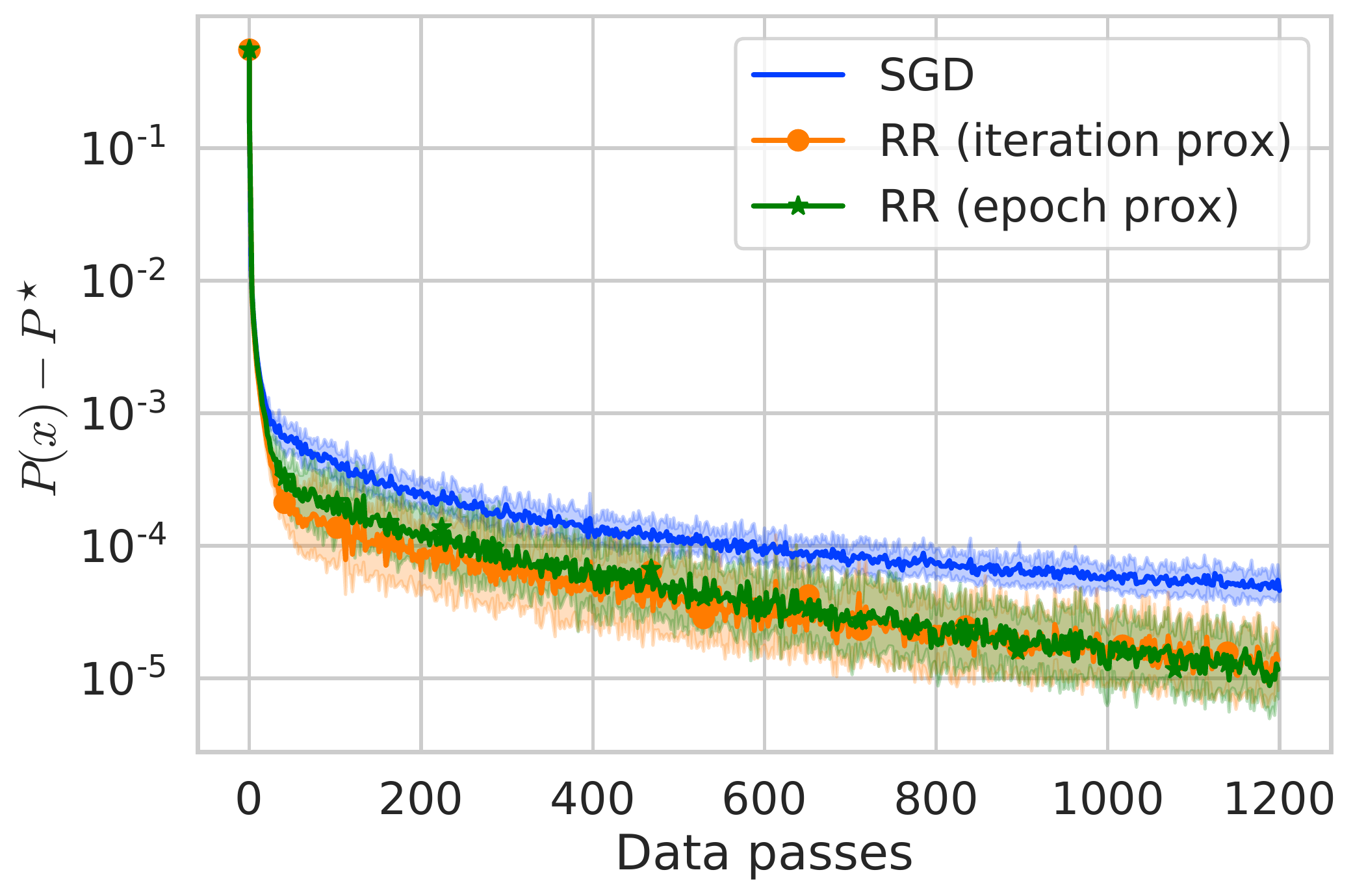}
 	\includegraphics[scale=0.23]{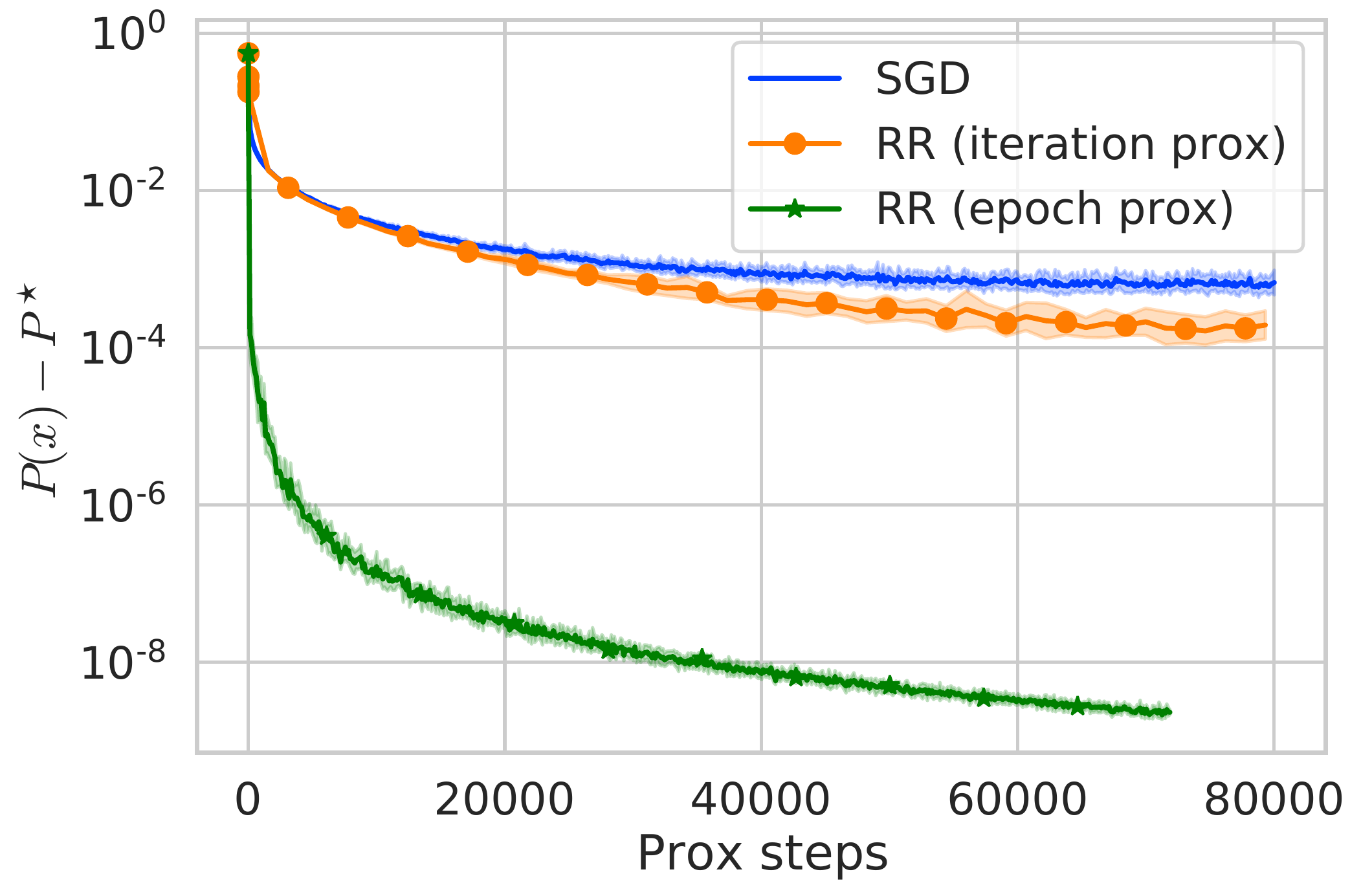}
 	\includegraphics[scale=0.23]{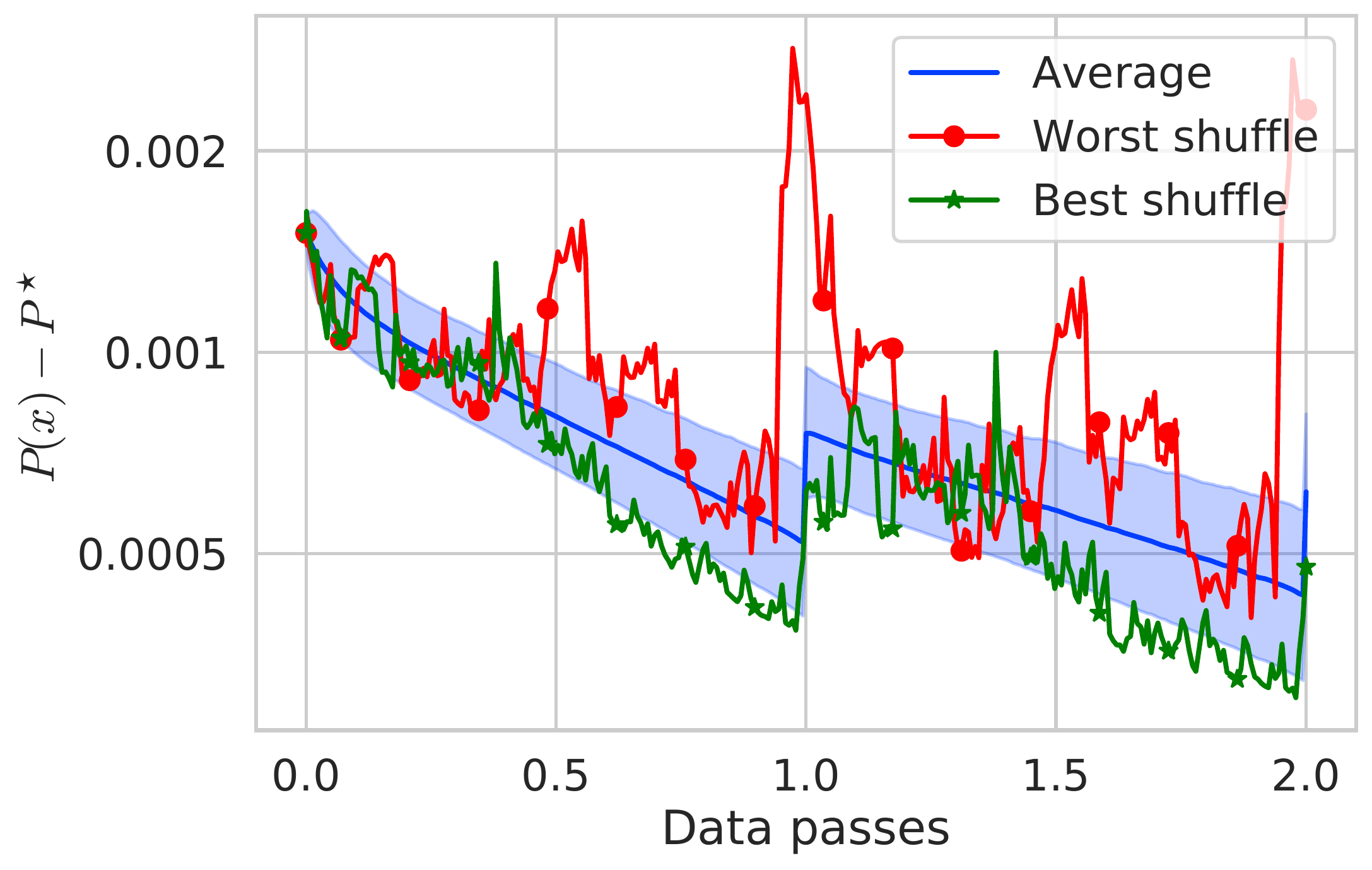}
 	\caption{Experimental results for problem \eqref{eq:exp_problem}. The first two plots show with average and confidence intervals estimated on 20 random seeds and clearly demonstrate that one can save a lot of proximal operator computations with our method. In the last plot, we show convergence oscillations of \algname{ProxSO} by sampling 20,000 permutations and choosing the best and the worst performing among them.}
 \end{figure}

{\bf Homogeneous data.}
For simplicity, in the homogeneous (i.e., i.i.d.) data case we provide guarantees without the proximal operator. Since then we have $F_1(x)=\dotsb = F_M(x)$, for any $m$ it holds $\nabla F_m(x^\ast)=0$, and thus $\sigma_{m, \ast}^2= \frac{1}{n}\sum_{j=1}^{n} \|\nabla f_{mj}(x^\ast)\|^2$. The full variance is then given by
\[
	\sum \limits_{m=1}^M \sigma_{m,\ast}^2 = \frac{1}{n}\sum \limits_{m=1}^M\sum \limits_{i=1}^n \|\nabla f_{mi}(x^\ast)\|^2
	= \frac{N}{n}\sigma_\ast^2
	= M\sigma_\ast^2,
\]
where $\sigma_\ast^2\eqdef \frac{1}{N}\sum_{i=1}^n\sum_{m=1}^M\|\nabla f_{mi}(x^\ast)\|^2$ is the variance of the gradients over all data.
\begin{theorem}\label{thm:fed_iid}
	Let $R(x)\equiv 0$ (no prox) and the data be i.i.d., that is $\nabla F_m(x^\ast)=0$ for any $m$, where $x^\ast$ is the solution of \eqref{eq:fed_learning}. If each $f_{mj}$ is $L_{\max}$-smooth and $\mu$-strongly convex, then the iterates of \Cref{alg:fed_rr} satisfy
	\[
		\ec{\|x^T - x^\ast\|^2}
		\le (1-\gamma\mu)^{nT}\|x^0-x^\ast\|^2 + \frac{\gamma^2L_{\max}N\sigma_\ast^2}{M\mu},
	\]
	where $\sigma_\ast^2\eqdef \frac{1}{N}\sum_{i=1}^n\sum_{m=1}^M\|\nabla f_{mi}(x^\ast)\|^2$.
\end{theorem}

The most important part of this result is that the last term in \Cref{thm:fed_iid} has a factor of $M$ in the denominator, meaning that the convergence bound improves with the number of devices involved.

\section[Experiments]{Experiments\footnote{Our code: \href{https://github.com/konstmish/rr_prox_fed}{https://github.com/konstmish/rr\_prox\_fed}}}

We look at the logistic regression loss with the elastic net regularization,
\begin{equation} \label{eq:exp_problem}
	\frac{1}{N}\sum \limits_{i=1}^N f_i (x) + \lambda_1\|x\|_1+ \frac{\lambda_2}{2}\|x\|^2,
\end{equation}
where each $f_i\colon \R^d \to \R$ is defined as
\[ f_i (x) \eqdef -\big(b_i \log \big(h(a_i^\top x)\big) + (1-b_i)\log\big(1-h(a_i^\top x)\big)\big) \]
and where $(a_i, b_i)\in \R^d\times \{0, 1\}$, $i=1,\dotsc, N$ are the data samples, $h\colon t\to1/(1+e^{-t})$ is the sigmoid function, and $\lambda_1, \lambda_2\ge 0$ are parameters.  We set mini-batch sizes to 1 for all methods and use theoretical stepsizes, without any tuning. We denote the version of \algname{RR} that performs proximal operator step after each iteration as `RR (iteration prox)'. We give more details in the supplementary. From the experiments, we can see that all methods behave more or less the same way. However, the algorithm that we propose needs only a small fraction of proximal operator evaluations, which gives it a huge advantage whenever the operator takes more time to compute than stochastic gradients.

\chapter{The First Adaptive Stepsize Rule for Gradient Descent that Provably Works}\label{chapter:adaptive}

\graphicspath{{adaptive/}}

\section{Introduction}
Since the early days of optimization it was evident that there is a
need for algorithms that are as independent from the
user as possible. First-order methods have proven to be versatile and
efficient in a wide range of applications, but one drawback has been
present all that time: the stepsize. Despite certain success
stories, line search procedures and adaptive online methods have not
removed the need to manually tune the optimization parameters. Even in
smooth convex optimization, which is often believed to be much simpler
than the non-convex counterpart, robust rules for stepsize selection
have been elusive. The purpose of this chapter is to remedy this deficiency.

The problem formulation that we consider is the basic unconstrained optimization problem
\begin{equation}
  \label{main}
  \min_{x\in \R^d}\ f(x),
\end{equation}
where $f\colon \R^d \to \R$ is  a differentiable
function. Throughout the chapter, we assume that \eqref{main} has a
solution and we denote its optimal value by $f^\star$.

The simplest and most known approach to this problem is the \algname{Gradient
Descent}method (\algname{GD}), whose origin can be traced back to
Cauchy~\cite{cauchy1847methode,lemarechal2012cauchy}. Although it is
probably the oldest optimization method, it continues to play a
central role in modern algorithmic theory and applications. Its
definition can be written in a mere one line,
\begin{equation}
    \label{eq:grad}
    x^{k+1} = x^k - \gamma \nabla f(x^k), \qquad k\geq 0,
\end{equation}
where $x^0\in \R^d$ is arbitrary and $\gamma >0$. Under  assumptions that
$ f$ is convex and $L$--smooth (equivalently, $\nabla f$ is $L$-Lipschitz)
that is
\begin{equation}\label{Lipschitz}
    \n{\nabla f(x)-\nabla f(y)}\leq L \n{x-y}, \quad \forall x,y,
\end{equation}
one can show that \algname{GD} with $\gamma \in (0, \frac 2 L)$ converges to an
optimal solution~\cite{polyak1963gradient}.  Moreover, with
$\gamma = \frac 1 L $ the convergence rate~\cite{drori2014performance} is
\begin{equation}
    \label{eq:grad_rate}
    f(x^k)-f^\star\leq \frac{L\n{x^0-x^\star}^2}{2(2k+1)},
\end{equation}
where $x^\star$ is any solution of \eqref{main}. Note that this bound is not improvable~\cite{drori2014performance}.

We identify four important challenges that limit the applications of
\algname{Gradient Descent} even in the convex case:
\begin{enumerate}
\itemsep0em
    \item \algname{GD} is not general: many functions do not satisfy \eqref{Lipschitz} globally.
    \item \algname{GD} is not a free lunch: one needs to guess $\gamma$,
    potentially trying many values before a success.
    \item \algname{GD} is not robust: failing to provide $\gamma < \frac{2}{L}$ may lead to divergence.
    \item \algname{GD} is slow: even if $L$ is finite, it might be arbitrarily
    larger than local smoothness.
\end{enumerate}

\section{Related Work}
Certain ways to address some of the issues above already exist in the literature. They include line search, adaptive \algname{Polyak's stepsize}, \algname{Mirror Descent}, \algname{Dual Preconditioning}, and stepsize estimation for subgradient methods. We discuss them one by one below, in a process reminiscent of cutting off Hydra's limbs: if one issue is fixed, two others take its place.

The most practical and generic solution to the aforementioned issues
is known as line search (or backtracking). This direction of research
started from the seminal works of Goldstein~\cite{goldstein1962cauchy}
and Armijo~\cite{armijo1966} and continues to attract attention, see
\cite{bello2016convergence,salzo2017variable} and references therein.
In general, at each iteration the line search executes another
subroutine with additional evaluations of $\nabla f$ and/or $f$ until
some condition is met. Obviously, this makes each iteration more
expensive.

At the same time, the famous \algname{Polyak's
stepsize}~\cite{polyak1969minimization} stands out as a very fast
alternative to \algname{Gradient Descent}. Furthermore, it does not depend on
the global smoothness constant and uses the current gradient to
estimate the geometry. The formula might look deceitfully simple,
$\gamma_k = \frac{f(x^k)-f^\star}{\n{\nabla f(x^k)}^2}$, but there is a
catch: it is rarely possible to know $f^\star$. This method, again,
requires the user to guess $f^\star$. What is more, with $\gamma$ it was fine to underestimate it by a factor of 10, but the guess for $f^\star$ must be tight, otherwise it has to be reestimated later~\cite{hazan2019revisiting}.

Seemingly no issue is present in the
\algname{Barzilai--Borwein} stepsize. Motivated by the quasi-Newton
schemes, Barzilai and Borwein~\cite{barzilai1988two} suggested using steps
\[
	\gamma_k = \frac{\lr{x^k-x^{k-1},\nabla f(x^k)-\nabla f(x^{k-1})}}{\n{\nabla f(x^k)-\nabla f(x^{k-1})}^2}.
\]
Alas, the convergence results regarding this choice of $\gamma_k$ are
very limited and the only known case where it provably works is
quadratic problems~\cite{raydan1993barzilai, dai2002r}. In general it
may not work even for smooth strongly convex functions, see the
counterexample in~\cite{burdakov2019stabilized}.

Other more interesting ways to deal with
non-Lipschitzness of $\nabla f$ use the problem structure. The first
method, proposed in~\cite{birnbaum2011distributed} and further
developed in~\cite{Bauschke2016}, shows that the \algname{Mirror Descent} method \cite{nemirovsky1983problem}, which is
another extension of \algname{GD}, can be used with a fixed stepsize, whenever
$f$ satisfies a certain generalization of~\eqref{Lipschitz}.
In addition, Maddison et al.~\cite{maddison2019dual} proposed the \algname{Dual
Preconditioning} method---another refined version of \algname{GD}. Similarly to
the former technique, it also goes beyond the standard smoothness
assumption of $f$, but in a different way. Unfortunately, these two simple and elegant
approaches cannot resolve all issues yet.
First, not many
functions fulfill respective generalized conditions. And secondly, both
methods still get us back to the problem of not knowing the allowed
range of stepsizes.

A whole branch of optimization considers adaptive extensions of \algname{GD} that deal with functions whose (sub)gradients are bounded. Probably the earliest work in that direction was written by Shor~\cite{shor1962application}. He showed that the method
\begin{align*}
        x^{k+1} = x^k - \gamma_k \frac{g^k}{\|g^k\|},
\end{align*}
where $g^k\in \partial f(x^k)$ is a subgradient, converges for
properly chosen sequences $(\gamma_k)_k$, see, e.g., Section 3.2.3 in~\cite{Nesterov2013}.
Moreover, $\gamma_k$ requires no knowledge about the function whatsoever.

 Similar methods that work in online setting such as \algname{Adagrad}~\cite{duchi2011adaptive,mcmahan2010adaptive} received a lot of attention in recent years and remain an active topic of research~\cite{ward19a}.
Methods similar to \algname{Adagrad}---\algname{Adam}~\cite{adam, adam2},
\algname{RMSprop}~\cite{Tieleman2012} and
\algname{Adadelta}~\cite{zeiler2012adadelta}---remain state-of-the-art for
training neural networks. The corresponding objective is usually
neither smooth nor convex, and the theory often assumes Lipschitzness
of the function rather than of the gradients. Therefore, this
direction of research is mostly orthogonal to ours, although we do
compare with some of these  methods in our neural networks experiment.

We also note that without momentum \algname{Adam} and \algname{RMSprop} reduce to \algname{signSGD}~\cite{bernstein2018signsgd}, which is known to be non-convergent for arbitrary stepsizes on a simple quadratic problem~\cite{karimireddy2019error}.

In a close relation to ours is the recent work~\cite{Malitsky2019},
where there was proposed an adaptive golden ratio algorithm for monotone
variational inequalities. As it solves a more general problem, it does not exploit the structure
of~\eqref{main} and, as most variational inequality methods, has a more conservative update. Although the method estimates the smoothness, it still requires an upper bound on the stepsize as input.

\paragraph{Contribution.} We propose a new version of \algname{GD} that at no
cost resolves all aforementioned issues. The idea is simple, and it is
surprising that it has not been yet discovered.  In each iteration we
choose $\gamma_k$ as a certain approximation of the inverse local Lipschitz
constant. With such a choice, we prove that convexity and local
smoothness of $f$ are sufficient for convergence of iterates with the
complexity  $\mathcal{O}(1/k)$ for $f(x^k)-f^\star$ in the worst case.

\paragraph{Discussion.} Let us now briefly discuss why we believe that proofs based on monotonicity and global smoothness lead to slower methods.

\algname{Gradient Descent} is by far not a recent method, so there have been
obtained optimal rates of convergence. However, we argue that adaptive
methods require rethinking optimality of the stepsizes. Take as an
example a simple quadratic problem,
 $f(x, y)=\frac{1}{2}x^2+\frac{\delta}{2}y^2$, where
$\delta \ll 1$. Clearly, the smoothness constant of this problem is
equal to $L=1$ and the strong convexity one is $\mu=\delta$. If we run
\algname{GD} from an arbitrary point $(x^0, y^0)$ with the ``optimal'' stepsize
$\gamma = \frac{1}{L} =1$, then one iteration of \algname{GD}  gives us $(x^1, y^1)
= (0, (1-\delta) y^0)$, and similarly $(x^k, y^k) = (0,
(1-\delta)^{k}y^0)$. Evidently for $\delta $ small enough it will take a long time to converge to the solution
$(0,0)$. Instead \algname{GD} would converge in two iterations if it adjusts its step after the first iteration to $\gamma = \frac{1}{\delta}$.

Nevertheless, all existing analyses of the \algname{Gradient Descent} with $L$-smooth $f$ use
stepsizes bounded by $2/L$. Besides, functional analysis gives
\[f(x^{k+1})\le f(x^k) - \gamma\Bigl(1 - \frac{\gamma
    L}{2}\Bigr)\|\nabla f(x^k)\|^2,\]
     from which $1/L$ can be seen as the
``optimal'' stepsize. Alternatively, we can assume that $f$ is
$\mu$-strongly convex, and the analysis in norms gives
\begin{align*}
\|x^{k+1}-x^\star\|^2\le \Bigl(1 - 2\frac{\gamma L \mu}{L + \mu} \Bigr)\|x^k-x^\star\|^2 - \gamma\Bigl(\frac{2}{L+\mu} - \gamma\Bigr)\|\nabla f(x^k) \|^2,
\end{align*}
whence the ``optimal'' step is $\frac{2}{L+\mu}$.

Finally,
line search procedures use some certain type of monotonicity, for
instance ensuring that $f(x^{k+1})\le f(x^k) - c\|\nabla f(x^k)\|^2$
for some $c>0$. We break with this tradition and merely ask for
convergence in the end.

\section{Convergence Theory}
\subsection{Local smoothness of $f$}\label{subs:main}
Recall that a mapping is locally Lipschitz if it is Lipschitz over
any compact set of its domain. A function $f$ with (locally) Lipschitz
gradient $\nabla f$ is called (locally) smooth.  It is natural to
ask whether some interesting functions are smooth locally, but not
globally.

It turns out there is no shortage of examples, most prominently among
highly nonlinear functions. In $\R$, they include
$x\mapsto \exp(x)$, $\log(x)$, $\tan(x)$, $x^p$, for $p > 2$,
etc. More generally, they include any twice differentiable $f$, since  $\nabla^2
f(x)$, as a continuous mapping, is  bounded over any bounded set
$\mathcal{C}$. In this case, we have that $\nabla f$ is Lipschitz on
$\mathcal{C}$, due to the mean value inequality
\[\n{\nabla f(x) - \nabla f(y)}\leq \max_{z\in \mathcal{C}}\n{\nabla^2 f(z)}\n{x-y},\quad \forall
    x,y\in \mathcal{C}.\]

Algorithm~\ref{alg:main} that we propose is just a slight modification of
 \algname{GD}. The quick explanation why local Lipschitzness of $\nabla f$
does not cause us any problems, unlike  most other methods, lies in
the way we prove its convergence. Whenever the stepsize $\gamma_k$
satisfies two inequalities\footnote{It can be shown that instead of the second condition it is enough to ask for $\gamma_k^2\le \frac{\n{x^{k}-x^{k-1}}^2}{[3\n{\nabla
        f(x^{k})}^2 - 4\<\nabla f(x^k), \nabla f(x^{k-1})>]_+}$, where $[a]_+\eqdef \max \{0, a\}$, but we prefer the option written in the main text for its simplicity.}
\begin{align*}
        \begin{cases}
  \gamma_k^2 & \leq (1+\th_{k-1})\gamma_{k-1}^2,\\ \gamma_k & \leq \frac{\n{x^{k}-x^{k-1}}}{2\n{\nabla
        f(x^{k})-\nabla f(x^{k-1})}},
        \end{cases}
\end{align*}
independently of the properties of $f$ (apart from convexity), we can show
that the iterates $(x^k)_k$ remain bounded. Here and everywhere else we
use the convention $1/0=+\infty$, so if $\nabla f(x^k) -\nabla
f(x^{k-1})=0$, the second inequality can be ignored. In the
first iteration, when $\theta_0=+\infty$, it might happen that $L_0=0$ and $\gamma_1 = \min\{+\infty, +\infty\}$. In
this case, we suppose that any choice of $\gamma_1>0$ is possible.
 \begin{algorithm}[t]
 \caption{\algname{Adaptive Gradient Descent} \rrbox{(\algname{AdGD})}}
 \label{alg:main}
 \begin{algorithmic}[1]
     \State \textbf{Input:} $x^0 \in \R^d$, $\gamma_0>0$,
     $\th_0=+\infty$, number of steps $K$
     \State  $x^1= x^0-\gamma_0\nabla f(x^0)$
        \For{$k = 1,2,\dotsc, K-1$}
        \State $\gamma_k = \min\Bigl\{
        \sqrt{1+\th_{k-1}}\gamma_{k-1},\frac{\n{x^{k}-x^{k-1}}}{2\n{\nabla
        f(x^{k})-\nabla f(x^{k-1})}}\Bigr\}$
		\State $x^{k+1} = x^k - \gamma_k \nabla f(x^k)$
		\State $\th_k = \frac{\gamma_k}{\gamma_{k-1}}$
        \EndFor
 \end{algorithmic}
 \end{algorithm}

 Although Algorithm~\ref{alg:main} needs $x^0$ and $\gamma_0$ as input,
 this is not an issue as one can simply fix $x^0=0$ and
 $\gamma_0=10^{-10}$. Equipped with a tiny $\gamma_0$, we ensure
 that $x^1$ will be close enough to $x^0$ and likely will give a good
 estimate for $\gamma_1$. Otherwise, this has no  influence on further steps.

\subsection{Analysis without descent}
It is now time to show our main contribution, the new analysis technique. The tools that we are going to use are the well-known Cauchy-Schwarz and convexity inequalities. In addition, our methods are related to potential functions~\cite{taylor19a}, which is a powerful tool for producing tight bounds for \algname{GD}.

Another divergence from the common practice is that our main lemma includes not only $x^{k+1}$ and $x^k$, but also $x^{k-1}$. This can be seen as a two-step analysis, while the majority of optimization methods have one-step bounds. However, as we want to adapt to the local geometry of our objective, it is rather natural to have two terms to capture the change in the gradients.

Now, it is time to derive a characteristic inequality for a specific
Lyapunov energy.
\begin{lemma}\label{lemma:energy}
Let $f\colon \R^d\to \R$ be convex and differential and let $x^\star$ be any solution of \eqref{main}.  Then for $(x^k)_k$ generated
by Algorithm~\ref{alg:main} it holds
\begin{multline}
  \label{eq:lemma_ineq}
\n{x^{k+1}-x^\star}^2+ \frac 1 2 \n{x^{k+1}-x^k}^2  + 2\gamma_{k}(1+\th_{k})
 (f(x^k)-f^\star)  \\ \leq \n{x^k-x^\star}^2  + \frac 1 2
 \n{x^k-x^{k-1}}^2    + 2\gamma_k \th_k (f(x^{k-1})-f^\star).
\end{multline}
\end{lemma}

\begin{proof}
    Let $k\geq 1$. We start from the standard way of analyzing \algname{GD}:
    \begin{align*}
      \|x^{k+1}- x^\star\|^2
      &= \|x^k - x^\star\|^2 + 2\<x^{k+1} - x^{k}, x^k-x^\star>+ \|x^{k+1} - x^{k}\|^2\\
      &= \|x^k - x^\star\|^2 + 2\gamma_k \<\nabla f(x^k), x^\star - x^k>  + \|x^{k+1} - x^{k}\|^2.
    \end{align*}
    As usually, we bound the scalar product by convexity of $f$:
    \begin{align}\label{eq:conv}
       2\gamma_k \<\nabla f(x^k), x^\star - x^k> \le 2\gamma_k (f^\star - f(x^k)),
    \end{align}
    which gives us
    \begin{align}\label{eq:norms}
          \|x^{k+1}- x^\star\|^2
          \le \|x^k - x^\star\|^2 - 2\gamma_k(f(x^k)-f^\star)  + \|x^{k+1} - x^{k}\|^2.
    \end{align}
    These two steps have been repeated thousands of times, but now we continue in a completely different manner. We have precisely one ``bad'' term in~\eqref{eq:norms}, which
    is $\n{x^{k+1}-x^k}^2$. We will bound it
    using the difference of gradients:
    \begin{align}\label{dif_x}
      \|x^{k+1} -x^k\|^2 & = 2 \n{x^{k+1}-x^k}^2 -
      \n{x^{k+1}-x^k}^2 \\
      & = -2\gamma_k \lr{\nabla f(x^k),
                               x^{k+1}-x^k}- \n{x^{k+1}-x^k}^2\notag
                               \\ &= 2\gamma_k \lr{\nabla f(x^k)-\nabla
                                   f(x^{k-1}), x^{k}-x^{k+1}}\notag\\
                               & \qquad \qquad +
                                   2\gamma_k\lr{\nabla f(x^{k-1}), x^k-x^{k+1}}
      -  \n{x^{k+1}-x^k}^2.
  \end{align}
Let us estimate the first two terms in the right-hand
    side above. First,   definition of $\gamma_k$, followed by
Cauchy-Schwarz and Young's inequalities, yields
\begin{align}\label{cs}
   2\gamma_k \lr{\nabla f(x^k) -\nabla f(x^{k-1}), x^k - x^{k+1}}  & \leq
2\gamma_k \n{\nabla f(x^k) -\nabla f(x^{k-1})} \n{x^k - x^{k+1}} \notag \\ &\leq
\n{x^k -x^{k-1}} \n{x^k - x^{k+1}} \notag \\ &\leq
  \frac 1 2 \n{x^{k}-x^{k-1}}^2 + \frac{1}{2}\n{x^{k+1}-x^k}^2.
\end{align}
Secondly, by convexity of $f$,
    \begin{align}    \label{eq:terrible_simple}
        2\gamma_k\lr{\nabla f(x^{k-1}), x^k-x^{k+1}}
      &=  \frac{2\gamma_k}{\gamma_{k-1}}\lr{x^{k-1} - x^{k},
        x^{k}-x^{k+1}}\notag  =  2\gamma_k\th_k \lr{x^{k-1}-x^{k}, \nabla
        f(x^k)}  \notag \\ & \leq  2\gamma_k\th_k
    (f(x^{k-1})-f(x^k)).
\end{align}
Plugging~\eqref{cs} and \eqref{eq:terrible_simple}
in~\eqref{dif_x}, we obtain
\begin{align*}
  \n{x^{k+1}-x^k}^2 \leq \frac 1 2 \n{x^k-x^{k-1}}^2 - \frac 1 2
  \n{x^{k+1}-x^k}^2 + 2\gamma_k\th_k(f(x^{k-1})-f(x^{k})).
\end{align*}
Finally, using the produced estimate for $\n{x^{k+1}-x^k}^2$ in
\eqref{eq:norms}, we deduce the desired inequality~\eqref{eq:lemma_ineq}.
\end{proof}

The above lemma already might give a good hint why our method
works. From inequality~\eqref{eq:lemma_ineq} together with condition
$\gamma_k^2\leq (1+\th_{k-1})\gamma_{k-1}^2$, we obtain that the Lyapunov
energy---the left-hand side of \eqref{eq:lemma_ineq}---is
decreasing. This gives us boundedness of $(x^k)_k$, which is often the
key ingredient for proving convergence. In the next theorem we
formally state our result.
\begin{theorem}\label{th:main}
    Suppose that $f\colon \R^d\to \R$ is convex with locally Lipschitz
    gradient $\nabla f$. Then $(x^k)_k$ generated by Algorithm~\ref{alg:main} converges to a solution of \eqref{main} and we have
    that
\[
	f(\hat x^K)-f^\star \leq \frac{D}{2S_k}=\mathcal{O}\Bigl(\frac{1}{K}\Bigr),
\]
where
\begin{align*}
	\hat x^K &=  \frac{\gamma_K(1+\th_K)x^K +
        \sum_{i=1}^{K-1}w_i
           x^i}{S_K},\\
    w_i &= \gamma_i(1+\th_i)-\gamma_{i+1}\th_{i+1},\\
  S_K &= \gamma_K(1+\th_K) + \sum_{i=1}^{K-1}w_i = \sum_{i=1}^K \gamma_i + \gamma_1\th_1,
\end{align*}
and $D$ is a constant that explicitly depends on the initial data and
the solution set, see \eqref{eq:telescope}.

\end{theorem}
Our proof will consist of two parts. The first one is a
straightforward application of Lemma~\ref{lemma:energy}, from which we
derive boundedness of $(x^k)_k$ and complexity result. Due to its
conciseness, we provide it directly after this remark. In the second
part, we prove that the whole sequence $(x^k)_k$ converges to a
solution. Surprisingly, this part is a bit more technical than
expected, and thus we postpone it to the appendix.

\begin{proof} \textit{(Boundedness and complexity result.)}

    Fix any $x^\star$ from the solution set of \cref{main}.  Telescoping
    inequality~\eqref{eq:lemma_ineq}, we deduce
    \begin{multline}\label{eq:telescope}
        \n{x^{k+1}-x^\star}^2+ \frac 1 2 \n{x^{k+1}-x^k}^2  + 2\gamma_k
  (1+\th_k) (f(x^k)-f^\star) \\+
  2\sum_{i=1}^{k-1}[\gamma_i(1+\th_i)-\gamma_{i+1}\th_{i+1}](f(x^i)-f^\star) \\
  \leq \, \n{x^1-x^\star}^2  + \frac 1 2 \n{x^1-x^{0}}^2   + 2\gamma_1 \th_1 [f(x^{0})-f^\star]\eqdef D.
\end{multline}
Note that by definition of $\gamma_k$, the second line above is always
nonnegative. Thus, the sequence $(x^k)_k$ is bounded. Since $\nabla f$
is locally Lipschitz, it is Lipschitz continuous on bounded sets. It
means that for the set
$\mathcal{C} = \clconv \{x^\star, x^0, x^1,\dots\}$, which is bounded as the convex hull of bounded points, there exists $L>0$ such that
\[\n{\nabla f(x)-\nabla f(y)} \leq L \n{x-y} \quad \forall
    x,y\in \mathcal{C}.\]
    Clearly, $\gamma_1=\frac{\n{x^1-x^0}}{2\n{\nabla
        f(x^1)-\nabla f(x^0)}}\geq
\frac{1}{2L}$, thus, by induction one can prove that
\(\gamma_k \geq \frac{1}{2L} \), in other words, the sequence
$(\gamma_k)_k$ is separated from zero.

Now we want to apply the Jensen's inequality for the sum of all terms
$f(x^i)-f^\star$ in the left-hand side of \eqref{eq:telescope}.  Notice, that the
total sum of coefficients at these terms is
\[\gamma_k(1+\th_k)
    +\sum_{i=1}^{k-1}[\gamma_i(1+\th_i)-\gamma_{i+1}\th_{i+1}]
    =\sum_{i=1}^k \gamma_i + \gamma_1\th_1=S_k\]
Thus, by Jensen's inequality,
\[\frac D 2 \geq \frac{\text{LHS of \eqref{eq:telescope}}}{2} \geq S_k (f(\hat
    x^k)-f^\star),\]
where $\hat x^k$ is given in the statement of the theorem.
By this, the first part of the proof is complete. Convergence of $(x^k)_k$ to a solution is provided in the appendix.
\end{proof}
As we have shown that $\gamma_i\geq \frac{1}{2L}$ for all $i$, we have a
theoretical upper bound $f(\hat x^k)-f^\star\leq \frac{D L}{k}$. Note that
in practice, however, $(\gamma_k)$ might be much larger than the
pessimistic lower bound $\frac{1}{2L}$, which we observe in our experiments
together with a faster convergence.

\subsection{$f$ is locally strongly convex}
Since one of our goals is to make optimization easy to use, we believe that
a good method should have state-of-the-art guarantees in various scenarios. For strongly convex functions, this means that we want to see linear convergence, which is not covered by normalized \algname{GD} or online methods. In
\cref{subs:main} we have shown that Algorithm~\ref{alg:main} matches
the $\mathcal{O}(1/\e)$ complexity of \algname{GD} on convex problems. Now we show that it also matches
$\mathcal{O}(\frac{L}{\mu}\log\frac{1}{\e})$ complexity of \algname{GD} when $f$
is locally strongly convex. Similarly to local smoothness, we
call $f$ locally strongly convex if it is strongly convex over any
compact set of its domain.

For proof simplicity, instead of using bound
$\gamma_k \leq \sqrt{1+\th_{k-1}}\gamma_{k-1}$ as in step~4 of
Algorithm~\ref{alg:main} we will use a more conservative bound
$\gamma_k\leq \sqrt{1+\frac{\th_{k-1}}{2}}\gamma_{k-1}$ (otherwise the
derivation would be too technical). It is clear that with such a
change \Cref{th:main} still holds true, so the sequence is bounded and we can rely on local smoothness and local strong convexity.
\begin{theorem}\label{th:strong}
    Suppose that $f\colon \R^d\to \R$ is locally strongly convex and
    $\nabla f$ is locally Lipschitz.  Then $(x^k)_k$ generated by
    Algorithm~\ref{alg:main} (with the modification mentioned above)
    converges to the solution $x^\star$ of \eqref{main}. The complexity to
    get $\|x^k - x^\star\|^2\le \varepsilon$ is
    $\mathcal{O}(\kappa \log\frac 1 \e)$, where
    $\kappa = \frac{L}{\mu}$ and $L,\mu$ are the smoothness and strong
    convexity constants of $f$ on the set
    $\mathcal{C}=\clconv\{x^\star, x^0, x^1, \dots\}$.
\end{theorem}

We want to highlight that in our rate $\kappa$ depends on the local
Lipschitz and strong convexity constants $L$ and $\mu$, which is
meaningful even when these properties are not satisfied
globally. Similarly, if $f$ is globally smooth and strongly convex,
our rate is still faster as it depends on the smaller local constants.

\section{Heuristics}
In this section, we describe several extensions of our method. We do
not have a full theory for them, but believe that they are of interest
in applications.
\subsection{Acceleration}
\begin{algorithm}[t]
 \caption{\algname{Adaptive Accelerated Gradient Descent} \rrbox{(\algname{AdGD-accel})}}
 \label{alg:accel}
 \begin{algorithmic}[1]
     \State \textbf{Input:} $x^0 \in \R^d$, $\gamma_0>0$, $\gamma_0>0$,
     $\th_0=\Theta_0 = +\infty$, number of steps $K$\\
     \State $y^1 = x^1= x^0-\gamma_0\nabla f(x^0)$
        \For{$k = 1,2,\dots, K-1$}
        \State $\gamma_k = \min\Bigl\{
        \sqrt{1+\frac{\th_{k-1}}{2}}\gamma_{k-1},\frac{\n{x^{k}-x^{k-1}}}{2\n{\nabla
                f(x^{k})-\nabla f(x^{k-1})}}\Bigr\}$
        \State $\mu_k = \min\Bigl\{\sqrt{1 +
                \frac{\Theta_{k-1}}{2}}\mu_{k-1}, \frac{\| \nabla f(x^k) -
                \nabla f(x^{k-1}) \|}{2\|x^k - x^{k-1}\|} \Bigr\}$
        \State $\beta_k =\frac{\sqrt{1/\gamma_k} - \sqrt{\mu_k}}{\sqrt{1/\gamma_k}+\sqrt{\mu_k}} $
        \State $y^{k+1} = x^k - \gamma_k \nabla f(x^k)$
        \State $x^{k+1}=y^{k+1} + \beta_k (y^{k+1}-y^k)$
\State $\th_k = \frac{\gamma_k}{\gamma_{k-1}}$, $\Theta_k = \frac{\mu_k}{\mu_{k-1}}$
        \EndFor
 \end{algorithmic}
 \end{algorithm}

Suppose that $f$ is $\mu$-strongly convex.
One version of the  accelerated gradient method proposed by
Nesterov~\cite{Nesterov2013} is
\begin{align*}
        y^{k+1} &= x^k  - \frac{1}{L}\nabla f(x^k),\\
        x^{k+1}&= y^{k+1} + \beta(y^{k+1} - y^{k}),
\end{align*}
where $\beta = \frac{\sqrt{L} - \sqrt{\mu}}{\sqrt{L}+\sqrt{\mu}}$.
\algname{Adaptive Gradient Descent} for
 strongly convex $f$ efficiently estimated $\frac{1}{2L}$ by
\begin{align*}
        \gamma_k = \min\biggl\{\sqrt{1 + \frac{\th_{k-1}}{2}}\gamma_{k-1}, \frac{\|x^k - x^{k-1}\|}{2\|\nabla f(x^k) -\nabla f(x^{k-1})\|}\biggr\}.
\end{align*}
What about the strong convexity constant $\mu$?  We know that it
 equals to the inverse smoothness constant of
the conjugate $f^*(y) \eqdef \sup_x\{\lr{x,y} - f(x) \}$. Thus, it is
tempting to estimate this inverse constant just as we estimated inverse
smoothness of $f$, i.e., by formula
\begin{align*}
        \mu_k = \min\biggl\{\sqrt{1 +\frac{\Theta_{k-1}}{2}} \mu_{k-1}, \frac{\| p^k - p^{k-1} \|}{2\|\nabla f^*(p^k) - \nabla f^*(p^{k-1})\|} \biggr\}
\end{align*}
where $p^k$ and $p^{k-1}$ are some elements of the dual space and
$\Theta_k = \frac{\mu_k}{\mu{k-1}}$. A natural choice then is $p^k = \nabla f(x^k)$ since it is an element of the dual space that we use. What is its value? It is well known that $\nabla f^*(\nabla f(x))=x$, so we come up with the update rule
\begin{align*}
        \mu_k = \min\biggl\{\sqrt{1 + \frac{\Theta_{k-1}}{2}}\mu{k-1}, \frac{\| \nabla f(x^k) - \nabla f(x^{k-1}) \|}{2\|x^k - x^{k-1}\|} \biggr\},
\end{align*}
and hence we can estimate $\beta$ by $
        \beta_k = \frac{\sqrt{1/\gamma_k} - \sqrt{\mu_k}}{\sqrt{1/\gamma_k}+\sqrt{\mu_k}}$.

We summarize our arguments in
Algorithm~\ref{alg:accel}. Unfortunately, we do not have any
theoretical guarantees for it.

 Estimating strong convexity parameter $\mu$ is important in
 practice. Most common approaches rely on restarting technique
 proposed by Nesterov~\cite{Nesterov2013a}, see
 also~\cite{fercoq2017adaptive} and references therein. Unlike
 Algorithm~\ref{alg:accel}, these works have theoretical guarantees,
 however, the methods themselves are more complicated and still
 require tuning of other unknown parameters.

 \subsection{Uniting our steps with stochastic gradients}
Here we would like to discuss applications of our method to the problem
\begin{align*}
        \min_x \E{f(x; \xi)},
\end{align*}
where $f(\cdot; \xi)$ is almost surely $L$-smooth and $\mu$-strongly convex. Assume that at each iteration we get sample $\xi^k$ to make a stochastic gradient step,
\begin{align*}
x^{k+1}=x^k - \gamma_k\nabla f(x^k; \xi^k).
\end{align*}
Then, we have two ways of incorporating our stepsize into \algname{SGD}.
The first is to reuse $\nabla f(x^k; \xi^k)$ to estimate $L_k=\frac{\n{\nabla f(x^{k}; \xi^k)-\nabla
        f(x^{k-1};\xi^k)}}{\n{x^{k}-x^{k-1}}}$, but this would make
        $\gamma_k\nabla f(x^k; \xi^k)$ biased. Alternatively, one can use an extra
        sample to estimate $L_k$, but this is less intuitive since our goal is to estimate
        the curvature of the function used in the update.

        We give a full description in~\Cref{alg:stoch}. We remark that
        the option with a biased estimate performed much better in our
        experiments with neural networks. The theorem below provides
        convergence guarantees for both cases, but with different
        assumptions.

\begin{algorithm}[t]
 \caption{\algname{Adaptive SGD} \rrbox{(\algname{AdSGD})}}
 \label{alg:stoch}
 \begin{algorithmic}[1]
     \State \textbf{Input:} $x^0 \in \R^d$, $\gamma_0>0$,
     $\th_0=+\infty$, $\xi^0$, $\alpha>0$, number of steps $K$
     \State $x^1 = x^0-\gamma_0\nabla f(x^0; \xi^0)$
        \For{$k = 1,2,\dotsc, K-1$}
        \State Sample $\xi^k$ and optionally $\zeta^k$
        \State Option I (biased):  $L_k = \frac{\n{\nabla
        f(x^{k}; \xi^k)-\nabla f(x^{k-1};\xi^k)}}{\n{x^{k}-x^{k-1}}}$
        \State Option II (unbiased):  $L_k = \frac{\n{\nabla
        f(x^{k};\zeta^k)-\nabla f(x^{k-1};\zeta^{k-1})}}{\n{x^{k}-x^{k-1}}}$
        \State $\gamma_k = \min\Bigl\{
        \sqrt{1+\th_{k-1}}\gamma_{k-1},\frac{\alpha}{ L_k}\Bigr\}$
\State $x^{k+1} = x^k - \gamma_k \nabla f(x^k; \xi^k)$
\State $\th_k = \frac{\gamma_k}{\gamma_{k-1}}$
        \EndFor
 \end{algorithmic}
 \end{algorithm}

\begin{theorem}
        Let $f(\cdot;\xi)$ be $L$-smooth and $\mu$-strongly convex almost surely. Assuming $\alpha\le \frac{1}{2\kappa}$ and estimating $L_k$ with $\nabla f(\cdot; {\zeta^k})$, the complexity to get $\E{\n{x^k- x^\star}^2}\le \varepsilon$ is not worse than $\mathcal{O}\left(\frac{\kappa^2}{\varepsilon}\log \frac{\kappa}{\varepsilon}\right)$. Furthermore, if the model is overparameterized, i.e., $\nabla f(x^\star; \xi)=0$ almost surely, then one can estimate $L_k$ with $\xi^k$ and the complexity is $\mathcal{O}\left(\kappa^2 \log \frac{1}{\varepsilon}\right)$.
\end{theorem}
Note that in both cases we match the known dependency on $\varepsilon$ up to logarithmic terms, but we get an extra $\kappa$ as the price for adaptive estimation of the stepsize.

Another potential application of our techniques is estimation of decreasing stepsizes in \algname{SGD}. The best known rates for \algname{SGD}~\cite{stich2019unified}, are obtained using $\gamma_k$ that evolves as $\mathcal{O}\left(\frac{1}{L+\mu k}\right)$. This requires estimates of both smoothness and strong convexity, which can be borrowed from the previous discussion. We leave rigorous proof of such schemes for future work.
\section{Experiments}
In the experiments\footnote{See
  \href{https://github.com/ymalitsky/adaptive_GD}{https://github.com/ymalitsky/adaptive\_gd}},
we compare our approach with the two most related methods: \algname{GD} and
\algname{Nesterov's Accelerated method} for convex
functions~\cite{Nesterov1983a}. Additionally, we consider line search,
\algname{Polyak} step, and \algname{Barzilai--Borwein} method. For neural networks we also
include a comparison with \algname{SGD}, \algname{SGDm} and \algname{Adam}.

\begin{figure*}[!t]
\centering
\begin{subfigure}[t]{0.32\textwidth}
    {\includegraphics[trim={3mm 0 3mm 0},clip,width=1\textwidth]{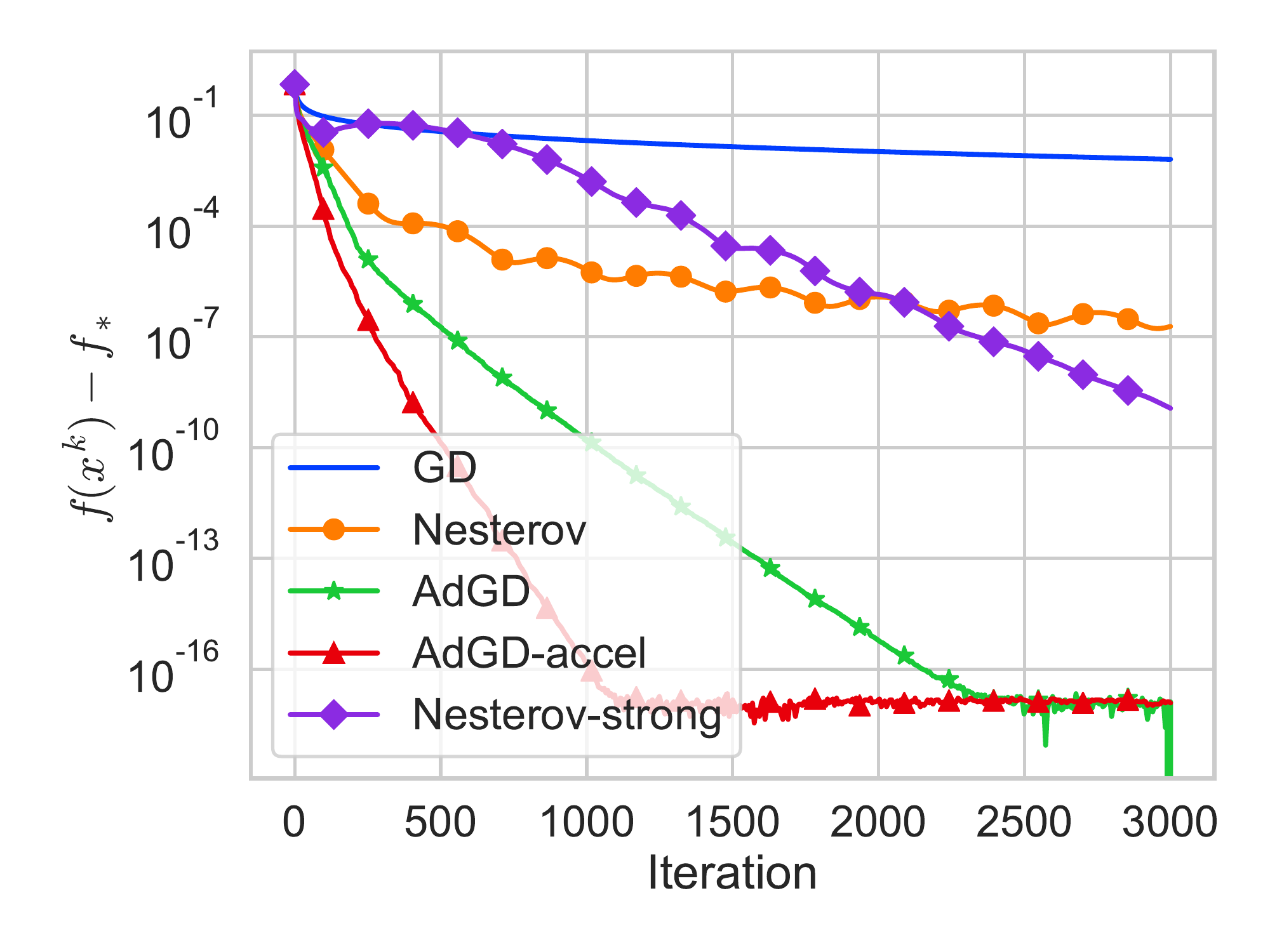}}
\caption{Mushrooms dataset, objective}
\end{subfigure}
\hfill
\begin{subfigure}[t]{0.32\textwidth}
{\includegraphics[trim={3mm 0 3mm 0},clip,width=1\textwidth]{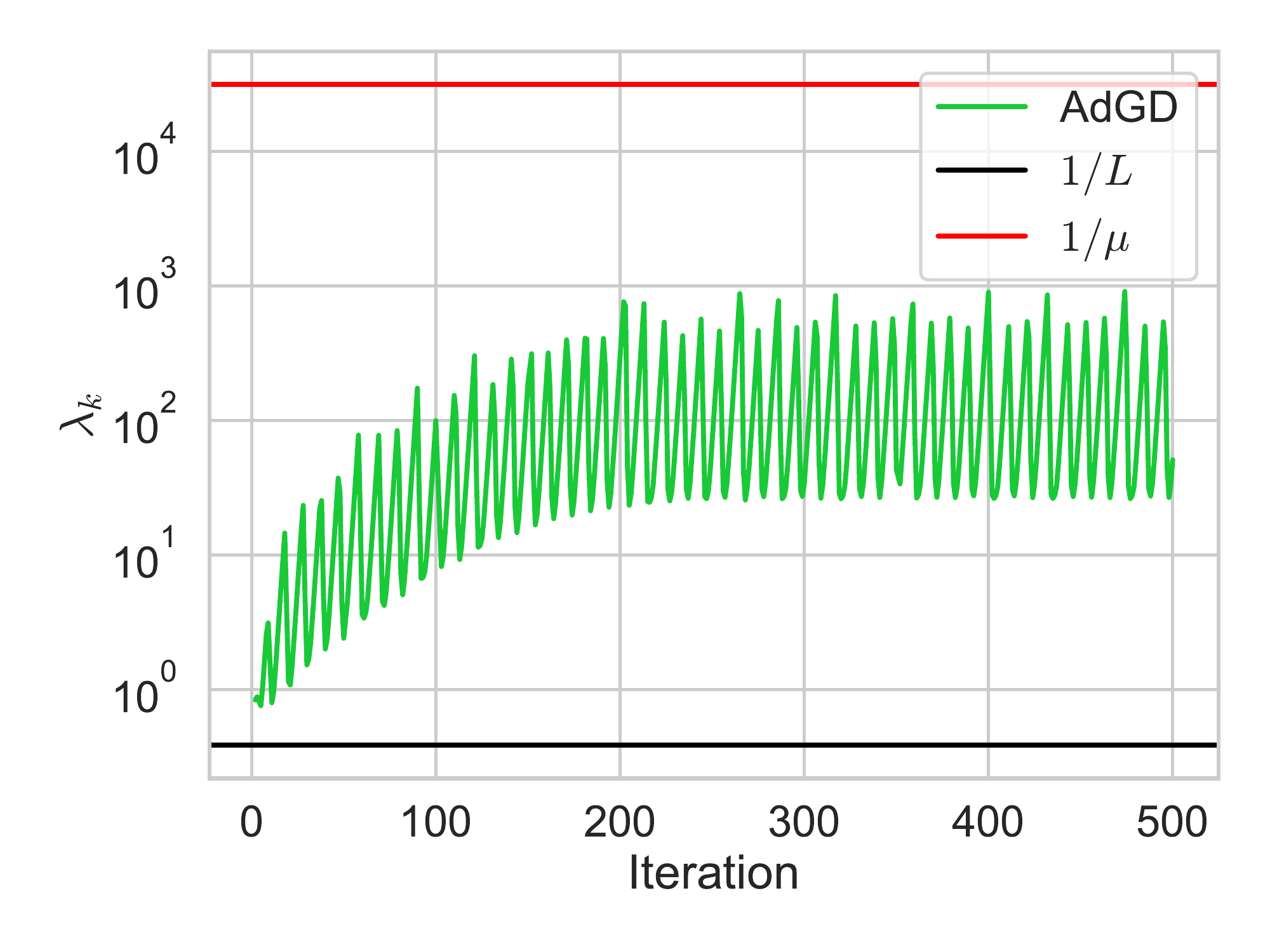}}
\caption{Mushrooms dataset, stepsize}
\end{subfigure}
\hfill
\begin{subfigure}[t]{0.32\textwidth}
{\includegraphics[trim={3mm 0 3mm 0},clip,width=1\textwidth]{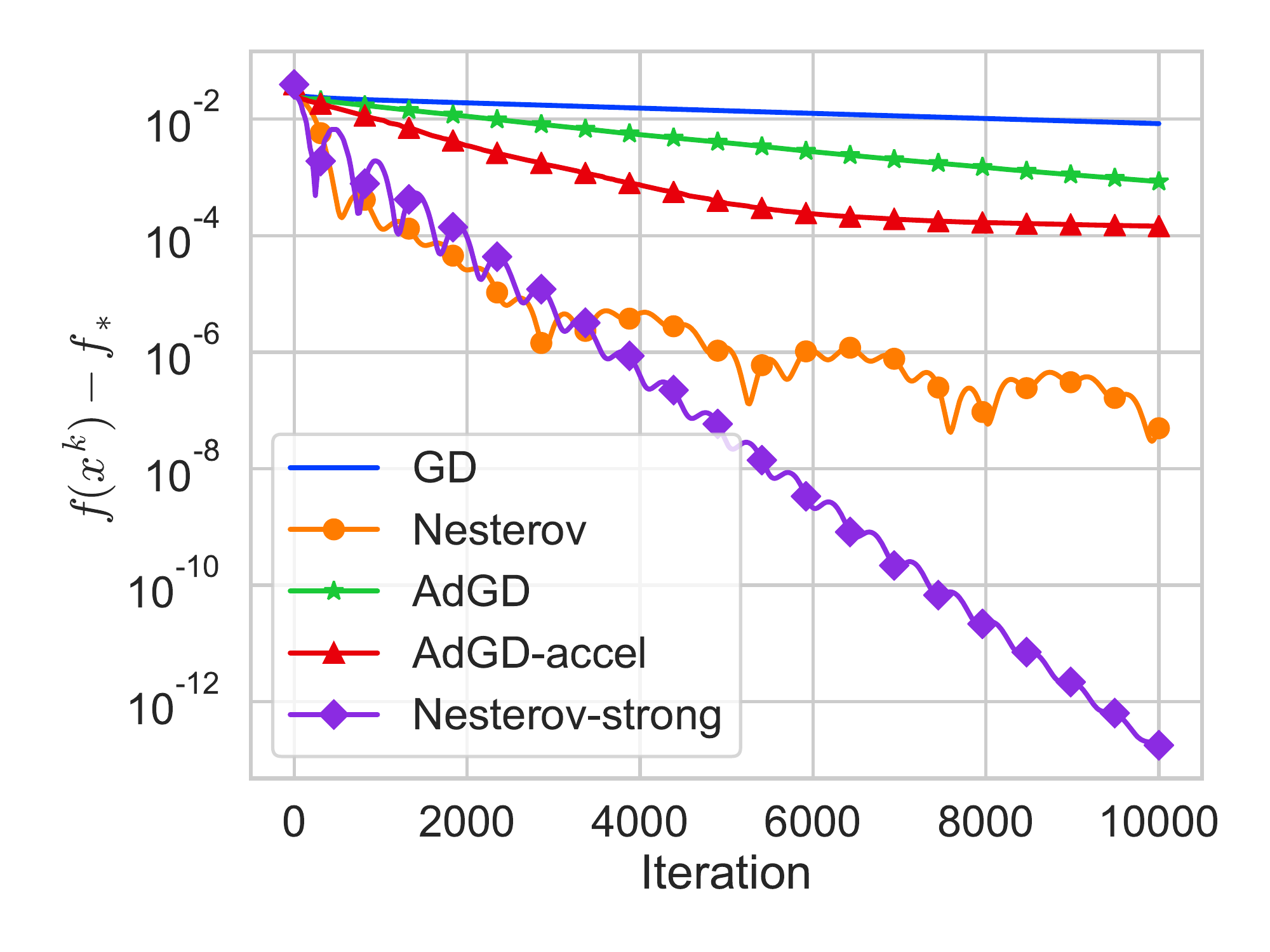}}
\caption{Covtype dataset, objective}
\end{subfigure}
\caption{Results for the logistic regression problem.}
\label{fig:logistic}
\end{figure*}

\paragraph{Logistic regression.}
The  logistic loss with $\ell_2$-regularization is given by
$\frac{1}{n}\sum_{i=1}^n \log(1 + \exp(-b_i a_i^\top x)) +
\frac{\lambda}{2}\|x\|^2$, where $n$ is the number of observations,
$\lambda>0$ is a regularization parameter, and
$(a_i, b_i)\in\R^{d}\times \R$, $i=1,\dots, n$, are the observations. We
use `mushrooms' and `covtype' datasets to run the experiments. We choose $\lambda$
proportionally to $\frac{1}{n}$ as often done in practice. Since we
have closed-form expressions to estimate $L=\frac{1}{4n}\|\mA\|^2+\lambda$, where $\mA=(a_1^\top, \dotsc, a_n^\top)^\top$, we used stepsize $\frac{1}{L}$ in \algname{GD} and its acceleration. The
results are provided in Figure~\ref{fig:logistic}.

\paragraph{Matrix factorization.}

\begin{figure*}[t]
\centering
\begin{subfigure}[t]{0.32\textwidth}
    {\includegraphics[trim={3mm 0 3mm 0},clip,width=1\textwidth]{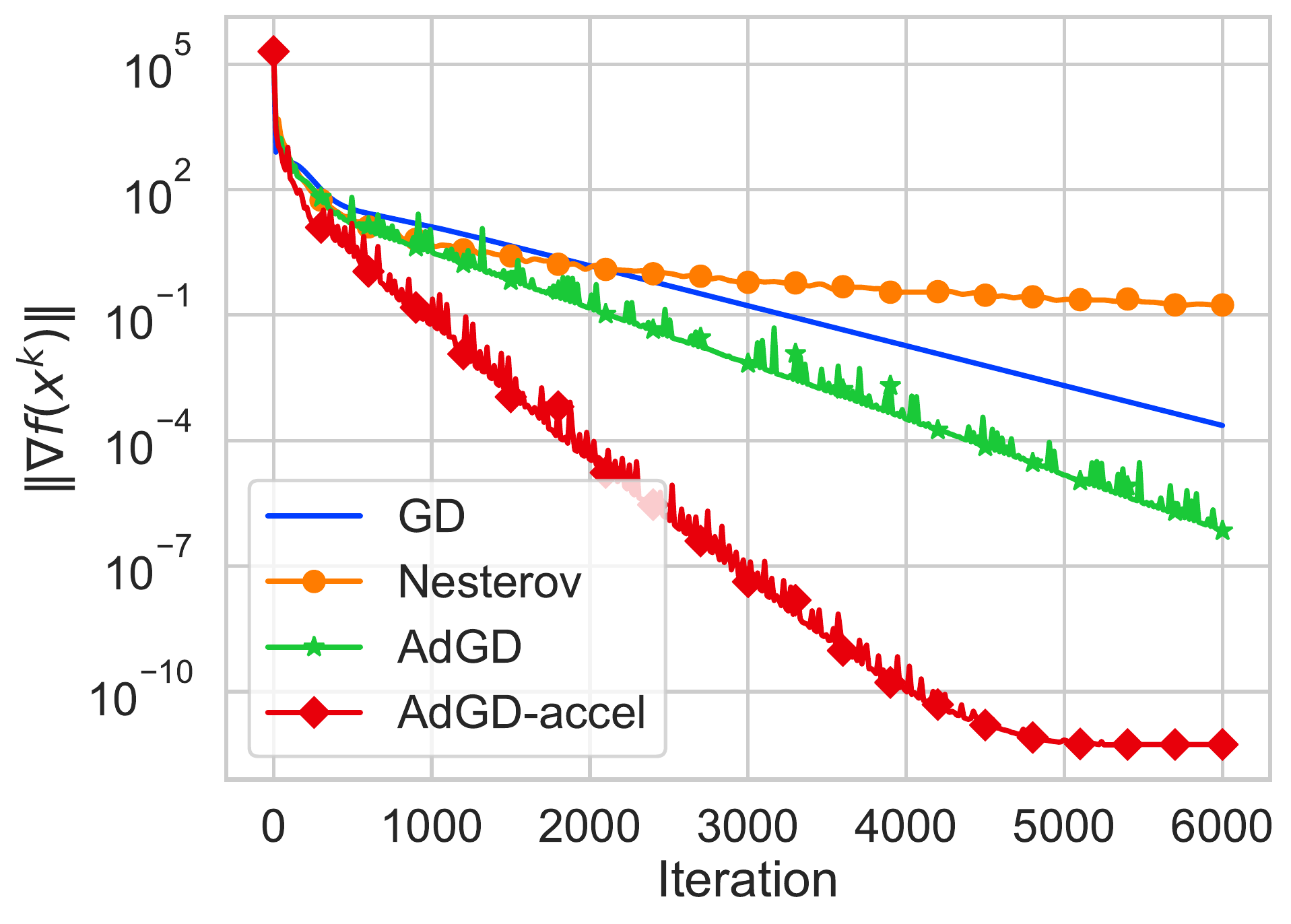}}
\caption{$r=10$}
\end{subfigure}
\hfill
\begin{subfigure}[t]{0.32\textwidth}
{\includegraphics[trim={3mm 0 3mm 0},clip,width=1\textwidth]{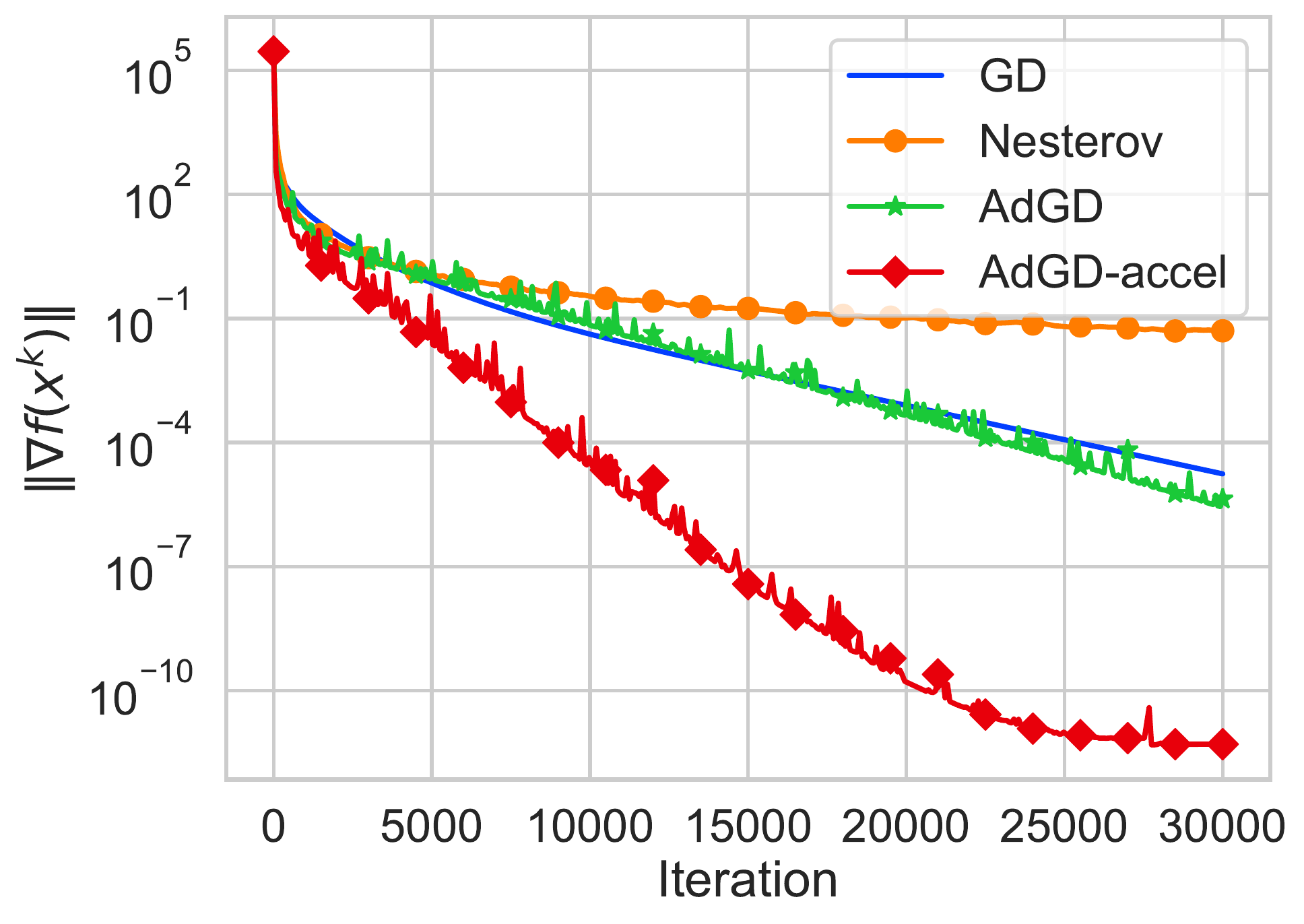}}
\caption{$r=20$}
\end{subfigure}
\hfill
\begin{subfigure}[t]{0.32\textwidth}
{\includegraphics[trim={3mm 0 3mm 0},clip,width=1\textwidth]{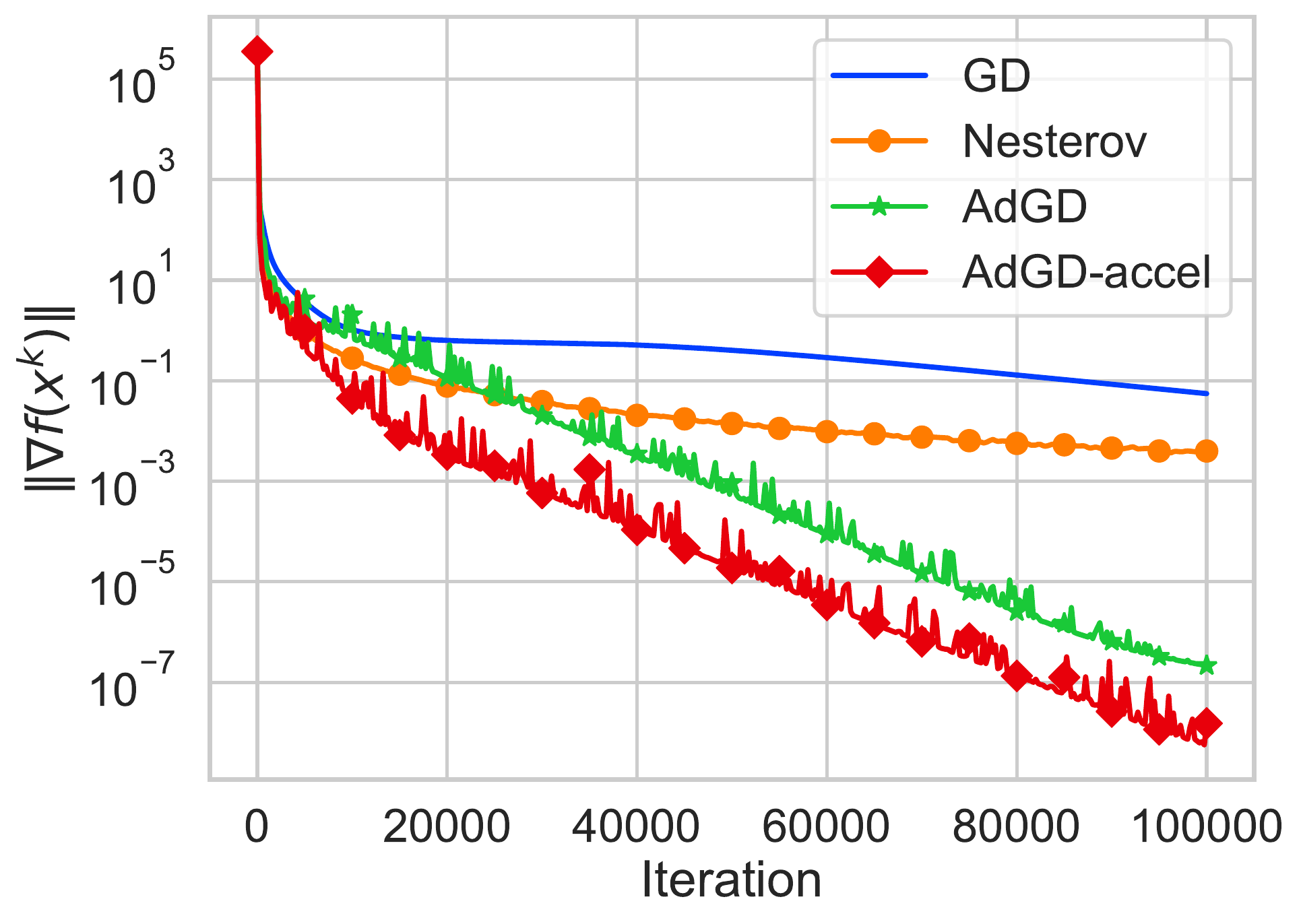}}
\caption{$r=30$}
\end{subfigure}
\caption{Results for matrix factorization. The objective is neither convex nor
    smooth.}
\label{fig:factorization}
\end{figure*}

Given a matrix $\mA\in \R^{m\times n}$ and $r<\min\{m,n\}$, we want to
solve $\min_{X=[\mU, \mV]} f(\mX)=f(\mU, \mV)=\frac 1 2 \n{\mU \mV^\top-\mA}^2_F$ for
$\mU\in \R^{m\times r}$ and $\mV\in \R^{n\times r}$.  It is a non-convex
problem, and the gradient $\nabla f$ is not globally
Lipschitz. With some tuning, one still can apply \algname{GD} and \algname{Nesterov's
accelerated method}, but---and we want to emphasize it---it was not a
trivial thing to find the steps in practice.  The steps we have chosen
were almost optimal, namely, the methods did not  converge if we
doubled the steps.  In contrast, our methods do
not require any tuning, so even in this regard they are much more
practical. For the experiments we used Movilens 100K
dataset~\cite{harper2016movielens} with more than million entries and
several values of $r=10,\ 20,\ 30$. All algorithms were initialized at
the same point, chosen randomly. The results are presented in
\Cref{fig:factorization}.

\paragraph{Cubic regularization.}

\begin{figure*}[t]
\centering
\begin{subfigure}[t]{0.32\textwidth}
    {\includegraphics[trim={3mm 0 3mm 0},clip,width=1\textwidth]{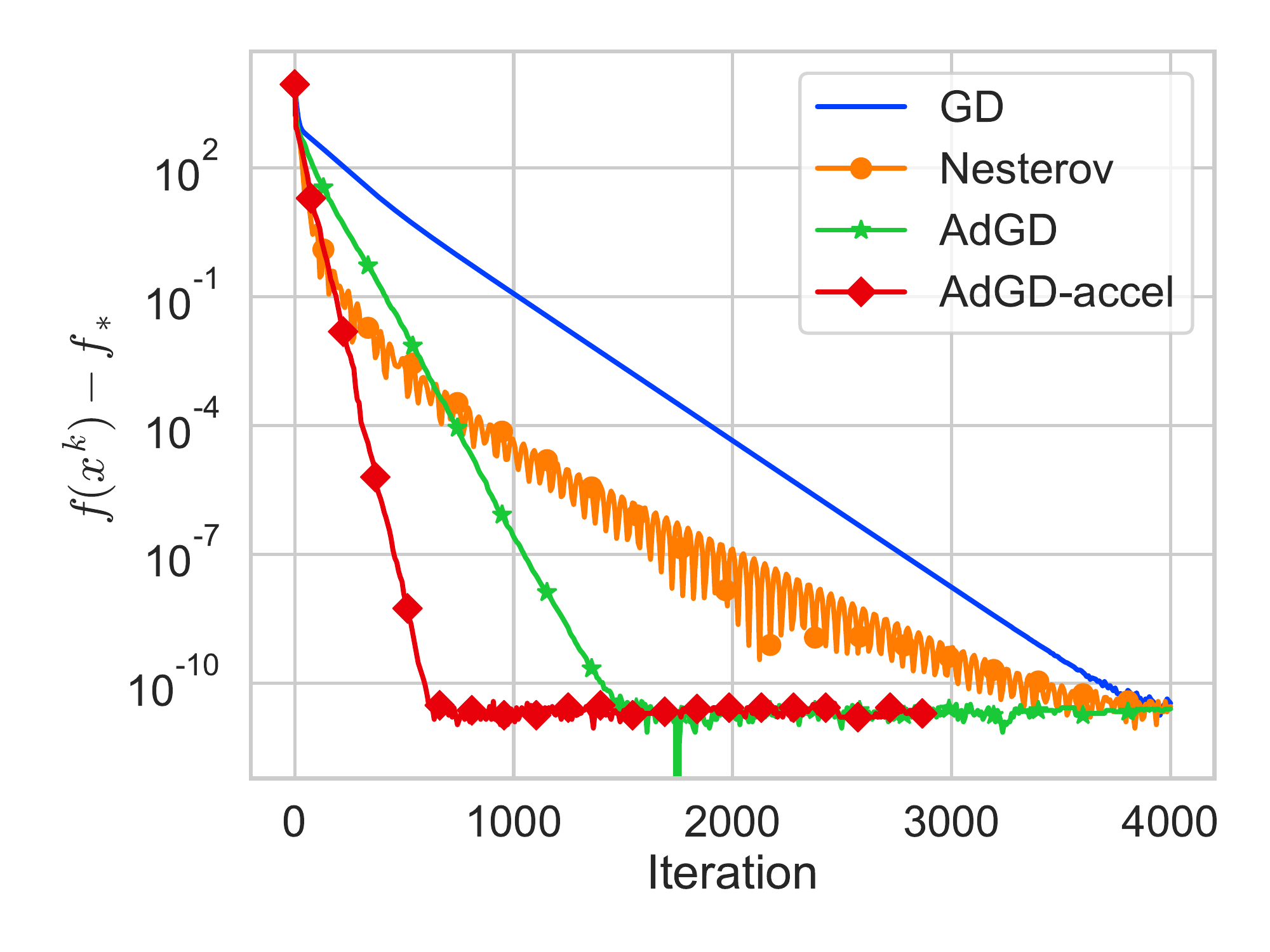}}
\caption{$M=10$}
\end{subfigure}
\hfill
\begin{subfigure}[t]{0.32\textwidth}
{\includegraphics[trim={3mm 0 3mm 0},clip,width=1\textwidth]{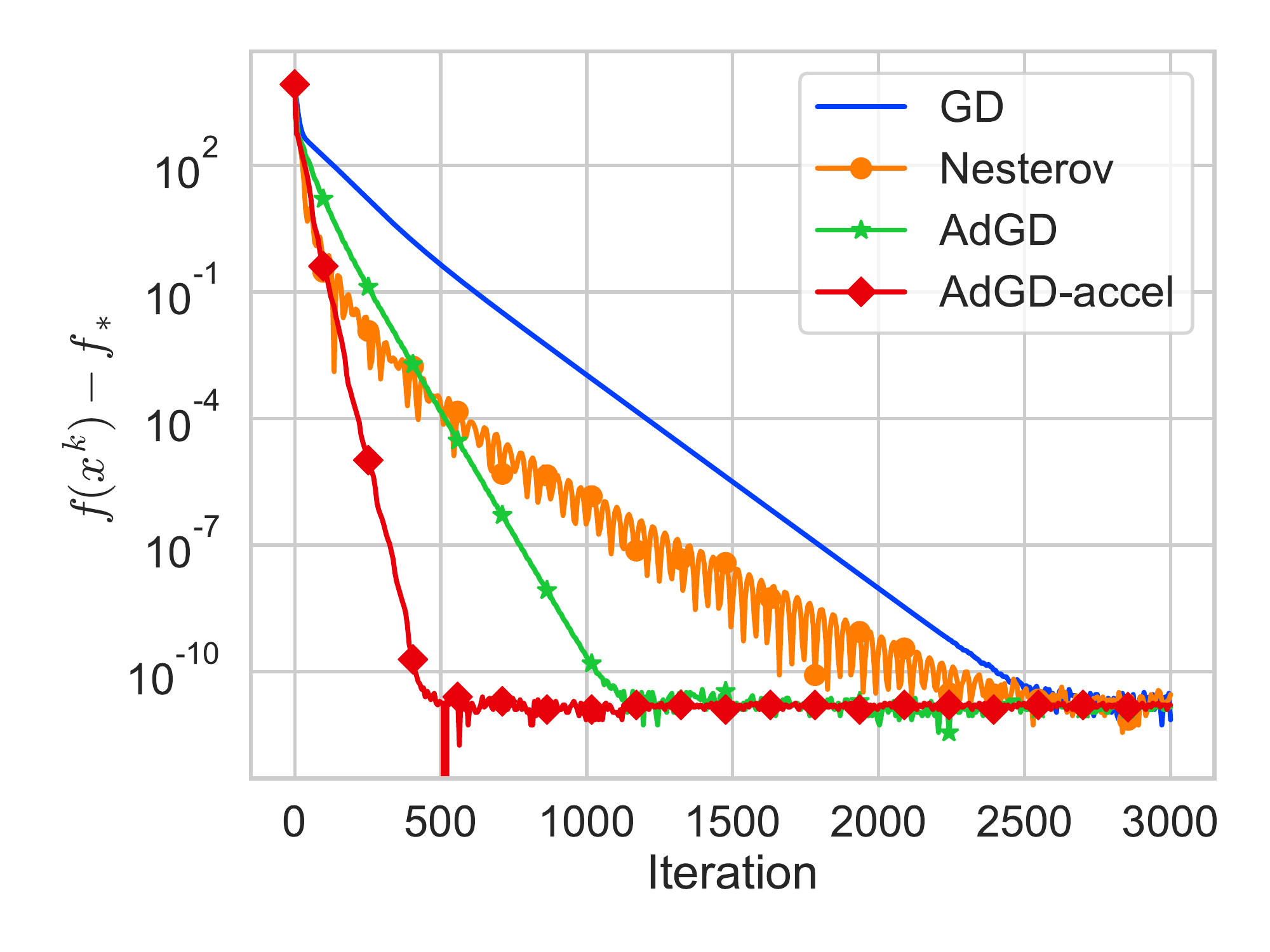}}
\caption{$M=20$}
\end{subfigure}
\hfill
\begin{subfigure}[t]{0.32\textwidth}
{\includegraphics[trim={3mm 0 3mm 0},clip,width=1\textwidth]{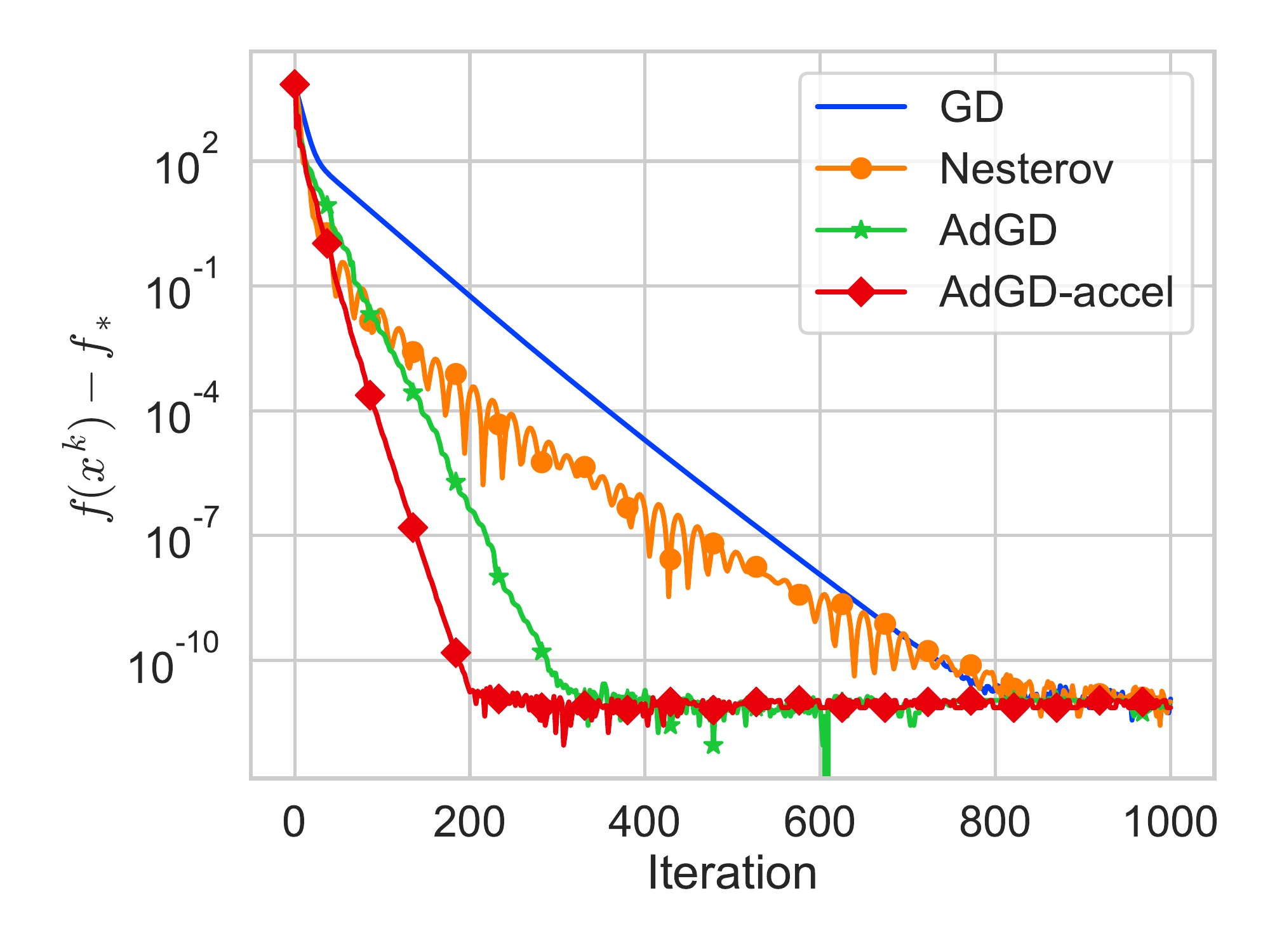}}
\caption{$M=100$}
\end{subfigure}
\caption{Results for the non-smooth subproblem from cubic regularization.}
\label{fig:cubic}
\end{figure*}

In cubic regularization of Newton method~\cite{nesterov2006cubic}, at
each iteration we need to  minimize $f(x)= g^\top x +
\frac{1}{2}x^\top H x + \frac{M}{6}\|x\|^3$, where $g\in \R^d, H\in
\R^{d\times d}$ and $M>0$ are given. This objective is smooth only
locally due to the cubic term, which is our motivation to consider
it. $g$ and $H$ were the gradient and the Hessian of the logistic loss
with the `covtype' dataset, evaluated at $x=0\in\R^d$. Although the values of
$M=10$, $20$, $100$ led to similar results, they also required different numbers of iterations, so we present
the corresponding results in \Cref{fig:cubic}.

\paragraph{\algname{Barzilai--Borwein}, \algname{Polyak} and line searches.}
\begin{figure*}[t]
\centering
\begin{subfigure}[t]{0.32\textwidth}
    {\includegraphics[trim={3mm 0 3mm 0},clip,width=1\textwidth]{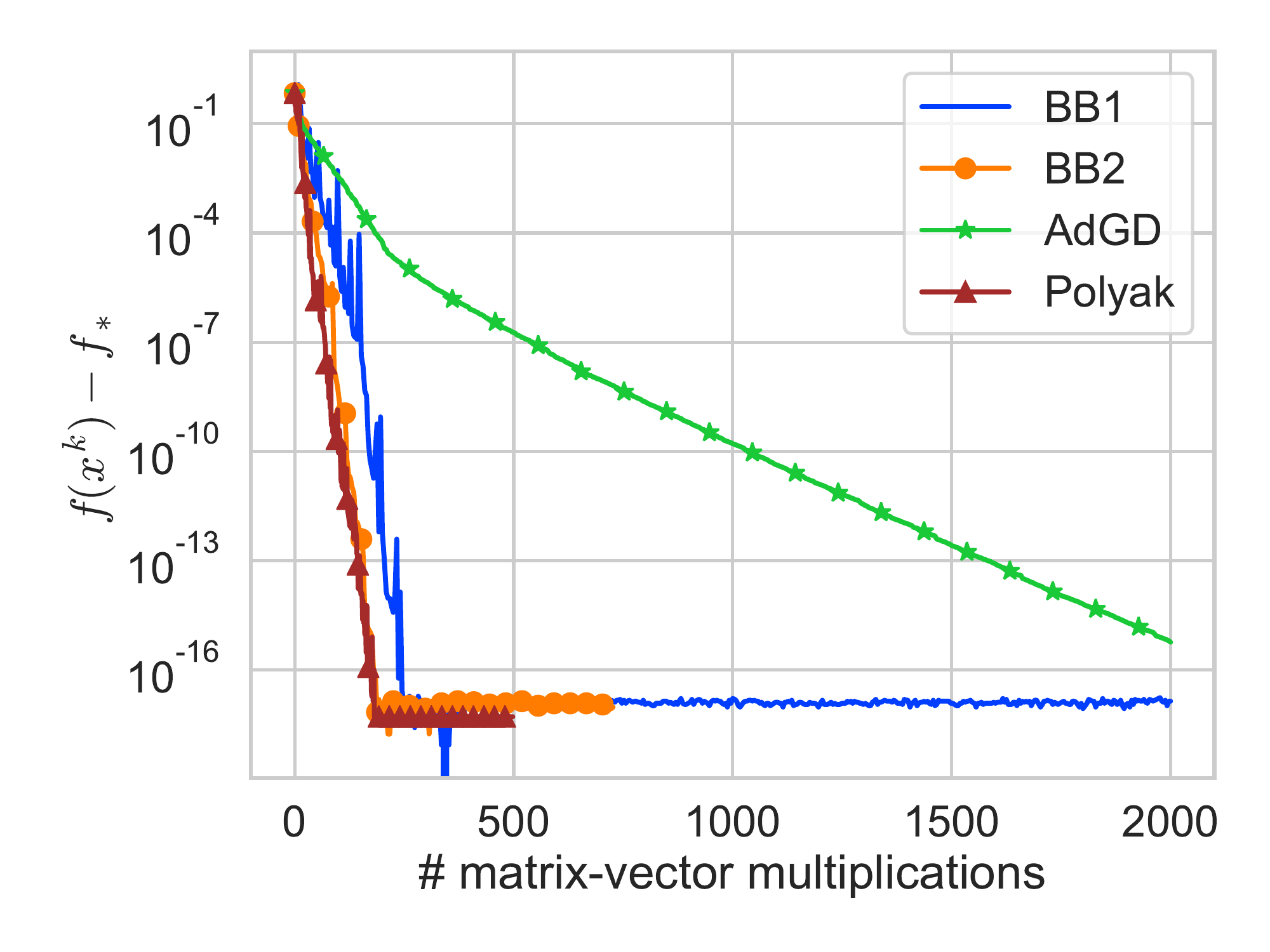}}
\caption{Mushrooms dataset, objective}
\end{subfigure}
\hfill
\begin{subfigure}[t]{0.32\textwidth}
{\includegraphics[trim={3mm 0 3mm 0},clip,width=1\textwidth]{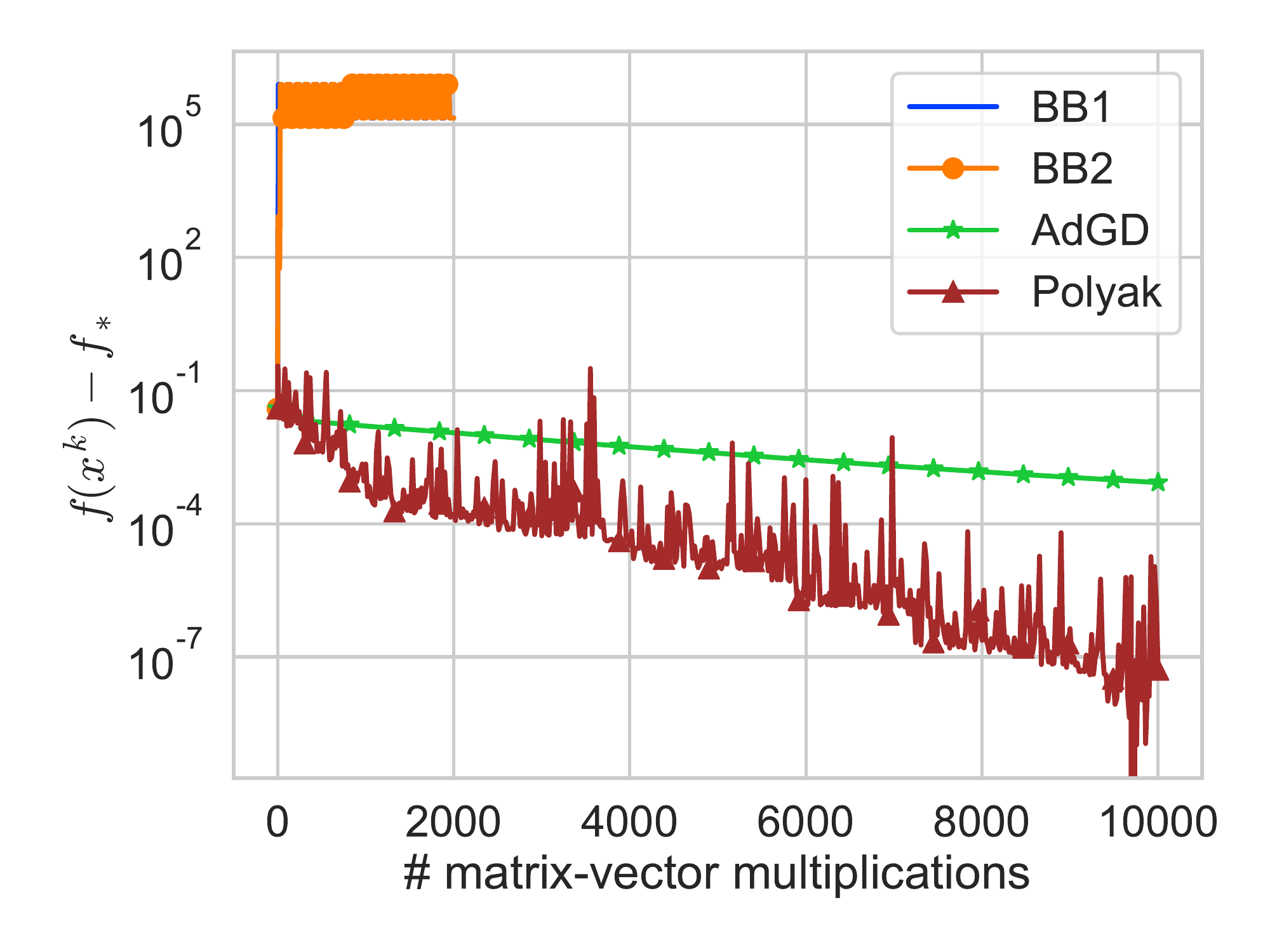}}
\caption{Covtype dataset, objective}
\end{subfigure}
\hfill
 \begin{subfigure}[t]{0.32\textwidth}
{\includegraphics[trim={3mm 0 3mm 0},clip,width=1\textwidth]{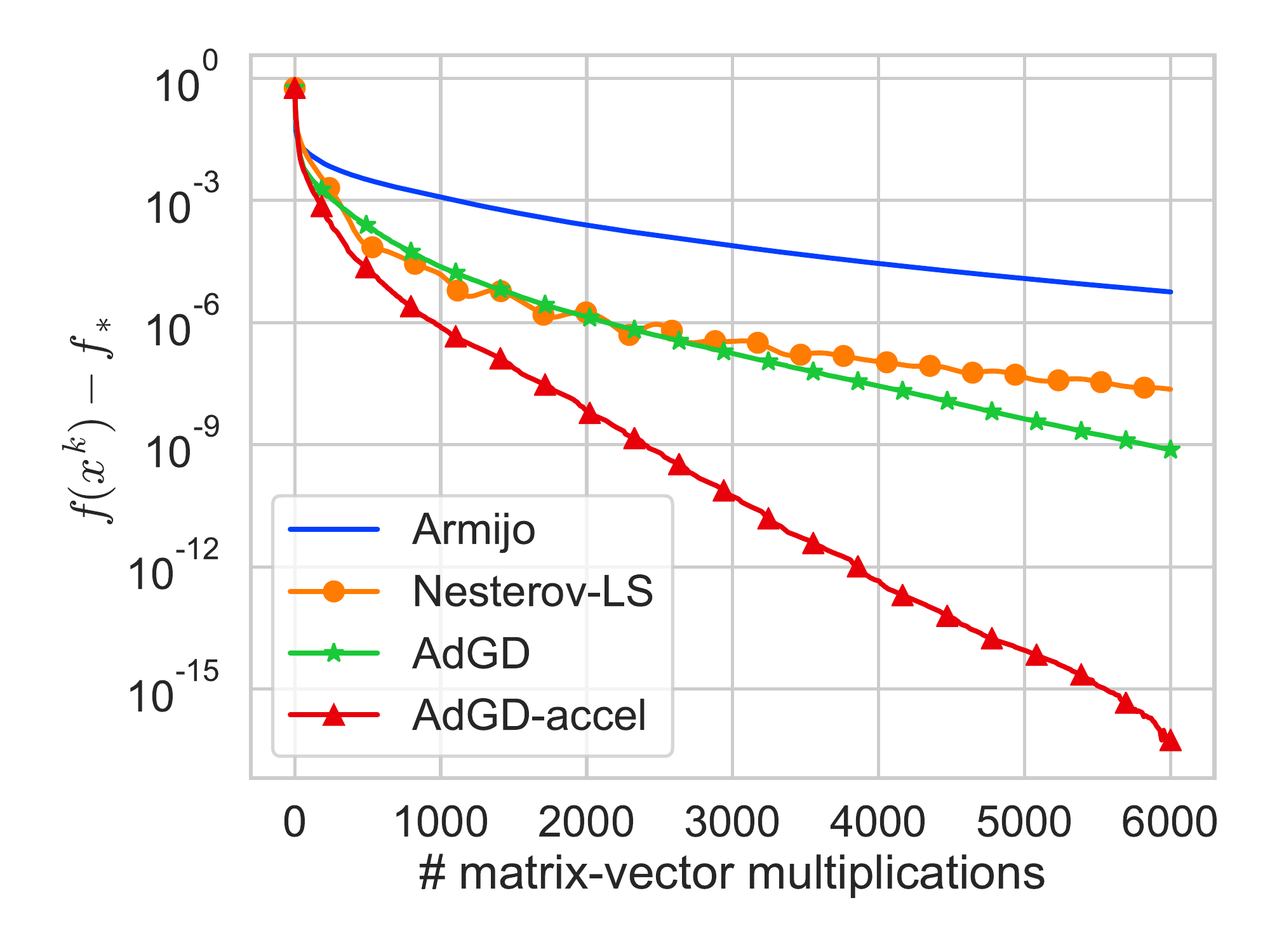}}
 \caption{W8a dataset, objective}
 \end{subfigure}
\caption{Additional results for the logistic regression problem.}
\label{fig:logistic_extra}
\end{figure*}

We have started this chapter with an overview of different approaches to
tackle the issue of a stepsize for \algname{GD}. Now, we
demonstrate some of those solutions. We again consider the
$\ell_2$-regularized logistic regression (same setting as before) with
`mushrooms', `covtype', and `w8a' datasets.

In~\Cref{fig:logistic_extra} (left) we see that the \algname{Barzilai--Borwein}
method can indeed be very fast.  However, as we said before, it lacks
a theoretical basis and \Cref{fig:logistic_extra} (middle) illustrates
this quite well. Just changing one dataset to another makes both
versions of this method to diverge on a strongly convex and smooth
problem. \algname{Polyak's method} consistently performs well (see
\Cref{fig:logistic_extra} (left and middle)), however, only after it
was fed with $f^\star$ that we found by running another method. Unfortunately,
for logistic regression there is no way to guess this value
beforehand.

Finally, line search for \algname{GD} (Armijo version) and \algname{Nesterov GD}
(implemented as in \cite{Nesterov2013a}) eliminates the need to know
the stepsize, but this comes with a higher price per iteration as
\Cref{fig:logistic_extra} (right) shows. Actually in all our
experiments for logistic regression with different datasets one
iteration of Armijo line search was approximately 2 times more
expensive than \algname{AdGD}, while line search for \algname{Nesterov GD} was 4 times
more expensive. We note that these observations are consistent with the
theoretical derivations in \cite{Nesterov2013a}.

\paragraph{Neural networks.}
\begin{figure*}[t]
\centering
\begin{subfigure}[t]{0.32\textwidth}
    {\includegraphics[trim={3mm 0 3mm 0},clip,width=1\textwidth]{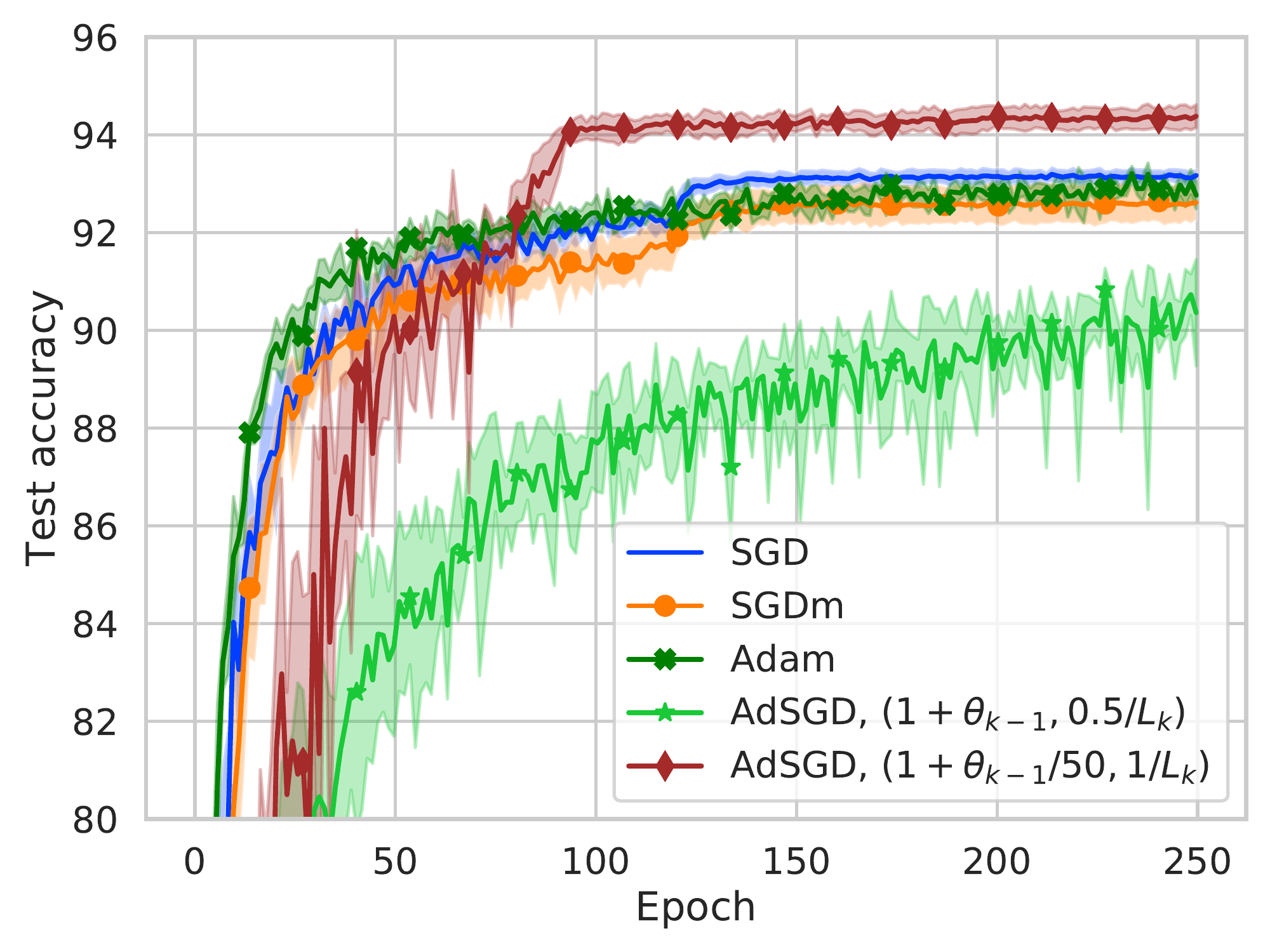}}
\caption{Test accuracy}
\end{subfigure}
\hfill
\begin{subfigure}[t]{0.32\textwidth}
{\includegraphics[trim={3mm 0 3mm 0},clip,width=1\textwidth]{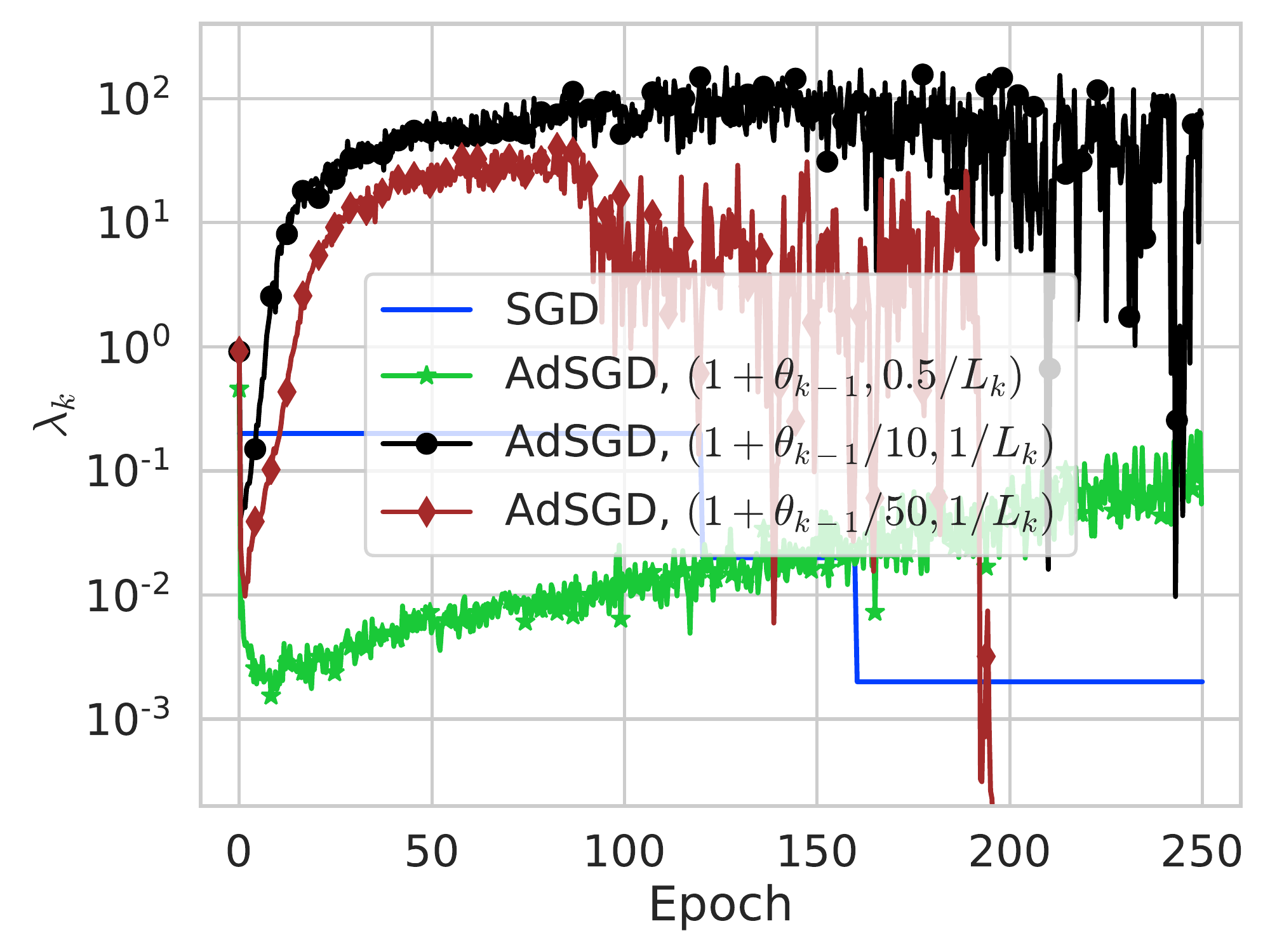}}
\caption{Stepsize}
\end{subfigure}
\hfill
\begin{subfigure}[t]{0.32\textwidth}
{\includegraphics[trim={3mm 0 3mm 0},clip,width=1\textwidth]{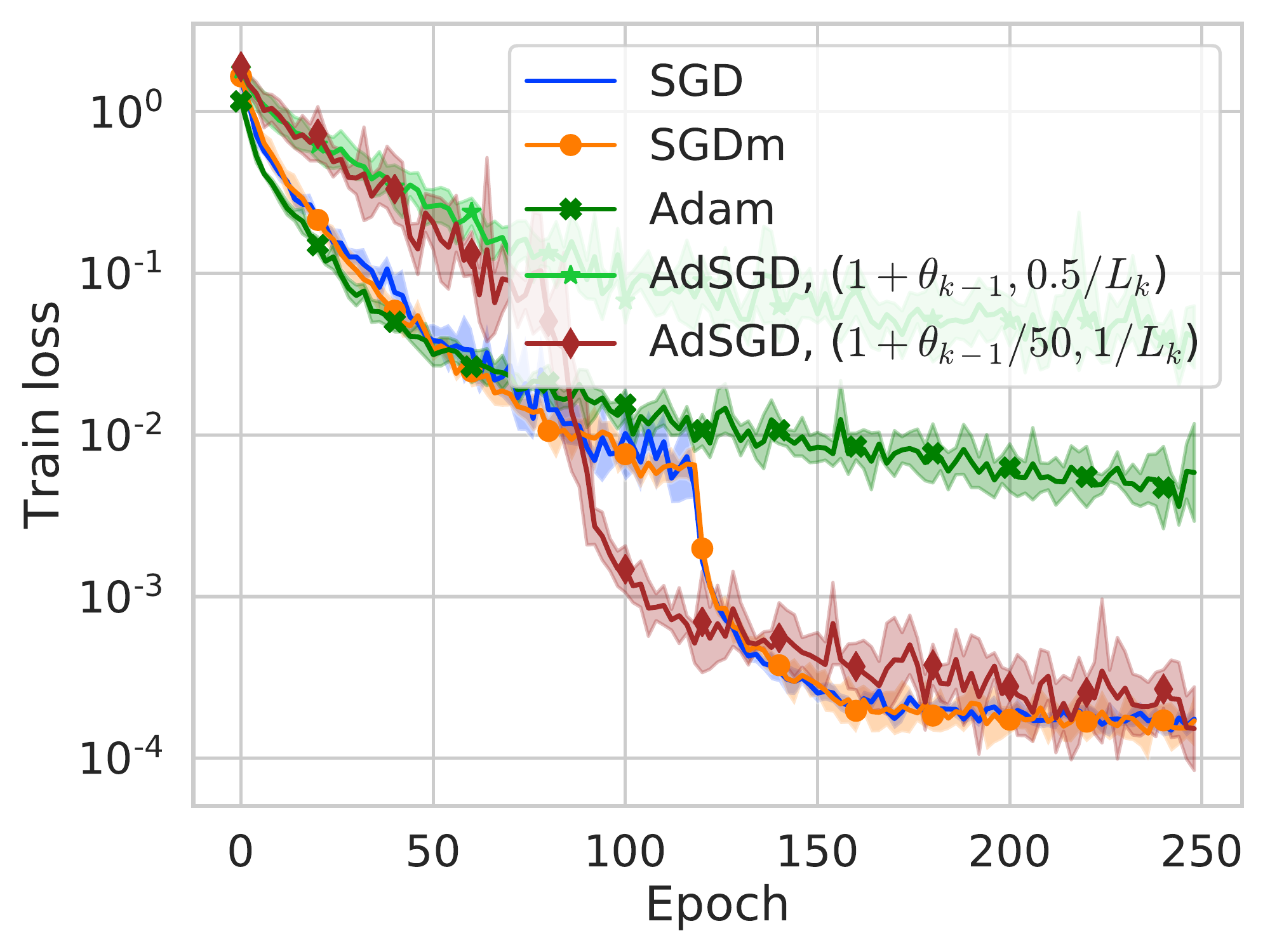}}
\caption{Train loss}
\end{subfigure}
\caption{Results for training ResNet-18 on CIFAR-10. Labels for \algname{AdGD} correspond to how $\gamma_k$ was estimated.}
\label{fig:resnet}
\end{figure*}

\begin{figure*}[t]
\centering
\begin{subfigure}[t]{0.32\textwidth}
    {\includegraphics[trim={3mm 0 3mm 0},clip,width=1\textwidth]{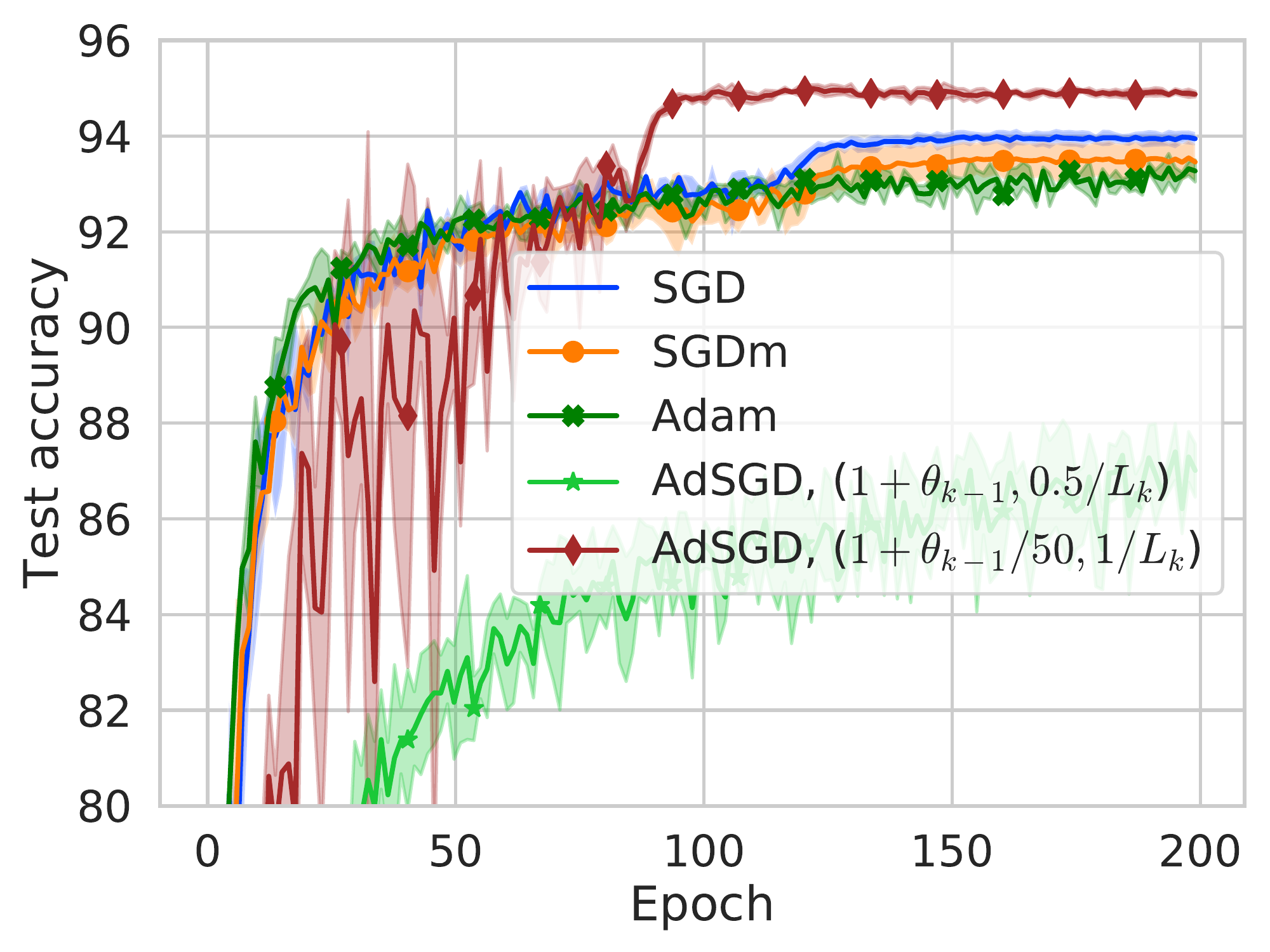}}
\caption{Test accuracy}
\end{subfigure}
\hfill
\begin{subfigure}[t]{0.32\textwidth}
{\includegraphics[trim={3mm 0 3mm 0},clip,width=1\textwidth]{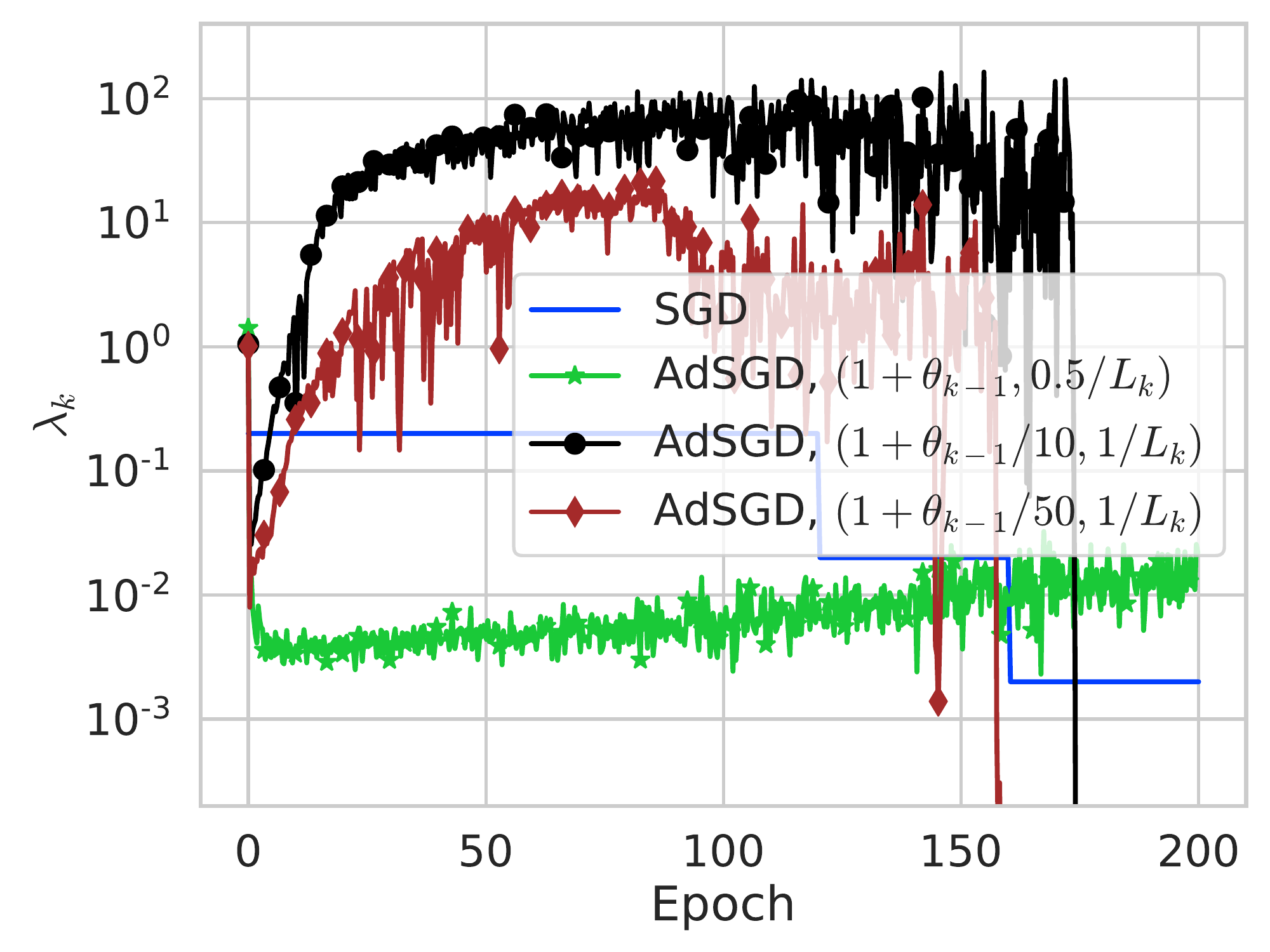}}
\caption{Stepsize}
\end{subfigure}
\hfill
\begin{subfigure}[t]{0.32\textwidth}
{\includegraphics[trim={3mm 0 3mm 0},clip,width=1\textwidth]{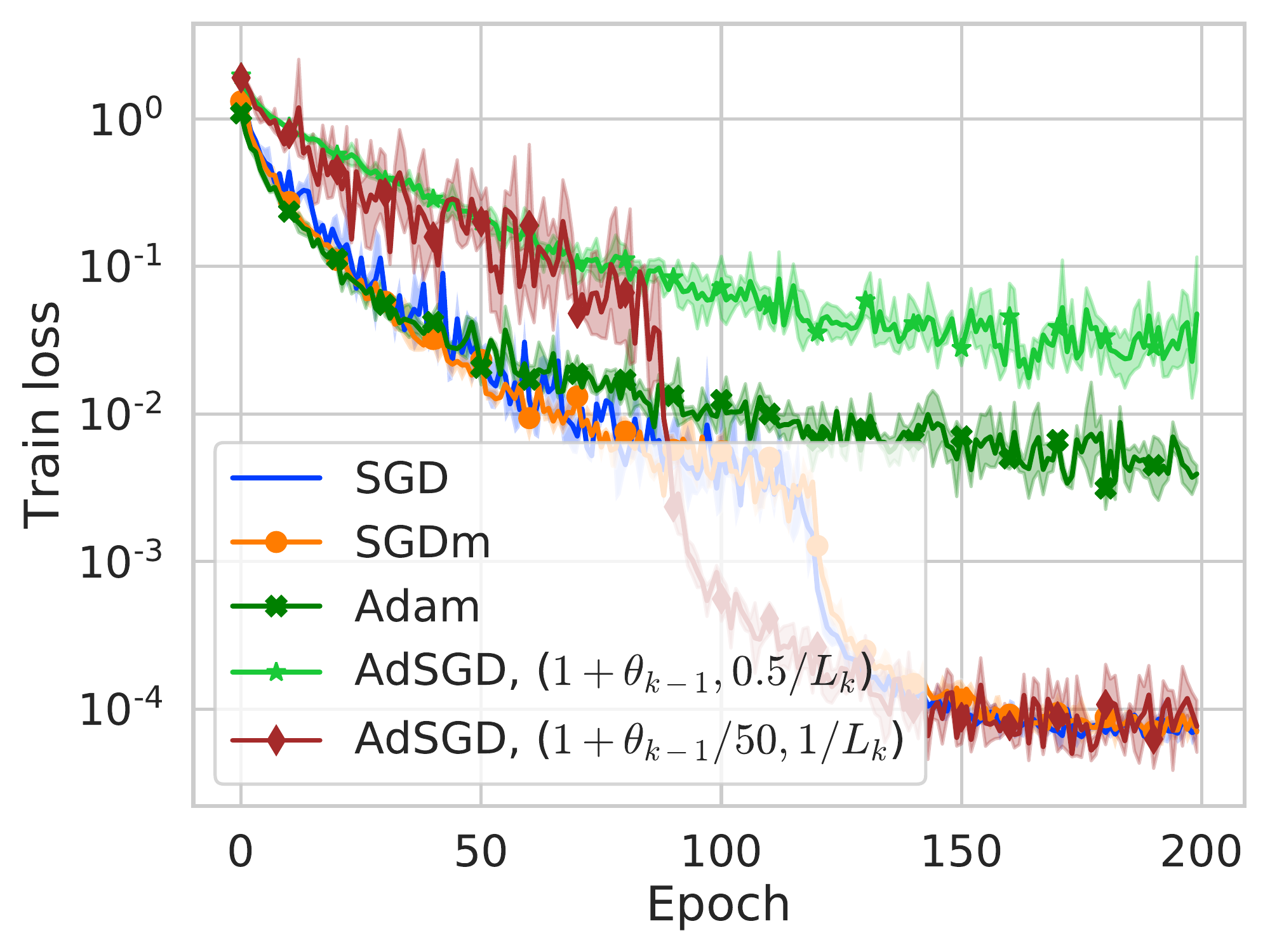}}
\caption{Train loss}
\end{subfigure}
\caption{Results for training DenseNet-121 on CIFAR-10.}
\label{fig:densenet}
\end{figure*}

We use standard ResNet-18 and DenseNet-121 architectures implemented in PyTorch~\cite{paszke2017automatic} and train them to classify images from the CIFAR-10 dataset~\cite{krizhevsky2009learning} with cross-entropy loss.

We use batch size 128 for all methods. For our method, we observed
that $\frac{1}{L_k}$ works better than $\frac{1}{2L_k}$. We ran it
with $\sqrt{1+\nu\th_k}$ in the other factor with values of
$\nu$ from $\{1, 0.1, 0.05, 0.02, \\ 0.01\}$ and $\nu=0.02$ performed the best. For reference, we provide the result for the theoretical estimate as well and value $\nu=0.1$ in the plot with estimated stepsizes.  The results are depicted in~\Cref{fig:resnet,fig:densenet} and other details are provided in~\cref{ap:exp_details}.

We can see that, surprisingly, our method achieves better test accuracy than \algname{SGD} despite having the same train loss. At the same time, our method is significantly slower at the early stage and the results are quite noisy for the first 75 epochs. Another observation is that the smoothness estimates are very non-uniform and $\gamma_k$ plummets once train loss becomes small.

\chapter{Achieving Fast Rates in Distributed Optimization with Quantization}\label{chapter:diana}
\graphicspath{{diana/}{diana/img/}}

\section{Introduction}

Big machine learning models are typically trained in a distributed fashion, with  the training data distributed across several workers, all of which compute in parallel an update to the model based on their local data. For instance, they can all perform a single step of Gradient Descent (\algname{GD} ) or Stochastic Gradient Descent (\algname{SGD} ). These updates are then sent to a parameter server which performs aggregation (typically this means just averaging of the updates) and then broadcasts  the aggregated updates back to the workers. 
The process is repeated until a good solution is found.

When doubling the amount of computational power, one usually expects to see the learning process finish in half time. If this is the case, the considered system is called to scale linearly. For various reasons, however, this does not happen, even to the extent that the system might become slower with more resources. At the same time, the surge of big data applications increased the demand for distributed optimization methods, often requiring new properties such as ability to find a sparse solution. It is, therefore, of great importance to design new methods that are versatile, efficient and scale linearly with the amount of available resources. In fact, the applications vary a lot in their desiderata. There is a rising interest in federated learning~\cite{FEDLEARN}, where the main concerns include the communication cost and ability to use  local data only in an attempt to provide a certain level of privacy. In high-dimensional machine learning problems, non-smooth $\ell_1$-penalty is often utilized, so one wants to have a support for proximable regularization. The efficiency of deep learning, in contrast, is dependent on heavy-ball momentum and non-convex convergence to criticality, while sampling from the full dataset might not be an issue. In this chapter, we try to address all of these questions.

{\bf Communication as the bottleneck.}
The  key aspects of distributed optimization efficiency are computational and communication complexity. In general, evaluating full gradients is intractable due to time and memory restrictions, so computation is made cheap by employing stochastic updates. On the other hand, in typical distributed computing architectures, communication is much slower (see Figure~\ref{fig:communication} for our experiments with communication cost of aggregating and broadcasting) than a stochastic update, and the design of a training algorithm needs to find a trade-off between them. 

\section{Related Work}

There have been considered several ways of dealing with the issue of slow communication. One of the early approaches is to have each worker perform a block descent step, which leads to the \algname{Hydra}  family of methods \cite{Hydra, Hydra2}. By choosing the size of the block, one directly chooses the amount of data that needs to be communicated. An alternative idea is for each worker to do more work between communication rounds (e.g., by employing a more powerful local solver, such as a second order method), so that computation roughly balances out with communication. The key methods in this sphere include \algname{CoCoA}  and its variants 
\cite{cocoa, COCOA+journal, cocoa-2018-JMLR}, \algname{DANE}  \cite{DANE}, \algname{DiSCO} \cite{DISCO, ma2015partitioning}, \algname{DANCE} \cite{jahani2018efficient} and \algname{AIDE} \cite{reddi2016aide}.

{\bf Update compression via randomized sparsification and/or quantization.} Practitioners suggested a number of heuristics to find a remedy for the communication botlleneck. Of special interest to this chapter is the idea of compressing \algname{SGD} updates, proposed in~\cite{seide20141}. Building off of this work, Alistarh et al.~\cite{alistarh2017qsgd} designed a variant of \algname{SGD} that guarantees convergence with compressed updates. Other works with \algname{SGD} update structure include~\cite{konecny2016randomized, bernstein2018signsgd, khirirat2018distributed}. Despite proving a convergence rate, Alistarh et al.~\cite{alistarh2017qsgd} also left many new questions open and introduced an additional, unexplained, heuristic of quantizing only vector blocks. Moreover, their analysis implicitly makes an assumption that all data should be available to each worker, which is hard and sometimes even impossible to satisfy. In a concurrent with~\cite{alistarh2017qsgd} work~\cite{wen2017terngrad}, the \algname{TernGrad} method was analyzed for stochastic updates that in expectation have positive correlation with the vector pointing to the solution. While giving more intuition about convergence of quantized methods, this work used $\ell_{\infty}$ norm for quantization, unlike $\ell_2$-quantization of~\cite{alistarh2017qsgd}. 

\section{Settings and Contributions}
{\bf The problem.} Let $f_i\colon \R^d\to \R$ be the loss of model $x$ obtained on data points belonging to distribution $\cD_i$, i.e., 
$
    f_i(x) \eqdef \mathbb{E}_{\xi\sim \cD_i} [f(x; \xi)]
$ and $\psi\colon\R^d\to \R\cup \{+\infty\}$ be a proper closed convex regularizer. In this chapter, we focus on the problem of training a machine learning model via regularized empirical risk minimization:
\begin{equation} \label{eq:main}
     \min_{x\in \R^d} \left[f(x) + \psi(x) \eqdef \frac{1}{M}\sum \limits_{i=1}^M f_i(x) + \psi(x)\right].
\end{equation}
We do not assume any kind of similarity between distributions $\cD_1,\dotsc, \cD_M$. 

We also need to introduce several key concepts and ingredients that come together to make the algorithm.  In each iteration $k$ of \algname{DIANA}, each node will sample an unbiased estimator of the local gradient.  We assume that these gradients have bounded variance.

\begin{assumption}[Stochastic gradients]\label{as:noise}
For every $i = 1,2,\dots,M$,
	$\mathbb{E} [g_i^k \;|\; x^k] =  \nabla f_i(x^k)$. Moreover, the variance is bounded:
	\begin{align}\label{eq:bounded_noise}
		\mathbb{E} \left[\|g_i^k - \nabla f_i(x^k)\|^2\right]
		\le \sigma_i^2.
	\end{align}

\end{assumption}

Note that $g^k \eqdef \frac{1}{M}\sumiM g_i^k$ is an unbiased estimator of $\nabla f(x^k)$: \begin{equation}\label{eq:hat_g_expectation}  \mathbb{E}[g^k \;|\; x^k] = \frac{1}{M} \sum \limits_{i=1}^M \nabla f_i(x^k) = \nabla f(x^k).\end{equation}
Let $\sigma^2 \eqdef \frac{1}{M}\sumiM \sigma_i^2$. By independence of the random vectors $(g_i^k - \nabla f_i (x^k))_{i=1}^M$, its variance is bounded above by
\begin{equation}\label{eq:bgud7t9gf}
	\mathbb{E}\left[ \|g^k -\nabla f(x^k)\|^2 \; | \; x^k \right] \leq \frac{\sigma^2}{M}.
\end{equation}

{\bf Notation.}  By $\sign (t)$ we denote the sign of $t\in \R$ (-1 if $t<0$, 0 if $t=0$ and $1$ if $t>0$). The $j$-th element of a vector $x\in \R^d$ is denoted as $x_{(j)}$. For $x = (x_{(1)},\dots,x_{(d)})\in \R^d$ and $p\geq 1$, we let $\|x\|_p = (\sum_i |x_{(i)}|^p)^{1/p}$. Note that $\|x\|_1 \geq \|x\|_p \geq \|x\|_\infty$ for all $x$. By $\|x\|_0$ we denote the number of nonzero elements of $x$. Detailed description of the notation is in Table~\ref{tbl:notation-table} in the appendix.  

\begin{table}
\begin{center}
\caption{Comparison of \algname{DIANA} and related methods. Here ``lin. rate'' means that  linear convergence either to a ball around the optimum or to the optimum was proved, ``loc. data'' describes whether or not authors assume that $f_i$ is available at node $i$ only, ``non-smooth'' means support for a non-smooth regularizer, ``momentum'' says whether or not authors consider momentum in their algorithm,  and ``block quant.'' means theoretical justification for using block quantization.
\label{tbl:comparison}
}
\footnotesize
\begin{tabular}{cccccc}
\toprule
Method & Linear rate & Local data & Non-smooth & Momentum & Block quant.\\
\midrule
\algname{DIANA} (New!) & \cmark & \cmark & \cmark & \cmark & \cmark\\
\hline
\algname{QSGD}  \cite{alistarh2017qsgd}  & \xmark & \xmark & \xmark & \xmark & \xmark\\
\algname{TernGrad} \cite{wen2017terngrad} & \xmark & \xmark & \xmark & \xmark & \xmark\\
\algname{DQGD}  \cite{khirirat2018distributed}  & \cmark & \cmark & \xmark & \xmark & \xmark\\
\algname{ QSVRG} \cite{alistarh2017qsgd}  & \cmark & \cmark & \xmark & \xmark & \xmark\\
\bottomrule
\end{tabular}
\end{center}
\vspace{-20pt}
\end{table}

 \textbf{\algname{DIANA}.} We develop a  distributed gradient-type method with compression of gradient differences, which we call \algname{DIANA} (Algorithm~\ref{alg:distributed1}).
 
  \textbf{Rate in the strongly convex case.} We show that when applied to the smooth strongly convex minimization problem with arbitrary closed convex regularizer, \algname{DIANA} has the iteration complexity 
  \[
  	\cO\left(\left(\sqrt{\frac{d}{m}} + \kappa\left(1 + \frac{1}{M}\sqrt{\frac{d}{m}}\right)\right)\log\frac{1}{\varepsilon}\right),
  \]  
  to a ball with center at the optimum (see Sec~\ref{sec:theory-strong-convex}, Theorem~\ref{thm:DIANA-strongly_convex} and Corollary~\ref{cor:DIANA-strong-convex} for the details). In the case of decreasing stepsize we show $\cO\left(\frac{1}{\varepsilon}\right)$ iteration complexity (see Sec~\ref{sec:DIANA-decreasing-stepsizes}, Theorem~\ref{th:str_cvx_decr_steps} and Corollary~\ref{cor:str_cvx_decr_step} for the details). Unlike  in \cite{khirirat2018distributed}, in a noiseless regime our method converges to the exact optimum, and at a  linear rate.
  
\textbf{Rate in the non-convex case.}  We prove that \algname{DIANA} also works for smooth non-convex problems and get the iteration complexity 
\[
	\cO\left(\frac{1}{\varepsilon^2}\left(\frac{L^2(f(x^0) - f^*)^2}{M^2\alpha_p^2} +  \frac{\sigma^4}{(1+M\alpha_p)^2}\right)\right)
\]
(see Sec~\ref{sec:DIANA-non-convex}, Theorem~\ref{thm:DIANA-non-convex} and Corollary~\ref{cor:DIANA-non-convex} for the details).

\textbf{\algname{DIANA} with momentum.} We study momentum version of \algname{DIANA} for the case of smooth non-convex objective with constant regularizer and $f_i = f$ (see Sec~\ref{sec:DIANA-momentum}, Theorem~\ref{thm:DIANA-momentum} and Corollary~\ref{cor:DIANA-momentum} for the details). We summarize a few key features of our complexity results established  in Table~\ref{tbl:comparison}.

\textbf{First rate for  \algname{TernGrad}.} We provide first convergence rate of \algname{TernGrad} and provide new tight analysis of 1-bit \algname{QSGD} under less restrictive assumptions for both smooth strongly convex objectives with arbitrary closed convex regularizer and non-convex objective (see Sec~\ref{sec:Alg} for the detailed comparison). Both of these methods are just special cases of our Algorithm~\ref{alg:terngrad} which is also a special case of Algorithm~\ref{alg:distributed1} with $\alpha = 0$ and $h_i^0 = 0$ for all $i$. We show that Algorithm~\ref{alg:terngrad} has $\cO\left(\frac{\kappa}{M\alpha_p}\right)$ iteration complexity of convergence to the ball with center at the optimum in the case of the smooth strongly convex minimization problem with arbitrary closed convex regularizer (see Sec~\ref{sec:Terngrad}, Theorem~\ref{thm:terngrad_strg_cvx_prox}) and 
	\[
		\cO\left(\frac{1}{\varepsilon^2}\left(\frac{L^2(f(x^0) - f(x^\star))^2}{M^2\alpha_p^2} + \frac{\sigma^4}{(1+M\alpha_p)^2}\right)\right)
	\]
	in the case of non-convex minimization problem (see  Theorem~\ref{thm:TernGrad-nonconvex} and Corollary~\ref{cor:TernGrad-nonconvex}).

\textbf{\algname{QSGD} and \algname{TernGrad} with momentum.}	  We study momentum version of \algname{DIANA} for $\alpha = 0, h_i^0 =0$ and, in particular, we propose momentum versions of (1-bit) \algname{QSGD} and \algname{TernGrad} the case of smooth non-convex objective with constant regularizer and $f_i = f$ (see Sec~\ref{sec:TernGrad-momentum}, Theorem~\ref{thm:TernGrad-momentum} and Corollary~\ref{cor:TernGrad-momentum}).

\textbf{Optimal norm power.} We find the answer for the following question: \textit{which $\ell_p$ norm to use for quantization in order to get the best iteration complexity of the algorithm?} It is easy to see that all the bounds that we propose depend on $\frac{1}{\alpha_p}$ where $\alpha_p$ is an increasing function of $1\le p \le \infty$ (see Lemma~\ref{lema:alpha_p} for the details). That is, \textit{for both Algorithm~\ref{alg:distributed1}~and~\ref{alg:terngrad} the iteration complexity reduces when $p$ is growing and the best iteration complexity for our algorithms is achieved for $p = \infty$.} This implies that \algname{TernGrad} has better iteration complexity than 1-bit \algname{QSGD}.

\textbf{First analysis for block-quantization.}	  We give a first analysis of block-quantization (i.e., bucket-quantization), which was mentioned in \cite{alistarh2017qsgd} as a useful heuristic.

\section{The Algorithm} \label{sec:Alg}

\begin{algorithm}[t]
   \caption{\rrbox{\algname{DIANA}} ($M$ nodes)}
   \label{alg:distributed1}
\begin{algorithmic}[1]
   \State \textbf{Input:}  learning rates $\alpha>0$ and $(\gamma_k)_{k\geq 0}$, initial vectors $x^0, h_1^0,\dotsc, h_M^0 \in \R^d$ and $h^0 = \frac{1}{M}\sum_{i=1}^M h_i^0$, quantization parameter $p \geq 1$, sizes of blocks $(d_l)_{l=1}^M$, momentum parameter $0\le \beta < 1$, number of steps $K$
   \State $v^0 = \nabla f(x^0)$
   \For{$k=0,1,\dotsc, K-1$}
	   \State Broadcast $x^{k}$ to all workers
        \For{$i=1,\dotsc,n$ in parallel}
	        \State Sample $g^{k}_i$ such that $\mathbb{E} [g^k_i \;|\; x^k]  =\nabla f_i(x^k)$ and let $\Delta^k_i = g^k_i - h^k_i$
\State Sample $\hat \Delta^k_i \sim {\rm Quant}_p(\Delta^k_i,(d_l)_{l=1}^m)$ and let $h_i^{k+1} = h_i^k + \alpha \hat \Delta_i^k$ and  $\hat g_i^k = h_i^k + \hat \Delta_i^k$
        \EndFor
       \State $\hat \Delta^k = \frac{1}{M}\sum_{i=1}^M \hat \Delta_i^k$
       \State $\hat g^k = \frac{1}{M}\sum_{i=1}^M \hat g_i^k = h^k + \hat \Delta^k $
       \State $v^k = \beta v^{k-1} + \hat g^k$
        \State $x^{k+1} = \prox_{\gamma_k \psi}\left(x^k - \gamma_kv^k \right)$
        \State  $h^{k+1}  = \frac{1}{M}\sum_{i=1}^M h_i^{k+1} = h^k + \alpha \hat \Delta^k$
   \EndFor
\end{algorithmic}
\end{algorithm}

In this section we describe our main method---\algname{DIANA}. 

{\bf Quantization.} \algname{DIANA} applies random compression (quantization) to gradient differences, which are then communicated  to a parameter server.  We now define the random quantization transformations used. Our first quantization operator transforms a vector $\Delta \in \R^d$ into a random vector $\hat{\Delta} \in \R^d$ whose entries belong to the set $\{-t,0,t\}$ for some $t>0$.

\begin{definition}[$p$-quantization]\label{def:p-quant}  Let $\Delta \in \R^d$ and let $p\geq 1$. If $\Delta=0$, we define $\widetilde{\Delta}=\Delta$. If $\Delta \neq 0$, we define $\widetilde{\Delta}$ by setting 
\begin{equation}\label{eq:quant-j}\widetilde{\Delta}_{(j)} = \|\Delta\|_p \sign(\Delta_{(j)}) \xi_{(j)}, \quad j=1,2,\dots,d,\end{equation}  where
 $\xi_{(j)}\sim {\rm Be}\left(|\Delta_{(j)}|/\|\Delta\|_p\right)$  are Bernoulli random variables\footnote{That is, $\xi_{(j)} = 1$ with probability $|\Delta_{(j)}|/\|\Delta\|_p$ (observe that this quantity is always upper bounded by 1) and $\xi_{(j)} = 0$ with probability $1-|\Delta_{(j)}|/\|\Delta\|_p$. }. Note that \begin{equation}\label{eq:quant}\widetilde{\Delta} = \|\Delta\|_p \; \sign(\Delta) \circ \xi,\end{equation} where $\sign$ is applied elementwise, and $\circ$ denotes the Hadamard (i.e., elementwise) product.  We say that $\widetilde{\Delta}$ is $p$-quantization of $\Delta$. When sampling $\widetilde{\Delta}$, we shall write $\widetilde{\Delta} \sim {\rm Quant}_{p}(\Delta)$.  
\end{definition}

In addition, we consider a block variant of $p$-quantization operators. These are defined and their properties studying in Section~ in the appendix.

{\bf Communication cost.} If $b$ bits are used to encode a float number, then at most $C(\hat \Delta) \eqdef \|\hat \Delta\|_0^{1/2}(\log \|\hat \Delta\|_0 + \log 2 + 1) + b$ bits are needed to communicate $\hat \Delta$ with Elias coding (see Theorem~3.3 in \cite{alistarh2017qsgd}). In our next result, we given an upper bound on the expected communication cost.

\begin{theorem}[Expected sparsity]\label{th:quantization_quality1} 
	Let $0\neq \Delta\in \R^{\tilde d}$ and ${\widetilde \Delta} \sim {\rm Quant}_{p}(\Delta)$ be its $p$-quantization. Then  
	\begin{equation} \label{eq:expected_comm_cost}
	 \mathbb{E} \|\widetilde \Delta\|_0 =  \frac{\|\Delta\|_1}{\|\Delta\|_{p}} \leq \|\Delta\|_0^{1-1/p} \leq \tilde{d}^{1-1/p},
	\end{equation}
 \begin{equation} 
 \label{eq:expected_comm_cost2} 
	C_p \eqdef 
	\mathbb{E}\left[ C(\widetilde \Delta)\right]
	\leq \frac{\|\Delta\|_1^{1/2}}{\|\Delta\|_{p}^{1/2}} (\log \tilde{d} + \log 2 + 1) + b.
 \end{equation}
 All expressions in \eqref{eq:expected_comm_cost} and  \eqref{eq:expected_comm_cost2} are  increasing functions of $p$. 

\end{theorem}

{\bf \algname{DIANA}.}    In \algname{DIANA}, each machine $i\in \{1,2,\dots,n\}$ first computes a stochastic gradient $g_i^k$ at current iterate $x^k$. We do not quantize this information and send it off to the parameter server as that approach would not converge for $\psi\neq 0$. Instead, we maintain memory $h^k_i$ at each node $i$ (initialized to arbitrary values), and {quantize the difference $\delta_i^k\eqdef g_i^k-h_i^k$ instead. Both the node and the parameter server update $h_i^k$ in an appropriate manner, and a proximal gradient descent step is taken with respect to direction $v^k = \beta v^{k-1} + \hat{g}^k$, where $0\leq \beta\leq 1$ is a momentum parameter, whereas $\hat{g}^k$ is an unbiased estimator of the full gradient, assembled from the memory $h_i^k$ and the transmitted quantized vectors. Note that we allows for block quantization for more flexibility. In practice, we want the transmitted quantized vectors to be much easier to communicate than the full dimensional vector in $\R^d$, which can be tuned by the choice of $p$ defining the quantization norm, and the choice of blocks.
    
{\bf Relation to \algname{QSGD} and \algname{TernGrad}.} If the initialization is done with $h^0=0$ and $\alpha=0$, our method reduces to either 1-bit \algname{QSGD} or \algname{TernGrad} with $p=2$ and $p=\infty$ respectively. We unify them in the Algorithm~\ref{alg:terngrad}. We analyze this algorithm (i.e., \algname{DIANA} with $\alpha = 0$ and $h_i^0 = 0$) in three cases: i) smooth strongly convex objective with arbitrary closed convex regularizer; ii) smooth non-convex objective with constant regularizer; iii) smooth non-convex objective with constant regularizer for the momentum version of the algorithm. We notice, that in the original paper \cite{wen2017terngrad}, the authors do not provide the rate of convergence for \algname{TernGrad}, and we get the convergence rate for the three aforementioned situations as a special case of our results. Moreover, we emphasize that our analysis is new even for 1-bit \algname{QSGD}, since in the original paper \cite{alistarh2017qsgd}, the authors consider only the case of bounded gradients ($\mathbb{E}\left[\|g^k\|^2\right] \le B^2$), which is a very restrictive assumption, and they do not provide rigorous analysis of block-quantization. In contrast, we consider more general case of block-quantization and assume only that the variance of the stochastic gradients is bounded, which is strictly less restrictive since the inequality $\mathbb{E}\left[\|g^k\|^2\right] \le B^2$ implies  $\mathbb{E}\left[\|g^k - \nabla f(x^k)\|^2\right] \le \mathbb{E}\left[\|g^k\|^2\right] \le B^2$.

We obtain the convergence rate for arbitrary $p\ge 1$ for the three aforementioned cases (see Theorems~\ref{thm:TernGrad-nonconvex},~\ref{thm:TernGrad-momentum},~\ref{thm:terngrad_strg_cvx_prox},~\ref{thm:TernGrad-decreasing-stepsizes} and Corollaries~\ref{cor:TernGrad-nonconvex},~\ref{cor:TernGrad-momentum},~\ref{cor:TernGrad-decreasing-stepsizes} for the details) and all obtained bounds becomes better when $p$ is growing, which means that \algname{TernGrad} has better iteration complexity than \algname{QSGD} and, more generally, the best iteration complexity attains for $\ell_\infty$ norm quantization. 

\section{Convergence Theory for Strongly Convex Losses} \label{sec:theory-strong-convex}
\begin{table}[t]
    \centering
    \caption{Summary of iteration complexity results. }
    \footnotesize
    \begin{tabular}{cccccccc}
    \toprule
    Block quant. & Local data & Non-convex & Str.\ convex & $\psi$ & Momentum & $\alpha>0$ & Theorem\\
	\midrule 
	\cmark & \cmark & \cmark & \xmark & \xmark & \xmark & \cmark &  \ref{thm:DIANA-non-convex} \\
	\cmark & \cmark & \cmark & \xmark & \xmark & \cmark & \cmark &  \ref{thm:DIANA-momentum}\\
	\cmark & \cmark & \xmark & \cmark & \cmark & \xmark & \cmark &  \ref{thm:DIANA-strongly_convex}, \ref{th:str_cvx_decr_steps} \\
	\cmark & \cmark & \cmark & \xmark & \xmark & \xmark & \xmark &  \ref{thm:TernGrad-nonconvex} \\
	\cmark & \cmark & \cmark & \xmark & \xmark & \cmark & \xmark &  \ref{thm:TernGrad-momentum}\\
	\cmark & \cmark & \xmark & \cmark & \xmark & \xmark & \xmark &  \ref{thm:terngrad_strg_cvx_prox}, \ref{thm:TernGrad-decreasing-stepsizes} \\
	\bottomrule 
    \end{tabular}
    \label{tab:results_diana}
    \vspace{-20pt}
\end{table}

Let us introduce two key assumptions of this section.

\begin{assumption}[$L$--smoothness]
    We say that a function $f$ is $L$-smooth if
    \begin{align}
        f(x) \le f(y) + \< \nabla f(y), x - y> + \frac{L}{2}\|x - y\|^2, \quad \forall x, y. \label{eq:smoothness_functional}
    \end{align}
\end{assumption}

\begin{assumption}[$\mu$-strong convexity]
	 $f$ is $\mu$-strongly convex, i.e., 
	\begin{align}
		f(x) \ge f(y) + \< \nabla f(y), x - y> + \frac{\mu}{2}\|x - y\|^2, \quad \forall x,y. \label{eq:strong_cvx_functional}
	\end{align}
\end{assumption}

For $1\leq p \leq +\infty$, define
\begin{equation}\label{eq:alpha_p} \alpha_p(d) \eqdef \inf_{x\neq 0,x\in\R^d} \frac{\|x\|^2}{\|x\|_1 \|x\|_p}. \end{equation}

\begin{lemma} \label{lema:alpha_p} $\alpha_p$ is  increasing as a function of $p$ and decreasing as a function of $d$. In particular,  $\alpha_1 \leq \alpha_2 \leq \alpha_\infty$, and moreover,
$\alpha_1(d) = \frac{1}{d}$, $\alpha_2(d) = \frac{1}{\sqrt{d}}$, $\alpha_{\infty}(d) = \frac{2}{(1+\sqrt{d})}$
and, as a consequence, for all positive integers $\widetilde{d}$ and $d$ the following relations holds 
$\alpha_1(\widetilde{d}) = \frac{\alpha_1(d) d}{\widetilde{d}}$, $\alpha_2(\widetilde{d}) = \alpha_2(d)\sqrt{\frac{d}{\widetilde{d}}},$ and
$\alpha_\infty(\widetilde{d}) = \frac{\alpha_\infty(d)(1+\sqrt{d})}{(1 + \sqrt{\widetilde{d}})}.$
\end{lemma}

\begin{theorem} \label{thm:DIANA-strongly_convex}  Assume the functions $f_1,\dots,f_M$ are $L$--smooth and $\mu$--strongly convex. Choose stepsizes $\alpha>0$ and $\gamma_k=\gamma>0$, block sizes $(d_l)_{l=1}^m$, where $\widetilde{d} = \max\limits_{l=1,\ldots,m}d_l$, and parameter $c>0$ satisfying the following relations:
\begin{equation}\label{eq:cond1} \frac{1 + M c \alpha^2}{1 + M c \alpha}   \leq \alpha_p \eqdef \alpha_p(\widetilde{d}),\end{equation}
\begin{equation}\label{eq:cond2}\gamma \leq \min\left\{\frac{\alpha}{\mu}, \frac{2}{(\mu+L)(1+c \alpha)} \right\}.\end{equation}

For any $k\ge 0$, define the Lyapunov function
\begin{equation}\label{eq:strong_convex_Lyapunov}   V^k \eqdef \|x^{k} - x^\star\|^2 + \frac{c\gamma^2}{M}\sum \limits_{i=1}^M \|h_i^{k} - h_i^*\|^2, \end{equation}
where $x^\star$ is the solution of \eqref{eq:main} and $h^* \eqdef \nabla f(x^\star)$. Then for all $k\geq 0$,
\begin{equation} \label{eq:strong_convex_rate} 
	\mathbb{E}\left[ V^k\right] 
	\le (1 - \gamma\mu)^k V^0 + \frac{\gamma}{\mu}(1+Mc\alpha)\frac{\sigma^2}{M}. 
\end{equation}
This implies that as long as $K \geq \frac{1}{\gamma \mu} \log \frac{V^0}{\varepsilon}$, we have $\mathbb{E}\left[ V^K\right] \leq \varepsilon + \frac{\gamma}{\mu}(1+Mc\alpha)\frac{\sigma^2}{M}.$
\end{theorem}

In particular, if we set $\gamma$ to be equal to the minimum in \eqref{eq:cond2}, then the leading term in the iteration complexity bound is $\frac{1}{\gamma \mu} = \max \left\{ \frac{1}{\alpha}, \frac{(\mu+L)(1+c\alpha)}{2\mu}\right\}.$

\begin{corollary} \label{cor:DIANA-strong-convex} Let $\kappa = \frac{L}{\mu}$, $\alpha = \frac{\alpha_p}{2}$, $c = \frac{4(1-\alpha_p)}{M\alpha_p^2}$, and $\gamma = \min\left\{\frac{\alpha}{\mu}, \frac{2}{(L+\mu)(1+c \alpha)}\right\}$. Then the conditions \eqref{eq:cond1} and \eqref{eq:cond2} are satisfied, and the leading  iteration complexity term  is equal to
\begin{equation}\label{eq:bu987gd9} \frac{1}{\gamma \mu} = \max\left\{\frac{2}{\alpha_p}, (\kappa+1)\left(\frac{1}{2} - \frac{1}{M} + \frac{1}{M\alpha_p}\right)\right\}.\end{equation} 
This is a decreasing function of $p$, and  hence $p=+\infty$ is the optimal choice.
\end{corollary}

In Table~\ref{tbl:complex_strong_conv} we calculate the leading term \eqref{eq:bu987gd9}  in the complexity of \algname{DIANA}  for $p\in \{1,2,+\infty\}$, each  for two condition number regimes: $n=\kappa$ (standard) and $n=\kappa^2$ (large).

\begin{table}
\caption{The leading term of the iteration complexity of \algname{DIANA} in the strongly convex case (Theorem~\ref{thm:DIANA-strongly_convex}, Corollary~\ref{cor:DIANA-strong-convex} and Lemma \ref{lema:alpha_p}). Logarithmic dependence on $1/\varepsilon$ is suppressed. Condition number: $\kappa= \frac{L}{\mu}$.}
\begin{center}
\footnotesize
\begin{tabular}{cccc}
\toprule
$p$ & Iteration complexity & $\kappa=\Theta(M)$ & $\kappa=\Theta(M^2)$\\
\midrule
1 &$\frac{2d}{m} + (\kappa+1)A$; \quad  $A = \left(\frac{1}{2}-\frac{1}{M}+\frac{d}{Mm}\right)$   & $\cO\left(M+\frac{d}{m}\right)$ & $\cO\left(M^2 + \frac{Md}{m}\right)$ \\
\hline
$2$ & $\frac{2\sqrt{d}}{\sqrt{m}} + (\kappa+1)B$; \quad $B = \left(\frac{1}{2}-\frac{1}{M}+\frac{\sqrt{d}}{M\sqrt{m}}\right)$  & $\cO\left(M+\sqrt{\frac{d}{m}}\right)$ & $\cO\left(M^2 + \frac{M\sqrt{d}}{\sqrt{m}}\right)$ \\
\hline
$\infty$ & $1+\sqrt{\frac{d}{m}} + (\kappa+1)C$; \quad $C = \left(\frac{1}{2}-\frac{1}{M}+\frac{1+\sqrt{\frac{d}{m}}}{2M}\right)$  & $\cO\left(M+\sqrt{\frac{d}{m}}\right)$ & $\cO\left(M^2 + \frac{M\sqrt{d}}{\sqrt{m}}\right)$ \\
\bottomrule
\end{tabular}
\label{tbl:complex_strong_conv}
\end{center}
\vspace{-20pt}
\end{table}

{\bf Matching the rate of gradient descent for quadratic size models.} Note that as long as the model size is not too big; in particular, when $d = \cO(\min\{\kappa^2,M^2\}), $
the linear rate of \algname{DIANA} with $p\geq 2$ is $\cO(\kappa \log (1/\varepsilon))$, which matches the rate of gradient descent.

{\bf Optimal block quantization.} If the dimension of the problem is large, it becomes reasonable to quantize vector's blocks, also called blocks. For example, if we had a vector which consists of 2 smaller blocks each of which is proportional to the vector of all ones, we can transmit just the blocks without any loss of information. In the real world, we have a similar situation when different parts of the parameter vector have different scale. A straightforward example is deep neural networks, layers of which have pairwise different scales. If we quantized the whole vector at once, we would zero most of the update for the layer with the smallest scale.

Our theory says that if we have $M$ workers, then the iteration complexity increase of quantization is about $\frac{\sqrt{d}}{M}$. However, if quantization is applied to a block of size $M^2$, then this number becomes 1, implying that the complexity remains the same. Therefore, if one uses about 100 workers and splits the parameter vector into parts of size about 10,000, the algorithm will work as fast as \algname{SGD}, while communicating bits instead of floats!

Some consideration related to the question of optimal number of nodes are included in Section~\ref{sec:opt_no_nodes}.


\subsection{Decreasing stepsizes}\label{sec:DIANA-decreasing-stepsizes}

We now provide a convergence result for \algname{DIANA} with decreasing stepsizes, obtaining a $\cO(1/k)$ rate.

\begin{theorem}\label{th:str_cvx_decr_steps}
    Assume that $f$ is $L$-smooth, $\mu$-strongly convex and we have access to its gradients with bounded noise. Set $\gamma_k = \frac{2}{\mu k + \theta}$ with some $\theta \ge 2\max\left\{\frac{\mu}{\alpha}, \frac{(\mu+L)(1+c\alpha)}{2} \right\}$ for some numbers $\alpha > 0$ and $c > 0$ satisfying $\frac{1+Mc\alpha^2}{1+Mc\alpha} \le \alpha_p$. After $K$ iterations of \algname{DIANA} we have
    \begin{align*}
        \mathbb{E}\left[ V^K\right]
        \le \frac{1}{\eta K+1}\max\left\{ V^0, 4\frac{(1+Mc\alpha)\sigma^2}{M\theta\mu} \right\},
    \end{align*}
    where $\eta\eqdef \frac{\mu}{\theta}$, $V^k=\|x^k - x^\star\|^2 + \frac{c\gamma_k}{M}\sumiM\|h_i^0 - h_i^*\|^2$ and $\sigma^2$ is the variance of the gradient noise.
\end{theorem}

\begin{corollary}\label{cor:str_cvx_decr_step}
	If we choose $\alpha = \frac{\alpha_p}{2}$, $c = \frac{4(1-\alpha_p)}{M\alpha_p^2}$,  $\theta=2\max\left\{\frac{\mu}{\alpha}, \frac{\left(\mu+L\right)\left(1 + c\alpha\right)}{2} \right\} = \frac{\mu}{\alpha_p}\max\left\{4, \frac{2(\kappa + 1)}{M} + \frac{(\kappa+1)(M-2)}{ M}\alpha_p\right\}$, then there are three regimes:
i) if $1 = \max\left\{1,\frac{\kappa}{M},\kappa\alpha_p\right\}$, then $\theta = \Theta\left(\frac{\mu}{\alpha_p}\right)$ and to achieve $\mathbb{E}\left[ V^k\right]\le \varepsilon$ we need at most 
\[
	\cO\left( \frac{1}{\alpha_p}\left(V^0 + \frac{(1-\alpha_p)\sigma^2}{M\mu^2} \right)\frac{1}{\varepsilon} \right)
\]
 iterations;
ii) if $\frac{\kappa}{M} = \max\left\{1,\frac{\kappa}{M},\kappa\alpha_p\right\}$, then $\theta = \Theta\left(\frac{L}{M\alpha_p}\right)$ and to achieve $\mathbb{E}\left[ V^k\right]\le \varepsilon$ we need at most 
\[
	\cO\left( \frac{\kappa}{M\alpha_p}\left(V^0 + \frac{(1-\alpha_p)\sigma^2}{\mu L} \right)\frac{1}{\varepsilon} \right)
\]
iterations;
iii) if $\kappa\alpha_p = \max\left\{1,\frac{\kappa}{M},\kappa\alpha_p\right\}$, then $\theta = \Theta\left(L\right)$ and to achieve $\mathbb{E}\left[ V^k\right]\le \varepsilon$ we need at most 
\[
	\cO\left( \kappa\left(V^0 + \frac{(1-\alpha_p)\sigma^2}{\mu Ln\alpha_p} \right)\frac{1}{\varepsilon} \right)
\]
iterations.		
\end{corollary}

\section{Convergence Theory for Non-Convex Losses}\label{sec:DIANA-non-convex}

In this section we consider the non-convex case under the following assumption which we call \textit{bounded data dissimilarity}.

\begin{assumption}[Bounded data dissimilarity]\label{as:almost_identical data}
	We assume that there exists a constant $\zeta \ge 0$ such that for all $x\in\R^d$
	\begin{equation}\label{eq:almost_identical_data}
		\frac{1}{M}\sum_{i=1}^M\|\nabla f_i(x) - \nabla f(x)\|^2 \le \zeta^2.
	\end{equation}
\end{assumption}
In particular, Assumption~\ref{as:almost_identical data} holds with $\zeta = 0$ when all $f_i$'s are the same up to some additive constant (i.e., each worker samples from one dataset). We note that it is also possible to extend our analysis to a more general assumption with extra $\cO(\|\nabla f(x)\|^2)$ term in the right-hand side of~\eqref{eq:almost_identical_data}. However, this would overcomplicate the theory without providing more insight.

\begin{theorem}\label{thm:DIANA-non-convex}
Assume that $\psi$ is constant and Assumption~\ref{as:almost_identical data} holds.    
    Also assume that $f$ is $L$-smooth, stepsizes $\alpha>0$ and $\gamma_k=\gamma>0$ and parameter $c>0$ satisfying $\frac{1 + M c \alpha^2}{1 + M c \alpha}   \leq \alpha_p,$ $\gamma \le \frac{2}{L(1+2c\alpha)}$ and $\overline x^K$ is chosen randomly from $\{x^0,\dotsc, x^{K-1} \}$. Then
    \begin{align*}
        \mathbb{E}\left[ \|\nabla f(\overline x^K)\|^2\right] 
        &\le \frac{2}{K}\frac{\Lambda^0}{\gamma(2 - L\gamma - 2c\alpha L \gamma)}+ \frac{(1+2cn\alpha)L\gamma}{2 - L\gamma - 2c\alpha L \gamma}\frac{\sigma^2}{M} + \frac{4c\alpha L\gamma \zeta^2}{2-L\gamma -2c\alpha L\gamma},
    \end{align*}
    where $\Lambda^k\eqdef  f(x^k) - f^* + c\frac{L\gamma^2}{2}\frac{1}{M}\sumiM \|h_i^{k}-h_i^*\|^2$ for any $k\ge 0$.
\end{theorem}

\begin{corollary}\label{cor:DIANA-non-convex}
	Set $\alpha = \frac{\alpha_p}{2}$, $c = \frac{4(1-\alpha_p)}{M\alpha_p^2}$, $\gamma = \frac{M\alpha_p}{L(4 + (M-4)\alpha_p)\sqrt{K}}$, $h^0 = 0$ and run the algorithm for $K$ iterations. Then, the final accuracy is at most 
	\[
		\frac{2}{\sqrt{K}} \frac{L(4+(M-4)\alpha_p)}{M\alpha_p} \Lambda^0 + \frac{1}{\sqrt{K}}\frac{(4-3\alpha_p)\sigma^2}{4+(M-4)\alpha_p} + \frac{8(1-\alpha_p)\zeta^2}{(4+(M-4)\alpha_p)\sqrt{K}}.
	\]
\end{corollary}

Moreover, if the first term in Corollary~\ref{cor:DIANA-non-convex} is leading and $\frac{1}{M} = \Omega(\alpha_p)$, the resulting complexity is $\cO(\frac{1}{\sqrt{K}})$, i.e., the same as that of \algname{SGD} . For instance, if sufficiently large mini-batches are used, the former condition holds, while for the latter it is enough to quantize vectors in blocks of size $\cO(M^2)$.

\section{Experiments}
Following advice from~\cite{alistarh2017qsgd}, we encourage the use of \textit{blocks} when quantizing large vectors. To this effect, a vector can decomposed into a number of blocks, each of which should  then be quantized separately. If coordinates have different scales, as is the case in deep learning, it will prevent undersampling of those with typically smaller values. Moreover, our theoretical results predict that applying quantization to blocks or layers will result in superlinear acceleration.

In our convex experiments, the optimal values of $\alpha$ were usually around $\min_i\frac{1}{\sqrt{d_i}}$, where the minimum is taken with respect to blocks and $d_i$ are their sizes.

Finally, higher mini-batch sizes make the sampled gradients less noisy, which in turn is favorable to more uniform differences $g_i^k - h_i^k$ and faster convergence.

Detailed description of the experiments can be found in Section~\ref{sec:A:detailsOfNumericalExperiments} as well as extra numerical results.

{\bf \algname{DIANA} with momentum works best.} We implement \algname{DIANA}, \algname{QSGD},  \algname{TernGrad} and \algname{DQGD}  in Python using MPI4PY for processes communication. This is then tested on a machine with 24 cores, each is Intel(R) Xeon(R) Gold 6146 CPU @ 3.20GHz. The problem considered is binary classification with logistic loss and $\ell_2$ penalty, chosen to be of order $1/N$, where $N$ is the total number of data points. We experiment with choices of $\alpha$, choice of norm type $p$, different number of workers and search for optimal block sizes. $h_i^0$ is always set to be zero vector for all $i$. We observe that for $\ell_{\infty}$-norm the optimal block size is significantly bigger than for $\ell_2$-norm. Here, however, we provide Figure~\ref{fig:diana_main} to show how vast the difference is with other methods.
\begin{figure}[h]
\centering

\includegraphics[scale=0.205]{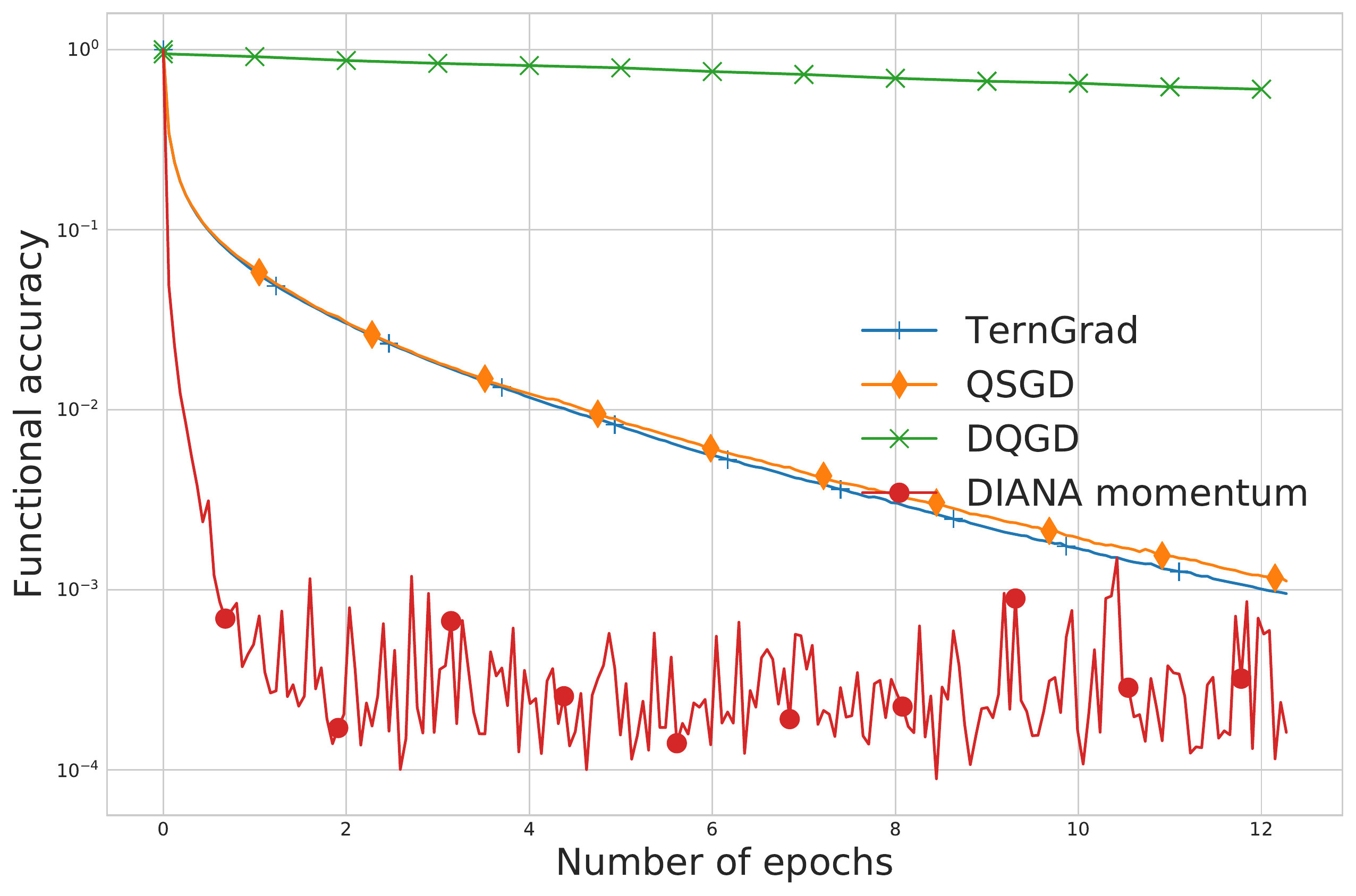}

\caption{Comparison of \algname{DIANA} ($\beta = 0.95$) with \algname{QSGD}, \algname{TernGrad} and \algname{DQGD}  on the logistic regression problem for the ``mushrooms'' dataset.}\label{fig:diana_main}
\end{figure}

{\bf \algname{DIANA} vs.\ MPI.} In Figure~\ref{fig:imagesPerSecond2}
we compare the performance of \algname{DIANA} vs.\ doing a MPI reduce operation with 32bit floats. The computing cluster had Cray Aries High Speed Network. However, for \algname{DIANA} we used 2bit per dimension and  have experienced a strange scaling behaviour, which was documented also in~\cite{parker2018performance}. In our case, this affected speed for alexnet and vgg\_a beyond 64 or 32 MPI processes respectively. For more detailed experiments, see Section~\ref{sec:A:detailsOfNumericalExperiments}.

{\bf Train and test accuracy on CIFAR-10.} In the next experiments, we run \algname{QSGD} \cite{alistarh2017qsgd}, \algname{TernGrad} \cite{wen2017terngrad}, \algname{SGD}  with momentum and \algname{DIANA} on CIFAR-10 dataset for 3 epochs. We have selected 8 workers and run each method for learning rate from $\{0.1, 0.2, 0.05\}$.
For \algname{QSGD}, \algname{DIANA} and \algname{TernGrad}, we also tried various quantization bucket sizes in $\{32, 128, 512\}$.
For \algname{QSGD} we have chosen $2,4,8$ quantization levels.
For \algname{DIANA} we have chosen $\alpha \in 
\{0, 1.0/\sqrt{\mbox{quantization bucket sizes }}\}$
and have selected initial $h^0 = 0$. 
For \algname{DIANA} and \algname{SGD}  we also run a momentum version, with a momentum parameter in $\{0, 0.95, 0.99\}$.
For \algname{DIANA} we also run with two choices of norm $\ell_2$ and $\ell_\infty$.
For each experiment we have selected softmax cross entropy loss. CIFAR-10-DNN is a convolutional DNN described here
\url{https://github.com/kuangliu/pytorch-cifar/blob/master/models/lenet.py}.

In Figure~\ref{fig:DNN:evolution2} we show the best runs over all the parameters for all the methods. We notice that \algname{DIANA} and \algname{SGD}   significantly outperform other methods.

\begin{figure}[h]

\centering 
\includegraphics[width=0.33\textwidth]{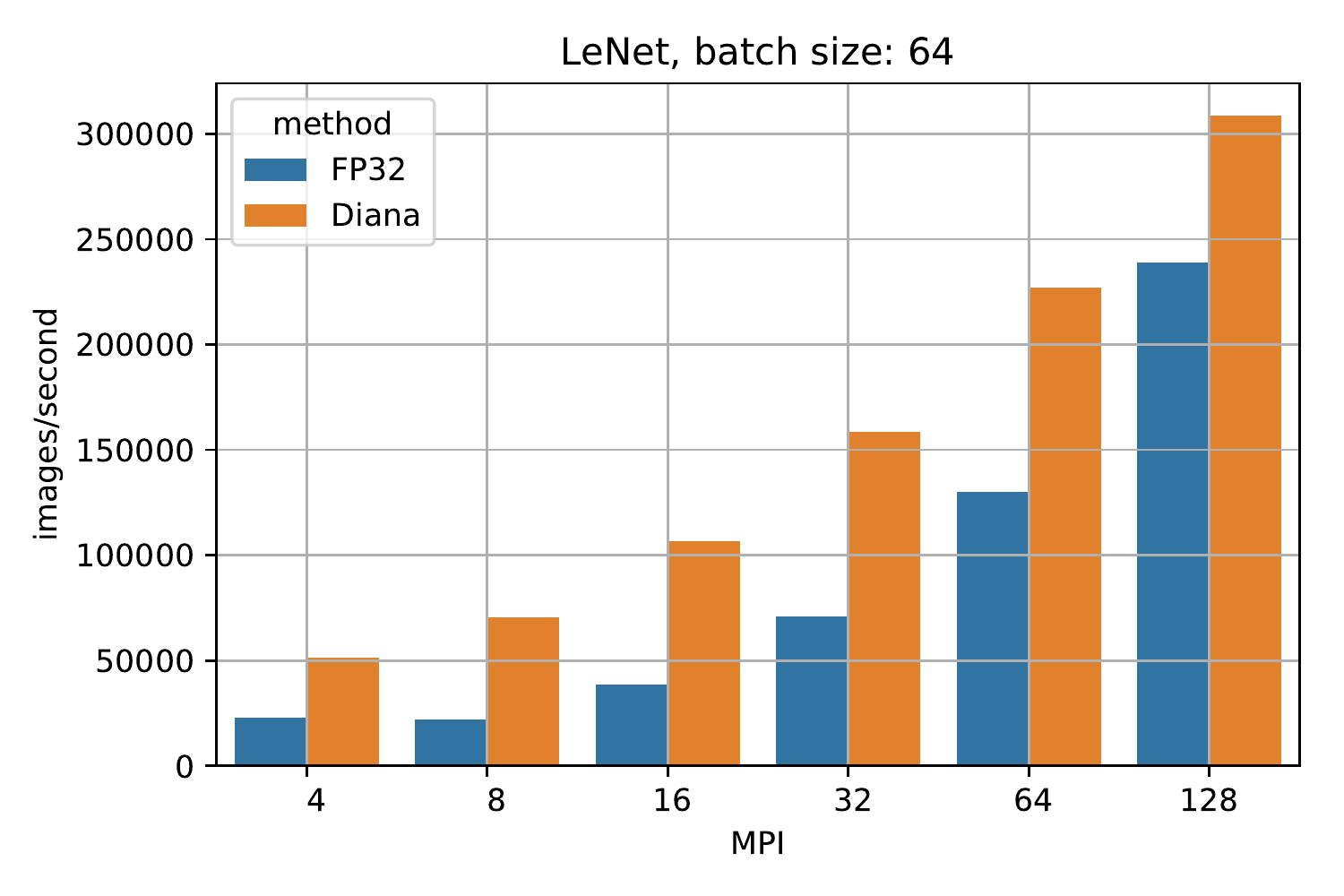}
\includegraphics[width=0.33\textwidth]{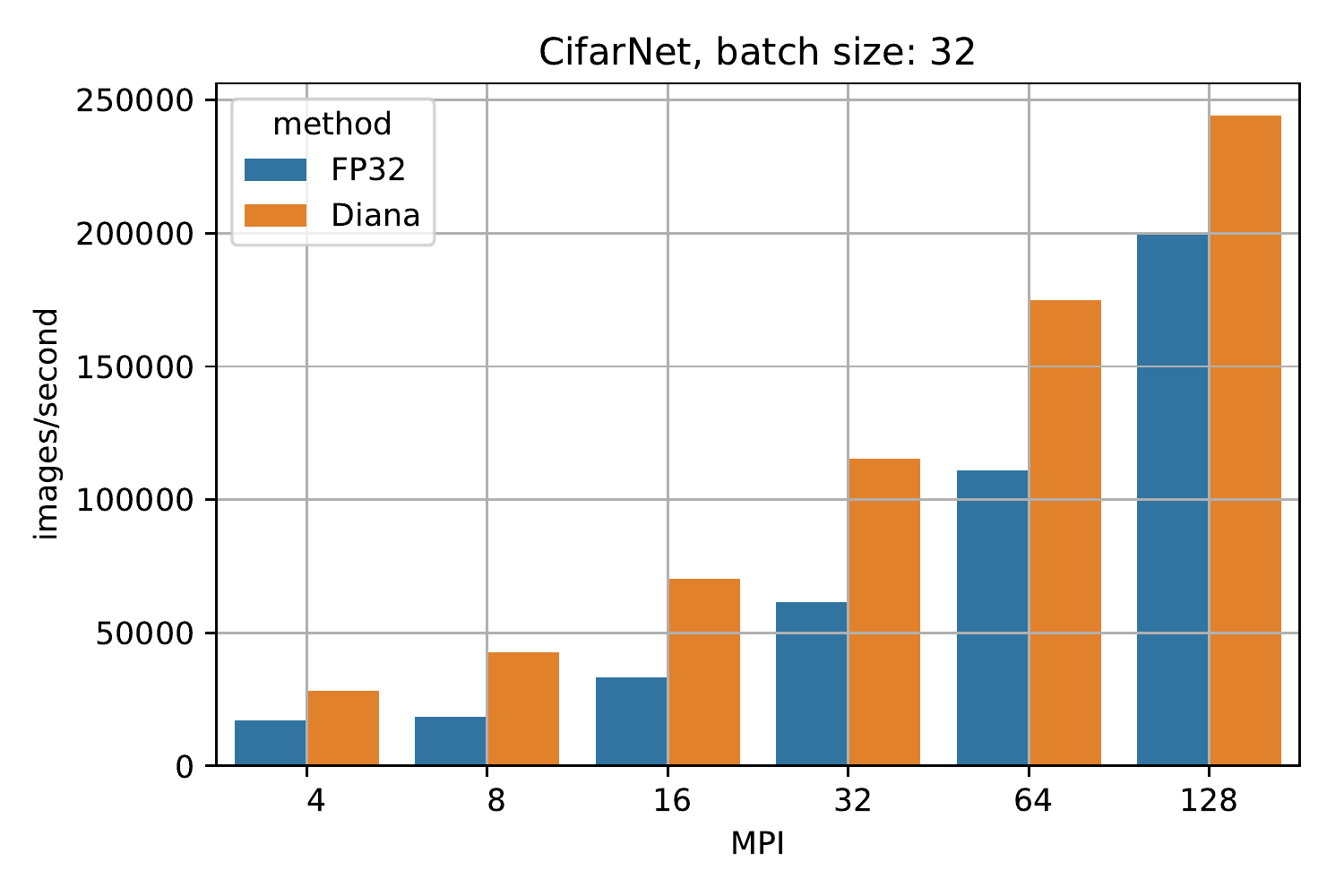}
\includegraphics[width=0.33\textwidth]{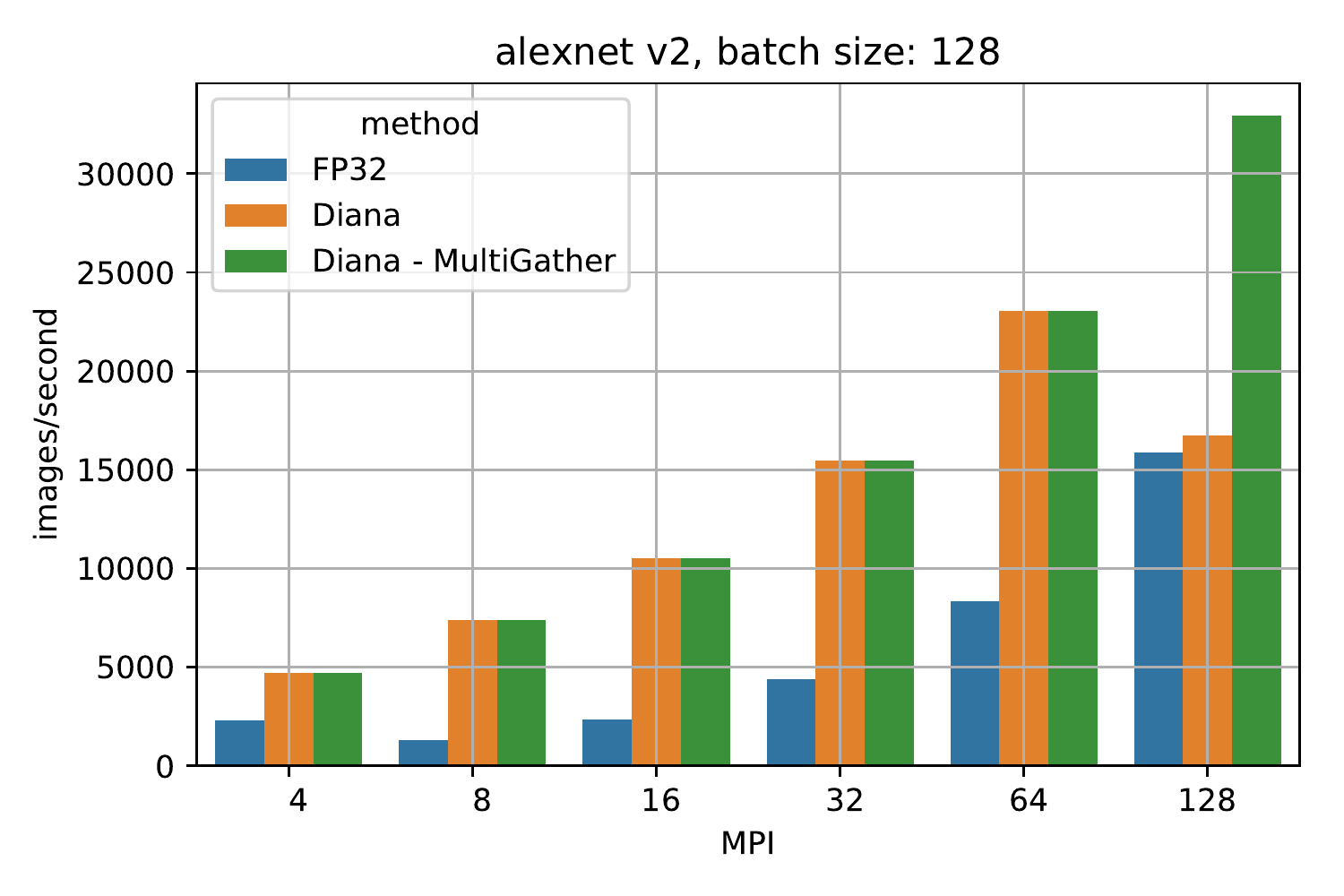}
\includegraphics[width=0.33\textwidth]{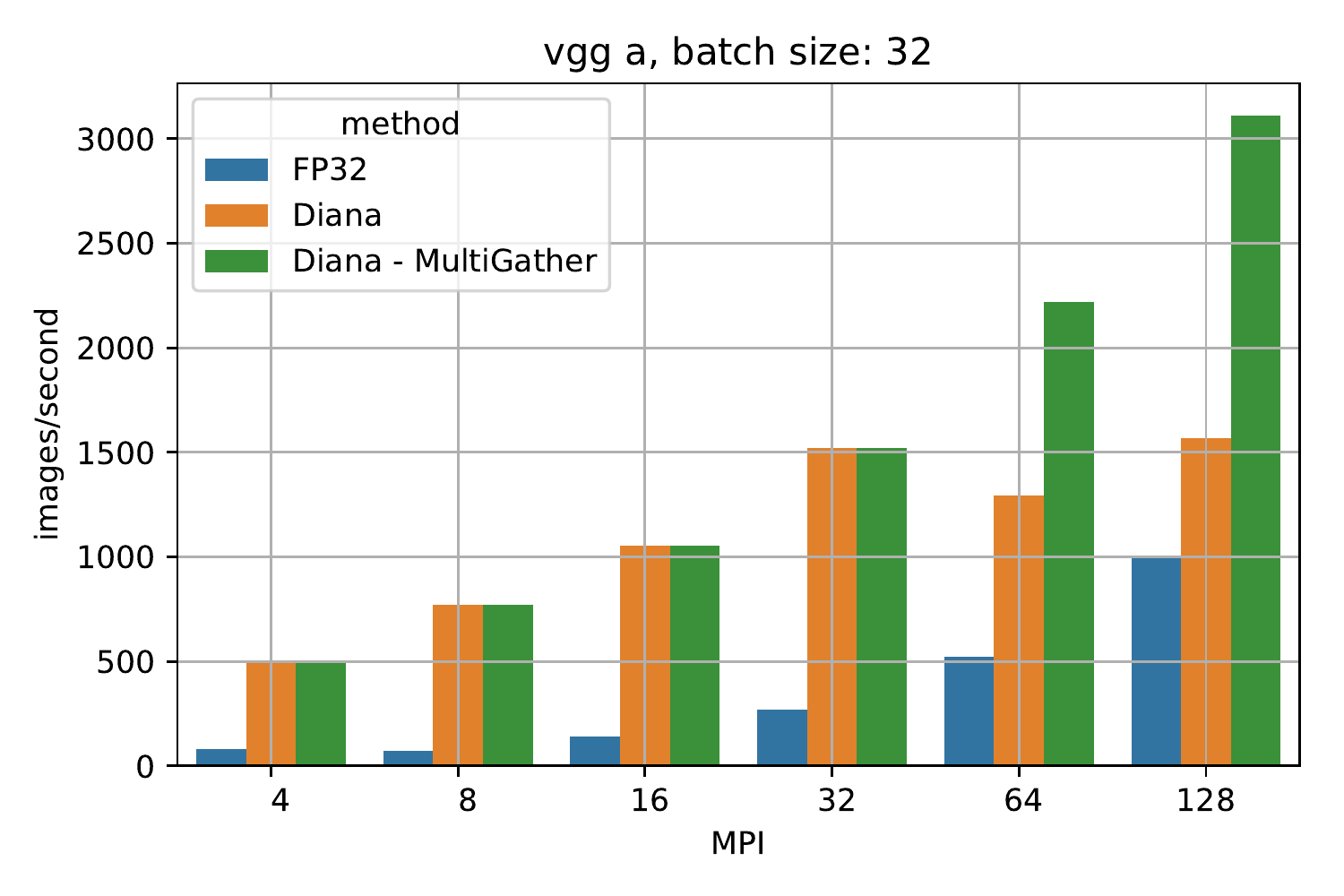}

\caption{Comparison of performance (images/second) for various number of GPUs/MPI processes and sparse communication \algname{DIANA} (2bit) vs.\ Reduce with 32bit float (FP32).}
\label{fig:imagesPerSecond2}

\end{figure}

\begin{figure}[h] 
\centering 
\includegraphics[width=0.4\textwidth]{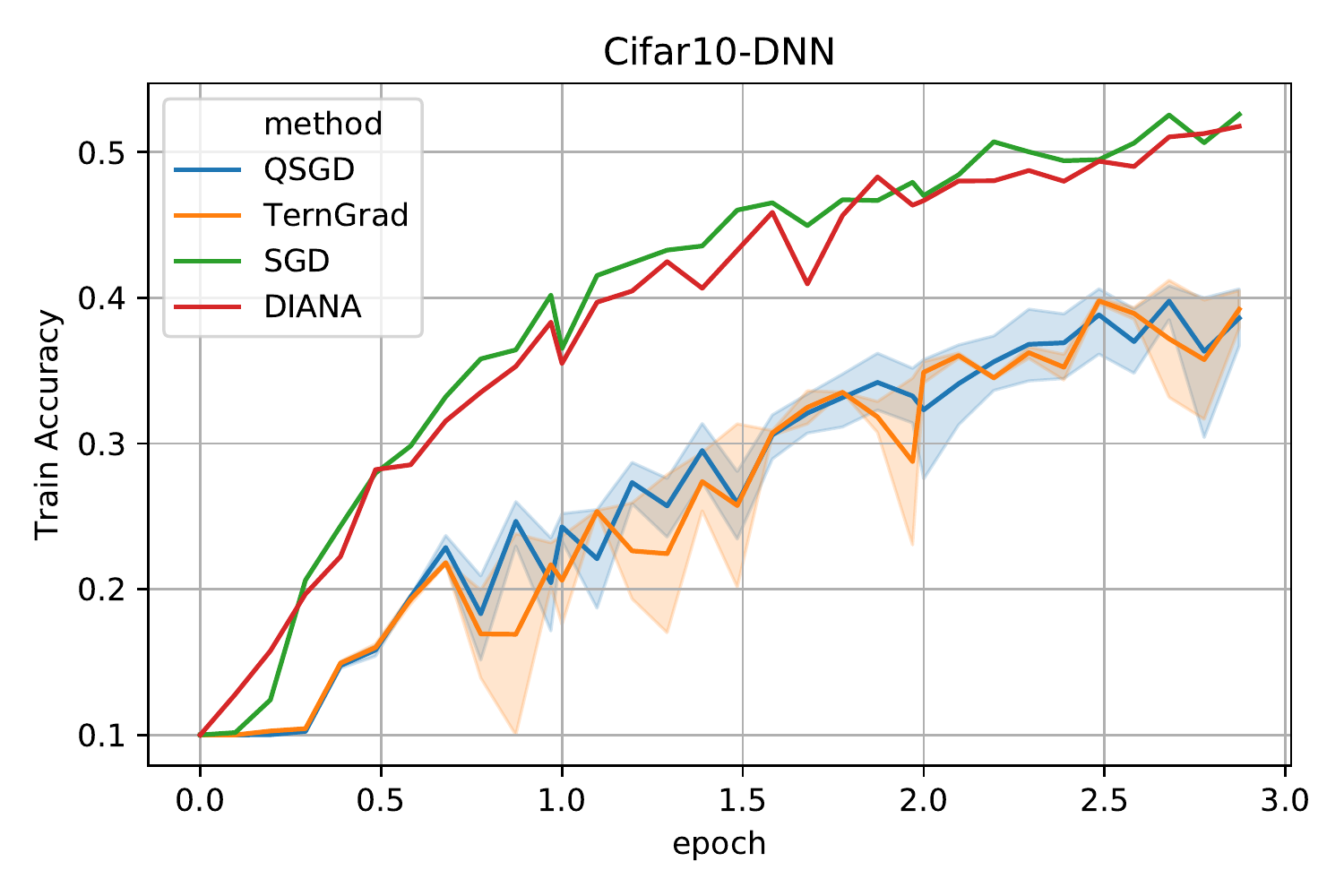}
\includegraphics[width=0.4\textwidth]{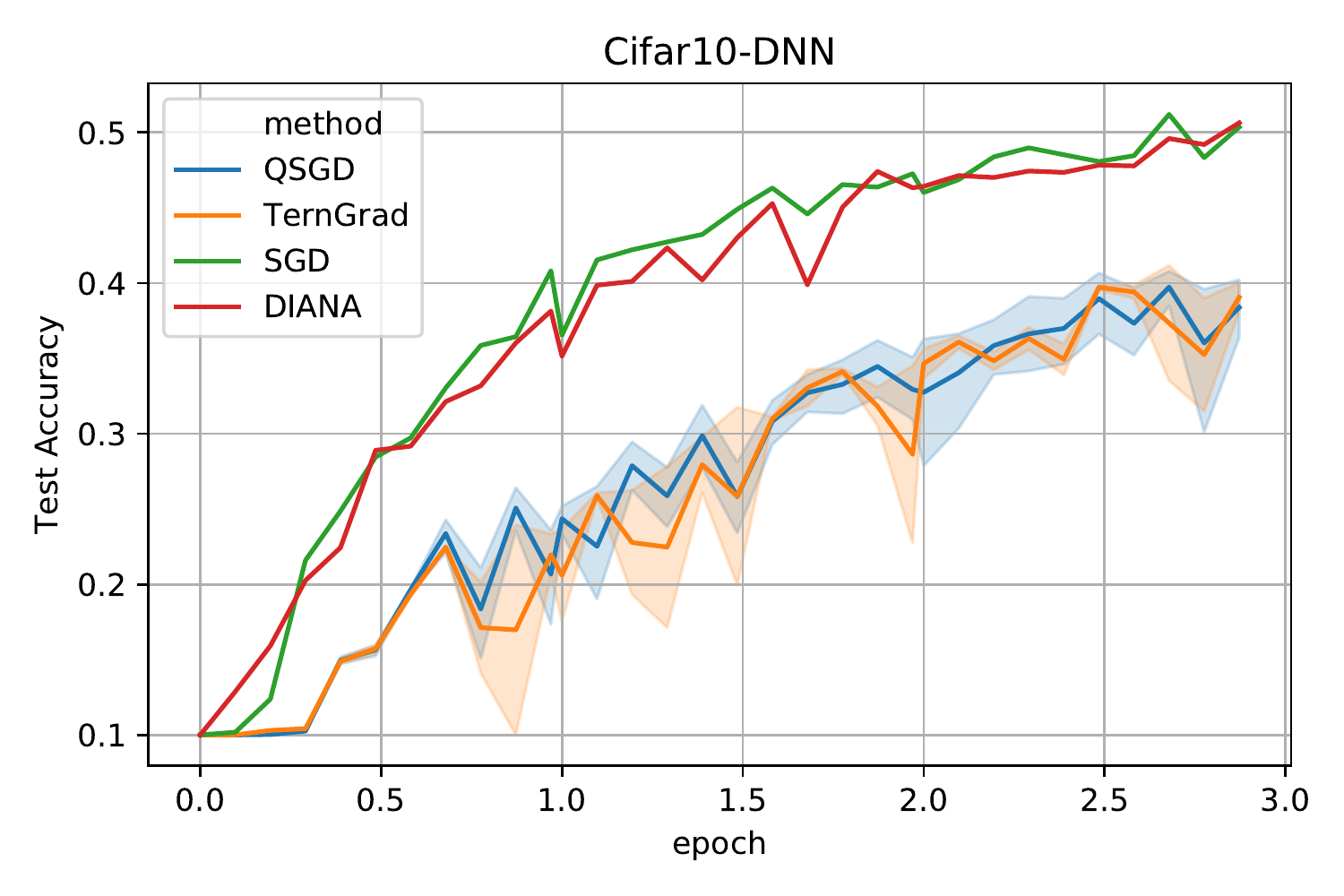}

\caption{Evolution of training (left) and testing (right) accuracy on CIFAR-10, using 4 algorithms: \algname{DIANA}, \algname{SGD}, \algname{QSGD} and \algname{TernGrad}. 
We have chosen the best runs over all tested hyper-parameters.}
\label{fig:DNN:evolution2}

\end{figure}

\chapter{Developing Variance Reduction for Proximal Operators}
\label{chapter:sdm}

\graphicspath{{sdm/}}

\section{Introduction}

In this chapter we address optimization problems of the form 
\begin{align}
	\min_{x\in\RR^d} \left[P(x) \eqdef f(x) + \frac{1}{m}\sum \limits_{j=1}^m g_j(x) + \psi(x)\right], \label{eq:pb_general}
\end{align}
where $f\colon\RR^d \to \RR$ is a smooth convex function,
and $\psi, g_1,\dotsc, g_m\colon \RR^d \to \RR\cup \{+\infty\}$ are  proper closed convex functions, admitting efficiently computable proximal operators.   We also assume throughout that $\dom(P) \eqdef \{x: P(x)<+\infty\} \neq \emptyset$ and, moreover, that the set of minimizers of \eqref{eq:pb_general}, $\cX^\star$, is non-empty.
    
    The main focus of this chapter is on how the difficult non-smooth term 
    \begin{equation}\label{eq:g_sum} 
	    g(x)\eqdef \frac{1}{m}\sum \limits_{j=1}^m g_j(x)    
    \end{equation} 
    should be treated in order to construct an efficient  algorithm for solving the problem.  We are specifically interested in the case when $m$ is very large, and when the proximal operators of  $g$ and $g+\psi$ are impossible or prohibitively difficult to evaluate. 
We thus need to rely on splitting approaches which make calls to proximal operators of functions $\{g_j\}$ and $\psi$ separately.

\section{Related Work}
Existing methods for solving problem~\eqref{eq:pb_general} can efficiently handle the case $m=1$ only \cite{alacaoglu2017smooth}. There were a few attempts to design methods capable of handling the general $m$ case, such as~\cite{allen2017katyusha, pedregosa2019proximal, ryu2017proximal} and~\cite{defazio2016simple}. None of the existing methods offer a linear rate for non-smooth problem except for random projection. In cases when sublinear rates are established, the assumptions on the functions $g_j$ are very restrictive. For instance, the results in cite{allen2017katyusha} are limited to Lipschitz continuous $g_j$ only, and Defazio~\cite{defazio2016simple} assumes $g_j$ to be strongly convex. This is very unfortunate because the majority of  problems appearing in popular data science and machine learning applications lack these properties. For instance, if we want to find a minimum of a smooth function over the intersection of $m$ convex sets, $g_j$ will be characteristic functions of sets, which are neither Lipschitz nor strongly convex. 

{\bf Applications.}  There is a long list of  applications of the non-smooth finite-sum problem \eqref{eq:pb_general}, including  convex feasibility~\cite{BauBor:96}, constrained optimization~\cite{Nocedal-Wright-book-2006}, decentralized optimization~\cite{Nedic2009distributed}, support vector machine~\cite{cortes1995support}, Dantzig selector~\cite{candes2007dantzig}, overlapping group Lasso~\cite{yuan2006model}, and fused Lasso. In Appendix~\ref{ap:applications} we elaborate in detail how these problems  can be mapped to the general problem \eqref{eq:pb_general} (in particular, see Table~\ref{tbl:summary_apps}).

{\bf Variance reduction.} Stochastic variance-reduction methods are a major breakthrough of the last decade, whose success started with  the \algname{Stochastic Dual Coordinate Ascent} (\algname{SDCA}) method~\cite{shalev2013stochastic} and the invention of the \algname{Stochastic Average Gradient} (\algname{SAG}) method~\cite{sag}. Variance reduction has attracted enormous attention and now its reach covers strongly convex, convex and non-convex~\cite{lei2017non} stochastic problems. Despite being originally developed for finite-sum problems,  variance reduction was shown to be applicable even to problems with $f$ expressed as a general expectation~\cite{lei2017less, nguyen2018inexact}. Further generalizations and extensions include variance reduction for minimax problems~\cite{palaniappan2016stochastic}, coordinate descent in the general $\psi$ case~\cite{hanzely2018sega}, and minimization with arbitrary sampling~\cite{gower2018stochastic}.  However, very little is known about variance reduction for non-smooth finite sum problems.

\section{Settings and Contributions}  

The departure point of our work is the observation that there is a class of non-smooth problems for which variance reduction is not required; these are the linear feasibility problems: given $\mA\in \RR^{m\times d}$ and $b\in \RR^m$, find $x\in \RR^d$ such that $\mA x=b$. Assuming the system is consistent, this problem can be cast as an instance of  \eqref{eq:pb_general}, with $\psi \equiv 0$ , $f(x)=\frac{1}{2}\|x\|^2$ and $g_j$ corresponding to the characteristic function of the $j$-th equation in the system. Efficient \algname{SGD} methods (or equivalently, randomized projection methods)  with linear convergence rates were recently developed for this problem \cite{gower2015randomized, richtarik2017stochastic, tu2017breaking}, as well as accelerated variants~\cite{tu2017breaking,richtarik2017stochastic, gower2018accelerated} whose linear rate yields a quadratic improvement in the iteration complexity. However, it is not known whether these or similar linear rates could be obtained when one considers $f$ to be an arbitrary smooth and strongly convex function. While our work was originally motivated by the quest to answer this question, and  we answer in the affirmative, we were able to build a much more general theory, as we explain below.

We now summarize some of the most important contributions of this chapter:

{\bf First variance reduction for $g$.}  We propose a variance-reduction strategy for progressively approximating the proximal operator of the average of a large number of non-smooth functions $g_j$ via only evaluating the proximal operator of a single function $g_j$ in each iteration. That is, unlike existing approaches, we are able to treat the difficult  term \eqref{eq:g_sum} for any $m$. Combined with a gradient-type step in $f$ (we allow for multiple ways in which the gradient estimator is built; more on that below), and a proximal step for $\psi$, this leads to a new and remarkably efficient method (Algorithm~\ref{alg:sdm}) for solving problem~\eqref{eq:pb_general}. 

{\bf Compatibility with any gradient estimator for $f$.}
Our variance-reduction scheme for the non-smooth term $g$ is  decoupled from the way we choose to construct gradient estimators for $f$.  This allows us to use the most efficient and suitable estimators depending on the structure of $f$. In this regard, two cases are of particular importance: i) $f(x) =\mathbb{E}_\xi [f(x; \xi)] $, where  $f(\cdot; \xi) \colon \RR^d\to \RR$ is almost surely convex and smooth, and ii) $f	= \frac{1}{n}\sum_i f_i$, where $\{f_i\}$ are convex and smooth.  In case i) one may consider the standard stochastic gradient estimator $\nabla f(x^k; {\xi^k})$, or a mini-batch variant thereof,  and in case ii) one may consider the batch gradient $\nabla f(x^k)$ if $n$ is small, or a variance-reduced gradient estimator, such as \algname{SVRG}~\cite{SVRG,kovalev2019don} or \algname{SAGA}~\cite{defazio2014saga, SAGA-AS}, if $n$ is large. Our general analysis allows for any estimator to be used as long as it satisfies a certain technical assumption (Assumption~\ref{as:method}). In particular, to illustrate the versatility of our approach, we show that this assumption holds for estimators used by \algname{Gradient Descent}, \algname{SVRG}, \algname{SAGA} and over-parameterized \algname{SGD}. 

{\bf Future-proof design.} Our analysis is compatible with a wide array of other estimators of the gradient of $f$ beyond the specific ones listed above. Therefore, new specific variants of our generic method for solving problem \eqref{eq:pb_general} can be obtained in the future by marrying any such new  estimators with our variance-reduction strategy for the non-smooth finite sum term $g$.

{\bf Special cases.} Special cases of our method include \algname{Randomized Kaczmarz method}~\cite{karczmarz1937angenaherte, RK}, \algname{Douglas--Rachford} splitting \cite{DRsplitting}, \algname{Forward--Backward} splitting~\cite{Nesterov_composite2013, FB2011}, a variant of \algname{SDCA}~\cite{shalev2013stochastic},  and \algname{Point--SAGA}~\cite{defazio2016simple}. Also, we obtain the first randomized variant of the famous \algname{Dykstra's algorithm}~\cite{dykstra1983algorithm} for projection onto the intersection of convex sets. These special cases are summarized in Table~\ref{tab:special_cases}.

\begin{table}[t]
    \caption{Selected special cases of our method. For \algname{Dykstra's algorithm}, $\cC_1, \dotsc, \cC_m$ are closed convex sets; and we  wish to find projection onto their intersection. \algname{Randomized Kaczmarz} is a special case for linear constraints (i.e., $\cC_j = \{x: a_j^\top x = b_j\}$). We do not prove convergence under the same assumptions as \algname{Point--SAGA} as they require strong convexity and smoothness of each $g_j$, but the algorithm is still a special case.}
   \centering    
    \footnotesize
    \begin{tabular}{cccccc}
        \toprule
         Method & $f$ & $g_j$ & $\psi$ & $\gamma$ & Citation  \\
         \midrule
         \algname{Forward--Backward} & $f_1=f$, $n=1$ & 0 & $\psi$ & $< \frac{2}{L}$ &   \cite{Nesterov_composite2013, FB2011}\\
                 \addlinespace
         \algname{Douglas--Rachford} & 0 & $g_1=g$, $m=1$ & $\psi$ & any & \cite{DRsplitting}\\ 
           \addlinespace                  
         \algname{Proximal SGD} & $\mathbb{E}_\xi [f(\cdot; \xi)] $ & 0 & $\psi$ & $\le \frac{1}{4L}$ & \cite{duchi2009efficient} \\ 
         \addlinespace
         \algname{Proximal SAGA} & $\frac{1}{n}  \sum_i f_i$ & 0 & $\psi$ & $\le \frac{1}{5L}$ & \cite{defazio2014saga} \\ 
         \addlinespace
         \algname{SDCA} & $\frac{1}{2}\|x - x^0\|^2$ & $g_j$ & 0 & $\gamma=\frac{1}{m}$ &  \cite{shalev2013stochastic} \\ 
         \addlinespace
         \algname{Randomized Dykstra's algorithm} & $\frac{1}{2}\|x - x^0\|^2$ & $\ind_{\cC_j}$ & 0 & $\gamma=\frac{1}{m}$ &  NEW \\ 
         \addlinespace
         \algname{Randomized Kaczmarz method} & $\frac{1}{2}\|x - x^0\|^2$ & $\ind_{\{x:a_j^\top x = b_j\}}$ & 0 & $\gamma=\frac{1}{m}$ &   \cite{karczmarz1937angenaherte, RK} \\
          \addlinespace
         \algname{Point--SAGA} & 0 & $g_j$ & 0 & any &  \cite{defazio2016simple} \\
          \addlinespace
         \algname{Condat--V{\~u} algorithm} & $f_1=f$, $n=1$ & $g_1=g$, $m=1$ & $\psi$ & $< \frac{2}{L}$ &   \cite{vu2013splitting, condat2013primal} \\ 
         \bottomrule 
    \end{tabular}
    \label{tab:special_cases}
\end{table}

{\bf Sublinear rates.} 
We first prove convergence of the iterates to the solution set in a Bregman sense, without quantifying the rate (see Appendix~\ref{ap:almost_sure}). Next, we establish $\cO\left(\frac{1}{K}\right)$ rate with constant stepsizes under no assumption on problem~\eqref{eq:pb_general} beyond the existence of a solution and a few technical assumptions (see Theorem~\ref{th:1_over_t_rate}). The rate improves to $\cO\left(\frac{1}{K^2}\right)$ once we assume strong convexity of $f$, and allow for carefully designed decreasing stepsizes (see Theorem~\ref{th:1_t2_rate}).

{\bf Linear rate in the non-smooth case with favourable data.} Consider the special case of \eqref{eq:pb_general} with $f$ being strongly convex, $\psi\equiv 0$ and $g_j(x) = \phi_j(\mA_j^\top x)$, where $\phi_j\colon\RR^{d_j}\to \RR\cup \{+\infty\}$ are proper closed convex functions, and $\mA_j\in \RR^{d\times d_j}$ are given (data) matrices: \begin{align}
    \min_{x\in\RR^d} f(x) + \frac{1}{m} \sum \limits_{j=1}^m \phi_j(\mA_j^\top x). \label{eq:pb_linear}
\end{align}
Let us define
\[
	\mA \eqdef [\mA_1, \dots, \mA_m]\in \RR^{d\times \sum_j d_j}.
\]
If the smallest eigenvalue of $\mA^\top \mA$ is positive, i.e., $\lambda_{\min}(\mA^\top \mA)>0$, then our method converges linearly (see Theorem~\ref{th:lin_conv_lin_model}; and note that this can only happen if $\sum_j d_j \leq d$). Moreover, picking $j$ with probability proportional to $\|\mA_j\|$ is optimal (Corollary~\ref{cor:imp_sampl}).
In the special case when $\phi_j(y)=\ind_{\{x\;:\;\mA_j^\top x = b_j\}}(x)$ for some vectors $b_1 \in \RR^{d_1},\dots, b_m\in \RR^{d_1}$, i.e., if we are minimizing a strongly convex function under a linear constraint,
\[\min_{x\in \RR^d} \left\{f(x) \;:\; \mA^\top x = b \right\},\]
then the rate is linear even if $\mA^\top \mA$ is not positive definite\footnote{By $\ind_{\cC}(x)$ we denote the characteristic function of the set $\cC$, defined as $\ind_{\cC}(x) = 0$ if $x\in \cC$ and $\ind_{\cC}(x) = +\infty$ if $x\notin \cC$}. The rate will depend on $\lambda_{\min}^+(\mA^\top \mA)$, i.e., the smallest positive eigenvalue (see Theorem~\ref{th:lin_constr}).

{\bf Linear and accelerated rate in the smooth case.}
If $g_1, \dotsc, g_m$ are smooth functions, the rate is linear (see Theorem~\ref{th:lin_conv_smooth}). If $m$ is big enough, then it is also accelerated (Corollary~\ref{cor:acc_in_g}).  A summary of our iteration complexity results is provided in Table~\ref{tab:results}.

\begin{table}[t]
    \caption{Summary of iteration complexity results that we proved. We assume by default that all functions are convex, but provide different rates based on whether $f$ is strongly convex (scvx) and whether $g_1, \dotsc, g_m$ are smooth functions, which is represented by the check marks.}
    \centering
    \footnotesize
    \begin{tabular}{cccccc}
        \toprule
         Problem & $f$ scvx & $g_j$ smooth & Oracle & Rate & Theorem \\
         \midrule
         \multirow{2}{*}{$\mathbb{E}\left[f(x; \xi)\right] + \frac{1}{m}\sum\limits_{j=1}^m g_j(x) + \psi(x)$}  & \centering {\large\red\xmark} & {\large\red\xmark} &
\multirow{3}{*}{\algname{SGD}}     & \large{$\cO\left( \frac{1}{\sqrt{K}}\right)$} & \ref{cor:sgd}\\[0.5ex]
         \cline{2-3} \cline {5-6}
         & \centering {\large\green\cmark} & {\large\red\xmark} &   &{\large$\cO\left(\frac{1}{K}\right)$}  & \ref{th:sgd_str_cvx} \\[0.5ex]
         \hline
         \multirow{4}{*}{$\frac{1}{n}\sum\limits_{i=1}^n f_i(x) + \frac{1}{m}\sum\limits_{j=1}^m g_j(x) + \psi(x)$}  & \centering {\large\red\xmark} & {\large\red\xmark} &   \multirow{6}{*}{\shortstack{\algname{GD},\\ \algname{SVRG}\\ and\\ SAGA}}  
         & {\large$\cO\left( \frac{1}{K} \right)$} & \ref{th:1_over_t_rate} \\[0.5ex]
         \cline{2-3} \cline {5-6}
         & \centering {\large\green\cmark} & {\large\red\xmark} &   & {\large$\cO\left( \frac{1}{K^2} \right)$} & \ref{th:1_t2_rate} \\[0.5ex]
         \cline{2-3} \cline {5-6}
         & \centering {\large\green\cmark} & {\large\green\cmark} &   & Linear & \ref{th:lin_conv_smooth} \\[0.5ex]
         \cline{1-3} \cline {5-6}
         $\frac{1}{n}\sum\limits_{i=1}^n f_i(x) + \frac{1}{m}\sum\limits_{j=1}^m \phi_j(\mA_j^\top x)$ & \centering {\large\green\cmark} & {\large\red\xmark} & & Linear & \ref{th:lin_conv_lin_model}, \ref{th:lin_constr}\\[0.5ex]
         \bottomrule
    \end{tabular}
    \label{tab:results}
\end{table}


{\bf Related work.} The problems that we consider recently received a lot of attention. However, we are the first to show linear convergence on non-smooth problems. $\cO\left(\frac{1}{K}\right)$ convergence with stochastic variance reduction was obtained in~\cite{ryu2017proximal} and~\cite{pedregosa2019proximal}, although both works do not have $\cO\left(\frac{1}{K^2}\right)$ rate as we do. On  the other hand, works such as~\cite{yurtsever2016stochastic,cevher2018stochastic} managed to prove $\cO\left(\frac{1}{K^2}\right)$ convergence, but only with all functions from $f$ and $g$ used at every iteration. Stochastic $\cO\left(\frac{1}{K^2}\right)$ for constrained minimization can be found in~\cite{mishchenko2018stochastic}. There is also a number of works that consider parallel~\cite{deng2017parallel} ($\cO\left(\frac{1}{K}\right)$ rate) and stochastic~\cite{zheng2016fast, liu2017accelerated} variants of \algname{ADMM}, which work with one non-smooth term composed with a linear transformation. To show linear convergence they require matrix in the transformation to be positive-definite. Variance-reduced \algname{ADMM} for compositions, which is an orthogonal direction to ours, was considered in~\cite{yu2017fast}. There is a method for non-smooth problems with $f\equiv 0$ and proximal operator preconditioning that was analyzed in detail in~\cite{chambolle2018stochastic}, we discuss the relation to it in Appendix~\ref{ap:spdhg}. Many methods were designed to work with non-smooth functions  in parallel only, and one can obtain more of them from three-operator splitting methods such as the \algname{Condat--V{\~u}} algorithm~\cite{vu2013splitting, condat2013primal}. Several works obtained linear convergence for smooth $g$~\cite{du2018linear, palaniappan2016stochastic}. Coordinate descent methods for two non-smooth functions were considered in~\cite{alacaoglu2017smooth}.

We will make the following assumption related to optimality conditions.
\begin{assumption}\label{as:optimality}
    There exists $x^\star\in \cX^\star$ and vectors $y_1^\star\in\partial g_1(x^\star), \dots, y_m^\star\in \partial g_m(x^\star)$ and $r^\star\in \partial \psi(x^\star)$ such that
$
        \nabla f(x^\star) + \avejm y_j^\star + r^\star = 0.
$
\end{assumption}
Throughout the chapter, we will assume that some $x^\star$ and $y_1^\star, \dotsc, y_m^\star$ satisfying Assumption~\ref{as:optimality} are fixed and all statements relate to these objects. We will denote $y^\star\eqdef \frac{1}{m}\sumjm y_j^\star$. 
A commentary and further details related to this assumption can be found in Appendix~\ref{sec:Opt_Cond}.

\section{The Algorithm}

\begin{algorithm}[t]
   \caption{Stochastic Decoupling Method \rrbox{(\algname{SDM})}}
   \label{alg:sdm}
\begin{algorithmic}[1]
   \Require Stepsize $\gamma$, initial vectors $x^0$, $y_1^0, \dotsc, y_m^0\in \RR^d$, probabilities $p_1,\dotsc, p_m$, oracle that gives gradient estimates, number of steps $K$
   \For{$k=0,1,\dotsc, K-1$}
	   \State Produce an estimate $v^k$ of $\nabla f(x^k)$, e.g., $v^k=\nabla f(x^k)$ 
	   \State $y^k = \frac{1}{m}\sum_{j=1}^m y_j^k$
	   \State $z^{k} = \proxR(x^k - \gamma v^k - \gamma y^k)$
	   \State Sample $j$ from $\{1,\dotsc, m\}$ with probabilities $\{p_1, \dotsc, p_m\}$ and set $\eta_j = \frac{\gamma}{mp_j}$
	   \State $x^{k+1} = \proxj(z^k + \eta_j y_j^k)$
	   \State $y_j^{t+1} = y_j^k + \frac{1}{\eta_j}(z^k - x^{k+1})$ \Comment{$y_j^{k+1}\in \partial g_j(x^{k+1})$}
   \EndFor
\end{algorithmic}
\end{algorithm}

Our method  is very general and can work with different types of gradient update. One only needs to have for each $x^k$ an estimate of the gradient $v^k$ such that $\mathbb{E}\left[v^k\right] = \nabla f(x^k)$ plus an additional assumption about its variance. We also maintain an estimate $y^k $ of full proximal step with respect to $g$, which allows us to make an intermediate step
$
    z^{k} 
    = \prox_{\gamma \psi}(x^k - \gamma v^k - \gamma y^k).
$
The key idea of this chapter is then to combine it with variance reduction in the non-smooth part. In fact, it mimics variance-reduction step from~\cite{defazio2016simple}, which was motivated by the \algname{SAGA} algorithm~\cite{defazio2014saga}. Essentially, the expression above for $z^k$ does not allow for update of $y^k$, so we do one more step,
\begin{align*}
	x^{k+1} =  \proxj(z^{k} + \eta_j y_j^k).
\end{align*}
This can additionally be rewritten using the identity $\prox_{\gamma g}(x) \in x - \gamma \partial g(\prox_{\gamma g}(x))$ as
\begin{align*}
    x^{k+1} 
    &\in x^k - \gamma (v^k + \partial \psi(z^k) + y^k) - \eta_j(\partial g_j(x^{k+1}) - y_j^k)  \approx \prox_{\gamma (\psi+g)}(x^k - \gamma \nabla f(x^k)).
\end{align*}
To make sure that the approximation works, we want to make $y_j^k$ be close to $\partial g_j(x^{k+1})$, which we do not know in advance. However, we do it in hindsight by updating $y_j^{k+1}$ with a particular subgradient from $\partial g_j(x^{k+1})$, namely 
\[
	y_j^{k+1}
	 = \frac{1}{\eta_j}(z^k + \eta_j y_j^k - \proxj(z^k + \eta_j y_j^k))\in \partial g_j(x^{k+1}).
\]

We also need to accurately estimate $\nabla f(x^k)$, and there several options for this. The simplest choice is setting $v^k = \nabla f(x^k)$. Often this is too expensive and one can instead construct $v^k$ using a variance-reduction technique, such as \algname{SAGA}~\cite{defazio2014saga} (see Algorithm~\ref{alg:v_saga}). To a reader  familiar with Fenchel duality, it might be of some interest that there is an explanation of our ideas using the dual.
Indeed, Problem~\eqref{eq:pb_general} can be recast into
\[
	\min_x \max_{y_1,\dotsc, y_m} f(x) + \psi(x) + \frac{1}{m}\sum_{j=1}^m x^\top y_j - \frac{1}{m} \sumjm g_j^\star\left(y_j\right),
\]
where $g_j^\star$ is the Fenchel conjugate of $g_j$. Then, the proximal gradient step in $x$ would be
\[
	z 
	= \prox_{\gamma \psi}\left(x - \gamma\nabla f(x) - \gamma \frac{1}{m} \sumjm y_j \right).
\]
In contrast, our update in $y_j$ is a proximal block-coordinate ascent step, so the overall process is akin to \algname{Proximal Alternating Gradient Descent-Ascent} (see~\cite{bianchi2015coordinate, combettes2015stochastic} for related ideas). However, this is neither how we developed nor analyze the method, so this should not be seen as a formal explanation.

\section{Gradient Estimators}

Since we want to have analysis that puts many different methods under the same umbrella, we need an assumption that is easy to satisfy. In particular, the following will fit our needs.
\begin{assumption}\label{as:method}
	Let $w^k \eqdef x^k - \gamma v^k$ and $w^\star \eqdef x^\star - \gamma \nabla f(x^\star)$. We assume that the oracle produces $v^k$ and (potentially) updates some other variables in such a way that for some constants $\gamma_{\max}>0$, $\omega > 0$ and nonnegative sequence $\{\cM^k\}_{k=0}^{+\infty}$, such that the following holds for any $\gamma \le \gamma_{\max}$:
	\begin{enumerate}[(a)]
		\item If $f$ is convex, then 
		\[
			\mathbb{E}\left[\|w^{k} - w^\star\|^2\right] + \cM^{k+1}
			\le \mathbb{E}\left[\|x^k - x^\star\|^2\right] - \omega\gamma\mathbb{E}\left[D_f(x^k, x^\star)\right] + \cM^k.
		\]
		\item If $f$ is $\mu$-strongly convex, then  either $\cM^k=0$ for all $k$ or there exists $\rho > 0$ such that
		\begin{align*}
			\mathbb{E}\left[\|w^{k} - w^\star\|^2\right] + \cM^{k+1}
			\le (1 -\omega\gamma\mu)\mathbb{E}\left[\|x^k - x^\star\|^2\right] + (1 - \rho)\cM^k.
		\end{align*}
	\end{enumerate}
\end{assumption}
We note that we could easily make a slightly different assumption to allow for a strongly convex $\psi$, but this would be at the cost of analysis clarity. Since the assumption above is already quite general, we choose to stick to it and claim without a proof that in the analysis it is possible to transfer strong convexity from $f$ to $\psi$.

Another observation is that part~(a) of Assumption~\ref{as:method} implies its part~(b) with $\omega/2$. However, to achieve tight bounds for \algname{Gradient Descent} we need to consider them separately.

\begin{lemma}[Proof in Appendix~\ref{ap:gd}]\label{lem:gd}
	If $f$ is convex, the \algname{Gradient Descent} estimate, $v^k=\nabla f(x^k)$, satisfies Assumption~\ref{as:method}(a) with any $\gamma_{\max} < \frac{2}{L}$, $\omega = 2 - \gamma_{\max} L$ and $\cM^k = 0$. If $f$ is $\mu$-strongly convex, \algname{Gradient Descent} satisfies Assumption~\ref{as:method}(b) with $\gamma_{\max} = \frac{2}{L + \mu}$, $\omega = 1$ and $\cM^k=0$.
\end{lemma}
Since $\cM^k=0$ for \algname{Gradient Descent}, one can ignore $\rho$ in the convergence results or treat it as $+\infty$.

\begin{lemma}[Proof in Appendix~\ref{ap:svrg_saga}]\label{lem:svrg_saga}
	In \algname{SVRG} and \algname{SAGA}, if $f_i$ is $L$-smooth and convex for all $i$, Assumption~\ref{as:method}(a) is satisfied with $\gamma_{\max} = \frac{1}{6L}$, $\omega = \frac{1}{3}$ and 
$	\cM^k 
		= \frac{3\gamma^2}{n} \sum_{i=1}^n \mathbb{E}\left[\|\nabla f_i(u_i^k) - \nabla f_i(x^\star)\|^2\right],
$
	where in \algname{SVRG} $u_i^k=u^k$ is the reference point of the current loop, and in \algname{SAGA} $u_i^k$ is the point whose gradient is stored in memory for function $f_i$. If $f$ is also strongly convex, then Assumption~\ref{as:method} holds with $\gamma_{\max} = \frac{1}{5L}$, $\omega = 1$, $\rho = \frac{1}{3n}$ and the same $\cM^k$.
\end{lemma}

\begin{algorithm}[t]
   \caption{\algname{SAGA} Oracle}
   \label{alg:v_saga}
\begin{algorithmic}[1]
	\Require $x^k$, table of past gradients $\nabla f_1(u_1^k), \dotsc, \nabla f_n(u_n^k)$ and their average $\alpha^k$
	\State Sample subset $S$ from $\{1,\dotsc, n\}$ of size $\tau$
	\State $v^k= \frac{1}{\tau}\sum_{i\in S}\left(\nabla f_i(x^k) - \nabla f_i(u_i^k)\right) + \alpha^k$
	\State For all $i\in S$ update  $\nabla f_i(u_i^{k+1})$ with $u_i^{k+1} = x^k$
	\\ \Return $v^k$
\end{algorithmic}
\end{algorithm}

\begin{lemma}[Proof in Appendix~\ref{ap:sgd}]\label{lem:sgd}
	Assume that at an optimum $x^\star$ the variance of stochastic gradients is finite, i.e., 
	\[
		\sigma_\star^2\eqdef \mathbb{E}_\xi \left[ \|\nabla f(x^\star; \xi) - \nabla f(x^\star)\|^2\right] < +\infty.
	\]
	Then, \algname{SGD} that terminates after at most $K$ iterations satisfies Assumption~\ref{as:method}(a) with $\gamma_{\max}=\frac{1}{4L}$, $\omega=1$ and $\rho=0$. In this case, sequence $\{\cM^k\}_{k=0}^{K}$ is given by
$
		\cM^k = 2\gamma^2(K - k)\sigma_\star^2.
$
	If $f$ is strongly convex and $\sigma_\star=0$, it satisfies Assumption~\ref{as:method}(b) with $\gamma_{\max} = \frac{1}{2L}$, $\omega=1$ and $\cM^k=0$.
\end{lemma}

There are two important cases for \algname{SGD}. If the model is overparameterized, i.e., $\sigma_\star\approx 0$,  we get almost the same guarantees for \algname{SGD} as for \algname{GD}. If,  $\sigma_\star\gg 0$, then one needs to choose $\gamma = \cO\left(\frac{1}{\sqrt{K} L}\right)$ in order to keep $\cM^0$ away from $+\infty$. This effectively changes the $\cO\left(\frac{1}{K}\right)$ rate to $\cO\left(\frac{1}{\sqrt{K}}\right)$, see Corollary~\ref{cor:sgd}. Moreover,  obtaining a $\cO\left(\frac{1}{K}\right)$ rate for strongly convex case requires a separate proof.

\section{Convergence Theory}
Let  
\[
	\nu \eqdef \min_{j=1,\dotsc,m} \frac{1}{\eta_j L_j} = \min_{j=1,\dotsc,m} m\frac{p_j}{\gamma L_j},
\]
where $L_j\in\RR\cup \{+\infty\}$ is the smoothness constant of $g_j$, in most cases giving $L_j=+\infty$ and $\nu=0$. Tho goal of our analysis is to show that with introducing new term in the Lyapunov function,
\[
	\cY^k
	\eqdef (1+\nu)\sumlm \eta_l^2\mathbb{E}\left[\|y_l^{k} - y_l^\star\|^2\right],
\]
the convergence is not significantly hurt. This term will be always incorporated in the full Lyapunov function defined as
\begin{align*}
        \cL^k
        \eqdef  \mathbb{E}\left[\|x^k - x^\star\|^2\right] + \cM^k + \cY^k,
\end{align*}
where $\cM^k$ is from Assumption~\ref{as:method}. In the proof of $\cO\left(\frac{1}{K^2}\right)$ rate we will use decreasing stepsizes and $\cY^k$ will be defined slightly differently, but except for this, it is going to be the same Lyapunov function everywhere.

\subsection{$\cO(\frac{1}{K})$ convergence for general convex problem}

\begin{theorem}[Proof in Appendix~\ref{ap:1_t_rate}]\label{th:1_over_t_rate}
    Assume $f$ is $L$-smooth and $\mu$-strongly convex, $g_1, \dotsc, g_m, \psi$ are proper, closed and convex. If we use a method for generating $v^k$ which satisfies Assumption~\ref{as:method} and $\gamma\le \gamma_{\max}$, then
\[
        \mathbb{E}\left[D_f(\overline x^K, x^\star)\right]
        \le \frac{1}{\omega \gamma K}\cL^0,
\]
    where $\cL^0\eqdef \|x^0 - x^\star\|^2 + \cM^0 + \sumlm  \eta_l^2 \|y_l^0 - y_l^\star\|^2$ and $\overline x^K \eqdef \frac{1}{K}\sum_{l=0}^{k-1} x^l$.
\end{theorem}

If $\psi\equiv 0$ and $g_j\equiv 0$ for all $j$, then this transforms into $\cO(\frac{1}{K})$ convergence of $\ec{f(x^K) - \min_x f(x)}$, which is the correct rate of \algname{SGD}.

The next result takes care of the case when \algname{SGD} is used, which requires special consideration.
\begin{corollary}\label{cor:sgd}
	If we use \algname{SGD} for $k$ iterations with constant stepsize, the method converges to a neighborhood of radius $\frac{\cM^0}{\gamma k}=2\gamma\sigma_\star^2$. If we choose the stepsize $\gamma = \Theta\left(\frac{1}{L\sqrt{K}} \right)$, then $2\gamma\sigma_\star^2=\cO(\frac{1}{\sqrt{K}})$, and we recover $\cO\left(\frac{1}{\sqrt{K}}\right)$ rate.
\end{corollary}

\subsection{$\cO(\frac{1}{K^2})$ convergence for strongly convex $f$}
In this section, we consider a variant of Algorithm~\ref{alg:sdm} with time-varying stepsizes,
\begin{align*}
	z^k 
	&= \prox_{\gamma_{k} \psi}(x^k - \gamma_{k} v^k - \gamma_{k} y^k), \qquad x^{k+1} 
	= \prox_{\eta_{k,j} g_j}(z^k + \eta_{k,j} y_j^k).	
\end{align*}

\begin{theorem}[Proof in Appendix~\ref{ap:1_t2_rate}]\label{th:1_t2_rate}
	Consider updates with time-varying stepsizes, $\gamma_k = \frac{2}{\mu\omega(a + k + 1)}$ and $\eta_{k,j} = \frac{\gamma_k}{m p_j}$ for $j=1,\dotsc, m$, where $a\ge 2\max\left\{\frac{1}{\omega\mu\gamma_{\max}}, \frac{1}{\rho} \right\}$. Then
\[ 
		\mathbb{E}\left[\|x^K - x^\star\|^2 \right]
		\le \frac{a^2}{(K+ a-1)^2}\cL^0,
\]
	where $\cL^0 = \|x^0 - x^\star\|^2 + \cM^0 + \sum_{j=1}^m \eta_{0,j}^2 \|y_j^0 - y_j^\star\|^2$.
\end{theorem}
This improves upon $\cO(\frac{1}{K})$ convergence proved in~\cite{defazio2016simple} under similar assumptions and matches the bound in~\cite{chambolle2018stochastic}.

In Corollary~\ref{cor:sgd} we obtained $\cO(\frac{1}{\sqrt{K}})$ rate for \algname{SGD} with $\sigma_\star\neq 0$. It is not surprising that the rate is worse as it is so even with $g\equiv 0$. For standard \algname{SGD} we are able to improve the guarantee above to $\cO(\frac{1}{K})$ when the objective is strongly convex.

\begin{theorem}[Proof in Appendix~\ref{ap:sgd_str_cvx}]\label{th:sgd_str_cvx}
	Assume $f$ is $\mu$-strongly convex, $f(\cdot; \xi)$ is almost surely convex and $L$-smooth. Let the update be produced by \algname{SGD}, i.e., $v^k = \nabla f(x^k; \xi^k)$, and let us use time-varying stepsizes $\gamma_{k-1} = \frac{2}{a + \mu k}$ with $a\ge 4L$. Then
\[ 
		\mathbb{E}\left[\|x^K - x^\star\|^2 \right]
		\le \frac{8\sigma_\star^2}{\mu(a + \mu K)} + \frac{a^2}{(a + \mu K)^2}\cL^0.
\]
\end{theorem}

\subsection{Linear convergence for linear non-smoothness} \label{sec:linear_non_smoothness}

We now provide two linear convergence rates in the case when $\psi\equiv 0$ and $g_j(x) = \phi_j(\mA_j^\top x)$. 

\begin{theorem}[Proof in Appendix~\ref{ap:lin_conv_lin_model}]\label{th:lin_conv_lin_model}
    Assume that $f$ is $\mu$-strongly convex, $\psi\equiv 0$, $g_j(x) = \phi_j(\mA_j^\top x)$ for $j=1,\dotsc, m$ and take a method satisfying Assumption~\ref{as:method} with $\rho>0$. Then, if $\gamma\le \gamma_{\max}$,
\[
        \mathbb{E}\left[\|x^K - x^\star\|^2 \right]
        \le \left(1 - \min\{\rho, \omega\gamma\mu, \rho_{A}\} \right)^K \cL^0,
\]
    where $\rho_{A} \eqdef \lambda_{\min}(\mA^\top\mA) \min_j \left(\frac{p_j}{\|\mA_j\|}\right)^2$, and $\cL^0\eqdef \|x^0 - x^\star\|^2 + \cM^0 + \sum_{j=1}^m  \gamma_j^2 \|y_j^0 - y_j^\star\|^2$.
\end{theorem}

\begin{corollary}\label{cor:imp_sampl}
    If oracle from Algorithm~\ref{alg:v_saga} (\algname{SAGA}) is used with probabilities $p_j \propto \|\mA_j\|$, then to get $\mathbb{E}\left[\|x^K - x^\star\|^2\right] \le \varepsilon$, it is enough to run it for 
    \[ 
    		K=\cO\left(\left(n+ \frac{L}{\mu}+ \frac{\|\mA\|_{2, 1}^2}{\lambda_{\min} (\mA^\top \mA)}\right)\log\frac{1}{\varepsilon}\right)
    \]
     iterations.
\end{corollary}
Now let us show that this can be improved to depend only on positive eigenvalues if the problem is linearly constrained.
\begin{theorem}[Proof in 	Appendix~\ref{ap:lin_constr}]\label{th:lin_constr}
	Under the same assumptions as in Theorem~\ref{th:lin_conv_lin_model} and assuming, in addition, that $g_j = \ind_{\{x: \mA_j^\top x = b_j\}}$ it holds 
	\[
		\mathbb{E}\left[\|x^k - x^\star\|^2\right] \le (1 - \min\{\rho, \omega\gamma\mu, \rho_{\mA}\})^k\cL^0
	\]
	 with $\rho_{\mA}=\lambda_{\min}^+(\mA^\top\mA)\min_{j}\left(\frac{p_j}{\|\mA_j\|}\right)^2$, i.e., $\rho_{\mA}$ depends only on  the smallest positive eigenvalue of $\mA^\top\mA$.
\end{theorem}
One implication of Theorem~\ref{th:lin_constr} is that just by taking a solver such as \algname{SVRG} we immediately obtain a method for decentralized optimization that will converge linearly. Furthermore, if the problem is ill-conditioned or the communication graph is well conditioned, the leading term is still $\frac{L}{\mu}$, meaning that the rate for decentralized method is the same as for centralized up to constant factors. In Appendix~\ref{ap:lin_constr}, we also give a version of our method specialized to the linearly constrained problem that requires only one extra vector, $y^k$.

\subsection{Linear convergence if all $g_j$ are smooth}
\begin{theorem}[Proof in Appendix~\ref{ap:lin_conv_smooth}]\label{th:lin_conv_smooth}
    Assume that $f$ is $L$-smooth and $\mu$-strongly convex, $g_j$ is $L_j$-smooth for all $j$, Assumption~\ref{as:method}(b) is satisfied and $\gamma\le \gamma_{\max}$. Then, Algorithm~\ref{alg:sdm} converges as
\[
		\mathbb{E}\left[\|x^K - x^\star\|^2\right]
		\le \left(1 - \min\left\{\omega\gamma\mu, \rho, \frac{\nu}{m(1+\nu)}\right\}\right)^K\cL^0,
\]   
    where $\nu \eqdef \min_{j=1,\dotsc,m} \frac{1}{\eta_j L_j}$.
\end{theorem}

Based on the theorem above, we suggest to choose probabilities $p_j$ to maximize $\nu$, which can be done by using $p_j\propto L_j$. If $p_j = \frac{L_j}{\sum_{k=1}^m L_k}$, then $\nu = \min_{j=1,\dotsc,m} \frac{m p_j}{\gamma L_j} = \frac{1}{\gamma \overline L}$ with $\overline L\eqdef \frac{1}{m}\sum_{j=1}^m L_j $.
\begin{corollary}[Proof in Appendix~\ref{ap:acc_in_g}]\label{cor:acc_in_g}
	Choose as solver for $f$ \algname{SVRG} or \algname{SAGA} without mini-batching, which satisfy Assumption~\ref{as:method} with $\gamma_{\max}=\frac{1}{5L}$ and $\rho = \frac{1}{3n}$, and consider for simplicity situation where $L_1= \dotsb = L_m \eqdef L_g$ and $p_1=\dotsb=p_m$. Define $\gamma_{\mathrm{best}} \eqdef (\omega \mu m L_g)^{-\frac{1}{2}}$,
	and set the stepsize to  $\gamma=\min\{\gamma_{\max}, \gamma_{\mathrm{best}}\}$. Then the complexity to get $\mathbb{E}\left[\|x^K - x^\star\|^2\right]\le \varepsilon$ is
\[
		\cO\left(\left(n + m + \frac{L}{\mu} + \sqrt{\frac{mL_g}{\mu}}\right) \log\frac{1}{\varepsilon}\right).
\]
\end{corollary}
Notably, the rate in Corollary~\ref{cor:acc_in_g} is accelerated in $g$, suggesting that the proposed update is in some cases optimal. Moreover, if $m$ becomes large, the last term is dominating everything else meaning that acceleration in $f$ might not be needed at all.

\section{Experiments}

{\bf Randomly generated linear system.}
In this experiment, we first generate a matrix with independent Gaussian entries of zero mean and scale $\frac{1}{\sqrt{d}}$, where $d=100$, and after that we set $\mW\in\RR^{d\times d}$ to be the product of the generated matrix with itself plus identity matrix with coefficient $10^{-2}$ to make sure $\mW$ is positive definite. We also generated a random vector $x^\star\in\RR^d$ and took $b = \mW x^\star$. The problem is to solve $\mW x=b$, or, equivalently, to minimize $\|\mW x - b\|^2$. We made this choice because it makes estimation of the parameters of accelerated Sketch-and-Project easier.

To run our method, we choose \[f(x) = \frac{1}{2}\|x\|^2\] and \[g_j(x) = \ind_{ \{x: w_j^\top x = b_j\}}(x), \quad j=1,\dotsc, d,\] where $\ind_{\{x \;:\; w_j^\top x=b_j\}}(x)$ is the characteristic function, whose value is $0$ if $w_j^\top x=b_j$ and $+\infty$ otherwise. Then, the proximal operator of $g_j$ is the projection operator onto the corresponding constraint. We found that the choice of stepsize is important for fast convergence and that the value approximately equal $1.3\cdot 10^{-4}\ll 1 = \frac{2}{(L + \mu)}$ led to the best performance for this matrix. 

We compare our method to the accelerated Sketch-and-Project method of~\cite{gower2018accelerated} using optimal parameters. The other method that we consider is classic Kaczmarz method that projects onto randomly chosen constraint. We run all methods with uniform sampling.

\textbf{Linear regression with linear constraints.} We took A9a dataset from LIBSVM and ran $\ell_2$-regularized linear regression, using first 50 observations of the dataset as tough constraints. We compare iteration complexity to precise projection onto all constraints and observe that it takes almost the same number of iterations, although stochastic iterations are significantly cheaper. For each method we chose mini-batch of size 20 and stepsizes of order $\frac{1}{L}$ for all methods.

More experiments are provided in Appendix~\ref{sec:add_experiments}.
\begin{figure}
\center
	\includegraphics[scale=0.25]{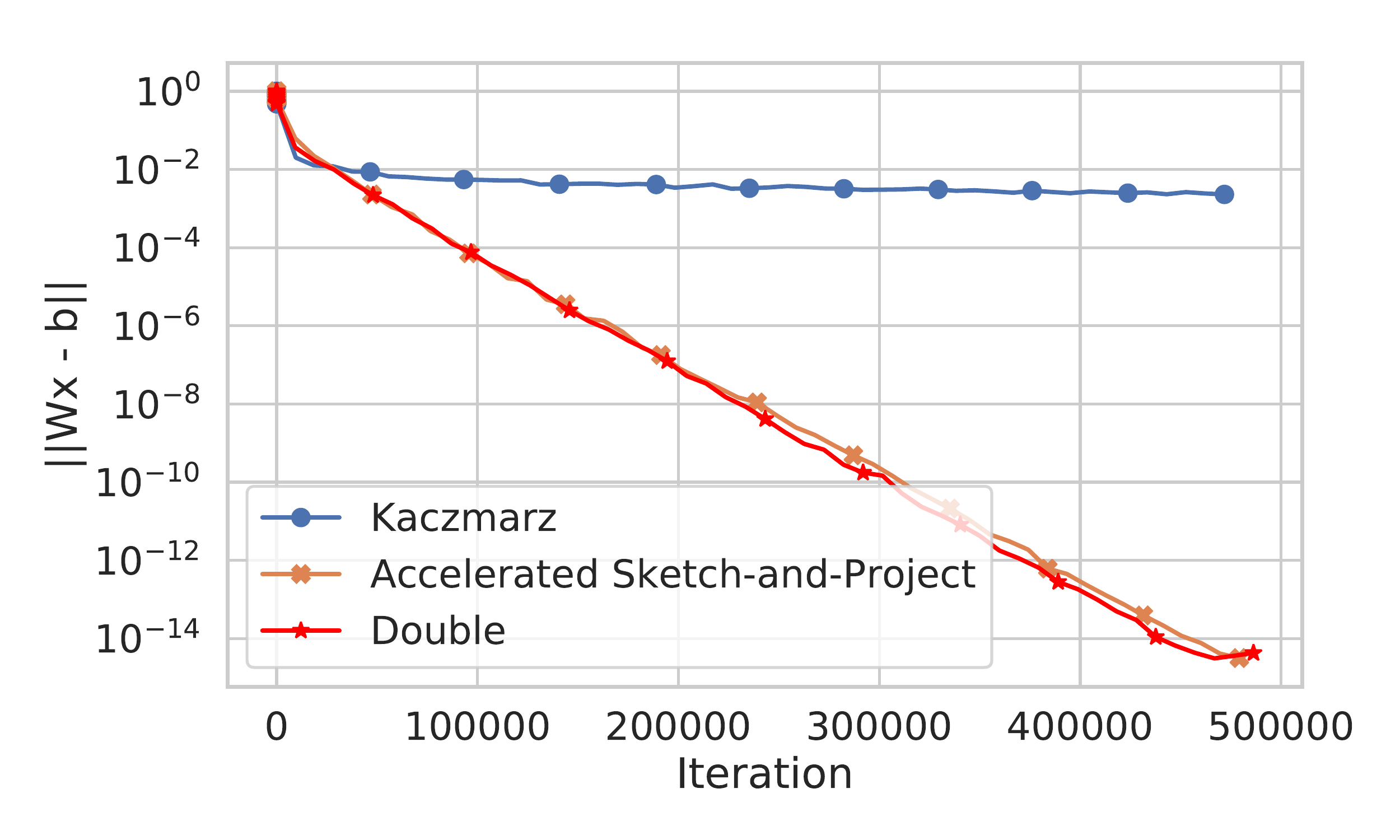}
	\includegraphics[scale=0.23]{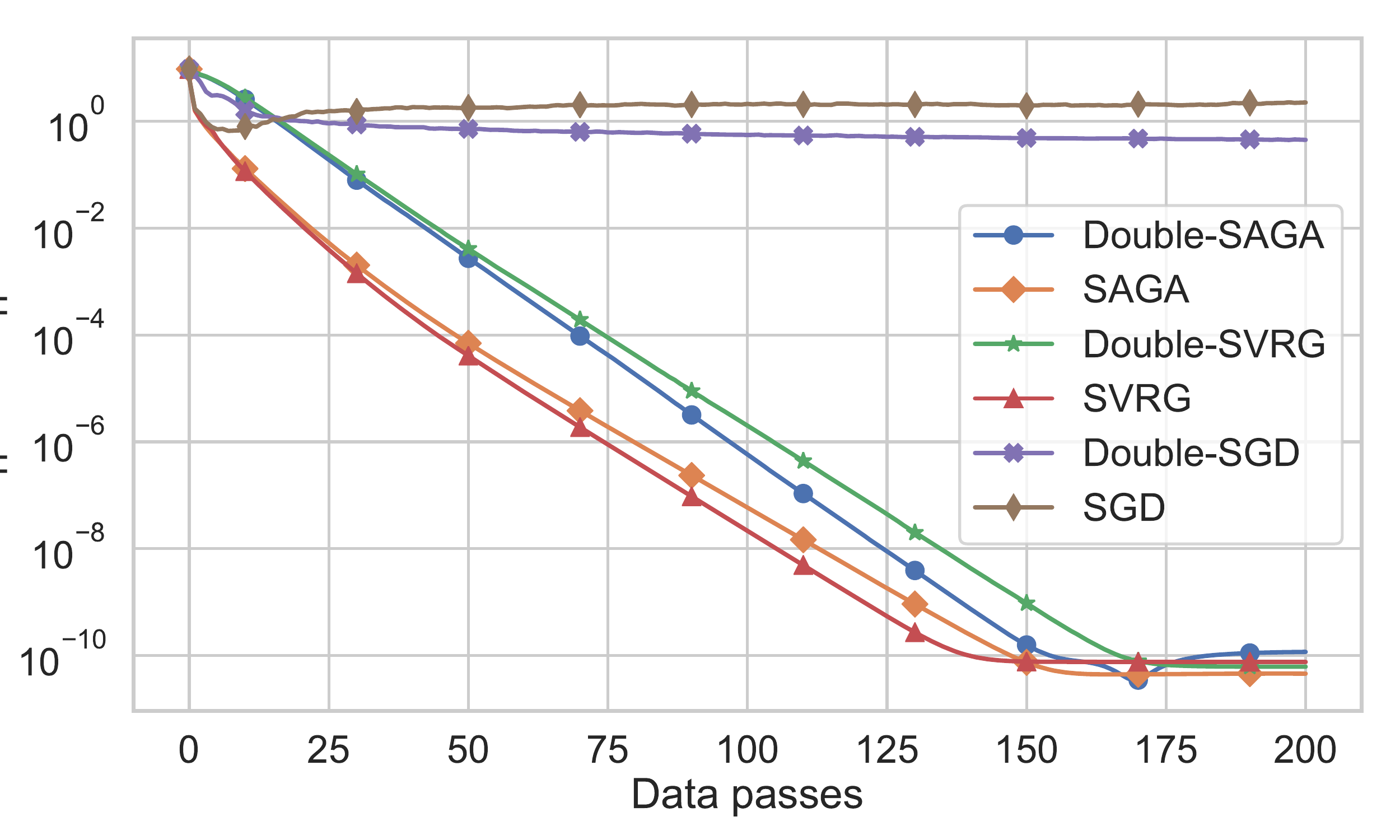}
	\caption{Left: convergence of the Stochastic Decoupling method, Kaczmarz and accelerated Kaczmarz of~\cite{gower2018accelerated} when solving $\mW x = b$ with random positive-definite $\mW\in\RR^{d\times d}$, where $d=100$. It is immediate to observe that the method we propose performs on a par with the accelerated Sketch-and-Project. Right: linear regression with A9a dataset from LIBSVM~\cite{chang2011libsvm} with first 50 observation used as linear constraints. We compare convergence of \algname{SVRG}, \algname{SAGA} and \algname{SGD} with full projections (labeled as 'SVRG', 'SAGA', 'SGD') to the same methods combined with Algorithm~\ref{alg:sdm} (labeled as 'Double-').}
\end{figure}

\chapter{Designing Variance-Reduced Algorithms with Splitting}
\label{chapter:pddy}

\graphicspath{{pddy/}}

\section{Introduction}
Many problems in statistics, machine learning or signal processing 
can be formulated as high-dimensional convex optimization problems~\cite{palomar2010convex, starck2010sparse, bach2012optimization, polson2015proximal, chambolle2016introduction, stathopoulos2016operator}. These optimization problems typically involve a smooth term $F$ and a non-smooth regularization $G$, and are often solved using a (variant of) the \algname{Proximal Stochastic Gradient Descent} (\algname{SGD})~\cite{atchade2017perturbed}. However, in many cases, $G$ is not proximable, i.e., its proximity operator does not admit a closed form expression. 

In particular, structured regularizations~\cite{chambolle2016introduction,cremers2011convex} like the total variation regularization over a graph~\cite{condat2017discrete,duran2016collaborative,bredies2010total,wan16} or the overlapping group Lasso~\cite{bach2012optimization} are known to have an expensive proximity operators~\cite{salim2019snake}. Another example is the case of linear constraints on the optimization problem. This corresponds to $G$ being an indicator function and the proximity operator of $G$ being the projection onto the constraints space. This projection requires the resolution of a high-dimensional linear system~\cite{bauschke2017convex} often intractable. The context of decentralized optimization~\cite{xu2020distributed}, in which a network of computing agents aims at jointly minimizing an objective function by performing local computations and exchanging information along the edges, is a particular case of the context of linearly constrained optimization. In this particular case, projecting onto the constraints space is equivalent to averaging across the network, which is prohibited. Finally, when $G$ is a sum of several regularizers, $G$ is not proximable even if the regularizers are proximable, because the proximity operator is not linear. 

Although in these examples $G$ is not proximable, $G$ takes the form $G = \psi + H \circ \mL$ where $\psi,H$ are proximable and $\mL$ is a linear operator\footnote{In these contexts, $H \circ \mL$ is not proximable as well (the symbol $\circ$ stands for the composition of functions).}. 
Therefore, in this chapter we study the problem 
\begin{equation}
    \label{eq:original-pb}
   \min_{x\in\mathcal{X}}\ f(x) + \psi(x)+H(\mL x),
\end{equation}
where $\mathcal{X}$ is a real Hilbert space, $f$ is a smooth convex function, $\psi,H$ are convex, possibly non-smooth, functions and $\mL$ is a linear operator. To solve Problem~\eqref{eq:original-pb}, we recast it as finding a zero of the sum of three operators which are monotone in a primal--dual product space, under a metric depending on the problem. Then, we apply \algname{Davis--Yin Splitting} (\algname{DYS})~\cite{davis2017three}, a generic method for this type of monotone inclusions. This way, we recover existing algorithms but we also discover a new one, which we call the \algname{Primal--Dual Davis--Yin} (\algname{PDDY}) algorithm. Moreover, this \algname{DYS} representation of our algorithms allow us to use an important inequality regarding \algname{DYS} for their analysis. More precisely, we can apply Lemma~\ref{lem:funda-DYS} below by instanciating the monotone operators and replacing the inner product by the inner product under which we apply \algname{DYS}. Thanks to this, we can afterwards painlessly replace the gradient $\nabla f$ by a stochastic variance-reduced estimator (see Section~\ref{sec:grad-estimators}), which can be much cheaper to evaluate. This machinery allows us to obtain new convergence rates for the stochastic primal--dual algorithms studied in this chapter, and opens the door to a new class of randomized proximal algorithms for large-scale convex non-smooth optimization.

\section{Related Work} 

\textit{Splitting algorithms}: Algorithms allowing to minimize a function involving several non-smooth proximable terms are called splitting algorithms. At the core of splitting algorithms is the \algname{Douglas--Rachford} (or \algname{ADMM}) algorithm~\cite{lions1979splitting,glowinski1975approximation} which is, under reasonable assumptions, the only splitting algorithm that can minimize the sum of two non-smooth functions $\psi+H$~\cite{ryu2020uniqueness}. To minimize $G = \psi +H \circ \mL$, the \algname{Douglas--Rachford} algorithm can be generalized to the \algname{Primal--Dual Hybrid Gradient} (\algname{PDHG}) algorithm, also called \algname{Chambolle--Pock} algorithm~\cite{chambolle2011first}. Behind the success of \algname{PDHG} is the ability to handle such a composite function $G$ and hence the regularizations mentioned above. However, in signal processing and machine learning applications, the objective function usually involves a smooth data fitting term $f$. In order to cover these applications, splitting algorithms like \algname{Condat--V\~u}~\cite{condat2013primal,vu2013splitting} and \algname{PD3O}~\cite{yan2018new} were proposed to solve the Problem~\eqref{eq:original-pb}. These algorithms are primal--dual in nature, i.e., their iterates take the form $(x^k,y^k) \in \cX \times \cY$, where $\cY$ is another real Hilbert space, $x^k$ converges to a solution of Problem~\eqref{eq:original-pb} and $y^k$ converges to a solution of a dual of Problem~\eqref{eq:original-pb}. All Hilbert spaces are supposed of finite dimension.

\textit{\algname{DYS} and monotone operators}: The notion of monotone operator~\cite{bauschke2017convex} generalizes the notion of subdifferential of a convex function. The problem of finding a zero of a monotone operator finds many applications in optimization, beyond the case of a subdifferential. Davis-Yin Splitting~\cite{davis2017three} is a method for finding a zero of the sum of three monotone operators. When the three monotone operators are subdifferentials, \algname{DYS} boils down to an optimization algorithm for solving Problem~\eqref{eq:original-pb} with $\mL = \mI$. When one of the three monotone operators is equal to zero, the \algname{DYS} boils down to the \algname{Forward--Backward Splitting} (\algname{FBS}). If the two monotone operators of the \algname{FBS} are subdifferentials, the \algname{FBS} boils down to the standard proximal gradient algorithm. But the \algname{FBS} goes beyond the case of subdifferentials. For instance, the \algname{Condat--V\~u} algorithms can be seen as instances of the \algname{FBS} involving monotone operators which are not subdifferentials. The \algname{FBS} representation of the \algname{Condat--V\~u} algorithms is the cornerstone of their analysis in~\cite{condat2013primal}.

\textit{Stochastic splitting algorithms}: In machine learning applications, the gradient of $f$ is often intractable and replaced by a cheaper stochastic gradients. These stochastic gradients can be classified in two classes: variance reduced (VR) stochastic gradients~\cite{SVRG,defazio2014saga,defazio2016simple,gorbunov2020unified} and generic stochastic gradients, see, e.g., \cite{moulines2011non,lan2020first}. VR stochastic gradients are stochastic gradient estimators of the full gradient that ensure convergence to an exact solution, as for deterministic algorithms. The variance reduction allows to speedup stochastic algorithms and eventually recover the convergence rates of their deterministic counterparts. In the case where $\mL = \mI$, Problem~\eqref{eq:original-pb} was considered with generic stochastic gradients in~\cite{yurtsever2016stochastic} and with VR stochastic gradients in~\cite{pedregosa2019proximal}. In the general case $\mL \neq \mI$ that is of interest in this chapter, the resolution of \eqref{eq:original-pb} was considered with a generic stochastic gradient in~\cite{zhao2018stochastic}.
\paragraph{Contributions and technical challenges.}
In this chapter we consider the resolution of Problem~\eqref{eq:original-pb} with VR stochastic gradients, and make several contributions w.r.t. the understanding of primal--dual algorithms. 

We first propose a new algorithm called \algname{Primal--Dual Davis--Yin} (\algname{PDDY}) to solve \eqref{eq:original-pb}. This algorithm is obtained as a carefully designed instance of the \algname{DYS} involving operators which are monotone under a metric depending on $\mL$. This \algname{DYS} representation allows us to prove convergence rates for \algname{PDDY}, a task which would be lengthy and technical if such a representation was not obtained prior to proving the convergence rates. More precisely, we analyze \algname{PDDY} with a deterministic gradient and with a variance reduced stochastic gradient. Both settings are cast into a single assumption which can be elegantly plugged into our analysis of \algname{PDDY}, thanks to the flexibility of our framework.   

We believe that representing primal--dual algorithms as "simpler" algorithms using monotone operators in a primal--dual product space is important for understanding the primal--dual algorithms, as illustrated above for \algname{PDDY}. Therefore, our second contribution is to show how the \algname{Condat--V\~u} algorithms and the \algname{PD3O} algorithm can be viewed as instances of the \algname{DYS} involving operators which are monotone under a metric depending on $\mL$\footnote{These monotone operators are not subdifferentials in general.}. Such representation was not known for the \algname{Condat--V\~u} algorithms. Finally, we use again this \algname{DYS} representation to study the \algname{PD3O} algorithm with a VR stochastic gradient. 

One byproduct of our results is the discovery of one of the first linearly converging algorithm for the minimization of a smooth strongly convex function under linear constraints~\cite{mishchenko2019stochastic}. In the particular case where a full gradient is used and $\mL$ is a gossip matrix~\cite{xu2020distributed}, this algorithm leads to a decentralized algorithm whose complexity competes with optimization algorithms designed specifically for the decentralized optimization problem, see~\cite{xu2020distributed}. 

In summary, our contributions are the following:
\begin{itemize}
    \item We propose a new primal--dual algorithm called \algname{PDDY} to solve Problem~\eqref{eq:original-pb}
    \item We propose new generalizations of \algname{PDDY} and \algname{PD3O} using a VR stochastic gradient. These are the first variance reduced algorithms to tackle Problem~\eqref{eq:original-pb}. We leverage the \algname{DYS} representation of \algname{PDDY} and \algname{PD3O} to prove convergence rates which are faster than the convergence rates of non variance reduced algorithms like~\cite{zhao2018stochastic}. 
    \item We propose a new stochastic algorithm for the smooth strongly convex decentralized optimization problem and prove convergence rates for this algorithm. In the particular case of full gradients, its convergence rates is competitive with the decentralized algorithms listed in the Table 1 of~\cite{xu2020distributed}.
    \item We show that the \algname{Condat--V\~u} algorithms are instances of the \algname{DYS} involving monotone operators under a new metric.
\end{itemize}
The choice of studying stochastic versions of \algname{PDDY} and \algname{PD3O} instead of the \algname{Condat--V\~u} algorithms is motivated by the fact that \algname{PDDY} and \algname{PD3O} allow for larger stepsizes. Stochastic versions of the \algname{Condat--V\~u} algorithms could also be analyzed within our framework.

The remainder is organized as follows. In the next section we present the primal--dual formulation of Problem~\eqref{eq:original-pb}, and the associated monotone operators. In Section~\ref{sec3} we recall the \algname{DYS} along with an important inequality regarding its analysis. Then, in Section~\ref{sec:pdalgos} we present the four primal--dual algorithms studied in this chapter and their \algname{DYS} representation. Their stochastic versions and associated convergence rates are presented in Section~\ref{sec:grad-estimators}. The numerical experiments, the convergence proofs and the decentralized application of our algorithms are postponed to the appendix.

\section{Primal--Dual Formulations and Optimality Conditions}\label{sec2}

The necessary notions and notations of convex analysis and operator theory are introduced in \Cref{chapter:intro}. 
Let $\mathcal{X}$ and $\mathcal{Y}$ be finite-dimensional real Hilbert spaces, $\mL\colon \mathcal{X}\rightarrow \mathcal{Y}$ be a linear operator, $f,\psi$ and $H$ be convex, closed and proper. We assume that $f$ is $\nu$-smooth, for some $\nu>0$. We assume, as usual, that there exists  $x^\star\in\mathcal{X}$ such that $0\in \nabla f(x^\star) +\partial \psi(x^\star) +\mL^* \partial H(\mL x^\star)$. Then $x^\star$ is solution to~\eqref{eq:original-pb}. For instance, a standard qualification constraint for this condition to hold is that $0$ belongs to the relative interior of $\dom(H) - \mL \dom(\psi)$~\cite{combettes2012primal}. Therefore, there exists $y^\star \in \mathcal{Y}$ such that $(x^\star,y^\star) \in \mathrm{zer}(M)$, where $M$ is the set-valued operator defined by
\begin{equation}
    \label{eq:M}
    M(x,y) \eqdef \begin{bmatrix}  \nabla f(x)+\partial \psi(x)\!\!\!\!\!\!\!\!\!\!\!\!\!\!&{}+\mL^* y\\-\mL x&\ \ \ \  \ \ \,{}+\partial H^*(y)\end{bmatrix}.
\end{equation}
In other words, there exist
$r^\star \in \partial \psi(x^\star)$ and $h^\star \in \partial H^*(y^\star)$ such that
\begin{equation}
    \label{eq:saddle0}
    \begin{bmatrix} 0 \\ 0\end{bmatrix} = \begin{bmatrix}  \nabla f(x^\star)+ r^\star\!\!\!\!\!&{}+\mL^* y^\star \\ -\mL  x^\star\ \ \ \ \ &\!\!\!\!\!\!{}+h^\star\end{bmatrix}.
\end{equation}
Conversely, for every solution $(x^\star,y^\star) \in \mathrm{zer}(M)$, $x^\star$ is a solution to \eqref{eq:original-pb}.
In the sequel, we let $(x^\star,y^\star) \in \mathrm{zer}(M)$ and $r^\star,h^\star$ be any elements such that Equation~\eqref{eq:saddle0} holds.

The inclusion~\eqref{eq:saddle0} characterizes the first-order optimality conditions associated with the convex--concave Lagrangian function defined as
\begin{equation}
\label{eq:lagrangian}
    \cL(x,y) \eqdef (f+\psi)(x) - H^*(y) + \ps{\mL x,y}.
\end{equation}
For every $x \in \mathcal{X}$, $y \in \mathcal{Y}$, we define the duality gap at $(x,y)$ as $\cL(x,y^\star) - \cL(x^\star,y)$. Then 
\begin{lemma}[Duality gap]
    \label{lem:duality-gap}
    For every $x \in \mathcal{X}$, $y \in \mathcal{Y}$, we have
    \begin{equation}
        \cL(x,y^\star) - \cL(x^\star,y) = D_f(x,x^\star)+D_\psi(x,x^\star)+D_{H^*}(y,y^\star),
    \end{equation}
    where the Bregman divergence of the smooth function $f$ between any two points $x,x$  is $D_f(x,x') \eqdef f(x) - f(x') - \ps{\nabla f(x'),x-x'}$, and
          $D_\psi(x,x^\star) \eqdef \psi(x) - \psi(x^\star) - \ps{r^\star,x-x^\star}$, 
        $D_{H^{*}}(y,y^\star) \eqdef H^*(y) - H^*(y^\star) - \ps{h^\star,y-y^\star}$.
    \end{lemma}
For every $x \in \mathcal{X}$, $y \in \mathcal{Y}$, Lemma~\ref{lem:duality-gap} and the convexity of $f,\psi,H^*$ imply that 
\begin{equation}
  \cL(x^\star,y) \leq \cL(x^\star,y^\star) \leq \cL(x,y^\star).
\end{equation}
 So, the duality gap $\cL(x,y^\star) - \cL(x^\star,y)$ is nonnegative, and it is zero 
if  $x$ is a solution to Problem~\eqref{eq:original-pb} and $y$ is a solution to
the dual problem $\min_{y \in \mathcal{Y}} (f+\psi)^*(-\mL^* y) + H^*(y)$, see Section 15.3 of~\cite{bauschke2017convex}. The converse is true under mild assumptions, for instance strict convexity of the functions around $x^\star$ and $y^\star$.

Finally, one can check that the operator $M$ defined in~\eqref{eq:M} is monotone. Moreover, we have
\begin{align}
    M(x,y) &= \begin{bmatrix} \partial \psi(x) \\ 0 \end{bmatrix} +  \begin{bmatrix} &   \mL^* y \\ -\mL  x\!\!\!\!\!\!& {}+\partial H^*(y)\end{bmatrix}+\begin{bmatrix} \nabla f(x)\\ 0\end{bmatrix} \label{eq:saddle}\\  
    &= \begin{bmatrix} 0 \\ \partial H^*(y) \end{bmatrix} +  \begin{bmatrix} \partial \psi(x)\!\!\!\!\!\!\!\!\!&{}+\mL^* y \\ -\mL  x\end{bmatrix}+\begin{bmatrix} \nabla f(x)\\ 0\end{bmatrix} \label{eq:saddle12}
    \end{align}
and each term at the right hand side of~\eqref{eq:saddle} or~\eqref{eq:saddle12} is maximal monotone, see Corollary 25.5 in~\cite{bauschke2017convex}.

\begin{figure*}[t]
\begin{minipage}{.48\textwidth}\begin{algorithm}[H]
 \caption*{Davis--Yin Splitting \rrbox{\algname{DYS}$(\tilde{A},\tilde{B},\tilde{C})$}~\cite{davis2017three}}
 \begin{algorithmic}[1]
    \State \textbf{Input:} $v^0\in \mathcal{Z}$, $\gamma>0$, number of steps $K$
 	\For{$k = 0,1,2,\dots, K-1$}
  \State $z^{k} = J_{\gamma \tilde{B}}(v^k)$
  \State $u^{k+1} = J_{\gamma \tilde{A}}(2 z^{k} - v^k - \gamma \tilde{C}(z^{k}))$
  \State $v^{k+1} = v^k + u^{k+1} - z^{k}$
 	\EndFor
 \end{algorithmic}
 \end{algorithm}
 \end{minipage}\ \ \ \ \   \begin{minipage}{.48\textwidth}\begin{algorithm}[H]
 \caption*{\rrbox{\algname{PriLiCoSGD}} \\$\big($deterministic version: $g^{k+1}=\nabla f(x^k)\big)$}
 \begin{algorithmic}[1]
    \State \textbf{Input:} $x^0\in \mathcal{X}$, $\gamma>0$, $\tau>0$, number of steps $K$
 	\For{$k = 0,1,2,\dots, K-1$}
	  \State $t^{k+1}=x^k - \gamma g^{k+1}$
	  \State $a^{k+1} = a^k + \tau \mW (t^{k+1}-\gamma a^k)-\tau c$
	  \State $x^{k+1} = t^{k+1}-\gamma a^{k+1}$
  	\EndFor
 \end{algorithmic}
 \end{algorithm}
 \end{minipage}
\begin{minipage}{.48\textwidth}\begin{algorithm}[H]
\caption*{\sobox{\algname{Stochastic PDDY}}\\$\big($deterministic version: $g^{k+1}=\nabla f(x^k)\big)$}
 \begin{algorithmic}[1]
    \State \textbf{Input:} $p^0 \in \cX, y^0 \in \cY$, $\gamma>0$, $\tau>0$, number of steps $K$
 	\For{$k = 0,1,2,\dots, K-1$}
 	\State $y^{k+1} = \prox_{\tau  H^*\!}\big(y^{k}+\tau \mL (p^k-\gamma \mL^*y^{k})\big)$%
 	\State $x^k=p^k -\gamma \mL^* y^{k+1}$
	\State $s^{k+1}=\prox_{\gamma \psi }\big(2x^k-p^k-\gamma g^{k+1}\big)$%
	\State $p^{k+1}=p^k+s^{k+1}-x^k$
 	\EndFor
 \end{algorithmic}
 \end{algorithm}\end{minipage}\ \ \ \ \ \ \begin{minipage}{.48\textwidth}\begin{algorithm}[H]
\caption*{\sobox{\algname{Stochastic PD3O}} \\$\big($deterministic version: $g^{k+1}=\nabla f(x^k)\big)$}
 \begin{algorithmic}[1]
    \State \textbf{Input:} $p^0 \in \cX, y^0 \in \cY $, $\gamma>0$, $\tau>0$, number of steps $K$
 	\For{$k = 0,1,2,\dots, K-1$}
 	\State $x^k=\prox_{\gamma \psi }(p^k)$
 	\State $w^k = 2x^k - p^k-\gamma g^{k+1}$
	\State $y^{k+1} \!= \prox_{\tau H^*}\!\big( y^k + \tau \mL (w^k-\gamma \mL^* y^k)\big)$%
	\State $p^{k+1} = x^k - \gamma g^{k+1} - \gamma \mL^* y^{k+1}$
 	\EndFor
 \end{algorithmic}
\end{algorithm}\end{minipage}
 \end{figure*}

\section{\algname{Davis--Yin Splitting}}\label{sec3}
Solving the optimization problem \eqref{eq:original-pb} boils down to finding a zero $(x^\star,y^\star)$ of the monotone operator $M$ defined in \eqref{eq:M}, which can be written as the sum of three monotone operators, like in \eqref{eq:saddle} or \eqref{eq:saddle12}.
 The \algname{Davis--Yin Splitting} (\algname{DYS}) algorithm~\cite{davis2017three}, is dedicated to this problem; that is, find a zero of the sum of three monotone operators, one of which is cocoercive.

Let $\mathcal{Z}$ be a real Hilbert space. Let $\tilde{A}, \tilde{B}, \tilde{C}$  be maximal monotone operators on $\mathcal{Z}$. We assume that $\tilde{C}$ is $\xi$-cocoercive, for some $\xi >0$.
The \algname{DYS} algorithm, denoted by $\algname{DYS}(\tilde{A},\tilde{B},\tilde{C})$ and shown above, aims at finding an element in  $\mathrm{zer}(\tilde{A}+\tilde{B}+\tilde{C})\neq\emptyset$.
The fixed points of $\algname{DYS}(\tilde{A},\tilde{B},\tilde{C})$ are the triplets $(v^\star, z^\star, u^\star)\in \mathcal{Z}^3$, such that
\begin{equation}
\label{eq:DYSfix}
	z^\star = J_{\gamma \tilde{B}}(v^\star),\quad u^\star = J_{\gamma \tilde{A}}\big(2 z^\star - v^\star - \gamma \tilde{C}(z^\star)\big),\quad u^\star = z^\star.
\end{equation}
These fixed points are related to the zeros of $\tilde{A}+\tilde{B}+\tilde{C}$ as follows, see Lemma 2.2 in~\cite{davis2017three}: for every $(v^\star, z^\star, u^\star)\in \mathcal{Z}^3$ satisfying \eqref{eq:DYSfix}, $z^\star \in \mathrm{zer}(\tilde{A}+\tilde{B}+\tilde{C})$. Conversely, for every $z^\star \in \mathrm{zer}(\tilde{A}+\tilde{B}+\tilde{C})$, there exists $(v^\star, u^\star) \in\mathcal{Z}^2$, such that  $(v^\star, z^\star, u^\star)$ satisfies  \eqref{eq:DYSfix}. We have~\cite{davis2017three}:

\begin{lemma}[Convergence of the \algname{DYS} algorithm]
    \label{lem:DYS-cv}
    Suppose that $\gamma\in (0,2\xi)$. Then the sequences  $(v^k)_{k}$, $(z^k)_{k}$, $(u^k)_{k}$ generated by $\algname{DYS}(\tilde{A},\tilde{B},\tilde{C})$ 
converge 
to some elements $v^\star$, $z^\star$, $u^\star$ in $\mathcal{Z}$, respectively. Moreover, $(v^\star, z^\star, u^\star)$ satisfies  \eqref{eq:DYSfix} and $u^\star = z^\star \in \mathrm{zer}(\tilde{A}+\tilde{B}+\tilde{C})$.
\end{lemma}
The following equality, proved in the Appendix, is at the heart of the convergence proofs: \begin{lemma}[Fundamental equality of the \algname{DYS} algorithm]
\label{lem:funda-DYS}
Let $(v^k,z^k,u^k) \in \cZ^3$ be the iterates of the \algname{DYS} algorithm, and $(v^\star, z^\star, u^\star)\in \cZ^3$ be such that~\eqref{eq:DYSfix} holds.
Then, for every $k \geq 0$, there exist $b^k \in \tilde{B}(z^k), b^\star \in \tilde{B}(z^\star), a^{k+1} \in \tilde{A}(u^{k+1})$ and $a^\star \in \tilde{A}(u^\star)$ such that
\begin{align}
    \label{eq:funda}
    \|v^{k+1} &-v^\star\|^2 = \|v^k - v^\star\|^2 -2\gamma\ps{b^k - b^\star,z^{k} - z^\star}-2\gamma\ps{\tilde{C}(z^{k}) - \tilde{C}(z^\star),z^{k} - z^\star}\\
    &{}-2\gamma\ps{a^{k+1} - a^\star,u^{k+1} - u^\star}-\gamma^2\|a^{k+1}+b^k - \left(a^{\star}+b^\star\right)\|^2+\gamma^2\|\tilde{C}(z^{k}) - \tilde{C}(z^\star)\|^2\notag.
\end{align}
\end{lemma}

\section{Primal--Dual Optimization Algorithms}
\label{sec:pdalgos}

We now set $\mathcal{Z}\eqdef\mathcal{X}\times\mathcal{Y}$, where $\mathcal{X}$ and $\mathcal{Y}$ are the spaces defined in Sect.~\ref{sec2}. To solve the primal--dual problem \eqref{eq:saddle} or \eqref{eq:saddle12}, which consists in finding a zero of the sum $A+B+C$ of 3 operators in $\mathcal{Z}$, of which $C$ is cocoercive, a natural idea is to apply the \algname{Davis--Yin} algorithm  \algname{DYS}$(A,B,C)$. But the resolvent of $\tilde{A}$ or $\tilde{B}$ is often intractable. In this section, we show that preconditioning is the solution; that is,  we exhibit a positive definite linear operator $\mP$, such that \algname{DYS}$(\mP^{-1}A,\mP^{-1}B,\mP^{-1}C)$  is tractable. Since
$\mP^{-1}A,\mP^{-1}B,\mP^{-1}C$ are monotone operators in $\mathcal{Z}_\mP$,
 the algorithm will converge to a  zero of $\mP^{-1}A + \mP^{-1}B + \mP^{-1}C$, or, equivalently, of $A+B+C$.

Let us apply this idea in four different ways.

\subsection{A new primal--dual algorithm: the \algname{PDDY} algorithm}

Let $\gamma >0$ and $\tau>0$ be real parameters. We introduce the four operators on $\mathcal{Z}$, using matrix-vector notations:
  \begin{align}
    A(x,y)\
    & \!=\!\begin{bmatrix}    \mL^* y \\ -\mL  x +\partial H^*(y)\end{bmatrix}\!,\ B(x,y)\!=\!\begin{bmatrix} \partial \psi(x) \\ 0 \end{bmatrix}\!, \notag\\
    C(x,y) &\!=\!\begin{bmatrix} \nabla f(x)\\ 0\end{bmatrix}\!,\ \mP\!=\!\begin{bmatrix} \mI&  0 \\ 0 & \!\frac{\gamma}{\tau}\mI - \gamma^2 \mL \mL^* \end{bmatrix}\!.
    \label{eq:ABCP}
  \end{align}

Note that $\mP$ is positive definite if and only if $\gamma\tau \|\mL\|^2 < 1$. Since $A$, $B$, $C$ are maximal monotone in $\mathcal{Z}$, $\mP ^{-1}A,\mP ^{-1}B,\mP ^{-1}C$ are maximal monotone in $\mathcal{Z}_\mP $. Moreover, $\mP ^{-1}C$ is $1/\nu$-cocoercive in $\mathcal{Z}_\mP $. Importantly, we have:
\begin{align}
    \mP ^{-1}C:(x,y) \mapsto \big(\nabla f(x),0\big),\ \ &\ \  J_{\gamma \mP ^{-1}B}:(x,y)\mapsto \big(\mathrm{prox}_{\gamma \psi}(x),y\big)\label{eq:PB},\\
    J_{\gamma \mP ^{-1}A}:(x,y)\mapsto (x',y'),\ \ &\ \mbox{ where} \label{algorithm_reslv}
    \left\lfloor
       \begin{array}{l}
       y'
   = \mathrm{prox}_{\tau H^*}\big(y+\tau \mL (x - \gamma \mL^* y)\big)\\
   x'=x -\gamma \mL^* y'.
   \end{array}\right.
\end{align}
We plug these explicit steps into the \algname{Davis--Yin} algorithm $(\mP ^{-1}B,\mP ^{-1}A,\mP ^{-1}C)$ and we identify the variables as  $v^k = (p^k,q^k)$, $z^k = (x^k,y^k)$, $u^k = (s^k,h^k)$.
After some simplifications, we obtain the new \algname{Primal--Dual Davis--Yin} (\algname{PDDY}) algorithm, shown above, for Problem~\eqref{eq:original-pb}. Note that it can be written with only one call to $\mL$ and $\mL^*$ per iteration.
Also, the \algname{PDDY} algorithm can be overrelaxed~\cite{con19}, since this possibility exists for the Davis--Yin algorithm. We have:

\begin{theorem}[Convergence of the \algname{PDDY} algorithm]
    \label{th:pddy-cv}
    Suppose that $\gamma\in (0,2/\nu)$ and that $\tau\gamma\|\mL\|^2<1$. Then the sequences $(x^k)_{k}$ and $(s^k)_{k}$ (resp.\ the sequence $(q^k)_{k\in\mathbb{N}}$) generated by the \algname{PDDY} algorithm converge to some 
    solution $x^\star$ to Problem~\eqref{eq:original-pb}
    (resp.\ some $y^\star \in \argmin (f+\psi)^* \circ (-\mL^*) + H^*$).
\end{theorem}
\begin{proof}
    Under the assumptions of Theorem~\ref{th:pddy-cv}, $\mP $ is positive definite. Then the result follows from Lemma~\ref{lem:DYS-cv} applied in $\mathcal{Z}_\mP $ and from the analysis in Sect.~\ref{sec2}.
\end{proof}

Whether the \algname{PDDY} algorithm is a good alternative in practice to the \algname{PD3O} or the \algname{Condat--V\~u} algorithm, which solve the same problems, should be considered on a case-by-case basis. In particular, their memory requirements can be different, and depend on which operations are performed in-place.

\subsection{The \algname{PD3O} algorithm}
We consider the same notations as in the previous section. We switch the roles of $A$ and $B$ and consider \algname{DYS}$(\mP ^{-1}A,\mP ^{-1}B,\mP ^{-1}C)$. Then we recover exactly the \algname{PD3O} algorithm proposed in~\cite{yan2018new}, shown above. Although it is not derived this way, its interpretation as a primal--dual Davis--Yin algorithm is mentioned by its author. Its convergence properties are the same as for the \algname{PDDY} algorithm, as stated in Theorem~\ref{th:pddy-cv}. 

We can note that in a recent work \cite{o2020equivalence}, the \algname{PD3O} algorithm has been shown to be an instance of the Davis--Yin algorithm, with a different reformulation, which does not involve duality. Whether this connection could yield different insights on the \algname{PD3O} algorithm is left for future investigation.

\subsection{The \algname{Condat--V\~u} algorithm}\label{seccv}

    Let $\gamma >0$ and $\tau>0$ be real parameters. We want to study the decomposition \eqref{eq:saddle12} instead of \eqref{eq:saddle}. For this, we define the operators
       \begin{equation}
        \label{eq:ABC2}
        \bar{A}(x,y)\!=\!\begin{bmatrix}\partial \psi(x)\!\!\!\!\!\!\!\!&{}  + \mL^* y \\ -\mL  x&\end{bmatrix}\!,\ \bar{B}(x,y)\!=\!\begin{bmatrix}  0\\ \partial H^*(y)\end{bmatrix}\!, \       C(x,y)\!=\!\begin{bmatrix} \nabla f(x)\\ 0\end{bmatrix}\!,\ \mQ=\begin{bmatrix} \mK&  0 \\ 0 & \mI \end{bmatrix}\!,
      \end{equation}
where $\mK \eqdef \frac{\gamma}{\tau}\mI - \gamma^2 \mL^* \mL$.  If $\gamma\tau \|\mL\|^2 < 1$, $\mK$ and $\mQ$ are positive definite. In that case, since $\bar{A}$, $\bar{B}$, $C$ are maximal monotone in $\mathcal{Z}=\mathcal{X}\times\mathcal{Y}$, $\mQ^{-1}\bar{A}$, $\mQ^{-1}\bar{B}$, $\mQ^{-1}C$ are maximal monotone in $\mathcal{Z}_{\mQ}$.
Moreover, we have:
\begin{align}
    \mQ^{-1}C:(x,y)\mapsto \big(\mK^{-1}\nabla f(x),0\big)\label{eq:QC},\ \ &\ \  J_{\gamma \mQ^{-1}\bar{B}}:(x,y)\mapsto \big(x,\mathrm{prox}_{\gamma H^*}(y)\big),\\
     J_{\gamma \mQ^{-1}\bar{A}}:(x,y)\mapsto (x',y'),\ \ &\ \mbox{ where} \label{algorithm_reslv-cv}
 \left\lfloor
    \begin{array}{l}
    x'
= \mathrm{prox}_{\tau \psi}\big((\mI-\tau\gamma \mL^* \mL) x-\tau \mL^* y \big)\\
y'=y +\gamma \mL x'.
\end{array}\right.
\end{align}
As proved in the Appendix, if we plug these explicit steps into the \algname{Davis--Yin} algorithm \linebreak \algname{DYS}$(\mQ^{-1}\bar{A}, \mQ^{-1}\bar{B}, \mQ^{-1}C)$ or \algname{DYS}$(\mQ^{-1}\bar{B}, \mQ^{-1}\bar{A}, \mQ^{-1}C)$, we recover the two forms of the \algname{Condat--V\~u} algorithm~\cite{condat2013primal,vu2013splitting}; that is, Algorithms 3.1 and  3.2 of~\cite{condat2013primal}, respectively. 
The \algname{Condat--V\~u} algorithm has the form of a  primal--dual forward--backward algorithm~\cite{he2012convergence,combettes2014forward,komodakis2015playing,con19}. But we have just seen that it can be viewed as a primal--dual Davis--Yin algorithm, with a different metric, as well. 
Hence, convergence follows from Lemma~\ref{lem:DYS-cv};
 a technical point is to determine the value of $\xi$, the cocoercivity constant of $\mQ^{-1}C$ in $\mathcal{Z}_{\mQ}$. We prove in the Appendix that we recover the same conditions on $\tau$ and $\gamma$ as in Theorem 3.1 of~\cite{condat2013primal}.

\section{Stochastic Primal--Dual Algorithms: Non-Asymptotic Analysis}\label{sec:grad-estimators}
We now introduce stochastic versions of the \algname{PD3O} and \algname{PDDY} algorithms; we omit the analysis of the stochastic version of the \algname{Condat--V\~u} algorithm, which is the same,  with added technicalities due to cocoercivity with respect to the metric induced by $\mQ$ in \eqref{eq:ABC2}. Moreover, linear convergence results under stronger assumptions are deferred to the 
Appendix. 
Our approach has a `plug-and-play' flavor: we show that we have all the ingredients to leverage the unified theory of stochastic gradient estimators recently presented in \cite{gorbunov2020unified}. 

In the stochastic versions of the algorithms, the gradient $\nabla f(x^k)$ is replaced by a stochastic gradient $g^{k+1}$. More precisely, we consider a filtered probability space $(\Omega,\cF,(\cF_k)_k,\bP)$, an $(\cF_k)_k$-adapted stochastic process $(g^k)_k$, we denote by $\bE$ the mathematical expectation and by $\bE_k$ the conditional expectation w.r.t.\ $\cF_k$. The following assumption is made on the process $(g^k)_{k}$.
\begin{assumption}
    \label{as:sto-grad}
    There exist $\alpha,\beta,\delta \geq 0$, $\rho \in (0,1]$ and a $(\cF_k)_k$-adapted stochastic process denoted by $(\sigma_k)_k$, such that, for every $k \in \mathbb{N}$ 
    \begin{align*}
     &\bE_k[g^{k+1}] = \nabla f(x^k), \\
     &\bE_k[\|g^{k+1} - \nabla f(x^\star)\|^2] \leq 2\alpha D_f(x^k,x^\star) + \beta\sigma_k^2\ ,\\
     & \bE_k[\sigma_{k+1}^2] \leq (1-\rho)\sigma_k^2 + 2\delta D_f(x^k,x^\star).
    \end{align*}
\end{assumption}
Assumption~\ref{as:sto-grad} is a consequence of the smoothness of $f$ and the choice of the stochastic gradient estimator, see~\cite{gorbunov2020unified}. Assumption~\ref{as:sto-grad} is satisfied by several stochastic gradient estimators used in machine learning, including some kinds of coordinate descent~\cite{hanzely2018sega}, variance reduction~\cite{defazio2014saga,hofmann2015variance,gower2018stochastic,kovalev2019don}, and also compressed gradients used to reduce the communication cost in distributed optimization~\cite{horvath2019stochastic}, see Table 1 in~\cite{gorbunov2020unified}. Also, the full gradient estimator defined by
    $g^{k+1} = \nabla f(x^k)$
satisfies Assumption~\ref{as:sto-grad} with $\alpha = \nu$, the smoothness constant of $f$, $\sigma_k \equiv 0$, $\rho = 1$, and $\delta = \beta = 0$, see Theorem 2.1.5 in~\cite{nesterov2018lectures}. The \algname{Loopless SVRG} estimator~\cite{hofmann2015variance,kovalev2019don} also satisfies Assumption~\ref{as:sto-grad}. \begin{proposition}[\algname{Loopless SVRG} estimator]
    Assume that $f$ is written as a finite sum $f = \frac{1}{n}\sum_{i = 1}^n f_i,$ where for every $i \in \{1,\ldots,n\}$, $f_i : \cX \to \bR$ is a $\nu_i$-smooth convex function. Let $p \in (0,1)$, and $(\Omega,\cF,\bP)$ be a probability space. On $(\Omega,\cF,\bP)$, consider:
    
    \noindent$\bullet\ \ $a sequence of i.i.d.\ random variables $(\theta^k)_k$ with Bernoulli distribution of parameter $p$,
    
     \noindent$\bullet\ \ $a sequence of i.i.d.\ random variables $(\xi^k)_k$ with uniform distribution over $\{1,\ldots,n\}$,

    \noindent$\bullet\ \ $the sigma-field $\cF_k$ generated by $(\theta^k,\xi^k)_{0 \leq j \leq k}$ and a $(\cF_k)_k$-adapted stochastic process $(x^k)_k$,

  \noindent$\bullet\ \ $a stochastic process $(\tilde{x}^k)_k$ defined by
        $\tilde{x}^{k+1} = \theta^{k+1} x^k + (1-\theta^{k+1})\tilde{x}^{k}$,
   
  \noindent$\bullet\ \ $a stochastic process $(g^{k})_k$ defined by
        $g^{k+1} = \nabla f(x^k; \xi^{k+1}) - \nabla f(\tilde{x}^{k}; \xi^{k+1}) + \nabla f(\tilde{x}^{k})$.
   
Then, the process $(g^k)_k$ satisfies Assumption~\ref{as:sto-grad} with $\alpha = 2\max_{i \in \{1,\ldots,n\}} \nu_i$, $\beta = 2$, $\rho = p$, $\delta = \alpha p /2$, and 
\[
	\sigma_k^2 = \frac{1}{n}\sum_{i=1}^n \bE_k \left[\|\nabla f_i(\tilde{x}^{k}) - \nabla f_i(x^\star)\|^2 \right].
\]
\end{proposition}
\begin{proof}
    The proof is the same as the proof of Lemma A.11 of~\cite{gorbunov2020unified}. Although this Lemma is only stated for $(x^k)$  generated by a specific algorithm, it remains true for any $(\cF_k)$-adapted stochastic process $(x^k)$.
\end{proof}

We can now exhibit our main results; the details are provided in the Appendix. In a nutshell,  $\mP ^{-1}C(z^k)$ is replaced by the stochastic outcome $\mP ^{-1}(g^{k+1},0)$
    and the last term of Equation~\eqref{eq:funda}, which is nonnegative,      is handled using Assumption~\ref{as:sto-grad}.

\subsection{The \algname{Stochastic PD3O} algorithm} 
We denote by $\|\cdot\|_\mP $ the norm induced by $\mP $ on $\cZ$. The \algname{Stochastic PD3O} algorithm, shown above,  has $\cO\left(\frac{1}{K}\right)$ ergodic convergence in the general case:

\begin{theorem}[Convergence of the \algname{Stochastic PD3O} algorithm]
\label{th:cvx:PD3O}
Suppose that Assumption~\ref{as:sto-grad} holds.
Let $\kappa \eqdef \beta/\rho$, $\gamma, \tau >0$ be such that $\gamma \leq 1/{2(\alpha+\kappa\delta)}$ and $\gamma\tau\|\mL\|^2 < 1$.
Set $V^0 \eqdef \|v^{0} - v^\star\|_\mP ^2 + \gamma^2 \kappa \sigma_{0}^2,$ where $v^0 = (p^0,y^0)$.
Then,
\begin{equation*}
     \mathbb{E}\left[\cL(\bar{x}^{K},y^\star) - \cL(x^\star,\bar{y}^{K+1})\right]
     \leq \frac{V^0}{K \gamma},
\end{equation*}
where $\bar{x}^{K} = \frac{1}{K} \sum_{j = 0}^{K-1} x^j$ and $\bar{y}^{K+1} = \frac{1}{K} \sum_{j = 1}^{K} y^j$.
\end{theorem}

 In the deterministic case $g^{k+1} = \nabla f(x^k)$, we recover the same rate as in  \cite[Theorem 2]{yan2018new}.

\begin{remark}[Primal--Dual gap]
    Deriving a similar bound on the stronger primal--dual gap $(f+\psi+H\circ \mL)(\bar{x}^{k})+((f+\psi)^*\circ (-\mL) +H^*)(\bar{y}^{k})$ 
    requires additional assumptions; for instance, even for the \algname{Chambolle--Pock} algorithm, which is the particular case of the \algname{PD3O}, \algname{PPDY} and \algname{Condat--V\~u} algorithm when $f\equiv 0$, the best available result \cite[Theorem 1]{chambolle2016ergodic} is not stronger than Theorem \ref{th:cvx:PD3O} 
 \end{remark}
\begin{remark}[Particular case of \algname{SGD}]
    In the case where $H =0$ and $\mL = 0$, the \algname{Stochastic PD3O} algorithm boils down to \algname{Proximal Stochastic Gradient Descent} (\algname{Proximal SGD}) and Theorem~\ref{th:cvx:PD3O} implies that
        \[
        		\ec{(f+\psi)(\bar{x}^{K}) - (f+\psi)(x^\star)} \leq \frac{V^0}{K \gamma}.
        \] 
   This $\cO\left(\frac{1}{K}\right)$ ergodic convergence rate 
   unifies known results on \algname{SGD} in the non-strongly-convex case, where the stochastic gradient satisfies Assumption~\ref{as:sto-grad}. This covers 
   coordinate descent and variance-reduced versions, as discussed previously. 
\end{remark}

\subsection{The  algorithm}
We now analyze the proposed \algname{Stochastic PDDY} algorithm, shown above. For it too, we have $\cO\left(\frac{1}{K}\right)$ ergodic convergence in the general case:

    \begin{theorem}[Convergence of the \algname{Stochastic PDDY} algorithm]
        \label{th:cvx:PDDY}
        Suppose that Assumption~\ref{as:sto-grad} holds.
        Let $\kappa \eqdef \beta/\rho$, $\gamma, \tau >0$ be such that $\gamma \leq 1/{2(\alpha+\kappa\delta)}$ and $\gamma\tau\|\mL\|^2 < 1$.
        Define $V^0 \eqdef \|v^{0} - v^\star\|_\mP ^2 + \gamma^2 \kappa \sigma_{0}^2,$ where $v^0 = (p^0,y^0)$.
        Then,
        \begin{equation*}
             \mathbb{E}\left[D_f(\bar{x}^{K},x^\star)+D_{H^*}(\bar{y}^{K+1},y^\star)+ D_\psi(\bar{s}^{K+1},s^\star)\right]
             \leq \frac{V^0}{K \gamma},
        \end{equation*}
        where $\bar{x}^{K} = \frac{1}{K} \sum_{j = 0}^{K-1} x^j$, $\bar{y}^{K+1} = \frac{1}{K} \sum_{j = 1}^{K} y^j$ and $\bar{s}^{K+1} = \frac{1}{K} \sum_{j = 1}^{K} s^j$.
\end{theorem}

\subsection{Linearly constrained or decentralized optimization}

In this section, we set $\psi=0$ and $H\colon y\mapsto (0$ if $y=b$, $+\infty$ else$)$, for some $b \in \mathrm{ran}(\mL)$.
    In this case, Problem~\eqref{eq:original-pb} boils down to $\min_x f(x)$ s.t.\ $\mL x=b$. The stochastic \algname{PD3O} and \algname{PDDY} algorithms both revert to the same algorithm, shown in the Appendix, which 
    we call \algname{Linearly Constrained Stochastic Gradient Descent} (\algname{LiCoSGD}). Note that it is fully split: it does not make use of projections  onto the affine space $\{x \in \cX, \mL x = b\}$ and only makes calls to $\mL$ and $\mL^*$. 
   
   \begin{theorem}[Linear convergence of \algname{LiCoSGD} with $f$ strongly convex]
   \label{th:LV0}
   Suppose that Assumption~\ref{as:sto-grad} holds, that $f$ is $\mu_f$-strongly convex, for some $\mu_f> 0$, and  that $y^0 \in \Range{\mL}$. 
Let $y^\star$ be the unique element of $\Range{\mL}$ such that $\nabla f(x^\star) + \mL^* y^{\star} = 0$, and $\lambda_{\min}^{+}(\mL^* \mL)>0$
 be the smallest positive eigenvalue of $\mL^* \mL$. For every $\kappa > \beta/\rho$ and every $\gamma, \tau >0$ such that $\gamma \leq \frac{1}{\alpha + \kappa\delta}$ and $\gamma\tau\|\mL\|^2 < 1$, we define
   \begin{equation}
       \label{eq:lyapunov-lin}
       V^k 
       \eqdef \|x^{k} - x^\star\|^2+ \left(1+\tau\gamma\lambda_{\min}^{+}(\mL^* \mL)\right)\|y^{k} - y^\star\|_{\gamma,\tau}^2 + \kappa\gamma^2\ec{ \sigma_{k}^2},
   \end{equation}
   and
   \begin{equation}
       \label{eq:rate-lin}
       r 
       \eqdef \max\left(1-\gamma\mu_f,1-\rho+\frac{\beta}{\kappa},\frac{1}{1+\tau\gamma\lambda_{\min}^{+}(\mL^* \mL)}\right)<1.
   \end{equation}
  Then, for every $k\geq 0$,
   \begin{equation}
   	   \bE \left[V^{k}\right] \leq r^k V^0.
   \end{equation}
   \end{theorem}
     
   Furthermore, \algname{LiCoSGD} can be written using $\mW=\mL^* \mL$, $c=\mL^*b$ and primal variables in $\mathcal{X}$ only; this version, called \algname{PriLiCoSGD}, is shown above. Now, consider that $f=\frac{1}{M}\sum_{m=1}^M f_m$ is a finite sum of functions,
   that $\mW$ is a gossip matrix of a network with $M$ nodes~\cite{xu2020distributed}, and that $c=0$. We obtain a new decentralized algorithm, shown and discussed in the Appendix. Theorem~\ref{th:LV0} applies and
   shows that, with the full gradient, 
  $\varepsilon$-accuracy is reached after
   $\cO\left( (\kappa + \kappa_\mW)\log\frac{1}{\varepsilon} \right)$ iterations, where $\kappa$ is the condition number of $f$ and $\kappa_\mW = \frac{\|\mW\|}{\lambda_{\min}^{+}(\mW)}$. 
   This rate is better or equivalent to the one of recently proposed decentralized algorithms, like \algname{EXTRA}, \algname{DIGing}, \algname{NIDS}, \algname{NEXT}, \algname{Harness}, \algname{Exact Diffusion}, see Table 1 of~\cite{xu2020distributed}, \cite[Theorem 1]{li2020revisiting} and \cite{alghunaim2019decentralized}. With a stochastic gradient, the rate of our algorithm is also better than~\cite[Equation 99]{mokhtari2016dsa}.


\chapter{Concluding Remarks}
In this concluding chapter, we summarize the obtained results and outline potential future directions and open problems.
\section{Summary}
In this thesis, we have addressed several issues arising when applying optimization to machine learning problems. Our particular interest was in stochastic first-order methods that scale best on problems where both the dimension and the number of data samples are large.

In \Cref{chapter:local_sgd}, we established several upper bounds for Local SGD under different settings of the noise and data heterogeneity. As a follow-up work of Woodworth et al.~\cite{woodworth2020minibatch} showed, our results are optimal in certain regimes, although there are some remaining gaps. Many other works tried to refine (Woodworth et al., \cite{woodworth2020minibatch}) or build on top of our results. For instance, Karimireddy et al.~\cite{karimireddy2019scaffold} proposed a variance-reduction technique to tackle data heterogeneity, and Malinovskiy et al.~\cite{malinovskiy2020local} designed new methods based on fixed-point iterations.

To improve the algorithms proposed in \Cref{chapter:local_sgd}, in \Cref{chapter:rr} we obtained tight results for Random Reshuffling, and then obtained a local RR algorithm in \Cref{chapter:proxrr}. Moreover, our results in \Cref{chapter:rr} are of independent interest as they close a number of open questions on convergence of incremental methods that date back to the eighties. 

In Chapters~\ref{chapter:sdm} and \ref{chapter:pddy}, we developed new variance-reduction and splitting methods for structured optimization. Our algorithm from both chapters were subsequently shown by Salim et al.~\cite{salim2021optimal} to be optimal in terms of the required number of matrix-vector multiplications when the conditioning of the smooth objective, $\kappa = \frac{L}{\mu}$, and of the constraints, $\kappa_\mA = \frac{\lambda_{\max}(\mA)}{\lambda_{\min}^+(\mA)}$, are of the same order, i.e., $\kappa = \Theta(\kappa_\mA)$. To the best of our knowledge, the results from these two chapters were also the first ones to show that one can obtain linear rates of convergence for arbitrary problems with linear constraints, provided that the objective is smooth and strongly convex. The prior literature provided such guarantees for the case $f(x)=\frac{1}{2}\|x-x^0\|^2$, while we only require $f$ to be upper- and lower-bounded by a quadratic. Our methods also generalize many classic splitting techniques and stochastic algorithms, such as the Kaczmarz method, PDDY, PDHG and \algname{Condat--V\~u} algorithm.

In \Cref{chapter:adaptive}, we addressed a question that is ubiquitous in optimization: stepsize estimation. In particular, we provided the first stepsize rule for gradient descent that converges with $\cO\left(\frac{1}{k}\right)$ rate on any convex problem that is locally smooth. Not only it does not require knowledge of the smoothness constant, it does not need any other parameters either, which is in contrast to all other stepsize rules for gradient descent. Unlike the classical line search strategies, our stepsize does not require any subroutines and can adapt to the local curvature arbitrarily many times.

In \Cref{chapter:diana}, we proposed the first technique for learning the gradients under quantized communication in distributed learning. In particular, we proposed to quantize the differences of gradients instead of gradients themselves and built estimators based on this information. This approach allowed us to improve the communication complexity of distributed first-order method and provide convergence guarantees under larger stepsizes. Moreover, a number of works developed extensions of our algorithm to various settings, see for instance (Horv{\'a}th et al., \cite{horvath2019stochastic}; Liu et al., \cite{liu2020double}; Philippenko and Dieuleveut, \cite{philippenko2020bidirectional}; Liu et al., \cite{li2020acceleration}; Gorbunov et al., \cite{gorbunov2020linearly}; Gorbunov et al., \cite{gorbunov2020unified}).

\section{Future Research Work}
As we discovered answers to some of the important challenges or their aspects, we also faced new challenges and can see new gaps between theory and practice. Below, we briefly provide a few directions which we personally consider to be important and challenging.
\subsection{Federated learning}
Despite all progress, many questions remain open even regarding convergence of Local SGD. First of all, in a work that built on top of our results by Woodworth et al.~\cite{woodworth2020minibatch}, it was pointed out that the current best upper bounds still do not match the available lower bounds. Unfortunately, it remains unknown whether the lower or the upper bounds are to be blamed for this gap. We believe that further improvements upon the bounds may lead to a better understanding and potential development of new methods that fix the technical issues of those works.

Another question that we did not address in the thesis is whether one can use Local SGD to perform meta learning. Since the number of devices in federated learning can be vast, it might be impossible to obtain a single model that performs well on all of those devices. This motivates us to study local methods in the context of learning multiple models with some level of personalization, as was studied, for instance, by Hanzely and Richt{\'a}rik~\cite{hanzely2020federated}. At the moment, this question has received substantially less attention, however, we believe that a proper study may lead to more breakthroughs in the area.

\subsection{Random Reshuffling}
As we mentioned in \Cref{chapter:rr}, our new analysis of Random Reshuffling under strong convexity is optimal since it yields complexity that matches the lower bounds of Safran and Shamir~\cite{Safran2020good} and Rajput et al.~\cite{Rajput2020}. At the same time, neither our analysis nor any other available in the literature provides guarantees showing that Random Reshuffling is better than Gradient Descent in the non-convex regime. Among the remaining questions, we want to particularly emphasize that no deterministic strategy to choose permutation is known to have matching convergence guarantees with Random Reshuffling, not to mention giving better guarantees. After our work was published, it was shown by Rajput et al.~\cite{rajput2021permutation} that in some cases one can achieve a better rate by reversing the previously sampled permutation and using this instead of sampling a completely new one. Unfortunately, no results are yet available for non-quadratic functions showing that this procedure or any other sampling strategy can beat Random Reshuffling. Since Random Reshuffling is widely used in practice, we expect that any improved version could have a huge impact on the applications, and hence we believe it is an important direction to pursue.
    
\subsection{Variance reduction and splitting}
One can mention that in \Cref{chapter:sdm}, we required access to the proximal operator of any term in the summation, and we did not split matrix multiplications. In \Cref{chapter:pddy}, in contrast, we managed to split the linear operator from the proximal operator. However, we required all proximal operators to be used at each iteration. Since both splitting and variance reduction are individually possible in this setting, it is reasonable to expect that one can combine them to obtain a faster method for applications with regularization of linear transformations, such as the PC-Lasso proposed by Tay et al.~\cite{tay2018principal}. Unfortunately, we are not aware of any method that achieves this.
\subsection{Adaptive methods}
We now suggest several promising directions for the extensions of our results.
    First and foremost, a great challenge for us is to
    obtain theoretical guarantees of the proposed method in the
    non-convex setting. Unfortunately, we are not aware of any generic first-order
    method for non-convex optimization that does not rely on the
    descent lemma (or its
    generalization), see, e.g., (Attouch et al., \cite{attouch2013convergence}). We hope to see more work on this side as non-convex problems pose a great challenge and often have more complicated structure than convex ones. 

    We also want to point out that our proposed stepsize requires the objective to be differentiable and smooth, which limits the potential applications of our method. Since the transition from smooth to composite
    minimization (Nesterov, \cite{Nesterov2013a}) in classical first-order methods is rather
    straightforward, one would expect that it is trivial to break this limit for our method too. Unfortunately, the proposed proof of
    \Cref{alg:main} does not seem to provide any route for
    generalization and we hope there is some way of resolving this
    issue.

    Finally, we note that the derived bounds for the stochastic case have a suboptimal dependency on $\kappa$. However, it is not clear to us whether one can extend the techniques from the deterministic analysis to improve the rate. We hope to see more progress on adaptive stochastic algorithms, which are particularly challenging because adaptivity can often break unbiasedness.

\subsection{Communication efficiency}
The theoretical advances in various aspects of communication efficiency that we discussed before have not yet led to a method that would completely eliminate the communication bottleneck. Even though the theory sometimes allows for free reduction of communication complexity, there is a discrepancy between the theoretical improvement and the practical performance of communication-efficient methods that was studied in detail by Dutta et al.~\cite{dutta2020discrepancy}. The primary reason for this discrepancy lies in the difficulty of implementing efficient communication primitives for the known quantization operators. For instance, the binary compression considered in \Cref{chapter:diana} requires Gather operation, which is not as efficient as Reduce used in standard SGD. Moreover, Dutta et al.~\cite{dutta2020discrepancy} point out that vector compression might take longer than communication without compression, so the advantage is lost for compressors with expensive computation. We believe that new compressors are required that are capable of running Reduce operations without any computational overheads. An example of such an approach is the recently proposed IntSGD compressor of Mishchenko et al.~\cite{mishchenko2021intsgd}; we hope to see even finer compressors in the future.

\begin{onehalfspacing}
\renewcommand*\bibname{\centerline{REFERENCES}} 
\addcontentsline{toc}{chapter}{References}
\newcommand{\BIBdecl}{\setlength{\itemsep}{0pt}}

\bibliographystyle{plain}

\hypersetup{linkcolor=magenta}
\bibliography{parts/thesis}

\end{onehalfspacing}

\appendix
\hypersetup{linkcolor=.}
\newpage

\begingroup
\let\clearpage\relax
\let\cleardoublepage\relax 
\begin{center}
\vspace*{2\baselineskip}
{ \textbf{{\large APPENDICES}}} 
\addcontentsline{toc}{chapter}{Appendices} 
\end{center}
\endgroup

\begingroup
\let\clearpage\relax
\chapter{Tables of Frequently Used Notation}\label{sec:table}
\endgroup

\renewcommand{\EE}{\mathbb{E}}

{
\footnotesize

 \begin{longtable}[H]{| p{.20\textwidth} | p{.80\textwidth}| } 
 \caption{Summary of frequently used notation.}\label{tbl:notation}\\
\hline
\multicolumn{2}{|c|}{{\bf \normalsize For all chapters} }\\
\hline
\multicolumn{2}{|c|}{ Basic }\\
\hline
$k$ & Iteration/epoch counter \\
$\E{\cdot}$& Expectation\\
$\EE_\xi\left[\cdot\right]$& Expectation with respect to the variable $\xi$ \\
$\EE_k\left[\cdot\right]$& Expectation conditioned on the information up to iteration $k$ \\
$\langle \cdot ,\cdot \rangle$,  $\| \cdot \|$ & Standard inner product and norm in $\R^d$ $\langle x ,y \rangle =\sum_{i=1}^d x_iy_i$; $\| x \| =\sqrt{\langle x ,x \rangle} $ \\
$\ell_p$ norm & $\|x\|_p= \left(\sum_{l=1}^d x_l^p \right)^{\frac{1}{p}}$, $p\in [1, +\infty]$ \\
$\lambda_{\min} (\cdot), \lambda_{\min}^+(\cdot)$ & The smallest and the smallest positive eigenvalues of a matrix   \\
$\lambda_{\max}(\cdot)$ & The largest eigenvalue of a matrix   \\
$\nabla h(x)$& Gradient of a differentiable function $h$ \\
$\nabla^2 h(x)$& The Hessian of a twice differentiable function $h$  \\
\hline
\multicolumn{2}{|c|}{ Objective }\\
\hline
	$P$ & Full objective if there is a non-smooth term\\
	$d$ & Dimension of the primal space $x\in \RR^d$  \\
	$f\colon \R^d \rightarrow \R$ & Smooth part of the objective \\
    $f_j\colon \R^d \rightarrow \R$ & Differentiable function ($1\leq j \leq n$) if finite-sum setting \\
    $f(\cdot; \xi)\colon \R^d \rightarrow \R$ & Almost-surely differentiable function in expectation setting \\
    $n$ & Number of samples in finite-sum setting \\
    $M$ & Number of nodes in parallel setting \\
	$\psi\colon \R^d \rightarrow \R \cup \{\infty \}$ & Non-smooth part of the objective (proper, closed, convex)\\
	$x^\star$ & Global optimum\\
	$L$ & (Global) smoothness constant of $f$ \\
	$\mu$ & Strong convexity constant  \\
	$\kappa$  & The condition number $\kappa = L/\mu$ for strongly convex objectives \\
  \hline
\multicolumn{2}{|c|}{ Other }\\
\hline
     $[n]$ & Set $\{1,2,\dots,n\}$ \\
     $D_h(x,y)$ & Bregman distance $D_h(x,y) = h(x) - h(y)-\langle \nabla h(y),x-y \rangle$ \\
	$\dom(h)$ & Domain of $h$, $\dom (h) = \{x\;:\; h(x)<+\infty\} $ \\
	$\partial h$ & Subdifferential of $h$, $\partial h(x)= \{s \;:\; h(y) \geq h(x) + \langle s, y-x \rangle\}$, $x\in \dom(h)$ \\
	$\prox_{\gamma h}(x)$ & Proximal operator of a function $h$, $\prox_{\gamma hi}(x) = \argmin_{u} \left\{h(u) + \frac{1}{2\gamma}\|u - x\|^2 \right\}$ \\ 
	$\gamma$ & Fixed stepsize \\
	$\gamma_k$ & Stepsize at iteration $k$ \\
\hline
\end{longtable}


\begin{longtable}{| p{.20\textwidth} | p{.80\textwidth}| } 
 \caption{Notation specific to Chapter~\ref{chapter:local_sgd}.}\label{tab:notation-summary}\\
\hline
\multicolumn{2}{|c|}{{\bf \normalsize Chapter~\ref{chapter:local_sgd}} }\\
\hline
    $g^k_m$ & Stochastic gradient at time $k$ on node $m$  \\ 
    $x^k_m$ & Local iterate at time $k$ on node $m$ \\ 
    $g^k$ & The average of stochastic gradients across nodes at time $k$  \\ 
    $\bar{g}^k$ & Expected value of $g^k$, $\ec{g^k} = \bar{g}^k$ \\ 
    $\hat{x}^k$ & The average of all local iterates at time $k$ \\ 
    $r^{k}$ & The deviation of the average iterate from the optimum $\hat{x}^k - x^\ast$ at time $k$ \\ 
    $\sigma^2$ & Global uniform bound on the variance of the stochastic gradients for identical data   \\ 
    $\sigmaopt^2$ & The variance of the stochastic gradients at the optimum  for identical data \\ 
    $\sigmaf^2$ & The variance of the stochastic gradients at the optimum  for heterogeneous data \\ 
    $k_1, k_2, \ldots, k_p$ & Timesteps at which synchronization happens in \algname{Local SGD} \\ 
    $H$ & Maximum number of local steps between communications, $H=\max_{p} \abs{k_p - k_{p+1}} $   \\ \hline
\end{longtable}


\begin{longtable}{| p{.20\textwidth} | p{.80\textwidth}| } 
 \caption{Notation specific to Chapters~\ref{chapter:rr} and \ref{chapter:proxrr}.}\label{tab:notation}\\
  \hline
  \multicolumn{2}{|c|}{{\bf \normalsize Chapters~\ref{chapter:rr} and~\ref{chapter:proxrr}} }\\
  \hline
  $\pi$ &
    \begin{tabular}[l]{@{}l@{}}A permutation $\pi = \br{ \prm{0}, \prm{1}, \ldots, \prm{n-1} }$ of $\{ 1, 2, \ldots, n \}$ \\ Fixed for \algname{Shuffle-Once} and resampled every epoch for \algname{Random Reshuffling}.\end{tabular} \\
  $x^k_i$              & The current iterate after $i$ steps in epoch $k$, for $0 \leq i \leq n$                                                                       \\ 
  $g^k$                & The sum of gradients used over epoch $k$ such that $x^{k+1}=x^k-\gamma g^k$                                                                                                      \\ 
  $\sigma_{k}^2$       & The variance of the individual loss gradients from the average loss at point $x^k$                                                            \\ 
  $A, B$ & Assumption~\ref{asm:2nd-moment} constants                                                                                                  \\
  $\delta_{k}$         & Functional suboptimality, $\delta_{k} = f(x^k) - f^\ast$, where $f^\ast= \inf_{x} f(x)$ \\ 
  $r_{k}$              & The squared iterate distance from an optimum for convex losses $r_{k} = \sqn{x^k - x^\ast}$       \\ \hline
\end{longtable}


\begin{longtable}{| p{.20\textwidth} | p{.80\textwidth}| } 
 \caption{Notation specific to Chapter~\ref{chapter:adaptive}.}\\
\hline
\multicolumn{2}{|c|}{{\bf \normalsize Chapter~\ref{chapter:adaptive}} }\\
\hline
$L_k $ & Estimated gradient Lipschitz constant \\
$\gamma_k $ & Adaptive stepsize \\
$\mu_k $ & Estimated strong convexity \\
$\theta_k $ & Stepsize ratio, $\theta_k = \frac{\gamma_k}{\gamma_{k-1}}$ \\
$\Theta_k $ & Strong convexity ratio, $\Theta_k = \frac{\mu_k}{\mu_{k-1}}$ \\
$\beta_k $ & Estimated momentum \\
$w_i $ & The weight of $x^i$ in the ergodic iterate $\hat x^K$ \\
$\Psi^k$ & Lyapunov function \\
\hline
\end{longtable}


\begin{longtable}{| p{.20\textwidth} | p{.80\textwidth}| } 
 \caption{Notation specific to Chapter~\ref{chapter:diana}.}\label{tbl:notation-table}\\
\hline
\multicolumn{2}{|c|}{{\bf \normalsize Chapter~\ref{chapter:diana}} }\\
\hline
$\sign(t)$ & The sign of $t$  ($-1$ if $t < 0$, $0$ if $t=0$ and $1$ if $t>1$) \\
$x_{(j)}$ & The $j$-th element of $x\in\R^d$\\
$g_i^k$ & Stochastic gradient at time $k$ on node $i$  \\
$g^k$ & The average of stochastic gradients across nodes at time $k$ \\ 
$\sigma_i^2$ & Variance of the stochastic gradient $g_i^k$  \\
$\sigma^2$ & Variance of the stochastic gradient $g^k$ \\ 
$h_i^k$ & Stochastic approximation of $\nabla f_i(x^\star)$  \\
$\Delta_i^k$ & $\Delta_i^k = g_i^k - h_i^k$ \\
$\text{Quant}_p(\Delta)$ & Full $p$-quantization of vector $\Delta$ \\
$d_l$ & Size of the $l$-th block for quantization \\ 
$m$ & Number of blocks for quantization\\
$\beta$ & Momentum parameter \\
$\hat \Delta_i^k$& Block $p$-quantization of $\Delta_i^k = g_i^k - h_i^k$ \\
$\hat \Delta$ & $\hat \Delta^k = \frac{1}{M}\sumiM\hat\Delta_i^k$  \\ 
$\hat g_i^k$ & Stochastic approximation of $\nabla f_i(x^k)$ \\
$\hat g^k$ & $g^k = \frac{1}{M}\sumiM g_i^k$ \\
$v^k$ & Stochastic gradient with momentum $v^k = \beta v^{k-1} + \hat g^k$\\
$\Psi_l(x)$ & Variance of the $l$-th block $\Psi_l(x) = \|x(l)\|_1\|x(l)\|_p - \|x(l)\|^2$ \\
$\Psi(x)$ & Variance of the quantized vector $\Psi(x) = \sum\limits_{l=1}^m\Psi_l(x)$ \\
$\widetilde d$ & $\widetilde d = \max\limits_{l=1,\ldots,m}d_l$ \\
$V^k$ & Lyapunov function $V^k = \|x^k - x^\star\|^2 + \frac{c\gamma^2}{M}\sumiM\|h_i^k - h_i^\star\|$  \\
$\zeta$ & Data dissimilarity, $\frac{1}{M}\sum_{i=1}^M\|\nabla f_i(x) - \nabla f(x)\|^2 \le \zeta^2$ \\
$\mathbb{E}_{Q^k} \left[ \cdot \right]$ & Expectation with respect to quantization  \\
\hline
\end{longtable}


\begin{longtable}{| p{.20\textwidth} | p{.80\textwidth}| } 
 \caption{Notation specific to Chapter~\ref{chapter:sdm}.}\label{tbl:notation_sdm}\\
\hline
\multicolumn{2}{|c|}{{\bf \normalsize Chapter~\ref{chapter:sdm}} }\\
\hline
$\cX^\star $ & Set of optimal solutions $\cX^\star = \{x^\star\in \RR^d \;:\; P(x) \geq P(x^\star) \; \forall x\in \RR^d\}$ \\
$\sigma_{\star}^2$ & Gradient noise at optimum $\sigma_{\star}^2= \EE_\xi \left[\|\nabla f(x^\star; \xi) - \nabla f(x^\star)\|^2 \right]$ \\
$g$ & $g(x) = \frac{1}{m}\sum \limits_{j=1}^m g_j(x)$  (proper, closed, convex)\\
$g_j$ & $g_j\colon\RR^d\to \RR\cup \{+\infty\}$ (proper, closed, convex)\\
$\phi_j$ & Structured form of $g_j$ $g_j(x) = \phi_j(\mA_j^\top x)$, $\mA_j\in \RR^{d\times d_j}$ \\
$v^k$ & Estimator of $\nabla f(x^k)$ \\
$y_1^k,\dots,y_j^k $ & Dual variables \\
$p_j$ & Probability of selecting index $j$ \\
$\gamma$ & Stepsize associated with $f$ and $\psi$ \\
$\eta_j$  & Stepsize associated with $g_j$,  $\eta_j = \frac{\gamma}{m p_j}$ \\
$\cL^k$ & Lyapunov function   $\cL^k
        =  \ec{\|x^k - x^\star\|^2} + \cM^k + \cY^k$ \\
$\cY^k$ & Dual distance to solution $ \cY^k
	= (1+\nu)\sumlm \eta_l^2\ec{\|y_l^{k} - y_l^\star\|^2}$\\      
$L_j$ & Smoothness constant of $g_j$, $L_j\in\RR\cup \{+\infty\}$ ($L_j=+\infty$ if $g_j$ is non-smooth)\\
$\nu$	& $\nu = \min_{j=1,\dots,m} \frac{1}{\eta_j L_j}$ ($\nu=0$ if any $g_j$ is non-smooth)  \\
$\ind_{\cC}(x)$ & Characteristic function of a set $\cC$, $\ind_{\cC}(x)= \begin{cases}0 & x\in \cC \\ +\infty & x\notin \cC\end{cases}$ \\
$\Pi_{\cC}(x)$ & Projection onto a set $\cC$, $\Pi_{\cC}(x) =  \prox_{\ind_{\cC}}(x)= \argmin_{u\in \cC} \|u-x\|$ \\
\hline
\end{longtable}


\begin{longtable}{| p{.20\textwidth} | p{.80\textwidth}| } 
 \caption{Notation specific to Chapter~\ref{chapter:pddy}.}\\
\hline
\multicolumn{2}{|c|}{{\bf \normalsize Chapter~\ref{chapter:pddy}} }\\
\hline
$\mL $ & Linear operator from the objective \\
$\cX $ & Primal space \\
$\cY $ & Dual space \\
$\cL\colon \cX\times \cY\to \R $ & Lagrangian function \\
$y^\star$ & Dual solution \\
$\tau$ & Dual stepsize \\
$\mP$ & Linear operator that defines the primal-dual metric \\
$(\Omega,\cF,(\cF_k)_k,\bP)$ & Probability space \\
$V^k$ & Lyapunov function \\
\hline
\end{longtable}

}

\renewcommand{\EE}{\mathbb{E}}

\chapter{Appendix for Chapter~\ref{chapter:local_sgd}}
\label{local_sgd_appendix}

\graphicspath{{local_sgd/}}

\section{Proofs for Identical Data under Assumption~\ref{asm:uniformly-bounded-variance}}
\subsection{Proof of Lemma~\ref{lemma:uniform-var-iterate-variance-bound}}
\begin{proof}
     Let $k \in \N$ be such that $k_p \leq k \leq k_{p+1} - 1$. Recall that for a time $k$ such that $k_p \leq k < k_{p+1}$ we have $x^{k+1}_m = x^k_m - \gamma g^{k}_m$ and $\hat{x}^{k+1} = \hat{x}^k - \gamma g^k$. Hence for the expectation conditional on $x^k_1, x^k_2, \ldots, x^k_M$ we have:
     \begin{align*}
         \ecn{x^{k+1}_m - \hat{x}^{k+1}} &= \sqn{x^k_m - \hat{x}^k} + \gamma^2 \ecn{g^k_m - g^k} - 2 \gamma \ec{\ev{x^k_m - \hat{x}^k, g^k_m - g^k}} \\
         &= \sqn{x^k_m - \hat{x}^k} + \gamma^2 \ecn{g^k_m - g^k} - 2 \gamma \ev{x^k_m - \hat{x}^k, \nabla f(x^{k}_m)} \\
         &\quad + 2 \gamma \ev{x^k_m - \hat{x}^k, \bar{g}^k}.
     \end{align*}
     Averaging both sides and letting $V^k = \frac{1}{M} \sum_{m} \sqn{x^{k}_m - \hat{x}^k}$, we have
     \begin{align}
         \ec{V^{k+1}} &= V^k + \frac{\gamma^2}{M} \sum_{m} \ecn{g^k_m - g^k} - \frac{2 \gamma}{M} \sum_{m} \ev{x^k_m - \hat{x}^k, \nabla f(x^k_m)} + 2 \gamma \underbrace{\ev{ \hat{x}^k - \hat{x}^k, \bar{g}^k}}_{=0} \nonumber \\
         \label{iterate-variance-bound-recursion}
         &= V^k + \frac{\gamma^2}{M} \sum_{m} \ecn{g^k_m - g^k} - \frac{2 \gamma}{M} \sum_{m} \ev{x^k_m - \hat{x}^k, \nabla f(x^k_m)}.
     \end{align}
     Now note that by expanding the square we have,
     \begin{align}
         \label{iterate-gradient-variance-bound-1}
         \ecn{g^k_m - g^k} 
         &= \ecn{g^k_m - \bar{g}^k} + \ecn{\bar{g}^k - g^k} + 2 \ec{\ev{g^k_m - \bar{g}^k, \bar{g}^k - g^k}}.
     \end{align}
     We decompose the first term in the last equality again by expanding the square,
     \begin{align*}
         \ecn{g^k_m - \bar{g}^k} &= \ecn{g^k_m - \bar{g}^k_m} + \sqn{\bar{g}^k_m - \bar{g}^k} + 2 \ec{\ev{ g^k_m - \bar{g}^k_m, \bar{g}^k_m - \bar{g}^k }} \\
         &= \ecn{g^k_m - \bar{g}^k_m} + \sqn{\bar{g}^k_m - \bar{g}^k} + 2 \underbrace{\ev{ \bar{g}^k_m - \bar{g}^k_m, \bar{g}^k_m - \bar{g}^k }}_{=0} \\
         &= \ecn{g^k_m - \bar{g}^k_m} + \sqn{\bar{g}^k_m - \bar{g}^k}.
     \end{align*}
     Plugging this into~\eqref{iterate-gradient-variance-bound-1} we have,
     \begin{align*}
         \ecn{g^k_m - g^k} &= \ecn{g^k_m - \bar{g}^k_m} + \sqn{\bar{g}^k_m - \bar{g}^k} + \ecn{\bar{g}^k - g^k} + 2 \ec{\ev{g^k_m - \bar{g}^k, \bar{g}^k - g^k}}.
     \end{align*}
     Now average over $m$:
     \begin{align*}
         \frac{1}{M} \sum_{m} \ecn{g^k_m - g^k} &= \frac{1}{M} \sum_{m} \ecn{g^k_m - \bar{g}^k_m} + \frac{1}{M} \sum_m \sqn{\bar{g}^k_m - \bar{g}^k} + \ecn{\bar{g}^k - g^k} \\
         &- 2 \ecn{\bar{g}^k - g^k},
     \end{align*}
     where we used that by definition $\avemm g^k_m = g^k$. Hence,
     \begin{align}
         \frac{1}{M} \sum_m \ecn{g^k_m - g^k} &= \frac{1}{M} \sum_{m} \ecn{g^k_m - \bar{g}^k_m} + \frac{1}{M} \sum_m \sqn{\bar{g}^k_m - \bar{g}^k} - \ecn{\bar{g}^k - g^k} \nonumber \\
         \label{iterate-gradient-variance-bound-2}
         &\leq \frac{1}{M} \sum_{m} \ecn{g^k_m - \bar{g}^k_m} + \frac{1}{M} \sum_{m} \sqn{\bar{g}^k_m - \bar{g}^k}.
     \end{align}
     Now note that for the first term in~\eqref{iterate-gradient-variance-bound-2} we have by Assumption~\ref{asm:uniformly-bounded-variance},
     \begin{align}
         \label{iterate-gradient-variance-bound-3}
         \ecn{g^k_m - \bar{g}^k_m} &= \ecn{g^k_m - \nabla f(x^k_m)} \leq \sigma^2.
     \end{align}
     For the second term in \eqref{iterate-gradient-variance-bound-2} we have
     \begin{align*}
         \sqn{\bar{g}^k_m - \bar{g}^k} &= \sqn{\bar{g}^k_m - \nabla f(\hat{x}^k)} + \sqn{ \nabla f(\hat{x}^k) - \bar{g}^k } + 2 \ev{ \bar{g}^k_m - \nabla f(\hat{x}^k), \nabla f(\hat{x}^k) - \bar{g}^k  }.
     \end{align*}
     Averaging over $m$,
     \begin{align*}
         &\frac{1}{M} \sum_{m=1}^{M} \sqn{ \bar{g}^k_m - \bar{g}^k }\\
         &= \frac{1}{M} \sum_{m} \sqn{\bar{g}^k_m - \nabla f(\hat{x}^k)} + \sqn{\nabla f(\hat{x}^k) - \bar{g}^k} + 2 \ev{ \bar{g}^k - \nabla f(\hat{x}^k), \nabla f(\hat{x}^k) - \bar{g}^k } \\
         &= \frac{1}{M} \sum_{m} \sqn{\bar{g}^k_m - \nabla f(\hat{x}^k)} + \sqn{\nabla f(\hat{x}^k) - \bar{g}^k} - 2 \sqn{\nabla f(\hat{x}^k) - \bar{g}^k} \\
         &= \frac{1}{M} \sum_{m} \sqn{\bar{g}^k_m - \nabla f(\hat{x}^k)} - \sqn{\nabla f(\hat{x}^k) - \bar{g}^k} \leq \frac{1}{M} \sum_{m} \sqn{\bar{g}^k_m - \nabla f(\hat{x}^k)},
     \end{align*}
     where we used the fact that $\frac{1}{M} \sum_{m} \bar{g}^k_m = \bar{g}^k$, which comes from the linearity of expectation. Now we bound $\sqn{\bar{g}^k_m - \nabla f(\hat{x}^k)}$ in the last inequality by smoothness and then use that Jensen's inequality implies $\sum_{m=1}^{M} (f(\hat{x}^k) - f(x^k_m)) \leq 0$,
     \begin{align}
         \frac{1}{M} \sum_{m} \sqn{ \bar{g}^k_m - \nabla f(\hat{x}^k) } &= \frac{1}{M} \sum_{m} \sqn{ \nabla f(x^{k}_m) - \nabla f(\hat{x}^k) } \nonumber \\
         &\overset{\eqref{eq:grad_dif_bregman}}{\le} \frac{1}{M} \sum_{m} 2L (f(\hat x^k) - f(x^k_m) - \ev{\hat x^k - x^k_m, \nabla f(x^k_m)}) \nonumber \\
         \label{eq:iterate-gradient-variance-bound-4}
         &\le \frac{2L}{M} \sum_{m} \ev{x^k_m - \hat x^k, \nabla f(x^k_m)}.
     \end{align}
     Plugging in \eqref{eq:iterate-gradient-variance-bound-4} and \eqref{iterate-gradient-variance-bound-3} into \eqref{iterate-gradient-variance-bound-2} we have,
     \begin{align}
         \label{iterate-gradient-variance-bound-5}
         \frac{1}{M} \sum_{m} \ecn{g^k_m - g^k} \leq \sigma^2 + \frac{2L}{M} \sum_{m} \ev{x^k_m - \hat x^k, \nabla f(x^k_m)}.
     \end{align}
     Plugging \eqref{iterate-gradient-variance-bound-5} into \eqref{iterate-variance-bound-recursion}, we get
     \begin{align}
         \ec{V^{k+1}} 
         &\le V^k + \gamma^2 \sigma^2 - \frac{2 \gamma}{M} \sum_m \ev{x^k_m - \hat{x}^k, \nabla f(x^k_m)} + \frac{2L \gamma^2}{M} \sum_{m} \ev{x^k_m - \hat x^k, \nabla f(x^k_m)} \nonumber \\
         \label{iterate-variance-recursion-pre-sc}
         &= V^k + \gamma^2 \sigma^2 - \frac{2 \gamma(1 - \gamma L)}{M} \sum_m \ev{x^k_m - \hat{x}^k, \nabla f(x^k_m)} \\
         &\overset{\eqref{eq:asm-strong-convexity}}{\le} \br{1 - \gamma(1 - \gamma L)\mu} V^k + \gamma^2\sigma^2. \nonumber
     \end{align}
     Using that $\gamma \leq \frac{1}{2L}$ we can conclude,
     \begin{align*}
         \ec{V^{k+1}}
         &\le \br{1 - \frac{\gamma\mu}{2}} V^k + \gamma^2\sigma^2 \\
         &\le V^k + \gamma^2 \sigma^2.
     \end{align*}
     Taking expectations and iterating the above inequality,
     \begin{align*}
         \ec{V^k} &\leq \ec{V^{k_p}} + \gamma^2 \sigma^2 \br{k - k_p} \\
         &\leq \ec{V^{k_p}} + \gamma^2 \sigma^2 \br{k_{p+1} - k_p - 1} \\
         &\leq \ec{V^{k_p}} + \gamma^2 \sigma^2 \br{H - 1}.
     \end{align*}
     It remains to notice that by assumption we have $V^{k_p} = 0$.
\end{proof}

\subsection{Two more lemmas}
\begin{lemma}{\cite{Stich2018}.}
     \label{lemma:iterate-one-recursion}
     Let $(x^k_m)_{k}$ be the iterates generated by Algorithm~\ref{alg:local_sgd} run with identical data. Suppose that $f$ satisfies Assumption~\ref{asm:convexity-and-smoothness} and that $\gamma \leq \frac{1}{2 L}$. Then,
     \begin{align}
         \begin{split}
             \ecn{\hat{x}^{k+1} - x^\ast} \leq (1 &- \gamma \mu) \ecn{\hat{x}^k - x^\ast} + \gamma^2 \ecn{g^k - \bar{g}^k} \\
             &- \frac{\gamma}{2} \ec{D_{f} (\hat{x}^k, x^\ast)} + 2 \gamma L \ec{V^k}.
         \end{split}
     \end{align}
 \end{lemma}
 \begin{proof}
     This is Lemma 3.1 in \cite{Stich2018}.
 \end{proof}
 
 \begin{lemma}
     \label{lemma:minibatch-variance-reduction}
     Suppose that Assumption~\ref{asm:uniformly-bounded-variance} holds. Then,
     \[ \ecn{g^k - \bar{g}^k} \leq \frac{\sigma^2}{M}. \]
 \end{lemma}
 \begin{proof}
     This is Lemma 3.2 in \cite{Stich2018}. Because the stochastic gradients $g^k_m$ are independent we have that the variance of their sum is the sum of their variances, hence
     \begin{align*}
         \ecn{g^k - \bar{g}^k} &= \frac{1}{M^2} \ecn{ \sum_{m=1}^{M} g^k_m - \bar{g}^k_m } = \frac{1}{M^2} \sum_{m=1}^{M} \ecn{g^k_m - \bar{g}^k_m} \leq \frac{\sigma^2}{M}.
     \end{align*}
 \end{proof}

\subsection{Proof of Theorem~\ref{thm:sc-convergence-theorem}}
\begin{proof}
     Combining Lemma~\ref{lemma:iterate-one-recursion} and Lemma~\ref{lemma:minibatch-variance-reduction}, we have
     \begin{equation}
         \label{sc-thm-proof-1}
         \ecn{\hat{x}^{k+1} - x^\ast} \leq (1 - \gamma \mu) \ecn{\hat{x}^k - x^\ast} + \frac{\gamma^2 \sigma^2}{M} - \frac{\gamma}{2} \ec{D_{f} (\hat{x}^k, x^\ast)} + 2 \gamma L \ec{V^k}.
     \end{equation}
     Using Lemma~\ref{lemma:uniform-var-iterate-variance-bound} we can upper bound the $\ec{V^k}$ term in $\eqref{sc-thm-proof-1}$:
     \begin{equation*}
         \ecn{\hat{x}^{k+1} - x^\ast} \leq (1 - \gamma \mu) \ecn{\hat{x}^k - x^\ast} + \frac{\gamma^2 \sigma^2}{M} - \frac{\gamma}{2} \ec{D_{f} (\hat{x}^k, x^\ast)} + 2 \gamma^3 L \br{H - 1} \sigma^2.
     \end{equation*}
     Letting $r^{k+1} = \hat{x}^{k+1} - x^\ast$, we have
     \begin{equation*}
         \ecn{r^{k+1}} \leq \br{1 - \gamma \mu} \ecn{r^k} + \frac{\gamma^2 \sigma^2}{M} + 2 \gamma^3 L \br{H - 1} \sigma^2.
     \end{equation*}
     Recursing the above inequality we have,
     \begin{equation*}
         \ecn{r^{K}} \leq \br{1 - \gamma \mu}^{K} \ecn{x^0-x^\ast} + \br{ \sum_{k=0}^{K-1} \br{1 - \gamma \mu}^k } \br{ \frac{\gamma^2 \sigma^2}{M} + 2 \gamma^3 L \br{H - 1} \sigma^2 }.
     \end{equation*}
     Using that $\sum_{k=0}^{K-1} \br{1 - \gamma \mu}^{k} \leq \sum_{k=0}^{\infty} \br{1 - \gamma \mu}^{k} = \frac{1}{\gamma \mu}$ we have,
     \begin{equation*}
         \ecn{r^{K}} \leq \br{1 - \gamma \mu}^{K} \ecn{x^0-x^\ast} + \frac{\gamma \sigma^2}{\mu M} + \frac{2 \gamma^2 L \br{H - 1} \sigma^2}{\mu},
     \end{equation*}
     which is the claim of this theorem.
 \end{proof}

\subsection{Proof of Theorem~\ref{thm:weakly-convex-thm}}
\begin{proof}
    Let $r^k = \hat{x}^k - x^\ast$, then putting $\mu = 0$ in Lemma~\ref{lemma:iterate-one-recursion} and combining it with Lemma~\ref{lemma:minibatch-variance-reduction}, we have
    \begin{equation*}
        \ecn{r^{k+1}} \leq \ecn{r^k} + \frac{\gamma^2 \sigma^2}{M} - \frac{\gamma}{2} \ec{D_{f} (\hat{x}^k, x^\ast)} + 2 \gamma L \ec{V^k}.
    \end{equation*}
    Further using Lemma~\ref{lemma:uniform-var-iterate-variance-bound},
    \begin{equation*}
        \ecn{r^{k+1}} \leq \ecn{r^k} + \frac{\gamma^2 \sigma^2}{M} - \frac{\gamma}{2} \ec{D_{f} (\hat{x}^k, x^\ast)} + 2 \gamma^3 L \br{H - 1} \sigma^2.
    \end{equation*}
    Rearranging we have,
    \begin{align*}
        \frac{\gamma}{2} \ec{D_{f} (\hat{x}^k, x^\ast)} \leq \ecn{r^k} - \ecn{r^{k+1}} + \frac{\gamma^2 \sigma^2}{M} + 2 \gamma^3 L \br{H - 1} \sigma^2.
    \end{align*}
    Averaging the above equation as $k$ varies between $0$ and $K-1$,
    \begin{align}
        \frac{\gamma}{2 K} \sum_{k=0}^{K-1} \ec{D_{f} (\hat{x}^k, x^\ast)} 
        &\leq \frac{1}{K} \sum_{k=0}^{K-1} \left(\ecn{r^k} - \ecn{r^{k+1}}\right) \nonumber\\
        &\qquad  + \frac{1}{K} \sum_{k=0}^{K-1} \br{ \frac{\gamma^2 \sigma^2}{M} + 2 \gamma^3 L \br{H - 1} \sigma^2} \nonumber \\
        &= \frac{\sqn{x^0-x^\ast} - \ecn{r^K}}{K} + \frac{\gamma^2 \sigma^2}{M} + 2 \gamma^3 L \br{H - 1} \sigma^2 \nonumber \\
        \label{wc-thm-proof-2}
        &\leq \frac{\sqn{x^0-x^\ast}}{K} + \frac{\gamma^2 \sigma^2}{M} + 2 \gamma^3 L \br{H - 1} \sigma^2.
    \end{align}
    By Jensen's inequality we have $D_{f} (\bar{x}^K, x^\ast) \leq \frac{1}{K} \sum_{k=0}^{K-1} D_{f} (\hat{x}^k, x^\ast)$. Using this in \eqref{wc-thm-proof-2} we have,
    \begin{align*}
        \frac{\gamma}{2} \ec{D_{f} (\bar{x}^K, x^\ast)} \leq \frac{\sqn{x^0-x^\ast}}{K} + \frac{\gamma^2 \sigma^2}{M} + 2 \gamma^3 L \br{H - 1} \sigma^2.
    \end{align*}
    Dividing both sides by $\gamma/2$ yields the theorem's claim.
\end{proof}

\section{Proofs for Identical Data under Assumption~\ref{asm:finite-sum-stochastic-gradients}}
\subsection{Preliminary lemmas}
\begin{lemma}[Individual gradient variance bound]
     \label{lemma:gradient-smoothness-bound}
     Assume that Assumption~\ref{asm:finite-sum-stochastic-gradients} holds with identical data, then for all $k \geq 0$ and  $m \in [M]$ we have
     \begin{equation}
         \label{eq:lma-gradient-smoothness-bound}
         \ecn{g^k_m} \leq 4 L D_{f} (x^k_m, x^\ast) + 2 \sigma_{m}^2,
     \end{equation}
    where $\sigma_m^2 \eqdef \ecn[\xi_m \sim \D_m]{\nabla f (x^\ast; \xi_m)}$ is the noise at the optimum on the $m$-th node.
\end{lemma}
\begin{proof}
     Using that $g^k_m = \nabla f (x^k_m; \xi_m)$ for some $\xi_m \sim \D_m$, 
     \begin{align*}
         \sqn{g^k_m} &= \sqn{\nabla f (x^k_m; \xi_m)} \\
         &\overset{\eqref{eq:sum_sqnorm}}{\leq} 2 \sqn{\nabla f (x^k_m; \xi_m) - \nabla f (x^\ast; \xi_m)} + 2 \sqn{\nabla f (x^\ast; \xi_m)} \\
         &\overset{\eqref{eq:grad_dif_bregman}}{\leq} 4 L \br{f (x^k_m; \xi_m) - f (x^\ast; \xi_m) - \ev{\nabla f (x^\ast; \xi_m), x^k_m - x^\ast}} + 2 \sqn{\nabla f (x^\ast; \xi_m)}.
     \end{align*}
     Taking expectations and using that $\ec{\nabla f (x^\ast; \xi)} = \nabla f(x^\ast) = 0$ we get,
     \begin{align*}
         \ecn{g^k_m} &\leq 4 L \br{f(x^k_m) - f(x^\ast)} + 2 \sigma_m^2 \\
         &= 4 L D_{f} (x^k_m, x^\ast) + 2 \sigma_m^2.
     \end{align*}
\end{proof}

\begin{lemma}[Average gradient variance reduction]
    \label{lemma:minibatch-gradient-noise-reduction}
    Assume that Assumption~\ref{asm:finite-sum-stochastic-gradients} holds with identical data, then for all $k \geq 0$ and for $M$ nodes we have,
    \begin{align}
        \label{eq:lma-mbns}
        \ecn{g^k - \bar{g}^k} \leq \frac{2 \sigmaopt^2}{M} + \frac{4 L}{M^2} \sum_{m=1}^{M} D_{f} (x^k_m, x^\ast).
    \end{align}
\end{lemma}
\begin{proof}
    Using the definition of $g^k$ and $\bar{g}^k$,
    \begin{align}
        \ecn{g^k - \bar{g}^k} &= \ecn{\frac{1}{M} \sum_{m=1}^{M} g^k_m - \nabla f(x^k_m)} \nonumber \\
        &= \frac{1}{M^2} \ecn{\sum_{m=1}^{M} \br{g^k_m - \nabla f(x^k_m)}}. \label{eq:lma-mgns-proof-1}
    \end{align}
    The sum in \eqref{eq:lma-mgns-proof-1} is the variance of a sum of independent random variables and hence can be decomposed into the sum of their individual variances which we can use Lemma~\ref{lemma:gradient-smoothness-bound} to bound:
    \begin{align*}
        \ecn{g^k - \bar{g}^k} &= \frac{1}{M^2} \sum_{m=1}^{M} \ecn{g^k_m - \nabla f(x^k_m)} \\
        &\overset{\eqref{eq:variance_sqnorm_upperbound}}{\leq} \frac{1}{M^2} \sum_{m=1}^{M} \ecn{g^k_m} \\
        &\overset{\eqref{eq:lma-gradient-smoothness-bound}}{\leq} \frac{1}{M^2} \sum_{m=1}^{M} \br{ 2 \sigma_m^2 + 4 L D_{f} (x^k_m, x^\ast) } \\
        &= \frac{2 \sigmaopt^2}{M} + \frac{4 L}{M^2} \sum_{m=1}^{M} D_f (x^k_m, x^\ast),
    \end{align*}
    where in the last equality we used that $\sigmaopt^2$ is by definition equal to $\sum_{m=1}^{M} \sigma_m^2/M$.
\end{proof}
The next Lemma bounds the optimality gap across one iteration when the descent step is $\hat{x}^{k+1} = \hat{x}^k - \frac{\gamma}{M} \sum_{m=1}^{M} \nabla f(x^k_m)$, i.e., when the expectation of the \algname{Local SGD} update is used.
\begin{lemma}[Perturbed iterate analysis]
    \label{lemma:perturbed-iterate-analysis}
     Suppose that Assumptions~\ref{asm:convexity-and-smoothness}, and \ref{asm:finite-sum-stochastic-gradients} hold with identical data. Then,
    \begin{align}
        \begin{split}
            \label{eq:lma-perturbed-iterate-analysis}
            \sqn{\hat{x}^{k+1} - \gamma \bar{g}^k - x^\ast} &\leq \sqn{\hat{x}^k - x^\ast} + 2 \gamma L V^k \\
            &\quad + \frac{2 \gamma}{M} \sum_{m=1}^{M} \br{\br{\gamma L - \frac{1}{2}} \br{f(x^k_m) - f(x^\ast)} - \frac{\mu}{2} \sqn{x^k_m - x^\ast}}.
        \end{split}
    \end{align}
\end{lemma}
\begin{proof}
    This is the first part of Lemma~3.1 in~\cite{Stich2018} and we reproduce it for completeness:
    \begin{align}
        &\sqn{\hat{x}^k - x^\ast - \gamma \bar{g}^k} \nonumber\\
        & = \sqn{\hat{x}^k - x^\ast} + \gamma^2 \sqn{\bar{g}^k} - 2 \gamma \ev{\hat{x}^k - x^\ast, \bar{g}^k} \nonumber \\
        &= \sqn{\hat{x}^k - x^\ast} + \gamma^2 \sqn{\bar{g}^k} - \frac{2 \gamma}{M} \sum_{m=1}^{M} \ev{\hat{x}^k - x^\ast, \nabla f(x^k_m)} \nonumber \\
        &\overset{\eqref{eq:sqnorm-jensen}}{\leq} \sqn{\hat{x}^k - x^\ast} + \frac{\gamma^2}{M} \sum_{m=1}^{M} \sqn{\nabla f(x^k_m)} - \frac{2 \gamma}{M} \sum_{m=1}^{M} \ev{\hat{x}^k - x^k_m + x^k_m - x^\ast, \nabla f(x^k_m)} \nonumber \\
            &= \sqn{\hat{x}^k - x^\ast} + \frac{\gamma^2}{M} \sum_{m=1}^{M} \sqn{\nabla f(x^k_m) - \nabla f(x^\ast)} - \frac{2 \gamma}{M} \sum_{m=1}^{M} \ev{x^k_m - x^\ast, \nabla f(x^k_m)} \notag \\
            &\quad - \frac{2 \gamma}{M} \sum_{m=1}^{M} \ev{\hat{x}^k - x^k_m, \nabla f(x^k_m)}
        \nonumber \\
           & \overset{\eqref{eq:grad_dif_bregman}}{\leq} \sqn{\hat{x}^k - x^\ast} + \frac{2 L \gamma^2}{M} \sum_{m=1}^{M} \br{f(x^k_m) - f(x^\ast)} - \frac{2 \gamma}{M} \sum_{m=1}^{M} \ev{x^k_m - x^\ast, \nabla f(x^k_m)} \notag \\
            &\quad- \frac{2 \gamma}{M} \sum_{m=1}^{M} \ev{\hat{x}^k - x^k_m, \nabla f(x^k_m)}
        \nonumber \\
            &\overset{\eqref{eq:asm-strong-convexity}}{\leq} \sqn{\hat{x}^k - x^\ast} + \frac{2 \gamma}{M} \sum_{m=1}^{M} \br{\br{ \gamma L - 1 } \br{f(\hat{x}^k_m) - f(x^\ast)} - \frac{\mu}{2} \sqn{x^k_m - x^\ast}} \notag \\
            &\quad- \frac{2 \gamma}{M} \sum_{m=1}^{M} \ev{\hat{x}^k - x^k_m, \nabla f(x^k_m)}. \label{eq:lma-pia-proof-1}
    \end{align}
    To bound the last term in (\ref{eq:lma-pia-proof-1}) we use the generalized Young's inequality $2 \ev{a, b} \leq \zeta \sqn{a} + \zeta^{-1} \sqn{b}$ with $\zeta = 2L$:
    \begin{align}
        - 2 \ev{\hat{x}^k - x^k_m, \nabla f(x^k_m)} &\overset{\eqref{eq:youngs-inequality}}{\leq} 2L \sqn{x^k_m - \hat{x}^k} + \frac{1}{2L} \sqn{\nabla f(x^k_m)} \nonumber \\
        &= 2L \sqn{x^k_m - \hat{x}^k} + \frac{1}{2L} \sqn{\nabla f(x^k_m) - f(x^\ast)} \nonumber \\
        &\overset{\eqref{eq:grad_dif_bregman}}{\leq} 2L \sqn{x^k_m - \hat{x}^k} + \br{f(x^k_m) - f(x^\ast)}. \label{eq:lma-pia-proof-2}
    \end{align}
    Finally, using \eqref{eq:lma-pia-proof-2} in \eqref{eq:lma-pia-proof-1} we get,
    \begin{align*}
        \begin{split}
            \sqn{\hat{x}^k - \gamma \bar{g}^k - x^\ast} &\overset{\eqref{eq:lma-pia-proof-1}, \eqref{eq:lma-pia-proof-2}}{\leq} \sqn{\hat{x}^k - x^\ast}+ \frac{2 \gamma L}{M} \sum_{m=1}^{M} \sqn{\hat{x}^k - x^k_m} \\
            &\qquad + \frac{2 \gamma}{M} \sum_{m=1}^{M} \br{\br{ \gamma L - \frac{1}{2} } \br{f(\hat{x}^k_m) - f(x^\ast)} - \frac{\mu}{2} \sqn{x^k_m - x^\ast}} .
        \end{split}
    \end{align*}
\end{proof}

\begin{lemma}[Single-iterate optimality gap analysis]
    \label{lemma:optimality-gap-contraction}
    Suppose that Assumptions~\ref{asm:convexity-and-smoothness} and \ref{asm:finite-sum-stochastic-gradients} hold with identical data. Choose a stepsize $\gamma > 0$ such that $\gamma \leq \frac{1}{4 L \br{1 + \frac{2}{M}}}$ where $M$ is the number of nodes, then for expectation conditional on $x^k_1, x^k_2, \ldots, x^k_M$ we have
    \begin{align}
        \label{eq:lma-optimality-gap-contraction}
        \ecn{\hat{x}^{k+1} - x^\ast} \leq \br{1 - \gamma \mu} \sqn{\hat{x}^k - x^\ast} &+ 2 \gamma L V^k + \frac{2 \gamma^2 \sigmaopt^2}{M} - \frac{\gamma}{2} \br{f(\hat{x}^k) - f(x^\ast)},
    \end{align}
    where $\hat{x}^k = \frac{1}{M} \sum_{m=1}^{M} x^k_m$ and $V^k \eqdef \frac{1}{M} \sum_{m=1}^{M} \sqn{x^k_m - \hat{x}^k}$ is the iterate variance across the different nodes from their mean at timestep $k$.
\end{lemma}
\begin{proof}
    This is a modification of Lemma~3.1 in~\cite{Stich2018}. For expectation conditional on $(x^k_m)_{m=1}^{M}$ and using Lemma~\ref{lemma:perturbed-iterate-analysis},
    \begin{align*}
        \ecn{\hat{x}^{k+1} - x^\ast} &\overset{\eqref{eq:variance_def}}{=} \sqn{\hat{x}^k - x^\ast - \gamma \bar{g}^k} + \gamma^2 \ecn{g^k - \bar{g}^k} \\
        &\overset{\eqref{eq:lma-perturbed-iterate-analysis}}{\leq} \sqn{\hat{x}^k - x^\ast} + 2 \gamma L V^k + \gamma^2 \ecn{g^k - \bar{g}^k}\\
        &\qquad + \frac{2 \gamma}{M} \sum_{m=1}^{M} \br{ \br{\gamma L - \frac{1}{2}} \br{f(x^k_m) - f(x^\ast)} - \frac{\mu}{2} \sqn{x^k_m - x^\ast} }.
    \end{align*}
    Now use Lemma~\ref{lemma:minibatch-gradient-noise-reduction} to bound $\sqn{g^k - \bar{g}^k}$:
    \begin{align}
        \label{eq:lma-ogc-proof-1}
        \begin{split}
            \ecn{\hat{x}^{k+1} - x^\ast} &\overset{\eqref{eq:lma-mbns}}{\leq} \sqn{\hat{x}^k - x^\ast} + 2 \gamma L V^k + \frac{2 \gamma^2 \sigmaopt^2}{M} \\
            &+ \frac{2 \gamma}{M} \sum_{m=1}^{M} \br{ \br{\gamma L + \frac{2 \gamma L}{M} - \frac{1}{2}} \br{f(x^k_m) - f(x^\ast)} - \frac{\mu}{2} \sqn{x^k_m - x^\ast} }.
        \end{split}
    \end{align}
    We now use that the stepsize $\gamma \leq \frac{1}{4 L \br{1 + \frac{2}{M}}}$ to bound the last term in \eqref{eq:lma-ogc-proof-1},
    \begin{align}
        \label{eq:lma-ogc-proof-2}
        \begin{split}
            \ecn{\hat{x}^{k+1} - x^\ast} &\leq \sqn{\hat{x}^k - x^\ast} + 2 \gamma L V^k + \frac{2 \gamma^2 \sigmaopt^2}{M} \\
            &\qquad + \frac{2 \gamma}{M} \sum_{m=1}^{M} \br{ -\frac{1}{4} \br{f(x^k_m) - f(x^\ast)} - \frac{\mu}{2} \sqn{x^k_m - x^\ast} }.
        \end{split}
    \end{align}
    Applying Jensen's inequality from Proposition~\ref{pr:jensen} to $\frac{1}{4} \br{f(x^k_m)- f(x^\ast)} + \frac{\mu}{2} \sqn{x^k_m - x^\ast}$, we obtain
    \begin{align}
        \label{eq:lma-ogc-proof-3}
        - \avemm \br{\frac{1}{4} \br{f(x^k_m)- f(x^\ast)} + \frac{\mu}{2} \sqn{x^k_m - x^\ast} }
        &\overset{\eqref{eq:jensen}}{\leq} - \frac{1}{4} \br{f(\hat{x}^k) - f(x^\ast)} \\
        &\qquad - \frac{\mu}{2} \sqn{\hat{x}^k - x^\ast}.
    \end{align}
    Plugging~\eqref{eq:lma-ogc-proof-3} in \eqref{eq:lma-ogc-proof-2}, we get
    \begin{eqnarray*}
        \ecn{\hat{x}^{k+1} - x^\ast} 
        &\overset{\eqref{eq:lma-ogc-proof-2}, \eqref{eq:lma-ogc-proof-3}}{\leq} \br{1 - \gamma \mu} \sqn{\hat{x}^k - x^\ast} + 2 \gamma L V^k + \frac{2 \gamma^2 \sigmaopt^2}{M}\\
        & - \frac{\gamma}{2} \br{f(\hat{x}^k) - f(x^\ast)},
    \end{eqnarray*}
    which is the claim of this lemma.
\end{proof}

\begin{lemma}[Bounding the deviation of the gradients from their average]
    \label{lemma:average-gradient-deviation}
    Under Assumptions~\ref{asm:convexity-and-smoothness} and \ref{asm:finite-sum-stochastic-gradients} for identical data we have for all $k \geq 0$,
    \begin{equation}
        \label{eq:lma-average-gradient-deviation}
        \avemm \ecn{g^k_m - \avemm g^k_m} \leq 2 \sigmaopt^2 + \frac{4L}{M} \sum_{m=1}^{M} D_{f} (x^k_m, x^\ast).
    \end{equation}
\end{lemma}
\begin{proof}
    We start by the variance bound,
    \begin{align*}
        \avemm \sqn{g^k_m - \avemm g^k_m} &\overset{\eqref{eq:variance_m}}{=} \avemm \sqn{g^k_m} - \sqn{\avemm g^k_m} \\
        &\leq \avemm \sqn{g^k_m}.
    \end{align*}
    We now take expectations and use Lemma \ref{lemma:gradient-smoothness-bound}:
    \begin{align*}
        \avemm \ecn{g^k_m - \avemm g^k_m} &\leq \avemm \ecn{g^k_m} \\
        &\leq \frac{1}{M} \sum_{m=1}^{M} \br{2 \sigma_m^2 + 4L D_{f} (x^k_m, x^\ast)} \\
        &= 2 \sigmaopt^2 + \frac{4 L}{M} \sum_{m=1}^{M} D_{f} (x^k_m, x^\ast).
    \end{align*}
\end{proof}

\begin{lemma}
    \label{lemma:average-iterate-deviation-recursion}
    Suppose that Assumptions~\ref{asm:convexity-and-smoothness}, and \ref{asm:finite-sum-stochastic-gradients} hold for identical data. Choose $\gamma \leq \frac{1}{2 L}$, then for all $k \geq 0$:
    \begin{equation}
        \label{eq:lma-average-iterate-deviation-recursion}
        \ec{V^{k+1}} \leq \br{1 - \gamma \mu} \ec{V^k} + 2 \gamma \ec{D_{f} (\hat{x}^k, x^\ast)} + 2 \gamma^2 \sigmaopt^2.
    \end{equation}
\end{lemma}
\begin{proof}
    If $k + 1 = k_p$ for some $p \in \N$ then the left hand side is zero and the above inequality is trivially satisfied. If not, then recall that $x^{k+1}_m = x^k_m - \gamma g^k_m$ and $\hat{x}^{k+1} = \hat{x}^k - \gamma g^k$ where $\ec{g^k_m} = \nabla f(x^k_m)$ and $g^k = \avemm g^k_m$. Hence, for the expectation conditional on $(x^k_m)_{m=1}^{M}$ we have
    \begin{align*}
        \ecn{x^{k+1}_m - \hat{x}^{k+1}} &= \ecn{x^k_m - \hat{x}^k - \gamma \br{g^k_m - g^k}} \\
        &= \ecn{x^k_m - \hat{x}^k_m} + \gamma^2 \ecn{g^k_m - g^k} - 2 \gamma \ec{\ev{x^k_m - \hat{x}^k, g^k_m - g^k}} \\
        &= \ecn{x^k_m - \hat{x}^k_m} + \gamma^2 \ecn{g^k_m - g^k} - 2 \gamma \ev{x^k_m - \hat{x}^k, \nabla f(x^k_m) - \bar{g}^k},
    \end{align*}
    where $\bar{g}^k = \ec{\avemm g^k_m} = \avemm \nabla f(x^k_m)$. Averaging over $m$ in the last equality,
    \begin{align}
        \ec{V^{k+1}} &= V^k + \frac{\gamma^2}{M} \sum_{m=1}^{M}\ecn{g^k_m - g^k} - \frac{2\gamma}{M} \sum_{m=1}^{M} \ev{x^k_m - \hat{x}^k, \nabla f(x^k_m)} \nonumber\\
        &\qquad + \frac{2 \gamma}{M} \sum_{m=1}^{M} \ev{x^k_m - \hat{x}^k, \bar{g}^k} \nonumber \\
        &= V^k + \frac{\gamma^2}{M} \sum_{m=1}^{M} \ecn{g^k_m - g^k} - \frac{2\gamma}{M} \sum_{m=1}^{M} \ev{x^k_m - \hat{x}^k, \nabla f(x^k_m)} \nonumber \\
        &\qquad + 2 \gamma \underbrace{\ev{\hat{x}^k - \hat{x}^k, \bar{g}^k}}_{=0} \nonumber \\
        &= V^k + \frac{\gamma^2}{M} \sum_{m=1}^{M} \ecn{g^k_m - g^k} - \frac{2 \gamma}{M} \sum_{m=1}^{M} \ev{x^k_m - \hat{x}^k, \nabla f(x^k_m)}. \label{eq:lma-aidr-proof-1}
    \end{align}
    We now use Lemma~\ref{lemma:average-gradient-deviation} to bound the second term in \eqref{eq:lma-aidr-proof-1},
    \begin{align}
        \label{eq:lma-aidr-proof-2}
        \ec{V^{k+1}} &\overset{\eqref{eq:lma-average-gradient-deviation}}{\leq} V^k + \frac{4 L \gamma^2}{M} \sum_{m=1}^{M} D_{f} (x^k_m, x^\ast) + 2 \gamma^2 \sigmaopt^2 - \frac{2 \gamma}{M} \sum_{m=1}^{M} \ev{x^k_m - \hat{x}^k, \nabla f(x^k_m)}.
    \end{align}
    We now use Assumption~\ref{asm:convexity-and-smoothness} to bound the last term in \eqref{eq:lma-aidr-proof-2}:
    \begin{align}
        \label{eq:lma-aidr-proof-3}
        \ev{\hat{x}^k - x^k_m, \nabla f(x^k_m)} &\overset{\eqref{eq:asm-strong-convexity}}{\leq} f(\hat{x}^k) - f(x^k_m) - \frac{\mu}{2} \sqn{x^k_m - \hat{x}^k}.
    \end{align}
    Plugging \eqref{eq:lma-aidr-proof-3} into \eqref{eq:lma-aidr-proof-2},
    \begin{align}
        \label{eq:lma-aidr-proof-4}
        \ec{V^{k+1}}
         &\leq \br{1 - \gamma \mu} V^k + 2 \gamma^2 \sigmaopt^2 + \frac{4 L \gamma^2}{M} \sum_{m=1}^{M} \br{f(x^k_m) - f(x^\ast)}\\
        &\qquad + \frac{2 \gamma}{M} \sum_{m=1}^{M} \br{f(\hat{x}^k) - f(x^k_m)}.
    \end{align}
    Using that $\gamma \leq \frac{1}{2L}$ in \eqref{eq:lma-aidr-proof-4},
    \begin{align*}
        \ec{V^{k+1}} &\leq \br{1 - \gamma \mu} V^k + 2 \gamma^2 \sigmaopt^2 + \frac{2 \gamma}{M} \sum_{m=1}^{M} \br{f(x^k_m) - f(x^\ast) + f(\hat{x}^k) - f(x^k_m)} \\
        &= \br{1 - \gamma \mu} V^k + 2 \gamma^2 \sigmaopt^2 + 2\gamma \br{f(\hat{x}^k) - f(x^\ast)}.
    \end{align*}
    Taking unconditional expectations and using the tower property yields the lemma's statement.
\end{proof}

\begin{lemma}[Epoch iterate deviation bound]
    \label{lemma:epoch-iterate-deviation-bound}
    Suppose that Assumptions~\ref{asm:convexity-and-smoothness}, and \ref{asm:finite-sum-stochastic-gradients} hold with identical data. Assume that Algorithm~\ref{alg:local_sgd} is run with stepsize $\gamma > 0$, let $p \in \N$ be such that $k_{p}$ is a synchronization point then for $v=k_{p+1} - 1$ we have for $\alpha \eqdef 1 - \gamma \mu$,
    \begin{align*}
        \sum_{k=k_p}^{v} \alpha^{v-k} \cdot \ec{V^k} 
        &\leq \frac{2 \gamma \br{H-1}}{\alpha} \sum_{k=k_p}^{v} \alpha^{v-k} \cdot \ec{D_f (\hat{x}^k, x^\ast)} + 2 \gamma^2 \sigmaopt^2 \br{H - 1} \sum_{k=k_p}^{v} \alpha^{v-k}. 
    \end{align*}
\end{lemma}
\begin{proof}
    We start with Lemma~\ref{lemma:average-iterate-deviation-recursion},
    \begin{align*}
        \ec{V^{k}} &\leq \br{1 - \gamma \mu} \ec{V^{k-1}} + 2 \gamma \ec{D_{f} (\hat{x}^{k-1}, x^\ast)} + 2 \gamma^2 \sigmaopt^2 \\
        &= \alpha \cdot \ec{V^{k-1}} + 2 \gamma \ec{D_{f} (\hat{x}^{k-1}, x^\ast)} + 2 \gamma^2 \sigmaopt^2.
    \end{align*}
    By assumption there is some synchronization point $p \in \N$ such that $k_p \leq k \leq k_{p+1} - 1$ and $k_{p+1} - k_p \leq H$, recursing the above inequality until $k_{p}$ and using that $V^{k_p} = 0$,
    \begin{align}
        \ec{V^{k}} &\leq \alpha^{k - k_p} \ec{V^{k_p}} + 2 \gamma \sum_{j=k_p}^{k-1} \alpha^{k - j - 1} \ec{D_{f} (\hat{x}^{j}, x^\ast)} + 2 \sigmaopt^2 \sum_{j=k_p}^{k-1} \gamma^2 \alpha^{k - 1 - j} \nonumber \\
        \label{eq:lma-eidb-1}
        &= \frac{2 \gamma}{\alpha} \sum_{j=k_p}^{k-1} \alpha^{k - j} \ec{D_{f} (\hat{x}^{j}, x^\ast)} + 2 \gamma^2 \sigmaopt^2 \sum_{j=k_p}^{t-1} \alpha^{k - 1 - j}.
    \end{align}
    The second term in \eqref{eq:lma-eidb-1} can be bounded as follows: because $\alpha \leq 1$ then $\alpha^{k - 1 - j} \leq 1$ for $j \leq k-1$, hence for $k \leq k_{p+1} - 1$
    \begin{align}
        2 \gamma^2 \sigmaopt^2 \sum_{j=k_p}^{t-1} \alpha^{k - 1 - j} &\leq 2 \gamma^2 \sigmaopt^2 \sum_{j=k_p}^{k-1} 1 \nonumber \\
        &= 2 \gamma^2 \sigmaopt^2 \br{k - k_p} \nonumber \\
        &\leq 2 \gamma^2 \sigmaopt^2 \br{k_{p+1} - k_p - 1} \nonumber \\
        \label{eq:lma-eidb-2}
        &\leq 2 \gamma^2 \sigmaopt^2 \br{H - 1}.
    \end{align}
    Using \eqref{eq:lma-eidb-2} in \eqref{eq:lma-eidb-1},
    \begin{align}
        \label{eq:lma-eidb-3}
        \ec{V^k} &\leq \frac{2 \gamma}{\alpha} \sum_{j=k_p}^{k-1} \alpha^{k - j} \ec{D_{f} (\hat{x}^{j}, x^\ast)} + 2 \gamma^2 \sigmaopt^2 \br{H - 1}.
    \end{align}
    Then summing up \eqref{eq:lma-eidb-3} weighted by $\alpha^{v-k}$ for $v = k_{p+1} - 1$,
    \begin{align}
        \label{eq:lma-eidb-4}
        \begin{split}
            \sum_{k=k_p}^{v} \alpha^{v-k} \ec{V^k} 
            &\leq \frac{2 \gamma}{\alpha} \sum_{k=k_p}^{v} \alpha^{v-k} \sum_{j=k_p}^{k-1} \alpha^{k-j} \ec{D_f (\hat{x}^j, x^\ast)} + 2 \gamma^2 \sum_{k=k_p}^{v} \alpha^{v-k} \sigmaopt^2 \br{H - 1}.
        \end{split}
    \end{align}
    We now bound the first term in \eqref{eq:lma-eidb-4} by adding more terms in the inner sum, since $D_{f} (\hat{x}^j, x^\ast) \geq 0$ and $k -1 \leq v - 1 \leq v$:
    \begin{align}
        \sum_{k=k_p}^{v} \alpha^{v-k} \sum_{j=k_p}^{k-1} \alpha^{k-j} \ec{D_f (\hat{x}^j, x^\ast)} 
        &\leq \sum_{k=k_p}^{v} \alpha^{v-k} \sum_{j=k_p}^{v} \alpha^{k-j} \ec{D_f (\hat{x}^j, x^\ast)} \nonumber \\
        &= \sum_{k=k_p}^{v} \sum_{j=k_p}^{v} \alpha^{v - j} \ec{D_f (\hat{x}^j, x^\ast)} \nonumber \\
        &= \br{v - k_p} \sum_{j=k_p}^{v} \alpha^{v - j} \ec{D_f (\hat{x}^j, x^\ast)} \nonumber \\
        &= \br{k_{p+1} - k_p - 1} \sum_{j=k_p}^{v} \alpha^{v - j} \ec{D_f (\hat{x}^j, x^\ast)} \nonumber \\
        \label{eq:lma-eidb-5}
        &\leq \br{H-1} \sum_{j=k_p}^{v} \alpha^{v - j} \ec{D_f (\hat{x}^j, x^\ast)}.
    \end{align}
    Combining \eqref{eq:lma-eidb-5} and \eqref{eq:lma-eidb-4} we have,
    \begin{align*}
        \sum_{k=k_p}^{v} \alpha^{v-k} \cdot \ec{V^k} &\leq \frac{2 \gamma \br{H-1}}{\alpha} \sum_{j=k_p}^{v} \alpha^{v-j} \cdot \ec{D_f (\hat{x}^j, x^\ast)} + 2 \gamma^2 \sigmaopt^2 \br{H - 1} \sum_{k=k_p}^{v} \alpha^{v-k}. 
    \end{align*}
    Finally, renaming the variable $j$ gives us the claim of this lemma.
\end{proof}

\subsection{Proof of Theorem~\ref{theorem:sc-unbounded-variance-iid}}
\begin{proof}
    Let $(k_{p})_{p}$ index all the times $k$ at which communication and averaging happen. Taking expectations in Lemma~\ref{lemma:optimality-gap-contraction} and letting $r^{k} = \hat{x}^k - x^\ast$,
    \begin{align}
        \label{smooth-case-thm-proof-3}
        \ecn{r^{k+1}} &\leq \br{1 - \gamma \mu} \ecn{r^k} + 2 \gamma L \ec{V^{k}} + \frac{2 \gamma^2 \sigmaopt^2}{M} - \frac{\gamma}{2} \ec{D_{f} (\hat{x}^k, x^\ast)} \\
        &= \br{1 - \gamma \mu} \ecn{r^k} + \br{2 \gamma L \ec{V^k} - \frac{\gamma}{2} D_{f} (\hat{x}^k, x^\ast)} + \frac{2 \gamma^2 \sigmaopt^2}{M}.
    \end{align}
    Let $K = k_p - 1$ for some $p \in \N$, then expanding out $\ecn{r^k}$ in \eqref{smooth-case-thm-proof-3},
    \begin{align}
        \ecn{r^{K+1}} &\leq \br{1 - \gamma \mu}^{K+1} \ecn{\hat{x}^{0} - x^\ast} + \sum_{k=0}^{K} \br{1 - \gamma \mu}^{K-i} \frac{2 \gamma^2 \sigmaopt^2}{M} \nonumber \\
        &\qquad + \sum_{i=0}^{K} \br{1 - \gamma \mu}^{K-i} \br{2 \gamma L \ec{V^i} - \frac{\gamma}{2} D_{f} (\hat{x}^i, x^\ast)} \nonumber \\
        &\leq \br{1 - \gamma \mu}^{K+1} \ecn{x^0 - x^\ast} + \frac{2 \gamma \sigmaopt^2}{\mu M} \notag\\
        &\qquad + \frac{\gamma}{2} \sum_{i=0}^{K} \br{1 - \gamma \mu}^{K-i} \ec{4 L V^i - D_{f} (\hat{x}^i, x^\ast)}. \label{smooth-case-thm-proof-4}
    \end{align}
    It remains to bound the last term in \eqref{smooth-case-thm-proof-4}. We have
    \begin{align}
        \sum_{i=0}^{K} \br{1 - \gamma \mu}^{K-i} &(4 L \ec{V^i} - D_{f} (\hat{x}^i, x^\ast)) \notag \\
        &= \sum_{l=1}^{p} \sum_{i=k_{l-1}}^{k_l - 1} \br{1 - \gamma \mu}^{K-i} \br{4 L \ec{V^i} - D_{f} \br{\hat{x}^i, x^\ast}} \nonumber \\
        \label{smooth-case-thm-proof-5}
        &= \sum_{l=1}^{p} \br{1 - \gamma \mu}^{K - \br{k_l - 1}} \sum_{k_{l-1}}^{k_l - 1} \br{1 - \gamma \mu}^{k_l - 1 - i} \ec{4 L V^i - D_{f} \br{\hat{x}^i, x^\ast}},
    \end{align}
    where in the first line we just count $i$ by decomposing it over all the communication intervals. Fix $l \in \N$ and let $v_l = k_{l} - 1$. Then by Lemma~\ref{lemma:epoch-iterate-deviation-bound} we have,
    \begin{align}
        \sum_{i=k_l}^{v_l} \br{1 - \gamma \mu}^{v_l - i} \ec{V^i} 
        &\leq \frac{2 \gamma (H - 1)}{\alpha} \sum_{i=k_l}^{v_l} \alpha^{v_l - i} \ec{D_{f} (\hat{x}^i, x^\ast)} + \sum_{i=k_l}^{v_l} \alpha^{v_l - i} 2 \gamma^2 \sigma^2 (H - 1),
        \label{smooth-case-thm-proof-2}
    \end{align}
    where $\alpha = 1 - \gamma \mu$. Using \eqref{smooth-case-thm-proof-2} in \eqref{smooth-case-thm-proof-5},
    \begin{align}
        4 L &\sum_{i = k_{l-1}}^{v_l} \br{1 - \gamma \mu}^{v_l - i} \ec{V^i} - \sum_{i=k_{l-1}}^{v_l} \br{1 - \gamma \mu}^{v_l-i} \ec{D_{f} (\hat{x}^i, x^\ast)} \nonumber \\
        &\leq 4 L  \frac{2 \gamma \br{H - 1}}{1 - \gamma \mu} \sum_{i=k_{l-1}}^{v_l} \br{1 - \gamma \mu}^{v_l-i} \ec{D_{f} (\hat{x}^i, x^\ast)} \nonumber\\
        &\quad+ 4L\sum_{i=k_{l-1}}^{v_l} \br{1 - \gamma \mu}^{v_l - i} 2 \gamma^2 \sigmaopt^2 \br{H - 1}  - \sum_{i=k_{l-1}}^{v_l} \br{1 - \gamma \mu}^{v_l-i} \ec{D_{f} (\hat{x}^i, x^\ast)} \nonumber \\
        &= \sum_{i=k_{l-1}}^{v_l} \br{1 - \gamma \mu}^{v_l - i} 8 \gamma^2 \sigmaopt^2 \br{H - 1} L \nonumber\\
        &\quad - \sum_{i=k_{l-1}}^{v_l} \br{1 - \frac{8 \gamma L \br{H - 1}}{1 - \gamma \mu}} \br{1 - \gamma \mu}^{v_l-i} \ec{D_{f} (\hat{x}^i, x^\ast)} \nonumber \\
        \label{smooth-case-thm-proof-6}
        &\leq \sum_{i=k_{l-1}}^{v_l} \br{1 - \gamma \mu}^{v_l - i} 8 \gamma^2 \sigmaopt^2 \br{H - 1} L,
    \end{align}
    where in the third line we used that our choice of $\gamma$ guarantees that $1 - \frac{8 \gamma L H}{1 - \gamma \mu} \geq 0$. Using \eqref{smooth-case-thm-proof-6} in \eqref{smooth-case-thm-proof-5},
    \begin{align}
        \sum_{i=0}^{K} \br{1 - \gamma \mu}^{K-i} &\ec{4 L V^i - D_{f} (\hat{x}^i, x^\ast)} \nonumber \\
        &\leq \sum_{l=1}^{p} \br{1 - \gamma \mu}^{K - \br{k_l - 1}} \sum_{i= k_{l-1}}^{k_l - 1} \br{1 - \gamma \mu}^{k_l - 1 - i} \ec{4 L V^i - D_{f} \br{\hat{x}^i, x^\ast}} \nonumber \\
        &\leq \sum_{l=1}^{p} \br{1 - \gamma \mu}^{K - \br{k_l - 1}} \sum_{i=k_{l-1}}^{k_l - 1} \br{1 - \gamma \mu}^{k_l - 1 - i} 8 \gamma^2 \sigmaopt^2 \br{H - 1} L \nonumber \\
        &= \sum_{kl=1}^{p} \sum_{i=k_{l-1}}^{k_l - 1} \br{1 - \gamma \mu}^{K - i} 8 \gamma^2 \sigma^2 \br{H - 1} L \nonumber \\
        &= \sum_{i=0}^{K} \br{1 - \gamma \mu}^{K - i} 8 \gamma^2 \sigmaopt^2 \br{H - 1} L \nonumber \\
        \label{smooth-case-thm-proof-7}
        &\leq \frac{8 \sigmaopt^2 \gamma \br{H - 1} L}{\mu}.
    \end{align}
    Using \eqref{smooth-case-thm-proof-7} in \eqref{smooth-case-thm-proof-4},
    \begin{align*}
        \ecn{r^{K+1}} 
        &\leq \br{1 - \gamma \mu}^{K+1} \ecn{x^0 - x^\ast} \\
        &\qquad + \frac{2 \gamma \sigmaopt^2}{\mu M} + \frac{\gamma}{2} \sum_{i=0}^{K} \br{1 - \gamma \mu}^{K-i} \ec{4 L V^i - D_{f} (\hat{x}^i, x^\ast)} \\
        &\leq \br{1 - \gamma \mu}^{K + 1} \ecn{x^0 - x^\ast} + \frac{2 \gamma \sigmaopt^2}{\mu M} + \frac{4 \sigmaopt^2 \gamma^2 \br{ H - 1} L}{\mu},
    \end{align*}
    which is the claim of the theorem.
\end{proof}

\subsection{Proof of Theorem~\ref{thm:wc-iid-unbounded-var}}
\begin{proof}
    Start with Lemma~\ref{lemma:optimality-gap-contraction} with $\mu = 0$, then the conditional expectations satisfies
    \begin{align*}
        \ecn{\hat{x}^{k+1} - x^\ast} &\overset{\eqref{eq:lma-optimality-gap-contraction}}{\leq} \sqn{\hat{x}^k - x^\ast} + 2 \gamma L V^k + \frac{2 \gamma^2 \sigmaopt^2}{M} - \frac{\gamma}{2} D_{f} (\hat{x}^k, x^\ast).
    \end{align*}
    Taking full expectations and rearranging,
    \begin{align*}
        \frac{\gamma}{2} \ec{D_{f} (\hat{x}^k, x^\ast)} &\leq \ecn{\hat{x}^k - x^\ast} - \ecn{\hat{x}^{k+1} - x^\ast} + 2 \gamma L \ec{V^k} + \frac{2 \gamma^2 \sigmaopt^2}{M}.
    \end{align*}
    Averaging as $k$ varies from $0$ to $K-1$,
    \begin{align}
        \frac{\gamma}{2 K} \sum_{k=0}^{K-1} \ec{D_{f} (\hat{x}^k, x^\ast)} 
        &\leq \frac{1}{K} \sum_{k=0}^{K-1} \br{\ecn{\hat{x}^k - x^\ast} - \ecn{\hat{x}^{k+1} - x^\ast}} \notag\\
        &\quad + \frac{2 \gamma L}{K} \sum_{k=0}^{K-1} \ec{V^k} + \frac{2 \gamma^2 \sigmaopt^2}{M} \nonumber \\
        &= \frac{1}{K} \br{\sqn{x^0 - x^\ast} - \ecn{\hat{x}^{K} - x^\ast}} + \frac{2 \gamma L}{K} \sum_{k=0}^{K-1} \ec{V^k} \nonumber\\
        &\quad + \frac{2 \gamma^2 \sigmaopt^2}{M} \nonumber \\
        \label{eq:thm-wcc-proof-1}
        &\leq \frac{\sqn{x^0 - x^\ast}}{K} + \frac{2 \gamma L}{K} \sum_{k=0}^{K-1} \ec{V^k} + \frac{2 \gamma^2 \sigmaopt^2}{M}.
    \end{align}
    To bound the sum of deviations in \eqref{eq:thm-wcc-proof-1}, we use Lemma~\ref{lemma:epoch-iterate-deviation-bound} with $\mu = 0$ (and noticing that because $\mu = 0$ we have $\alpha = 1$),
    \begin{align}
        \label{eq:thm-wcc-proof-2}
        \sum_{k=k_p}^{k_{p+1} - 1} \ec{V^k} &\leq \sum_{k=k_p}^{k_{p+1} - 1} \br{2 \gamma (H - 1) \ec{D_f (\hat{x}^k, x^\ast)} + 2 \gamma^2 \sigmaopt^2 (H - 1) }.
    \end{align}
    Since by assumption $K$ is a synchronization point, there is some $l \in \N$ such that $K = k_{l}$. To estimate the sum of deviations in \eqref{eq:thm-wcc-proof-1} we use double counting to decompose it over each epoch, use \eqref{eq:thm-wcc-proof-2}, and then use double counting again:
    \begin{align}
        \sum_{k=0}^{K - 1} \ec{V^k} 
        &= \sum_{p=0}^{l-1} \sum_{k=k_{p}}^{k_{p+1} - 1} \ec{V^k} \nonumber \\
        &\overset{\eqref{eq:thm-wcc-proof-2}}{\leq} \sum_{p=0}^{l-1} \sum_{k=k_{p}}^{k_{p+1} - 1} \br{2 \gamma (H-1) \ec{D_{f} (\hat{x}^k, x^\ast)} + 2 \gamma^2 \sigmaopt^2 \br{H - 1}} \nonumber \\
        \label{eq:thm-wcc-proof-3}
        &= \sum_{k=0}^{K-1} \br{2 \gamma \br{H - 1} \ec{D_{f} (\hat{x}^k, x^\ast)} + 2 \gamma^2 \sigmaopt^2 (H - 1)}.
    \end{align}
    Using \eqref{eq:thm-wcc-proof-3} in \eqref{eq:thm-wcc-proof-1} and rearranging we get,
    \begin{align*}
        \frac{\gamma}{2K} \sum_{k=0}^{K-1} \ec{D_{f} (\hat{x}^k, x^\ast)} 
        &\leq \frac{\sqn{x^0 - x^\ast}}{K}  + \frac{2 \gamma^2 \sigmaopt^2}{M}\\
        &\quad + \frac{2 \gamma L}{K} \sum_{k=0}^{K-1} \br{2 \gamma (H-1) \ec{D_{f} (\hat{x}^k, x^\ast)} + 2 \gamma^2 \sigmaopt^2 (H-1)}.
    \end{align*}
    Therefore,
    \begin{align*}
        \frac{\gamma}{2K} \sum_{k=0}^{K-1} \br{1 - 8 \gamma (H-1) L} \ec{D_{f} (\hat{x}^k, x^\ast)} 
        &\leq \frac{\sqn{x^0 - x^\ast}}{K} + \frac{2 \gamma^2 \sigmaopt^2}{M} + 4 \gamma^3 L \sigmaopt^2 (H-1).	
    \end{align*}
    By our choice of $\gamma$ we have that $1 - 8 \gamma L (H-1) \geq \frac{2}{10}$. Using this with some algebra, we get
    \begin{align*}
        \frac{\gamma}{10 K} \sum_{k=0}^{K-1} \ec{D_{f} (\hat{x}^k, x^\ast)} 
        &\leq \frac{\sqn{x^0 - x^\ast}}{K} + \frac{2 \gamma^2 \sigmaopt^2}{M} + 4 \gamma^3 L \sigmaopt^2 (H-1).
    \end{align*}
    Dividing both sides by $\gamma/10$ and using Jensen's inequality yields the theorem's claim.
\end{proof}

\section{Proofs for Heterogeneous Data}
\subsection{Preliminary lemmas}
\begin{lemma}
    \label{lemma:average-gradient-bound}
    Suppose that Assumptions~\ref{asm:convexity-and-smoothness} and \ref{asm:finite-sum-stochastic-gradients} hold with $\mu \geq 0$ for heterogeneous data. Then for expectation conditional on $x^k_1, x^k_2, \ldots, x^k_m$ and for $M \geq 2$, we have
    \begin{align}
        \label{eq:lma-average-gradient-bound}
        \ecn{g^k} \leq 2 L^2 V^k + 8 L D_{f} (\hat{x}^k, x^\ast) + \frac{4 \sigmaf^2}{M}.
    \end{align}
\end{lemma}
\begin{proof}
    Starting with the left-hand side,
    \begin{align}
        \label{eq:lma-agb-proof-1}
        \ecn{g^k} 
        &\overset{\eqref{eq:sum_sqnorm}}{\leq} 2 \mathbb{E}\left[\biggl\|g^k - \frac{1}{M} \sum_{m=1}^{M} \nabla f_m (\hat{x}^k; \xi_m)\biggr\| \right] + 2 \mathbb{E}\left[\biggl\| \frac{1}{M} \sum_{m=1}^{n} \nabla f_m (\hat{x}^k; \xi_m) \biggr\| \right].
    \end{align}
    To bound the first term in \eqref{eq:lma-agb-proof-1}, we use $L$-smoothness of $f_m (\cdot; \xi_m)$ to obtain
    \begin{align}
        2 \ecn{g^k - \frac{1}{M} \sum_{m=1}^{M} \nabla f_m (\hat{x}^k; \xi_m) } &= 2 \ecn{ \frac{1}{M} \sum_{m=1}^{M}  (\nabla f_m (x^k_m; \xi_m) - \nabla f_m (\hat{x}^k; \xi_m))} \nonumber
        \\
        &\leq \frac{2}{M} \sum_{m=1}^{M} \ecn{\nabla f_m (x^k_m; \xi_m) - \nabla f_m (\hat{x}^k; \xi_m)} \nonumber \\
        \label{eq:lma-agb-proof-2}
        &\leq  \frac{2 L^2}{M} \sum_{m=1}^{M} \sqn{x^k_m - \hat{x}^k},
    \end{align}
    where in the second inequality we have used Jensen's inequality and the convexity of the map $x\mapsto \|x\|^2$. For the second term in \eqref{eq:lma-agb-proof-1}, we have
    \begin{align}
            \ecn{ \frac{1}{M} \sum_{m=1}^{M} \nabla f_m (\hat{x}^k; \xi_m) } &\overset{\eqref{eq:variance_def}}{=} \ecn{ \frac{1}{M} \sum_{m=1}^{M} \nabla f_m (\hat{x}^k; \xi_m) - \frac{1}{M} \sum_{m=1}^{M} \nabla f_m (\hat{x}^k)} \nonumber \\
            &\qquad+ \sqn{ \frac{1}{M} \sum_{m=1}^{M} \nabla f_m (\hat{x}^k)}. 
        \label{eq:lma-agb-proof-2-2}
    \end{align}
    For the first term in \eqref{eq:lma-agb-proof-2-2} we have by the independence of $\xi_1, \xi_2, \ldots, \xi_m$,
    \begin{align*}
        &\ecn{ \frac{1}{M} \sum_{m=1}^{M} \nabla f_m (\hat{x}^k; \xi_m) - \frac{1}{M} \sum_{m=1}^{M} \nabla f_m (\hat{x}^k)}\\
        &\quad = \frac{1}{M^2} \sum_{m=1}^{M} \ecn{\nabla f_m (\hat{x}^k; \xi_m) - \nabla f_m (\hat{x}^k)} \\
        &\quad \overset{\eqref{eq:variance_sqnorm_upperbound}}{\leq} \frac{1}{M^2} \sum_{m=1}^{M} \ecn{\nabla f_m (\hat{x}^k; \xi_m)} \\
        &\quad \overset{\eqref{eq:sum_sqnorm}}{\leq} \frac{2}{M^2} \sum_{m=1}^{M} \ecn{\nabla f_m (\hat{x}^k; \xi_m) - \nabla f_m (x^\ast; \xi_m)} + \frac{2}{M^2} \sum_{m=1}^{M} \ecn{\nabla f_m (x^\ast; \xi_m)} \\
        &\quad \overset{\eqref{eq:grad_dif_bregman}}{\leq} \frac{4 L}{M^2} \sum_{m=1}^{M} D_{f_m} (\hat{x}^k, x^\ast) + \frac{2 \sigmaf^2}{M}\\
        &\quad = \frac{4L}{M} D_{f} (\hat{x}^k, x^\ast) + \frac{2 \sigmaf^2}{M}.
    \end{align*}
    Using this in \eqref{eq:lma-agb-proof-2-2} we have,
    \begin{align*}
         \ecn{ \frac{1}{M} \sum_{m=1}^{M} \nabla f_m (\hat{x}^k; \xi_m) } &\leq \frac{4L}{M} D_{f} (\hat{x}^k, x^\ast) + \frac{2 \sigmaf^2}{M} + \ecn{ \frac{1}{M} \sum_{m=1}^{M} \nabla f_m (\hat{x}^k)} \\
         &= \frac{4L}{M} D_{f} (\hat{x}^k, x^\ast) + \frac{2 \sigmaf^2}{M} + \sqn{ \nabla f(\hat{x}^k)}.
    \end{align*}
    Now notice that
    \[ \sqn{\nabla f(\hat{x}^k)} = \sqn{\nabla f(\hat{x}^k) - \nabla f(x^\ast)} \leq 2 L D_{f} (\hat{x}^k, x^\ast).  \]
    Using this in the previous inequality we get
    \begin{align*}
        \ecn{\frac{1}{M} \sum_{m=1}^{M} \nabla f_m (\hat{x}^k; \xi_m)} \leq 2L \br{1 + \frac{2}{M}} D_{f} (\hat{x}^k, x^\ast) + \frac{2 \sigmaf^2}{M}.
    \end{align*}
    Because $M \geq 2$ we have $1 + \frac{2}{M} \leq 2$, hence
    \begin{align}
        \label{eq:lma-agb-proof-3}
        \ecn{\frac{1}{M} \sum_{m=1}^{M} \nabla f_m (\hat{x}^k; \xi_m)} \leq 4 L D_{f} (\hat{x}^k, x^\ast) + \frac{2 \sigmaf^2}{M}.
    \end{align}
    Combining \eqref{eq:lma-agb-proof-2} and \eqref{eq:lma-agb-proof-3} in \eqref{eq:lma-agb-proof-1} we have,
    \begin{align*}
        \ecn{g^k} \leq 2 L^2 V^k + 8 L D_{f} (\hat{x}^k, x^\ast) + \frac{4 \sigmaf^2}{M}.
    \end{align*}
\end{proof}

\begin{lemma}
    \label{lemma:inner-product-bound}
    Suppose that Assumption~\ref{asm:convexity-and-smoothness} holds with $\mu \geq 0$ for heterogeneous data (i.e., it holds for each $f_m$ for $m = 1, 2, \ldots, M$). Then we have,
    \begin{equation}
        \label{eq:lma-inner-product-bound}
        -\frac{2}{M} \sum_{m=1}^{M} \ev{\hat{x}^k - x^\ast, \nabla f_m (x^k_m)} \leq - 2  D_{f} (\hat{x}^k, x^\ast) -  \mu \sqn{ \hat{x}^k - x^\ast} +  L V^k.
    \end{equation}
\end{lemma}
\begin{proof}
    Starting with the left-hand side,
    \begin{align}
        \label{eq:lma-inner-prod-proof-1}
        - 2  \ev{\hat{x}^k - x^\ast, \nabla f_m (x^k_m)} 
        &= - 2  \ev{x^k_m - x^\ast, \nabla f_m (x^k_m)} - 2  \ev{\hat{x}^k - x^k_m, \nabla f_m (x^k_m)}.
    \end{align}
    The first term in \eqref{eq:lma-inner-prod-proof-1} is bounded by strong convexity:
    \begin{align}
        \label{eq:lma-inner-prod-proof-2}
        - \ev{x^k_m - x^\ast, \nabla f_m (x^k_m)} &\leq f_m (x^\ast) - f_m (x^k_m) - \frac{\mu}{2} \sqn{x^k_m - x^\ast}.
    \end{align}
    For the second term, we use $L$-smoothness,
    \begin{align}
        \label{eq:lma-inner-prod-proof-3}
        - \ev{\hat{x}^k  - x^k_m, \nabla f_m (x^k_m)} \leq f_m (x^k_m) - f_m (\hat{x}^k) + \frac{L}{2} \sqn{x^k_m - \hat{x}^k}.
    \end{align}
    Combining \eqref{eq:lma-inner-prod-proof-3} and \eqref{eq:lma-inner-prod-proof-2} in \eqref{eq:lma-inner-prod-proof-1},
    \begin{align*}
        - 2  \ev{\hat{x}^k - x^\ast, \nabla f_m (x^k_m)} &\leq 2  \br{f_m (x^\ast) - f_m (x^k_m) - \frac{\mu}{2} \sqn{x^k_m - x^\ast}} \\
        &+ 2  \br{f_m (x^k_m) - f_m (\hat{x}^k) + \frac{L}{2} \sqn{x^k_m - \hat{x}^k}} \\
        &= 2  \br{f_m (x^\ast) - f_m (\hat{x}^k) - \frac{\mu}{2} \sqn{x^k_m - x^\ast} + \frac{L}{2} \sqn{x^k_m - \hat{x}^k}}.
    \end{align*}
    Averaging over $m$,
    \begin{align*}
        -\frac{2 }{M} \sum_{m=1}^{M} \ev{\hat{x}^k - x^\ast, \nabla f_m (x^k_m)} 
        &\leq - 2  \br{f(\hat{x}^k) - f(x^\ast)} - \frac{ \mu}{M} \sum_{m=1}^{M} \sqn{x^k_m - x^\ast}\\
        &\quad + \frac{ L}{M} \sum_{m=1}^{M} \sqn{x^k_m - \hat{x}^k}.
    \end{align*}
    Note that the first term is the Bregman divergence $D_{f} (\hat{x}^k, x^\ast)$, and using Jensen's inequality we have $- \frac{1}{M} \sum_{m=1}^{M} \sqn{x^k_m - x^\ast} \leq - \sqn{\hat{x}^k - x^\ast}$, hence
    \begin{align*}
        -\frac{2 }{M} \sum_{m=1}^{M} \ev{\hat{x}^k - x^\ast, \nabla f_m (x^k_m)} &\leq - 2  D_{f} (\hat{x}^k, x^\ast) -  \mu \sqn{\hat{x}^k - x^\ast} +  L V^k,
    \end{align*}
    which is the claim of this lemma.
\end{proof}

\begin{lemma}
    \label{lemma:iterate-deviation-epoch}
    Suppose that Assumptions~\ref{asm:convexity-and-smoothness} and \ref{asm:finite-sum-stochastic-gradients} hold for Algorithm~\ref{alg:local_sgd} with heterogeneous data and with $\sup_{p} \abs{k_p - k_{p+1}} \leq H$. Let $p \in \N$, then for $v = k_{p+1} - 1$ and $\gamma \leq \frac{1}{4 L \br{H - 1}}$ we have,
    \begin{equation}
        \label{eq:lma-iterate-deviation-epoch}
        \sum_{k=k_p}^{v} \ec{V^k} \leq 8 L \gamma^2 \br{H-1}^2 \sum_{k=k_p}^{v} \ec{D_{f} (\hat{x}^k, x^\ast)} + 4 \gamma^2 \br{H-1}^2 \sum_{k=k_p}^{v} \sigmaf^2.
    \end{equation}
\end{lemma}
\begin{proof}
    Let $k$ be such that $k_p \leq k \leq k_{p+1} - 1 = v$. From the definition of $V^k$,
    \begin{align*}
        \ec{V^k} &= \frac{1}{M} \sum_{m=1}^{M} \ecn{x^k_m - \hat{x}^k} \\
        &= \frac{1}{M} \sum_{m=1}^{M} \ecn{ \br{x_{k_p}^{m} - \gamma \sum_{i=k_p}^{k-1} g^i_m } - \br{x_{k_p} - \gamma \sum_{i=k_p}^{k-1} g^i } }.
    \end{align*}
    Using that $x_{k_p} = x_{k_p}^{m}$ for all $m$ we have,
    \begin{align*}
        \ec{V^k} &= \frac{\gamma^2}{M} \sum_{m=1}^{M} \ecn{ \sum_{i=k_p}^{k-1} g^i_m - g^i } \\
        &\overset{\eqref{eq:sqnorm-jensen}}{\leq} \frac{\gamma^2 \br{k - k_p}}{M} \sum_{m=1}^{M} \sum_{i=k_p}^{k-1} \ecn{g^i_m - g^i} \\
        &\overset{\eqref{eq:variance_sqnorm_upperbound}}{\leq} \frac{\gamma^2 \br{k - k_p}}{M} \sum_{m=1}^{M} \sum_{i=k_p}^{t-1} \ecn{g^i_m} \\
        &\leq \frac{\gamma^2 \br{H - 1}}{M} \sum_{m=1}^{M}  \sum_{i=k_p}^{k-1} \ecn{g^i_m},
    \end{align*}
    where in the third line we used that $g^i = \avemm g^i_m$ and $\avemm \ecn{g^i_m - g^i} \leq \avemm \ecn{g^i_m}$, while in the fourth line we used that $k - k_p \leq k_{p+1} - k_p - 1 \leq H - 1$.  Summing up as $k$ varies from $k_p$ to $v$,
    \begin{align*}
        \sum_{k=k_p}^{v} \ec{V^k} 
        &\leq \sum_{k=k_p}^{v} \frac{\gamma^2 (H-1)}{M} \sum_{m=1}^{M}  \sum_{i=k_p}^{k-1} \ecn{g^i_m}.
    \end{align*}
    Because $k - 1 \leq v - 1 \leq v$ we can upper bound the inner sum as follows
    \begin{align}
        \sum_{k=k_p}^{v} \ec{V^k} 
        &\leq \sum_{k=k_p}^{v} \frac{\gamma^2 (H-1)}{M} \sum_{m=1}^{M}  \sum_{i=k_p}^{k-1} \ecn{g^i_m} \nonumber \\
        &\leq \sum_{k=k_p}^{v} \frac{\gamma^2 (H-1)}{M} \sum_{m=1}^{M}  \sum_{i=k_p}^{v-1} \ecn{g^i_m} \nonumber \\
        &= \frac{\gamma^2 (H-1) \br{v - k_p}}{M} \sum_{m=1}^{M} \sum_{k=k_p}^{v - 1} \ecn{g^k_m} \nonumber \\
        &\leq \frac{\gamma^2 (H-1)^2}{M} \sum_{m=1}^{M} \sum_{k=k_p}^{v-1} \ecn{g^k_m} \nonumber \\
        \label{eq:lma-ide-proof-1}
        &\leq \frac{\gamma^2 (H-1)^2}{M} \sum_{m=1}^{M} \sum_{i=k_p}^{v} \ecn{g^i_m}.
    \end{align}
    To bound the gradient norm term in \eqref{eq:lma-ide-proof-1}, we have
    \begin{align}
        \label{eq:lma-ide-proof-2}
        \begin{split}
            \ecn{g^i_m} \leq 3 \ecn{g^i_m - \nabla f_m (\hat{x}^i; \xi_m)} &+ 3 \ecn{\nabla f_m (\hat{x}^i; \xi_m) - \nabla f_m (x^\ast; \xi_m)} \\
            &+ 3 \ecn{\nabla f_m (x^\ast; \xi_m)}.
        \end{split}
    \end{align}
    The first term in \eqref{eq:lma-ide-proof-2} can be bounded by smoothness:
    \begin{align}
        \ecn{g^i_m - \nabla f_m (\hat{x}^i; \xi_m)} 
        &= \ecn{\nabla f_m (x^i_m; \xi_m) - \nabla f_m (\hat{x}^i; \xi_m)} \\
        &\leq L^2 \ecn{x^i_m - \hat{x}^i}.
        \label{eq:lma-ide-proof-3}
    \end{align}
    The second term in \eqref{eq:lma-ide-proof-2} can be bounded by smoothness and convexity:
    \begin{equation}
        \label{eq:lma-ide-proof-4}
        \ecn{\nabla f_m (\hat{x}^i; \xi_m) - \nabla f_m (x^\ast; \xi_m)} 
        \overset{\eqref{eq:grad_dif_bregman}}{\leq} 2 L \ec{D_{f_m} (\hat{x}^i, x^\ast)}.
    \end{equation}
    Using \eqref{eq:lma-ide-proof-4} and \eqref{eq:lma-ide-proof-3} in \eqref{eq:lma-ide-proof-2} and averaging with respect to $m$,
    \begin{align}
        \label{eq:lma-ide-proof-5}
        \avemm \ecn{g^i_m} 
        &\leq \frac{3L^2}{M} \sum_{m=1}^{M} \ecn{x^i_m - \hat{x}^i} + 6 L D_{f} (\hat{x}^i, x^\ast) + 3 \sigmaf^2 \nonumber \\
        &= 3 L^2 \ec{V^i} + 6 L \ec{D_{f} (\hat{x}^i, x^\ast)} + 3 \sigmaf^2.
    \end{align}
    Using \eqref{eq:lma-ide-proof-5} in \eqref{eq:lma-ide-proof-2},
    \begin{align*}
        \sum_{k=k_p}^{v} \ec{V^k} 
        &\leq \gamma^2 \br{H-1}^2 \sum_{k=k_p}^{v} \ec{3 L^2 V^k + 6 L D_{f} (\hat{x}^k, x^\ast) + 3 \sigmaf^2}.
    \end{align*}
    Noticing that the sum $\sum_{k=k_p}^{v} \ec{V^k}$ appears in both sides, we can rearrange
    \begin{align*}
        \br{1 - 3 \gamma^2 \br{H-1}^2 L^2} \sum_{k=k_p}^{v} \ec{V^k} 
        &\leq 6 L \gamma^2 \br{H-1}^2 \sum_{k=k_p}^{v} \ec{D_{f} (\hat{x}^k, x^\ast)}\\
        &\qquad + 3 \gamma^2 \br{H-1}^2 \sum_{k=k_p}^{v} \sigmaf^2.
    \end{align*}
    Finally using that our choice $\gamma$ implies that $1 - 3 \gamma^2 \br{H-1}^2 L^2 \geq \frac{3}{4}$ we have,
    \begin{align*}
        \sum_{k=k_p}^{v} \ec{V^k} 
        \leq 8 L \gamma^2 \br{H-1}^2 \sum_{k=k_p}^{v} \ec{D_{f} (\hat{x}^k, x^\ast)} + 4 \gamma^2 \br{H-1}^2 \sum_{k=k_p}^{v} \sigmaf^2.
    \end{align*}
\end{proof}

\begin{lemma}[Optimality gap single recursion] 
    \label{lemma:optimality-gap-single-recursion}
    Suppose that Assumptions~\ref{asm:convexity-and-smoothness} and \ref{asm:finite-sum-stochastic-gradients} hold for Algorithm~\ref{alg:local_sgd} with heterogeneous data and with $M \geq 2$. Then for any $\gamma \geq 0$ we have for expectation conditional on $x^k_1, x^k_2, \ldots, x^k_m$,
    \begin{equation}
        \label{eq:9f8gff}
        \ecn{r^{k+1}} \leq \br{1 - \gamma \mu} \sqn{r^k} + \gamma L \br{1 + 2 \gamma L} V^k - 2 \gamma \br{1 - 4 \gamma L} D_{f} (\hat{x}^k, x^\ast) + \frac{4 \gamma^2 \sigmaf^2}{M},
     \end{equation}   
    where $r^{k} \eqdef \hat{x}^k - x^\ast$. In particular, if $\gamma \leq \frac{1}{8L}$, then
    \begin{equation}
        \ecn{r^{k+1}} \leq \br{1 - \gamma \mu} \sqn{r^k} + \frac{5}{4} \gamma L V^k - \frac{\gamma}{2} D_{f} (\hat{x}^k, x^\ast) + \frac{4 \gamma^2 \sigmaf^2}{M},
    \end{equation}
\end{lemma}
\begin{proof}
    First note that $\hat{x}^{k+1} = \hat{x}^k - \gamma g^k$ is always true (regardless of whether or not synchronization happens), hence
    \begin{align*}
        \sqn{\hat{x}^{k+1} - x^\ast} &= \sqn{\hat{x}^k - \gamma g^k - x^\ast} \\
        &= \sqn{\hat{x}^k - x^\ast} + \gamma^2 \sqn{g^k} - 2 \gamma \ev{\hat{x}^k - x^\ast, g^k} \\
        &= \sqn{\hat{x}^k - x^\ast} + \gamma^2 \sqn{g^k} - \frac{2 \gamma}{M} \sum_{m=1}^{M} \ev{\hat{x}^k - x^\ast, g^k_m}.
    \end{align*}
    Taking conditional expectations then using Lemmas~\ref{lemma:average-gradient-bound} and \ref{lemma:inner-product-bound},
    \begin{align*}
        \ecn{r^{k+1}} 
        &\leq \sqn{r^k} + \gamma^2 \ecn{g^k} - \frac{2 \gamma}{M} \sum_{m=1}^{M} \ev{\hat{x}^k - x^\ast, \nabla f_m (x^k_m)} \\
        &\overset{\eqref{eq:lma-average-gradient-bound}}{\leq} \sqn{r^k} + \gamma^2 \br{ 2 L^2 V^k + 8 L D_{f} (\hat{x}^k, x^\ast) + \frac{4 \sigmaf^2}{M}}\\
        &\qquad - \frac{2 \gamma}{M} \sum_{m=1}^{M} \ev{\hat{x}^k - x^\ast, \nabla f_m (x^k_m)}  \\
        &\overset{\eqref{eq:lma-inner-product-bound}}{\leq} \br{1 - \gamma \mu} \sqn{r^k} + \gamma L \br{1 + 2 \gamma L} V^k - 2 \gamma \br{1 - 4 \gamma L} D_{f} (\hat{x}^k, x^\ast)\\
        &\qquad + \frac{4 \gamma^2 \sigmaf^2}{M}.
    \end{align*}
    If $\gamma \leq \frac{1}{8L}$, then $1 - 4 \gamma L \geq \frac{1}{2}$ and $1 + 2 \gamma L \leq \frac{5}{4}$, and hence
    \begin{align*}
        \ecn{r^{k+1}} &\leq \br{1 - \gamma \mu} \sqn{r^k} + \frac{5}{4} \gamma L V^k - \frac{\gamma}{2} D_{f} (\hat{x}^k, x^\ast) + \frac{4 \gamma^2 \sigmaf^2}{M},
    \end{align*}
    as claimed.
\end{proof}

\subsection{Proof of Theorem~\ref{thm:wc-noniid-unbounded-var}}
\begin{proof}
    Start with Lemma~\ref{lemma:optimality-gap-single-recursion} with $\mu = 0$,
    \begin{align*}
        \ecn{r^{k+1}} &\leq \sqn{r^k} + \frac{\gamma}{2} \br{\frac{5}{2} L V^k - D_{f} (\hat{x}^k, x^\ast)} + \frac{4 \gamma^2 \sigmaf^2}{M}.
    \end{align*}
    Taking unconditional expectations and summing up,
    \begin{align}
        \label{eq:thm-lsgd-dd-proof-1}
        \sum_{i=1}^{K} \ecn{r^k} &\leq \sum_{i=0}^{K-1} \ecn{r^k} + \frac{\gamma}{2} \sum_{i=0}^{K-1} \ec{\frac{5}{2} L V^i - D_{f} (\hat{x}^i, x^\ast)} + \sum_{i=0}^{K-1} \frac{4 \gamma^2 \sigmaf^2}{M}.
    \end{align}
    Using that $K = k_p$ for some $p \in \N$, we can decompose the second term by double counting and bound it by Lemma~\ref{lemma:iterate-deviation-epoch},
    \begin{align*}
        \sum_{i=0}^{K-1} \ec{\frac{5}{2} L V^i - D_{f} (\hat{x}^i, x^\ast)} 
        &= \sum_{l=1}^{p} \sum_{i=k_{l-1}}^{k_{l} - 1} \ec{\frac{5}{2} L V^i - D_{f} (\hat{x}^i, x^\ast)} \\
        &\leq \sum_{l=1}^{p} \sum_{i=k_{l-1}}^{k_l - 1} \br{20 L^2 \gamma^2 (H-1)^2 - 1} \ec{D_f (\hat{x}^i, x^\ast}\\
        &\qquad + \sum_{l=1}^{p} \sum_{i=k_{l-1}}^{k_l - 1} 10 L \gamma^2 (H-1)^2 \sigmaf^2.
    \end{align*}
    By assumption on $\gamma$ we have that $20 L^2 \gamma^2 (H-1)^2 - 1 \leq -\frac{1}{2}$, using this and then using double counting again we have,
    \begin{align*}
        \sum_{i=0}^{K-1} \ec{\frac{5}{2} L V^i - D_{f} (\hat{x}^i, x^\ast)} 
        &\leq -\frac{1}{2} \sum_{l=1}^{p} \sum_{i=k_{l-1}}^{k_l - 1} \ec{D_f (\hat{x}^i, x^\ast)} + \sum_{l=1}^{p} \sum_{i=k_{l-1}}^{k_l - 1} 10 L \gamma^2 (H-1)^2 \sigmaf^2 \\
        &= -\frac{1}{2} \sum_{i=0}^{K-1} \ec{D_f (\hat{x}^i, x^\ast)} + \sum_{i=0}^{K-1} 10 L \gamma^2 (H-1)^2 \sigmaf^2.
    \end{align*}
    Using this in \eqref{eq:thm-lsgd-dd-proof-1},
    \begin{align*}
        \sum_{i=1}^{K} \ecn{r^k} 
        &\leq \sum_{i=0}^{K-1} \ecn{r^k} - \frac{\gamma}{4} \sum_{i=0}^{K-1} \ec{D_{f} (\hat{x}^i, x^\ast)}\\
        &\qquad + \sum_{i=0}^{K-1} \br{5 L \gamma^3 (H-1)^2 \sigmaf^2 + \frac{4 \gamma^2 \sigmaf^2}{M}}.
    \end{align*}
    Rearranging, we get
    \begin{align*}
        \frac{\gamma}{4} \sum_{i=0}^{K-1} \ec{D_f (\hat{x}^i, x^\ast)} 
        &\leq \sum_{i=0}^{K-1} \ecn{r^k} - \sum_{i=1}^{K} \ecn{r^k}\\
        &\qquad + \sum_{i=0}^{K-1} \br{5 L \gamma^3 (H-1)^2 \sigmaf^2 + \frac{4 \gamma^2 \sigmaf^2}{M}} \\
        &= \sqn{x^0-x^\ast} - \ecn{r^k} + \sum_{i=0}^{K-1} \br{5 L \gamma^3 (H-1)^2 \sigmaf^2 + \frac{4 \gamma^2 \sigmaf^2}{M}} \\
        &\leq \sqn{x^0-x^\ast} + K \br{5 L \gamma^3 (H-1)^2 \sigmaf^2 + \frac{4 \gamma^2 \sigmaf^2}{M}}.
    \end{align*}
    Dividing both sides by $\gamma K / 4$, we get
    \begin{align*}
        \frac{1}{K} \sum_{i=0}^{K-1} \ec{D_f (\hat{x}^i, x^\ast)} 
        &\leq \frac{4 \sqn{x^0-x^\ast}}{\gamma K} + \frac{20 \gamma \sigmaf^2}{M} + 16 \gamma^2 L (H-1)^2 \sigmaf^2.
    \end{align*}
    Finally, using Jensen's inequality and the convexity of $f$ we get the required claim.
\end{proof}

\section{Extra Experiments}
Figure~\ref{fig:a5a_same_data_001} shows experiments done with identical data and Figure~\ref{fig:mushrooms_different_H} shows experiments done with heterogeneous data in the same setting as described in the main text but with different datasets.

\begin{figure}
	\centering
	\includegraphics[scale=0.32]{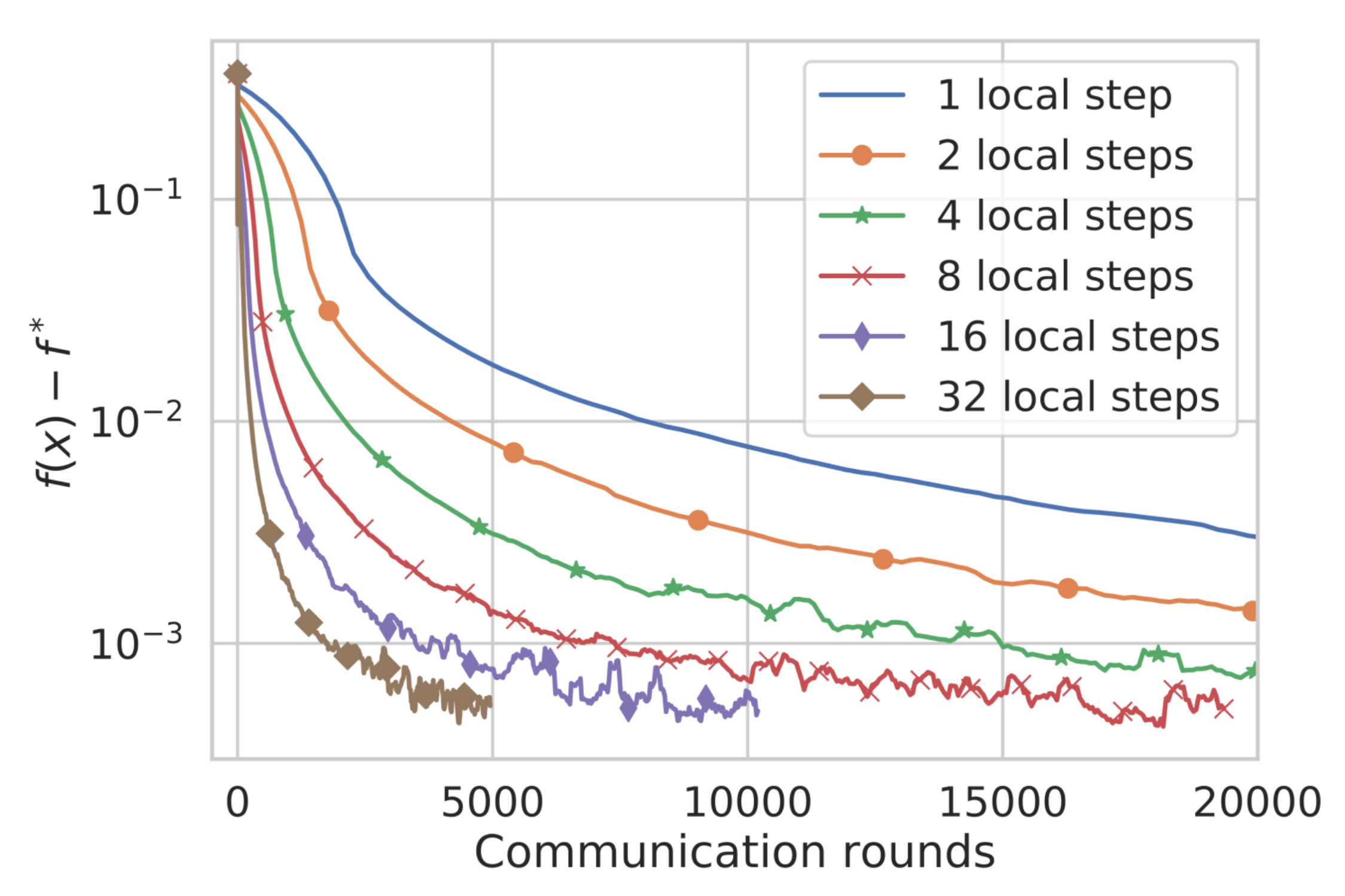}
	\includegraphics[scale=0.32]{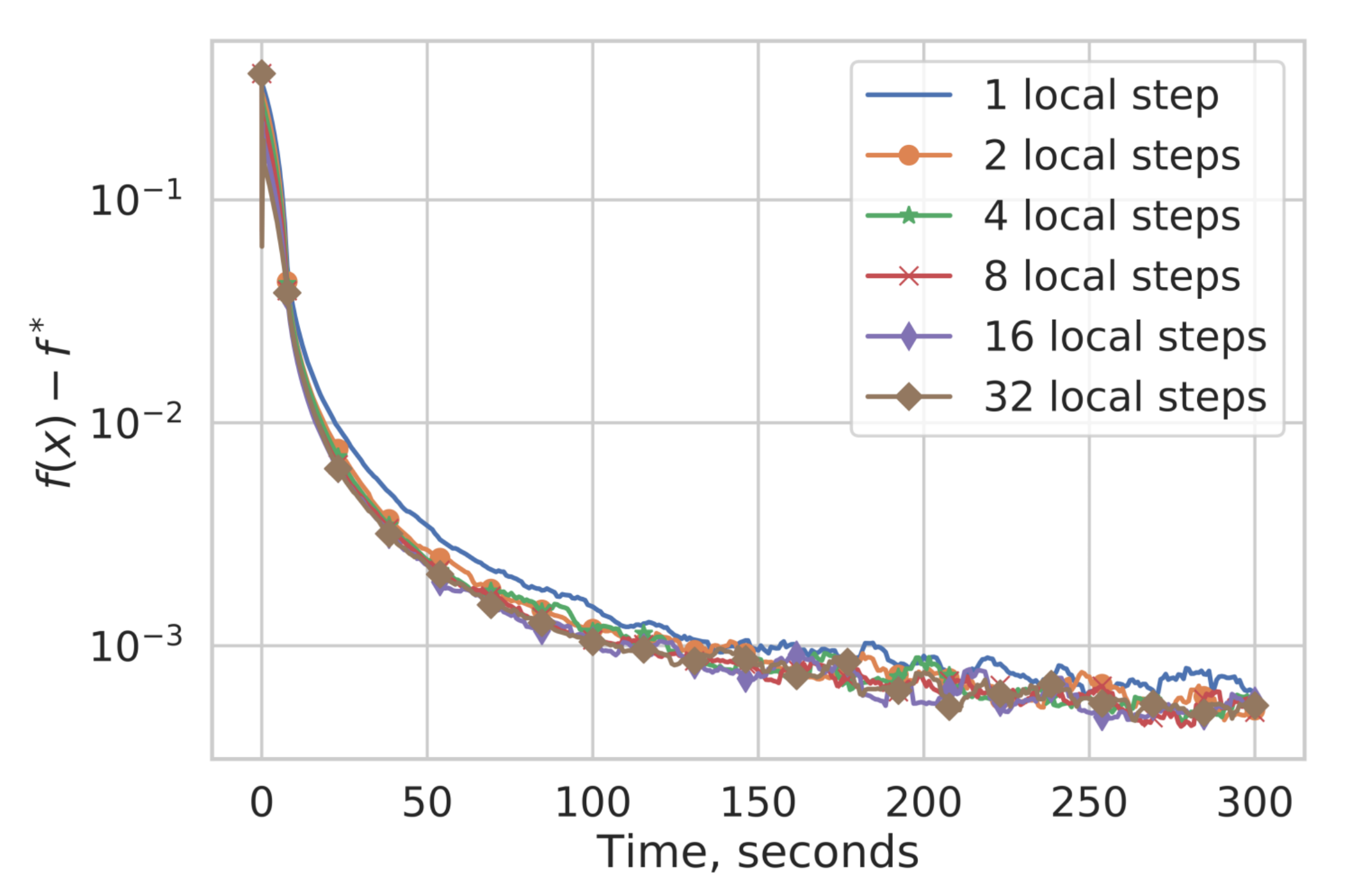}
	\caption{Results on regularized logistic regression with shared data, `a9a' dataset, with stepsize $\frac{0.05}{L}$. With more local iterations, fewer communication rounds are required to get to a neighborhood of the solution.}
	\label{fig:a5a_same_data_001}
\end{figure}

\begin{figure}[!b]
	\centering
	\includegraphics[scale=0.32]{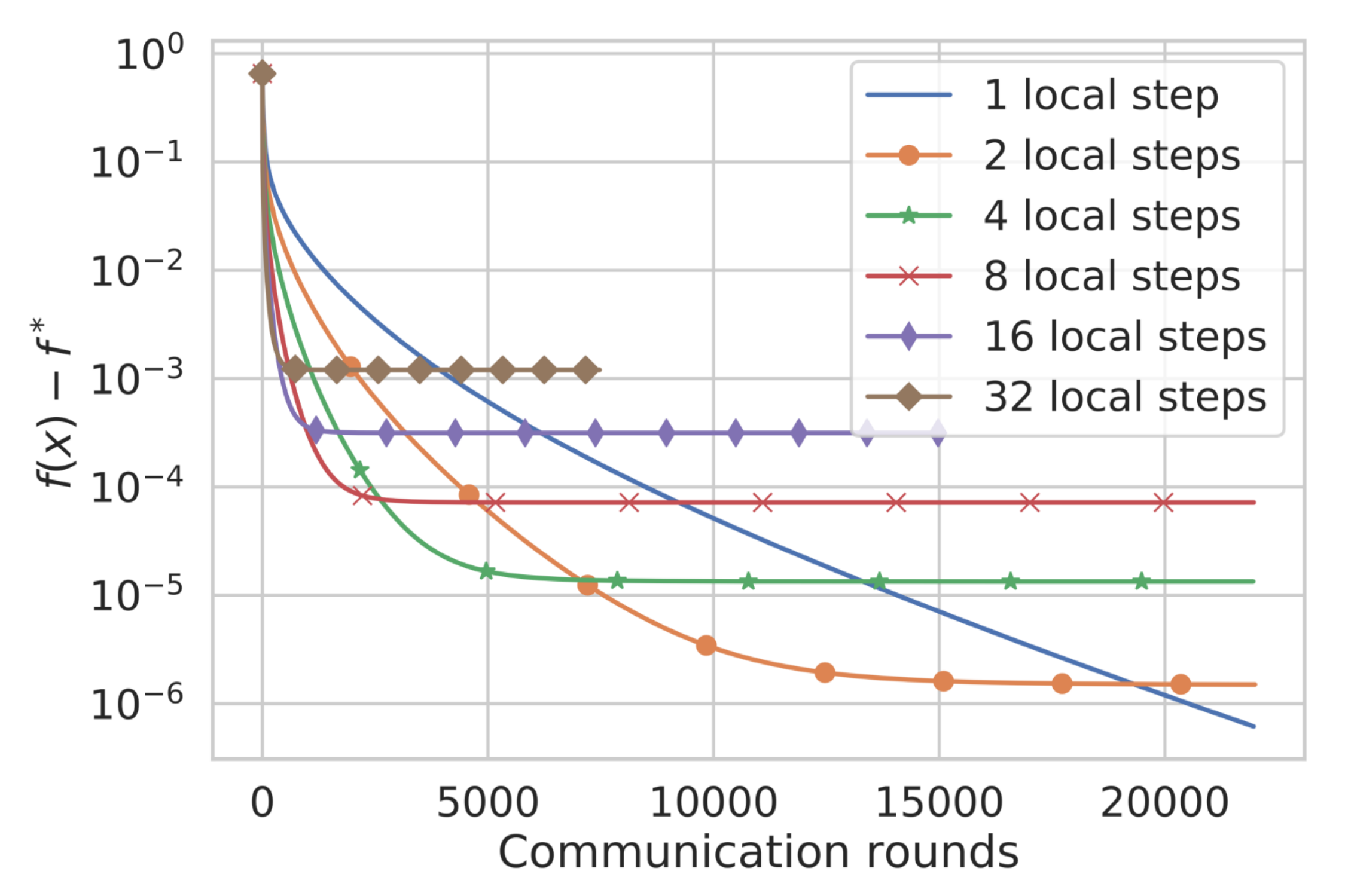}
	\includegraphics[scale=0.32]{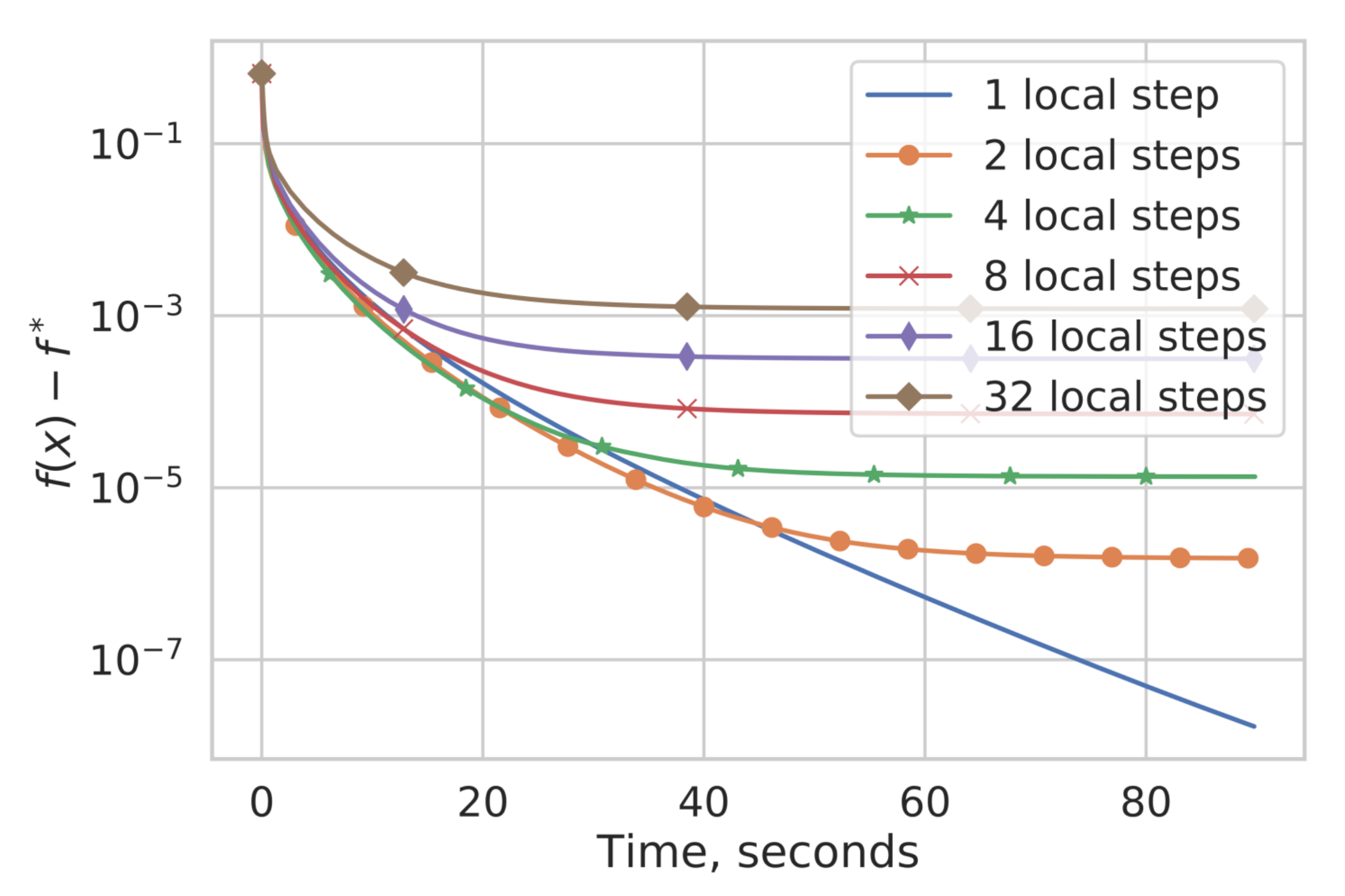}
	\includegraphics[scale=0.32]{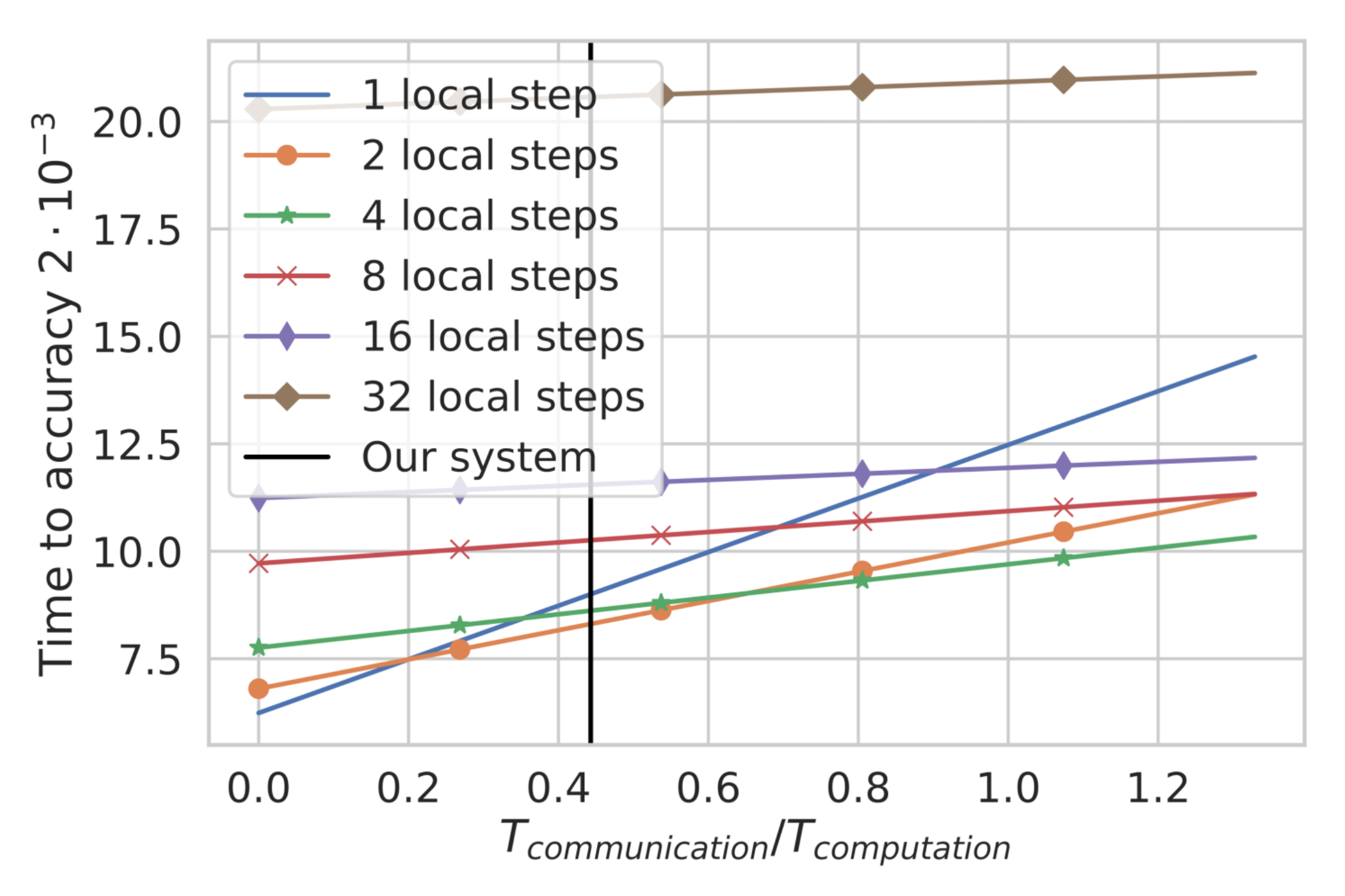}
	\caption{Same experiment as in Figure~\ref{fig:a5a_different_H}, performed on the `mushrooms' dataset.} 
	\label{fig:mushrooms_different_H}
\end{figure}

\section{Discussion of Dieuleveut and Patel (2019)}
\label{sec:patel-discuss}
An analysis of \algname{Local SGD} for identical data under strong convexity, Lipschitzness of $\nabla f$, uniformly bounded variance, and Lipschitzness of $\nabla^2 f$ is given in \cite{Patel19}, where they obtain a similar communication complexity to \cite{Stich2018} without bounded gradients. However, in the proof of their result for general non-quadratic functions (Proposition S20) they make the following assumption, rewritten in our notation:
\[ G = \sup_{p} \br{ 1 + M L_H \gamma \sum_{k=k_{p}}^{k_{p+1} - 1} \sqn{\hat{x}^k - x^\ast} } < \infty,  \]
where $L_H$ is the Lipschitz constant of the Hessian of $f$ (assumed thrice differentiable). Their discussion of $G$ speculates on the behaviour of iterate distances, e.g., saying that if they are bounded, then the guarantee is good. Unfortunately, assuming this quantity bounded implies that gradients are bounded as well, making the improvement over \cite{Stich2018} unclear to us. Furthermore, as $G$ depends on the algorithm's convergence (it is the distance from the optimum evaluated at various points), assuming it is bounded to prove convergence to a compact set results in a possibly circular argument. Since $G$ is also used as an upper bound on $H$ in their analysis, it is not possible to calculate the communication complexity.

\chapter{Appendix for Chapter~\ref{chapter:rr}}
\label{rr_appendix}

\graphicspath{{rr/}}

\section{Additional Experiment Details}
\textbf{Objective properties.} To better correspond to the theoretical setting of our main result, we use $\ell_2$ regularization in every element of the finite-sum. To obtain mini-batches for the \algname{RR}, \algname{SO} and \algname{IG} we permute the dataset and then split it into $n=\lceil \frac{N}{\tau}\rceil$ groups of sizes $\tau, \dotsc, \tau, N-\tau \br{\lceil\frac{N}{\tau}\rceil - 1}$. In other words, the first $n-1$ groups are of size $\tau$ and the remaining samples go to the last group. For \algname{SO} and \algname{IG}, we split the data only once, and for \algname{RR}, we do this at the beginning of each epoch. The permutation of samples used in \algname{IG} is the one in which the datasets are stored online. The smoothness constant of the sum of logistic regression losses admits a closed form expression $L_f=\frac{1}{4N}\|\mA\|^2+\lambda$. The individual losses are $L_{\max}$-smooth with $L_{\max}=\max_{i=1,\dotsc, n} \| a_i\|^2+\lambda$. 

\textbf{Stepsizes.} For all methods in Figure~\ref{fig:conv_plots}, we keep the stepsize equal to $\frac{1}{L}$ for the first $k_0=\lfloor K/40 \rfloor$ iterations, where $K$ is the total number of stochastic steps. This is important to ensure that there is an exponential convergence before the methods reach their convergence neighborhoods~\cite{stich2019unified}. After the initial $k_0$ iterations, the stepsizes used for \algname{RR}, \algname{SO} and \algname{IG} were chosen as $\gamma_{k}=\min\left\{\frac{1}{L}, \frac{3}{\mu \max\{1, k-k_0\}}\right\}$ and for \algname{SGD} as $\gamma_{k}=\min\left\{\frac{1}{L}, \frac{2}{\mu \max\{1, k-k_0\}}\right\}$. Although these stepsizes for \algname{RR} are commonly used in practice~\cite{Bottou2009}, we do not analyze them and leave decreasing-stepsize analysis for future work. We also note that although $L$ is generally not available, it can be estimated using empirically observed gradients \cite{malitsky2019adaptive}. For our experiments, we estimate $L$ of mini-batches of size $\tau$ using the closed-form expressions from Proposition~3.8 in \cite{Gower2019} as $L\le \frac{n(\tau -1)}{\tau(n-1)}L_f + \frac{n-\tau}{\tau(n-1)}L_{\max}$. The confidence intervals in Figure~\ref{fig:conv_plots} are estimated using 20 random seeds.

For the experiments in Figure~\ref{fig:variance}, we estimate the expectation from~\eqref{eq:bregman-div-noise} with 20 permutations, which provides sufficiently stable estimates. In addition, we use $L=L_{\max}$ (instead of using the batch smoothness of \cite{Gower2019}) as the plots in this figure use different mini-batch sizes and we want to isolate the effect of reducing the variance by mini-batching from the effect of changing $L$.

\textbf{\algname{SGD} implementation.} For \algname{SGD}, we used two approaches to mini-batching. In the first, we sampled $\tau$ indices from $\{1,\dotsc, N\}$ and used them to form the mini-batch, where $N$ is the total number of data samples. In the second approach, we permuted the data once and then at each iteration, we only sampled one index $i$ and formed the mini-batch from indices $i, (i+1) \mod N,\dotsc, (i+\tau-1) \mod N$. The latter approach is much more cash-friendly and runs significantly faster, while the iteration convergence was the same in our experiments. Thus, we used the latter option to produce the final plots.

For all plots and methods, we use zero initialization, $x^0=(0,\dotsc, 0)^\top \in \R^d$. We obtain the optimum, $x^\ast$, by running Nesterov's accelerated gradient method until it reaches machine precision. The plots in the right column in Figure~\ref{fig:variance} were obtained by initializing the methods at an intermediate iterate of Nesterov's method, and we found the average, best and worst results by sampling 1,000 permutations.

\section{Notation}

We define the epoch total gradient $g^k$ as
\[ g^k \eqdef \sum_{i=0}^{n-1} \nabla f_{\prm{i}} (x^k_i). \]
We define the variance of the local gradients from their average at a point $x^k$ as
\[ \sigma_{k}^2 \eqdef \frac{1}{n} \sum_{j=1}^{n} \sqn{\nabla f_{j} (x^k) - \nabla f(x^k)}. \]
By $\et{\cdot}$ we denote the expectation conditional on all information prior to iteration $k$, including $x^k$. To avoid issues with the special case $n=1$, we use the convention $0/0=0$.
A summary of the notation used in this chapter is given in Table~\ref{tab:notation}.

\section{A Lemma for Sampling without Replacement}
The following algorithm-independent lemma characterizes the variance of sampling a number of vectors from a finite set of vectors, without replacement. It is a key ingredient in our results on the convergence of the \algname{RR} and \algname{SO} methods.

\begin{lemma}\label{lem:sampling_wo_replacement}
	Let $X_1,\dotsc, X_n\in \R^d$ be fixed vectors, $\overline X\eqdef \frac{1}{n}\sum_{i=1}^n X_i$ be their average and $\sigma^2 \eqdef \frac{1}{n}\sum_{i=1}^n \norm{X_i-\overline X}^2$ be the population variance. Fix any $m\in\{1,\dotsc, n\}$, let $X_{\pi_1}, \dotsc X_{\pi_k}$ be sampled uniformly without replacement from $\{X_1,\dotsc, X_n\}$ and $\overline X_\pi$ be their average. Then, the sample average and variance are given by
	\begin{align}
		\ec{\overline X_\pi}=\overline X, && \ec{\norm{\overline X_{\pi} - \overline X}^2}= \frac{n-m}{m(n-1)}\sigma^2. \label{eq:sampling_wo_replacement}
	\end{align}
\end{lemma}
\begin{proof}
	The first claim follows by linearity of expectation and uniformity of  sampling:
	\[
		\ec{\overline X_\pi} 
		= \frac{1}{m}\sum_{i=1}^m \ec{X_{\pi_i}}
		= \frac{1}{m}\sum_{i=1}^m \overline X
		= \overline X.
	\]
	To prove the second claim, let us first establish that the identity $\mathrm{cov}(X_{\pi_i}, X_{\pi_j})=-\frac{\sigma^2}{n-1}$ holds  for any $i\neq j$. Indeed,
	\begin{align*}
		\mathrm{cov}(X_{\pi_i}, X_{\pi_j})
		&= \ec{ \ev{X_{\pi_i} - \overline X, X_{\pi_j} - \overline X}}
		= \frac{1}{n(n-1)}\sum_{l=1}^n\sum_{p=1,p\neq l}^n\ev{X_l - \overline X, X_p - \overline X} \\
		&= \frac{1}{n(n-1)}\sum_{l=1}^n\sum_{p=1}^n\ev{X_l - \overline X, X_p - \overline X} - \frac{1}{n(n-1)}\sum_{l=1}^n \norm{X_l - \overline X}^2 \\
		&= \frac{1}{n(n-1)}\sum_{l=1}^n \ev{X_l - \overline X, \sum_{p=1}^n(X_p - \overline X)} - \frac{\sigma^2}{n-1} \\
		&=-\frac{\sigma^2}{n-1}.
	\end{align*}
This identity helps us to establish the formula for sample variance:
	\begin{align*}
		\ecn{\overline X_{\pi} - \overline X}
		&= \frac{1}{k^2} \sum_{i=1}^k\sum_{j=1}^k \mathrm{cov}(X_{\pi_i}, X_{\pi_j}) \\
    &= \frac{1}{m^2}\ec{\sum_{i=1}^k \sqn{X_{\pi_i} - \overline X}} + \sum_{i=1}^k\sum_{j=1,j\neq i}^{m} \mathrm{cov}(X_{\pi_i}, X_{\pi_j})  \\
		&=\frac{1}{m^2}\left(m\sigma^2 - m(m-1)\frac{\sigma^2}{n-1}\right)
    = \frac{n-m}{m(n-1)}\sigma^2.
	\end{align*}
\end{proof}

\section{Proofs for Convex Objectives (Sections~\ref{sec:strongly-convex} and~\ref{sec:weakly-convex})}

\subsection{Proof of Proposition~\ref{prop:shuffling-variance-normal-variance-bound}}
\begin{proof}
  Let us start with the upper bound. Fixing any $i$ such that $1 \leq i \leq n-1$, we have $i(n-i)\le \frac{n^2}{4}\le \frac{n(n-1)}{2}$ and using smoothness and Lemma~\ref{lem:sampling_wo_replacement} leads to
	\begin{eqnarray*}
		\ec{D_{f_{\pi_i}}(x^\ast_i, x^\ast)}
		\overset{\eqref{eq:L-smooth-intro}}{\le}  \frac{L}{2}\ec{\|x^\ast_i-x^\ast\|^2} 
		&=& \frac{L}{2}\ec{\sqn{\sum_{k=0}^{i-1}\gamma \nabla f_{\pi_k}(x^\ast)}} \\
		&\overset{\eqref{eq:sampling_wo_replacement}}{=} & \frac{\gamma^2 L i(n-i)}{2(n-1)}\sigma_\star^2 \\
    &\le & \frac{\gamma^2 L n}{4}\sigma_\star^2.
  \end{eqnarray*}
  To obtain the upper bound, it remains to take  maximum with respect to $i$ on both sides and divide by $\gamma$.  To prove the lower bound, we use strong convexity and the fact that  $\max_i i(n-i)\ge \frac{n(n-1)}{4}$ holds for any integer $n$. Together, this leads to
  \[
  	\max_i\ec{D_{f_{\pi_i}}(x^\ast_i, x^\ast)}
		\overset{\eqref{eq:asm-strong-convexity}}{\ge} \max_i\frac{\mu}{2}\ec{\|x^\ast_i-x^\ast\|^2} 
    = \max_i\frac{\gamma^2 \mu i(n-i)}{2(n-1)}\sigma_\star^2 \ge \frac{\gamma^2 \mu n}{8}\sigma_\star^2,
  \]
  as desired.
\end{proof}

\subsection{Proof remainder for Theorem~\ref{thm:all-sc-rr-conv}}
\label{sec:proof-of-thm-1}
\begin{proof}
  We start from \eqref{eq:thm_str_cvx_main-proof-1} proved in the main text:
	\[
		\ec{\|x^k_{i+1}-x^\ast_{i+1}\|^2}
		\le (1-\gamma\mu)\ec{\|x^k_i-x^\ast_i\|^2}+ 2 \gamma^2 \sigmass.
	\]
	Since $x^{k+1}-x^\ast=x^k_n-x^\ast_n$ and $x^k-x^\ast=x^k_0-x^\ast_0$, we can unroll the recursion, obtaining the epoch level recursion
  \[
    \ecn{x^{k+1} - x^\ast} \leq \br{1 - \gamma \mu}^{n} \ecn{x^k - x^\ast} + 2 \gamma^2 \sigmass \br{ \sum_{i=0}^{n-1} \br{1 - \gamma \mu}^{i} }.
  \]
  Unrolling this recursion across $T$ epochs, we obtain
  \begin{equation}
    \label{eq:thm_str_cvx_main-proof-2}
    \ecn{x^{T} - x^\ast} \leq \br{1 - \gamma \mu}^{nT} \sqn{x^0 - x^\ast} + 2 \gamma^2 \sigmass \br{\sum_{i=0}^{n-1} \br{1 - \gamma \mu}^{i}} \Biggl( \sum_{j=0}^{T-1} \br{1 - \gamma \mu}^{nj} \Biggr).
  \end{equation}
 The product of the two sums in \eqref{eq:thm_str_cvx_main-proof-2} can be bounded by  reparameterizing  the summation as follows:
  \begin{align*}
    \br{ \sum_{j=0}^{T-1} \br{1 - \gamma \mu}^{nj} } \br{\sum_{i=0}^{n-1} \br{1 - \gamma \mu}^{i}} &= \sum_{j=0}^{T-1} \sum_{i=0}^{n-1} \br{1 - \gamma \mu}^{nj + i} \\
    &= \sum_{k=0}^{nT-1} \br{1 - \gamma \mu}^{k} \leq \sum_{k=0}^{\infty} \br{1 - \gamma \mu}^{k} = \frac{1}{\gamma \mu}.
  \end{align*}
  Plugging this bound back into \eqref{eq:thm_str_cvx_main-proof-2}, we finally obtain the bound
  \[
    \ecn{x^{T} - x^\ast} \leq \br{1 - \gamma \mu}^{nT} \sqn{x^0 - x^\ast} + 2\gamma\frac{\sigmass}{\mu}.
  \]
\end{proof}

\subsection{Proof of complexity}
In this subsection, we show how we get from Theorem~\ref{thm:all-sc-rr-conv} the complexity for strongly convex functions.
\begin{corollary}
  \label{corr:all-sc-rr-conv}
  Under the same conditions as those in Theorem~\ref{thm:all-sc-rr-conv}, we choose stepsize 
  \[ \gamma = \min \pbr{ \frac{1}{L}, \frac{2}{\mu n T} \log\br{ \frac{\norm{x^0 - x^\ast} \mu T \sqrt{n}}{\sqrt{\kappa} \sigmaesc} } }. \]
The final iterate $x^T$ then satisfies
  \[ \ecn{x^{T} - x^\ast} = \mathcal{\tilde{O}} \br{ \exp\left( - \frac{\mu n T}{L}\right) \sqn{x^0 - x^\ast} + \frac{\kappa \sigmaesc^2}{\mu^2 n T^2} }, \]
  where $\tilde{O}(\cdot)$ denotes ignoring absolute constants and polylogarithmic factors. Thus, in order to obtain error (in squared distance to the optimum) less than $\e$, we require that the total number of iterations $n T$ satisfies
  \[ n T = \tilde{\Omega}\br{ \kappa + \frac{\sqrt{\kappa n} \sigmaesc}{\mu \sqrt{\e}} }. \]
\end{corollary}
\begin{proof}
  Applying \Cref{thm:all-sc-rr-conv}, the final iterate generated by Algorithms~\ref{alg:rr} or~\ref{alg:so} after $T$ epochs satisfies
  \[ \ecn{x^{T} - x^\ast} \leq \br{1 - \gamma \mu}^{n T} \sqn{x^0 - x^\ast} + 2\gamma\frac{\sigmass}{\mu}. \]
  Using \Cref{prop:shuffling-variance-normal-variance-bound} to bound $\sigmass$, we get
  \begin{equation}
    \label{eq:corr-all-sc-proof-1}
    \ecn{x^T - x^\ast} \leq \br{1 - \gamma \mu}^{n T} \sqn{x^0 - x^\ast} + \gamma^2 \kappa n \sigmaesc^2.
  \end{equation}
  We now have two cases: 
  \begin{itemize}[leftmargin=0.15in,itemsep=0.01in]
    \item \textbf{Case 1}: If $\frac{1}{L} \leq \frac{2}{\mu n T} \log\br{\frac{\norm{x^0 - x^\ast} \mu T \sqrt{n}}{\sqrt{\kappa} \sigmaesc} }$, then using $\gamma = \frac{1}{L}$ in \eqref{eq:corr-all-sc-proof-1} we have
    \begin{align*}
      \ecn{x^T - x^\ast} &\leq \br{1 - \frac{\mu}{L}}^{n T} \sqn{x^0 - x^\ast} + \frac{\kappa n \sigmaesc^2}{L^2} \\
      &\leq \br{1 - \frac{\mu}{L}}^{n T} \sqn{x^0 - x^\ast} + \frac{4 \kappa \sigmaesc^2}{\mu^2 n T^2} \log^2 \br{ \frac{\norm{x^0 - x^\ast} \mu T \sqrt{n}}{\sqrt{\kappa} \sigmaesc} }.
    \end{align*}
    Using that $1 - x \leq \exp(-x)$ in the previous inequality, we get
    \begin{equation}
      \label{eq:corr-all-sc-proof-2}
      \ecn{x^T - x^\ast} = \ctO\br{ \exp\br{- \frac{\mu n T}{L}} \sqn{x^0 - x^\ast} + \frac{\kappa \sigmaesc^2}{\mu^2 n T^2} },
    \end{equation}
    where $\ctO(\cdot)$ denotes ignoring polylogarithmic factors and absolute (non-problem specific) constants.
    \item \textbf{Case 2}: If $\frac{2}{\mu n T} \log\br{\frac{\norm{x^0 - x^\ast} \mu T \sqrt{n}}{\sqrt{\kappa} \sigmaesc} } < \frac{1}{L}$, recall that by \Cref{thm:all-sc-rr-conv},
    \begin{equation}
      \label{eq:corr-all-sc-proof-3}
      \ecn{x^T - x^\ast} \leq \br{1 - \gamma \mu}^{n T} \sqn{x^0 - x^\ast} + \gamma^2 \kappa n \sigmaesc^2.
    \end{equation}
    Plugging in $\gamma = \frac{2}{\mu n T} \log\br{\frac{\norm{x^0 - x^\ast} \mu T \sqrt{n}}{\sqrt{\kappa} \sigmaesc} }$,  the first term in \eqref{eq:corr-all-sc-proof-3} satisfies
    \begin{align}
      \br{1 - \gamma \mu}^{n T} \sqn{x^0 - x^\ast} &\leq \exp\br{- \gamma \mu n T} \sqn{x^0 - x^\ast} \nonumber \\
      &= \exp\br{ - 2 \log \br{ \frac{\norm{x^0 - x^\ast} \mu T \sqrt{n}}{\sqrt{\kappa} \sigmaesc} } } \sqn{x^0 - x^\ast} \nonumber \\
      &= \frac{\kappa \sigmaesc^2}{\mu^2 n T^2}.
      \label{eq:corr-all-sc-proof-4}
    \end{align}
    Furthermore, the second term in \eqref{eq:corr-all-sc-proof-3} satisfies
    \begin{equation}
      \gamma^2 \kappa n \sigmaesc^2 = \frac{4 \kappa \sigmaesc^2}{\mu^2 n T^2} \log^2 \br{ \frac{\norm{x^0 - x^\ast} \mu T \sqrt{n}}{\sqrt{\kappa} \sigmaesc} }.
      \label{eq:corr-all-sc-proof-5}
    \end{equation}
Substituting~\eqref{eq:corr-all-sc-proof-4} and~\eqref{eq:corr-all-sc-proof-5} into~\eqref{eq:corr-all-sc-proof-3}, we get
    \begin{equation}
      \label{eq:corr-all-sc-proof-6}
      \ecn{x^T - x^\ast} = \ctO\br{\frac{\kappa \sigmaesc^2}{\mu^2 n T^2} }.
    \end{equation}
    This concludes the second case.
  \end{itemize}
  It remains to take the maximum of~\eqref{eq:corr-all-sc-proof-2} from the first case and~\eqref{eq:corr-all-sc-proof-6} from the second case.
\end{proof}

\subsection{Two lemmas for Theorems~\ref{thm:only-f-sc-rr-conv} and~\ref{thm:weakly-convex-f-rr-conv}}

In order to prove Theorems~\ref{thm:only-f-sc-rr-conv} and~\ref{thm:weakly-convex-f-rr-conv}, it will be useful to define the following quantity.

\begin{definition}
  \label{def:cVt}
  Let $x^k_0, x^k_1, \ldots, x^k_n$ be iterates generated by Algorithms~\ref{alg:rr} or~\ref{alg:so}. We define the forward per-epoch deviation over the $k$-th epoch $\cV^k$ as
  \begin{equation}
    \label{eq:cVt_def}
    \cV^k\eqdef \sum_{i=0}^{n-1} \norm{x^k_i - x^{k+1}}^2.
  \end{equation}
\end{definition}

We will now establish two lemmas. First, we will show that $\cV^k$ can be efficiently upper bounded using Bregman divergences and the variance at the optimum. Subsequently use this bound to establish the convergence of \algname{RR}/\algname{SO}.

\subsubsection{Bounding the forward per-epoch deviation}
\begin{lemma}
  \label{lem:convex_rr_deviation}
	Consider the iterates of \algname{Random Reshuffling} (Algorithm~\ref{alg:rr}) or \algname{Shuffle-Once} (Algorithm~\ref{alg:so}). If the functions $f_1, \ldots, f_n$ are convex and Assumption~\ref{asm:f-smoothness} is satisfied, then 
  \begin{equation}
    \label{eq:lma-cVt-bound}
		\ec{\cV^k}
		\le 4\gamma^2n^2 L\sum_{i=0}^{n-1}\ec{D_{f_{\pi_i}}(x^\ast, x^k_i)}
		+ \frac{1}{2}\gamma^2n^2\sigma_\star^2,
  \end{equation}
  where $\cV^k$ is defined as in Definition~\ref{def:cVt}, and $\sigmaesc^2$ is the variance at the optimum given by $\sigmaesc^2 \eqdef \frac{1}{n} \sum_{i=1}^{n} \sqn{\nabla f_{i} (x^\ast)}$.
\end{lemma}
\begin{proof}
	For any fixed $m\in\{0,\dotsc, n-1\}$, by definition of $x^k_m$ and $x^{k+1}$ we get the decomposition
	\[
		x^k_m - x^{k+1}
		=\gamma\sum_{i=m}^{n-1}\nabla f_{\pi_i}(x^k_i)
		= \gamma\sum_{i=m}^{n-1}(\nabla f_{\pi_i}(x^k_i) - \nabla f_{\pi_i}(x^\ast))+\gamma\sum_{i=m}^{n-1}\nabla f_{\pi_i}(x^\ast).
	\]
	Applying Young's inequality to the sums above yields
	\begin{eqnarray*}
		\|x^k_m - x^{k+1}\|^2
		&\overset{\eqref{eq:sum_sqnorm}}{\le}& 2\gamma^2\sqn{\sum_{i=m}^{n-1}(\nabla f_{\pi_i}(x^k_i) - \nabla f_{\pi_i}(x^\ast))} + 2\gamma^2\sqn{\sum_{i=m}^{n-1}\nabla f_{\pi_i}(x^\ast)} \\
		&\overset{\eqref{eq:sqnorm-sum-bound}}{\le}& 2\gamma^2n\sum_{i=m}^{n-1}\sqn{\nabla f_{\pi_i}(x^k_i) - \nabla f_{\pi_i}(x^\ast)} + 2\gamma^2\sqn{\sum_{i=m}^{n-1}\nabla f_{\pi_i}(x^\ast)}\\
		&\overset{\eqref{eq:grad_dif_bregman}}{\le}& 4\gamma^2Ln\sum_{i=m}^{n-1}D_{f_{\pi_i}}(x^\ast, x^k_i) + 2\gamma^2\sqn{\sum_{i=m}^{n-1}\nabla f_{\pi_i}(x^\ast)} \\
		&\le& 4\gamma^2Ln\sum_{i=0}^{n-1}D_{f_{\pi_i}}(x^\ast, x^k_i) + 2\gamma^2\sqn{\sum_{i=m}^{n-1}\nabla f_{\pi_i}(x^\ast)} .
	\end{eqnarray*}
  Summing up and taking expectations leads to
  \begin{equation}
    \label{eq:cVt-bound-proof-1}
    \sum_{m=0}^{n-1} \ecn{x^k_m - x^{k+1}} \leq 4\gamma^2Ln^2\sum_{i=0}^{n-1}\ec{D_{f_{\pi_i}}(x^\ast, x^k_i)} + 2\gamma^2 \sum_{k=0}^{n-1} \ecn{\sum_{i=k}^{n-1}\nabla f_{\pi_i}(x^\ast)} .
  \end{equation}
  We now bound the second term in the right-hand side of \eqref{eq:cVt-bound-proof-1}. First,  using \Cref{lem:sampling_wo_replacement}, we get
	\begin{align*}
		\ecn{\sum_{i=k}^{n-1}\nabla f_{\pi_i}(x^\ast)}
		&= (n-k)^2\ecn{\frac{1}{n-k}\sum_{i=k}^{n-1}\nabla f_{\pi_i}(x^\ast)} \\
		&= (n-k)^2\frac{k}{(n-k)(n-1)}\sigma_\star^2 \\
		&=\frac{k(n-k)}{n-1}\sigma_\star^2.
	\end{align*}
	Next, by summing this for $k$ from 0 to $n-1$, we obtain
	\[
		\sum_{k=0}^{n-1} \ec{\sqn{\sum_{i=k}^{n-1}\nabla f_{\pi_i}(x^\ast)}} 
		= \sum_{k=0}^{n-1}\frac{k(n-k)}{n-1} \sigmaesc^2
		= \frac{1}{6}n(n+1) \sigmaesc^2 \le \frac{n^2 \sigmaesc^2}{4},
	\]
	where in the last step we also used $n\ge 2$. The result follows.
\end{proof}

\subsubsection{Finding a per-epoch recursion}
\begin{lemma}\label{lem:improved_convex}
	Assume that functions $f_1,\dotsc,f_n$ are  convex and that Assumption~\ref{asm:f-smoothness} is satisfied. If \algname{Random Reshuffling} (Algorithm~\ref{alg:rr}) or \algname{Shuffle-Once} (Algorithm~\ref{alg:so}) is run with a stepsize satisfying $\gamma\le \frac{1}{\sqrt{2}Ln}$, then
	\[
		\ec{\|x^{k+1}-x^\ast\|^2}
		\le \ec{\|x^k-x^\ast\|^2} - 2\gamma n\ec{f(x^{k+1})-f^\star} + \frac{\gamma^3 L n^2\sigma_\star^2}{2}.
	\]
\end{lemma}
\begin{proof}
  Define the sum of gradients used in the $k$-th epoch as $g^{k} \eqdef \sum_{i=0}^{n-1} \nabla f_{\pi_{i}} (x^k_i)$. We will use $g^k$ to relate the iterates $x^k$ and $x^{k+1}$. By definition of $x^{k+1}$, we can write
  \[ x^{k+1} = x^k_{n} = x^k_{n-1} - \gamma \nabla f_{\pi_{n-1}} (x^k_{n-1}) = \cdots = x^k_{0} - \gamma \sum_{i=0}^{n-1} \nabla f_{\pi_{i}} (x^k_i).  \]
Further, since $x^k_0 = x^k$, we see that $x^{k+1} = x^k - \gamma g^k$, which leads to
	\begin{align*}
		\|x^{k}-x^\ast\|^2
		=\|x^{k+1}+\gamma g^k-x^\ast\|^2
		&=\|x^{k+1}-x^\ast\|^2+2\gamma\ev{g^k, x^{k+1}-x^\ast}+\gamma^2\|g^k\|^2 \\
		&\ge \|x^{k+1}-x^\ast\|^2+2\gamma\ev{g^k, x^{k+1}-x^\ast} \\
		&= \|x^{k+1}-x^\ast\|^2+2\gamma\sum_{i=0}^{n-1}\ev{\nabla f_{\pi_i}(x^k_i), x^{k+1}-x^\ast}.
	\end{align*}
	Observe that for any $i$, we have the following decomposition
	\begin{align}
		\ev{\nabla f_{\pi_i}(x^k_i), x^{k+1}-x^\ast}
		&= [f_{\pi_i}(x^{k+1})-f_{\pi_i}(x^\ast)]+ [f_{\pi_i}(x^\ast) - f_{\pi_i}(x^k_i) - \<\nabla f_{\pi_i}(x^k_i), x^\ast - x^{k}_i>] \notag\\
		& \quad - [f_{\pi_i}(x^{k+1}) - f_{\pi_i}(x^k_i)-\<\nabla f_{\pi_i}(x^k_i), x^{k+1}-x^k_i>] \notag\\
		&= [f_{\pi_i}(x^{k+1})-f_{\pi_i}(x^\ast)] + D_{f_{\pi_i}}(x^\ast, x^k_i) - D_{f_{\pi_i}}(x^{k+1}, x^k_i). \label{eq:rr_decomposition}
	\end{align}
	Summing the first quantity in~\eqref{eq:rr_decomposition} over $i$ from $0$ to $n-1$ gives
	\[
		\sum_{i=0}^{n-1}[f_{\pi_i}(x^{k+1})-f_{\pi_i}(x^\ast)]
		= n(f(x^{k+1})-f^\star).
	\]
	Now, we can bound the third term in the decomposition~\eqref{eq:rr_decomposition} using $L$-smoothness as follows:
	\[
		D_{f_{\pi_i}}(x^{k+1}, x^k_i)
		\le \frac{L}{2}\|x^{k+1}-x^k_i\|^2.
	\]
	By summing the right-hand side over $i$ from $0$ to $n-1$ we get the forward deviation over an epoch $\cV^k$, which we bound by Lemma~\ref{lem:convex_rr_deviation} to get
	\[
		\sum_{i=0}^{n-1}\ec{D_{f_{\pi_i}}(x^{k+1}, x^k_i)}
    \overset{\eqref{eq:cVt_def}}{\leq} \frac{L}{2} \ec{\cV^k} 
		\overset{\eqref{eq:lma-cVt-bound}}{\le}  2\gamma^2 L^2n^2\sum_{i=0}^{n-1}\ec{D_{f_{\pi_i}}(x^\ast, x^k_i)}
		+ \frac{\gamma^2 L n^2\sigma_\star^2}{4}.
	\]
	Therefore, we can lower-bound the sum of the second and the third term in~\eqref{eq:rr_decomposition} as
	\begin{align*}
		\sum_{i=0}^{n-1}\ec{D_{f_{\pi_i}}(x^\ast, x^k_i) - D_{f_{\pi_i}}(x^{k+1}, x^k_i)}
		&\ge \sum_{i=0}^{n-1}\ec{D_{f_{\pi_i}}(x^\ast, x^k_i)}\\
		&\qquad  - 2\gamma^2L^2n^2\sum_{i=0}^{n-1}\ec{D_{f_{\pi_i}}(x^\ast, x^k_i)} - \frac{\gamma^2 L n^2\sigma_\star^2}{4} \\
		&\ge (1-2\gamma^2L^2n^2)\sum_{i=0}^{n-1}\ec{D_{f_{\pi_i}}(x^\ast, x^k_i)} - \frac{\gamma^2 L n^2\sigma_\star^2}{4} \\
		&\ge  - \frac{\gamma^2 L n^2\sigma_\star^2}{4},
	\end{align*}
	where in the third inequality we used that $\gamma \leq \frac{1}{\sqrt{2} L n}$ and that $D_{f_{\pi_{i}}} (x^\ast, x^k_i)$ is nonnegative. Plugging this back into the lower-bound on $\|x^k-x^\ast\|^2$ yields
	\[
		\ecn{x^k - x^\ast}
		\ge \ec{\|x^{k+1}-x^\ast\|^2} + 2\gamma n\ec{f(x^{k+1})-f^\star}- \frac{\gamma^3 L n^2\sigma_\star^2}{2}.
	\]
	Rearranging the terms gives the result.
\end{proof}

\subsection{Proof of Theorem~\ref{thm:only-f-sc-rr-conv}}
\begin{proof}
	We can use \Cref{lem:improved_convex} and strong convexity to obtain
	\begin{align*}
		\ecn{x^{k+1}-x^\ast}
		&\le \ecn{x^k-x^\ast} - 2\gamma n\ec{f(x^{k+1})-f^\star} + \frac{\gamma^3 L n^2\sigma_\star^2}{2}\\
		&\overset{\eqref{eq:asm-strong-convexity}}{\le}  \ecn{x^k-x^\ast} - \gamma n\mu\ecn{x^{k+1}-x^\ast} + \frac{\gamma^3 L n^2\sigma_\star^2}{2},
	\end{align*}
	whence
	\begin{align*}
		\ecn{x^{k+1}-x^\ast}
    &\le\frac{1}{1+\gamma\mu n}\left( \ecn{x^k-x^\ast} + \frac{\gamma^3 L n^2\sigma_\star^2}{2} \right) \\
    &= \frac{1}{1 + \gamma \mu n} \ecn{x^k - x^\ast} +  \frac{1}{1 + \gamma \mu n} \frac{\gamma^3 L n^2\sigma_\star^2}{2} \\
		&\le \left(1 - \frac{\gamma\mu n}{2}\right) \ecn{x^k - x^\ast} + \frac{\gamma^3 L n^2 \sigmaesc^2}{2} .
  \end{align*}  
  Recursing for $T$ iterations, we get that the final iterate satisfies
  \begin{align*}
    \ecn{x^T - x^\ast} &\leq \br{ 1 - \frac{\gamma \mu n}{2} }^{T} \sqn{x^0 - x^\ast} + \frac{\gamma^3 L n^2 \sigmaesc^2 }{2} \br{ \sum_{j=0}^{T-1} \br{1 - \frac{\gamma \mu n}{2}}^{j} } \\
    &\leq \br{ 1 - \frac{\gamma \mu n}{2} }^{T} \sqn{x^0 - x^\ast} + \frac{\gamma^3 L n^2 \sigmaesc^2 }{2} \br{ \frac{2}{\gamma \mu n} } \\
    &= \br{ 1 - \frac{\gamma \mu n}{2} }^T \sqn{x^0 - x^\ast} + \gamma^2 \kappa n \sigmaesc^2.
  \end{align*}
\end{proof}

\subsection{Proof of Theorem~\ref{thm:weakly-convex-f-rr-conv}}
\begin{proof}
  We start with \Cref{lem:improved_convex}, which states that the following inequality holds:
  \[ \ec{\sqn{x^{k+1} - x^\ast}} \leq \ecn{x^k - x^\ast} - 2 \gamma n \ec{f(x^{k+1}) - f (x^\ast)} + \frac{\gamma^3 L n^2 \sigmaesc^2}{2} . \]
  Rearranging the result leads to
  \[ 2 \gamma n \ec{f(x^{k+1}) - f(x^\ast)} \leq \ecn{x^{k} - x^\ast} - \ecn{x^{k+1} - x^\ast} +\frac{\gamma^3 L n^2 \sigmaesc^2}{2} .  \]
  Summing these inequalities for $k=0,1,\dots, T-1$ gives
  \begin{align*}
    2 \gamma n \sum_{k=0}^{T-1} \ec{f(x^{k+1}) - f(x^\ast)} 
    &\leq \sum_{k=0}^{T-1} \br{ \ecn{x^k - x^\ast} - \ecn{x^{k+1} - x^\ast} }\\
    &\qquad + \frac{\gamma^3 L n^2 \sigmaesc^2 T}{2}  \\
    &= \sqn{x_{0} - x^\ast} - \ecn{x^{T} - x^\ast} +\frac{\gamma^3 L n^2 \sigmaesc^2 T}{2} \\
    &\leq \sqn{x^0 - x^\ast} + \frac{\gamma^3 L n^2 \sigmaesc^2 T}{2} ,
  \end{align*}
  and dividing both sides by $2 \gamma n T$, we get
  \[
    \frac{1}{T} \sum_{k=0}^{T-1} \ec{f(x^{k+1}) - f(x^\ast)} 
    \leq \frac{\sqn{x^0 - x^\ast}}{2 \gamma n T} + \frac{ \gamma^2 L  n \sigmaesc^2}{4}. 
  \]
  Finally, using convexity of $f$,  the average iterate $\hat{x}_{T} \eqdef \frac{1}{T} \sum_{k=1}^{T} x^k$ satisfies
	\[
		\ec{f(\hat{x}^T) - f(x^\ast)} \leq \frac{1}{T} \sum_{k=1}^{T} \ec{f(x^k)-f(x^\ast)}
    \le \frac{\sqn{x^0 - x^\ast}}{2 \gamma n T} + \frac{ \gamma^2 L  n \sigmaesc^2}{4}. 
	\]
\end{proof}

\subsection{Proof of complexity}
\begin{corollary}
  \label{corr:weakly-convex-f}
  Under the same conditions as Theorem~\ref{thm:weakly-convex-f-rr-conv}, choose the stepsize 
  \[ \gamma = \min \pbr{ \frac{1}{\sqrt{2} L n}, \br{ \frac{\sqn{x^0 - x^\ast}}{L n^2 T \sigmaesc^2} }^{1/3} }. \]
  Then
  \[ \ec{f(\hat{x}^T) - f(x^\ast)} \leq \frac{L \sqn{x^0 - x^\ast}}{\sqrt{2} T} + \frac{3 L^{1/3} \norm{x^0 - x^\ast}^{4/3} \sigmaesc^{2/3}}{4 n^{1/3} T^{2/3}}. \]
We can guarantee $\ec{f(\hat{x}^T) - f(x^\ast)} \leq \e^2$ provided that the total number of iterations satisfies
  \[ Tn \geq  \frac{2\sqn{x^0 - x^\ast} \sqrt{Ln}}{\e^2} \max \pbr{ \sqrt{2 L n}, \frac{\sigmaesc}{\e} }. \]
\end{corollary}
\begin{proof}
  We start with the guarantee of Theorem~\ref{thm:weakly-convex-f-rr-conv}:
  \begin{equation}
    \label{eq:corr-wc-proof-1}
    \ec{f(\hat{x}^T) - f(x^\ast)} \leq \frac{\sqn{x^0 - x^\ast}}{2 \gamma n T} + \frac{\gamma^2 L n \sigmaesc^2}{4} .
  \end{equation}
  We now have two cases depending on the stepsize: 
  \begin{itemize}[leftmargin=0.15in,itemsep=0.01in]
    \item \textbf{Case 1}: If $\gamma = \frac{1}{\sqrt{2} L n} \leq \br{ \frac{\sqn{x^0 - x^\ast}}{L n^2 T \sigmaesc^2} }^{1/3}$, then plugging this $\gamma$ into \eqref{eq:corr-wc-proof-1} gives
    \begin{align}
      \ec{f(\hat{x}^T) - f(x^\ast)} &\leq \frac{L \sqn{x^0 - x^\ast}}{\sqrt{2} T} + \frac{\gamma^2 L n \sigmaesc^2 }{4}  \nonumber \\
      &\leq \frac{L \sqn{x^0 - x^\ast}}{\sqrt{2} T} +  \br{ \frac{\sqn{x^0 - x^\ast}}{L n^2 T \sigmaesc^2} }^{2/3} \frac{Ln \sigmaesc^2}{4}  \nonumber \\
      &= \frac{L \sqn{x^0 - x^\ast}}{\sqrt{2} T} + \frac{L^{1/3} \sigmaesc^{2/3} \norm{x^0 - x^\ast}^{4/3}}{4 n^{1/3} T^{2/3}}.
      \label{eq:corr-wc-proof-2}
    \end{align}
    \item \textbf{Case 2}: If $\gamma = \br{ \frac{\sqn{x^0 - x^\ast}}{L n^2 T \sigmaesc^2} }^{1/3} \leq \frac{1}{\sqrt{2} L n}$, then plugging this $\gamma$ into \eqref{eq:corr-wc-proof-1} gives
    \begin{align}
      \ec{f(\hat{x}^T) - f(x^\ast)} &\leq  \frac{L^{1/3} \norm{x^0 - x^\ast}^{4/3} \sigmaesc^{2/3}}{2 n^{1/3} T^{2/3}} + \frac{L^{1/3} \sigmaesc^{2/3} \norm{x^0 - x^\ast}^{4/3}}{4 n^{1/3} T^{2/3}} \nonumber \\
      &= \frac{3 L^{1/3} \norm{x^0 - x^\ast}^{4/3} \sigmaesc^{2/3}}{4 n^{1/3} T^{2/3}}.
      \label{eq:corr-wc-proof-3}
    \end{align}
  \end{itemize}
  Combining \eqref{eq:corr-wc-proof-2} and \eqref{eq:corr-wc-proof-3}, we see that in both cases we have
  \[ \ec{f(\hat{x}^T) - f(x^\ast)} \leq \frac{L \sqn{x^0 - x^\ast}}{\sqrt{2} T} + \frac{3 L^{1/3} \norm{x^0 - x^\ast}^{4/3} \sigmaesc^{2/3}}{4 n^{1/3} T^{2/3}}. \]
  Translating this to sample complexity, we can guarantee that $\ec{f(\hat{x}^T) - f(x^\ast)} \leq \e^2$ provided
  \[
    n T \geq \frac{2 \sqn{x^0 - x^\ast} \sqrt{L n}}{\e^2} \max \pbr{ \sqrt{L n}, \frac{\sigmaesc}{\e} }. 
  \]
\end{proof}

\section{Proofs for Non-Convex Objectives (Section~\ref{sec:non-convex})}

\subsection{Proof of Proposition~\ref{prop:2nd-moment-bound}}
\begin{proof}
  This proposition is a special case of Lemma~3 in \cite{Khaled2020} and we prove it here for completeness. Let $x \in \R^d$. We start with \eqref{eq:grad-bound} (which does not require convexity) applied to each $f_i$:
  \[ \sqn{\nabla f_{i} (x)} \leq 2 L \br{f_{i} (x) - f_i^\ast}. \]
  Averaging, we derive
  \begin{align*}
    \frac{1}{n} \sum_{i=1}^{n} \sqn{\nabla f_{i} (x)} &\leq 2 L \br{f(x) - \frac{1}{n} \sum_{i=1}^{n} f_i^\ast} \\
    &= 2 L \br{f(x) - f^\ast} + 2 L \br{f^\ast - \frac{1}{n} \sum_{i=1}^{n} f_i^\ast}.
  \end{align*}
  Note that because $f^\ast$ is the infimum of $f(\cdot)$ and $\frac{1}{n} \sum_{i=1}^{n} f_i^\ast$ is a lower bound on $f$ then $f^\ast - \frac{1}{n} \sum_{i=1}^{n} f_i^\ast \geq 0$. We may now use the variance decomposition
  \begin{eqnarray*}
    \frac{1}{n} \sum_{i=1}^{n} \sqn{\nabla f_{i} (x) - \nabla f(x)} &\overset{\eqref{eq:variance_m}}{=}& \frac{1}{n} \sum_{i=1}^{n} \sqn{\nabla f_i (x)} - \sqn{\nabla f(x)} \\
    &\leq& \frac{1}{n} \sum_{i=1}^{n} \sqn{\nabla f_{i} (x)} \\
    &\leq& 2 L \br{f(x) - f^\ast} + 2 L \br{f^\ast - \frac{1}{n} \sum_{i=1}^{n} f_i^\ast}.
  \end{eqnarray*}
  It follows that Assumption~\ref{asm:2nd-moment} holds with $A = L$ and $B^2 = 2 L \br{f^\ast - \frac{1}{n} \sum_{i=1}^{n} f_i^\ast}$.
\end{proof}

\subsection{Finding a per-epoch recursion}
For this subsection and the rest of this section, we need to define the following quantity:
\begin{definition}
  \label{def:vt}
  For Algorithm~\ref{alg:rr} we define the backward per-epoch deviation at timestep $k$ by
  \[ V^k \eqdef \frac{1}{n} \sum_{i=1}^{n} \sqn{x^k_i - x^k}. \]
\end{definition}
We will study the convergence of Algorithm~\ref{alg:rr} for non-convex objectives as follows: we first derive a per-epoch recursion that involves $V^k$ in Lemma~\ref{lemma:epoch-recursion-non-convex}, then we show that $V^k$ can be bounded using smoothness and probability theory in Lemma~\ref{lemma:vt_rr}, and finally combine these two to prove Theorem~\ref{thm:rr-nonconvex}.

\begin{lemma}
  \label{lemma:epoch-recursion-non-convex}
  Suppose that Assumption~\ref{asm:f-smoothness} holds. Then for iterates $x^k$ generated by Algorithm~\ref{alg:rr} with stepsize $\gamma \leq \frac{1}{L n}$, we have
  \begin{equation}
    \label{eq:epoch-recursion-non-convex}
    f(x^{k+1}) \leq f(x^k) - \frac{\gamma n}{2} \sqn{\nabla f(x^k)} + \frac{\gamma L^2}{2} V^k,
  \end{equation}
  where $V^k$ is defined as in \Cref{def:vt}.
\end{lemma}
\begin{proof}
Our approach for establishing this lemma is similar to that of \cite[Theorem~1]{Nguyen2020}, which we became aware of in the course of preparing this manuscript. Recall that $x^{k+1} = x^k - \gamma g^k$, where $g^k = \sum_{i=0}^{n-1} \nabla f_{\pi_{i}} (x^k_i)$. Using $L$-smoothness of $f$, we get
  \begin{eqnarray}
    f(x^{k+1}) &\overset{\eqref{eq:L-smooth-intro}}{\leq}& f(x^{k}) + \ev{\nabla f(x^k), x^{k+1} - x^k} + \frac{L}{2} \sqn{x^{k+1} - x^k} \nonumber \\
    &=& f(x^k) - \gamma n \ev{\nabla f(x^k), \frac{g^k}{n}} + \frac{ \gamma^2 L n^2}{2} \sqn{ \frac{g^k}{n} } \nonumber \\
    &\overset{\eqref{eq:square-decompos}}{=}& f(x^k) - \frac{\gamma n}{2} \br{ \sqn{\nabla f(x^k)} + \sqn{\frac{g^k}{n}} - \sqn{ \nabla f(x^k) - \frac{g^k}{n} } } + \frac{\gamma^2 L  n^2}{2} \sqn{ \frac{g^k}{n} } \nonumber \\
    &=& f(x^k) - \frac{\gamma n}{2} \sqn{\nabla f(x^k)} - \frac{\gamma n}{2} \br{1 - L \gamma n} \sqn{\frac{g^k}{n}} \notag \\
    && + \frac{\gamma n}{2} \sqn{\nabla f(x^k) - \frac{g^k}{n}}.
    \label{eq:ernc-1}
  \end{eqnarray}
  By assumption, we have $\gamma \leq \frac{1}{Ln}$, and hence $1 - L \gamma n \geq 0$. Using this in \eqref{eq:ernc-1}, we get
  \begin{equation}
    \label{eq:ernc-2}
    f(x^{k+1}) \leq f(x^k) - \frac{\gamma n}{2} \sqn{\nabla f(x^k)} + \frac{\gamma n}{2} \sqn{\nabla f(x^k) - \frac{g^k}{n}}. 
  \end{equation}
  For the last term in \eqref{eq:ernc-2}, we note
  \begin{eqnarray}
    \sqn{\nabla f(x^k) - \frac{g^k}{n}} &=& \sqn{ \frac{1}{n} \sum_{i=0}^{n-1} \left [ \nabla f_{\prm{i}} (x^k) - \nabla f_{\prm{i}} (x^k_i) \right ] } \nonumber \\
    &\overset{\eqref{eq:sqnorm-jensen}}{\leq}&  \frac{1}{n} \sum_{i=0}^{n-1} \sqn{\nabla f_{\prm{i}} (x^k) - \nabla f_{\prm{i}} (x^k_i) } \nonumber \\
    &\overset{\eqref{eq:nabla-Lip}}{\leq}&  \frac{1}{n} \sum_{i=0}^{n-1} L^2 \sqn{x^k - x^k_i} = \frac{L^2}{n} V^k.
    \label{eq:ernc-3}
  \end{eqnarray}
  Plugging in~\eqref{eq:ernc-3} into~\eqref{eq:ernc-2} yields the lemma's claim.
\end{proof}

\subsection{Bounding the backward per-epoch deviation}
\begin{lemma}
  \label{lemma:vt_rr}
  Suppose that Assumption~\ref{asm:f-smoothness} holds (with each $f_{i}$ possibly non-convex) and that Algorithm~\ref{alg:rr} is used with a stepsize $\gamma \le \frac{1}{2 L n}$. Then
  \begin{equation}
    \label{eq:vt_bound_rr}
		\et{V^k}
		\le \gamma^2n^3\sqn{\nabla f(x^k)} + \gamma^2n^2\sigma_{k}^2,
	\end{equation}
  where $V^k$ is defined as in \Cref{def:vt} and $\sigma_k^2 \eqdef \frac{1}{n} \sum_{j=1}^{n} \sqn{\nabla f_{j} (x^k) - \nabla f(x^k)}$. 
\end{lemma}
\begin{proof}
  Let us fix any $m\in[1, n-1]$ and find an upper bound for $\et{\norm{x^k_m - x^k}^2}$. First, note that
	\[
		x^k_m = x^k - \gamma\sum_{i=0}^{m-1}\nabla f_{\pi_i}(x^k_i).
	\]
	Therefore, by Young's inequality, Jensen's inequality and gradient Lipschitzness
	\begin{eqnarray*}
		&&\et{\|x^k_m - x^k\|^2}\\
		&=& \gamma^2\et{\left\lVert\sum_{i=0}^{m-1}\nabla f_{\pi_i}(x^k_i) \right\rVert^2} \\
		&\overset{\eqref{eq:sum_sqnorm}}{\le}& 2\gamma^2\et{\left\lVert\sum_{i=0}^{m-1}\left(\nabla f_{\pi_i}(x^k_i)-\nabla f_{\pi_i}(x^k)\right)  \right\rVert^2} + 2\gamma^2\et{\left\lVert\sum_{i=0}^{m-1}\nabla f_{\pi_i}(x^k) \right\rVert^2} \\
		&\overset{\eqref{eq:sqnorm-sum-bound}}{\le}& 2\gamma^2k\sum_{i=0}^{m-1}\et{\norm{\nabla f_{\pi_i}(x^k_i)-\nabla f_{\pi_i}(x^k)}^2} + 2\gamma^2\et{\left\lVert\sum_{i=0}^{m-1}\nabla f_{\pi_i}(x^k) \right\rVert^2} \\
		&\overset{\eqref{eq:nabla-Lip}}{\le}& 2\gamma^2 L^2 k\sum_{i=0}^{m-1}\et{\|x^k_i-x^k\|^2} + 2\gamma^2\et{\left\lVert\sum_{i=0}^{m-1}\nabla f_{\pi_i}(x^k) \right\rVert^2}.
	\end{eqnarray*}
	Let us bound the second term. For any $i$ we have $\et{\nabla f_{\pi_i}(x^k)}=\nabla f(x^k)$, so using \Cref{lem:sampling_wo_replacement} (with vectors $\nabla f_{\pi_0} (x^k), \nabla f_{\pi_{1}} (x^k), \ldots, \nabla f_{\pi_{k-1}} (x^k)$) we obtain
	\begin{eqnarray*}
		\et{\left\lVert\sum_{i=0}^{m-1}\nabla f_{\pi_i}(x^k) \right\rVert^2}
		&\overset{\eqref{eq:variance_m}}{=}& m^2\sqn{\nabla f(x^k)} + m^2\et{\left\lVert\frac{1}{m}\sum_{i=0}^{m-1}  (\nabla f_{\pi_i}(x^k) - \nabla f(x^k))\right\rVert^2} \\
		&\overset{\eqref{eq:sampling_wo_replacement}}{\le}& m^2\sqn{\nabla f(x^k)}+ \frac{m(n-m)}{n-1}\sigma_k^2.
	\end{eqnarray*}
	where $\sigma_k^2 \eqdef \frac{1}{n} \sum_{j=1}^{n} \sqn{\nabla f_{j} (x^k) - \nabla f(x^k)}$. Combining the produced bounds yields 
	\begin{align*}
		\et{\norm{x^k_m - x^k}^2}
		&\le 2\gamma^2 L^2 k \sum_{i=0}^{m-1}\et{\norm{x^k_i-x^k}^2} + 2\gamma^2 k^2\sqn{\nabla f(x^k)} + 2\gamma^2\frac{m(n-m)}{n-1}\sigma_k^2 \\
		&\le 2\gamma^2 L^2 k \ec{V^k} + 2\gamma^2 k^2\sqn{\nabla f(x^k)} + 2\gamma^2\frac{m(n-m)}{n-1}\sigma_k^2,
	\end{align*}
	whence
	\begin{align*}
		\ec{V^k}
		& = \sum_{m=0}^{n-1} \et{\|x^k_m - x^k\|^2} \\
    & \le \gamma^2 L^2 n(n-1) \ec{V^k} + \frac{1}{3}\gamma^2(n-1)n(2n-1)\sqn{\nabla f(x^k)} 
    + \frac{1}{3}\gamma^2 n(n+1)\sigma_k^2.
  \end{align*}
  Since $\ec{V^k}$ appears in both sides of the equation, we rearrange and use that $\gamma\le \frac{1}{2Ln}$ by assumption, which leads to
	\begin{align*}
		\ec{V^k} 
		&\le  \frac{4}{3}(1 - \gamma^2L^2n(n-1)) \ec{V^k} \\
		&\le  \frac{4}{9}\gamma^2(n-1)n(2n-1)\sqn{\nabla f(x^k)} + \frac{4}{9}\gamma^2 n(n+1)\sigma_k^2 \\
    &\le  \gamma^2n^3\sqn{\nabla f(x^k)} + \gamma^2n^2\sigma_k^2.
	\end{align*}
\end{proof}

\subsection{A lemma for solving the non-convex recursion}
\begin{lemma}
  \label{lemma:noncvx-recursion-solution}
  Suppose that there exist constants $a, b, c \geq 0$ and nonnegative sequences $(s_{k})_{k=0}^{T}, (q_{k})_{k=0}^{T}$ such that for any $k$ satisfying $0 \leq k \leq T$ we have the recursion
  \begin{equation}
    \label{eq:nonconvex-recursion-init}
    s_{k+1} 
    \leq \br{1 + a} s_{k} - b q_{k} + c.
  \end{equation}
  Then, the following holds:
  \begin{equation}
    \label{eq:nonconvex-recursion-soln}
    \min_{k=0, \ldots, T-1} q_{k} \leq \frac{\br{1+a}^T}{b T} s_{0} + \frac{c}{b}.
  \end{equation}
\end{lemma}
\begin{proof}
  The first part of the proof (for $a > 0$) is a distillation of the recursion solution in Lemma~2 of \cite{Khaled2020} and we closely follow their proof. Define 
  \[ 
  	w_{k} \eqdef \frac{1}{\br{1+a}^{k+1}}. 
  \]
  Note that $w_{k} \br{1+a} = w_{k-1}$ for all $k$. Multiplying both sides of \eqref{eq:nonconvex-recursion-init} by $w_{k}$,
  \[
    w_{k} s_{k+1} \leq \br{1 + a} w_k s_k - b w_k q_k + c w_k = w_{k-1} s_{k} - b w_{k} q_k + c w_k.
  \]
  Rearranging, we get
$
    b w_{k} q_{k} \leq w_{k-1} s_{k} - w_{k} s_{k+1} + c w_{k}.
$
  Summing up as $k$ varies from $0$ to $T-1$ and noting that the sum telescopes leads to
  \begin{align*}
    \sum_{k=0}^{T-1} b w_{k} q_{k} &\leq \sum_{k=0}^{T-1} \br{w_{k-1} s_{k} - w_{k} s_{k+1}} + c \sum_{k=0}^{T-1} w_{k} \\
    &= w_{0} s_{0} - w_{T-1} s_{T} + c \sum_{k=0}^{T-1} w_{k} \\
    &\leq w_{0} s_{0} + c \sum_{k=0}^{T-1} w_{k}.
  \end{align*}
  Let $W_{T} = \sum_{k=0}^{T-1} w_{k}$. Dividing both sides by $W_{T}$, we get
  \begin{align}
    \frac{1}{W_{T}} \sum_{k=0}^{T-1} b w_{k} q_{k} \leq \frac{w_0 s_0}{W_T} + c.
    \label{eq:nc-rec-1}
  \end{align}
  Note that the left-hand side of \eqref{eq:nc-rec-1} satisfies
  \begin{equation}
    b \min_{k=0, \ldots, T-1} q_k \leq \frac{1}{W_T} \sum_{k=0}^{T-1} b w_k q_k.
    \label{eq:nc-rec-2}
  \end{equation}
 For the right-hand side of \eqref{eq:nc-rec-1}, we have
  \begin{equation}
    W_{T} = \sum_{k=0}^{T-1} w_{k} \geq T \min_{k=0, \ldots, T-1} w_{k} = T w_{T-1} = \frac{T}{\br{1+a}^{T}}.
    \label{eq:nc-rec-3}
  \end{equation}
  Substituting with \eqref{eq:nc-rec-3} in \eqref{eq:nc-rec-2} and dividing both sides by $b$, we finally get
  \[
    \min_{k=0, \ldots, T-1} q_{k} \leq \frac{\br{1 + a}^T}{bT} s_0 + \frac{c}{b}.
  \]
\end{proof}

\subsection{Proof of Theorem~\ref{thm:rr-nonconvex}}
\begin{proof}
  \textbf{Without PL.} Taking expectation in \Cref{lemma:epoch-recursion-non-convex} and then using \Cref{lemma:vt_rr}, we have that for any $t \in \{ 0, 1, \ldots, T-1 \}$,
  \begin{eqnarray*}
    \et{f(x^{k+1})} &\overset{\eqref{eq:epoch-recursion-non-convex}}{\leq} & f(x^k) - \frac{\gamma n}{2} \sqn{\nabla f(x^k)} + \frac{\gamma L^2}{2} \et{V^k} \\
    &\overset{\eqref{eq:vt_bound_rr}}{\leq} & f(x^k) - \frac{\gamma n}{2} \sqn{\nabla f(x^k)} + \frac{\gamma L^2}{2} \br{\gamma^2 n^3 \sqn{\nabla f(x^k)} + \gamma^2 n^2 \sigma_{k}^2} \\
    &= & f(x^k) - \frac{\gamma n}{2} \br{1 - \gamma^2 L^2 n^2} \sqn{\nabla f(x^k)} + \frac{ \gamma^3 L^2 n^2 \sigma_{k}^2}{2} .
  \end{eqnarray*}
  Let $\delta_{k} \eqdef f(x^k) - f^\ast$. Adding $-f^\ast$ to both sides and using Assumption~\ref{asm:2nd-moment},
  \begin{align*}
    \et{\delta_{k+1}}  & \leq  \delta_{k} - \frac{\gamma n}{2} \br{1 - \gamma^2 L^2 n^2} \sqn{\nabla f(x^k)} + \frac{\gamma^3 L^2  n^2 \sigma_{k}^2}{2}  \\
    & \leq  \br{1 + \gamma^3  A L^2 n^2} \delta_{k} - \frac{\gamma n}{2} \br{1 - \gamma^2 L^2 n^2} \sqn{\nabla f(x^k)} + \frac{ \gamma^3 L^2 n^2 B^2}{2} .
  \end{align*}
  Taking unconditional expectations in the last inequality and using that by assumption on $\gamma$ we have $1 - \gamma^2 L^2 n^2 \geq \frac{1}{2}$, we get the estimate
  \begin{equation}
    \label{eq:thm-rr-nc-1}
    \ec{\delta_{k+1}} \leq \br{1 + \gamma^3 A L^2  n^2} \ec{\delta_{k}} - \frac{\gamma n}{4} \ecn{\nabla f(x^k)} + \frac{ \gamma^3 L^2 n^2 B^2}{2} .
  \end{equation}
  Comparing \eqref{eq:nonconvex-recursion-init} with \eqref{eq:thm-rr-nc-1} verifies that the conditions of \Cref{lemma:noncvx-recursion-solution} are readily satisfied. Applying the lemma, we get
  \[
    \min_{k=0, \ldots, T-1} \ecn{\nabla f(x^k)} \leq \frac{4 \br{1 + \gamma^3 A L^2 n^2}^{T}}{\gamma n T} \br{f(x^0) - f^\ast} + 2 \gamma^2 L^2  n B^2.
  \]
  Using that $1 + x \leq \exp(x)$ and that the stepsize $\gamma$ satisfies $\gamma \leq \br{A L^2 n^2 T}^{-1/3}$, we have
  \[ \br{1 + \gamma^3 A L^2  n^2}^{T} \leq \exp\br{\gamma^3 A L^2  n^2 T} \leq \exp\br{1} \leq 3. \]
  Using this in the previous bound, we finally obtain
  \[
    \min_{k=0, \ldots, T-1} \ecn{\nabla f(x^k)} \leq \frac{12 \br{f(x^0) - f^\ast}}{\gamma n T} + 2  \gamma^2 L^2 n B^2.
  \]
  
  \textbf{With PL.} Now we additionally assume that $A=0$ and that $\frac{1}{2}\|\nabla f(x)\|^2\ge \mu(f(x)-f^\star)$. Then, \eqref{eq:thm-rr-nc-1} yields
  \begin{align*}
  	\ec{\delta_{k+1}}
  	&\leq \ec{\delta_{k}} - \frac{\gamma n}{4} \ecn{\nabla f(x^k)} + \frac{ \gamma^3 L^2 n^2 B^2}{2} \\
  	&\le \ec{\delta_{k}} - \frac{\gamma \mu n}{2} \ec{f(x^k)-f^\star} + \frac{ \gamma^3 L^2 n^2 B^2}{2} \\
  	&= \br{1-\frac{\gamma\mu n}{2}}\ec{\delta_k} +  \frac{ \gamma^3 L^2 n^2 B^2}{2} .
  \end{align*}
  As in the proof of Theorem~\ref{thm:only-f-sc-rr-conv}, we recurse this bound to $x^0$:
  \begin{align*}
	  \ec{\delta_{T}}
	  &\le \br{1-\frac{\gamma\mu n}{2}}^T\delta_0 + \frac{ \gamma^3 L^2 n^2 B^2}{2} \sum_{j=0}^{T-1}\br{1-\frac{\gamma\mu n}{2}}^j \\
	  &\le \br{1-\frac{\gamma\mu n}{2}}^T\delta_0 + \frac{ \gamma^3 L^2 n^2 B^2}{2} \frac{2}{\gamma\mu n} \\
	  &= \br{1-\frac{\gamma\mu n}{2}}^T\delta_0 + \gamma^2\kappa L n B^2.
  \end{align*}
\end{proof}

\subsection{Proof of complexity}
\begin{corollary}
  \label{corr:rr-nonconvex}
  Choose the stepsize $\gamma$ as
  \[ \gamma = \min \pbr{ \frac{1}{2 L n}, \frac{1}{A^{1/3} L^{2/3} n^{2/3} T^{1/3} }, \frac{\e}{2 L \sqrt{n} B} }. \]
  Then the minimum gradient norm satisfies $$\min \limits_{k=0, \ldots, T-1} \ecn{ \nabla f(x^k)} \leq \e^2$$ provided the total number of iterations satisfies
  \[ Tn \geq \frac{48 \delta_{0} L \sqrt{n}}{\e^2} \max \pbr{ \sqrt{n}, \frac{\sqrt{6 \delta_{0} A}}{\e}, \frac{B}{\e} }.  \]
\end{corollary}
\begin{proof}
  From \Cref{thm:rr-nonconvex}
  \[ \min_{k=0, \ldots, T-1} \ecn{\nabla f(x^k)} \leq \frac{12 \br{f(x^0) - f^\ast}}{\gamma n T} + 2 \gamma^2 L^2  n B^2.  \]
  Note that by condition on the stepsize $2 L^2 \gamma^2 n B^2 \leq \e^2/2$, hence 
  \[ \min_{k=0, \ldots, T-1} \ecn{\nabla f(x^k)} \leq \frac{12 \br{f(x^0) - f^\ast}}{\gamma n T} + \frac{\e^2}{2}. \]
  Thus, to make the squared gradient norm smaller than $\e^2$ we require
  \[ \frac{12 \br{f(x^0) - f^\ast}}{\gamma n T} \leq \frac{\e^2}{2}, \]
  or equivalently
  \begin{equation}
    \label{corr:rr-nonconvex-proof-1}
    n T \geq \frac{24 \br{f(x^0) - f^\ast}}{\e^2 \gamma} = \frac{24 \delta_{0}}{\e^2} \max \pbr{ 2 L n, \br{A L^2 n^2 T}^{1/3}, \frac{2 L \sqrt{n} B}{\e} }, 
  \end{equation}
  where $\delta_{0} \eqdef f(x^0) - f^\ast$ and where we plugged in the value of the stepsize $\gamma$ we use. Note that $nT$ appears on both sides in the second term in the maximum in \eqref{corr:rr-nonconvex-proof-1}, hence we can cancel out and simplify:
  \[ n T \geq \frac{24 \delta_{0}}{\e^2} (A L^2 n^2 T)^{1/3} \Longleftrightarrow n T \geq \frac{(24 \delta_{0})^{3/2} L \sqrt{An}}{\e^3}.   \]
  Using this simplified bound in \eqref{corr:rr-nonconvex-proof-1} we obtain that $\min_{k=0, \ldots, T-1} \ecn{\nabla f(x^k)} \leq \e^2$ provided
  \[
    n T \geq \frac{48 \delta_{0} L \sqrt{n}}{\e^2} \max \pbr{ \sqrt{n}, \frac{\sqrt{6 \delta_{0} A}}{\e}, \frac{B}{\e} }.
  \]
\end{proof}

\section{Convergence Results for \algname{IG}}
\label{sec:ig-convergence}
In this section we present results that are extremely similar to the previously obtained bounds for \algname{RR} and \algname{SO}. For completeness, we also provide a full description of \algname{IG} in Algorithm~\ref{alg:ig}.
\begin{algorithm}[t]
    \caption{\algname{Incremental Gradient} \rrbox{(\algname{IG})}}
    \label{alg:ig}
 \begin{algorithmic}[1]
   \Require Stepsize $\gamma > 0$, initial vector $x^0 = x_0^0 \in \R^d$, number of epochs $T$
    \For{epochs $k=0,1,\dotsc,T-1$}
       \For{$i=0, 1, \ldots, n-1$}
          \State $x^{k}_{i+1} = x^k_{i} - \gamma \nabla f_{i+1} (x^k_i)$
       \EndFor
       \State $x^{k+1} = x^{k}_{n}$
    \EndFor
 \end{algorithmic}
 \end{algorithm}

\begin{theorem}
  \label{thm:ig}
  Suppose that Assumption~\ref{asm:f-smoothness} is satisfied. Then we have the following results for the Incremental Gradient method:
    
\begin{itemize}[leftmargin=0.15in,itemsep=0.01in,topsep=0pt]
  \item {\emone If each $f_i$ is $\mu$-strongly convex}: if $\gamma \leq \frac{1}{L}$, then
  \[ \sqn{x^{T} - x^\ast} \leq \br{1 - \gamma \mu}^{nT} \sqn{x^0 - x^\ast} + \frac{\gamma^2 L n^2  \sigmaesc^2}{\mu}. \]
  By carefully choosing the stepsize as in \Cref{corr:all-sc-rr-conv}, we see that this result implies that \algname{IG} has sample complexity  $\ctO\br{\kappa + \frac{\sqrt{\kappa} n \sigmaesc}{\mu \sqrt{\e}}}$ in order to reach a point $\tilde{x}$ with $\norm{\tilde{x} - x^\ast}^2 \leq \e$.
  \item {\emone If $f$ is $\mu$-strongly convex and each $f_i$ is convex}: if $\gamma \leq \frac{1}{\sqrt{2} n L}$, then
  \[ \sqn{x^T - x^\ast} \leq \br{1 - \frac{\gamma \mu n}{2}}^{T} \sqn{x^0 - x^\ast} + 2 \gamma^2 \kappa  n^2 \sigmaesc^2.  \]
  Using the same approach for choosing the stepsize as \Cref{corr:all-sc-rr-conv}, we see that \algname{IG} in this setting reaches an $\e$-accurate solution after $\ctO\br{ n \kappa + \frac{\sqrt{\kappa} n \sigmaesc}{\mu \sqrt{\e}} }$ individual gradient accesses.
  \item {\emone If each $f_i$ is convex}: if $\gamma \leq \frac{1}{\sqrt{2} n L}$, then
  \[ f(\hat{x}^T) - f(x^\ast) \leq \frac{\sqn{x^0 - x^\ast}}{2 \gamma n T} + \frac{\gamma^2 L n^2 \sigmaesc^2}{2} , \]
  where $\hat{x}^T \eqdef \frac{1}{T} \sum_{k=1}^{T} x^k$. Choosing the stepsize $\gamma = \min\pbr{\frac{1}{\sqrt{2} n L}, \frac{\sqrt{\e}}{\sqrt{L} n \sigmaesc} }$, then the average of iterate generated by \algname{IG} is an $\e$-accurate solution (i.e., $f(\hat{x}^T) - f(x^\ast) \leq \e$) provided that the total number of iterations satisfies
  \[ n T \geq \frac{\sqn{x^0 - x^\ast}}{\e} \max \pbr{ \sqrt{8} n L, \frac{\sqrt{L} \sigmaesc n}{\sqrt{\e}} }. \]
  \item {\emone If each $f_i$ is possibly non-convex}: if Assumption~\ref{asm:2nd-moment} holds with constants $A, B \geq 0$ and $\gamma \leq \min \pbr{ \frac{1}{ \sqrt{8} n L}, \frac{1}{(4 L^2 n^3 A T)^{1/3}}}$, then
  \[ 
    \min_{k=0, \ldots, T-1} \sqn{\nabla f(x^k)} \leq \frac{12 \br{f(x^0) - f^\ast}}{\gamma n T} + 8 \gamma^2 L^2 n^2 B^2. \]
  Using an approach similar to \Cref{corr:rr-nonconvex}, we can establish that \algname{IG} reaches a point with gradient norm less than $\e$ provided that the total number of iterations exceeds
  \[ n T \geq \frac{48 \br{f(x^0) - f^\ast} L n}{\e^2} \max \pbr{ \sqrt{2}, \frac{\sqrt{24 \br{f(x^0) - f^\ast}A}}{\e}, \frac{2 B}{\e}  }. \]
\end{itemize}
\end{theorem}

The proof of Thoerem~\ref{thm:ig} is given in the rest of the section, but first we briefly discuss the convergence rates and the relation of the result on strongly convex objectives to the lower bound of \cite{Safran2020good}.

\paragraph{Discussion of the convergence rates.} A brief comparison between the sample complexities given for \algname{IG} in Theorem~\ref{thm:ig} and those given for \algname{RR} (in Table~\ref{tab:conv-rates}) reveals that \algname{IG} has similar rates to \algname{RR} but with a worse dependence on $n$ in the variance term (the term associated with $\sigmaesc$ in the convex case and $B$ in the non-convex case), in particular \algname{IG} is worse by a factor of $\sqrt{n}$. This difference is significant in the large-scale machine learning regime, where the number of data points $n$ can be on the order of thousands to millions.

\paragraph{Discussion of existing lower bounds.} \cite{Safran2020good} give the lower bound (in a problem with $\kappa = 1$) 
\[ \sqn{x^{T} - x^\ast} = \Omega\br{ \frac{\sigmaesc^2}{\mu^2 T^2} }. \]

This implies a sample complexity of $\cO\br{ \frac{n \sigmaesc}{\mu \sqrt{\e}}}$, which matches our upper bound (up to an extra iteration and log factors) in the case each $f_i$ is strongly convex and $\kappa = 1$.

\subsection{Preliminary lemmas for Theorem~\ref{thm:ig}}
\subsubsection{Two lemmas for convex objectives}
\begin{lemma}
  \label{lem:convex_ig_deviation}
	Consider the iterates of Incremental Gradient (Algorithm~\ref{alg:ig}). Suppose that functions $f_1, \ldots, f_n$ are convex and that Assumption~\ref{asm:f-smoothness} is satisfied. Then it holds
  \begin{equation}
    \label{eq:lma-ig-cVt-bound}
    \sum_{k=0}^{n-1} \sqn{x^k_m - x^{k+1}} \leq 4 \gamma^2 L n^2 \sum_{i=0}^{n-1} D_{f_{i+1}} (x^\ast, x^k_i) + 2 \gamma^2 n^3 \sigmaesc^2,
  \end{equation}
  where $\sigmaesc^2$ is the variance at the optimum given by $\sigmaesc^2 \eqdef \frac{1}{n} \sum_{i=1}^{n} \sqn{\nabla f_{i} (x^\ast)}$.
\end{lemma}
\begin{proof}
	The proof of this Lemma is similar to that of Lemma~\ref{lem:convex_rr_deviation} but with a worse dependence on the variance term, since there is no randomness in \algname{IG}. Fix any $m\in\{0,\dotsc, n-1\}$. It holds by definition
  \[
      x^k_m - x^{k+1}
      =\gamma\sum_{i=k}^{n-1}\nabla f_{i+1}(x^k_i)
      = \gamma\sum_{i=k}^{n-1}(\nabla f_{i+1}(x^k_i) - \nabla f_{i+1}(x^\ast))+\gamma\sum_{i=k}^{n-1}\nabla f_{i+1}(x^\ast).
  \]
  Applying Young's inequality to the sums above yields
  \begin{eqnarray*}
      \|x^k_m - x^{k+1}\|^2
      &\overset{\eqref{eq:sum_sqnorm}}{\le}& 2\gamma^2\sqn{\sum_{i=k}^{n-1}(\nabla f_{i+1}(x^k_i) - \nabla f_{i+1}(x^\star))} + 2\gamma^2\sqn{\sum_{i=k}^{n-1}\nabla f_{i+1}(x^\star)} \\
      &\overset{\eqref{eq:sqnorm-sum-bound}}{\le}& 2\gamma^2n\sum_{i=k}^{n-1}\sqn{\nabla f_{i+1}(x^k_i) - \nabla f_{i+1}(x^\star)} + 2\gamma^2\sqn{\sum_{i=k}^{n-1}\nabla f_{i+1}(x^\star)}\\
      &\overset{\eqref{eq:grad_dif_bregman}}{\le}& 4\gamma^2Ln\sum_{i=k}^{n-1}D_{f_{i+1}}(x^\star, x^k_i) + 2\gamma^2\sqn{\sum_{i=k}^{n-1}\nabla f_{i+1}(x^\star)} \\
      &\le& 4\gamma^2Ln\sum_{i=0}^{n-1}D_{f_{i+1}}(x^\star, x^k_i) + 2\gamma^2\sqn{\sum_{i=k}^{n-1}\nabla f_{i+1}(x^\star)} .
  \end{eqnarray*}
  Summing up,
  \begin{equation}
    \label{eq:ig-cVt-bound-proof-1}
    \sum_{m=0}^{n-1} \sqn{x^k_m - x^{k+1}} \leq 4\gamma^2Ln^2\sum_{i=0}^{n-1} D_{f_{i+1}} (x^\ast, x^k_i) + 2\gamma^2 \sum_{m=0}^{n-1} \sqn{\sum_{i=k}^{n-1}\nabla f_{i+1}(x^\star)} .
  \end{equation}
  We now bound the second term in \eqref{eq:ig-cVt-bound-proof-1}. We have
  \begin{eqnarray}
    \sum_{k=0}^{n-1} \sqn{\sum_{i=k}^{n-1} \nabla f_{i+1} (x^\ast)} &\overset{\eqref{eq:sqnorm-sum-bound}}{\leq}& \sum_{k=0}^{n-1} (n-k) \sum_{i=k}^{n-1} \sqn{\nabla f_{i+1} (x^\ast)} \nonumber \\
    &\leq& \sum_{k=0}^{n-1} \br{n-k} \sum_{i=0}^{n-1} \sqn{\nabla f_{i+1} (x^\ast)} \nonumber \\
    &=& \sum_{k=0}^{n-1} \br{n-k} n \sigmaesc^2 = \frac{n^2 (n+1)}{2} \sigmaesc^2 \leq n^3 \sigmaesc^2.
    \label{eq:ig-cVt-bound-proof-2}
  \end{eqnarray} 
  Using \eqref{eq:ig-cVt-bound-proof-2} in \eqref{eq:ig-cVt-bound-proof-1}, we derive
  \[
    \sum_{k=0}^{n-1} \sqn{x^k_m - x^{k+1}} \leq 4 \gamma^2 L n^2 \sum_{i=0}^{n-1} D_{f_{i+1}} (x^\ast, x^k_i) + 2 \gamma^2 n^3 \sigmaesc^2.
  \]
\end{proof}

\begin{lemma}
  \label{lem:convex-ig-recursion}
  Assume the functions $f_1, \ldots, f_n$ are convex and that Assumption~\ref{asm:f-smoothness} is satisfied. If Algorithm~\ref{alg:ig} is run with a stepsize $\gamma \leq \frac{1}{\sqrt{2} L n}$, then
  \[ \sqn{x^{k+1} - x^\ast} \leq \sqn{x^k - x^\ast} - 2 \gamma n \br{f(x^{k+1}) - f(x^\ast)} +  \gamma^3 L n^3 \sigmaesc^2. \]
\end{lemma}
\begin{proof}
  The proof for this lemma is identical to \Cref{lem:improved_convex} but with the estimate of \Cref{lem:convex_ig_deviation} used for $\sum_{i=0}^{n-1} \sqn{x^k_i - x^{k+1}}$ instead of \Cref{lem:convex_rr_deviation}. We only include it for completeness. Define the sum of gradients used in the $k$-th epoch as $g^{k} \eqdef \sum_{i=0}^{n-1} \nabla f_{i+1} (x^k_i)$. By definition of $x^{k+1}$, we have $x^{k+1} = x^k - \gamma g^k$. Using this,
	\begin{align*}
		\|x^{k}-x^\star\|^2
		=\|x^{k+1}+\gamma g^k-x^\star\|^2
		&=\|x^{k+1}-x^\star\|^2+2\gamma\ev{g^k, x^{k+1}-x^\star}+\gamma^2\|g^k\|^2 \\
		&\ge \|x^{k+1}-x^\star\|^2+2\gamma\ev{g^k, x^{k+1}-x^\star} \\
		&= \|x^{k+1}-x^\star\|^2+2\gamma\sum_{i=0}^{n-1}\ev{\nabla f_{i+1}(x^k_i), x^{k+1}-x^\star}.
	\end{align*}
	For any $i$ we have the following decomposition
	\begin{align}
		\ev{\nabla f_{i+1}(x^k_i), x^{k+1}-x^\star}
    &= [f_{i+1}(x^{k+1})-f_{i+1}(x^\star)]\\
    & \quad + [f_{i+1}(x^\star) - f_{i+1}(x^k_i) - \<\nabla f_{i+1}(x^k_i), x^{k}_i-x^\star>] \notag\\
		& \quad - [f_{i+1}(x^{k+1}) - f_{i+1}(x^k_i)-\<\nabla f_{i+1}(x^k_i), x^{k+1}-x^k_i>] \notag\\
		&= [f_{i+1}(x^{k+1})-f_{i+1}(x^\star)] + D_{f_{i+1}}(x^\star, x^k_i) - D_{f_{i+1}}(x^{k+1}, x^k_i). \label{eq:ig_decomposition}
	\end{align}
	Summing the first quantity in~\eqref{eq:ig_decomposition} over $i$ from $0$ to $n-1$ gives
	\begin{align*}
		\sum_{i=0}^{n-1}[f_{i+1}(x^{k+1})-f_{i+1}(x^\star)]
		= n(f(x^{k+1})-f^\star).
	\end{align*}
	Now let us work out the third term in the decomposition~\eqref{eq:ig_decomposition} using $L$-smoothness,
	\begin{align*}
		D_{f_{i+1}}(x^{k+1}, x^k_i)
		\le \frac{L}{2}\|x^{k+1}-x^k_i\|^2.
	\end{align*}
	We next sum the right-hand side over $i$ from $0$ to $n-1$ and use \Cref{lem:convex_ig_deviation}
	\begin{eqnarray*}
		\sum_{i=0}^{n-1} D_{f_{i+1}}(x^{k+1}, x^k_i)
    &\leq& \frac{L}{2} \sum_{i=0}^{n-1} \sqn{x^{k+1} - x^k_i} \\ 
		&\overset{\eqref{eq:lma-ig-cVt-bound}}{\le}&  2\gamma^2 L^2n^2\sum_{i=0}^{n-1} D_{f_{i+1}}(x^\star, x^k_i)
		+\gamma^2  L n^3\sigma_\star^2.
	\end{eqnarray*}
	Therefore, we can lower-bound the sum of the second and the third term in~\eqref{eq:ig_decomposition} as
	\begin{align*}
		\sum_{i=0}^{n-1}(D_{f_{i+1}}(x^\star, x^k_i) - D_{f_{i+1}}(x^{k+1}, x^k_i))
    &\ge \sum_{i=0}^{n-1} D_{f_{i+1}}(x^\star, x^k_i)  \\
    &\qquad - \br{2\gamma^2L^2n^2\sum_{i=0}^{n-1} D_{f_{i+1}}(x^\star, x^k_i) - \gamma^2 L n^3\sigma_\star^2} \\
		&= (1-2\gamma^2L^2n^2)\sum_{i=0}^{n-1} D_{f_{i+1}}(x^\star, x^k_i) - \gamma^2 L n^3\sigma_\star^2 \\
		&\ge  - \gamma^2 L n^3\sigma_\star^2,
	\end{align*}
	where in the third inequality we used that $\gamma \leq \frac{1}{\sqrt{2} L n}$ and that $D_{f_{i+1}} (x^\ast, x^k_i)$ is nonnegative. Plugging this back into the lower-bound on $\|x^k-x^\star\|^2$ yields
	\[
		\|x^k-x^\star\|^2
		\ge \sqn{x^{k+1}-x^\star} + 2\gamma n\br{f(x^{k+1})-f^\star}- \gamma^3 L n^3\sigma_\star^2.
	\]
	Rearranging the terms gives the result.
\end{proof}

\subsubsection{A lemma for non-convex objectives}

\begin{lemma}
  \label{lemma:vt-bound-ig}
  Suppose that Assumption~\ref{asm:f-smoothness} holds. Suppose that Algorithm~\ref{alg:ig} is used with a stepsize $\gamma > 0$ such that $\gamma \leq \frac{1}{2 L n}$. Then we have,
  \begin{equation}
    \label{eq:vt-bound-ig}
    \sum_{i=1}^{n} \sqn{x^k_i - x^k} \leq 4 \gamma^2 n^3 \sqn{\nabla f(x^k)} + 4 \gamma^2 n^3 \sigma_{k}^2,
  \end{equation}
  where $\sigma_{k}^2 \eqdef \frac{1}{n} \sum_{j=1}^{n} \sqn{\nabla f_{j} (x^k) - \nabla f(x^k)}$.
\end{lemma}
\begin{proof}
  Let $i \in \{ 1, 2, \ldots, n \}$. Then we can bound the deviation of a single iterate as,
  \begin{eqnarray*}
      \sqn{x^k_i - x^k} = \sqn{x^k_0 - \gamma \sum_{j=0}^{i-1} \nabla f_{j+1} (x^k_{j}) - x^k} &=& \gamma^2 \sqn{\sum_{j=0}^{i-1} \nabla f_{j+1} (x^k_{j})} \\
      &\overset{\eqref{eq:sqnorm-sum-bound}}{\leq}& \gamma^2 i \sum_{j=0}^{i-1} \sqn{\nabla f_{j+1} (x^k_{j})}.
  \end{eqnarray*}
  Because $i \leq n$, we have
  \begin{equation}
      \sqn{x^k_i - x^k} \leq \gamma^2 i \sum_{j=0}^{i-1} \sqn{\nabla f_{i+1} (x^k_{j})} \leq \gamma^2 n \sum_{j=0}^{i-1} \sqn{\nabla f_{i+1} (x^k_{j})} \leq \gamma^2 n \sum_{j=0}^{n-1} \sqn{\nabla f_{i+1} (x^k_{j})}.
      \label{eq:vt-ig-1}
  \end{equation}
  Summing up allows us to estimate $V^k$:
  \begin{eqnarray}
      V^k &=& \sum_{i=1}^{n} \sqn{x^k_i - x^k} \nonumber \\
      &\overset{\eqref{eq:vt-ig-1}}{\leq} & \sum_{i=1}^{n} \br{\gamma^2 n \sum_{j=0}^{n-1} \sqn{\nabla f_{j+1} (x^k_j)}} \nonumber \\ & = &\gamma^2 n^2 \sum_{j=0}^{n-1} \sqn{\nabla f_{j+1} (x^k_j)} \nonumber \\
      &\overset{\eqref{eq:sum_sqnorm}}{\leq} & 2 \gamma^2 n^2 \sum_{j=0}^{n-1} \br{ \sqn{\nabla f_{j+1} (x^k_j) - \nabla f_{i+1} (x^k) } + \sqn{\nabla f_{j+1} (x^k) } } \nonumber \\
      &=& 2 \gamma^2 n^2 \sum_{j=0}^{n-1} \sqn{\nabla f_{j+1} (x^k_j) - \nabla f_{j+1} (x^k) } + 2 \gamma^2 n^2 \sum_{j=0}^{n-1} \sqn{\nabla f_{j+1} (x^k)}.
      \label{eq:vt-ig-2}
  \end{eqnarray}
  For the first term in \eqref{eq:vt-ig-2} we can use the smoothness of individual losses and that $x^k_0 = x^k$:
  \begin{align}
      \sum_{j=0}^{n-1} \sqn{\nabla f_{j+1} (x^k_j) - \nabla f_{j+1} (x^k)} &\overset{\eqref{eq:nabla-Lip}}{\leq}  L^2 \sum_{j=0}^{n-1}  \sqn{x^k_j - x^k} = L^2 \sum_{j=1}^{n-1} \sqn{x^k_j - x^k} = L^2 V^k.
      \label{eq:vt-ig-3}
  \end{align}
  The second term in \eqref{eq:vt-ig-2} is a sum over all the individual gradient evaluated at the same point $x^k$. Hence, we can drop the permutation subscript and then use the variance decomposition:
  \begin{eqnarray}
    \sum_{j=0}^{n-1} \sqn{\nabla f_{i+1} (x^k)} &=& \sum_{j=1}^{n} \sqn{\nabla f_{j} (x^k)} \nonumber \\
    &\overset{\eqref{eq:variance_m}}{=}&  n \sqn{\nabla f(x^k)} + \sum_{j=1}^{n} \sqn{\nabla f_{j} (x^k) - \nabla f(x^k)} \nonumber \\
    &=& n \sqn{\nabla f(x^k)} + n \sigma_{k}^2.
    \label{eq:vt-ig-4}
  \end{eqnarray}
  We can then use~\eqref{eq:vt-ig-3} and~\eqref{eq:vt-ig-4} in~\eqref{eq:vt-ig-2},
  \[
   V^k \leq 2 \gamma^2 L^2 n^2 V^k + 2 \gamma^2 n^3 \sqn{\nabla f(x^k)} + 2 \gamma^2 n^3 \sigma_{k}^2.
  \]
Since  $V^k$ shows up in both sides of the equation, we can rearrange to obtain
  \[
    \br{1 - 2 \gamma^2 L^2 n^2 } V^k \leq 2 \gamma^2 n^3 \sqn{\nabla f(x^k)} + 2 \gamma^2 n^3 \sigma_{k}^2.
  \]
  If $\gamma \leq \frac{1}{2 L n}$, then $1 - 2 \gamma^2 L^2 n^2  \geq \frac{1}{2}$ and hence
  \[
    V^k \leq 4 \gamma^2 n^3 \sqn{\nabla f(x^k)} + 4 \gamma^2 n^3 \sigma_{k}^2.
  \]
\end{proof}

\subsection{Proof of Theorem~\ref{thm:ig}}
\begin{proof}

  \begin{itemize}[leftmargin=0.15in,itemsep=0.01in,topsep=0pt]
    \item {\emone If each $f_i$ is $\mu$-strongly convex}: The proof follows that of \Cref{thm:all-sc-rr-conv}. Define $$x^\ast_i = x^\ast - \gamma \sum_{j=0}^{i-1} \nabla f_{j+1} (x^\ast).$$ First, we have
    \begin{align*}
      &\|x^k_{i+1}-x^\ast_{i+1}\|^2\\
      &=\|x^k_{i}-x^\ast_{i}\|^2-2\gamma\<\nabla f_{i+1}(x^k_i)-\nabla f_{i+1}(x^\star), x^k_i - x^\ast_i>+\gamma^2\|\nabla f_{i+1}(x^k_i) - \nabla f_{i+1}(x^\star)\|^2.
    \end{align*}
    Using the same three-point decomposition as \Cref{thm:all-sc-rr-conv} and strong convexity, we have
    \begin{align*}
      - \ev{\nabla f_{i+1} (x^k_i) - \nabla f_{i+1} (x^\ast), x^k_i - x^\ast_i} &= - D_{f_{i+1}} (x^\ast_i, x^k_i) - D_{f_{i+1}} (x^k_i, x^\ast) + D_{f_{i+1}} (x^\ast_i, x^\ast) \\
      &\leq - \frac{\mu}{2} \sqn{x^k_i - x^\ast_i} - D_{f_{i+1}} (x^k_i, x^\ast) + D_{f_{i+1}} (x^\ast_i, x^\ast).
    \end{align*}
    Using smoothness and convexity
    \[ \frac{1}{2L} \sqn{\nabla f_{i+1} (x^k_i) - \nabla f_{i+1} (x^\ast)} \leq D_{f_{i+1}} (x^k_i, x^\ast). \]
    Plugging in the last two inequalities into the recursion, we get
    \begin{align}
      \sqn{x^k_{i+1} - x^\ast_{i+1}} &\leq \br{1 - \gamma \mu} \sqn{x^k_i - x^\ast_i} - 2 \gamma \br{1 - \gamma L} D_{f_{i+1}} (x^k_i, x^\ast) + 2 \gamma D_{f_{i+1}} (x^\ast_i, x^\ast). \nonumber \\
      &\le \br{1 - \gamma \mu} \sqn{x^k_i - x^\ast_i} + 2 \gamma D_{f_{i+1}} (x^\ast_i, x^\ast).
      \label{eq:ig-proof-allsc-1}
    \end{align}
    For the last Bregman divergence, we have
    \begin{eqnarray*}
      D_{f_{i+1}} (x^\ast_i, x^\ast) & \overset{\eqref{eq:L-smooth-intro}}{\leq} & \frac{L}{2} \sqn{x^\ast_i - x^\ast} \nonumber \\ &=& \frac{\gamma^2 L }{2} \sqn{\sum_{j=0}^{i-1} \nabla f_{j+1} (x^\ast)} \\
      &\overset{\eqref{eq:sqnorm-sum-bound}}{\leq}& \frac{\gamma^2 L  i}{2} \sum_{j=0}^{i-1} \sqn{\nabla f_{j+1} (x^\ast)} \\
      &=& \frac{\gamma^2 L  i n}{2} \sigmaesc^2 \leq \frac{\gamma^2 L  n^2}{2} \sigmaesc^2.
    \end{eqnarray*}
    Plugging this into \eqref{eq:ig-proof-allsc-1}, we get
    \[ \sqn{x^k_{i+1} - x^\ast_{i+1}} \leq \br{1 - \gamma \mu} \sqn{x^k_i - x^\ast_i} +  \gamma^3 L n^2 \sigmaesc^2. \]
    We recurse and then use that $x^\ast_n = x^\ast$, $x^{k+1} = x^k_n$, and that $x^\ast_0 = x^\ast$, obtaining
    \begin{align*}
      \sqn{x^{k+1} - x^\ast} = \sqn{x^k_n - x^\ast_n} &\leq \br{1 - \gamma \mu}^{n} \sqn{x^k_0 - x^\ast_0} + \gamma^3 L  n^2 \sigmaesc^2 \sum_{j=0}^{n-1} \br{1 - \gamma \mu}^{j} \\
      &= \br{1 - \gamma \mu}^{n} \sqn{x^k - x^\ast} +  \gamma^3 L n^2 \sigmaesc^2 \sum_{j=0}^{n-1} \br{1 - \gamma \mu}^{j}.
    \end{align*}
    Recursing again,
    \begin{align*}
      \sqn{x^{T} - x^\ast} &\leq \br{1 - \gamma \mu}^{n T} \sqn{x^0 - x^\ast} + \gamma^3 L  n^2 \sigmaesc^2 \sum_{j=0}^{n-1} \br{1 - \gamma \mu}^{j} \sum_{k=0}^{T-1} \br{1 - \gamma \mu}^{nt} \\
      &=  \br{1 - \gamma \mu}^{n T} \sqn{x^0 - x^\ast} +  \gamma^3 L n^2 \sigmaesc^2 \sum_{k=0}^{nT-1} \br{1 - \gamma \mu}^{k} \\
      &\leq \br{1 - \gamma \mu}^{n T} \sqn{x^0 - x^\ast} + \frac{\gamma^3 L  n^2 \sigmaesc^2}{\gamma \mu} \\
      &= \br{1 - \gamma \mu}^{n T} \sqn{x^0 - x^\ast} +\gamma^2  \kappa  n^2 \sigmaesc^2.
    \end{align*}
    \item {\emone If $f$ is $\mu$-strongly convex and each $f_i$ is convex}: the proof is identical to that of Theorem~\ref{thm:only-f-sc-rr-conv} but using Lemma~\ref{lem:convex-ig-recursion} instead of Lemma~\ref{lem:improved_convex}, and we omit it for brevity.
    \item {\emone If each $f_i$ is convex}: the proof is identical to that of Theorem~\ref{thm:weakly-convex-f-rr-conv} but using Lemma~\ref{lem:convex-ig-recursion} instead of Lemma~\ref{lem:improved_convex}, and we omit it for brevity.
    \item {\emone If each $f_i$ is possibly non-convex}: note that Lemma~\ref{lemma:epoch-recursion-non-convex} also applies to \algname{IG} without change, hence if $\gamma \leq \frac{1}{Ln}$ we have
    \[ f(x^{k+1}) \leq f(x^k) - \frac{\gamma n}{2} \sqn{\nabla f(x^k)} + \frac{\gamma L^2}{2} \sum_{i=1}^{n} \sqn{x^k - x^k_i}. \]
    We may then apply \Cref{lemma:vt-bound-ig} to get for $\gamma \leq \frac{1}{2 L n}$
    \begin{align*}
      f(x^{k+1}) &\leq f(x^k) - \frac{\gamma n}{2} \sqn{\nabla f(x^k)} + \frac{\gamma L^2}{2} \br{4 \gamma^2 n^3 \sqn{\nabla f(x^k)} + 4 \gamma^2 n^3 \sigma_{k}^2} \\
      &= f(x^k) - \frac{\gamma n}{2} \br{1 - 4 \gamma^2 L^2 n^2} \sqn{\nabla f(x^k)} + 2 \gamma^3 L^2 n^3 \sigma_{k}^2.
    \end{align*}
    Using that $\gamma \leq \frac{1}{\sqrt{8} L n}$ and subtracting $f^\ast$ from both sides, we derive
    \[
      f(x^{k+1}) - f^\ast \leq \br{f(x^k) - f^\ast} - \frac{\gamma n}{4} \sqn{\nabla f(x^k)} + 2  \gamma^3 L^2 n^3 \sigma_{k}^2.
    \]
    Using Assumption~\ref{asm:2nd-moment}, we get
    \begin{equation}
      f(x^{k+1}) - f^\ast \leq \br{1 + 4 \gamma^3 L^2 A n^3 } \br{f(x^k) - f^\ast} - \frac{\gamma n}{4} \sqn{\nabla f(x^k)} + 2 \gamma^3 L^2 n^3 B^2.
      \label{eq:ig-nonconvex-thm-1}
    \end{equation}
    Applying Lemma~\ref{lemma:noncvx-recursion-solution} to \eqref{eq:ig-nonconvex-thm-1}, thus, gives
    \begin{equation}
      \min_{k=0, \ldots, T-1} \sqn{\nabla f(x^k)} \leq \frac{4 \br{1 + 4 \gamma^3 L^2 A n^3 }^{T}}{\gamma n T} \br{f(x^0) - f^\ast} + 8 \gamma^2 L^2 n^2 B^2.
      \label{eq:ig-nonconvex-thm-2}
    \end{equation}
    Note that by our assumption on the stepsize, $4 \gamma^3 L^2 A n^3  T \leq 1$, hence,
    \[ \br{1 + 4 \gamma^3 L^2 A n^3 }^{T} \leq \exp\br{4 \gamma^3 L^2 A n^3  T} \leq \exp\br{1} \leq 3. \]
    It remains to use this in \eqref{eq:ig-nonconvex-thm-2}.
  \end{itemize}
\end{proof}

\chapter{Appendix for Chapter~\ref{chapter:proxrr}}
\graphicspath{{rr_prox/}}

\section{Proof of \Cref{thm:shuffling-radius-bound}}
\begin{proof}
	By the $L_{i}$-smoothness of $f_i$ and the definition of $x^\ast_i$, we can replace the Bregman divergence in \eqref{eq:bregman-div-rad}  with the bound
		\begin{align}
		\ec{D_{f_{\pi_{i}}} (x^\ast_i, x^\ast)} \; \overset{\eqref{eq:L-smooth-intro}}{\le} \; \ec{\frac{L_{\pi_i}}{2}\norm{x^\ast_i - x^\ast}^2}  &\le   \frac{L_{\max}}{2}\ec{\|x^\ast_i - x^\ast\|^2} \notag \\ 
		& \overset{\eqref{eq:x_ast_i}}{=}   \frac{\gamma^2 L_{\max}}{2} \mathbb{E}\Biggl[\biggl\| \sum_{j=0}^{i-1}\nabla f_{\pi_j}(x^\ast) \biggr\|^2\Biggr]\notag  \\
		&=  \frac{\gamma^2 L_{\max} i^2}{2} \mathbb{E}\Biggl[\biggl\|\frac{1}{i} \sum_{j=0}^{i-1}\nabla f_{\pi_j}(x^\ast) \biggr\|^2\Biggr] \notag \\
		&= \frac{\gamma^2 L_{\max} i^2}{2}  \ecn{\bar{X}_\pi },\label{eq:8g9fd8gf9d}
	\end{align}
	where $\bar{X}_\pi = \frac{1}{j}\sum_{j=0}^{i-1} X_{\pi_j}$ with  $X_j \eqdef  \nabla f_{j}(x^\ast)$ for $j=1,2,\dots,n$.  Since $\bar{X} = \nabla f(x^\ast)$,  by applying Lemma~\ref{lemma:sampling-without-replacement} we get
    \begin{equation} 
        \ecn{\bar{X}_\pi } \;= \;   \norm{\bar{X}}^2 + \ecn{\bar{X}_\pi - \bar{X}}
        \; \overset{\eqref{eq:b97fg07gdf_08yf8d} + \eqref{eq:UG(*G(DG(*DGg87gf7ff}}{=} \; \norm{\nabla f(x^\ast)}^2 +\frac{n-i}{i (n-1)} \sigmaesc^2. \label{eq:08gfd_898fd*(*gdJS}
    \end{equation}
    It remains to combine  \eqref{eq:8g9fd8gf9d}  and \eqref{eq:08gfd_898fd*(*gdJS}, use the bounds  $i^2\leq n^2$ and
    $i(n-i)\le \frac{n(n-1)}{2}$, which holds for all $i\in \{0,1,  \dots, n-1\}$, and divide both sides of the resulting inequality by $\gamma^2$.
\end{proof}

\section{Main Convergence Proofs}
\subsection{A key lemma for shuffling-based methods}
The intermediate limit points $x^\ast_i$ are extremely important for showing tight convergence guarantees for Random Reshuffling even without proximal operator. The following lemma illustrates that by giving a simple recursion, whose derivation follows \cite[Proof of Theorem 1]{MKR2020rr}. The proof is included for completeness.

\begin{lemma}[Theorem 1 in \cite{MKR2020rr}]
	\label{lemma:inner-loop}
	Suppose that each $f_i$ is $L_i$-smooth and $\lambda$-strongly convex (where $\lambda = 0$ means each $f_i$ is just convex). Then the inner iterates generated by Algorithm~\ref{alg:proxrr} satisfy
	\begin{align}
		\ecn{x^k_{i+1} - x^\ast_{i+1}} 
		&\leq \br{1 - \gamma \lambda} \ecn{x^k_i - x^\ast_i} - 2 \gamma \br{1 - \gamma L_{\max}} \ec{D_{f_{\pi_i}} (x^k_i, x^\ast)} \notag\\
		&\qquad + 2 \gamma^3 \sigmarr^2,		\label{eq:inner-loop-recursion}
	\end{align}
	where $x^\ast_i $ is as in \eqref{eq:x_ast_i}, $i=0,1,\dots,n-1$, and $x^\ast$ is any minimizer of $P$.
\end{lemma}
\begin{proof}
	By definition of $x^{k}_{i+1}$ and $x^\ast_{i+1}$, we have
	\begin{align}
		\begin{split}
			\ecn{x^k_{i+1}-x^\ast_{i+1}}  &= \ecn{x^k_{i}-x^\ast_{i}}-2 \gamma \ec{\<\nabla f_{\pi_i}(x^k_i)-\nabla f_{\pi_i}(x^\ast), x^k_i - x^\ast_i> } \\
			& \qquad + \gamma^2 \ecn{\nabla f_{\pi_i}(x^k_i) - \nabla f_{\pi_i}(x^\ast)}. 		
		\end{split}\label{eq:inner-loop-proof-1}
	\end{align}
	Note that the third term in \eqref{eq:inner-loop-proof-1} can be bounded as
	\begin{equation}
		\label{eq:inner-loop-proof-2}
		\sqn{\nabla f_{\pi_i} (x^k_i) - \nabla f_{\pi_i} (x^\ast)} \leq 2 L_{\max} \cdot D_{f_{\pi_i}} (x^k_i, x^\ast).
	\end{equation}
	We may rewrite the second term in \eqref{eq:inner-loop-proof-1} using the three-point identity \cite[Lemma 3.1]{Chen1993} as
	\begin{equation}
		\label{eq:inner-loop-proof-3}
		\ev{\nabla f_{\pi_i} (x^k_i) - \nabla f_{\pi_{i}} (x^\ast), x^k_i - x^\ast_i } = D_{f_{\pi_i}} (x^\ast_i, x^k_i) + D_{f_{\pi_i}} (x^k_i, x^\ast) - D_{f_{\pi_i}} (x^\ast_i, x^\ast).
	\end{equation}
	Combining \eqref{eq:inner-loop-proof-1}, \eqref{eq:inner-loop-proof-2}, and \eqref{eq:inner-loop-proof-3} we obtain
	\begin{align}
		\begin{split}
			\ecn{x^k_{i+1} - x^\ast_{i+1}} \leq \ecn{x^k_i - x^\ast_i} &- 2 \gamma \cdot \ec{D_{f_{\pi_i}} (x^\ast_i, x^k_i)} + 2 \gamma \cdot \ec{D_{f_{\pi_i}} (x^\ast_i, x^\ast)} \\ 
			&- 2 \gamma \br{1 - \gamma L_{\max}} \ec{D_{f_{\pi_i}} (x^k_i, x^\ast)}.	
		\end{split}\label{eq:inner-loop-proof-4}
	\end{align}
	Using  $\lambda$-strong convexity of $f_{\pi_i}$, we derive
	\begin{equation}
		\label{eq:inner-loop-proof-5}
		\frac{\lambda}{2} \sqn{x^k_i - x^\ast_i} \leq D_{f_{\pi_i}} (x^\ast_i, x^k_i).
	\end{equation}
	Furthermore, by the definition of shuffling radius (Definition~\ref{def:bregman-div-rad}), we have
	\begin{equation}
		\label{eq:inner-loop-proof-6}
		\ec{D_{f_{\pi_i}} (x^\ast_i, x^\ast)} \leq \max_{i=0, \ldots, n-1}  \ec{D_{f_{\pi_i}} (x^\ast_i, x^\ast) }  = \gamma^2 \sigmarr^2.
	\end{equation}
	Using~\eqref{eq:inner-loop-proof-5} and~\eqref{eq:inner-loop-proof-6} in \eqref{eq:inner-loop-proof-4} yields \eqref{eq:inner-loop-recursion}.
\end{proof}

\subsection{Proof of \Cref{thm:f-strongly-convex-psi-convex}}
\begin{proof}
	Starting with Lemma~\ref{lemma:inner-loop} with $\lambda = \mu$, we have
	\begin{align*}
		\ecn{x^k_{i+1} - x^\ast_{i+1}} 
		&\leq \br{1 - \gamma \mu} \ecn{x^k_i - x^\ast_i} - 2 \gamma \br{ 1 - \gamma L_{\max}} \ec{D_{f_{\pi_i}} (x^k_i, x^\ast) }\\
		&\qquad + 2 \gamma^3 \sigmarr^2.
	\end{align*}
	Since $D_{f_{\pi}} (x^k_i, x^\ast)$ is a Bregman divergence of a convex function, it is nonnegative. Combining this with the fact that the stepsize satisfies $\gamma \leq \frac{1}{L_{\max}}$, we have
	\[ \ecn{x^k_{i+1} - x^\ast_{i+1}} \leq \br{1 - \gamma \mu} \ecn{x^k_i - x^\ast_i} + 2 \gamma^3 \sigmarr^2. \]
	Unrolling this recursion for $n$ steps, we get
	\begin{align}
		\ecn{x^k_{n} - x^\ast_n} &\leq \br{1 - \gamma \mu}^{n} \ecn{x^k_0 - x^\ast_0} + 2 \gamma^3 \sigmarr^2 \br{ \sum_{j=0}^{n-1} \br{1 - \gamma \mu}^{j} } \notag \\
		&= \br{1 - \gamma \mu}^{n} \ecn{x^k - x^\ast} + 2 \gamma^3 \sigmarr^2 \br{ \sum_{j=0}^{n-1} \br{1 - \gamma \mu}^j },		\label{eq:thm-f-sc-proof-1}
	\end{align}
	where we used the fact that $x^k_0 - x^\ast_0 = x^k - x^\ast$. Since $x^\ast$ minimizes $P$, we have by \Cref{prop:fixed-point} that
	\[ 
		x^\ast = \prox_{\gamma n \psi} \br{ x^\ast - \gamma \sum_{i=0}^{n-1} \nabla f_{\pi_i} (x^\ast) } = \prox_{\gamma n \psi} \br{x^\ast_n}. 
	\] 
	Moreover, by Lemma~\ref{prop:prox-contraction} we obtain that 
	\[ 
		\sqn{x^{k+1} - x^\ast} = \sqn{\prox_{\gamma n \psi} (x^k_n) - \prox_{\gamma n \psi} (x^\ast_n) } \leq \sqn{x^k_n - x^\ast_n}. 
	\]
	Using this in~\eqref{eq:thm-f-sc-proof-1} yields
	\[
		\ecn{x^{k+1} - x^\ast} \leq \br{1 - \gamma \mu}^n \ecn{x^k - x^\ast} + 2 \gamma^3 \sigmarr^2 \br{ \sum_{j=0}^{n-1} \br{1 - \gamma \mu}^j }.
	\]
	We now unroll this recursion again for $T$ steps
	\begin{align}
		\label{eq:thm-f-sc-proof-2}
			\ecn{x^{T} - x^\ast} 
			&\leq \br{1 - \gamma \mu}^{n T} \ecn{x^0 - x^\ast} \notag \\
			&\qquad + 2 \gamma^3 \sigmarr^2 \br{ \sum_{j=0}^{n-1} \br{1 - \gamma \mu}^j } \br{ \sum_{i=0}^{T-1} \br{1 - \gamma \mu}^{n i} }.
	\end{align}
	Following the analysis of Chapter~\ref{chapter:rr}, we rewrite and bound the product in the last term as
	\begin{align*}
		\br{ \sum_{j=0}^{n-1} \br{1 - \gamma \mu}^j } \br{ \sum_{i=0}^{T-1} \br{1 - \gamma \mu}^{n i} } &= \sum_{j=0}^{n-1} \sum_{i=0}^{T-1} \br{1 - \gamma \mu}^{ni + j} \\
		&= \sum_{k=0}^{nT-1} \br{1 - \gamma \mu}^k \\
		&\leq \sum_{k=0}^{\infty} \br{1 - \gamma \mu}^{k} \quad = \frac{1}{\gamma \mu}.
	\end{align*}
	It remains to plug this bound into~\eqref{eq:thm-f-sc-proof-2}.
\end{proof}

\subsection{Proof of \Cref{thm:psi-strongly-convex-f-convex}}
\begin{proof}
	Starting with Lemma~\ref{lemma:inner-loop} with $\lambda = 0$, we have
	\[ \ecn{x^k_{i+1} - x^\ast_{i+1}} \leq \ecn{x^k_i - x^\ast_i} - 2 \gamma \br{1 - \gamma L_{\max} } \ec{D_{f_{\pi_i}} (x^k_i, x^\ast) } + 2 \gamma^3 \sigmarr^2. \]
	Since $\gamma \leq \frac{1}{L_{\max}}$ and $D_{f_{\pi}} (x^k_i, x^\ast)$ is nonnegative we may simplify this to
	\[ \ecn{x^k_{i+1} - x^\ast_{i+1}} \leq \ecn{x^k_i - x^\ast_i} + 2 \gamma^3 \sigmarr^2. \]
	Unrolling this recursion over an epoch we have
	\begin{equation}
		\label{eq:thm-psi-sc-proof-1}
		\ecn{x^k_{n} - x^\ast_n} \leq \ecn{x^k_0 - x^\ast_0} + 2 \gamma^3 \sigmarr^2 n = \ecn{x^k - x^\ast} + 2 \gamma^3 \sigmarr^2 n.
	\end{equation}
	Since $x^\ast$ minimizes $P$, we have by \Cref{prop:fixed-point} that
	\[ x^\ast = \prox_{\gamma n \psi} \br{ x^\ast - \gamma \sum_{i=0}^{n-1} \nabla f_{\pi_i} (x^\ast) } = \prox_{\gamma n \psi} \br{x^\ast_n}. \] 
	Hence, $x^{k+1} - x^\ast = \prox_{\gamma n \psi} (x^k_n) - \prox_{\gamma n \psi} (x^\ast_n)$. We may now use Lemma~\ref{prop:prox-contraction} to get
	\[ \br{1 + 2 \gamma \mu n} \ecn{x^{k+1} - x^\ast} \leq \ecn{x^k_n - x^\ast_n}. \]
	Combining this with \eqref{eq:thm-psi-sc-proof-1}, we obtain
	\[
		\ecn{x^{k+1} - x^\ast} \leq \frac{1}{1 +2 \gamma \mu n} \ecn{x^k - x^\ast} + \frac{2 \gamma^3 \sigmarr^2 n}{1 + 2\gamma \mu n}.
	\]
	We may unroll this recursion again, this time for $T$ steps, and then use that $\sum_{j=1}^{T-1} \br{1 +2 \gamma \mu n}^{-j} \leq \sum_{j=1}^{\infty} \br{1 + 2\gamma \mu n}^{-j} = 1/(2\gamma \mu n)$:
	\begin{align*}
		\ecn{x^{T} - x^\ast} &\leq \br{1 + 2\gamma \mu n}^{-T} \ecn{x^0 - x^\ast} + \frac{2 \gamma^3 \sigmarr^2 n}{1 +2 \gamma \mu n} \Biggl(\sum_{j=0}^{T-1} \br{1 +2 \gamma \mu n}^{-j}\Biggr) \\
		&= \br{1 +2 \gamma \mu n}^{-T} \ecn{x^0 - x^\ast} + 2 \gamma^3 \sigmarr^2 n \Biggl(\sum_{j=1}^{T} \br{1 +2 \gamma \mu n}^{-j} \Biggr) \\
		&\leq \br{1 + 2\gamma \mu n}^{-T} \ecn{x^0 - x^\ast} + 2 \gamma^3 \sigmarr^2 n \frac{1}{2\gamma \mu n} \\
		&= \br{1 +2 \gamma \mu n}^{-T} \ecn{x^0 - x^\ast} + \frac{ \gamma^2 \sigmarr^2}{\mu}.
	\end{align*}
\end{proof}

\section{Convergence of \algname{SGD} (Proof of \Cref{thm:conv-prox-sgd})}
\begin{proof}
	We will prove the case when $\psi$ is $\mu$-strongly convex. The other result follows as a straightforward special case of \cite[Theorem 4.1]{gorbunov2020unified}. We start by analyzing one step of \algname{SGD} with stepsize $\gamma_k = \gamma$ and using \Cref{prop:fixed-point}
	\begin{align}
		\sqn{x^{k+1} - x^\ast} &= \sqn{ \prox_{\gamma \psi} (x^k - \gamma \nabla f (x^k; \xi)) - \prox_{\gamma \psi} (x^\ast - \gamma \nabla f (x^\ast))} \nonumber \\
		&\leq \frac{1}{1 + 2 \gamma \mu} \sqn{x^k - \gamma \nabla f (x^k; \xi) - (x^\ast - \gamma \nabla f (x^\ast))}. \label{eq:sgd-proof-1}
	\end{align}
	We may write the squared norm term in \eqref{eq:sgd-proof-1} as
	\begin{align}
			&\sqn{x^k - \gamma \nabla f (x^k; \xi) - (x^\ast - \gamma \nabla f (x^\ast))} \notag\\
			&= \sqn{x^k - x^\ast} - 2 \gamma \ev{x^k - x^\ast, \nabla f (x^k; \xi) - \nabla f (x^\ast)} + \gamma^2 \sqn{\nabla f (x^k; \xi) - \nabla f (x^\ast)}. \label{eq:sgd-proof-2}
	\end{align}
	Recall that we denote by $\ec[k]{\cdot}$ expectation conditional on $x^k$. Note that the gradient estimate is conditionally unbiased, i.e., that $\ec[k]{\nabla f (x^k; \xi)} = \frac{1}{n} \sum_{i=1}^{n} \nabla f_i (x^k) = \nabla f(x^k)$. Hence, taking conditional expectation in \eqref{eq:sgd-proof-2} and using unbiasedness we have
	\begin{align}
		&\ecn[k]{x^k - \gamma \nabla f (x^k; \xi) - (x^\ast - \gamma \nabla f (x^\ast))} \notag\\
		&= \sqn{x^k - x^\ast} - 2 \gamma \ev{x^k - x^\ast, \nabla f (x^k) - \nabla f (x^\ast)} + \gamma^2 \ecn[k]{\nabla f (x^k; \xi) - \nabla f (x^\ast)}. \label{eq:sgd-proof-3}
	\end{align}
	By the convexity of $f$ we have
	\[
		\ev{x^k - x^\ast, \nabla f(x^k) - \nabla f(x^\ast)} \geq D_{f} (x^k, x^\ast).
	\]
	Furthermore, we may estimate the third term in \eqref{eq:sgd-proof-3} by first using the fact that $\sqn{x + y} \leq 2 \sqn{x} + 2 \sqn {y}$ for any two vectors $x, y \in \R^d$
	\begin{align*}
		\ec[k]{\sqn{\nabla f (x^k; \xi) - \nabla f(x^\ast)}} &\leq 2 \ecn[k]{\nabla f (x^k; \xi) - \nabla f (x^\ast; \xi)}\\
		&\qquad + 2 \ecn[k]{\nabla f (x^\ast; \xi) - \nabla f(x^\ast)} \nonumber \\
		&= 2 \ecn[k]{\nabla f (x^k; \xi) - \nabla f (x^\ast; \xi)} + 2 \sigmaesc^2.
	\end{align*}
	We now use that by the $L_{\max}$-smoothness of $f_i$ we have that 
	\[ \sqn{\nabla f_i (x^k) - \nabla f_{i} (x^\ast)} \leq 2 L_{\max} \cdot D_{f_i} (x^k, x^\ast). \] 
	Hence
	\begin{align}
		\ecn[k]{\nabla f (x^k; \xi) - \nabla f (x^\ast; \xi)} &= \frac{1}{n} \sum_{i=1}^{n} \sqn{\nabla f_i (x^k) - \nabla f_{i} (x^\ast)} \nonumber \\
		&\leq \frac{2 L_{\max}}{n} \sum_{i=1}^{n} \left [ f_i (x^k) - f_i (x^\ast) - \ev{\nabla f_i (x^\ast), x^k - x^\ast} \right ] \nonumber \\
		&= 2 L_{\max} \left [ f(x^k) - f(x^\ast) - \ev{\nabla f(x^\ast), x^k - x^\ast} \right ] \nonumber \\
		&= 2 L_{\max} D_{f} (x^k, x^\ast).
		\label{eq:sgd-proof-6}
	\end{align}
	Combining equations \eqref{eq:sgd-proof-3}--\eqref{eq:sgd-proof-6} we obtain
	\begin{align*}
			\ecn[k]{x^k - \gamma \nabla f (x^k; \xi) - (x^\ast - \gamma \nabla f(x^\ast))} &\leq \sqn{x^k - x^\ast} - 2 \gamma \br{1 - 2 \gamma L_{\max}} D_{f} (x^k, x^\ast)\\
			&\qquad + 2 \gamma^2 \sigmaesc^2.
	\end{align*}
	Since $\gamma \leq \frac{1}{2 L_{\max}}$ by assumption we have that $1 - 2 \gamma L_{\max} \geq 0$. Since $D_{f} (x^k, x^\ast) \geq 0$ by the convexity of $f$ we arrive at
	\[
		\ecn[k]{x^k - \gamma \nabla f (x^k; \xi) - (x^\ast - \gamma \nabla f(x^\ast))} \leq \sqn{x^k - x^\ast} + 2 \gamma^2 \sigmaesc^2.
	\]
	Taking unconditional expectation and combining \eqref{eq:dec-step-1} with the last equation we have
	\begin{align*}
		\ecn{x^{k+1} - x^\ast} &\leq \frac{1}{1 + 2 \gamma \mu} \br{ \ecn{x^k - x^\ast} + 2 \gamma^2 \sigmaesc^2 } \\
		&= \frac{1}{1 + 2 \gamma \mu} \ecn{x^k - x^\ast} + \frac{2 \gamma^2 \sigmaesc^2}{1 + 2 \gamma \mu} \\
		&\leq \frac{1}{1 + 2 \gamma \mu} \ecn{x^k - x^\ast} + 2 \gamma^2 \sigmaesc^2.
	\end{align*}
	To simplify this further, we use that for any $x \leq \frac{1}{2}$ we have that $\frac{1}{1 + 2x} \leq 1-x$ and that $\gamma \mu \leq \frac{\mu}{2 L_{\max}} \leq \frac{1}{2}$, hence
	\begin{align*}
		\ecn{x^{k+1} - x^\ast} 
		\leq \br{1 - \gamma \mu} \ecn{x^k - x^\ast} + 2 \gamma^2 \sigmaesc^2.
	\end{align*}
	Recursing the above inequality for $K$ steps yields
	\begin{align*}
		\ecn{x^K - x^\ast} 
		&\leq \br{1 - \gamma \mu}^K \sqn{x^0 - x^\ast} + 2 \gamma^2 \sigmaesc^2 \br{ \sum_{k=0}^{K-1} \br{1 - \gamma \mu}^k } \\
		&\leq \br{1 - \gamma \mu}^K \sqn{x^0 - x^\ast} + 2 \gamma^2 \sigmaesc^2 \br{ \sum_{k=0}^{\infty} \br{1 - \gamma \mu}^k } \\
		&= \br{1 - \gamma \mu}^K \sqn{x^0 - x^\ast} + \frac{2 \gamma \sigmaesc^2}{\mu}.
	\end{align*}
\end{proof}

\section{Proofs for Decreasing Stepsize}

We first state and prove the following algorithm-independent lemma. This lemma plays a key role in the proof of \Cref{thm:psi-strongly-cvx-dec-stepsizes} and is heavily inspired by the stepsize schemes of \cite{stich2019unified} and \cite{Khaled2020} and their proofs. 

\begin{lemma}
	\label{lemma:decreasing-stepsizes-recursion-solution}
	Suppose that there exist constants $a, b, c \geq 0$ such that for all $\gamma_k \leq \frac{1}{b}$ we have
	\begin{equation} 
		\label{eq:dec-stepsizes-recursion-init}
		\br{1 + \gamma_k a n} r^{k+1} 
		\leq r^{k} + \gamma_k^3 c.
	\end{equation}
	Fix $T > 0$. Let $k_0 = \ceil{\frac{T}{2}}$. Then choosing stepsizes $\gamma_k > 0$ by
	\[
		\gamma_{k} = 
		\begin{cases}
			\frac{1}{b}, & \text { if } k \leq k_0 \text { or } T \leq \frac{b}{a n}, \\
			\frac{7}{a n \br{s + k - k_0}} & \text{ if } k > k_0 \text { and } T > \frac{b}{a n},
		\end{cases}	
	\]
	where $s = \frac{7b}{2 an}$. Then
	\[ r^{T} \leq \exp\br{-\frac{n T}{2 \br{b/a + n}}} r^0 + \frac{1421 c}{a^3 n^3 T^2}. \]
\end{lemma}
\begin{proof}
	If $T \leq \frac{7b}{an}$, then we have $\gamma_k = \gamma = \frac{1}{b}$ for all $k$. Hence recursing we have,
	\begin{align*}
		r^{T} &\leq \br{1 + \gamma a n}^{-T} r^0 + \frac{\gamma^3 c}{\gamma a n} = \br{1 + \gamma a n}^{-T} r^0 + \frac{\gamma^2 c}{a n}.
	\end{align*}
	Note that $\frac{1}{1+x} \leq \exp(-\frac{x}{1+x})$ for all $x$, hence
	\begin{align*}
		r^{T} &\leq \exp\br{ - \frac{\gamma a n T}{1 + \gamma a n} } r_{0} + \frac{\gamma^2 c}{an}
	\end{align*}
	Substituting for $\gamma$ yields
	\begin{align*}
		r^{T} &\leq \exp\br{-\frac{nT}{b/a + n}} r_{0} + \frac{c}{b^2 a n}.
	\end{align*}
	Note that by assumption we have $\frac{1}{b} \leq \frac{7}{T a n}$, hence
	\begin{align}
		\label{eq:dec-step-1}
		r^{T} &\leq \exp\br{-\frac{nT}{b/a + n}} r_{0} + \frac{49 c}{T^2 a^3 n^3}.
	\end{align}
	If $T > \frac{7 b}{an}$, then we have for the first phase when $k \leq k_0$ with stepsize $\gamma_k = \frac{1}{b}$ that
	\begin{align}
		\label{eq:dec-step-2}
		r^{k_0} &\leq \exp\br{-\frac{nk_0}{b/a + n}} r^0 + \frac{c}{b^2 a n} \leq \exp\br{-\frac{nT}{2(b/a + n)}} r_{0} + \frac{c}{b^2 a n}.
	\end{align}
	Then for $k > k_0$ we have
	\[
		\br{1 + \gamma_k a n} r^{k+1} \leq r^{T} + \gamma_k^3 c = r^{T} + \frac{7^3 c}{a^3 n^3 \br{s + k - k_0}^3}.
	\]
	Multiplying both sides by $(s+k - k_0)^3$ yields 
	\begin{equation}
		\label{eq:dec-step-3}
		\br{s + k - k_0}^3 \br{1 + \gamma_k a n} r^{k+1} \leq \br{s + k - k_0}^3 r^{T} + \frac{7^3 c}{a^3 n^3}.
	\end{equation}
	Note that because $k$ and $k_0$ are integers and $k > k_0$, we have that $k - k_0 \geq 1$ and therefore $s + k - k_0 \geq 1$. We may use this to lower bound the multiplicative factor in the left hand side of \eqref{eq:dec-step-3} as
	\begin{align}
		\br{s + k - k_0}^3 \br{1 + \gamma_k a n} &= \br{s + k - k_0}^3 \br{1 + \frac{7}{s + k - k_0}} \nonumber \\
		&= \br{s + k - k_0}^3 + 7 \br{s + k - k_0}^2 \nonumber \\
		&= \br{s + k - k_0}^3 + 3 \br{s + k - k_0}^2 + 3 \br{s + k - k_0}^2 \notag\\
		&\qquad + \br{s + k - k_0}^2 \nonumber \\
		&\geq \br{s + k - k_0}^3 + 3 \br{s + k - k_0}^2 + 3 \br{s + k - k_0} + 1 \nonumber \\
		\label{eq:dec-step-4}
		&= \br{s + k + 1 - k_0}^3.
	\end{align}
	Using \eqref{eq:dec-step-4} in \eqref{eq:dec-step-3} we obtain 
	\[
		\br{s + k + 1 - k_0}^3 r^{k+1} \leq \br{s + k - k_0}^3 r^k + \frac{7^3 c}{a^3 n^3}.
	\]
	Let $w^{k} = \br{s + k - k_0}^3$. Then we can rewrite the last inequality as
	\[
		w^{k+1} r^{k+1} - w^k r^k \leq \frac{7^3 c}{a^3 n^3}.
	\]
	Summing up and telescoping from $k=k_0$ to $T$ yields
	\[
		w^{T} r^{T} \leq w^{k_0} r^{k_0} + \frac{7^3 c}{a^3 n^3} \br{T - k_0}.
	\]
	Note that $w^{k_0} = s^3$ and $w^{T} = \br{s + T - k_0}^3$. Hence,
	\begin{align*}
		r^{T} &\leq \frac{s^3}{\br{s + T - k_0}^3} r^{k_0} + \frac{7^3 c}{a^3 n^3 \br{s + T - k_0}^2} \frac{T - k_0}{s + T - k_0} \\
		&\leq \frac{s^3}{\br{s + T - k_0}^3} r^{k_0} + \frac{7^3 c}{a^3 n^3 \br{s + T - k_0}^2}.
	\end{align*}
	Since we have $s + T - k_0 \geq T - k_0 \geq T/2$, it holds
	\begin{equation}
		\label{eq:dec-step-5}
		r^{T} \leq \frac{8 s^3}{T^3} r^{k_0} + \frac{4 \cdot 7^3 c}{a^3 n^3 T^2}.
	\end{equation}
	The bound in \eqref{eq:dec-step-2} can be rewritten as
	\[
		\frac{s^3}{T^3} r^{k_0} \leq \frac{s^3}{T^3} \exp\br{-\frac{n T}{2 \br{b/a + n}}} r_{0} + \frac{s^3 c}{b^2 a n T^3}.
	\]
	We now rewrite the last inequality, use that $T > 2s$ and further use the fact that $s = \frac{7b}{2an}$:
	\begin{align}
		\frac{s^3}{T^3} r^{k_0} &\leq \underbrace{\br{\frac{s}{T}}^3}_{\leq 1/8} \exp\br{-\frac{n T}{2 \br{b/a + n}}} r_{0} + \frac{s^2 c}{b^2 a n T^2} \underbrace{\br{\frac{s}{T}}}_{\leq 1/2} \nonumber \\
		&\leq \frac{1}{8} \exp\br{-\frac{n T}{2 \br{b/a + n}}} r_{0} + \frac{s^2 c}{2 b^2 a n T^2} \nonumber \\
		\label{eq:dec-step-6}
		&= \frac{1}{8} \exp\br{-\frac{n T}{2 \br{b/a + n}}} r_{0} + \frac{7^2 c}{8 a^3 n^3 T^2}.
	\end{align}
	Plugging in the estimate of \eqref{eq:dec-step-6} into \eqref{eq:dec-step-5} we obtain
	\begin{align}
		r^{T} &\leq \exp\br{-\frac{n T}{2 \br{b/a + n}}} r^0 + \frac{7^2 c}{a^3 n^3 T^2} + \frac{4 \cdot 7^3 c}{a^3 n^3 T^2} \nonumber \\
		\label{eq:dec-step-7}
		&= \exp\br{-\frac{n T}{2 \br{b/a + n}}} r^0 + \frac{1421 c}{a^3 n^3 T^2}.
	\end{align}
	Taking the maximum of \eqref{eq:dec-step-1} and \eqref{eq:dec-step-7} we see that for any $T > 0$ we have
	\[
		r^{T} \leq \exp\br{-\frac{n T}{2 \br{b/a + n}}} r^0 + \frac{1421 c}{a^3 n^3 T^2}.
	\]
\end{proof}

\subsection{Proof of \Cref{thm:psi-strongly-cvx-dec-stepsizes}}
\begin{proof}
	Start with Lemma~\ref{lemma:inner-loop} with $\lambda = 0$, $L=L_{\max}$, and $\gamma=\gamma_k$,
	\[ \ecn{x^k_{i+1} - x^\ast_{i+1}} \leq \ecn{x^k_i - x^\ast_i} - 2 \gamma \br{1 - \gamma L_{\max}} \ec{D_{f_{\pi_i}}(x^k_i, x^\ast) } + 2 \gamma_k^3 \sigmarr^2. \]
	Since $\gamma \leq \frac{1}{L_{\max}}$ and $D_{f_{\pi}} (x^k_i, x^\ast)$ is nonnegative we may simplify this to
	\[ \ecn{x^k_{i+1} - x^\ast_{i+1}} \leq \ecn{x^k_i - x^\ast_i} + 2 \gamma_k^3 \sigmarr^2. \]
	Unrolling this recursion for $n$ steps we get
	\[ 
		\ecn{x^k_{n} - x^\ast_n} 
		\leq \ecn{x^k_0 - x^\ast_0} + 2 n \gamma_k^3 \sigmarr^2. 
	\]
	By Lemma~\ref{prop:prox-contraction} and a similar reasoning to Theorem~\ref{thm:psi-strongly-convex-f-convex} we have
	\[ \br{1 + 2 \gamma_k \mu n} \ecn{x^{k+1} - x^\ast} \leq \ecn{x^k - x^\ast} + 2 \gamma_k^3 \sigmarr^2. \]
	We may then use Lemma~\ref{lemma:decreasing-stepsizes-recursion-solution} to obtain that
	\begin{align*}
		\ecn{x^T - x^\ast} &\leq \exp\br{-\frac{n T}{2(L_{\max}/\mu + n)}} \sqn{x^0 - x^\ast} + \frac{356 \sigmarr^2}{\mu^3 n^2 T^2} \\
		&= \mathcal{O}\br{ \exp\br{-\frac{n T}{\kappa + 2n}} \sqn{x^0 - x^\ast} + \frac{\sigmarr^2}{\mu^3 n^2 T^2} }.
	\end{align*}
\end{proof}

\section{Proof of \Cref{thm:IS} for Importance Resampling}
\begin{proof}
	We show that $N\le 2n$ as the rest of the theorem's claim trivially follows from \Cref{thm:psi-strongly-convex-f-convex}. Firstly, note that for any number $a\in\mathbb{R}$ we have $\lceil a \rceil \le a+1$. Therefore,
	\[
		N
		= \sum_{i=1}^n \left\lceil\frac{L_i}{\Lave}\right \rceil
		\le \sum_{i=1}^n \left(\frac{L_i}{\Lave}+1\right)
		= n + \sum_{i=1}^n \frac{L_i}{\Lave}
		= 2n.
	\]
\end{proof}

\section{Proofs for Federated Learning}
\subsection{Lemma for the extended proximal operator}
\begin{lemma}
    \label{lem:ext-proximal-operator}
    Let $\psi_C$ be the consensus constraint and $R$ be a closed convex proximable function. Suppose that $x_1, x_2, \ldots, x_M$ are all in $\R^d$. Then,
    \[ \prox_{\gamma (R + \psi_C)} (x_1, \ldots, x_M) = \prox_{\frac{\gamma}{M}R}(\overline x), \]
    where $\overline x = \frac{1}{M}\sum_{m=1}^M x_m$.
\end{lemma}
\begin{proof}
    We have,
    \[
        \prox_{\gamma(R+\psi_C)}(x_1, \dotsc, x_M)
        = \begin{pmatrix}
            \prox_{\frac{\gamma}{M}R}(\overline x)\\
            \vdots \\
            \prox_{\frac{\gamma}{M}R}(\overline x)
        \end{pmatrix}
        \quad \text{with} \quad \overline x = \frac{1}{M}\sum_{m=1}^M x_m.
    \]
    This is a simple consequence of the definition of the proximal operator. Indeed, the result of $\prox_{\gamma (R+\psi_C)}$ must have blocks equal to some vector $z$ such that
    \begin{align*}
        z
        &= \argmin_x \left\{\gamma R(x)+ \frac{1}{2}\sum_{m=1}^M\|x-x_m\|^2 \right\}\\
        &= \argmin_x \left\{\gamma R(x)+ \frac{1}{2}\sum_{m=1}^M \bigl(\|x-\overline x\|^2 + 2\<x-\overline x, \overline x - x_m>) + \|\overline x - x_m\|^2 \bigr)\right\}\\
        &= \argmin_x \left\{\gamma R(x)+ \frac{1}{2}M\|x-\overline x\|^2 \right\}
        \quad = \prox_{\frac{\gamma}{M}R}(\overline x).
    \end{align*}
\end{proof}

\subsection{Proof of \Cref{lem:fed_reform_properties}}
\begin{proof}
	Given some vectors $\xx, \yy\in\R^{d\cdot M}$, let us use their block representation $\xx=(x_1^\top,\dotsc, x_M^\top)^\top$, $\yy=(y_1^\top,\dotsc, y_M^\top)^\top$. Since we use the Euclidean norm, we have
	\[
		\|\nabla f_i(\xx)- \nabla f_i(\yy)\|^2
		= \sum_{m=1}^M \|\nabla f_{mi}(x_m)-\nabla f_{mi}(y_m)\|^2
		\le \sum_{m=1}^M L_i^2\|x_m-y_m\|^2
		= L_i^2\|\xx - \yy\|^2.
	\]
	We can obtain a lower bound by doing the same derivation and applying strong convexity instead of smoothness:
	\[
		\sum_{m=1}^M\|\nabla f_{mi}(x_m)-\nabla f_{mi}(y_m)\|^2
		\ge \mu^2 \sum_{m=1}^M\|x_m - y_m\|^2 
		= \mu^2 \|\xx-\yy\|^2.
	\]
	Thus, we have $\mu\|\xx-\yy\|\le \|\nabla f_i(\xx)- \nabla f_i(\yy)\|\le L_i\|\xx-\yy\|$, which is exactly $\mu$-strong convexity and $L_i$-smoothness of $f_i$.
\end{proof}

\subsection{Proof of \Cref{lem:fed_sigma}}
\begin{proof}
	By \Cref{thm:shuffling-radius-bound} we have
	\[
		\sigmarr^2 \le \frac{L_{\max}}{2}\Bigl(n^2\|\nabla f(\xx^\ast)\|^2 + \frac{n}{2}\sigmaesc^2\Bigr).
	\]
	Due to the separable structure of $f$, we have for the variance term
	\[
		n\sigmaesc^2 
		\eqdef  \sum_{i=1}^{n} \sqn{\nabla f_{i} (\xx^\ast) - \nabla f(\xx^\ast)}
		= \sum_{i=1}^{n}\sum_{m=1}^M \sqn{\nabla f_{mi} (x^\ast) - \frac{1}{n}\nabla F_m(x^\ast)}.
	\]
	The expression inside the summation is not exactly the variance due to the different normalization: $\frac{1}{n}$ instead of $\frac{1}{N_m}$. Nevertheless, we can expand the norm and try to get the actual variance:
	\begin{align*}
	     &\sum_{i=1}^{n}\sqn{\nabla f_{mi} (x^\ast) - \frac{1}{n}\nabla F_m(x^\ast)} \\
	     &=  \sum_{i=1}^{N_m}\biggl( \sqn{\nabla f_{mi} (x^\ast) - \frac{1}{N_m}\nabla F_m(x^\ast)} +\Bigl(\frac{1}{N_m} - \frac{1}{n}\Bigr)^2\sqn{\nabla F_m(x^\ast)}   \biggr)\\
	     &\quad + 2\sum_{i=1}^{N_m}\bigl\langle\nabla f_{mi}(x^\ast) - \frac{1}{N_m}\nabla F_m(x^\ast), \Bigl(\frac{1}{N_m} - \frac{1}{n}\Bigr)\nabla F_m(x^\ast) \bigr\rangle \\
	     &=  N_m\sigma_{m, \ast}^2 + N_m\Bigl(\frac{1}{N_m} - \frac{1}{n}\Bigr)^2\sqn{\nabla F_m(x^\ast)} \\
	     &\le n\sigma_{m, \ast}^2 + \sqn{\nabla F_m(x^\ast)}.
	\end{align*}
	Moreover, the gradient term has the same block structure, so
	\[
		n^2\|\nabla f(\xx^\ast)\|^2
		= n^2\biggl\|\frac{1}{n}\sum_{i=1}^n \nabla f_i(\xx^\ast) \biggr\|^2
		= \sum_{m=1}^M \sqn{\sum_{i=1}^n \nabla f_{mi}(x^\ast)}
		= \sum_{m=1}^M \|\nabla F_m(x^\ast)\|^2.
	\]
	Plugging the last two bounds back inside the upper bound on $\sigmarr^2$, we deduce the lemma's statement.
\end{proof}

\subsection{Proof of \Cref{thm:fed_hetero}}
\begin{proof}
	Since we assume that $N_1=\dotsb=N_M=n$, we have $\frac{N}{M}=n$ and the strong convexity constant of $\psi=\frac{N}{n}(R+\psi_C)$ is equal to $\frac{N}{n}\cdot \frac{\mu}{M}=\mu$. By applying \Cref{thm:psi-strongly-convex-f-convex} we obtain
	\[
		\ecn{\xx^T - \xx^\ast} \leq \br{1 + 2\gamma \mu n}^{-T} \sqn{\xx^0 - \xx^\ast} + \frac{ \gamma^2 \sigmarr^2}{\mu}.
	\]
	Since $\xx^T = \prox_{\gamma N(R+\psi_C)}(\xx^{T-1}_n)$, we have $\xx^T \in C$, i.e., all of its blocks are equal to each other and we have $\xx^T=((x^T)^\top,\dotsc, (x^T)^\top)^\top$. Since we use the Euclidean norm, it also implies
	\[
		\ecn{\xx^T - \xx^\ast}
		= M \|x^T - x^\ast\|^2.
	\]
	The same is true for $\xx^0$, so we need to divide both sides of  the upper bound on $\|\xx^T -\xx^\ast\|^2$ by $M$. Doing so together with applying \Cref{lem:fed_sigma} yields
	\begin{align*}
		\ecn{x^T - x^\ast} 
		&\leq \br{1 + 2\gamma \mu n}^{-T} \sqn{x^0 - x^\ast}  + \frac{ \gamma^2 \sigmarr^2}{M\mu} \\
		&\le \br{1 + 2\gamma \mu n}^{-T} \sqn{x^0 - x^\ast}  + \frac{ \gamma^2 L_{\max} }{M\mu} \sum_{m=1}^M\Bigl( \|\nabla F_m(x^\ast)\|^2 + \frac{n}{4}\sigma_{m, \ast}^2\Bigr) \\
		&= \br{1 + 2\gamma \mu n}^{-T} \sqn{x^0 - x^\ast}  + \frac{ \gamma^2 L_{\max}}{M\mu} \sum_{m=1}^M\Bigl(  \|\nabla F_m(x^\ast)\|^2 + \frac{N}{4M}\sigma_{m, \ast}^2\Bigr).
	\end{align*}
\end{proof}

\subsection{Proof of \Cref{thm:fed_iid}}
\begin{proof}
	According to~\Cref{lem:fed_reform_properties}, each $f_i$ is $\mu$-strongly convex and $L_{\max}$-smooth, so we obtain the result by trivially applying \Cref{thm:f-strongly-convex-psi-convex} and upper bounding $\sigmarr^2$ the same way as in the proof of~\Cref{thm:fed_hetero}.
\end{proof}

\section{Federated Experiments and Experimental Details}
We also compare the performance of \algname{FedRR} and \algname{Local SGD} on homogeneous (i.e., i.i.d.) data. Since \algname{Local SGD} requires smaller stepsizes to converge, it is significantly slower at initialization, as can be seen in Figure~\ref{fig:fed_rr}. \algname{FedRR}, however, does not need small initial stepsize and very quickly converges to a noisy neighborhood of the solution. The advantage is clear both from the perspective of the number of communication rounds and data passes.

To illustrate the severe impact of the number of local steps in \algname{Local SGD} we show results with different number of local steps. The blue line shows \algname{Local SGD} that takes the number of steps equivalent to full pass over the data by each node. The orange line takes 5 times fewer local steps. Clearly, the latter performs better in terms of communication rounds and local steps, making it clear that \algname{Local SGD} scales worse with the number of local steps. This phenomenon is well-understood and has been in discussed by~\cite{khaled2020tighter}.
\begin{figure}[t]
\centering
	\includegraphics[scale=0.38]{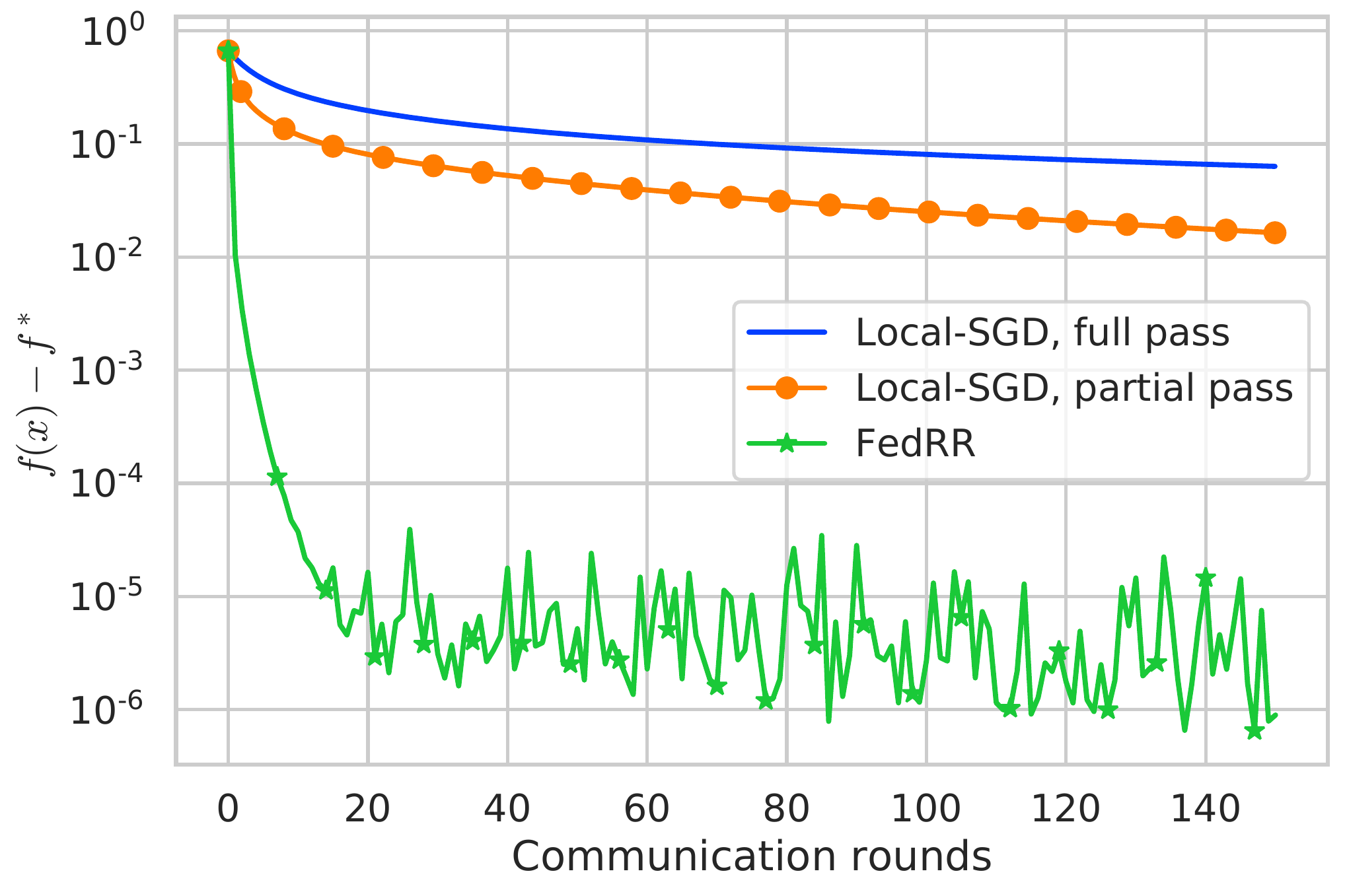}
	\includegraphics[scale=0.38]{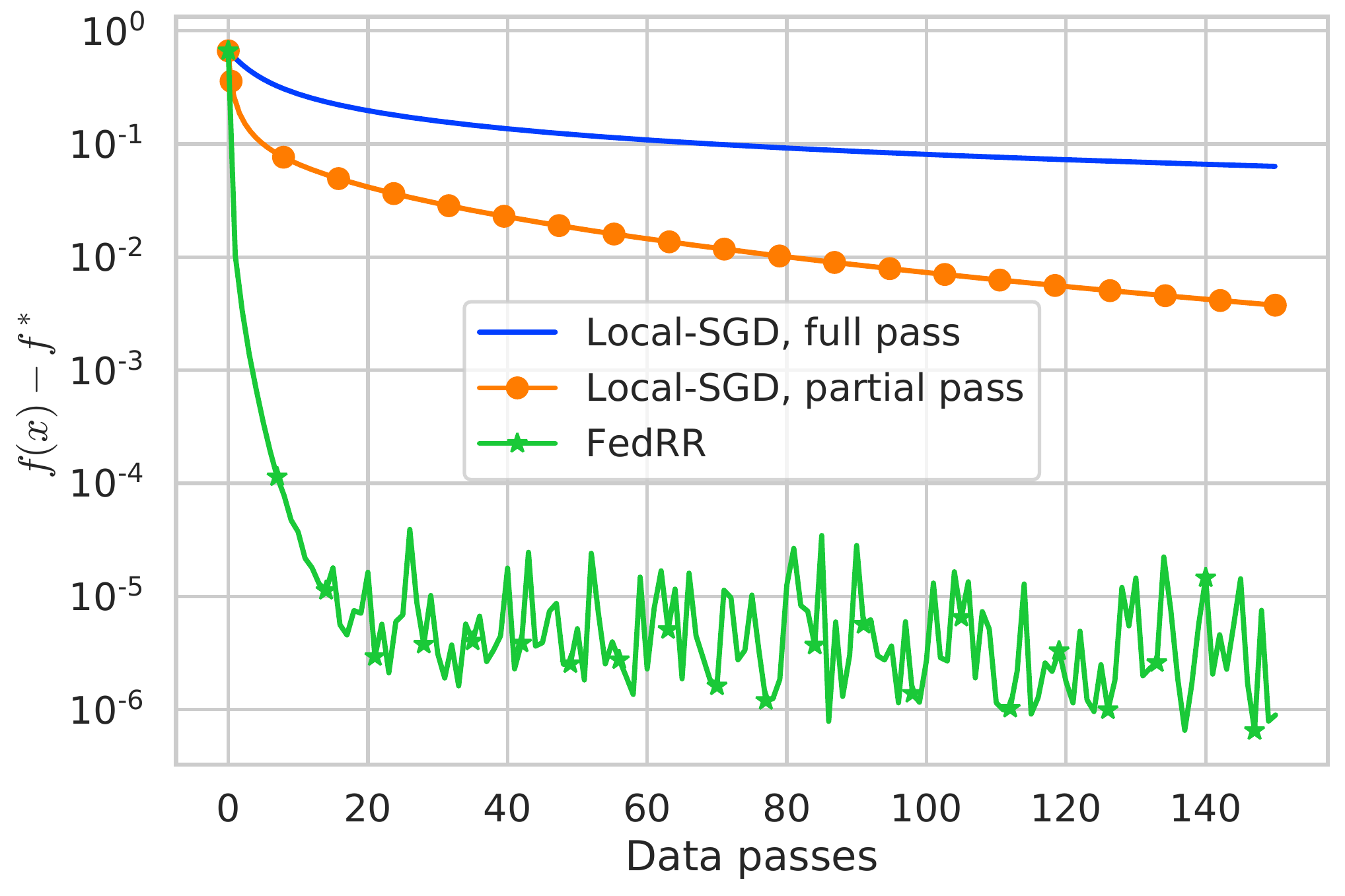}
\caption{Experimental results for parallel training. Left: comparison in terms of communication rounds, right: in terms of data passes}
\label{fig:fed_rr}
\end{figure}

\textbf{Implementation details.} For each $i$, we have $L_i=\frac{1}{4}\|a_i\|$. We set $\lambda_2=\frac{L}{N}$ and tune $\lambda_1$ to obtain a solution with less than 50\% coordinates (exact values are provided in the code). We use stepsizes decreasing as $\cO(\frac{1}{k})$ for all methods.  We use the `a1a' dataset for the experiment with $\ell_1$ regularization.

The experiment for the comparison of \algname{FedRR} and \algname{Local SGD} uses no $\ell_1$ regularization and $\lambda_2=\frac{L}{N}$. We choose the stepsizes according to the theory of \algname{Local SGD} and \algname{FedRR}. As per Theorem~3 in~\cite{khaled2020tighter}, the stepsizes for \algname{Local SGD} must satisfy $\gamma_k = \cO(1 / (LH))$, where $H$ is the number of local steps.  The parallelization of local runs is done using the Ray package\footnote{\href{https://ray.io/}{https://ray.io/}}. We use the `mushrooms' dataset for this experiment.

\textbf{Proximal operator calculation.} It is well-known (see, for instance, \cite{parikh2014proximal}) that the proximal operator for $\psi(x)=\lambda_1\|x\|_1 + \frac{\lambda_2}{2}\|x\|^2$ is given by
\[
	\prox_{\gamma\psi}(x) = \frac{1}{1+\gamma\lambda_2}\prox_{\gamma\lambda_1\|\cdot\|_1}(x),
\]
where the $j$-th coordinate of $\prox_{\gamma\lambda_1\|\cdot\|_1}(x)$ is
\[
	[\prox_{\gamma\lambda_1\|\cdot\|_1}(x)]_j
	= \begin{cases}
		\mathrm{sign}([x]_j)(|[x]_j|-\gamma\lambda_1), & \text{if } |[x]_j|\ge \gamma\lambda_1,\\
		0, &\text{otherwise}.
	\end{cases}
\]

\chapter{Appendix for Chapter~\ref{chapter:adaptive}}
\graphicspath{{adaptive/}}
\section{Missing Proofs}
Recall that in the proof of \Cref{th:main} we only showed boundedness of the iterates and complexity for minimizing $f(x)$. It remains to show that sequence $(x^k)_k$ converges to a solution.
For this, we  need some variation of the Opial lemma.
\begin{lemma}\label{opial-like}
    Let  $(x^k)_k$ and $(a_k)_k$ be two sequences in
    $\R^d$ and $\R_+$ respectively. Suppose that $(x^k)_k$ is bounded,
    its cluster points belong to $\cX \subset \R^d$ and it
    also holds that
\begin{equation}\label{fejer}
    \n{x^{k+1}-x}^2 + a_{k+1}\leq \n{x^k-x}^2 + a_k, \qquad \forall
    x\in \cX.
\end{equation}
Then $(x^k)_k$ converges to some element in $\cX$.
\end{lemma}
\begin{proof}
Let $\x_1$, $\x_2$ be any cluster points of $(x^k)_k$. Thus, there
    exist two subsequences $(x^{k_i})_i$ and $(x^{k_j})_j$ such that
    $x^{k_i}\to \x_1$ and $x^{k_j}\to \x_2$.
    Since $\n{x^k-x}^2 + a_k$ is nonnegative and bounded,  $\lim_{k\to \infty}(\n{x^k-x}^2 + a_k)$ exists for any
    $x\in \cX$. Let $x=\x_1$. This yields
\begin{align*}
    \lim_{k\to \infty}(\n{x^k-\x_1}^2 + a_k)&=\lim_{i\to
            \infty}(\n{x^{k_i}-\x_1}^2 + a_{k_i})=\lim_{i\to
                                            \infty}a_{k_i}\\
  & = \lim_{j\to
            \infty}(\n{x^{k_j}-\x_1}^2 + a_{k_j})=\n{\x_2-\x_1}^2 + \lim_{j\to
                                            \infty}a_{k_j}.
\end{align*}
Hence, $\lim_{i\to \infty} a_{k_i} = \lim_{j\to \infty}a_{k_j} +
\n{\x_1-\x_2}^2$. Doing the same with $x=\x_2$ instead of $x=\x_1$,
yields $\lim_{j\to \infty} a_{k_j} = \lim_{i\to \infty}a_{k_i} +
\n{\x_1-\x_2}^2$. Thus, we obtain that $\x_1=\x_2$, which
finishes the proof.
\end{proof}
Another statement that we need here is the following tightening of the convexity property.
\begin{lemma}[Theorem 2.1.5, \cite{Nesterov2013}]
    \label{lemma:coco}
    Let $\mathcal{C}$ be a closed convex set in  $\R^d$. If $f\colon \mathcal{C}\to \R $
    is convex and $L$-smooth, then $\forall x,y\in \mathcal{C}$ it holds
    \begin{align}\label{eq:smooth_and_convex}
    f(x)-f(y)-\lr{\nabla f(y),x-y}\geq \frac{1}{2L}\n{\nabla
            f(x)-\nabla f(y)}^2.
    \end{align}
\end{lemma}

\begin{proof}[\textbf{Proof of
    Theorem~\ref{th:main}}](\textit{Convergence of $(x^k)_k$})

Note that in the first part we have already proved that $(x^k)_k$ is bounded and that
$\nabla f$ is $L$-Lipschitz on  $\mathcal{C}=\clconv\{x^\star, x^0, x^1, \dots\}$. Invoking Lemma~\ref{lemma:coco}, we deduce that
\begin{equation}
    \label{eq:better_with_Lip}
    \gamma_k(f(x^\star) - f(x^k)) \geq  \gamma_k \lr{\nabla f(x^k), x^\star-x^{k}} +
    \frac{\gamma_k}{2L}\n{\nabla f(x^k)}^2.
\end{equation}
This indicates that instead of using
inequality~\eqref{eq:conv} in the proof of
Lemma~\ref{lemma:energy}, we could use a better
estimate~\eqref{eq:better_with_Lip}. However, we want to emphasize
that we did not assume that $\nabla f$ is globally Lipschitz, but rather
obtained Lipschitzness on $\mathcal{C}$ as an artifact of our
analysis. Clearly, in the end this improvement gives us an additional term
$\frac{\gamma_k}{L}\n{\nabla f(x^k)}^2$ in the left-hand side of
\eqref{eq:lemma_ineq}, that is
\begin{multline}
  \label{eq:lemma_ineq_appendix}
\n{x^{k+1}-x^\star}^2+ \frac 1 2 \n{x^{k+1}-x^k}^2  + 2\gamma_{k}(1+\th_{k})
 (f(x^k)-f^\star) + \frac{\gamma_k}{L}\n{\nabla f(x^k)}^2 \\ \leq \n{x^k-x^\star}^2  + \frac 1 2
 \n{x^k-x^{k-1}}^2    + 2\gamma_k \th_k (f(x^{k-1})-f^\star).
\end{multline}
Thus, telescoping \eqref{eq:lemma_ineq_appendix}, one
obtains that $\sum_{i=1}^k\frac{\gamma_k}{L}\n{\nabla f(x^k)}^2\leq D$.
As $\gamma_k\geq \frac{1}{2L}$, one has that $\nabla f(x^k)\to 0$. Now we
might conclude that all cluster points of $(x^k)_k$ are solutions of~\eqref{main}.

Let $\cX$ be the solution set of \eqref{main} and
$a_k = \frac 1 2 \n{x^k-x^{k-1}}^2 + 2\gamma_k \th_k(f(x^{k-1})-f^\star)$.
We want to finish the proof applying Lemma~\ref{opial-like}. To
this end, notice that inequality~\eqref{eq:lemma_ineq}
 yields \eqref{fejer}, since $\gamma_{k+1}\th_{k+1}\leq
 (1+\th_k)\gamma_k$. This completes the proof.
\end{proof}

\begin{proof}[\textbf{Proof of \Cref{th:strong}}]~\\
    First of all, we note that using the stricter inequality
    $\gamma_k\leq \sqrt{1+\frac{\th_{k-1}}{2}}\gamma_{k-1}$ does not change
    the statement of \Cref{th:main}. Hence, $x^k\to x^\star$ and there
    exist $\mu,L >0$ such that $f$ is $\mu$-strongly convex and
    $\nabla f$ is $L$-Lipschitz on
    $\mathcal{C}=\clconv\{x^\star, x^0, x^1, \dots\}$. Secondly, due to
    local strong convexity,
    $\n{\nabla f(x^k)-\nabla f(x^{k-1})}\geq \mu \n{x^k-x^{k-1}}$, and
    hence $\gamma_k \leq \frac{1}{2\mu}$ for $k\geq 1$.

    Now we tighten some steps in the analysis to
    improve bound~\eqref{eq:conv}. By strong convexity,
\begin{align*}
   \gamma_k \lr{\nabla f(x^k), x^\star-x^{k}} & \leq \gamma_k(f(x^\star) - f(x^k)) - \gamma_k\frac{\mu}{2}\|x^\star - x^k\|^2.
\end{align*}
By $L$-smoothness and bound $\gamma_k\le \frac{1}{2\mu}$,
\begin{align*}
   \gamma_k \lr{\nabla f(x^k), x^\star-x^{k}} & \leq \gamma_k(f(x^\star) - f(x^k)) - \gamma_k\frac{1}{2L}\|\nabla f(x^k)\|^2 \\
   &= \gamma_k(f^\star - f(x^k)) - \frac{1}{2L\gamma_k}\|x^{k+1}-x^k\|^2 \\
   &\le \gamma_k(f^\star - f(x^k)) - \frac{\mu}{L}\|x^{k+1}-x^k\|^2.
\end{align*}
Together, these two bounds give us
\begin{align*}
   \gamma_k \lr{\nabla f(x^k), x^\star-x^{k}} & \leq \gamma_k(f^\star - f(x^k)) -
                                         \gamma_k\frac{\mu}{4}\|x^k - x^\star\|^2 - \frac{\mu}{2L}\|x^{k+1}-x^k\|^2.
\end{align*}

We keep inequality~\eqref{eq:terrible_simple} and the rest of the proof as is.
Then the strengthen analog of~\eqref{eq:lemma_ineq} will be
\begin{align}
  &\n{x^{k+1}-x^\star}^2+ \frac 1 2\left(1 + \frac{2\mu}{L}\right) \n{x^{k+1}-x^k}^2 + 2\gamma_k
  (1+\th_k) (f(x^k)-f^\star) \nonumber \\
  \leq \, & \left(1 - \frac{\gamma_k\mu}{2}\right)\n{x^k-x^\star}^2  + \frac 1 2
            \n{x^k-x^{k-1}}^2   + 2\gamma_k \th_k (f(x^{k-1})-f^\star) \nonumber\\
   \leq \, & \left(1 - \frac{\gamma_k\mu}{2}\right)\n{x^k-x^\star}^2  + \frac 1 2
 \n{x^k-x^{k-1}}^2   + 2\gamma_{k-1} \left(1+\frac{\th_{k-1}}{2}\right) (f(x^{k-1})-f^\star),
            \label{eq:contraction}
\end{align}
where in the last inequality we used our new condition on $\gamma_k$.
Under the new update we have contraction in every term:
$1-\frac{\gamma_k\mu}{2}$ in the first,
$\frac{1}{1+2\mu/L}=1-\frac{2\mu}{L+2\mu}$ in the second and
$\frac{1+\th_{k-1}/2}{1+\th_{k-1}}=1-\frac{\th_{k-1}}{2(1+\th_{k-1})}$
in the last one.

To further bound the last contraction, recall that
$\gamma_k \in \left[\frac{1}{2L}, \frac{1}{2\mu}\right]$ for $k\ge 1$.
Therefore, $\th_k = \frac{\gamma_k}{\gamma_{k-1}} \ge \frac{1}{\kappa}$ for any $k>1$, where
$\kappa\eqdef \frac{L}{\mu}$. Since the function $\th\mapsto \frac{\th}{1+\th}$
monotonically increases with $\th>0$,  this implies
$\frac{\th_{k-1}}{2(1+\th_{k-1})}\ge \frac{1}{2(\kappa+1)}$ when
$k>2$. Thus, for the full energy $\Psi^{k+1}$ (the left-hand side
of~\eqref{eq:contraction}) we have
\begin{align*}
        \Psi^{k+1}
        &\le \left(1  - \min\left\{\frac{\gamma_k \mu}{2}, \frac{1}{2(\kappa+1)}, \frac{2\mu}{L+2\mu}\right\}\right) \Psi^k.
\end{align*}
Using simple bounds  $\frac{\gamma_k\mu}{2}\geq \frac{1}{4\kappa}$,
$\frac{2\mu}{L+2\mu}=\frac{2}{\kappa +2}\geq \frac{1}{4\kappa}$, and
$\frac{1}{2(\kappa +1)}\geq \frac{1}{4\kappa}$, we obtain $\Psi^{k+1}\le (1 - \frac{1}{4\kappa})\Psi^k$ for $k>2$.
This gives $\mathcal{O}\left(\kappa \log \frac{1}{\varepsilon}\right)$ convergence rate.
\end{proof}
\section{Extensions}
\subsection{More general update}
One may wonder how flexible the update for $\gamma_k$ in \Cref{alg:main}
is? For example, is it necessary to upper bound the stepsize with
$\sqrt{1+\th_k}\gamma_{k-1}$ and put $2$ in the denominator of
$\frac{\n{x^k-x^{k-1}}}{2\n{\nabla f(x^k)-\nabla f(x^{k-1})}}$?
\Cref{alg:general_update} that we present here  partially answers
 this question.
 \begin{algorithm}[t]
 \caption{\algname{Adaptive Gradient Descent} (general update)}
 \label{alg:general_update}
 \begin{algorithmic}[1]
     \State \textbf{Input:} $x^0 \in \R^d$, $\gamma_0>0$,
     $\th_0=+\infty$, $\a\in (0,1)$, $\b = \frac{1}{2(1-\a)}$, number of steps $K$
     \State  $x^1= x^0-\gamma_0\nabla f(x^0)$
        \For{$k = 1,2,\dots, K-1$}
        \State $\gamma_k = \min\Bigl\{
        \sqrt{\frac{1}{\beta}+\th_{k-1}}\gamma_{k-1},\frac{\a\n{x^{k}-x^{k-1}}}{\n{\nabla
        f(x^{k})-\nabla f(x^{k-1})}}\Bigr\}$
\State $x^{k+1} = x^k - \gamma_k \nabla f(x^k)$
\State $\th_k = \frac{\gamma_k}{\gamma_{k-1}}$
        \EndFor
 \end{algorithmic}
 \end{algorithm}
Obviously, \Cref{alg:main} is a particular case of
\Cref{alg:general_update} with $\a=\frac 12$ and $\b =1$.

\begin{theorem}\label{th:gen_update}
    Suppose that $f\colon \R^d\to \R$ is convex with locally Lipschitz
    gradient $\nabla f$. Then $(x^k)_k$ generated by~\Cref{alg:general_update} converges to a solution of \eqref{main}
    and we have that
\[
	f(\hat x^K)-f^\star 
	\leq \frac{D}{2S_K}
	=\mathcal{O}\Bigl(\frac{1}{K}\Bigr),\]
where
\begin{align*}
	\hat x^K &=  \frac{\gamma_K(1+\th_K\b)x^K +
        \sum_{i=1}^{K-1}w_i
           x^i}{S_K},\\
    w_i &= \gamma_i(1+\th_i\b)-\gamma_{i+1}\th_{i+1}\b,\\
  S_K &= \gamma_K(1+\th_K\b) + \sum_{i=1}^{K-1}w_i 
  = \sum_{i=1}^K \gamma_i + \gamma_1\th_1\b,
\end{align*}
and $D$ is a constant
that explicitly depends on the initial data and the solution set.
\end{theorem}
\begin{proof}
    Let $x^\star$ be arbitrary solution of \eqref{main}.  We note
    that equations~\eqref{eq:norms} and~\eqref{dif_x} hold for any
    variant of GD, independently of $\gamma_k$, $\a$, $\b$.  With the new
    rule for  $\gamma_k$, from \eqref{dif_x} it follows
\begin{align*}
        \|x^{k+1}-x^k\|^2
        &\le 2\gamma_k\th_k(f(x^{k-1})-f(x^k))-\|x^{k+1}-x^k\|^2\\
        &\qquad +2\gamma_k\|\nabla f(x^k)-\nabla f(x^{k-1})\|\| x^k-x^{k+1}\| \\
        &\le 2\gamma_k\th_k(f(x^{k-1})-f(x^k))-(1-\alpha)\|x^{k+1}-x^k\|^2+\alpha\|x^k-x^{k-1}\|^2,
\end{align*}
which, after multiplication by $\beta$ and reshuffling the terms, becomes
\begin{align*}
        \beta(2 - \alpha)\|x^{k+1}-x^k\|^2+2\beta\gamma_k\th_k (f(x^k)-f^\star)
        \le \alpha\beta\|x^k-x^{k-1}\|^2+ 2\beta\gamma_k\th_k (f(x^{k-1})-f^\star).
\end{align*}
Adding \eqref{eq:norms} and the latter inequality gives us
\begin{align*}
        &\|x^{k+1}-x^\star\|^2 + 2\gamma_k(1+\th_k\beta)(f(x^k)-f^\star) + (2\beta-\alpha\beta-1)\|x^{k+1}-x^k\|^2\\
        &\le \|x^k-x^\star\|^2+2\gamma_k\th_k\beta(f(x^{k-1})-f^\star) +\alpha\beta\|x^{k+1}-x^k\|^2.
\end{align*}
Notice that by $\b = \frac{1}{2(1-\a)}$, we have $2\b-\a\b -1= \a\b$ and
hence,
\begin{align*}
        &\|x^{k+1}-x^\star\|^2 + 2\gamma_k(1+\th_k\beta)(f(x^k)-f^\star) +\a\b\|x^{k+1}-x^k\|^2\\
        &\le \|x^k-x^\star\|^2+2\gamma_k\th_k\beta(f(x^{k-1})-f^\star) +\alpha\beta\|x^{k+1}-x^k\|^2.
\end{align*}
As a sanity check, we can see that with $\a=\frac{1}{2}$ and $\b
=1$, the above inequality coincides with~\eqref{eq:lemma_ineq}.

Telescoping this inequality, we deduce
\begin{align}\label{eq:telescope2}
        \n{x^{k+1}-x^\star}^2&+ \a\b \n{x^{k+1}-x^k}^2  + 2\gamma_k
  (1+\th_k\b) (f(x^k)-f^\star) \notag \\ &\qquad +
  2\sum_{i=1}^{k-1}[\gamma_i(1+\th_i\b)-\gamma_{i+1}\th_{i+1}\b](f(x^i)-f^\star) \notag\\
  &\leq  \n{x^1-x^\star}^2  + \a\b \n{x^1-x^{0}}^2   + 2\gamma_1 \th_1\b [f(x^{0})-f^\star]\eqdef D.
\end{align}
Note that because of the way we defined stepsize,
$\gamma_i(1+\th_i\b)-\gamma_{i+1}\th_{i+1}\b\geq 0$. Thus, the sequence
$(x^k)_k$ is bounded. Since $\nabla f$ is locally Lipschitz, it is
Lipschitz continuous on bounded sets. Let $L$ be a Lipschitz constant
of $\nabla f$ on a bounded set
$\mathcal{C}=\clconv\{x^\star,x^1,x^2,\dots\}$.

If $\a \leq \frac 1 2 $, then $\frac{1}{\b}>1$ and similarly
to~\Cref{th:main} we might conclude that $\gamma_k\geq\frac{\a}{L}$ for
all $k$. However, for the case $\a > \frac 1 2$ we cannot do
this. Instead, we prove  that $\gamma_k\ge \frac{2\a(1-\a)}{L}$, which suffices
 for our purposes.

 Let $m, n\in \N$ be the smallest numbers such that
 $\b^{-\frac{1}{2^m}}\ge 1 - \frac{1}{2\b}$ and
 $(1+\frac{1}{\b})^{\frac n 2}\ge \b$. We want
 to prove that for any $k$ it holds $\gamma_k\ge \frac{2\a(1-\a)}{L}$ and among every $m+n+1$ consecutive elements
 $\gamma_k,\gamma_{k+1},\dots,\gamma_{k+m+n}$ at least one is no less than
 $\frac{\a}{L}$. We shall prove this by induction. First, note that
 the second bound always satisfies
 $\frac{\a\n{x^k-x^{k-1}}}{\n{\nabla f(x^k)-\nabla f(x^{k-1})}}\geq
 \frac{\a}{L}$ for all $k$,  which also
 implies that $\gamma_1\geq \frac{\a}{L}$. If for all $k$ we have $\gamma_k\ge \frac{\a}{L}$, then we are done. Now assume that $\gamma_{k-1}\geq
 \frac{\a}{L}$ and $\gamma_k<\frac{\a}{L}$ for some $k$. Choose the largest $j$ (possibly infinite) such that the second bound is not
 active for $\gamma_k, \dotsc, \gamma_{k+j-1}$, i.e.,  $\gamma_{k+i}=\sqrt{\frac{1}{\b} + \th_{k+i-i}}\gamma_{k+i-1}$ for $i<j$.

Let us prove that $\gamma_k, \dotsc, \gamma_{k+j-1}\ge \frac{2\a(1-\a)}{L}$.  The definition of $j$ yields $\th_{k+i}=\sqrt{\frac{1}{\b} +
     \th_{k+i-1}}$ for all
 $i=0,\dots, j-1$. Recall that $\b > 1$, and
 thus,
 \[\th_k\geq  \b^{-\frac 1 2},\quad \th_{k+1}\geq \sqrt{\frac{1}{\b}
         + \sqrt{\frac{1}{\b}}}\geq \b^{-\frac{1}{4}},\quad \dotsc,
     \quad \th_{k+i}\geq \sqrt{\frac{1}{\b} + \b^{-\frac{1}{2^i}}} \geq
    \b^{-\frac{1}{2^{i+1}}}\]
 for all $i<j$. Now it remains to notice that for any $i<j$
 \[\frac{\gamma_{k+i}}{\gamma_{k-1}} = \th_{k}\th_{k+1}\dots
 \th_{k+i} \geq \beta^{-\frac{1}{2}-\frac{1}{4}-\dotsb -\frac{1}{2^{i+1}} }\geq \frac{1}{\beta}=2(1-\a)\]
and hence $\gamma_{k+i}\geq 2(1-\a)\gamma_{k-1}\geq \frac{2\a(1-\a)}{L}$. If
$j\leq m+n$, then at $(k+j)$-th iteration the second bound is active,
i.e., $\gamma_{k+j}\ge \frac{\a}{L}$, and we are done with the other claim as well. Otherwise, note
\[
        \th_{k+m-1} \ge \b^{-\frac{1}{2^m}} \ge 1 - \frac{1}{2\b},
\]
 so $\th_{k+m}=\sqrt{\frac{1}{\b}+\th_{k+m-1}}\ge \sqrt{\frac{1}{\b}+1 - \frac{1}{2\b}}=\sqrt{1 + \frac{1}{2\b}}$ and for any $i\in [m, j-2]$ we have $\th_{k+i+1}= \sqrt{\frac{1}{\b}+\th_{k+i}}\ge \sqrt{\frac{1}{\b}+1}$. Thus,
 \[
 \gamma_{k+m+n} =
 \gamma_{k-1}\Bigl(\prod_{l=k}^{k+m-1}\th_l\Bigr)\Bigl(\prod_{l=k+m}^{k+m+n}\th_l\Bigr)\ge
 \gamma_{k-1}\frac{1}{\b}\sqrt{1+\frac{1}{2\b}}\Bigl(1+\frac{1}{\b}\Bigr)^{\frac
 n 2}\ge \gamma_{k-1} \ge \frac{\a}{L},
 \]
 so we have shown the second claim too.

 To conclude, in both cases
$\a\leq \frac 12$ and $\a>\frac 12$, we have $S_k = \Omega(k)$.

Applying the Jensen inequality for the sum of all terms
$f(x^i)-f^\star$ in the left-hand side of~\eqref{eq:telescope2}, we obtain
\[\frac D 2 \geq \frac{\text{LHS of \eqref{eq:telescope2}}}{2} \geq S_K (f(\hat
    x^K)-f^\star),\]
where $\hat x^K$ is defined in the statement of the
theorem. Finally, convergence of $(x^k)_k$ can be proved in a similar
way as \Cref{th:main}.
\end{proof}

\subsection{$f$ is $L$-smooth}
Often, it is known that $f$ is smooth and even some estimate for the
Lipschitz constant $L$ of $\nabla f$ is available. In this case, we
can use slightly larger steps, since instead of just convexity the stronger inequality in \Cref{lemma:coco} holds. To take advantage of it, we present a modified version of \Cref{alg:main} in
\Cref{alg:smooth}. Note that we have chosen to modify \Cref{alg:main} and
not its more general variant~\Cref{alg:general_update} only for simplicity.

\begin{algorithm}[t]
 \caption{\algname{Adaptive Gradient Descent} ($L$ is known)}
 \label{alg:smooth}
 \begin{algorithmic}[1]
     \State \textbf{Input:} $x^0 \in \R^d$, $\gamma_0 = \frac{1}{L}$,
     $\th_0=+\infty$, number of steps $K$
     \State $x^1 = x^0-\gamma_0\nabla f(x^0)$
     \For{$k = 1,2,\dots, K-1$}
        \State $L_k = \frac{\n{\nabla
        f(x^{k})-\nabla f(x^{k-1})}}{\n{x^{k}-x^{k-1}}}$
     \State  $\gamma_k = \min\left\{
        \sqrt{1+\th_{k-1}}\gamma_{k-1},\frac{1}{\gamma_{k-1}L^2}+\frac{1}{2L_k} \right\}$
\State $x^{k+1} = x^k - \gamma_k \nabla f(x^k)$
\State $\th_k = \frac{\gamma_k}{\gamma_{k-1}}$
 \EndFor
 \end{algorithmic}
 \end{algorithm}
\begin{theorem}\label{theorem:energy-L}
Let $f$ be convex and $L$-smooth. Then for $(x^k)_k$ generated
by Algorithm~\ref{alg:smooth} inequality~\eqref{eq:lemma_ineq} holds. As a corollary, it holds for some ergodic vector $\hat x^K$ that $f(\hat x^K)-f^\star=\mathcal{O}\left(\frac{1}{K}\right)$.
\end{theorem}
\begin{proof}
    Proceeding similarly as in Lemma~\ref{lemma:energy}, we have
  \begin{equation}
      \label{eq:2_simple-L}
          \|x^{k+1}- x^\star\|^2 = \|x^k - x^\star\|^2 - 2\gamma_k \<\nabla
            f(x^k), x^k-x^\star> + \|x^{k+1} - x^{k}\|^2.
\end{equation}

By convexity of $f$ and Lemma~\ref{lemma:coco},
\begin{align}\label{eq:conv2_simple-L}
   2\gamma_k \lr{\nabla f(x^k), x^\star-x^{k}} & \overset{\eqref{eq:smooth_and_convex}}{\leq} 2\gamma_k(f(x^\star) -
                                         f(x^k)-\frac{1}{2L}\n{\nabla
                                         f(x^k)}^2) \notag \\&= 2\gamma_k(f^\star -
                                         f(x^k))-\frac{1}{\gamma_k L}\n{x^{k+1}-x^k}^2.
\end{align}
As in~\eqref{dif_x}, we have
 \begin{equation}\label{dif_x-L}
      \|x^{k+1} -x^k\|^2  =  2\gamma_k \lr{\nabla
                                   f(x^k)-f(x^{k-1}), x^{k}-x^{k+1}} +
                                   2\gamma_k\lr{f(x^{k-1}), x^k-x^{k+1}}
      -  \n{x^{k+1}-x^k}^2.
\end{equation}
Again, instead of  using merely convexity of $f$, we combine it with
\Cref{lemma:coco}. This gives
    \begin{align}    \label{eq:terrible_simple-L}
          2\gamma_k\lr{\nabla f(x^{k-1}), x^k-x^{k+1}}
      &=  \frac{2\gamma_k}{\gamma_{k-1}}\lr{x^{k-1} - x^{k}, x^{k}-x^{k+1}}\notag
      \\ &=  2\gamma_k\th_k \lr{x^{k-1}-x^{k}, \nabla
        f(x^k)} \notag \\ & \overset{\eqref{eq:smooth_and_convex}}{\leq}
      2\gamma_k\th_k
    (f(x^{k-1})-f(x^k)) - \frac{\gamma_k\th_k}{L}\n{\nabla f(x^k)-\nabla f(x^{k-1})}^2.
\end{align}
Since now we have two additional terms
$\frac{1}{\gamma_kL}\n{x^{k+1}-x^k}^2$ and $\frac{\gamma_k\th_k}{L}\n{\nabla
    f(x^{k}) - \nabla f(x^{k-1})}^2$, we can do better than~\eqref{cs}. But first we need a simple, yet a bit tedious fact.
By our choice  of $\gamma_k$, in every iteration $\gamma_k\leq \frac{1}{\gamma_{k-1}L^2} +
\frac{1}{2L_k}$ with $L_k= \frac{\n{\nabla f(x^k)-\nabla f(x^{k-1})}}{\n{x^k-x^{k-1}}}
$. We want to show that it implies
\begin{equation}\label{so_much_trouble}
    2\left(\gamma_k - \frac{\sqrt{\th_k}}{L} \right)\leq \frac{1}{L_k},
\end{equation}
which is equivalent to $\gamma_k-\frac{\sqrt{\gamma_k}}{\sqrt{\gamma_{k-1}}L}-\frac{1}{2L_k}\le 0$. Nonnegative solutions of the quadratic inequality $t^2 -
\frac{t}{\sqrt{\gamma_{k-1}}L}-\frac{1}{2L_k}\leq 0$ are
\[0\leq t\leq \frac{1}{2\sqrt{\gamma_{k-1}}L}+\frac{1}{2}\sqrt{\frac{1}{\gamma_{k-1}L^2}+\frac{2}{L_k}}
    = \frac{1}{2\sqrt{\gamma_{k-1}}L}\left(1 + \sqrt{1 + \frac{2\gamma_{k-1}L^2}{L_k}}\right).\]
Let us prove that $\sqrt{\frac{1}{\gamma_{k-1}L^2} +
\frac{1}{2L_k}}$ falls into this segment and, hence, $\sqrt{\gamma_k}$ does as well. Using a simple inequality $4 + a\leq (1 + \sqrt{1+a})^2 $, for
$a>0$, we obtain
\[\frac{1}{\gamma_{k-1}L^2} + \frac{1}{2L_k} =
    \frac{1}{4\gamma_{k-1}L^2} \bigl(4 +
    \frac{2\gamma_{k-1}L^2}{L_k}\bigr)\leq \frac{1}{4{\gamma_{k-1}}L^2}\left(1 + \sqrt{1 + \frac{2\gamma_{k-1}L^2}{L_k}}\right)^2.\]
This confirms that~\eqref{so_much_trouble} is true.
Thus, by
Cauchy-Schwarz and Young's inequalities, one has
\begin{equation}\label{cs-L}
\begin{aligned}
  & 2\gamma_k \lr{\nabla f(x^k) -\nabla f(x^{k-1}), x^k - x^{k+1}} \\
  &\leq 2\gamma_k \n{\nabla f(x^k) -\nabla f(x^{k-1})} \n{x^k - x^{k+1}} \\ 
  &= 2\left(\gamma_k - \frac{\sqrt{\th_k}}{L}\right) \n{\nabla f(x^k) -\nabla f(x^{k-1})}\n{x^k -
      x^{k+1}} \\
    &\qquad + \frac{2\sqrt{\th_k}}{L}\n{\nabla f(x^k) -\nabla f(x^{k-1})}\n{x^k - x^{k+1}}\\
  & \overset{\eqref{so_much_trouble}}{\leq} \frac{1}{L_k}
 \n{\nabla f(x^k) -\nabla f(x^{k-1})} \n{x^k -
    x^{k+1}} + \frac{\gamma_k\th_k}{L}\n{\nabla f(x^k)-\nabla f(x^{k-1})}^2\\
    &\qquad + \frac{1}{\gamma_k L}\n{x^{k+1}-x^k}^2 \\
& = \n{x^k -x^{k-1}} \n{x^k - x^{k+1}} + \frac{\gamma_k\th_k}{L}\n{\nabla f(x^k)-\nabla f(x^{k-1})}^2 +
\frac{1}{\gamma_k L}\n{x^{k+1}-x^k}^2 \\ & \leq
  \frac 1 2 \n{x^{k}-x^{k-1}}^2 + \frac{1}{2}\n{x^{k+1}-x^k}^2 + \frac{\gamma_k\th_k}{L}\n{\nabla f(x^k)-\nabla f(x^{k-1})}^2 +
\frac{1}{\gamma_k L}\n{x^{k+1}-x^k}^2.
\end{aligned}
\end{equation}
Combining everything together,  we obtain the statement of the  theorem.
\end{proof}

\section{Stochastic Analysis}\label{ap:stoch}
\subsection{Different samples}
Consider the following version of \algname{SGD}, in which we have two samples at each iteration, $\xi^k$ and $\zeta^k$ to compute
\begin{align*}
        \gamma_k &= \min\left\{\sqrt{1+\th_k} \gamma_{k-1}, \frac{\alpha\|x^k - x^{k-1}\|}{\|\nabla f(x^k; \zeta^k) - \nabla f(x^{k-1};\zeta^k)\|} \right\}, \\
        x^{k+1} &= x^k - \gamma_k \nabla f(x^k; \xi^k).
\end{align*}
As before, we assume that $\theta_0=+\infty$, so $\gamma_1 = \frac{\alpha\|x^1 - x^{0}\|}{\|\nabla f(x^1; \zeta^1) - \nabla f(x^{0}; \zeta^1)\|}$.
\begin{lemma}
        Let $f(\cdot; \xi)$ be $L$-smooth $\mu$-strongly convex almost surely. It holds for $\gamma_k$ produced by the rule above
        \begin{align}
                \frac{\alpha}{L}
                \le \gamma_k
                \le \frac{\alpha}{\mu} \quad \text{a.s.} \label{eq:stochastic_la}
        \end{align}
\end{lemma}
\begin{proof}
        Let us start with the upper bound. Strong convexity of
        $f(\cdot; \zeta^k)$ implies that $\|x-y\|\le\frac{1}{\mu}\|\nabla
        f(x; \zeta^k) - \nabla f(y; \zeta^k)\|$ for any $x,
        y$. Therefore, $\gamma_k\le \min\left\{\sqrt{1+\th_k}\gamma_{k-1},
            \alpha/\mu \right\}\le \alpha/\mu$ a.s.

        On the other hand, $L$-smoothness gives
        \[
			\gamma_k\ge \min\left\{\sqrt{1+\th_k} \gamma_{k-1},
            \alpha/L\right\}\ge \min\left\{ \gamma_{k-1},
            \alpha/L\right\}
        \]
        almost surely. Iterating this inequality, we
        obtain the stated lower bound.
    \end{proof}
\begin{proposition}
        Denote $\sigma^2\eqdef \E{\|\nabla f (x^\star; \xi)\|^2}$ and assume $f$ to be almost surely $L$-smooth and convex. Then it holds for any $x$
        \begin{align}
                \E{\|\nabla f (x; \xi) \|^2}
                \le 4L (f(x) - f^\star) + 2\sigma^2. \label{eq:sgd_variance}
        \end{align}
\end{proposition}
Another fact that we will use is a strong convexity bound, which states for any $x,y$
\begin{align}
        \<\nabla f(x), x - y>\ge \frac{\mu}{2}\|x-y\|^2 + f(x) -f(y). \label{eq:grad_dist_bound}
\end{align}
\begin{theorem}
        Let $f(\cdot; \xi)$ be $L$-smooth and $\mu$-strongly convex almost surely. If we choose some $\alpha \le \frac{\mu}{2L}$, then
        \begin{align*}
                \E{\|x^k - x^\star\|^2}
                \le \exp\left(-k\frac{\mu\alpha}{L}\right)C_0 + \alpha\frac{\sigma^2}{\mu^2},
        \end{align*}
        where $C_0\eqdef 2(1+2\gamma_0^2L^2)\|x^0-x^\star\|^2 + 4\gamma_0^2\sigma^2$ and $\sigma^2\eqdef \E{\|\nabla fi(x^\star; \xi)\|^2}$.
    \end{theorem}
\begin{proof}
        Under our assumptions on $\alpha\le \frac{\mu}{2L}$, we have $\gamma_k\le \frac{\alpha}{\mu}\le  \frac{1}{2L}$. Since $\gamma_k$ is independent of $\xi^k$, we have $\E{\gamma_k \nabla f(x^k; \xi^k)} = \E{\gamma_k}\E{\nabla f(x^k)}$ and
        \begin{align*}
                &\E{\|x^{k+1} - x^\star\|^2} \\
                &= \E{\|x^k - x^\star\|^2} - 2 \E{\gamma_k\<\nabla f(x^k), x^k - x^\star>} + \E{\gamma_k^2}\E{\|\nabla f(x^k; \xi^k)\|^2} \\
                &\overset{\eqref{eq:grad_dist_bound}}{\le} \E{(1-\gamma_k\mu)\|x^k - x^\star\|^2} - 2\E{\gamma_k (f(x^k) - f^\star)} + \E{\gamma_k^2}\E{\|\nabla f(x^k; \xi^k)\|^2} \\
                &\overset{\eqref{eq:sgd_variance}}{\le} \E{(1-\gamma_k\mu)\|x^k - x^\star\|^2} - 2 \mathbb{E}\Bigl[\gamma_k\underbrace{(1 - 2\gamma_k L)}_{\ge 0}(f(x^k) - f^\star)\Bigr] + \E{\gamma_k^2}\sigma^2 \\
                &\overset{\eqref{eq:stochastic_la}}{\le} \E{1-\gamma_k\mu}\E{\|x^k - x^\star\|^2} + \alpha\frac{\E{\gamma_k}\sigma^2}{\mu}.
        \end{align*}
        Therefore, if we subtract $\alpha\frac{\sigma^2}{\mu^2}$ from both sides, we obtain
        \begin{align*}
                \E{\|x^{k+1} - x^\star\|^2 - \alpha\frac{\sigma^2}{\mu^2}}
                \le \E{1 - \gamma_k\mu}\E{\|x^k - x^\star\|^2 - \alpha\frac{\sigma^2}{\mu^2}}.
        \end{align*}
        If $\E{\|x^k-x^\star\|^2}\le \alpha\frac{\sigma^2}{\mu^2}$ for some $k$, it follows that $\E{\|x^{t}-x^\star\|^2}\le \alpha\frac{\sigma^2}{\mu^2}$ for any $t\ge k$. Otherwise, we can reuse the produced bound to obtain
        \begin{align*}
                \E{\|x^{k+1} - x^\star\|^2}
                \le \prod_{t=1}^k\E{1 - \gamma_t\mu}\|x^1 - x^\star\|^2 + \alpha\frac{\sigma^2}{\mu^2}.
        \end{align*}
        By inequality $1-x\le e^{-x}$, we have $\prod_{t=1}^k\E{1 - \gamma_t\mu}\le \exp\left(-\mu\sum_{t=0}^k\gamma_t \right)$.
        In addition, recall that in accordance with~\eqref{eq:stochastic_la}  we have $\gamma_k \ge \frac{\alpha}{L} $. Thus,
        \begin{align*}
                \E{\|x^{k+1} - x^\star\|^2}
                \le \exp\left(-k\alpha\frac{\mu}{L}\right)\E{\|x^1 - x^\star\|^2} + \alpha\frac{\sigma^2}{\mu^2}.
        \end{align*}
        It remains to mention that
        \begin{eqnarray*}
                \E{\|x^1-x^\star\|^2} 
                &\le& 2\|x^0-x^\star\|^2 + 2\gamma_0^2\E{\|\nabla f(x^0; \xi^0)\|^2}  \\
                &\overset{\eqref{eq:sgd_variance}}{\le}& 2\|x^0-x^\star\|^2 + 2\gamma_0^2\left(2L^2\|x^0-x^\star\|^2 + 2\sigma^2\right).
        \end{eqnarray*}
\end{proof}
This gives the following corollary.
\begin{corollary}
        Choose $\alpha = \eta\frac{\mu}{2 L}$ with $\eta\le 1$. Then, to achieve $\E{\|x^k - x^\star\|^2}= \cO(\varepsilon + \eta\sigma^2)$ we need only $k=\cO\left(\frac{L^2}{\eta\mu^2}\log \frac{1+\gamma_0^2}{\varepsilon} \right)$ iterations. If we choose $\eta$ proportionally to $\varepsilon$, it implies $\cO\left(\frac{1}{\varepsilon}\log \frac{1}{\varepsilon}\right)$ complexity.
    \end{corollary}
\subsection{Same sample: overparameterized models}
Assume additionally that the model is overparameterized, i.e., $\nabla f(x^\star; \xi)=0$ with probability one.
In that case, we can prove that one can use the same stochastic sample to compute the stepsize and to move the iterate. The update becomes
\begin{align*}
        \gamma_k &= \min\left\{\sqrt{1+\th_k} \gamma_{k-1}, \frac{\alpha\|x^k - x^{k-1}\|}{\|\nabla f(x^k; \xi^k) - \nabla f(x^{k-1}; \xi^k)\|} \right\}, \\
        x^{k+1} &= x^k - \gamma_k \nabla f(x^k; \xi^k).
\end{align*}
\begin{theorem}
        Let $f(\cdot; \xi)$ be $L$-smooth, $\mu$-strongly convex and satisfy $\nabla f(x^\star; \xi)=0$  with probability one. If we choose $\alpha \le \frac{\mu}{L}$, then
        \begin{align*}
                \E{\|x^K - x^\star\|^2}
                \le \exp\left(-K\alpha\frac{\mu}{L}\right)C_0,
        \end{align*}
        where $C_0\eqdef 2(1+\gamma_0^2L^2)\|x^0-x^\star\|^2$.
\end{theorem}
\begin{proof}
        Now $\gamma_k$ depends on $\xi^k$, so we do not have an unbiased update anymore. However, under the new assumption, $\nabla f(x^\star; \xi^k)=0$, so we can write
        \begin{align*}
                \<\nabla f(x^k; \xi^k), x^k-x^\star>
                \overset{\eqref{eq:grad_dist_bound}}{\ge} \frac{\mu}{2}\|x^k -x^\star\|^2 + f(x^k; \xi^k) - f(x^\star; \xi^k).
        \end{align*}
        In addition, $L$-smoothness and convexity of $f(\cdot; \xi^k)$ give
        \begin{align*}
                \|\nabla f(x^k; \xi^k)\|^2
                \le 2L (f(x^k; \xi^k) - f(x^\star; \xi^k)).
        \end{align*}
        Since our choice of $\alpha$ implies $\gamma_k\le \frac{1}{L}$, we conclude that
        \begin{align*}
                \|x^{k+1}-x^\star\|^2
                &= \|x^k - x^\star\|^2 - 2\gamma_k\<\nabla f(x^k; \xi^k), x^k- x^\star> + \gamma_k^2\|\nabla f(x^k; \xi^k)\|^2 \\
                &\le (1 - \gamma_k\mu)\|x^k - x^\star\|^2 - 2\gamma_k(1 - \gamma_k L)(f(x^k; \xi^k) - f(x^\star; \xi^k)) \\
                &\le (1 - \gamma_k\mu)\|x^k - x^\star\|^2.
        \end{align*}
        Furthermore, as $\|\nabla f(x^0; \xi^0)\|=\|\nabla f(x^0; \xi^0)- \nabla f(x^\star; \xi^0)\|\le L\|x^0-x^\star\|$, we also get a better bound on $\E{\|x^1-x^\star\|^2}$, namely
        \begin{align*}
                \E{\|x^1-x^\star\|^2}\le 2\|x^0-x^\star\| + 2\gamma_0^2\E{\|\nabla f(x^0; \xi^0)\|^2} \le 2(1+\gamma_0^2L^2)\|x^0-x^\star\|.
        \end{align*}
    \end{proof}

\section{Experiments Details}\label{ap:exp_details}
Here we provide some omitted details of the experiments with neural
networks. We took the implementation of neural networks from a
publicly available repository\footnote{\href{https://github.com/kuangliu/pytorch-cifar/blob/master/models/resnet.py}{https://github.com/kuangliu/pytorch-cifar/blob/master/models/resnet.py}}. All methods were run with standard data augmentation and no weight decay. The confidence intervals for ResNet-18 are obtained from 5 different random seeds and for DenseNet-121 from 3 seeds.

In our ResNet-18 experiments, we used the default parameters for
\algname{Adam}. \algname{SGD} was used with a stepsize divided by 10 at epochs 120 and
160 when the loss plateaus. Log grid search with a factor of 2 was used to tune the initial stepsize of \algname{SGD} and the best initial value was 0.2. Tuning was done by running \algname{SGD} 3 times and comparing the average of test accuracies over the runs at epoch 200. For the momentum version (\algname{SGDm}) we used the standard values of momentum and initial stepsize for training residual networks, 0.9 and 0.1 correspondingly. We used the same parameters for DenseNet-121 without extra tuning.

For our method we used the variant of \algname{SGD} $$x^{k+1}=x^k - \gamma_k \nabla f(x^k; \xi^k),$$ with $\gamma_k$ computed using $\xi^k$ as well (biased option). We did not test stepsizes that use values other than $\frac{1}{L_k}$ and $\frac{1}{2L_k}$, so it is possible that other options will perform better. Moreover, the coefficient before $\th_{k-1}$ might be suboptimal too.

\chapter{Appendix for Chapter~\ref{chapter:diana}}
\graphicspath{{diana/}{diana/img/}}

\section{Block $p$-Quantization}

We now introduce a block version of $p$-quantization. We found these quantization operators to have better properties in practice.

\begin{definition}[block-$p$-quantization]\label{def:block-p-quant} Let $\Delta = (\Delta(1), \Delta(2),\ldots, \Delta(m)) \in \R^d$, where $\Delta(1)\in\R^{d_1},\ldots, \Delta(m)\in\R^{d_m}$, $d_1+\ldots+d_m = d$ and $d_l > 1$ for all $l=1,\ldots,m$. We say that $\hat\Delta$ is $p$-quantization of $\Delta$ with sizes of blocks $(d_l)_{l=1}^m$ and write $\hat\Delta \sim {\rm Quant}_{p}(\Delta,(d_l)_{l=1}^m)$ if $\hat\Delta(l) \sim {\rm Quant}_{p}(\Delta)$ for all $l=1,\ldots,m$.
\end{definition}

In other words, we quantize blocks called \textit{blocks} of the initial vector. Note that in the special case when $m=1$ we get full quantization: ${\rm Quant}_{p}(\Delta,(d_l)_{l=1}^m) = {\rm Quant}_{p}(\Delta)$. Note that we do not assume independence of the quantization of blocks or independence of $\xi_{(j)}$. Lemma~\ref{lem:moments1} in the appendix states that $\hat{\Delta}$ is an unbiased estimator of  $\Delta$, and gives a formula for its variance.

Next we show that the block $p$-quantization operator $\hat{\Delta}$ introduced in Definition \ref{def:block-p-quant} is an unbiased estimator of  $\Delta$, and give a formula for its variance.

\begin{lemma}\label{lem:moments1}
	Let $\Delta \in \R^d$ and $\hat{\Delta} \sim {\rm Quant}_p(\Delta)$. Then for $l=1,\ldots,m$ 
	\begin{equation}
		\mathbb{E} \left[ \hat \Delta(l) \right]   = \Delta(l), \qquad
        \mathbb{E} \left[\|\hat \Delta(l) - \Delta(l)\|^2\right]  = \Psi_l(\Delta), \label{eq:tilde_v_moments1}
        \end{equation}	
	\begin{equation}
		\mathbb{E}  \left[\hat \Delta \right]  = \Delta, \qquad
        \mathbb{E} \left[\|\hat \Delta - \Delta\|^2 \right] = \Psi(\Delta), \label{eq:hat_v_moments1} 
            \end{equation}
        where
      $x = (x(1),x(2),\ldots, x(m))$, 
       $ \Psi_l(x) \eqdef \|x(l)\|_1 \|x(l)\|_p - \|x(l)\|^2 \geq 0,$
		and $\Psi(x) \eqdef \sum\limits_{l=1}^m\Psi_l(x) \ge 0.$
 Thus,  $\hat{\Delta}$ is an unbiased estimator of $\Delta$. Moreover, the variance of $\hat{\Delta}$ is a decreasing function of $p$, and is minimized for $p=\infty$.
\end{lemma}

\begin{proof}
Note that the first part of \eqref{eq:hat_v_moments1} follows from the first part of \eqref{eq:tilde_v_moments1} and the second part of \eqref{eq:hat_v_moments1} follows from the second part of \eqref{eq:tilde_v_moments1} and \[\|\hat{\Delta}-\Delta\|^2 = \sum\limits_{l=1}^m \|\hat\Delta(l)-\Delta(l)\|^2.\] Therefore, it is sufficient to prove \eqref{eq:tilde_v_moments1}. If $\Delta(l)=0$, the statements follow trivially. Assume $\Delta(l)\neq 0$.
In view of \eqref{eq:quant-j}, we have 
\[
	\ec{\hat{\Delta}_{(j)}(l)}
	 = \|\Delta(l)\|_p \sign(\Delta_{(j)}(l)) \ec{ \xi_{(j)}} =  \|\Delta(l)\|_p \sign(\Delta_{(j)}(l)) |\Delta_{(j)}(l)|/ \|\Delta(l)\|_p  = \Delta_{(j)}(l),
	 \] 
	 which establishes the first claim. We can write
\begin{eqnarray*}
	\mathbb{E} \left[\|\hat \Delta(l) - \Delta(l)\|^2 \right]
	&=& \ec{ \sum_j (\hat{\Delta}_{(j)}(l) - \Delta_{(j)}(l))^2 }\\
&=& \ec{ \sum_j (\hat{\Delta}_{(j)}(l) - \ec{ \hat{\Delta}_{(j)}(l))^2}} \\
& \overset{\eqref{eq:quant-j}}{=} &  \|\Delta(l)\|_p^2   \sum_j \sign^2(\Delta_{(j)}(l)) \ec{ (\xi_{(j)} - \ec{ \xi_{(j)}})^2} \\
&=& \|\Delta(l)\|_p^2 \sum_j \sign^2(\Delta_{(j)}(l)) \frac{|\Delta_{(j)}(l)|}{\|\Delta(l)\|_p} (1 - \frac{|\Delta_{(j)}(l)|}{\|\Delta(l)\|_p})\\
& =&  \sum_j |\Delta_{(j)}(l)| (\|\Delta(l)\|_p - |\Delta_{(j)}(l)|) \\
&= & \|\Delta(l)\|_1 \|\Delta(l)\|_p - \|\Delta(l)\|.
\end{eqnarray*}
\end{proof}

\section{Proof of Theorem~\ref{th:quantization_quality1}}

Let $1_{[\cdot]}$ denote the indicator random variable of an event.  In view of \eqref{eq:quant-j}, $\hat \Delta_{(j)} = \|\Delta\|_p \sign(\Delta_{(j)}) \xi_{(j)}$, where $\xi_{(j)}\sim {\rm Be}(|\Delta_{(j)}|/\|\Delta\|_p)$. Therefore,
\[
\|\hat \Delta\|_0 = \sum_{j=1}^d 1_{[\hat \Delta_{(j)} \neq 0]} = \sum_{j \;:\; \Delta_{(j)} \neq 0}^d 1_{[  \xi_{(j)}=1]} ,
\]
which implies that
\[
	\mathbb{E}  \left[\|\hat \Delta\|_0  \right]
= \mathbb{E} \left[ \sum_{j \;:\; \Delta_{(j)} \neq 0}^d 1_{[  \xi_{(j)}=1]}\right]  = \sum_{j \;:\; \Delta_{(j)} \neq 0}^d  \mathbb{E}   \left[1_{[  \xi_{(j)}=1]}\right]
 = \sum_{j \;:\; \Delta_{(j)} \neq 0}^d  \frac{|\Delta_{(j)}|}{\|\Delta\|_p} = \frac{\|\Delta\|_1}{\|\Delta\|_p}.
\]
To establish the first clam, it remains to recall that for all $x\in \R^d$ and $1 \leq q \leq p \leq +\infty$, one has the bound
\[ \|x\|_p \leq \|x\|_q \leq \|x\|_0^{1/q-1/p} \|x\|_p, \]
and apply it with $q  =1$.

The proof of the second claim follows the same pattern, but uses the concavity of $t \mapsto \sqrt{t}$ and Jensen's inequality in one step.

\section{Proof of Lemma~\ref{lema:alpha_p}}

$\alpha_p(d)$ is  increasing as a function of $p$ because $\|\cdot\|_p$ is decreasing as a function of $p$. Moreover, $\alpha_p(d)$ is decreasing as a function of $d$ since if we have $d < b$ then
\[
	\alpha_p(b) = \inf\limits_{x\neq 0, x\in\R^b}\frac{\|x\|^2}{\|x\|_1\|x\|_p} \leqslant \inf\limits_{x\neq 0, x\in\R_d^b}\frac{\|x\|^2}{\|x\|_1\|x\|_p} = \inf\limits_{x\neq 0, x\in\R^d}\frac{\|x\|^2}{\|x\|_1\|x\|_p},
\]
where $\R_d^b \eqdef \{x\in\R^b: x_{(d+1)} = \ldots = x_{(b)} = 0\}$.
It is known that $\frac{\|x\|}{\|x\|_1}\geq \frac{1}{\sqrt{d}}$, and that this bound is tight. Therefore, \[\alpha_1(d) = \inf_{x\neq 0,x\in\R^d} \frac{\|x\|^2}{\|x\|_1^2} = \frac{1}{d}\] 
and
\[\alpha_2(d) =  \inf_{x\neq 0,x\in\R^d} \frac{\|x\|}{\|x\|_1} = \frac{1}{\sqrt{d}}.\]

Let us now establish that $\alpha_\infty(d) = \frac{2}{1+\sqrt{d}}$. Note that
\[
\frac{\|x\|^2}{\|x\|_1\|x\|_\infty} = \frac{\left\|\frac{x}{\|x\|_\infty}\right\|^2}{\left\|\frac{x}{\|x\|_\infty}\right\|_1\left\|\frac{x}{\|x\|_\infty}\right\|_\infty} = \frac{\left\|\frac{x}{\|x\|_\infty}\right\|^2}{\left\|\frac{x}{\|x\|_\infty}\right\|_1}.
\]
Therefore, w.l.o.g.\ one can assume that $\|x\|_\infty = 1$. Moreover, signs of coordinates of vector $x$ do not influence aforementioned quantity either, so one can consider only $x\in\R_+^d$. In addition, since $\|x\|_\infty = 1$, one can assume that $x_{(1)} = 1$. Thus, our goal now is to show that the minimal value of the function
\[
	f(x) = \frac{1 + x_{(2)}^2 + \ldots + x_{(d)}^2}{1 + x_{(2)} + \ldots + x_{(d)}}
\]
on the set $M = \{x\in\R^d \mid x_{(1)} = 1, 0\le x_{(j)} \le 1, j=2,\ldots,d\}$ is equal to $\frac{2}{1+\sqrt{d}}$. 
By Cauchy-Schwartz inequality: $x_{(2)}^2 + \ldots + x_{(d)}^2 \ge \frac{(x_{(2)} + \ldots + x_{(d)})^2}{d-1}$ and it becomes equality if and only if all $x_{(j)},j=2,\ldots,d$ are equal. It means that if we fix $x_{(j)} = a$ for $j=2,\ldots,d$ and some $0\le a \le 1$ than the minimal value of the function
\[
	g(a) = \frac{1 + \frac{((d-1)a)^2}{d-1}}{1 + (d-1)a} = \frac{1 + (d-1)a^2}{1 + (d-1)a}
\]
on $[0,1]$ coincides with minimal value of $f$ on $M$. The derivative
\[
	g'(a) = \frac{2(d-1)a}{1+(d-1)a} - \frac{(d-1)(1+(d-1)a^2)}{(1+(d-1)a)^2}
\]
has the same $\sign$ on $[0,1]$ as the difference $a - \frac{1 + (d-1)a^2}{2(1+ (d-1)a)}$, which implies that $g$ attains its minimal value on $[0,1]$ at such $a$ that $a = \frac{1 + (d-1)a^2}{2(1+ (d-1)a)}$. It remains to find $a\in[0,1]$ which satisfies
\[
	a = \frac{1 + (d-1)a^2}{2(1 + (d-1)a)},\quad a\in[0,1] \Longleftrightarrow (d-1)a^2 + 2a - 1 = 0,\quad a\in[0,1].
\]
This quadratic equation has unique positive solution $a^\star = \frac{-1 + \sqrt{d}}{d-1} = \frac{1}{1 + \sqrt{d}} < 1$. Direct calculations show that $g(a^\star) = \frac{2}{1+\sqrt{d}}$. It implies that $\alpha(d) = \frac{2}{1+\sqrt{d}}$.

\section{Strongly Convex Case: Optimal Number of Nodes} \label{sec:opt_no_nodes}

In practice one has access to a finite dataset, consisting of $N$ data points, where $N$ is very large, and wishes to solve an empirical risk minimization (``finite-sum'') of the form
\begin{equation}\label{eq:N} 
	\min_{x\in \R^d} f(x) = \frac{1}{N} \sum \limits_{i=1}^N \phi_i(x) + \psi(x),
\end{equation}
where each $\phi_i$ is $L$--smooth and $\mu$-strongly convex. If $M \leq N$ compute nodes of a distributed system  are available, one may partition the $N$ functions into $M$ groups, $G_1,\dots,G_M$, each of size $|G_i| = N/n$, and define
$f_m(x) = \frac{M}{N} \sum_{i \in G_m} \phi_i(x). $
Note that it still holds that 
\[
	f(x) = \frac{1}{M} \sumiM f_i(x) + \psi(x).
\]
Moreover, each $f_m$ is also $L$--smooth and $\mu$--strongly convex.

This way, we  have  fit the original (and large) problem \eqref{eq:N} into our framework. One may now ask the question: How many many nodes $M$ should we use (other things equal)? If what we care about is iteration complexity, then insights can be gained by investigating~\eqref{eq:bu987gd9}. For instance, if $p=2$, then the complexity is $W(M) \eqdef \max\left\{\frac{2\sqrt{d}}{\sqrt{m}},(\kappa+1)\left(\frac{1}{2}-\frac{1}{M}+\frac{\sqrt{d}}{M\sqrt{m}}\right)\right\}.$ The optimal choice is to choose $M$ so that the term $-\frac{1}{M} + \frac{\sqrt{d}}{M\sqrt{m}}$ becomes (roughly) equal to $\frac{1}{2}$:
$-\frac{1}{M} + \frac{\sqrt{d}}{M\sqrt{m}} = \frac{1}{2}.$ This gives the formula for the optimal number of nodes
$M^\star = M(d) \eqdef 2\left(\sqrt{\frac{d}{m}} - 1\right),$
and the resulting iteration complexity is $W(M^\star) = \max\left\{\frac{2\sqrt{d}}{\sqrt{m}}, \kappa + 1\right\}$.  Note that $M(d)$ is increasing in $d$. Hence, it makes sense to use more nodes for larger models (big $d$).

\section{Quantization Lemmas}

Consider iteration $k$ of the \algname{DIANA} method (Algorithm~\ref{alg:distributed1}).  Let $\mathbb{E}_{Q^k}$ be the expectation with respect to the randomness inherent in the quantization steps $\hat \Delta_i^k \sim {\rm Quant}_p(\Delta_i^k,(d_l)_{l=1}^m)$ for $i=1,2,\dots,n$ (i.e., we condition on everything else).

\begin{lemma}\label{lem:3in1} For all iterations $k\geq 0$  of \algname{DIANA} and $i=1,2,\dots,n$ we have the identities
\begin{equation}\label{eq:hat_gi_moments1} 
	\mathbb{E}_{Q^k} \left[ \hat g_i^k \right]  = g_i^k, 
	\qquad \mathbb{E}_{Q^k} \left[ \| \hat g_i^k - g_i^k\|^2 \right] 
	=\Psi(\Delta_i^k) 
\end{equation}
and
\begin{equation} 
	\label{eq:distr_hat_g_moments1} \mathbb{E}_{Q^k} \left[ \hat g^k  \right]
	 = g^k 
	 \eqdef \frac{1}{M}\sum_{i=1}^M g_i^k, 
	 \qquad \mathbb{E}_{Q^k}  \left[ \| \hat g^k - g^k\|^2   \right]
	 = \frac{1}{M^2}\sum_{i=1}^M\Psi(\Delta_i^k) .
\end{equation}
Furthermore, letting $h^\star = \nabla f(x^\star)$, and invoking Assumption~\ref{as:noise}, we have
\begin{align}
	\mathbb{E} \left[\hat g^k \right]
	&= \nabla f(x^k), \\
	\mathbb{E} \left[ \|\hat g^k -h^\star\|^2\right]
	 &\leq \mathbb{E}\left[\|\nabla f(x^k) - h^\star\|^2\right] + \frac{1}{M^2}\sum_{i=1}^M\mathbb{E} \left[\Psi(\Delta_i^k)\right] + \frac{\sigma^2}{M} . \label{eq:full_variance_of_mean_g1}
\end{align}
\end{lemma}
\begin{proof}

\begin{itemize}
\item[(i)] Since $\hat g_i^k = h_i^k + \hat \Delta_i^k$ and $\Delta_i^k = g_i^k - h_i^k$, we can apply Lemma~\ref{lem:moments1} and obtain
\[
	\mathbb{E}_{Q^k}  \left[\hat g_i^k  \right]
	= h_i^k + \mathbb{E}_{Q^k}  \left[\hat \Delta_i^k \right]
	\overset{\eqref{eq:hat_v_moments1}}{=}  h_i^k + \Delta_i^k 
	= g_i^k.
\]
Since
$\hat g_i^k - g_i^k = \hat \Delta_i^k - \Delta_i^k $, applying the second part of Lemma 1 gives the second identity in \eqref{eq:hat_gi_moments1}.

\item [(ii)] The first part of \eqref{eq:distr_hat_g_moments1}  follows directly from the first part of \eqref{eq:hat_gi_moments1}: 
\[
	\mathbb{E}_{Q^k}  \left[\hat g^k \right]
	= \mathbb{E}_{Q^k} \left[ \frac{1}{M} \sum_{i=1}^M \hat g_i^k \right] 
	=  \frac{1}{M} \sum_{i=1}^M \mathbb{E}_{Q^k}  \left[\hat g_i^k \right]
	\overset{\eqref{eq:hat_gi_moments1}}{=} \frac{1}{M} \sum_{i=1}^M  g_i^k 
	\overset{\eqref{eq:distr_hat_g_moments1} }{=}  g^k.
\]
The second part in \eqref{eq:distr_hat_g_moments1} follows from the second part of \eqref{eq:hat_gi_moments1} and independence of $\hat g_1^k,\dotsc, \hat g_M^k$.

\item [(iii)] The first part of \eqref{eq:full_variance_of_mean_g1} follows directly from the first part of \eqref{eq:distr_hat_g_moments1} and the assumption that  $g_i^k$ is  and unbiased estimate of $\nabla f_i(x^k)$.
 It remains to establish the second part of  \eqref{eq:full_variance_of_mean_g1}. First, we shall decompose
 \begin{eqnarray*}
    	\mathbb{E}_{Q^k} \left[\|\hat g^k - h^\star\|^2 \right]
    	&\overset{\eqref{eq:second_moment_decomposition}}{=}& \mathbb{E}_{Q^k} \left[\|\hat g^k - \mathbb{E}_{Q^k} \left[\hat g^k\right]\|^2 \right]+ \| \mathbb{E}_{Q^k}  \left[\hat g^k\right] - h^\star\|^2\\
    	&\overset{\eqref{eq:distr_hat_g_moments1} }{=}&\mathbb{E}_{Q^k} \left[\|\hat g^k - g^k\|^2\right] + \|  g^k - h^\star\|^2\\
    &\overset{\eqref{eq:distr_hat_g_moments1} }{=}& \frac{1}{M^2}\sum_{i=1}^M \Psi(\Delta_i^k) +  \|  g^k - h^\star\|^2.
\end{eqnarray*}

Further, applying variance decomposition \eqref{eq:second_moment_decomposition}, we get
  \begin{eqnarray*}
    	\mathbb{E}\left[  \|g^k - h^\star\|^2 \;|\; x^k \right]
    	&\overset{\eqref{eq:second_moment_decomposition}}{=}&\mathbb{E} \left[ \|g^k - \mathbb{E} [g^k\;|\; x^k]\|^2 \; | \; x^k \right] +  \|\mathbb{E} [g^k \;|\; x^k ] - h^\star\|^2 \\
    	&\overset{\eqref{eq:hat_g_expectation}}{=}&\mathbb{E} \left[ \|g^k - \nabla f(x^k)\|^2 \; | \; x^k \right] +  \| \nabla f(x^k)- h^\star\|^2 \\
    	& \overset{\eqref{eq:bgud7t9gf}}{\leq} & \frac{\sigma^2}{M} + \|\nabla f(x^k) - h^\star\|^2 .
    \end{eqnarray*}

Combining the two results, we get
\begin{eqnarray*}
    \mathbb{E}\left[	\mathbb{E}_{Q^k}[\|\hat g^k - h^\star\|^2] \;|\; x^k \right]]
    	\leq  \frac{1}{M^2} \sum_{i=1}^M \mathbb{E}\left[ \Psi(\Delta_i^k) \;|\; x^k \right] + \frac{\sigma^2}{M} + \|\nabla f(x^k) - h^\star\|^2.
\end{eqnarray*}
After applying full expectation, and using tower property, we get the result.
\end{itemize}

\end{proof}

\begin{lemma}\label{lem:distr_h_recurrence1}
    Let $x^\star$ be a solution of \eqref{eq:main} and let $h_i^\star = \nabla f_i(x^\star)$ for $i=1,2,\dots,d$. For every $i$, we can estimate the first two moments of $h_i^{k+1}$ as
    \begin{align}
        \mathbb{E}_{Q^k} \left[h_i^{k+1} \right]
        &= (1 - \alpha)h^k_i + \alpha g_i^k, \nonumber\\
        \mathbb{E}_{Q^k} \left[\|h_i^{k+1} - h_i^\star\|^2\right]
        & = (1 - \alpha)\|h_i^k - h_i^\star\|^2 + \alpha \|g_i^k - h_i^\star\|^2 \notag\\
        &\qquad - \alpha \left( \|\Delta_i^k\|^2 - \alpha \sum\limits_{l=1}^m\|\Delta_i^k(l)\|_1 \|\Delta_i^k(l)\|_p  \right) . \label{eq:distr_h^(k+1)le}
    \end{align}
\end{lemma}

\begin{proof}
    Since 
    \begin{equation}\label{eq:b87f9h8hf9}
    		h_i^{k+1} 
    		= h_i^k + \alpha \hat \Delta_i^k
    	\end{equation} 
    	and $\Delta_i^k = g_i^k - h_i^k$, in view of Lemma~\ref{lem:moments1} we have
\begin{equation}\label{eq:09h80hdf}
        \mathbb{E}_{Q^k} \left[h_i^{k+1}\right]  \overset{\eqref{eq:b87f9h8hf9}}{=}   h_i^k + \alpha \mathbb{E}_{Q^k} \left[\hat \Delta_i^k \right]
         \overset{\eqref{eq:hat_v_moments1}}{=}  h_i^k + \alpha \Delta_i^k   
           =  (1 - \alpha)h_i^k + \alpha g_i^k,
\end{equation}
which establishes the first claim. Further, using $\|\Delta_i^k\|^2 = \sum\limits_{l=1}^m\|\Delta_i^k(l)\|^2$ we obtain
    \begin{eqnarray*}
        \mathbb{E}_{Q^k} \left[\|h_i^{k+1} - h_i^\star\|^2  \right]
        &\overset{\eqref{eq:second_moment_decomposition}}{=}& \|\mathbb{E}_{Q^k} [h_i^{k+1} - h_i^\star]\|^2 + \mathbb{E}_{Q^k} \left[\| h_i^{k + 1} - \mathbb{E}_{Q^k} \left[ h_i^{k+1} \right]\|^2\right] \nonumber \\
        &\overset{\eqref{eq:09h80hdf}+
        \eqref{eq:b87f9h8hf9}      }{=}& \|(1 - \alpha) h_i^k + \alpha g_i^k - h_i^\star\|^2 + \alpha^2 \mathbb{E}_{Q^k}\left[\| \hat \Delta_i^k - \mathbb{E}_{Q^k} \left[ \hat \Delta_i^k\right] \|^2\right] \nonumber \\
        &\overset{\eqref{eq:hat_v_moments1}}{=}& \|(1 - \alpha) (h_i^k - h_i^\star) + \alpha (g_i^k - h_i^\star)\|^2\\
        &&\qquad + \alpha^2 \sum\limits_{l=1}^m(\|\Delta_i^k(l)\|_1 \|\Delta_i^k(l)\|_p - \|\Delta_i^k(l)\|^2) \nonumber \\
        &\overset{\eqref{eq:sqaured_norm_of_lin_combination} }{=}& (1 - \alpha)\| h_i^k - h_i^\star\|^2 + \alpha \| g_i^k - h_i^\star\|^2 - \alpha(1 - \alpha)\| \Delta_i^k\|^2\nonumber\\
        &&\quad + \alpha^2 \sum\limits_{l=1}^m(\|\Delta_i^k(l)\|_1 \|\Delta_i^k(l)\|_p) - \alpha^2\|\Delta_i^k\|^2\nonumber \\
&=&        (1 - \alpha)\| h_i^k - h_i^\star\|^2 + \alpha \| g_i^k - h_i^\star\|^2\\
	&&\quad + \alpha^2\sum\limits_{l=1}^m( \|\Delta_i^k(l)\|_1 \|\Delta_i^k(l)\|_p) - \alpha \|\Delta_i^k\|^2.\nonumber \\
    \end{eqnarray*}

\end{proof}

\begin{lemma}
	We have
	\begin{align}
		\mathbb{E}\left[\|\hat g^k - h^\star\|^2 \mid x^k\right] \le \|\nabla f(x^k) - h^\star\|^2 + \left(\frac{1}{\alpha_p} - 1\right)\frac{1}{M^2}\sumiM \|\nabla f_i(x^k) - h_i^k\|^2 + \frac{\sigma^2}{\alpha_p M}. \label{eq:full_second_moment_of_hat_g}
	\end{align}
\end{lemma}
\begin{proof}
	Since $\alpha_p = \alpha_p(\max\limits_{l=1,\ldots,m}d_l)$ and $\alpha_p(d_l) = \inf\limits_{x\neq 0,x\in\R^{d_l}}\frac{\|x\|^2}{\|x\|_1\|x\|_p} $, we have for a particular choice of $x=\Delta_i^k(l)$ that $\alpha_p\le\alpha_p(d_l) \le \frac{\|\Delta_i^k(l)\|^2}{\|\Delta_i^k(l)\|_1\|\Delta_i^k(l)\|_p}$. Therefore,
	\begin{align*}
		\Psi(\Delta_i^k) = \sum\limits_{l=1}^m\Psi_l(\Delta_i^k) 
		&= \sum\limits_{l=1}^m(\|\Delta_i^k(l)\|_1\|\Delta_i^k(l)\|_\infty - \|\Delta_i^k(l)\|) \\
		&\le \sum\limits_{l=1}^m \left(\frac{1}{\alpha_p}-1\right)\|\Delta_i^k(l)\|^2 = \left(\frac{1}{\alpha_p}-1\right)\|\Delta_i^k\|^2.
	\end{align*}
	 This can be applied to~\eqref{eq:full_variance_of_mean_g1} in order to obtain
	\begin{align*}
		\mathbb{E}\left[\|\hat g^k - h^\star \|^2 \mid x^k\right] 
		&\le  \|\nabla f(x^k) - h^\star\|^2 + \frac{1}{M^2}\sumiM \mathbb{E}\left[\Psi(\Delta_i^k)\mid x^k \right]  + \frac{\sigma^2}{M}\\
		&\le  \|\nabla f(x^k) - h^\star\|^2 + \frac{1}{M^2}\sumiM \left(\frac{1}{\alpha_p} - 1\right)\mathbb{E}\left[\|\Delta_i^k\|^2\mid x^k \right]  + \frac{\sigma^2}{M}.
	\end{align*} 	
 	Note that for every $i$ we have $\mathbb{E}\left[\Delta_i^k\mid x^k\right] = \mathbb{E}\left[g_i^k - h_i^k\mid x\right] = \nabla f_i(x^k) - h_i^k$, so
 	\begin{eqnarray*}
 		\mathbb{E} \left[\|\Delta_i^k\|^2 \mid x^k\right] 
 		&\overset{\eqref{eq:second_moment_decomposition}}{=} & \|\nabla f_i(x^k) - h_i^k\|^2 +  \mathbb{E}\left[\|g_i^k - \nabla f_i(x^k)\|^2 \mid x^k\right]  \\
 		&\le & \|\nabla f_i(x^k) - h_i^k\|^2 +  \sigma_i^2 .
 	\end{eqnarray*}
 	Summing the produced bounds, we get the claim.
\end{proof}

\section{Proof of Theorem~\ref{thm:DIANA-strongly_convex}}

\begin{proof}
Note that $x^\star$ is a solution of \eqref{eq:main} if and only if $x^\star = \prox_{\gamma R}(x^\star-\gamma h^\star)$ (this holds for any $\gamma>0$).  Using this identity together with the nonexpansiveness of the proximal operator, we shall  bound the first term of the Lyapunov function:
    \begin{eqnarray*}
        \mathbb{E}_{Q^k} \left[\|x^{k+1} - x^\star\|^2 \right]
        &=& \mathbb{E}_{Q^k} \left[\|\proxR(x^k - \gamma \hat g^k) - \proxR(x^\star - \gamma h^\star)\|^2\right] \\
        &\overset{\eqref{eq:nonexpansive}}{\le}& \mathbb{E}_{Q^k}\left[\|x^k - \gamma \hat g^k - (x^\star - \gamma h^\star)\|^2\right] \\
        &=& \|x^k - x^\star\|^2 - 2\gamma \mathbb{E}_{Q^k} \left[\<\hat g^k - h^\star, x^k - x^\star>\right] + \gamma^2\mathbb{E}_{Q^k}\left[\|\hat g^k - h^\star\|^2\right] \\
        &\overset{\eqref{eq:distr_hat_g_moments1}}{=}& \|x^k - x^\star\|^2 - 2\gamma \<g^k - h^\star, x^k - x^\star> + \gamma^2\mathbb{E}_{Q^k} \left[\|\hat g^k - h^\star\|^2\right].
	\end{eqnarray*}
Next, taking conditional expectation on both sides of the above inequality, and using \eqref{eq:hat_g_expectation}, we get

\begin{align*}
       \mathbb{E}\left[ \mathbb{E}_{Q^k}\left[\|x^{k+1} - x^\star\|^2\right] \;|\; x^k \right]
        &\leq  \|x^k - x^\star\|^2 - 2\gamma \< \nabla f(x^k) - h^\star, x^k - x^\star>\\
        &\qquad + \gamma^2\mathbb{E}\left[ \mathbb{E}_{Q^k} [\|\hat g^k - h^\star\|^2] \;|\; x^k\right].
	\end{align*}
	Taking full expectation on both sides of the above inequality, and  applying the tower property and Lemma~\ref{lem:3in1} leads to
	\begin{eqnarray}
       \ec{ \|x^{k+1} - x^\star\|^2 }
        &\leq & \ec{ \|x^k - x^\star\|^2} - 2\gamma \ec{ \< \nabla f(x^k) - h^\star, x^k - x^\star>} + \gamma^2\ec{  \|\hat g^k - h^\star\|^2} \notag \\
        &\overset{\eqref{eq:full_variance_of_mean_g1}}{\leq} & \ec{ \|x^k - x^\star\|^2} - 2\gamma \ec{ \< \nabla f(x^k) - h^\star, x^k - x^\star>} \notag \\
        && \ + \gamma^2 \ec{\|\nabla f(x^k) - h^\star\|^2} + \frac{\gamma^2}{M^2}\sum_{i=1}^M \ec{ \Psi(\Delta_i^k)} + \frac{\gamma^2 \sigma^2}{M} \notag \\
        &\leq & 
        \ec{ \|x^k - x^\star\|^2} - 2\gamma \ec{\< \nabla f(x^k) - h^\star, x^k - x^\star>} \notag \\
        && + \frac{\gamma^2}{M} \sumiM \ec{ \|\nabla f_i(x^k) - h_i^\star\|^2} + \frac{\gamma^2}{M^2}\sum_{i=1}^M \ec{ \Psi(\Delta_i^k)} + \frac{\gamma^2 \sigma^2}{M},
        \label{eq:buf89gh38bf98}
	\end{eqnarray}
	where the last inequality follows from the identities $\nabla f(x^k) = \frac{1}{M}\sumiM f_i(x^k)$, $h^\star = \frac{1}{M}\sumiM h_i^\star$ and an application of Jensen's inequality.

Averaging over the identities \eqref{eq:distr_h^(k+1)le} for $i=1,2,\dots,n$ in Lemma~\ref{lem:distr_h_recurrence1}, we get
\begin{align*}
       \frac{1}{M}\sumiM \mathbb{E}_{Q^k} \left[\|h_i^{k+1} - h_i^\star\|^2 \right]
         &= \frac{1 - \alpha}{M}\sumiM\|h_i^k - h_i^\star\|^2 +  \frac{\alpha}{M}\sumiM\|g_i^k - h_i^\star\|^2\\
         &\qquad -  \frac{\alpha}{M} \sumiM \left(\|\Delta_i^k\|^2 -\alpha \sum\limits_{l=1}^m \|\Delta_i^k(l)\|_1 \|\Delta_i^k(l)\|_p  \right) . 
\end{align*}
    Applying expectation to both sides, and using the tower property, we get
   \begin{align}
       \frac{1}{M}\sumiM \ec{ \|h_i^{k+1} - h_i^\star\|^2}
         &= \frac{1 - \alpha}{M}\sumiM \ec{ \|h_i^k - h_i^\star\|^2} +  \frac{\alpha}{M}\sumiM \ec{ \|g_i^k - h_i^\star\|^2} \notag\\
         &\quad -  \frac{\alpha}{M} \sumiM \mathbb{E} \left[\|\Delta_i^k\|^2 -\alpha \sum\limits_{l=1}^m \|\Delta_i^k(l)\|_1 \|\Delta_i^k(l)\|_p  \right] . \label{eq:nbi987fg98bf9hbf}
\end{align}

Since
\begin{equation*}
		\mathbb{E} [ \|g_i^k - h_i^\star\|^2 \;|\; x^k]
		\overset{\eqref{eq:second_moment_decomposition}}{=} \|\nabla f_i(x^k) - h_i^\star\|^2 + \mathbb{E} [\|g_i^k - \nabla f_i(x^k)\|^2\;|\; x^k ] \overset{\eqref{eq:bounded_noise}}{\le} \|\nabla f_i(x^k) - h_i^\star\|^2 + \sigma_i^2,
	\end{equation*}
the second term on the right-hand side of \eqref{eq:nbi987fg98bf9hbf} can be bounded above as
	\begin{equation}\label{eq:b897fg98fss}
		\ec{ \|g_i^k - h_i^\star\|^2 }
		\leq  \ec{ \|\nabla f_i(x^k) - h_i^\star\|^2} + \sigma_i^2.
	\end{equation}
Plugging \eqref{eq:b897fg98fss} into \eqref{eq:nbi987fg98bf9hbf} leads to the estimate
   \begin{align}
       \frac{1}{M}\sumiM \ec{ \|h_i^{k+1} - h_i^\star\|^2}
         &\leq \frac{1 - \alpha}{M}\sumiM \ec{ \|h_i^k - h_i^\star\|^2} +  \frac{\alpha}{M}\sumiM \ec{ \|\nabla f_i(x^k) - h_i^\star\|^2} \notag \\
         & + \alpha \sigma^2    -  \frac{\alpha}{M} \sumiM \mathbb{E} \left[\|\Delta_i^k\|^2 -\alpha \sum\limits_{l=1}^m \|\Delta_i^k(l)\|_1 \|\Delta_i^k(l)\|_p  \right]  . \label{eq:nbi9bh98sgs}
\end{align}
	
Adding \eqref{eq:buf89gh38bf98} with the $c\gamma^2$ multiple of \eqref{eq:nbi9bh98sgs}, we get an upper bound one the Lyapunov function:
\begin{align}
	\mathbb{E} \left[V^{k+1}\right]
	& \leq  \ec{\|x^k-x^\star\|^2}  +  \frac{(1-\alpha) c\gamma^2}{M}\sumiM \ec{\|h_i^k - h_i^\star\|^2} \notag \\ 
	& \qquad +  \frac{\gamma^2 ( 1 + \alpha c) }{M}\sum_{i=1}^M\ec{\|\nabla f_i(x^k) - h_i^\star\|^2}   - 2\gamma \ec{\<\nabla f(x^k) - h^\star, x^k - x^\star> } \notag \\
	& \qquad  + \frac{\gamma^2}{M^2} \sumiM\sum\limits_{l=1}^m \ec{ T_i^k(l)} + (Mc\alpha + 1)\frac{\gamma^2\sigma^2}{M}, \label{eq:bd7g98gdh8d}
\end{align}
where
\[
	T_i^k(l)
	\eqdef \left[   \left(\|\Delta_i^k(l)\|_1 \|\Delta_i^k(l)\|_p - \|\Delta_i^k(l)\|^2\right)  - n \alpha c  \left(\|\Delta_i^k(l)\|^2 -\alpha  \|\Delta_i^k(l)\|_1 \|\Delta_i^k(l)\|_p  \right)  \right].
\]
We now claim that due to our choice of $\alpha$ and $c$, we have $ T_i^k(l)\leq 0$ for all $\Delta_i^k(l)\in \R^{d_l}$, which means that we can bound this term away by zero. Indeed, note that $T_k^i(l) = 0$ for $\Delta_i^k(l) = 0$. If $\Delta_i^k(l) \neq 0$, then $T_k^i(l) \leq 0$ can be equivalently written as 
\[\frac{1 + Mc \alpha^2}{1 + Mc \alpha} \leq \frac{\|\Delta_i^k(l)\|^2}{\|\Delta_i^k(l)\|_1\|\Delta_i^k(l)\|_p}.\] 
However, this inequality holds since in view of the first inequality in \eqref{eq:cond1} and the definitions of $\alpha_p$ and $\alpha_p(d_l)$, we have
\[\frac{1 + Mc \alpha^2}{1 + Mc \alpha} \overset{\eqref{eq:cond1} }{\leq} \alpha_p \le \alpha_p(d_l) \overset{\eqref{eq:alpha_p}}{=} \inf_{x\neq 0,x\in\R^{d_l}} \frac{\|x\|^2}{\|x\|_1\|x\|_p} \leq  \frac{\|\Delta_i^k(l)\|^2}{\|\Delta_i^k(l)\|_1\|\Delta_i^k(l)\|_p}.\] 

Therefore, from \eqref{eq:bd7g98gdh8d} we get

\begin{align}
	\ec{ V^{k+1}} 
	& \leq  \ec{ \|x^k-x^\star\|^2}  +  \frac{(1-\alpha) c\gamma^2}{M}\sumiM \ec{\|h_i^k - h_i^\star\|^2} \notag \\ 
	& \qquad +  \frac{\gamma^2 ( 1 + \alpha c) }{M}\sum_{i=1}^M\ec{\|\nabla f_i(x^k) - h_i^\star\|^2}   \notag \\
	& \qquad  - 2\gamma \ec{\<\nabla f(x^k) - h^\star, x^k - x^\star>} + (Mc\alpha + 1)\frac{\gamma^2\sigma^2}{M}. \label{eq:bd7g98gdh8d11}
\end{align}

The next trick is to split $\nabla f(x^k)$ into the average of $\nabla f_i(x^k)$ in order to apply strong convexity of each term:
	\begin{eqnarray}
		&&\ec{\<\nabla f(x^k) - h^\star, x^k - x^\star> } \\
		&=& \frac{1}{M}\sumiM \ec{\<\nabla f_i(x^k) - h_i^\star, x^k - x^\star>} \notag\\
		&\overset{\eqref{eq:scal_prod_tight_str_cvx}}{\ge} &  \frac{1}{M} \sumiM \mathbb{E} \left[ \frac{\mu L}{\mu + L} \|x^k - x^\star\|^2 +   \frac{1}{\mu + L} \|\nabla f_i(x^k) - h_i^\star\|^2 \right]\notag\\
		&=& \frac{\mu L}{\mu + L}\ec{\|x^k - x^\star\|^2} + \frac{1}{\mu + L}\frac{1}{M}\sumiM \ec{\|\nabla f_i(x^k) - h_i^\star\|^2} \label{eq:inner_product_splitting}.
	\end{eqnarray}

Plugging these estimates into \eqref{eq:bd7g98gdh8d11}, we obtain 
\begin{align}
	\ec{ V^{k+1}}
	& \leq  \left(1 - \frac{2\gamma \mu L}{\mu+L}\right) \ec{ \|x^k-x^\star\|^2}  +  \frac{(1-\alpha) c\gamma^2}{M}\sumiM \ec{\|h_i^k - h_i^\star\|^2} \notag \\ 
	& \quad +  \left(\gamma^2 ( 1 + \alpha c) -\frac{2\gamma}{\mu+L} \right)\frac{1}{M}\sum_{i=1}^M\ec{\|\nabla f_i(x^k) - h_i^\star\|^2} + (Mc\alpha + 1)\frac{\gamma^2\sigma^2}{M}. \label{eq:bd7nb98sgsdd}
\end{align}

Notice that in view of the second inequality in \eqref{eq:cond2}, we have $ \gamma^2 ( 1 + \alpha c) -\frac{2\gamma}{\mu+L} \leq 0$. Moreover, since $f_i$ is $\mu$--strongly convex, we have
$\mu \|x^k-x^\star\|^2 \leq \langle \nabla f_i(x^k) - h_i^\star, x^k -x^\star \rangle$. Applying the Cauchy-Schwarz inequality to further  bound the right-hand side, we get the inequality $\mu \|x^k-x^\star\| \leq \|\nabla f_i(x^k) - h_i^\star\|$. Using these observations, we can get rid of the term on the second line of \eqref{eq:bd7nb98sgsdd} and absorb it with the first term, obtaining
\begin{align}
	\ec{V^{k+1}}  &\leq  \left(1 - 2\gamma \mu + \gamma^2 \mu^2 + c\alpha \gamma^2 \mu^2 \right) \ec{ \|x^k-x^\star\|^2} \notag\\
	&\quad  +  \frac{(1-\alpha) c\gamma^2}{M}\sumiM \ec{\|h_i^k - h_i^\star\|^2}  + (Mc\alpha + 1)\frac{\gamma^2\sigma^2}{M}. \label{eq:bd7n98sh90shdw}
\end{align}
It follows from the second inequality in \eqref{eq:cond2} that 
\[
	1 - 2\gamma \mu + \gamma^2 \mu^2 + c\alpha \gamma^2 \mu^2 \leq 1 - \gamma \mu.
\] 
Moreover, the first inequality in \eqref{eq:cond2} implies that $1-\alpha \leq 1-\gamma \mu$. Consequently, from \eqref{eq:bd7n98sh90shdw} we obtain the recursion
\[
	\ec{ V^{k+1}}  \leq (1-\gamma \mu) \ec{ V^k} +(Mc\alpha + 1)\frac{\gamma^2\sigma^2}{M}. 
\] 
Finally, unrolling the recurrence leads to
    \begin{align*}
    		\ec{ V^k }
    		&\le (1 - \gamma\mu)^k V^0 + \sum\limits_{l=0}^{k-1}(1-\gamma\mu)^l\gamma^2(1+Mc\alpha)\frac{\sigma^2}{M} \\
    		&\le (1 - \gamma\mu)^k V^0 + \sum\limits_{l=0}^{\infty}(1-\gamma\mu)^l\gamma^2(1+Mc\alpha)\frac{\sigma^2}{M} \\
    		&= (1 - \gamma\mu)^k V^0 + \frac{\gamma}{\mu}(1+Mc\alpha)\frac{\sigma^2}{M}.
    \end{align*}
\end{proof}

\section{Proof of Corollary~\ref{cor:DIANA-strong-convex}}

\begin{corollary} Let $\kappa = \frac{L}{\mu}$, $\alpha = \frac{\alpha_p}{2}$, $c = \frac{4(1-\alpha_p)}{M\alpha_p^2}$, and $\gamma = \min\left\{\frac{\alpha}{\mu}, \frac{2}{(L+\mu)(1+c \alpha)}\right\}$. Then the conditions \eqref{eq:cond1} and \eqref{eq:cond2} are satisfied, and the leading term in the iteration complexity bound is equal to
\begin{equation} \frac{1}{\gamma \mu} = \max\left\{\frac{2}{\alpha_p}, (\kappa+1)\left(\frac{1}{2} - \frac{1}{M} + \frac{1}{M\alpha_p}\right)\right\}.\end{equation} 
This is a decreasing function of $p$.  Hence, from iteration complexity perspective, $p=+\infty$ is the optimal choice.
\end{corollary}
\begin{proof}
	Condition \eqref{eq:cond2} is satisfied since $\gamma = \min\left\{\frac{\alpha}{\mu}, \frac{2}{(L+\mu)(1+c \alpha)}\right\}$. Now we check that \eqref{eq:cond1} is also satisfied:
	\begin{eqnarray*}
		\frac{1+Mc\alpha^2}{1+Mc\alpha} \frac{1}{\alpha_p} &=& \frac{1 + M\cdot\frac{4(1-\alpha_p)}{M\alpha_p^2}\cdot\frac{\alpha_p^2}{4}}{1 + M\cdot\frac{4(1-\alpha_p)}{M\alpha_p^2}\cdot\frac{\alpha_p}{2}}\cdot\frac{1}{\alpha_p}\\
		&=& \frac{2 - \alpha_p}{\alpha_p + 2(1-\alpha_p)}\\
		&=& 1.
	\end{eqnarray*}
	Since $\alpha = \frac{\alpha_p}{2}$ and $c = \frac{4(1-\alpha_p)}{M\alpha_p^2}$ we have
	\[
		1+\alpha c = 1 + \frac{2(1-\alpha_p)}{M\alpha_p} = 1 - \frac{2}{M} + \frac{2}{M\alpha_p}
	\]
	and, therefore,
	\[
		\frac{1}{\gamma\mu} = \max\left\{\frac{1}{\alpha}, \frac{L+\mu}{2\mu}(1+c\alpha)\right\} = \max\left\{\frac{2}{\alpha_p}, (\kappa+1)\left(\frac{1}{2} - \frac{1}{M} + \frac{1}{M\alpha_p}\right)\right\},	
	\]
	which is a decreasing function of $p$, because $\alpha_p$ increases when $p$ increases.
\end{proof}

\section{Strongly Convex Case: Decreasing Stepsize}\label{sec:DIANA-decreasing-stepsize-appendix}

\begin{lemma}
\label{lem:sgd_recurrence}
	Let a sequence $(a^k)_k$ satisfy inequality $a^{k+1}\le (1 - \gamma_k\mu) a^k + \gamma_k^2 \nu$ for any positive $\gamma_k \le \gamma_0$ with some constants $\mu > 0, \nu>0, \gamma_0 > 0$. Further, let $\theta \ge \frac{2}{\gamma_0} $ and take $C$ such that $\nu\le \frac{\mu\theta}{4}C$ and $a_0\le C$. Then, it holds
	\begin{align*}
		a^k \le \frac{ C}{\frac{\mu}{\theta} k + 1}
	\end{align*}
	if we set $\gamma_k=\frac{2}{\mu k + \theta}$.
\end{lemma}
\begin{proof}
	We will show the inequality for $a^k$ by induction. Since inequality $a_0\le C$ is one of our assumptions, we have the initial step of the induction. To prove the inductive step, consider
	\begin{align*}
		a^{k+1} 
		\le (1 - \gamma_k\mu) a^k + \gamma_k^2 \nu
		\le \left(1 - \frac{2\mu}{\mu k + \theta} \right) \frac{\theta C}{\mu k + \theta} + \theta\mu\frac{C}{(\mu k + \theta)^2}.
	\end{align*}
	To show that the right-hand side is upper bounded by $\frac{\theta C}{\mu(k + 1) + \theta}$, one needs to have, after multiplying both sides by $(\mu k + \theta)(\mu k + \mu + \theta)(\theta C)^{-1}$,
	\begin{align*}
		\left(1 - \frac{2\mu}{\mu k + \theta} \right) (\mu k + \mu + \theta) + \mu\frac{\mu k + \mu + \theta}{\mu k + \theta} \le \mu k + \theta,
	\end{align*}
	which is equivalent to
	\begin{align*}
		\mu - \mu\frac{\mu k + \mu + \theta}{\mu k + \theta} \le 0.
	\end{align*}
	The last inequality is trivially satisfied for all $k \ge 0$.
\end{proof}
We are not ready to prove \Cref{th:str_cvx_decr_steps}. Below we repeat its statement and prove the proof right afterwards.
\begin{theorem}
    Assume that $f$ is $L$-smooth, $\mu$-strongly convex and we have access to its gradients with bounded noise. Set $\gamma_k = \frac{2}{\mu k + \theta}$ with some $\theta \ge 2\max\left\{\frac{\mu}{\alpha}, \frac{(\mu+L)(1+c\alpha)}{2} \right\}$ for some numbers $\alpha > 0$ and $c > 0$ satisfying $\frac{1+Mc\alpha^2}{1+Mc\alpha} \le \alpha_p$. After $k$ iterations of \algname{DIANA} we have
    \begin{align*}
        \ec{ V^k} 
        \le \frac{1}{\eta k+1}\max\left\{ V^0, 4\frac{(1+Mc\alpha)\sigma^2}{M\theta\mu} \right\},
    \end{align*}
    where $\eta\eqdef \frac{\mu}{\theta}$, $V^k=\|x^k - x^\star\|^2 + \frac{c\gamma_k}{M}\sumiM\|h_i^0 - h_i^\star\|^2$ and $\sigma$ is the standard deviation of the gradient noise.
\end{theorem}
\begin{proof}
	To get a recurrence, let us recall an upper bound we have proved before:
    \[
        \ec{ V^{k+1}}
        \le (1 - \gamma_k\mu)\ec{ V^k} + (\gamma_k)^2(1+Mc\alpha)\frac{\sigma^2}{M}.
    \]
    Having that, we can apply Lemma~\ref{lem:sgd_recurrence} to the sequence $\ec{ V^k}$. The constants for the lemma are: $\nu = (1 + Mc\alpha)\frac{\sigma^2}{M}$, $C=\max\left\{V^0, 4\frac{(1+Mc\alpha)\sigma^2}{M\theta\mu}\right\}$, and $\mu$ is the strong convexity constant.
\end{proof}
\begin{corollary}
	If we choose $\alpha = \frac{\alpha_p}{2}$, $c = \frac{4(1-\alpha_p)}{M\alpha_p^2}$,  $\theta=2\max\left\{\frac{\mu}{\alpha}, \frac{\left(\mu+L\right)\left(1 + c\alpha\right)}{2} \right\} = \frac{\mu}{\alpha_p}\max\left\{4, \frac{2(\kappa + 1)}{M} + \frac{(\kappa+1)(M-2)}{ M}\alpha_p\right\}$, then there are three regimes:
	\begin{enumerate}
		\item[1)] if $1 = \max\left\{1,\frac{\kappa}{M},\kappa\alpha_p\right\}$, then $\theta = \Theta\left(\frac{\mu}{\alpha_p}\right)$ and to achieve $\ec{V^K}\le \varepsilon$ we need at most 
		\[
			K=\cO\left( \frac{1}{\alpha_p}\left(V^0+ \frac{(1-\alpha_p)\sigma^2}{M\mu^2} \right)\frac{1}{\varepsilon} \right)
		\] 
		iterations;
		\item[2)] if $\frac{\kappa}{M} = \max\left\{1,\frac{\kappa}{M},\kappa\alpha_p\right\}$, then $\theta = \Theta\left(\frac{L}{M\alpha_p}\right)$ and to achieve $\ec{V^K}\le \varepsilon$ we need at most 
		\[
			K=\cO\left( \frac{\kappa}{M\alpha_p}\left(V^0 + \frac{(1-\alpha_p)\sigma^2}{\mu L} \right)\frac{1}{\varepsilon} \right)
		\] 
		iterations;
		\item[3)] if $\kappa\alpha_p = \max\left\{1,\frac{\kappa}{M},\kappa\alpha_p\right\}$, then $\theta = \Theta\left(L\right)$ and to achieve $\ec{V^K}\le \varepsilon$ we need at most 
		\[
			K=\cO\left( \kappa\left(V^0 + \frac{(1-\alpha_p)\sigma^2}{\mu LM\alpha_p} \right)\frac{1}{\varepsilon} \right)
		\]
		iterations.
	\end{enumerate}		
\end{corollary}
\begin{proof}
	First of all, let us show that $c = \frac{4(1-\alpha_p)}{M\alpha_p^2}$ and $\alpha$ satisfy inequality $\frac{1+Mc\alpha^2}{1+Mc\alpha} \le \alpha_p$:
	\begin{align*}
		\frac{1+Mc\alpha^2}{1+Mc\alpha} \frac{1}{\alpha_p} &= \frac{1 + M\cdot\frac{4(1-\alpha_p)}{M\alpha_p^2}\cdot\frac{\alpha_p^2}{4}}{1 + M\cdot\frac{4(1-\alpha_p)}{M\alpha_p^2}\cdot\frac{\alpha_p}{2}}\cdot\frac{1}{\alpha_p}\\
		&= \frac{2 - \alpha_p}{\alpha_p + 2(1-\alpha_p)}\\
		&= 1.
	\end{align*}
	Moreover, since
	\[
		1+c\alpha = 1 + \frac{2(1-\alpha_p)}{M\alpha_p} = \frac{2 + (M-2)\alpha_P}{M\alpha_p}
	\]
	we can simplify the definition of $\theta$:
	\begin{align*}
		\theta &= 2\max\left\{\frac{\mu}{\alpha}, \frac{\left(\mu+L\right)\left(1 + c\alpha\right)}{2} \right\}\\
		 &= \frac{\mu}{\alpha_p}\max\left\{4, \frac{2(\kappa + 1)}{M} + \frac{(\kappa+1)(M-2)}{ M}\alpha_p\right\}\\
		 &= \Theta\left(\frac{\mu}{\alpha_p}\max\left\{1,\frac{\kappa}{M},\kappa\alpha_p\right\}\right).
	\end{align*}
	Using Theorem~\ref{th:str_cvx_decr_steps}, we get in the case:
	\begin{enumerate}
		\item[1)] if $1 = \max\left\{1,\frac{\kappa}{M},\kappa\alpha_p\right\}$, then $\theta = \Theta\left(\frac{\mu}{\alpha_p}\right)$, $\eta = \Theta\left(\alpha_p\right)$, $\frac{4(1+Mc\alpha)\sigma^2}{M\theta\mu} = \Theta\left(\frac{(1-\alpha_p)\sigma^2}{M\mu^2}\right)$ and to achieve $\ec{V^K}\le \varepsilon$ we need at most 
		\[
			K=\cO\left( \frac{1}{\alpha_p}\max\left(V^0+ \frac{(1-\alpha_p)\sigma^2}{M\mu^2} \right)\frac{1}{\varepsilon} \right)
		\] 
		iterations;
		\item[2)] if $\frac{\kappa}{M} = \max\left\{1,\frac{\kappa}{M},\kappa\alpha_p\right\}$, then $\theta = \Theta\left(\frac{L}{M\alpha_p}\right)$, $\eta = \Theta\left(\frac{\alpha_pM}{\kappa}\right)$, $\frac{4(1+Mc\alpha)\sigma^2}{M\theta\mu} = \Theta\left(\frac{(1-\alpha_p)\sigma^2}{\mu L}\right)$ and to achieve $\ec{V^k}\le \varepsilon$ we need at most 
		\[
			K=\cO\left( \frac{\kappa}{M\alpha_p}\left(V^0+ \frac{(1-\alpha_p)\sigma^2}{\mu L} \right)\frac{1}{\varepsilon} \right)
		\] 
		iterations;
		\item[3)] if $\kappa\alpha_p = \max\left\{1,\frac{\kappa}{M},\kappa\alpha_p\right\}$, then $\theta = \Theta\left(L\right)$, $\eta = \Theta\left(\frac{1}{\kappa}\right)$, $\frac{4(1+Mc\alpha)\sigma^2}{M\theta\mu} = \Theta\left(\frac{(1-\alpha_p)\sigma^2}{\mu L M\alpha_p}\right)$ and to achieve $\ec{V^K}\le \varepsilon$ we need at most 
		\[
			K=\cO\left( \kappa\left(V^0 + \frac{(1-\alpha_p)\sigma^2}{\mu LM\alpha_p} \right)\frac{1}{\varepsilon} \right)
		\] 
	iterations.
	\end{enumerate}		
\end{proof}

\section{Non-Convex Analysis}

\begin{theorem}
	Assume that $\psi$ is constant  and Assumption~\ref{as:almost_identical data} holds.    
    Also assume that $f$ is $L$-smooth, stepsizes $\alpha>0$ and $\gamma_k=\gamma>0$ and parameter $c>0$ satisfying $\frac{1 + M c \alpha^2}{1 + M c \alpha}   \leq \alpha_p,$ $\gamma \le \frac{2}{L(1+2c\alpha)}$ and $\overline x^K$ is chosen randomly from $\{x^0,\dotsc, x^{K-1} \}$. Then
    \begin{align*}
        \ec{ \|\nabla f(\overline x^K)\|^2} \le \frac{2}{K}\frac{\Lambda^0}{\gamma(2 - L\gamma - 2c\alpha L \gamma)} + \frac{(1+2cM\alpha)L\gamma}{2 - L\gamma - 2c\alpha L \gamma}\frac{\sigma^2}{M} + \frac{4c\alpha L\gamma \zeta^2}{2-L\gamma -2c\alpha L\gamma},
    \end{align*}
    where $\Lambda^k\eqdef  f(x^k) - f^\star + c\frac{L\gamma^2}{2}\frac{1}{M}\sumiM \|h_i^{k}-h_i^\star\|^2$.
\end{theorem}
\begin{proof}
	The assumption that $\psi$ is constant implies that $x^{k+1} = x^k - \gamma\hat{g}^k$ and $h^\star = 0$.    
    Moreover, by smoothness of $f$
    \begin{eqnarray*}
        \ec{ f(x^{k+1}) }
        &\overset{\eqref{eq:smoothness_functional}}{\le}& \ec{ f(x^k)} + \ec{\< \nabla f(x^k), x^{k+1} - x^k> } + \frac{L}{2}\ec{\|x^{k+1} - x^k\|^2} \\
        &\overset{\eqref{eq:distr_hat_g_moments1}}{\le}& \ec{ f(x^k)} - \gamma \ec{\|\nabla f(x^k)\|^2} + \frac{L\gamma^2}{2}\ec{\|\hat g^k\|^2} \\
        &\overset{\eqref{eq:full_variance_of_mean_g1}}{\le}& \ec{ f(x^k)} - \left(\gamma - \frac{L\gamma^2}{2}\right)\ec{\|\nabla f(x^k)\|^2}\\
        &&\quad + \frac{L \gamma^2}{2} \frac{1}{M^2}\sumiM \mathbb{E}\left[\Psi(\Delta_i^k)\right] + \frac{L\gamma^2}{2M^2}\sumiM\sigma_i^2.
    \end{eqnarray*}
    Denote $\Lambda^k\eqdef  f(x^k) - f^\star + c\frac{L\gamma^2}{2}\frac{1}{M}\sumiM \|h_i^{k}-h_i^\star\|^2$. Due to Assumption~\ref{as:almost_identical data} we can rewrite the equation~\eqref{eq:distr_h^(k+1)le} after summing it up for $i=1,\ldots,n$ in the following form
	\begin{eqnarray*}
		&& \frac{1}{M}\sum\limits_{i=1}^M\mathbb{E}\left[\|h_i^{k+1}-h_i^\star\|^2 \mid x^k\right]  \\
		&\le& \frac{1-\alpha}{M}\sum\limits_{i=1}^M\|h_i^k - h_i^\star\| + \frac{\alpha}{M}\sum\limits_{i=1}^M\mathbb{E}\left[\|g_i^k - h_i^\star\|^2\mid x^k\right]\\
		&&\quad - \alpha\left(\|\Delta_i^k\|^2 - \alpha\sum\limits_{l=1}^m\|\Delta_i^k(l)\|_1\|\Delta_i^k(l)\|_p\right)\\
		&\le& \frac{1-\alpha}{M}\sum\limits_{i=1}^M\|h_i^k - h_i^\star\| + \frac{2\alpha}{M}\sum\limits_{i=1}^M\mathbb{E}\left[\|g_i^k\|^2\mid x^k\right] + \frac{2\alpha}{M}\sum\limits_{i=1}^M\|h_i^\star - \underbrace{h^\star}_{0}\|^2\\
		&&\quad - \alpha\left(\|\Delta_i^k\|^2 - \alpha\sum\limits_{l=1}^m\|\Delta_i^k(l)\|_1\|\Delta_i^k(l)\|_p\right)\\
		&\overset{\eqref{as:almost_identical data}}{\le}& \frac{1-\alpha}{M}\sum\limits_{i=1}^M\|h_i^k - h_i^\star\| + \frac{2\alpha}{M}\sum\limits_{i=1}^M\|\nabla f_i(x^k)\|^2 + \frac{2\alpha}{M}\sum\limits_{i=1}^M\sigma_i^2 + 2\alpha\zeta^2\\
		&&\quad - \alpha\left(\|\Delta_i^k\|^2 - \alpha\sum\limits_{l=1}^m\|\Delta_i^k(l)\|_1\|\Delta_i^k(l)\|_p\right)\\
		&\overset{\eqref{as:almost_identical data} + \eqref{eq:second_moment_decomposition}}{\le}& \frac{1-\alpha}{M}\sum\limits_{i=1}^M\|h_i^k - h_i^\star\| + 2\alpha\|\nabla f(x^k)\|^2 + 2\alpha\sigma^2 + 4\alpha\zeta^2\\
		&&\quad - \alpha\left(\|\Delta_i^k\|^2 - \alpha\sum\limits_{l=1}^m\|\Delta_i^k(l)\|_1\|\Delta_i^k(l)\|_p\right),
	\end{eqnarray*}	    
	If we add it to the bound above, we  get
    \begin{align*}
        \ec{\Lambda^{k+1}}
        &= \ec{ f(x^{k+1}) - f^\star} + c\frac{L\gamma^2}{2}\frac{1}{M}\sumiM\ec{ \|h_i^{k+1}-h_i^\star\|^2}\\
        &\le \ec{ f(x^k) - f^\star} - \gamma \left( 1 - \frac{L\gamma}{2} - c\alpha L\gamma \right)\ec{ \|\nabla f(x^k)\|^2}\\
        &\quad + (1 - \alpha)c\frac{L\gamma^2}{2}\frac{1}{M}\sumiM \ec{\|h_i^{k}-h_i^\star\|^2}  + \frac{L \gamma^2}{2} \frac{1}{M^2}\sumiM\sum_{l=1}^m \ec{ T_i^k(l)}\\
        &\quad + (1+2cM\alpha)\frac{L\gamma^2}{2}\frac{\sigma^2}{M} + 2c\alpha L\gamma^2\zeta^2,
    \end{align*}
    where
    \[
    		T_i^k(l)
    		\eqdef \left[(\|\Delta_i^k(l)\|_1\|\Delta_i^k(l)\|_p - \|\Delta_i^k(l)\|^2) - nc\alpha(\|\Delta_i^k(l)\|^2 - \alpha\|\Delta_i^k(l)\|_1\|\Delta_i^k(l)\|_p)\right].
    	\]
    As we have shown before, we have $ T_i^k(l)\leq 0$ for all $\Delta_i^k(l)\in \R^{d_l}$.
	Putting all together we have
    \begin{align*}
    	\ec{\Lambda^{k+1}} &\le \ec{ f(x^k) - f^\star} + c\frac{L\gamma^2}{2}\frac{1}{M}\sumiM \ec{\|h_i^{k}-h_i^\star\|^2}+ 2c\alpha L\gamma^2\zeta^2 \\
    	&\qquad + (1+2cM\alpha)\frac{L\gamma^2}{2}\frac{\sigma^2}{M} - \gamma \left( 1 - \frac{L\gamma}{2} - c\alpha L\gamma\right) \ec{\|\nabla f(x^k)\|^2}.
    \end{align*}
    Due to $\gamma \le \frac{2}{L(1+2c\alpha)}$ the coefficient before $\|\nabla f(x^k)\|^2$ is positive. Therefore, we can rearrange the terms and rewrite the last bound as
    \begin{align*}
        \mathbb{E}[\|\nabla f(x^k)\|^2] \le 2\frac{\ec{\Lambda^{k}} - \ec{\Lambda^{k+1}}}{\gamma(2 - L\gamma - 2c\alpha L \gamma)} + \frac{(1+2cM\alpha)L\gamma}{2 - L\gamma - 2c\alpha L \gamma}\frac{\sigma^2}{M} + \frac{4c\alpha L\gamma \zeta^2}{2-L\gamma -2c\alpha L\gamma}.
    \end{align*}
    Summing from $0$ to $k-1$ results in telescoping of the right-hand side, giving
    \begin{align*}
        \sum_{l=0}^{k-1}\mathbb{E}[\|\nabla f(x^l)\|^2] &\le 2\frac{\Lambda^{0} - \ec{\Lambda^{k}}}{\gamma(2 - L\gamma - 2c\alpha L \gamma)} + k\frac{(1+2cM\alpha)L\gamma}{2 - L\gamma - 2c\alpha L \gamma}\frac{\sigma^2}{M} + k\frac{4c\alpha L\gamma \zeta^2}{2-L\gamma -2c\alpha L\gamma}.
    \end{align*}
    Note that $\ec{\Lambda^k}$ is nonnegative and, thus, can be dropped. After that, it suffices to divide both sides by $k$ and rewrite the left-hand side as $\ec{\|\nabla f(\overline x^k)\|^2 }$ where expectation is taken w.r.t.\ all randomness.
\end{proof}

\begin{corollary}\label{cor:ncvx}
	Set $\alpha = \frac{\alpha_p}{2}$, $c = \frac{4(1-\alpha_p)}{M\alpha_p^2}$, $\gamma = \frac{M\alpha_p}{L(4 + (M-4)\alpha_p)\sqrt{K}}$, $h^0 = 0$ and run the algorithm for $K$ iterations. Then, the final accuracy is at most $\frac{2}{\sqrt{K}} \frac{L(4+(M-4)\alpha_p)}{M\alpha_p} \Lambda^0 + \frac{1}{\sqrt{K}}\frac{(4-3\alpha_p)\sigma^2}{4+(M-4)\alpha_p} + \frac{8(1-\alpha_p)\zeta^2}{(4+(M-4)\alpha_p)\sqrt{K}}$.
\end{corollary}
\begin{proof}
	Our choice of $\alpha$ and $c$ implies
	\[
		c\alpha = \frac{2 ( 1 - \alpha_p )}{ n \alpha_p},\quad 1 + 2c\alpha = \frac{4+(M-4)\alpha_p}{M\alpha_p},\quad  1+2cM\alpha = \frac{4-3\alpha_p}{\alpha_p}.
	\]
	Using this and the inequality $\gamma = \frac{M\alpha_p}{L(4 + (M-4)\alpha_p)\sqrt{K}} \le \frac{M\alpha_p}{L(4 + (M-4)\alpha_p)}$ we get $2 - L\gamma - 2c\alpha L \gamma = 2 - (1+2c\alpha)L\gamma \ge 1$. Putting all together we obtain 
	\begin{align*}
		&\frac{2}{K}\frac{\Lambda^0}{\gamma\left(2 - L\gamma - 2c\alpha L \gamma\right)} + (1+2cM\alpha)\frac{L\gamma}{2 - L\gamma - 2c\alpha L \gamma}\frac{\sigma^2}{M} +\frac{4c\alpha L\gamma \zeta^2}{2-L\gamma -2c\alpha L\gamma} \\
		&\le \frac{2}{\sqrt{K}} \frac{L(4+(M-4)\alpha_p)}{M\alpha_p} \Lambda^0 + \frac{1}{\sqrt{K}}\frac{(4-3\alpha_p)\sigma^2}{4+(M-4)\alpha_p} + \frac{8(1-\alpha_p)\zeta^2}{(4+(M-4)\alpha_p)\sqrt{K}}.
	\end{align*}
\end{proof}


\section{Momentum Version of \algname{DIANA}}\label{sec:DIANA-momentum}
\begin{theorem}\label{thm:DIANA-momentum}
	Assume that $f$ is $L$-smooth, $\psi\equiv \text{const}$, $h_i^0 = 0$ and Assumption~\ref{as:almost_identical data} holds. Choose $0\le \alpha < \alpha_p$, $\beta < 1-\alpha$ and $\gamma < \frac{1-\beta^2}{2L\left(2\omega - 1\right)}$, such that $\frac{\beta^2}{(1-\beta)^2\alpha} \le \frac{1-\beta^2-2L\gamma\left(2\omega -1\right)}{\gamma^2L^2\delta}$, where $\delta \eqdef 1 + \frac{2}{M}\left(\frac{1}{\alpha_p}-1\right)\left(1+\frac{\alpha}{1-\alpha-\beta}\right)$ and $\omega \eqdef \frac{M-1}{M} + \frac{1}{M\alpha_p}$, and sample $\overline x^K$ uniformly from $\{x^0, \dotsc, x^{K-1}\}$. Then
	\begin{align*}
		\ec{ \|\nabla f(\overline x^K)\|^2} 
		&\le \frac{4(f(z^0) - f^\star)}{\gamma K} + 2\gamma\frac{L \sigma^2}{(1-\beta)^2  M}\left(\frac{3}{\alpha_p}-2\right) + 2\gamma^2\frac{ L^2\beta^2\sigma^2}{(1 - \beta)^5M}\left(\frac{3}{\alpha_p}-2\right)\\
		&\quad + 3\gamma^2\frac{L^2\beta^2\zeta^2}{(1-\beta)^5M}\left(\frac{1}{\alpha_p}-1\right).
	\end{align*}
\end{theorem}
\begin{proof}
	The main idea of the proof is to find virtual iterates $z^k$ whose recursion would satisfy $z^{k+1} = z^k - \frac{\gamma}{1-\beta} \hat g^k$. Having found it, we can prove convergence by writing a recursion on $f(z^k)$. One possible choice is defined below:
	\begin{equation}
		z^k \eqdef x^k - \frac{\gamma \beta}{1 - \beta} v^{k-1}, \label{eq:def_zk}
	\end{equation}
	where for the edge case $k=0$ we simply set $v^{-1}=0$ and $z^0=x^0$.
	Although $z^k$ is just a slight perturbation of $x^k$, applying smoothness inequality~\eqref{eq:smoothness_functional} to it produces a more convenient bound than the one we would have if used $x^k$. But first of all, let us check that we have the desired recursion for $z^{k+1}$:
	\begin{eqnarray*}
		z^{k+1} 
		&\overset{\eqref{eq:def_zk}}{=}& x^{k+1} -  \frac{\gamma \beta}{1 - \beta} v^{k}  \\
		&{=}& x^k -  \frac{\gamma}{1 - \beta} v^{k} \\
		&{=}& x^k -  \frac{\gamma \beta}{1 - \beta} v^{k-1} -  \frac{\gamma}{1 - \beta} \hat g^k \\
		&\overset{\eqref{eq:def_zk}}{=}& z^k - \frac{\gamma}{1 - \beta} \hat g^k.
	\end{eqnarray*}
	Now, it is time to apply smoothness of $f$:
	\begin{eqnarray}
		\ec{ f(z^{k+1}) }
		&\le& \mathbb{E} \left[f(z^k) + \< \nabla f(z^k), z^{k+1} - z^k> + \frac{L}{2}\|z^{k+1} - z^k\|^2 \right] \nonumber\\
		&\overset{\eqref{eq:def_zk}}{=}& \mathbb{E} \left[f(z^k) - \frac{\gamma}{1 - \beta} \< \nabla f(z^k), \hat g^k> + \frac{L\gamma^2}{2(1-\beta)^2}\|\hat g^k\|^2 \right] . \label{eq:technical2}
	\end{eqnarray}
	The scalar product in~\eqref{eq:technical2} can be bounded using the fact that for any vectors $a$ and $b$ one has $-\< a, b> = \frac{1}{2}(\|a - b\|^2 - \|a\|^2 - \|b\|^2)$. In particular,
	\begin{align*}
		 -  \< \nabla f(z^k), \nabla f(x^k)> 
		 &= \frac{1}{2}\left(\|\nabla f(x^k) - \nabla f(z^k)\|^2 - \|\nabla f(x^k)\|^2 - \|\nabla f(z^k)\|^2 \right) \\
		 &\le  \frac{1}{2}\left(L^2\|x^k - z^k\|^2 - \|\nabla f(x^k)\|^2\right) \\
		 &= \frac{\gamma^2L^2\beta^2}{2(1 - \beta)^2}\|v^{k-1}\|^2 - \frac{1}{2}\|\nabla f(x^k)\|^2.
	\end{align*}
	The next step is to come up with an inequality for $\ec{\|v^k\|^2}$. Since we initialize $v^{-1}=0$, one can show by induction that 
	\begin{equation*}
		v^k = \sum_{l=0}^{k}\beta^{l} \hat g^{k - l}.
	\end{equation*}
	Define $B \eqdef \sum_{l=0}^k \beta^l = \frac{1 - \beta^{k+1}}{1 - \beta}$. Then, by Jensen's inequality
	\begin{align*}
		\ec{\|v^k\|^2} 
		&= B^2\ec{\biggl\|\sum_{l=0}^{k}\frac{\beta^{l}}{B} \hat g^{k - l} \biggr\|^2 }
		\le B^2 \sum_{l=0}^{k}\frac{\beta^{l}}{B} \ec{\|\hat g^{k - l}\|^2}.
	\end{align*}
	Since $\alpha < \alpha_p \le \alpha_p(d_l) \le \frac{\|\Delta_i^k(l)\|^2}{\|\Delta_i^k(l)\|_1\|\Delta_i^k(l)\|_p}$ for all $i,k$ and $l$, we have
	\[
		\|\Delta_i^k(l)\|^2 - \alpha\|\Delta_i^k(l)\|_1\|\Delta_i^k(l)\|_p \ge (\alpha_p - \alpha)\|\Delta_i^k(l)\|_1\|\Delta_i^k(l)\|_p \ge 0
	\]
	for the case when $\Delta_i^k(l)\neq 0$. When $\Delta_i^k(l) = 0$ we simply have $\|\Delta_i^k(l)\|^2 - \alpha\|\Delta_i^k(l)\|_1\|\Delta_i^k(l)\|_p = 0$. Taking into account this and the following equality
	\[
		\|\Delta_i^k\|^2 - \alpha\sum\limits_{l=1}^m\|\Delta_i^k(l)\|_1\|\Delta_i^k(l)\|_p = \sum\limits_{l=1}^m\left(\|\Delta_i^k(l)\|^2 - \alpha\|\Delta_i^k(l)\|_1\|\Delta_i^k(l)\|_p\right),	
	\]
	we get from \eqref{eq:distr_h^(k+1)le}
	\begin{align*}
		\ec{\|h_i^k\|^2} 
		&\le (1-\alpha)\mathbb{E}\left[\|h_i^{k-1}\|^2\right] + \alpha\mathbb{E}\left[\|g_i^{k-1}\|^2\right]\\
		&\le (1-\alpha)^2\mathbb{E}\left[\|h_i^{k-2}\|^2\right] + \alpha(1-\alpha)\mathbb{E}\left[\|g_i^{k-2}\|^2\right] + \alpha\mathbb{E}\left[\|g_i^{k-1}\|^2\right]\\
		&\le \ldots \le (1-\alpha)^k\underbrace{\|h_i^0\|^2}_0 + \alpha\sum\limits_{j=0}^{k-1}(1-\alpha)^{j}\mathbb{E}\left[\|g_i^{k-1-j}\|^2\right].
	\end{align*}
	Next, let us use bounded variance of the stochastic gradients to improve the upper bound to
	\begin{align*}
		\ec{\|h_i^k\|^2} 
		&\le \alpha\sum\limits_{j=0}^{k-1}(1-\alpha)^j\ec{\|\nabla f_i(x^{k-1-j})\|^2} + \alpha\sum\limits_{j=0}^{k-1}(1-\alpha)^j\sigma_i^2\\
		&\le \alpha\sum\limits_{j=0}^{k-1}(1-\alpha)^j\ec{\|\nabla f_i(x^{k-1-j})\|^2} + \alpha\cdot\frac{\sigma_i^2}{1-(1-\alpha)}\\
		&= \alpha\sum\limits_{j=0}^{k-1}(1-\alpha)^j\ec{\|\nabla f_i(x^{k-1-j})\|^2} + \sigma_i^2.
	\end{align*}
	Under our special assumptions, inequality~\eqref{eq:full_second_moment_of_hat_g} gives us
	\begin{eqnarray*}
		\mathbb{E}\left[\|\hat g^k\|^2\right] 
		&\le & \ec{\|\nabla f(x^k)\|^2} + \left(\frac{1}{\alpha_p} - 1\right)\frac{1}{M^2}\sumiM\mathbb{E}\underbrace{\left[\|\nabla f_i(x^k) - h_i^k\|^2\right]}_{\le 2\|\nabla f_i(x^k)\|^2 + 2\|h_i^k\|^2} + \frac{\sigma^2}{\alpha_p M}\\
		&\le & \ec{\|\nabla f(x^k)\|^2} + \frac{2}{M^2}\left(\frac{1}{\alpha_p}-1\right)\sum\limits_{i=1}^M\|\nabla f_i(x^k)\|^2\\
		&&\quad + \frac{2}{M^2}\left(\frac{1}{\alpha_p}-1\right)\sumiM \ec{\|h_i^k\|^2} + \frac{\sigma^2}{\alpha_pM}\\		
		&\overset{\eqref{eq:almost_identical_data}}{\le} & \left(1+\frac{2}{M}\left(\frac{1}{\alpha_p}-1\right)\right)\ec{\|\nabla f(x^k)\|^2} + \frac{2}{M}\left(\frac{1}{\alpha_p}-1\right)\zeta^2\\
		&&\quad  + \frac{2}{M^2}\left(\frac{1}{\alpha_p}-1\right)\sumiM \ec{\|h_i^k\|^2} + \frac{\sigma^2}{\alpha_pM}\\
		&\le & \left(1+\frac{2}{M}\left(\frac{1}{\alpha_p}-1\right)\right)\ec{\|\nabla f(x^k)\|^2}\\
		&&\quad + \frac{2\alpha}{M^2}\left(\frac{1}{\alpha_p}-1\right)\sum\limits_{i=1}^M\sum\limits_{j=0}^{k-1}(1-\alpha)^j\ec{\|\nabla f_i(x^{k-1-j})\|^2}\\
		&&\quad + \left(\frac{1}{\alpha_p}-1\right)\frac{2\sigma^2+2\zeta^2}{M} + \frac{\sigma^2}{\alpha_pM}\\
		&\overset{\eqref{eq:almost_identical_data}}{\le}&  \left(1+\frac{2}{M}\left(\frac{1}{\alpha_p}-1\right)\right)\ec{\|\nabla f(x^k)\|^2} \\
		&&\quad+ \frac{2\alpha}{M}\left(\frac{1}{\alpha_p}-1\right)\sum\limits_{j=0}^{k-1}(1-\alpha)^j\ec{\|\nabla f(x^{k-1-j})\|^2}\\
		&&\quad + \left(\frac{1}{\alpha_p}-1\right)\frac{2\sigma^2+2\zeta^2}{M} + \frac{\sigma^2}{\alpha_pM}  + \frac{2\alpha}{M}\left(\frac{1}{\alpha_p}-1\right)\sum\limits_{j=0}^{k-1}(1-\alpha)^j\zeta^2\\
		&\le& \left(1+\frac{2}{M}\left(\frac{1}{\alpha_p}-1\right)\right)\ec{\|\nabla f(x^k)\|^2}\\
		&&\quad + \frac{2\alpha}{M}\left(\frac{1}{\alpha_p}-1\right)\sum\limits_{j=0}^{k-1}(1-\alpha)^j\ec{\|\nabla f(x^{k-1-j})\|^2}\\
		&&\quad + \left(\frac{1}{\alpha_p}-1\right)\frac{2\sigma^2+3\zeta^2}{M} + \frac{\sigma^2}{\alpha_pM}.
	\end{eqnarray*}
	Using this, we continue our evaluation of $\ec{\|v^k\|^2}$:
	\begin{eqnarray*}
		\ec{\|v^k\|^2}
		 &\le& B\sum\limits_{l=0}^k\beta^l\left(1+\frac{2}{M}\left(\frac{1}{\alpha_p}-1\right)\right)\ec{\|\nabla f(x^{k-l})\|^2}\\
		&&\quad + B\left(\frac{1}{\alpha_p}-1\right)\frac{2\alpha}{M}\sum\limits_{l=0}^k\sum\limits_{j=0}^{k-l-1}\beta^l(1-\alpha)^j\ec{\|\nabla f(x^{k-l-1-j})\|^2}\\
		&&\quad + B\sum\limits_{l=0}^k\beta^l\left(\left(\frac{1}{\alpha_p}-1\right)\frac{2\sigma^2+3\zeta^2}{M} + \frac{\sigma^2}{\alpha_pM}\right).
	\end{eqnarray*}
	Now we are going to simplify the double summation:
	\begin{align*}
		\sum\limits_{l=0}^k\sum\limits_{j=0}^{k-l-1}\beta^l(1-\alpha)^j\ec{\|\nabla f(x^{k-l-1-j})\|^2} &=\sum\limits_{l=0}^k\sum\limits_{j=0}^{k-l-1}\beta^l(1-\alpha)^{k-l-1-j}\ec{\|\nabla f(x^{j})\|^2}\\
		 &= \sum\limits_{j=0}^{k-1}\ec{\|\nabla f(x^{j})\|^2}\sum\limits_{l=0}^{k-j-1}\beta^l(1-\alpha)^{k-l-1-j}\\
		 &= \sum\limits_{j=0}^{k-1}\ec{\|\nabla f(x^{j})\|^2}\cdot \frac{(1-\alpha)^{k-j} - \beta^{k-j}}{1-\alpha-\beta}\\
		 &\le \sum\limits_{j=0}^{k}\ec{\|\nabla f(x^{j})\|^2}\cdot \frac{(1-\alpha)^{k-j}}{1-\alpha-\beta}\\
		 &= \frac{1}{1-\alpha - \beta}\sum\limits_{j=0}^{k}(1-\alpha)^{j}\ec{\|\nabla f(x^{k-j})\|^2}.
	\end{align*}
	Note that $B \eqdef \sum\limits_{l=0}^k\beta^l \le \frac{1}{1-\beta}$. Putting all together we get
	\begin{eqnarray*}
		\ec{\|v^k\|^2} &\le& \frac{\delta}{1-\beta}\sum\limits_{l=0}^k(1-\alpha)^l\ec{\|\nabla f(x^{k-l})\|^2} + \frac{\sigma^2}{M(1-\beta)^2}\left(\frac{3}{\alpha_p}-2\right)\\
		&&\quad + \frac{3\zeta^2}{M(1-\beta)^2}\left(\frac{1}{\alpha_p}-1\right),
	\end{eqnarray*}
	where $\delta \eqdef 1 + \frac{2}{M}\left(\frac{1}{\alpha_p}-1\right)\left(1+\frac{\alpha}{1-\alpha-\beta}\right)$, and as a result
	\begin{align*}
		\frac{\gamma^3L^2\beta^2}{2(1-\beta)^3}\ec{\|v^{k-1}\|^2}
		 &\le \frac{\gamma^3L^2\beta^2\delta}{2(1-\beta)^4}\sum\limits_{l=0}^{k-1}(1-\alpha)^{k-1-l}\ec{\|\nabla f(x^{l})\|^2} \\
		 &\quad + \frac{\gamma^3L^2\beta^2\sigma^2}{2M(1-\beta)^5}\left(\frac{3}{\alpha_p}-2\right) + \frac{3\gamma^3L^2\beta^2\zeta^2}{2M(1-\beta)^5}\left(\frac{1}{\alpha_p}-1\right).
	\end{align*}
	To sum up, we have
	\begin{eqnarray*}
		\mathbb{E}\left[f(z^{k+1})\right] &\le& \mathbb{E}\left[f(z^k)\right] -\frac{\gamma}{2(1-\beta)}\left(1 - \frac{L\gamma\omega}{1-\beta}\right)\ec{\|\nabla f(x^k)\|^2}\\
		&&\quad + \left(\frac{L\gamma^2\alpha(\omega-1)}{2(1-\beta)^2}+\frac{\gamma^3L^2\beta^2\delta}{2(1-\beta)^4}\right)\sum\limits_{l=0}^{k-1}(1-\alpha)^{k-1-l}\ec{\|\nabla f(x^l)\|^2}\\
		&&\quad + \frac{\sigma^2}{M}\left(\frac{3}{\alpha_p}-2\right)\left(\frac{L\gamma^2}{2(1-\beta)^2}+\frac{\gamma^3L^2\beta^2}{2(1-\beta)^5}\right) + \frac{3\gamma^3L^2\beta^2\zeta^2}{2M(1-\beta)^5}\left(\frac{1}{\alpha_p}-1\right).
	\end{eqnarray*}
	Telescoping this inequality from 0 to $K-1$, we get
	\begin{align*}
		&\ec{ f(z^K) - f(z^0)}\\
		&\le K\frac{\sigma^2}{M}\left(\frac{3}{\alpha_p}-2\right)\left(\frac{L\gamma^2}{2(1-\beta)^2}+\frac{\gamma^3L^2\beta^2}{2(1-\beta)^5}\right) + K\frac{3\gamma^3L^2\beta^2\zeta^2}{2M(1-\beta)^5}\left(\frac{1}{\alpha_p}-1\right)\\
		&\quad + \frac{\gamma}{2}\sum\limits_{l=0}^{K-2}\left(\left(\frac{L\gamma\alpha(\omega-1)}{(1-\beta)^2}+\frac{\gamma^2L^2\beta^2\delta}{(1-\beta)^4}\right)\sum\limits_{k'=l+1}^{K-1}(1-\alpha)^{k'-1-l}\right)\ec{\|\nabla f(x^l)\|^2}\\
		&\quad + \frac{\gamma}{2}\sum\limits_{l=0}^{K-2}\left(\frac{L\gamma\omega}{(1-\beta)^2} - \frac{1}{1-\beta}\right)\ec{\|\nabla f(x^l)\|^2}\\
		&\quad + \frac{\gamma}{2}\left(\frac{L\gamma\omega}{(1-\beta)^2} - \frac{1}{1-\beta}\right)\ec{\|\nabla f(x^{K-1})\|^2}\\
		&\le K\frac{\sigma^2}{M}\left(\frac{3}{\alpha_p}-2\right)\left(\frac{L\gamma^2}{2(1-\beta)^2}+\frac{\gamma^3L^2\beta^2}{2(1-\beta)^5}\right)+ K\frac{3\gamma^3L^2\beta^2\zeta^2}{2M(1-\beta)^5}\left(\frac{1}{\alpha_p}-1\right)\\
		&\quad + \frac{\gamma}{2}\sum\limits_{l=0}^{K-1}\left(\frac{\gamma^2L^2\beta^2\delta}{(1-\beta)^4\alpha} + \frac{L\gamma}{(1-\beta)^2}\left(2\omega-1\right) - \frac{1}{1-\beta}\right)\ec{\|\nabla f(x^l)\|^2}
	\end{align*}
	It holds $f^\star\le f(z^K)$ and our assumption on $\beta$ implies that $\frac{\gamma^2L^2\beta^2\delta}{(1-\beta)^4\alpha} + \frac{L\gamma}{(1-\beta)^2}\left(2\omega -1\right) - \frac{1}{1-\beta} \le -\frac{1}{2}$, so it all results in
	\begin{align*}
		\frac{1}{K}\sum_{l=0}^{K-1} \|\nabla f(x^{l})\|^2 &\le  \frac{4(f(z^0) - f^\star)}{\gamma K} + 2\gamma\frac{L \sigma^2}{(1-\beta)^2  M}\left(\frac{3}{\alpha_p}-2\right)\\
		&\quad + 2\gamma^2\frac{ L^2\beta^2\sigma^2}{(1 - \beta)^5M}\left(\frac{3}{\alpha_p}-2\right) + \frac{3\gamma^2L^2\beta^2\zeta^2}{M(1-\beta)^5}\left(\frac{1}{\alpha_p}-1\right).
	\end{align*}
	Since $\overline x^K$ is sampled uniformly from $\{x^0, \dotsc, x^{K-1}\}$, the left-hand side is equal to $\ec{ \|\nabla f(\overline x^K)\|^2}$. Also note that $z^0=x^0$.
\end{proof}
\begin{corollary}\label{cor:DIANA-momentum}
		If we set $\gamma=\frac{1-\beta^2}{2\sqrt{K}L\left(2\omega -1\right)}$ and $\beta$ such that $\frac{\beta^2}{(1 - \beta)^2\alpha}\le \frac{4K\left(2\omega -1\right)}{\delta}$ with $K>1$, then the error after $K$ iterations is at most 
		\begin{align*}
			&\frac{1}{\sqrt{K}}\left(\frac{8L(2\omega -1)(f(x^0)-f^\star)}{1-\beta^2} + \frac{(1+\beta)\sigma^2}{(2\omega -1)\alpha_p M(1-\beta)}\left(\frac{3}{\alpha_p}-2\right)\right)\\ 
			&+ \frac{1}{K}\frac{(1+\beta)^4\beta^2\sigma^2}{2(1 - \beta)(2\omega -1)\alpha_pM}\left(\frac{3}{\alpha_p}-2\right)
			+\frac{1}{K}\frac{3(1+\beta)^4\beta^2\zeta^2}{4(1 - \beta)(2\omega -1)\alpha_pM}\left(\frac{1}{\alpha_p}-1\right).
		\end{align*}
\end{corollary}
\begin{proof}
	Our choice of $\gamma = \frac{1-\beta^2}{2\sqrt{K}L\left(2\omega -1\right)}$ implies that
	\[
		\frac{\beta^2}{(1 - \beta)^2\alpha}
		\le \frac{4K\left(2\omega -1\right)}{\delta} \Longleftrightarrow \frac{\beta^2}{(1-\beta)^2\alpha} \le \frac{1-\beta^2-2L\gamma\left(2\omega -1\right)}{\gamma^2L^2\delta}.
	\]
	After that, it remains to plug-in $\gamma = \frac{1-\beta^2}{2\sqrt{K}L\left(2\omega -1\right)}$ in 
	\[
		\frac{4(f(z^0) - f^\star)}{\gamma K} + \frac{2\gamma L \sigma^2}{(1-\beta)^2  M}\left(\frac{3}{\alpha_p}-2\right) + \frac{2\gamma^2 L^2\beta^2\sigma^2}{(1 - \beta)^5M}\left(\frac{3}{\alpha_p}-2\right) + \frac{3\gamma^2L^2\beta^2\zeta^2}{M(1-\beta)^5}\left(\frac{1}{\alpha_p}-1\right)
	\] 
	to get the desired result. 
\end{proof}

\section{Analysis of \algname{DIANA} with $\alpha = 0$ and $h_i^0 = 0$}
\label{sec:Terngrad}

\subsection{Convergence Rate of \algname{TernGrad}}

Here we give the convergence guarantees for \algname{TernGrad} and provide upper bounds for this method. The method coincides with Algorithm~\ref{alg:terngrad} for the case when $p = \infty$. In the original paper \cite{wen2017terngrad} no convergence rate was given and we close this gap.

To maintain consistent notation we rewrite the \algname{TernGrad} in notation which is close to the notation we used for \algname{DIANA}. Using our notation it is easy to see that \algname{TernGrad} is \algname{DIANA} with $h_1^0 = h_2^0 = \ldots = h_M^0 = 0, \alpha = 0$ and $p=\infty$. Firstly, it means that $h_i^k = 0$ for all $i=1,2,\ldots,n$ and $k\ge 1$. What is more, this observation tells us that Lemma \ref{lem:3in1} holds for the iterates of \algname{TernGrad} too. What is more, in the original paper \cite{wen2017terngrad} the quantization parameter $p$ was chosen as $\infty$. We generalize the method and we don't restrict our analysis only on the case of $\ell_\infty$ sampling.

As it was in the analysis of \algname{DIANA} our proofs for \algname{TernGrad} work under Assumption~\ref{as:noise}.

\begin{algorithm}[ht]
   \caption{\algname{DIANA} with $\alpha=0$ and $h_i^0 = 0$; equivalent to \algname{QSGD} for $p=2$ (1-bit)/ \algname{TernGrad} for $p=\infty$ (\algname{SGD})}
   \label{alg:terngrad}
\begin{algorithmic}[1]
   \State \textbf{Input:}  stepsizes $(\gamma_k)_{k}$, initial vector $x^0$, quantization parameter $p \geq 1$, sizes of blocks $(d_l)_{l=1}^m$, momentum parameter $0\le \beta < 1$, number of steps $K$
   \State $v^0 = \nabla f(x^0)$
   \For{$k=1,2,\dotsc, K-1$}
	   \State Broadcast $x^{k}$ to all workers
        \For{$i=1,\dotsc,n$ do in parallel}
			\State Sample $g^{k}_i$ such that $\mathbb{E} [g^k_i \;|\; x^k]  =\nabla f_i(x^k)$
			\State Sample $\hat g^k_i \sim {\rm Quant}_{p}(g^k_i,(d_l)_{l=1}^m)$
        \EndFor
        \State $\hat g^k = \frac{1}{M}\sum_{i=1}^M \hat g_i^k$
        \State  $v^k = \beta v^{k-1} + \hat g^k$
        \State $x^{k+1} = \prox_{\gamma_k \psi}\left(x^k - \gamma_kv^k \right)$
   \EndFor
\end{algorithmic}
\end{algorithm}

\subsection{Technical lemmas}

First of all, we notice that since \algname{TernGrad} coincides with \algname{DIANA}, having $h_i^k = 0, i,k \ge 1$, $\alpha = 0$ and $p = \infty$, all inequalities from Lemma~\ref{lem:3in1} holds for the iterates of \algname{TernGrad} as well because $\Delta_i^k = g_i^k$ and $\hat \Delta_i^k = \hat g_i^k$.

\begin{lemma}\label{lem:gamma_choice_terngrad_special}
	Assume $\gamma \le \frac{M\alpha_p}{L((M-1)\alpha_p+1)}$. Then
	\begin{align}
		2\gamma\mu\left(1 - \frac{\gamma L((M-1)\alpha_p+1)}{2M\alpha_p}\right) \ge \gamma\mu.\label{eq:conseq_gamma_choice_terngrad_special}
	\end{align}
\end{lemma}
\begin{proof}
	Since $\gamma \le \frac{M\alpha_p}{L((M-1)\alpha_p+1)}$ we have
	$$
		2\gamma\mu\left(1 - \frac{\gamma L((M-1)\alpha_p+1)}{2M\alpha_p}\right) \ge 2\gamma\mu \left(1-\frac{1}{2}\right) = \gamma\mu.
	$$
\end{proof}

\begin{lemma}\label{lem:gamma_choice_terngrad}
	Assume $\gamma \le \frac{1}{L(1+\kappa(1 - \alpha_p)/(M\alpha_p))}$, where $\kappa \eqdef \frac{L}{\mu}$ is the condition number of $f$. Then
	\begin{align}
		r \ge \gamma\mu\label{eq:conseq_gamma_choice_terngrad},
	\end{align}
	where $r = 2\mu\gamma - \gamma^2 \left(\mu L + \frac{L^2(1-\alpha_p)}{ M\alpha_p}\right)$.
\end{lemma}
\begin{proof}
	Since $\gamma \le \frac{1}{L(1+\kappa(1-\alpha_p)/(M\alpha_p))} = \frac{\mu M\alpha_p}{\mu M\alpha_p L + L^2(1-\alpha_p)}$ we have
	$$
		M\alpha_p r = \gamma\left(2\mu M\alpha_p - \gamma\left(\mu M\alpha_p L + L^2(1-\alpha_p)\right)\right) \ge \gamma\mu M\alpha_p,
	$$
	whence $r \ge \gamma\mu$.
\end{proof}

\begin{lemma}\label{lem:gamma_choice_terngrad_prox}
	Assume $\gamma \le \frac{2M\alpha_p}{(\mu + L)(2+(M-2)\alpha_p)}$. Then
	\begin{eqnarray}\label{eq:conseq_gamma_choice_terngrad_prox}
		2\gamma\mu - \gamma^2\mu^2\left(1 + \frac{2(1-\alpha_p)}{M\alpha_p}\right) &\ge & \gamma\mu.
	\end{eqnarray}
\end{lemma}
\begin{proof}
	Since $\gamma \le \frac{2M\alpha_p}{(\mu + L)(2+(M-2)\alpha_p)}$ we have
	\[
		\gamma\mu \le  \frac{2\mu M\alpha_p}{(\mu + L)(2+(M-2)\alpha_p)} \le  \frac{(\mu + L)M\alpha_p}{(\mu + L)(2+(M-2)\alpha_p)} =  \frac{M\alpha_p}{2+(M-2)\alpha_p},
	\]
	whence
	\begin{align*}
		2\gamma\mu - \gamma^2\mu^2\left(1 + \frac{2(1-\alpha_p)}{M\alpha_p}\right) 
		&\ge 2\gamma\mu - \gamma\mu \frac{M\alpha_p}{2+(M-2)\alpha_p} \left(1 + \frac{2(1-\alpha_p)}{M\alpha_p}\right)\\
		&= 2\gamma\mu - \gamma\mu = \gamma\mu.
	\end{align*}
\end{proof}

\begin{lemma}\label{lem:L_smoothness_consequense}
	Assume that each function $f_i$ is $L$-smooth and $\psi$ is a constant function. Then for the iterates of Algorithm $\ref{alg:terngrad}$ with $\gamma_k = \gamma$ we have
	\begin{align}
    		\ec{\Theta^{k+1}} 
    		&\le  \ec{\Theta^k} + \left(\frac{\gamma^2 L}{2} - \gamma\right)\mathbb{E}\left[\|\nabla f(x^k)\|^2\right] \notag \\
    		&\qquad + \frac{\gamma^2 L}{2M^2}\left(\frac{1}{\alpha_p}-1\right)\sumiM\mathbb{E}\left[\|g_i^k\|^2\right] + \frac{\gamma^2 L\sigma^2}{2M},\label{eq:L_smoothness_consequense}
    \end{align}
    where $\Theta^k = f(x^k) - f(x^\star)$ and $\sigma^2\eqdef \frac{1}{M}\sumiM\sigma_i^2$.
\end{lemma}
\begin{proof}
	Since $\psi$ is a constant function we have $x^{k+1} = x^k - \gamma \hat g^k$. Moreover, from the $L$-smoothness of $f$ we have
	\begin{eqnarray*}
		\ec{\Theta^{k+1}} &\le & \ec{\Theta^k} + \mathbb{E}\left[\langle \nabla f(x^k), x^{k+1} - x^k \rangle\right] + \frac{L}{2}\|x^{k+1}-x^k\|^2\\
		&= & \ec{\Theta^k} - \gamma\mathbb{E}\left[\|\nabla f(x^k)\|^2\right] + \frac{\gamma^2 L}{2}\mathbb{E}\left[\left\|\hat g^k\right\|^2\right]\\
		&\overset{\eqref{eq:full_variance_of_mean_g1}}{\le}& \ec{\Theta^k} + \left(\frac{\gamma^2 L}{2} - \gamma\right)\mathbb{E}\left[\|\nabla f(x^k)\|^2\right]\\
		&&\quad + \frac{\gamma^2 L}{2M^2}\sumiM\sum\limits_{l=1}^m \mathbb{E}\left[\|g_i^k(l)\|_1\|g_i^k(l)\|_p - \|g_i^k(l)\|^2 \right] + \frac{\gamma^2 L}{2M^2}\sumiM\sigma_i^2,
	\end{eqnarray*}
	where the first equality follows from $x^{k+1} - x^k = \hat g^k$, $\mathbb{E}\left[\hat g^k\mid x^k\right] = \nabla f(x^k)$ and the tower property of mathematical expectation. 
	By definition $\alpha_p(d_l) = \inf\limits_{x\neq 0,x\in\R^{d_l}}\frac{\|x\|^2}{\|x\|_1\|x\|_p} = \left(\sup\limits_{x\neq 0,x\in\R^{d_l}}\frac{\|x\|_1\|x\|_p}{\|x\|^2}\right)^{-1}$ and $\alpha_p = \alpha_p(\max\limits_{l=1,\ldots,m}d_l)$ which implies
	\begin{align*}
		\mathbb{E}\left[\|g_i^k(l)\|_1\|g_i^k(l)\|_p - \|g_i^k(l)\|^2 \right] 
		&= \mathbb{E}\left[\|g_i^k(l)\|^2\left(\frac{\|g_i^k(l)\|_1\|g_i^k(l)\|_p}{\|g_i^k(l)\|^2} - 1\right) \right] \\
		&\le \left(\frac{1}{\alpha_p(d_l)}-1\right)\ec{\|g_i^k(l)\|^2}\\
		&\le \left(\frac{1}{\alpha_p}-1\right)\ec{\|g_i^k(l)\|^2}.
	\end{align*}
 	Since $\|g_i^k\|^2 = \sum\limits_{l=1}^m\|g_i^k(l)\|^2$ we have
	\begin{align*}
		\ec{\Theta^{k+1}}
		&\le  \ec{\Theta^k} + \left(\frac{\gamma^2 L}{2} - \gamma\right)\mathbb{E}\left[\|\nabla f(x^k)\|^2\right]\\
		&\qquad + \frac{\gamma^2 L}{2M^2}\left(\frac{1}{\alpha_p}-1\right)\sumiM\mathbb{E}\left[\|g_i^k\|^2\right] + \frac{\gamma^2 L\sigma^2}{2M},
	\end{align*}
	where $\sigma^2 = \frac{1}{M}\sumiM\sigma_i^2$.
\end{proof}

\subsection{Non-convex analysis}\label{sec:TernGrad-nonconvex}
\begin{theorem}\label{thm:TernGrad-nonconvex}
    Assume that $\psi$ is constant  and Assumption~\ref{as:almost_identical data} holds.    
    Also assume that $f$ is $L$-smooth, $\gamma \le \frac{M\alpha_p}{L((M-1)\alpha_p+1)}$ and $\overline x^K$ is chosen randomly from $\{x^0,\dotsc, x^{K-1} \}$. Then
    \begin{align*}
        \ec{ \|\nabla f(\overline x^K)\|^2} 
        \le \frac{2}{K} \frac{f(x^0) - f(x^\star)}{\gamma\left(2 - \gamma\frac{L((M-1)\alpha_p+1)}{M\alpha_p}\right)} + \frac{\gamma L\left(\sigma^2+(1-\alpha_p)\zeta^2\right)}{M\alpha_p}.
    \end{align*}
\end{theorem}
\begin{proof}
	Recall that we defined $\Theta^k$ as $f(x^k) - f(x^\star)$ in Lemma~\ref{lem:L_smoothness_consequense}. From \eqref{eq:L_smoothness_consequense} we have
	\begin{eqnarray*}
		\mathbb{E}[\Theta^{k+1}] 
		&\le &  \mathbb{E} [\Theta^k] + \left(\frac{\gamma^2 L}{2} - \gamma\right)\mathbb{E}\left[\|\nabla f(x^k)\|^2\right]\\
		&&\quad + \frac{\gamma^2 L}{2M^2}\left(\frac{1}{\alpha_p}-1\right)\sumiM\mathbb{E}\left[\|g_i^k\|^2\right] + \frac{\gamma^2 L\sigma^2}{2M}.
	\end{eqnarray*}
	Using variance decomposition
	\begin{equation*}
		\mathbb{E}\left[\|g_i^k\|^2\right] 
		= \mathbb{E}\left[\|\nabla f_i(x^k)\|^2\right] + \mathbb{E}\left[\|g_i^k - \nabla f(x^k)\|^2\right] 
		\le \mathbb{E}\left[\|\nabla f_i(x^k)\|^2\right] + \sigma_i^2,
	\end{equation*}
	we get
	\begin{eqnarray*}
		\frac{1}{M}\sumiM \mathbb{E}\left[\|g_i^k\|^2\right] &\le & \frac{1}{M}\sumiM\mathbb{E}\left[\|\nabla f_i(x^k)\|^2\right] + \sigma^2\\
		&\overset{\eqref{eq:second_moment_decomposition}}{=}& \ec{\|\nabla f(x^k)\|^2} + \frac{1}{M}\sumiM\mathbb{E}\left[\|\nabla f_i(x^k)-\nabla f(x^k)\|^2\right] + \sigma^2\\
		&\overset{\eqref{eq:almost_identical_data}}{\le}& \ec{\|\nabla f(x^k)\|^2} + \zeta^2 + \sigma^2.
	\end{eqnarray*}
	Putting all together we obtain
	\begin{eqnarray}
		\mathbb{E}[\Theta^{k+1}]
		&\le &  \mathbb{E}[\Theta^k] + \left(\frac{\gamma^2L}{2}\cdot \frac{(M-1)\alpha_p+1}{M\alpha_p} - \gamma\right)\mathbb{E}\left[\|\nabla f(x^k)\|^2\right] +  \frac{\gamma^2 L\sigma^2}{2M\alpha_p}\notag\\
		&&\quad + \frac{\gamma^2L\zeta^2(1-\alpha_p)}{2M\alpha_p} \label{eq:terngrad_key_estimation_special}.
	\end{eqnarray}
	Since $\gamma \le \frac{M\alpha_p}{L((M-1)\alpha_p+1)}$ the factor $\left(\frac{\gamma^2L}{2}\cdot \frac{(M-1)\alpha_p+1}{M\alpha_p} - \gamma\right)$ is negative and therefore
	\begin{eqnarray*}
		\mathbb{E}\left[\|\nabla f(x^k)\|^2\right] &\le & \frac{\ec{\Theta^k} - \ec{\Theta^{k+1}}}{\gamma\left(1 - \gamma\frac{L((M-1)\alpha_p+1)}{2M\alpha_p}\right)} + \frac{\gamma L\left(\sigma^2+(1-\alpha_p)\zeta^2\right)}{2M\alpha_p - \gamma L((M-1)\alpha_p+1)}.
	\end{eqnarray*}
	Telescoping the previous inequality from $0$ to $K-1$ and using $\gamma\le \frac{M\alpha_p}{L((M-1)\alpha_p+1)}$ we obtain
	\begin{align*}
		\frac{1}{K}\sum\limits_{l=0}^{K-1}\mathbb{E}\left[\|\nabla f(x^l)\|^2\right] 
		&\le \frac{2}{K} \frac{\ec{\Theta^0} - \ec{\Theta^{K}}}{\gamma\left(2 - \gamma\frac{L((M-1)\alpha_p+1)}{M\alpha_p}\right)} + \frac{\gamma L\left(\sigma^2+(1-\alpha_p)\zeta^2\right)}{M\alpha_p}.
	\end{align*}
	It remains to notice that the left-hand side is just $\mathbb{E}\left[\|\nabla f(\overline x^K)\|^2\right]$, $\Theta^K \ge 0$ and $\ec{\Theta^0} = f(x^0) - f(x^\star)$.
\end{proof}

\begin{corollary}\label{cor:TernGrad-nonconvex}
	If we choose $\gamma = \frac{M\alpha_p}{L((M-1)\alpha_p+1)\sqrt{K}}$ then the rate we get is 
	\[
		\frac{2}{\sqrt{K}}L\frac{(M-1)\alpha_p+1}{M\alpha_p}\left(f(x^0) - f(x^\star)\right) + \frac{1}{\sqrt{K}}\frac{\sigma^2+(1-\alpha_p)\zeta^2}{(1+(M-1)\alpha_p)}.
	\]
\end{corollary}
\begin{proof}
	Our choice of $\gamma = \frac{M\alpha_p}{L((M-1)\alpha_p+1)\sqrt{K}} \le \frac{M\alpha_p}{L((M-1)\alpha_p+1)}$ implies that $2 - \gamma\frac{L((M-1)\alpha_p+1)}{M\alpha_p} \ge 1$. After that it remains to notice that for our choice of $\gamma = \frac{M\alpha_p}{L((M-1)\alpha_p+1)\sqrt{K}}$ we have 
	\begin{align*}
		&\frac{2}{K} \frac{f(x^0) - f(x^\star)}{\gamma\left(2 - \gamma\frac{L((M-1)\alpha_p+1)}{M\alpha_p}\right)} + \frac{\gamma L\left(\sigma^2+(1-\alpha_p)\zeta^2\right)}{M\alpha_p} \\
		&\le \frac{2}{\sqrt{K}}L\frac{(M-1)\alpha_p+1}{M\alpha_p}\left(f(x^0) - f(x^\star)\right)  + \frac{1}{\sqrt{K}}\frac{\sigma^2+(1-\alpha_p)\zeta^2}{(1+(M-1)\alpha_p)}.
	\end{align*}
\end{proof}

\subsection{Momentum version}\label{sec:TernGrad-momentum}
\begin{theorem}\label{thm:TernGrad-momentum}
	Assume that $f$ is $L$-smooth, $\psi\equiv \text{const}$, $\alpha=0$, $h_i=0$ and Assumption~\ref{as:almost_identical data} holds. Choose $\beta<1$ and $\gamma < \frac{1-\beta^2}{2L\omega}$ such that $\frac{\beta^2}{(1 - \beta)^3}\le \frac{1 - \beta^2 - 2L\gamma\omega}{\gamma^2 L^2\omega}$, where $\omega \eqdef \frac{M-1}{M} + \frac{1}{M\alpha_p}$ and sample $\overline x^K$ uniformly from $\{x^0, \dotsc, x^{K-1}\}$. Then
	\begin{align*}
		\ec{ \|\nabla f(\overline x^K)\|^2}
		&\le \frac{4(f(z^0) - f^\star)}{\gamma K} + 2\gamma\frac{L \sigma^2}{\alpha_p  M(1-\beta)^2}+ 2\gamma^2\frac{ L^2\beta^2\sigma^2}{(1 - \beta)^5\alpha_p M}\\
		&\qquad + 2\gamma^2\frac{L^2\beta^2(1-\alpha_p)\zeta^2}{2(1 - \beta)^5\alpha_p M}.
	\end{align*}
\end{theorem}
\begin{proof}
	The main idea of the proof is to find virtual iterates $z^k$ whose recursion would satisfy $z^{k+1} = z^k - \frac{\gamma}{1-\beta} \hat g^k$. Having found it, we can prove convergence by writing a recursion on $f(z^k)$. One possible choice is defined below:
	\begin{align}
		z^k \eqdef x^k - \frac{\gamma \beta}{1 - \beta} v^{k-1}, \label{eq:def_zk_tg}
	\end{align}
	where for the edge case $k=0$ we simply set $v^{-1}=0$ and $z^0=x^0$.
	Although $z^k$ is just a slight perturbation of $x^k$, applying smoothness inequality~\eqref{eq:smoothness_functional} to it produces a more convenient bound than the one we would have if used $x^k$. But first of all, let us check that we have the desired recursion for $z^{k+1}$:
	\begin{eqnarray*}
		z^{k+1} 
		&\overset{\eqref{eq:def_zk_tg}}{=}& x^{k+1} -  \frac{\gamma \beta}{1 - \beta} v^{k}  \\
		&{=}& x^k -  \frac{\gamma}{1 - \beta} v^{k} \\
		&{=}& x^k -  \frac{\gamma \beta}{1 - \beta} v^{k-1} -  \frac{\gamma}{1 - \beta} \hat g^k \\
		&\overset{\eqref{eq:def_zk_tg}}{=}& z^k - \frac{\gamma}{1 - \beta} \hat g^k.
	\end{eqnarray*}
	Now, it is time to apply smoothness of $f$:
	\begin{eqnarray}
		\mathbb{E}[f(z^{k+1})] 
		&\le& \mathbb{E} \left[f(z^k) + \< \nabla f(z^k), z^{k+1} - z^k> + \frac{L}{2}\|z^{k+1} - z^k\|^2 \right] \nonumber\\
		&\overset{\eqref{eq:def_zk_tg}}{=}& \mathbb{E} \left[f(z^k) - \frac{\gamma}{1 - \beta} \< \nabla f(z^k), \hat g^k> + \frac{L\gamma^2}{2(1-\beta)^2}\|\hat g^k\|^2 \right] . \label{eq:technical3}
	\end{eqnarray}
	Under our special assumption inequality~\eqref{eq:full_second_moment_of_hat_g} simplifies to
	\begin{eqnarray*}
		\mathbb{E}\left[\|\hat g^k\|^2 \mid x^k\right] 
		&\le& \|\nabla f(x^k)\|^2 + \left(\frac{1}{\alpha_p} - 1\right)\frac{1}{M^2}\sum\limits_{i=1}^M \|\nabla f_i(x^k)\|^2 + \frac{\sigma^2}{\alpha_p M}\\
		&\overset{\eqref{eq:almost_identical_data}}{\le}& \|\nabla f(x^k)\|^2 + \left(\frac{1}{\alpha_p} - 1\right)\frac{1}{M}\|\nabla f(x^k)\|^2 + \left(\frac{1}{\alpha_p} - 1\right)\frac{\zeta^2}{M} + \frac{\sigma^2}{\alpha_p M}.
	\end{eqnarray*}
	The scalar product in~\eqref{eq:technical3} can be bounded using the fact that for any vectors $a$ and $b$ one has $-2\< a, b> = \|a - b\|^2 - \|a\|^2 - \|b\|^2$. In particular,
	\begin{align*}
		 - \frac{2\gamma}{1 - \beta} \< \nabla f(z^k), \nabla f(x^k)> 
		 &= \frac{\gamma}{1-\beta}\left(\|\nabla f(x^k) - \nabla f(z^k)\|^2 - \|\nabla f(x^k)\|^2 - \|\nabla f(z^k)\|^2 \right) \\
		 &\le  \frac{\gamma}{1-\beta}\left(L^2\|x^k - z^k\|^2 - \|\nabla f(x^k)\|^2\right) \\
		 &= \frac{\gamma^3L^2\beta^2}{(1 - \beta)^3}\|v^{k-1}\|^2 - \frac{\gamma}{1-\beta}\|\nabla f(x^k)\|^2.
	\end{align*}
	The next step is to come up with an inequality for $\ec{\|v^k\|^2}$. Since we initialize $v^{-1}=0$, one can show by induction that 
	\begin{equation*}
		v^k = \sum_{l=0}^{k}\beta^{l} \hat g^{k - l}.
	\end{equation*}
	Define $B \eqdef \sum_{l=0}^k \beta^l = \frac{1 - \beta^{k+1}}{1 - \beta}$. Then, by Jensen's inequality
	\begin{align*}
		\mathbb{E}[\|v^k\|^2 ]
		&= B^2\ec{\biggl\|\sum_{l=0}^{k}\frac{\beta^{l}}{B} \hat g^{k - l} \biggr\|^2 }
		\le B^2 \sum_{l=0}^{k}\frac{\beta^{l}}{B} \ec{\|\hat g^{k - l}\|^2}\\
		&\le B \sum_{l=0}^{k}\beta^{l} \left(\left(\frac{M-1}{M} + \frac{1}{M\alpha_p}\right)\ec{\|\nabla f(x^{k-l})\|^2} + \left(\frac{1}{\alpha_p} - 1\right)\frac{\zeta^2}{M} + \frac{\sigma^2}{\alpha_p M}\right).
	\end{align*}
	Note that $B\le \frac{1}{1-\beta}$, so
	\begin{align*}
		\frac{\gamma^3L^2\beta^2}{2(1 - \beta)^3}\ec{\|v^{k-1}\|^2} 
		&\le \frac{\gamma^3 L^2\beta^2}{2(1 - \beta)^5} \frac{\sigma^2}{\alpha_p M} + \frac{\gamma^3 L^2\beta^2}{2(1 - \beta)^5} \frac{(1-\alpha_p)\zeta^2}{\alpha_p M}\\
		&\quad + \omega\frac{\gamma^3 L^2\beta^2}{2(1 - \beta)^4}\sum_{l=0}^{k-1}\beta^{k-1-l}\ec{\|\nabla f(x^{l})\|^2}
	\end{align*}
	with $\omega\eqdef \frac{M-1}{M} + \frac{1}{M\alpha_p}$.
	We, thus, obtain
	\begin{align*}
		\mathbb{E} [f(z^{k+1})]
		 &\le \mathbb{E} [f(z^k)] - \frac{\gamma}{2(1-\beta)}\left(1 - \frac{L \gamma\omega}{1-\beta} \right)\ec{\|\nabla f(x^k)\|^2}\\
		 &\quad + \frac{L\gamma^2 \sigma^2}{2M\alpha_p(1-\beta)^2} + \frac{\gamma^3 L^2\beta^2\sigma^2}{2(1 - \beta)^5\alpha_p M}\\
		&\quad + \frac{\gamma^3 L^2\beta^2(1-\alpha_p)\zeta^2}{2(1 - \beta)^5\alpha_p M} + \omega\frac{\gamma^3 L^2\beta^2}{2(1 - \beta)^4} \sum_{l=0}^{k-1}\beta^{k-1-l}\ec{\|\nabla f(x^{l})\|^2}.
	\end{align*}
	Telescoping this inequality from 0 to $K-1$, we get
	\begin{align*}
		\mathbb{E}[ f(z^K) - f(z^0)]
		&\le K\left( \frac{L\gamma^2 \sigma^2}{2\alpha_p M(1-\beta)^2} + \frac{\gamma^3 L^2\beta^2\sigma^2}{2(1 - \beta)^5\alpha_p M} + \frac{\gamma^3 L^2\beta^2(1-\alpha_p)\zeta^2}{2(1 - \beta)^5\alpha_p M} \right)\\
		& + \frac{\gamma}{2}\sum_{l=0}^{K-2}\left(\omega\frac{\gamma^2 L^2\beta^2}{(1 - \beta)^4}\sum_{k'=l+1}^{k-1}\beta^{k'-1-l} + \frac{L\gamma\omega}{(1-\beta)^2} - \frac{1}{1-\beta}\right)\ec{\|\nabla f(x^{l})\|^2} \\
		& + \frac{\gamma}{2}\left(\frac{L\gamma\omega}{(1-\beta)^2} - \frac{1}{1-\beta}\right)\ec{\|\nabla f(x^{K-1})\|^2}\\
		&\le K\left( \frac{L\gamma^2 \sigma^2}{2\alpha_p M(1-\beta)^2} + \frac{\gamma^3 L^2\beta^2\sigma^2}{2(1 - \beta)^5\alpha_p M} + \frac{\gamma^3 L^2\beta^2(1-\alpha_p)\zeta^2}{2(1 - \beta)^5\alpha_p M} \right)\\
		& + \frac{\gamma}{2}\sum_{l=0}^{K-1}\left(\omega\frac{\gamma^2 L^2\beta^2}{(1 - \beta)^5} + \frac{L\gamma\omega}{(1-\beta)^2} - \frac{1}{1-\beta}\right)\ec{\|\nabla f(x^{l})\|^2 }.
	\end{align*}
	It holds $f^\star\le f(z^K)$ and our assumption on $\beta$ implies that $\omega\frac{\gamma^2 L^2\beta^2}{(1 - \beta)^5} + \frac{L\gamma\omega}{(1-\beta)^2} - \frac{1}{1-\beta} \le -\frac{1}{2}$, so it all results in
	\begin{align*}
		\frac{1}{K}\sum_{l=0}^{K-1} \ec{\|\nabla f(x^{l})\|^2} 
		&\le  \frac{4(f(z^0) - f^\star)}{\gamma K} + 2\gamma\frac{L \sigma^2}{\alpha_p  M(1-\beta)^2} + 2\gamma^2\frac{ L^2\beta^2\sigma^2}{(1 - \beta)^5\alpha_p M}\\
		&\quad + 2\gamma^2\frac{L^2\beta^2(1-\alpha_p)\zeta^2}{2(1 - \beta)^5\alpha_p M}.
	\end{align*}
	Since $\overline x^K$ is sampled uniformly from $\{x^0, \dotsc, x^{K-1}\}$, the left-hand side is equal to $\ec{ \|\nabla f(\overline x^K)\|^2}$. Finally, note that $z^0=x^0$ and $f(z^0)=f(x^0)$.
\end{proof}
\begin{corollary}\label{cor:TernGrad-momentum}
		If we set $\gamma=\frac{1-\beta^2}{2\sqrt{K}L\omega}$, where $\omega = \frac{M-1}{M} + \frac{1}{M\alpha_p}$, and $\beta$ such that $\frac{\beta^2}{(1 - \beta)^3}\le 4K\omega$ with $K>1$, then the error after $K$ iterations is at most 
		\[
			\frac{1}{\sqrt{K}}\left(\frac{8L\omega(f(x^0)-f^\star)}{1-\beta^2} + \frac{(1+\beta)\sigma^2}{\omega\alpha_p M(1-\beta)}\right) + \frac{1}{K}\frac{(1+\beta)^4\beta^2\sigma^2}{2(1 - \beta)\omega\alpha_pM} + \frac{1}{K}\frac{(1+\beta)^4\beta^2(1-\alpha_p)\zeta^2}{2(1 - \beta)\omega\alpha_pM}.
		\]
\end{corollary}
\begin{proof}
	Our choice of $\gamma = \frac{1-\beta^2}{2\sqrt{K}L\omega}$ implies that
	\[
		\frac{\beta^2}{(1 - \beta)^3}\le \frac{1 - \beta^2 - 2L\gamma\omega}{\gamma^2 L^2\omega} \Longleftrightarrow \frac{\beta^2}{(1 - \beta)^3}\le 4k\omega.
	\]
	Now it remains to plug-in $\gamma = \frac{1-\beta^2}{2\sqrt{k}L\omega}$ in the expression $\frac{4(f(z^0) - f^\star)}{\gamma K} + 2\gamma\frac{L \sigma^2}{\alpha_p  M(1-\beta)^2} + 2\gamma^2\frac{ L^2\beta^2\sigma^2}{(1 - \beta)^5\alpha_p M} + \frac{1}{K}\frac{(1+\beta)^4\beta^2(1-\alpha_p)\zeta^2}{2(1 - \beta)\omega\alpha_pM}$ to get the desired result.
\end{proof}

\subsection{Strongly convex analysis}\label{sec:TernGrad-strongly-convex}
\begin{theorem}\label{thm:terngrad_strg_cvx_prox}
	Assume that each function $f_i$ is $\mu$-strongly convex and $L$-smooth. Choose stepsizes $\gamma_k = \gamma > 0$ satisfying
\begin{equation}\label{eq:gamma_cond_terngrad}\gamma \le \frac{2M\alpha_p}{(\mu + L)(2+(M-2)\alpha_p)}.\end{equation}
If we run Algorithm $\ref{alg:terngrad}$ for $K$ iterations with $\gamma_k=\gamma$, then
	\[
    	\mathbb{E}\left[\|x^K - x^\star\|^2\right] 
    	\le (1 - \gamma\mu)^K \|x^0-x^\star\|^2 +  \frac{\gamma}{\mu}\left(\frac{\sigma^2}{M\alpha_p} + \frac{2(1-\alpha_p)}{M^2\alpha_p}\sumiM\|h_i^\star\|^2\right),
    \]
    where $\sigma^2\eqdef \frac{1}{M}\sumiM\sigma_i^2$ and $h_i^\star = \nabla f_i(x^\star)$.
\end{theorem}
\begin{proof}
	In the similar way as we did in the proof of Theorem~\ref{thm:DIANA-strongly_convex} one can derive inequality~\eqref{eq:buf89gh38bf98} for the iterates of \algname{TernGrad}:
	\begin{align*}
       \ec{ \|x^{k+1} - x^\star\|^2} 
        &\le  
        \ec{ \|x^k - x^\star\|^2} - 2\gamma \ec{ \< \nabla f(x^k) - h^\star, x^k - x^\star>} \notag \\
        & \qquad + \frac{\gamma^2}{M} \sumiM \ec{ \|\nabla f_i(x^k) - h_i^\star\|^2} + \frac{\gamma^2}{M^2}\sum_{i=1}^M\ec {\Psi(g_i^k)} + \frac{\gamma^2 \sigma^2}{M}.
	\end{align*}
	By definition $\alpha_p(d_l) = \inf\limits_{x\neq 0,x\in\R^{d_l}}\frac{\|x\|^2}{\|x\|_1\|x\|_p} = \left(\sup\limits_{x\neq 0,x\in\R^{d_l}}\frac{\|x\|_1\|x\|_p}{\|x\|^2}\right)^{-1}$ and $\alpha_p = \alpha_p(\max\limits_{l=1,\ldots,m}d_l)$ which implies
	\begin{align*}
		\mathbb{E}\left[\Psi_l(g_i^k) \right] 
		&= \mathbb{E}\left[\|g_i^k(l)\|_1\|g_i^k(l)\|_p - \|g_i^k(l)\|^2 \right] 
		= \mathbb{E}\left[\|g_i^k(l)\|^2\left(\frac{\|g_i^k(l)\|_1\|g_i^k(l)\|_p}{\|g_i^k(l)\|^2} - 1\right) \right]\\
		&\le \left(\frac{1}{\alpha_p(d_l)}-1\right)\ec{\|g_i^k(l)\|^2} \le \left(\frac{1}{\alpha_p}-1\right)\ec{\|g_i^k(l)\|^2}.
	\end{align*}
	 Moreover, 
	\[
		\|g_i^k\|^2 = \sum_{l=1}^m\|g_i^k(l)\|^2,\quad \Psi(g_i^k) = \sum\limits_{l=1}^m\Psi_l(g_i^k).	
	\]
	This helps us to get the following inequality
	\begin{align*}
       \ec{ \|x^{k+1} - x^\star\|^2} 
        &\le 
        \ec{ \|x^k - x^\star\|^2} - 2\gamma \ec{ \< \nabla f(x^k) - h^\star, x^k - x^\star>} + \frac{\gamma^2 \sigma^2}{M} \notag \\
        & \quad + \frac{\gamma^2}{M} \sumiM \ec{ \|\nabla f_i(x^k) - h_i^\star\|^2 }+ \frac{\gamma^2}{M^2}\left(\frac{1}{\alpha_p}-1\right)\sum_{i=1}^M \mathbb{E}\left[\|g_i^k\|^2\right].
    \end{align*}
	Using tower property of mathematical expectation and the fact that 
	\[
		\mathbb{E}\left[\|g_i^k \|^2\mid x^k\right] = \mathbb{E}\left[\|g_i^k - \nabla f_i(x^k)\|^2\mid x^k\right] + \|\nabla f_i(x^k)\|^2 \le \sigma_i^2 + \|\nabla f_i(x^k)\|^2,
	\] 
	we obtain
	\[
		\mathbb{E}[\|g_i^k\|^2]
		\le \mathbb{E}[\|\nabla f_i(x^k)\|^2] + \sigma_i^2 \le 2\mathbb{E}\left[\|\nabla f_i(x^k) - h_i^\star\|^2\right] + 2\|h_i^\star\|^2 + \sigma_i^2,	
	\]
	where the last inequality follows from the fact that for all $x,y\in\R^M$ the inequality $\|x+y\|^2 \le 2\left(\|x\|^2 + \|y\|^2\right)$ holds.
	Putting all together we have
	\begin{eqnarray*}
		\mathbb{E}\left[\|x^{k+1}-x^\star\|^2\right]
		&\le & \mathbb{E} \left[\|x^k - x^\star\|^2\right] - 2\gamma \mathbb{E} \left[\< \nabla f(x^k) - h^\star, x^k - x^\star>\right] \notag \\
        && \qquad + \frac{\gamma^2}{M}\left(1 + \frac{2(1-\alpha_p)}{M\alpha_p}\right) \sumiM \ec{ \|\nabla f_i(x^k) - h_i^\star\|^2 } \notag\\
        &&\qquad + \frac{2\gamma^2(1-\alpha_p)}{M^2\alpha_p}\sum_{i=1}^M \|h_i^\star\|^2 + \frac{\gamma^2 \sigma^2}{M\alpha_p}.
	\end{eqnarray*}
	Using the splitting trick \eqref{eq:inner_product_splitting} we get
	\begin{eqnarray}
		\ec{\|x^{k+1}-x^\star\|^2} &\le & \left(1-\frac{2\gamma \mu L}{\mu+L}\right)\ec{\|x^k - x^\star\|^2}\notag\\
		&&\quad + \frac{1}{M}\left(\gamma^2\left(1 + \frac{2(1-\alpha_p)}{M\alpha_p}\right) - \frac{2\gamma}{\mu+L}\right)\sumiM \ec{\|\nabla f_i(x^k)-h_i^\star\|^2 } \notag\\
		&& \quad + \frac{2\gamma^2(1-\alpha_p)}{M^2\alpha_p}\sum_{i=1}^M \|h_i^\star\|^2 + \frac{\gamma^2 \sigma^2}{M\alpha_p}\label{eq:trngrd_str_cvx_pre_final}.
	\end{eqnarray}
	Since $\gamma \le \frac{2M\alpha_p}{(\mu+L)(2+(M-2)\alpha_p)}$ the term $\left(\gamma^2\left(1 + \frac{2(1-\alpha_p)}{M\alpha_p}\right) - \frac{2\gamma}{\mu+L}\right)$ is nonnegative. Moreover, since $f_i$ is $\mu$--strongly convex, we have
$\mu \|x^k-x^\star\|^2 \leq \langle \nabla f_i(x^k) - h_i^\star, x^k -x^\star \rangle$. Applying the Cauchy-Schwarz inequality to further bound the right-hand side, we get the inequality $\mu \|x^k-x^\star\| \leq \|\nabla f_i(x^k) - h_i^\star\|$. Using these observations, we can get rid of the second term in the \eqref{eq:trngrd_str_cvx_pre_final} and absorb it with the first term, obtaining
	\begin{eqnarray*}
		\ec{\|x^{k+1}-x^\star\|^2} &\le & \left(1 - 2\gamma\mu + \gamma^2\mu^2\left(1 + \frac{2(1-\alpha_p)}{M\alpha_p}\right)\right)\ec{\|x^k - x^\star\|^2}\\
		&&\qquad + \frac{2\gamma^2(1-\alpha_p)}{M^2\alpha_p}\sum_{i=1}^M \|h_i^\star\|^2 + \frac{\gamma^2 \sigma^2}{M\alpha_p}\\
		&\overset{\eqref{eq:conseq_gamma_choice_terngrad_prox}}{\le} & (1-\gamma\mu)\ec{\|x^k - x^\star\|^2} + \gamma^2\left(\frac{\sigma^2}{M\alpha_p} + \frac{2(1-\alpha_p)}{M^2\alpha_p}\sumiM\|h_i^\star\|^2\right).
	\end{eqnarray*}
Finally, unrolling the recurrence leads to
    \begin{align*}
    		\ec{\|x^k-x^\star\|^2 }
    		&\le (1 - \gamma\mu)^k \|x^0-x^\star\|^2 + \sum\limits_{l=0}^{k-1}(1-\gamma\mu)^l\gamma^2\left(\frac{\sigma^2}{M\alpha_p} + \frac{2(1-\alpha_p)}{M^2\alpha_p}\sumiM\|h_i^\star\|^2\right) \\
    		&\le (1 - \gamma\mu)^k \|x^0-x^\star\|^2 + \sum\limits_{l=0}^{\infty}(1-\gamma\mu)^l\gamma^2\left(\frac{\sigma^2}{M\alpha_p} + \frac{2(1-\alpha_p)}{M^2\alpha_p}\sumiM\|h_i^\star\|^2\right) \\
    		&= (1 - \gamma\mu)^k \|x^0-x^\star\|^2 + \frac{\gamma}{\mu}\left(\frac{\sigma^2}{M\alpha_p} + \frac{2(1-\alpha_p)}{M^2\alpha_p}\sumiM\|h_i^\star\|^2\right).
    \end{align*}
\end{proof}

\subsection{Decreasing stepsize}\label{sec:TernGrad-decreasing-stepsizes}

\begin{theorem}\label{thm:TernGrad-decreasing-stepsizes}
    Assume that $f$ is $L$-smooth, $\mu$-strongly convex and we have access to its gradients with bounded noise. Set $\gamma_k = \frac{2}{\mu k + \theta}$ with some $\theta \ge \frac{(\mu+L)(2+(M-2)\alpha_p)}{2M\alpha_p}$. After $K$ iterations of Algorithm~\ref{alg:terngrad} we have
    \begin{align*}
        \ec{ \|x^K - x^\star\|^2}
        \le \frac{1}{\eta K+1}\max\left\{ \|x^0-x^\star\|^2, \frac{4}{\mu\theta}\left(\frac{\sigma^2}{M\alpha_p} + \frac{2(1-\alpha_p)}{M^2\alpha_p}\sumiM\|h_i^\star\|^2\right) \right\},
    \end{align*}
    where $\eta\eqdef \frac{\mu}{\theta}$, $\sigma^2 = \frac{1}{M}\sumiM\sigma_i^2$ and $h_i^\star = \nabla f_i(x^\star)$.
\end{theorem}
\begin{proof}
	To get a recurrence, let us recall an upper bound we have proved before in Theorem~\ref{thm:terngrad_strg_cvx_prox}:
    \[
        \ec{\|x^{k+1}-x^\star\|^2}
        \le (1 - \gamma_k\mu)\ec{\|x^k-x^\star\|^2 }+ \gamma_k^2\left(\frac{\sigma^2}{M\alpha_p} + \frac{2(1-\alpha_p)}{M^2\alpha_p}\sumiM\|h_i^\star\|^2\right).
    \]
    Having that, we can apply Lemma~\ref{lem:sgd_recurrence} to the sequence $\ec{\|x^k-x^\star\|^2 }$. The constants for the Lemma are: $\nu = \left(\frac{\sigma^2}{M\alpha_p} + \frac{2(1-\alpha_p)}{M^2\alpha_p}\sumiM\|h_i^\star\|^2\right)$ and 
    \[
    		C=\max\left\{ \|x^0-x^\star\|^2, \frac{4}{\mu\theta}\left(\frac{\sigma^2}{M\alpha_p} + \frac{2(1-\alpha_p)}{M^2\alpha_p}\sumiM\|h_i^\star\|^2\right) \right\}.
    \]
\end{proof}
\begin{corollary}\label{cor:TernGrad-decreasing-stepsizes}
	If we choose $\theta=\frac{(\mu+L)(2+(M-2)\alpha_p)}{2M\alpha_p}$, then to achieve $\ec{ \|x^K-x^\star\|^2 }\le \varepsilon$ we need at most 
	\[
		\cO\left( \frac{\kappa(1+M\alpha_p)}{M\alpha_p}\max\left\{ \|x^0-x^\star\|^2, \frac{M\alpha_p}{(1+M\alpha_p)\mu L}\left(\frac{\sigma^2}{M\alpha_p} + \frac{1-\alpha_p}{M^2\alpha_p}\sumiM\|h_i^\star\|^2\right) \right\}\frac{1}{\varepsilon} \right)
	\]
	iterations, where $\kappa \eqdef \frac{L}{\mu}$ is the condition number of $f$.
\end{corollary}
\begin{proof}
	If $\theta=\frac{(\mu+L)(2+(M-2)\alpha_p)}{2M\alpha_p} = \Theta\left(\frac{L(1+M\alpha_p)}{M\alpha_p}\right)$, then $\eta = \Theta\left(\frac{M\alpha_p}{\kappa(1+M\alpha_p)}\right)$ and $\frac{1}{\mu\theta} = \Theta\left(\frac{M\alpha_p}{\mu L(1+M\alpha_p)}\right)$. Putting all together and using the bound from Theorem~\ref{thm:TernGrad-decreasing-stepsizes} we get the desired result.
\end{proof}


\label{sec:A:detailsOfNumericalExperiments}

\subsection{Performance of \algname{DIANA}, \algname{QSGD} and \algname{TernGrad} on the Rosenbrock function}

In Figure~\ref{fig:rosen} we illustrate the workings of \algname{DIANA}, \algname{QSGD} and \algname{TernGrad} with 2 workers on the 2-dimensional (non-convex) Rosenbrock function: \[f(x,y)=(x-1)^2 + 10(y - x^2)^2,\] decomposed into average of $f_1=(x + 16)^2 + 10(y - x^2)^2 + 16y$ and $f_2= (x - 18)^2 + 10(y - x^2)^2 - 16y + \mathrm{const}$.
Each worker  has access to its own piece of the Rosenbrock function with parameter $a=1$ and $b=10$. The gradients used are not stochastic, and we use 1-bit version of \algname{QSGD}, so it also coincides with \algname{QGD} in that situation. For all methods, its parameters were carefully tuned except for momentum and $\alpha$, which were simply set  to $0.9$ and $0.5$ correspondingly. We see that \algname{DIANA} vastly outperforms the competing methods.

\begin{figure}[h!]
  \centering
    \includegraphics[scale=0.5 ]{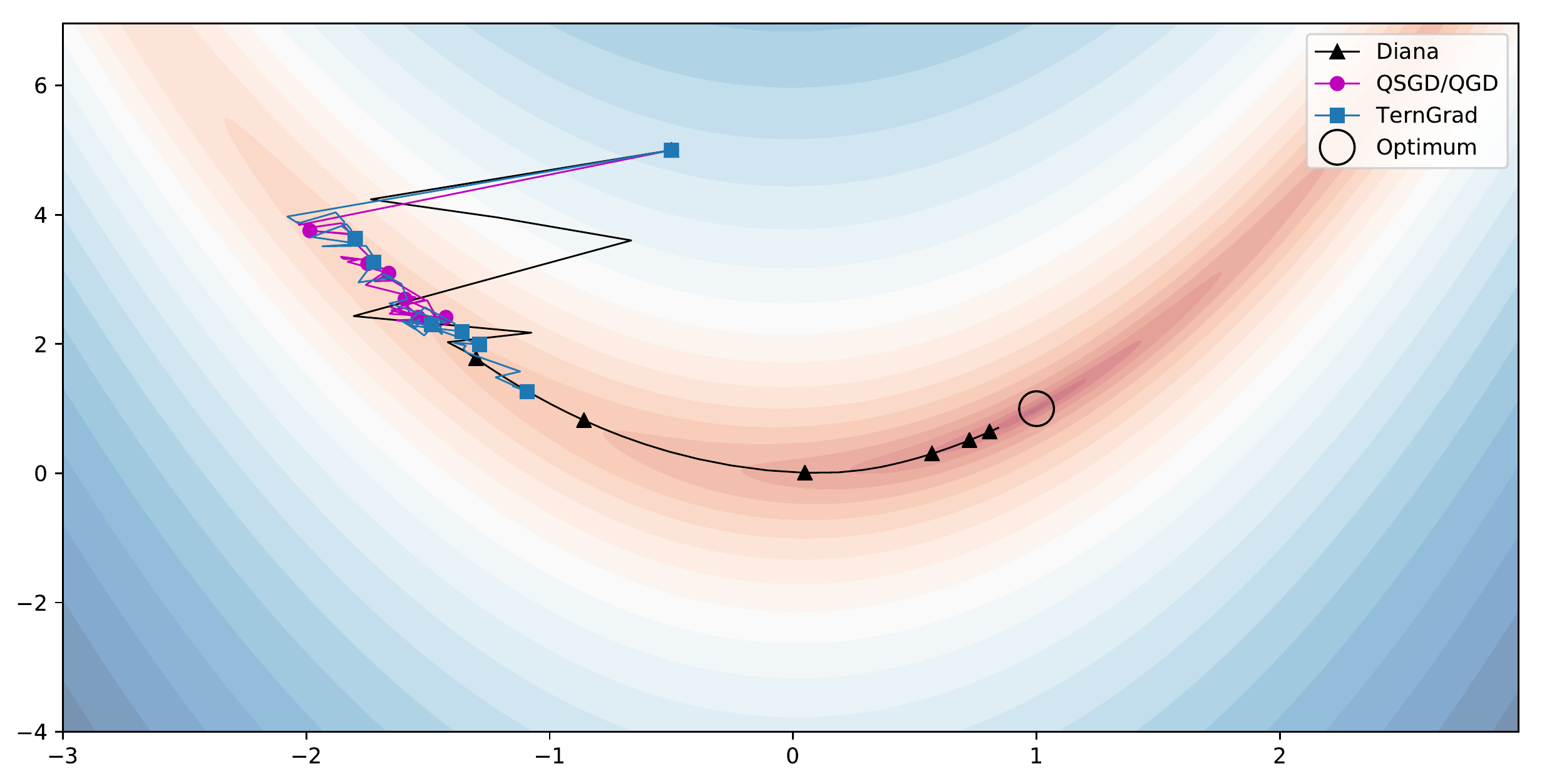}
  \caption{Illustration of the workings of \algname{DIANA}, \algname{QSGD} and \algname{TernGrad} on the Rosenbrock function.}
  \label{fig:rosen}
\end{figure}

\subsection{Logistic regression}\label{sec:log_reg}
We consider the logistic regression problem with $\ell_2$ and $\ell_1$ penalties for mushrooms dataset from LIBSVM. In our experiments we use $\ell_1$-penalty coefficient $l_1 = 2\cdot 10^{-3}$ and $\ell_2$-penalty coefficient $l_2 = \frac{L}{M}$. The coefficient $l_1$ is adjusted in order to have sparse enough solution ($\approx 20\%$ nonzero values). The main goal of this series of experiment is to compare the optimal parameters for $\ell_2$ and $\ell_\infty$ quantization.

\subsubsection{What $\alpha$ is better to choose?}
We run \algname{DIANA} with zero momentum ($\beta=0$) and obtain in our experiments that, actually, it is not important what $\alpha$ to choose for both $\ell_2$ and $\ell_\infty$ quantization. The only thing that we need to control is that $\alpha$ is small enough. 

\subsubsection{What is the optimal block-size?}
Since $\alpha$ is not so important, we run \algname{DIANA} with $\alpha = 10^{-3}$ and zero momentum ($\beta=0$) for different block sizes (see Figure~\ref{fig:block_tuning}). For the choice of $\ell_\infty$ quantization in our experiments it is always better to use full quantization. In the case of $\ell_2$ quantization it depends on the regularization: if the regularization is big then optimal block-size $\approx 25$ (dimension of the full vector of parameters is $d=112$), but if the regularization is small it is better to use small block sizes.
\begin{figure}[h!]
\centering

\includegraphics[scale=0.25]{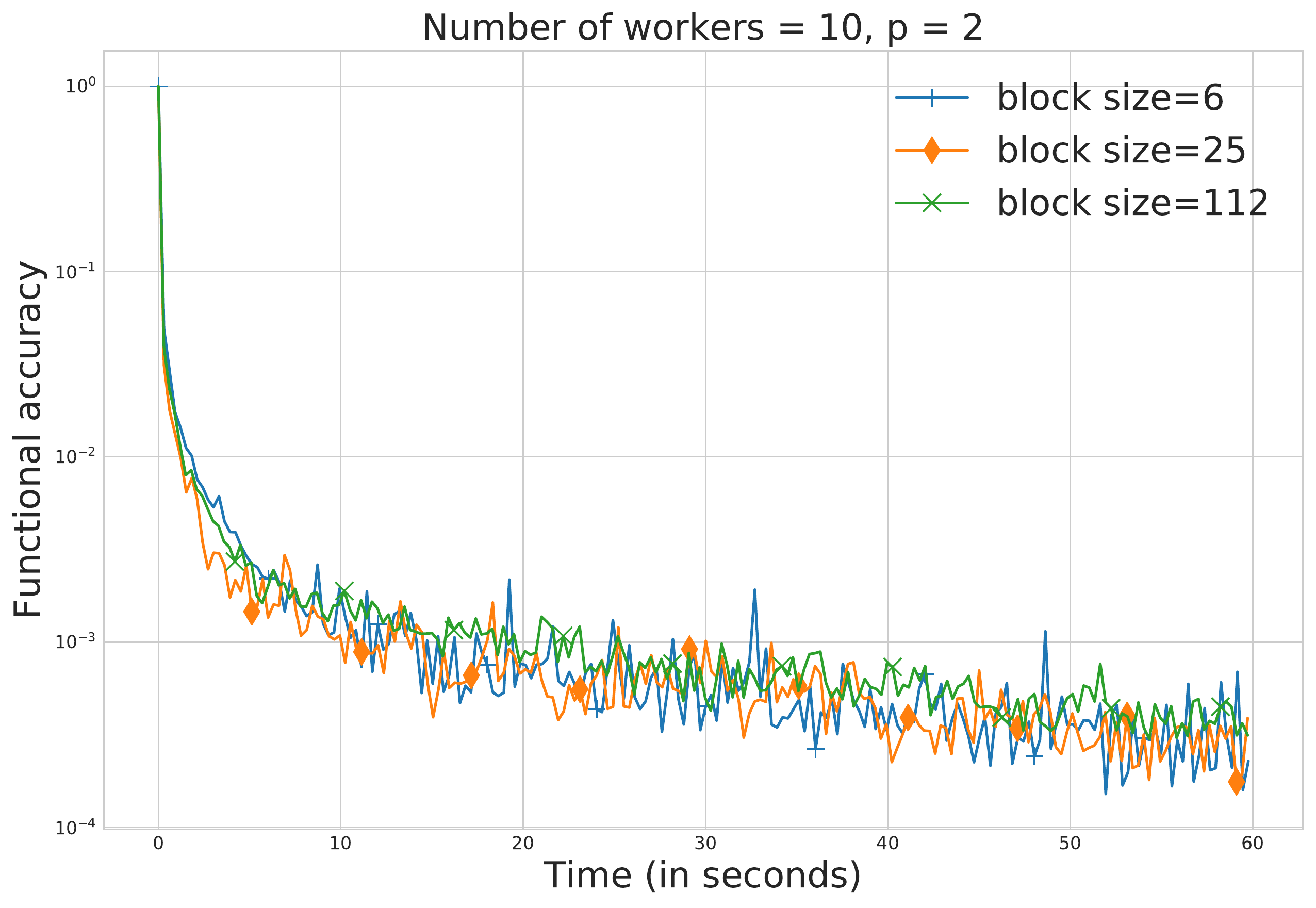}
\includegraphics[scale=0.25]{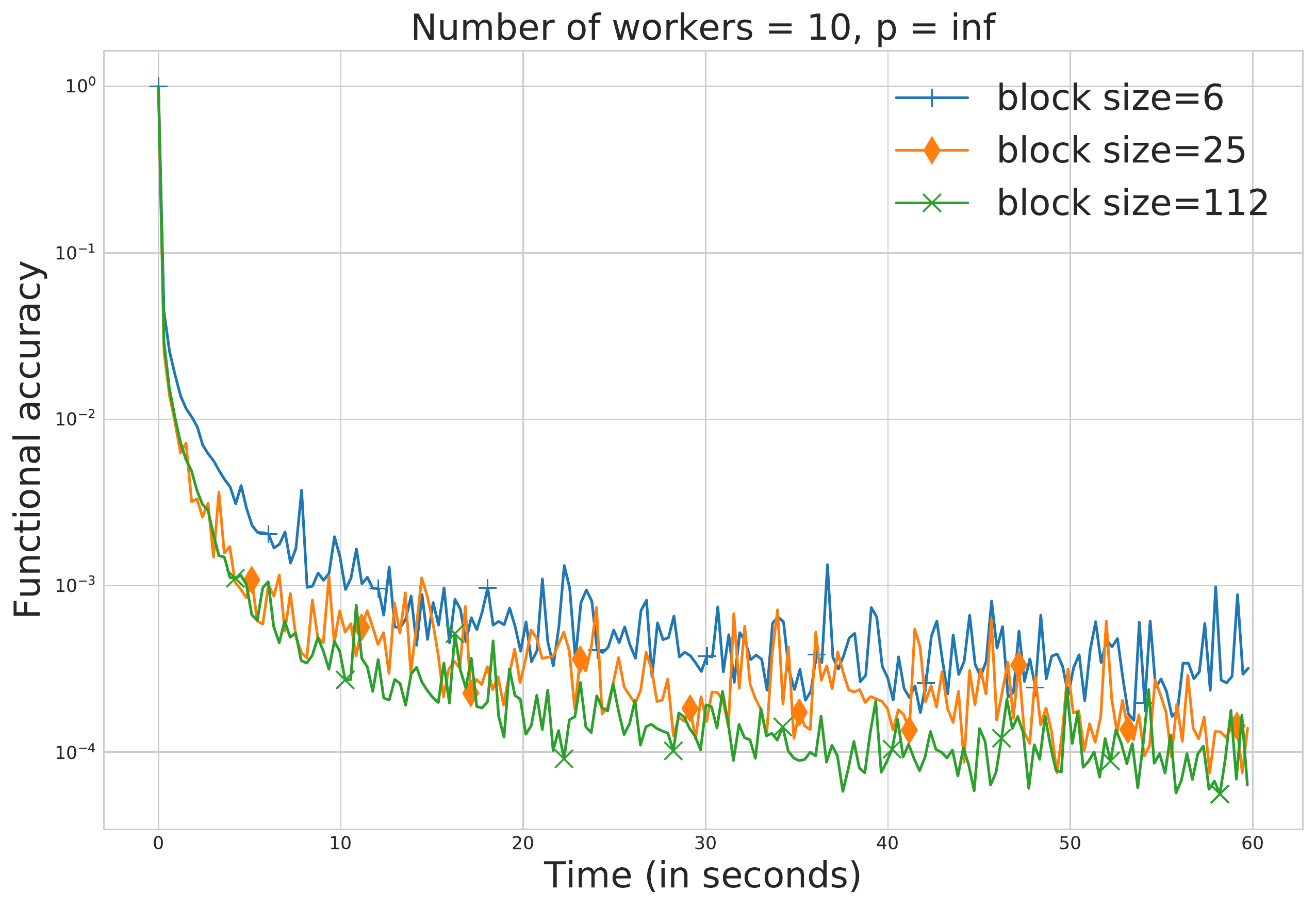}
\includegraphics[scale=0.25]{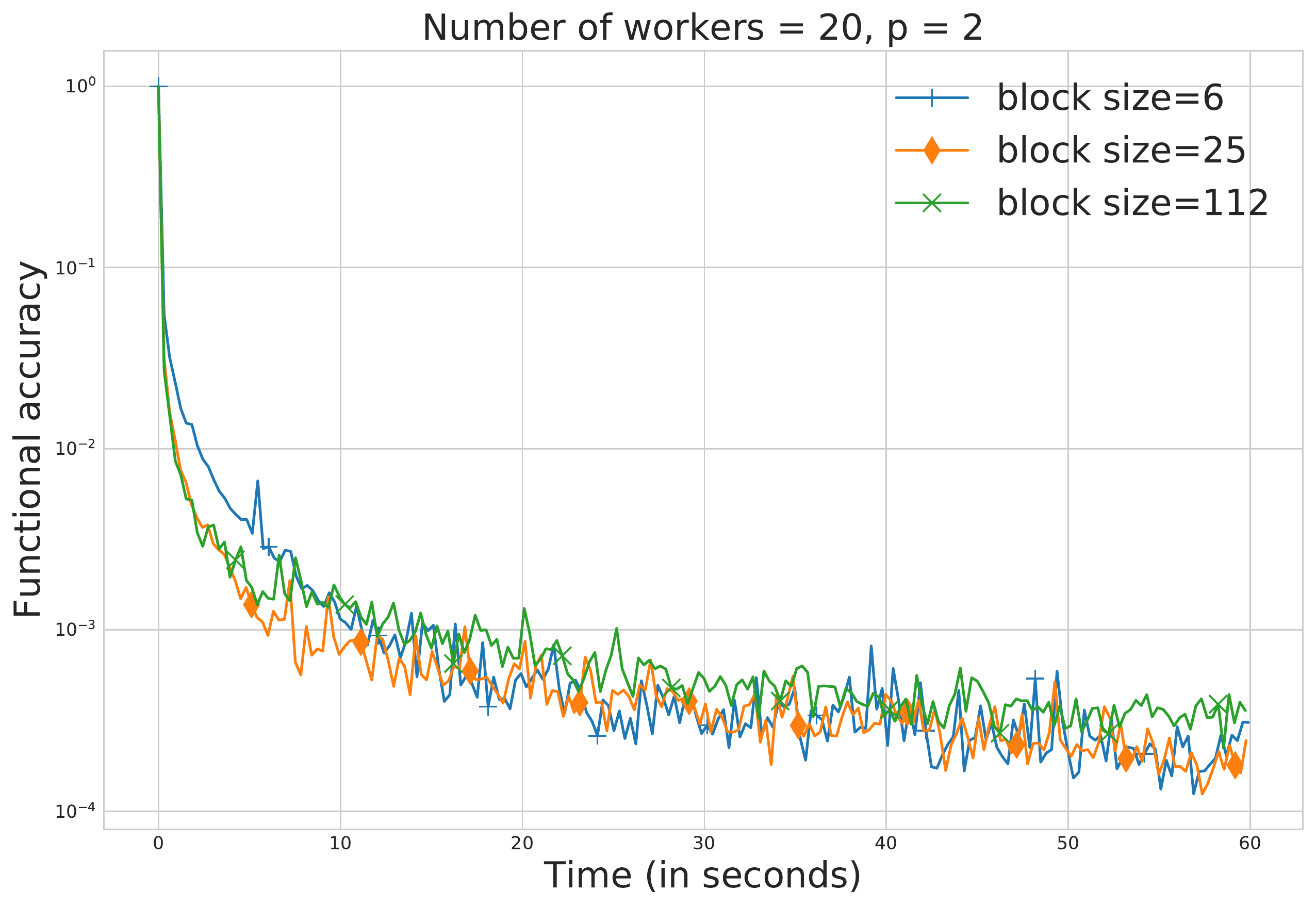}
\includegraphics[scale=0.25]{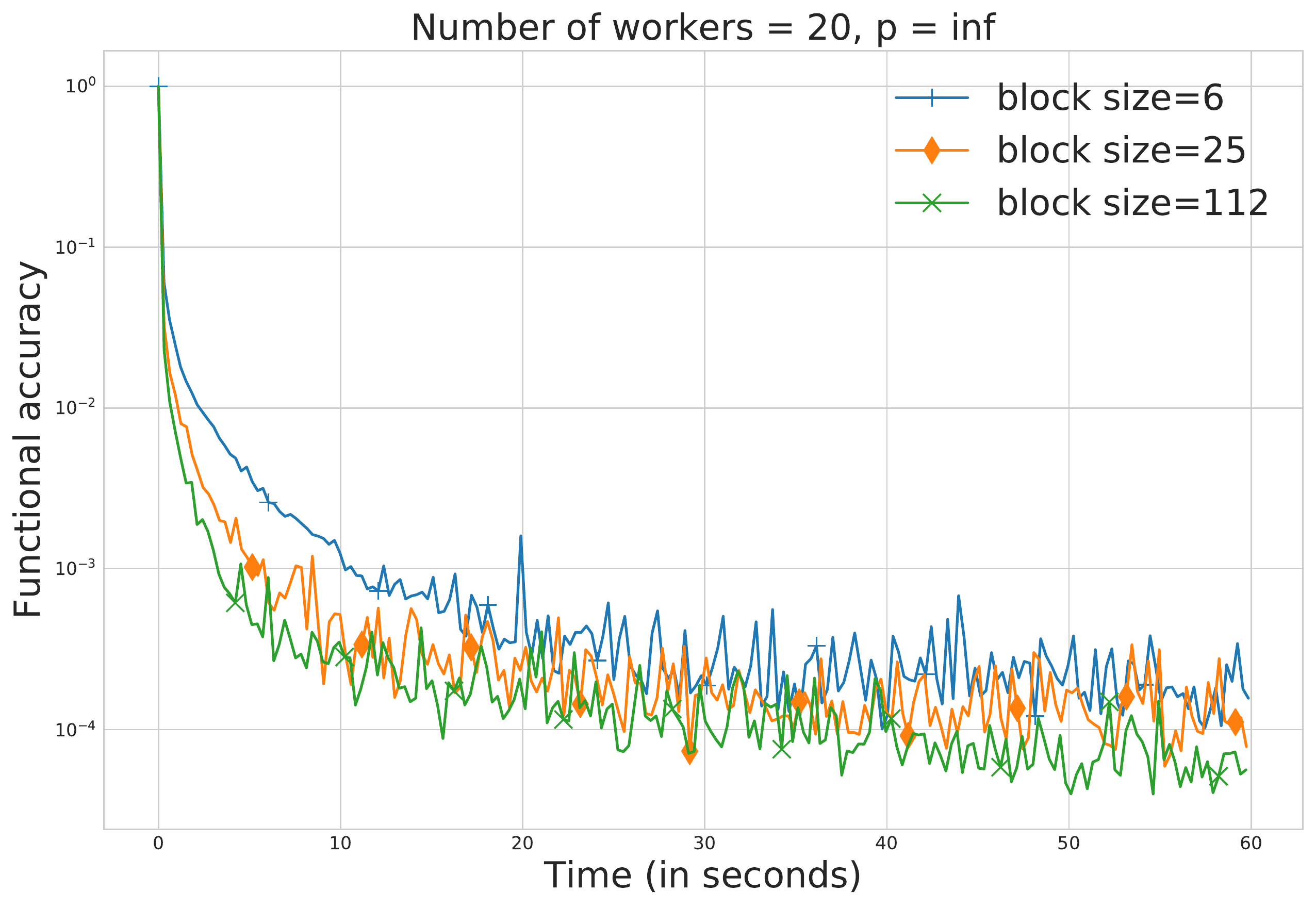}
\includegraphics[scale=0.25]{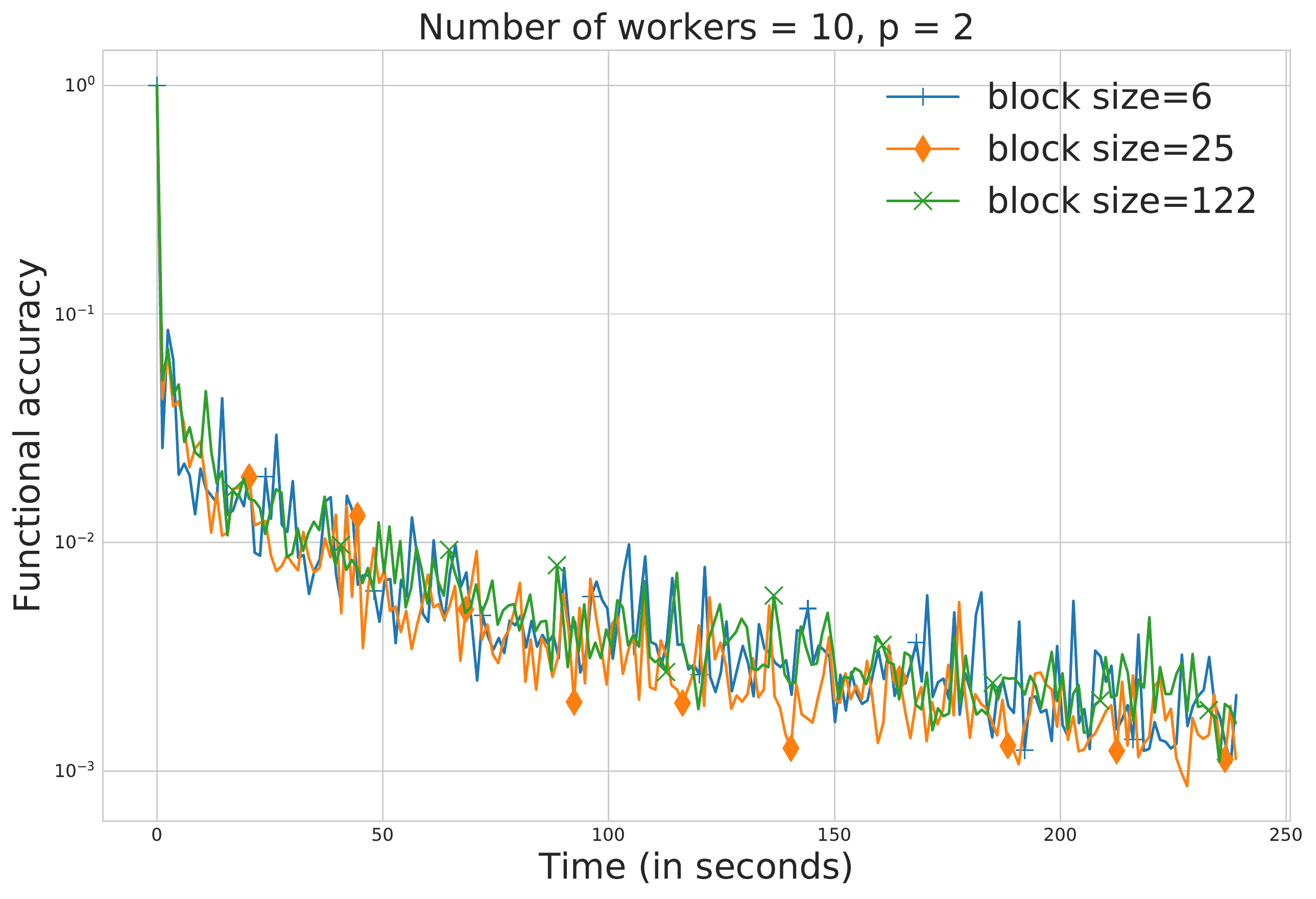}
\includegraphics[scale=0.25]{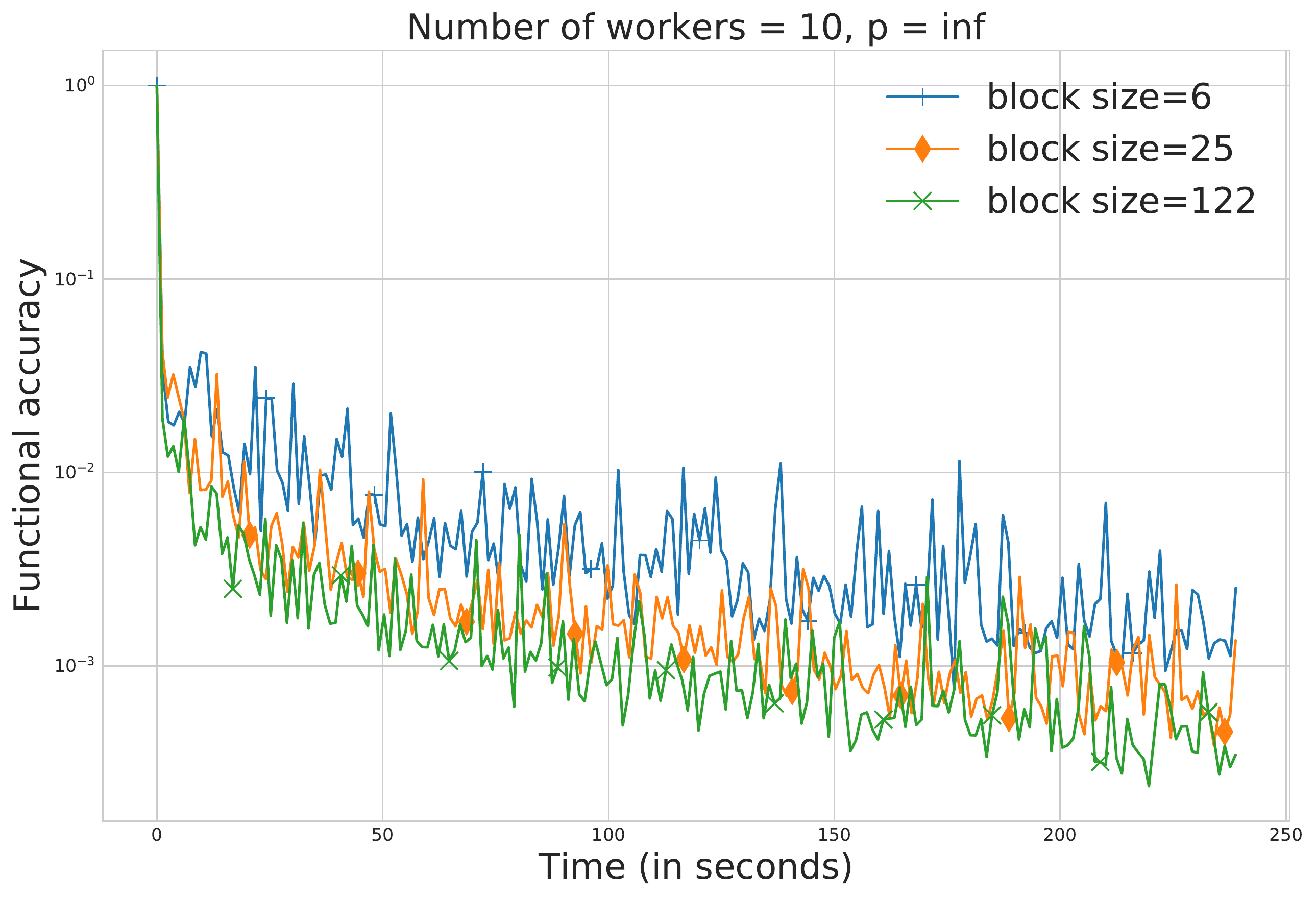}
\includegraphics[scale=0.25]{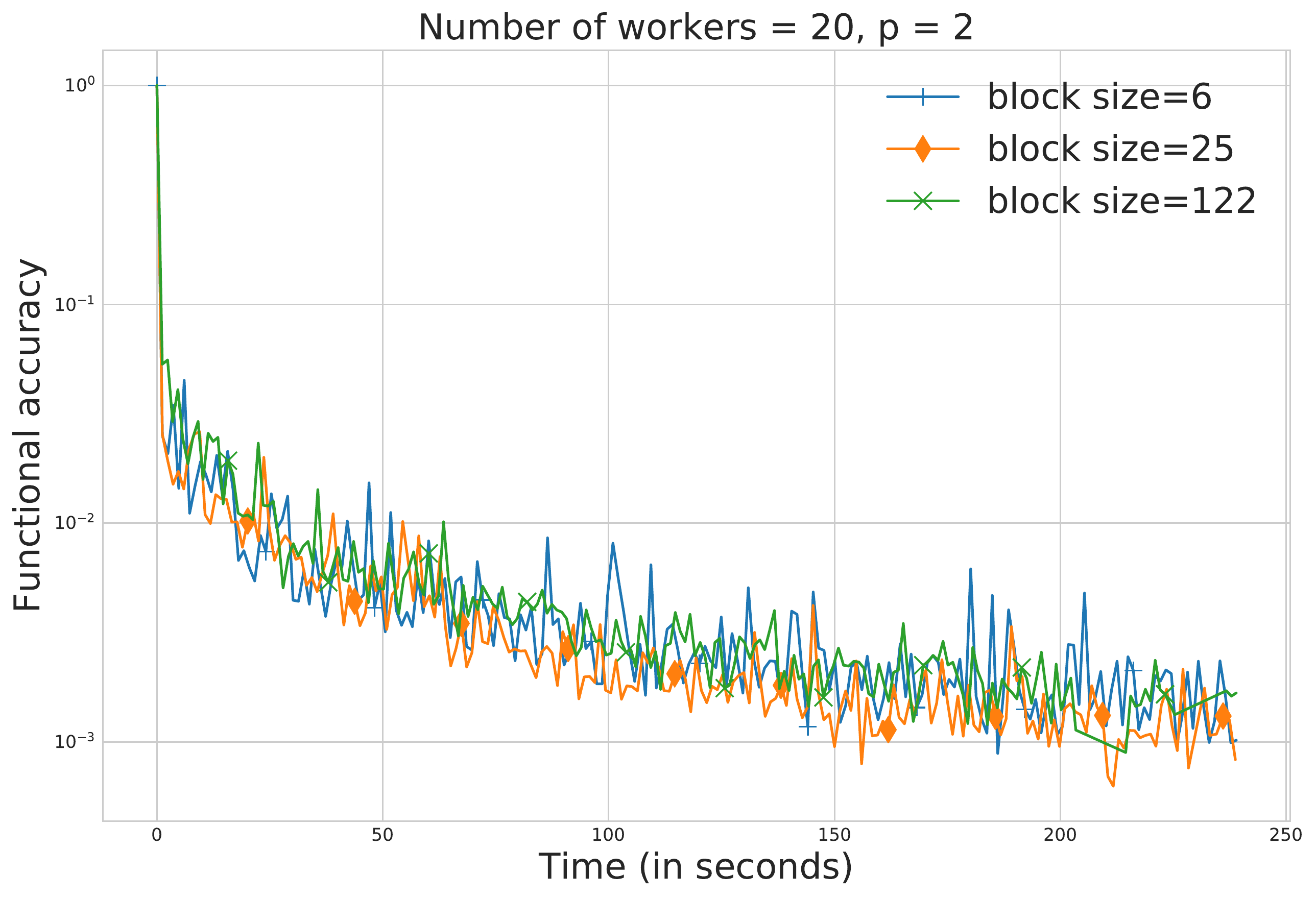}
\includegraphics[scale=0.25]{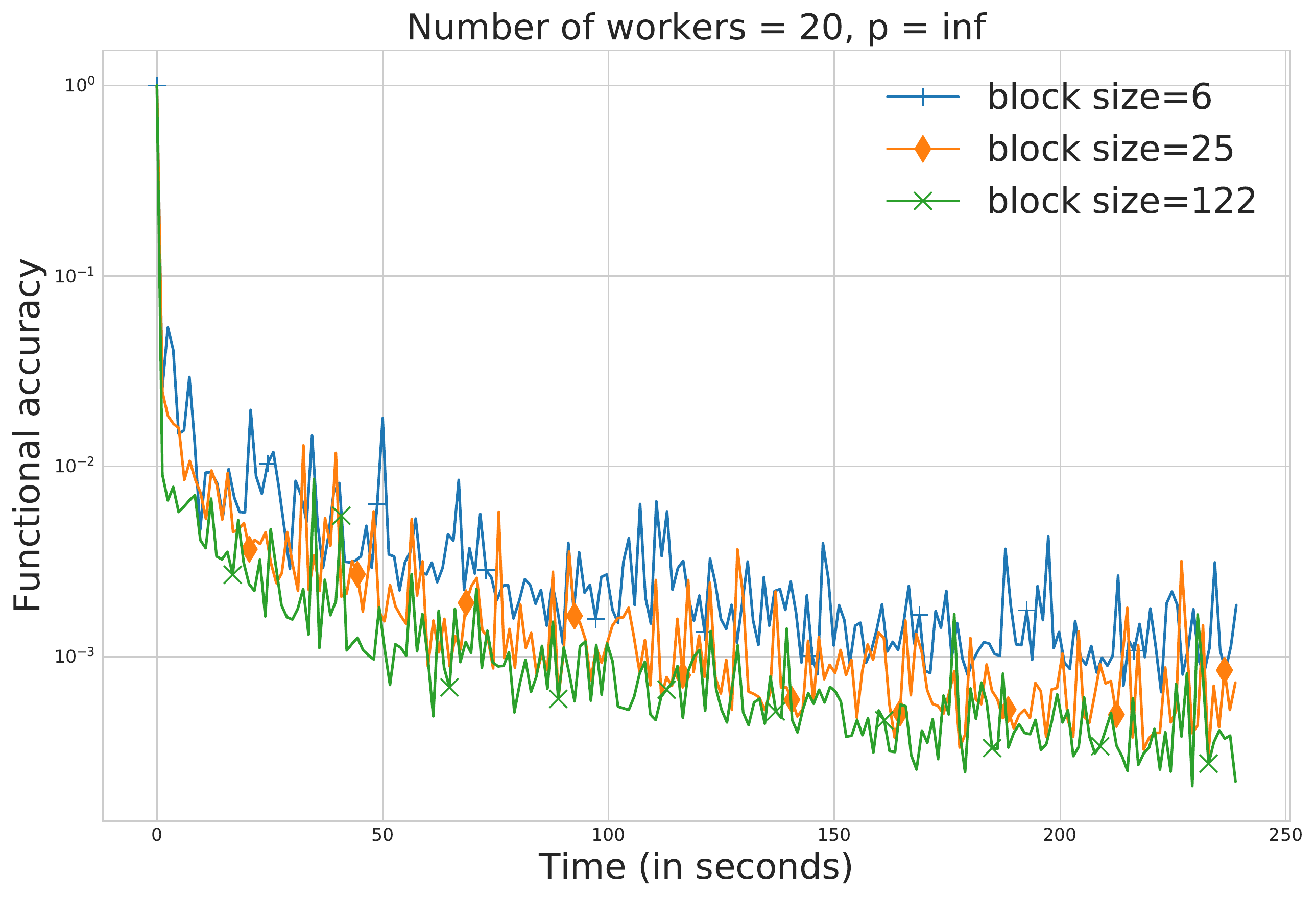}

\caption{Comparison of the influence of the block sizes on convergence for "mushrooms" (first row), "a5a" (second row)  datasets.} 
\label{fig:block_tuning}
\end{figure}

\begin{table}[ht!]
\begin{center}
\caption{Approximate optimal number of blocks for different dataset and configurations. Momentum equals zero for all experiments.
\label{tbl:opt_block}
}
\begin{tabular}{cccccc}
\toprule
Dataset & $N$ & $d$ & Number of workers & Quantization &  Optimal block size\\
\midrule
mushrooms & $8124$ & $112$ & $10$ & $\ell_2$ & $25$\\ 
mushrooms & $8124$ & $112$ & $10$ & $\ell_\infty$ & $112$\\ 
mushrooms & $8124$ & $112$ & $20$ & $\ell_2$  & $25$\\ 
mushrooms & $8124$ & $112$ & $20$ & $\ell_\infty$ & $112$\\ 
\hline 
a5a & $6414$ & $122$ & $10$ & $\ell_2$ & $25$\\ 
a5a & $6414$ & $122$ & $10$ & $\ell_\infty$ & $112$\\ 
a5a & $6414$ & $122$ & $20$ & $\ell_2$  & $25$\\ 
a5a & $6414$ & $122$ & $20$ & $\ell_\infty$ & $112$\\ 
\bottomrule

\end{tabular}
\end{center}
\end{table}

\subsubsection{\algname{DIANA} vs.\ \algname{QSGD} vs.\ \algname{TernGrad} vs.\ \algname{DQGD}}

%
%

\subsection{MPI: broadcast, reduce and gather}
\label{sec:A:MPI}

\begin{figure}[ht]
\centering
\includegraphics[scale=.5]{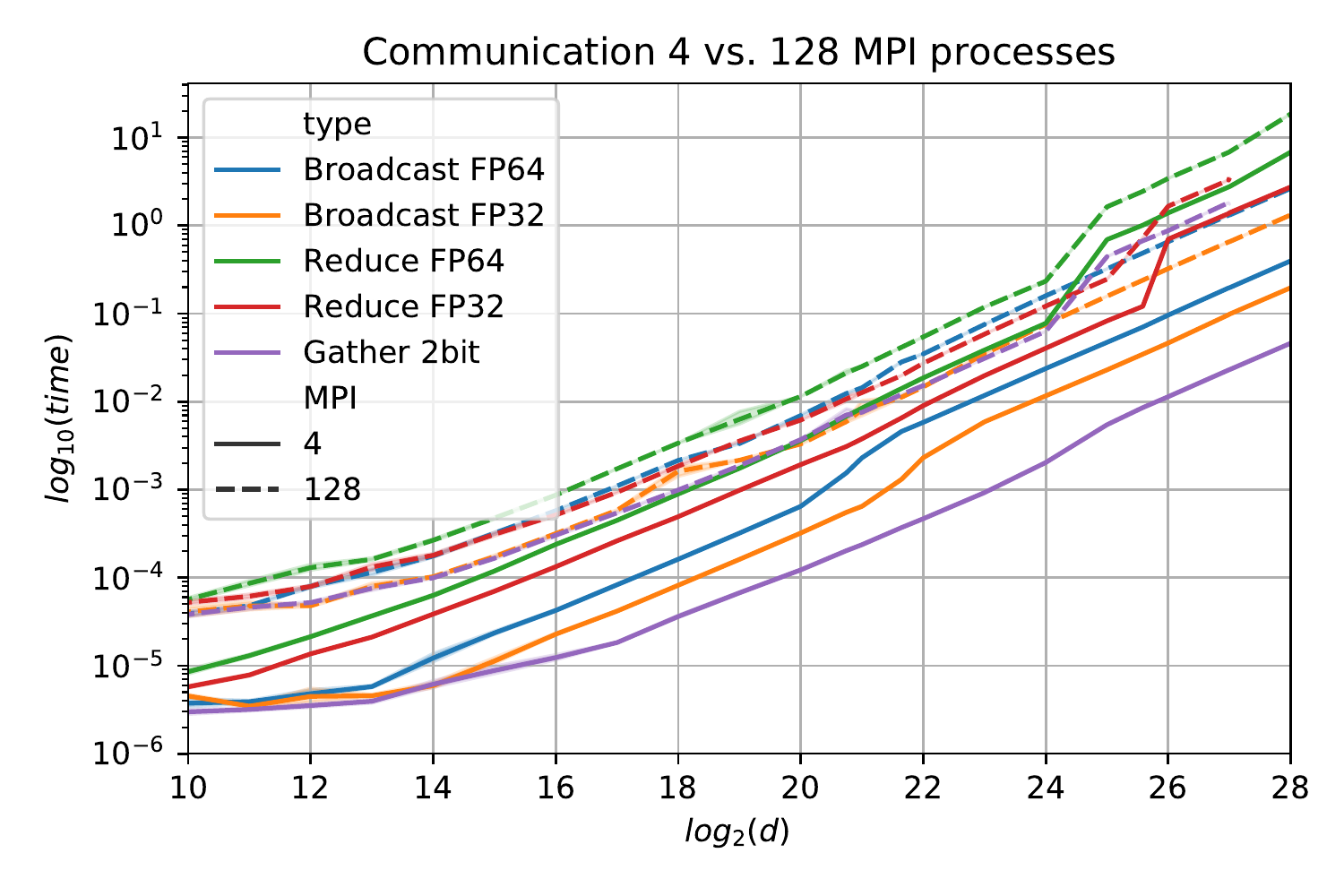}
\caption{Typical communication cost using  broadcast, reduce and gather  for 64 and 32 FP using 4 (solid) resp 128 (dashed) MPI processes.
See suppl.\ material for details about the network.
}
\label{fig:communication}
\end{figure}

In our experiments, we are running 4 MPI processes per physical node.
Nodes are connected by Cray Aries High Speed Network.

We utilize 3 MPI collective operations, Broadcast, Reduce and Gather. When implementing \algname{DIANA}, we could use P2P communication, but based on our experiments, we found that using Gather to collect data from workers significantly outperformed P2P communications.

In Figure~\ref{fig:communication:scale}
we show the duration of different communications for various MPI processes and message length.
Note that Gather 2bit do not scale linearly (as would be expected). It turns out, we are not the only one who observed such a weird behavior when using cray MPI implementation (see \cite{parker2018performance} for a nice study obtained by a team from Argonne National Laboratory).
\begin{figure}[h!]

\includegraphics[scale=0.46]{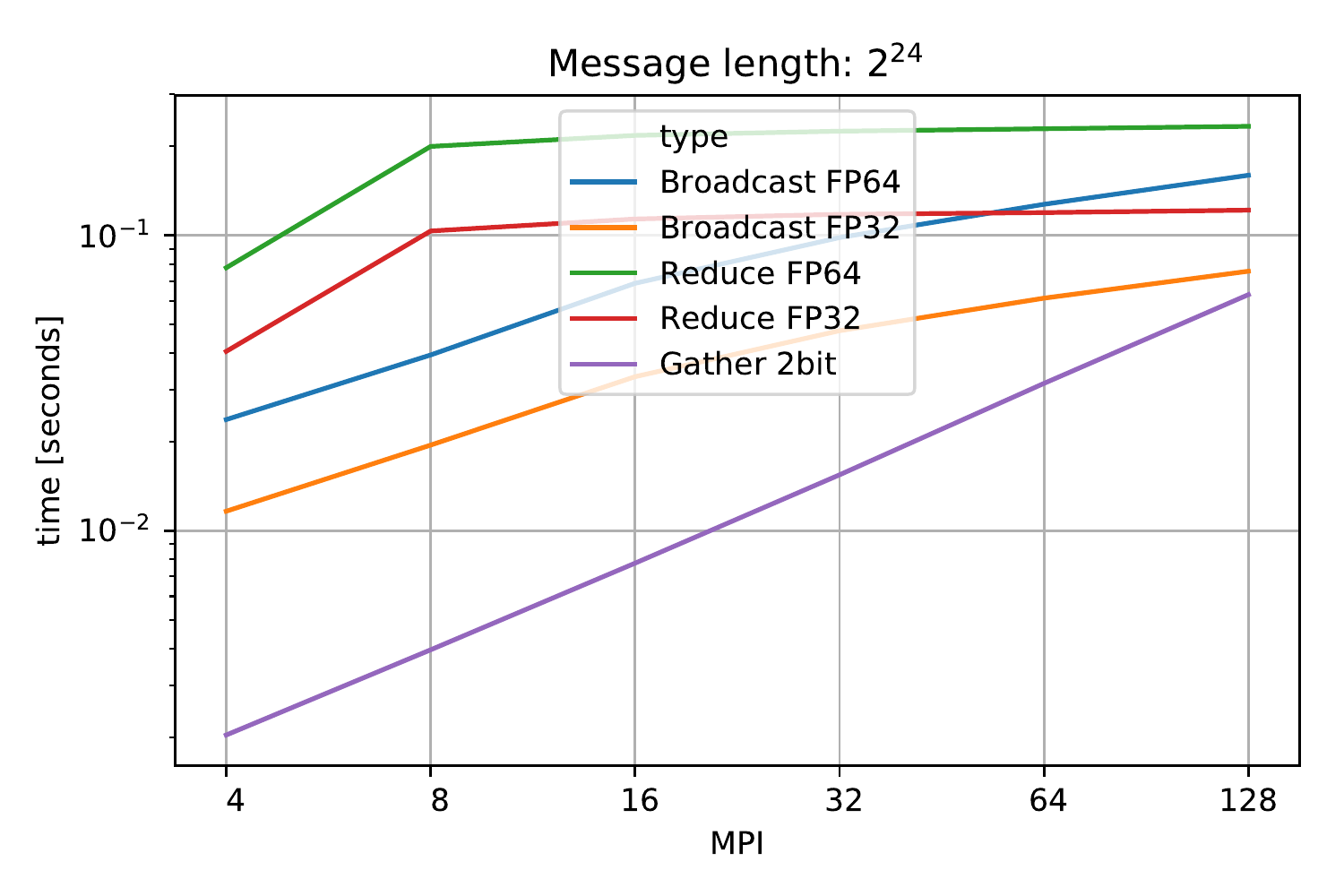}
\includegraphics[scale=0.46]{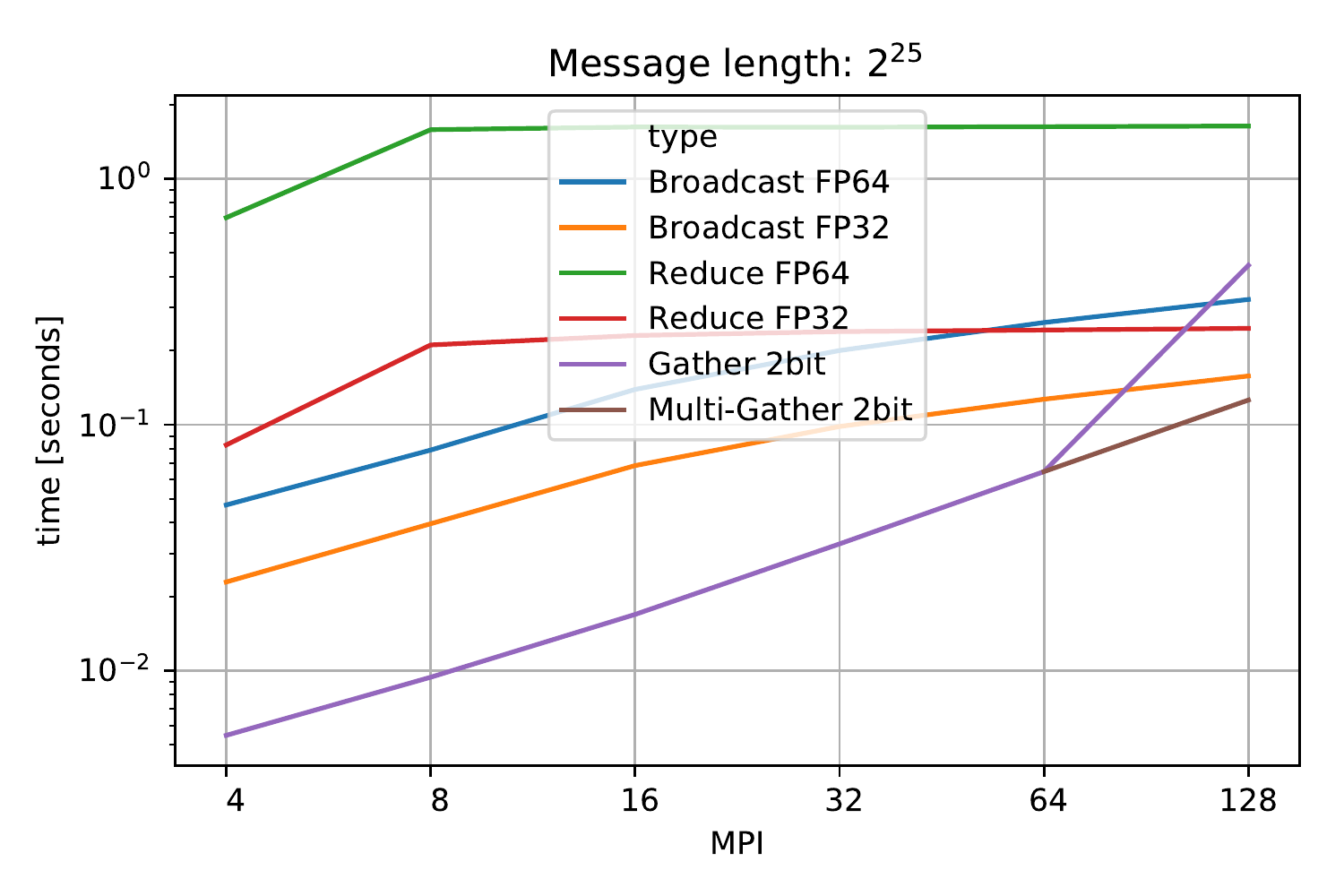}

\includegraphics[scale=0.46]{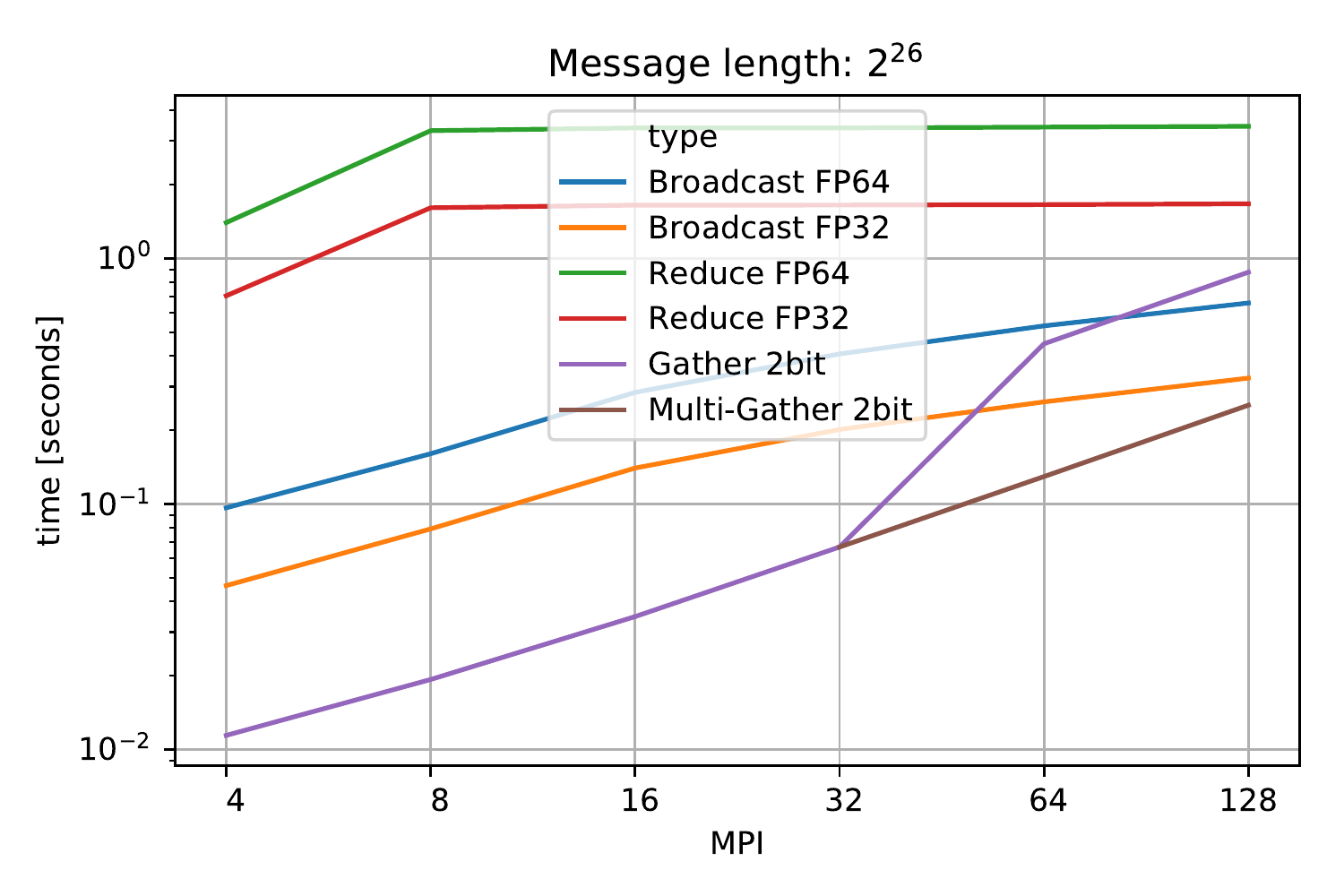}
\includegraphics[scale=0.46]{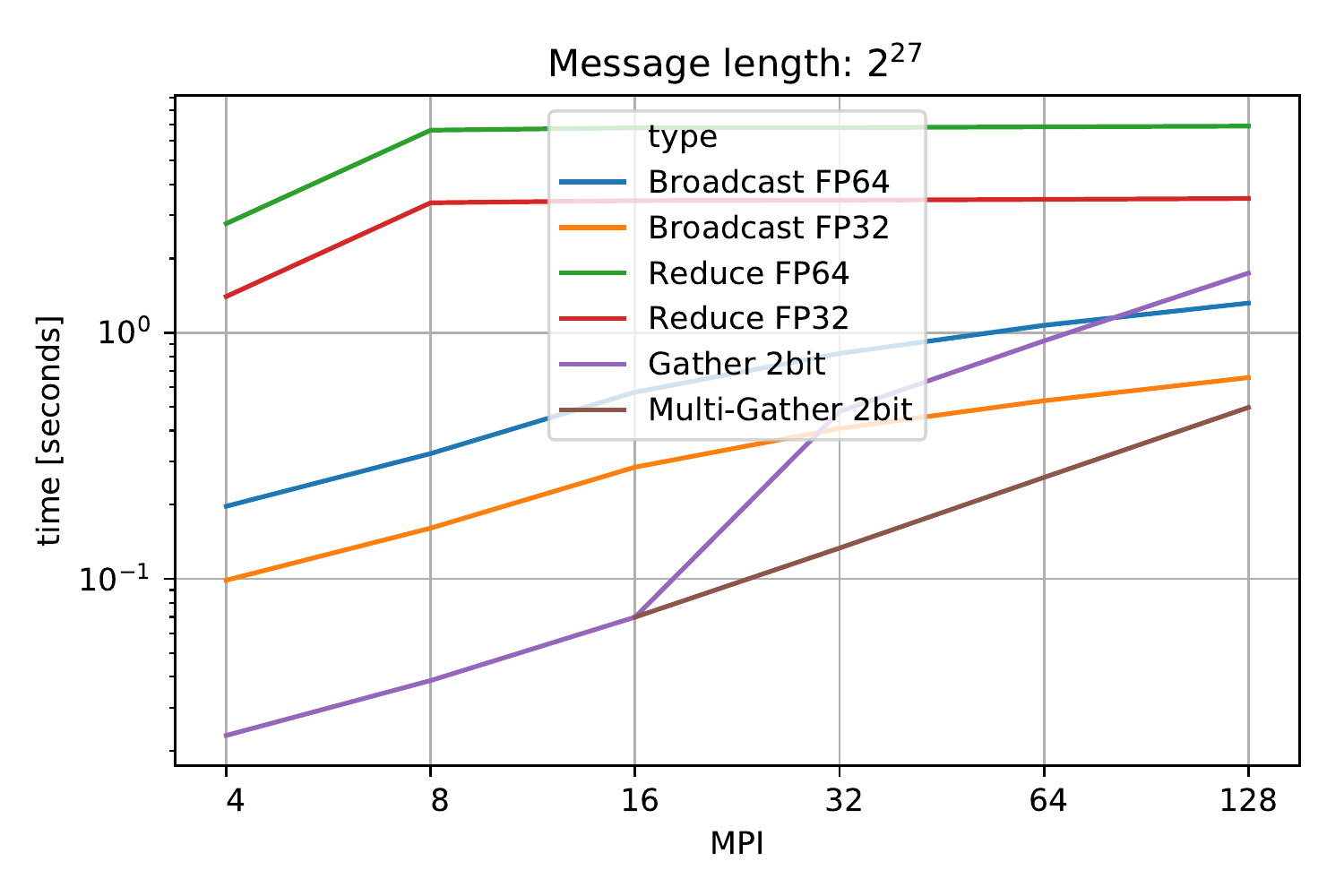}
\caption{
Time to communicate a vectors with different lengths for different methods as a function of \# of MPI processes. One can observe that Gather 2bit is not having nice scaling. We also show that the proposed Multi-Gather communication still achieves a nice scaling when more MPI processes are used.}

\label{fig:communication:scale}
\end{figure}
To correct for the unexpected behavior, 
we have performed MPI Gather multiple times on shorter vectors, such that the master node obtained all data, but in much faster time (see {\it Multi-Gather 2bit}).

\begin{figure}[b]
\centering

\includegraphics[scale=0.46]{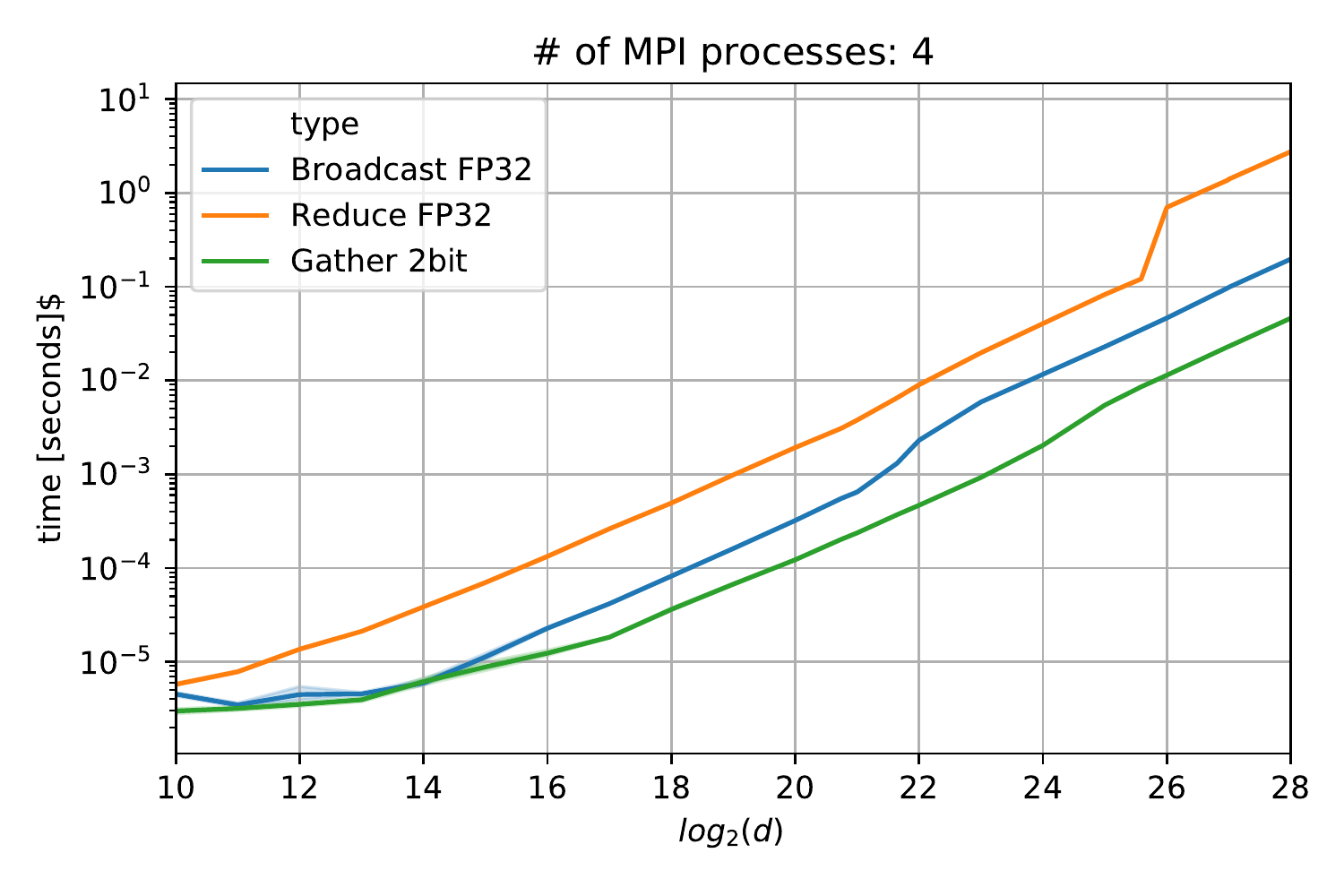}
\includegraphics[scale=0.46]{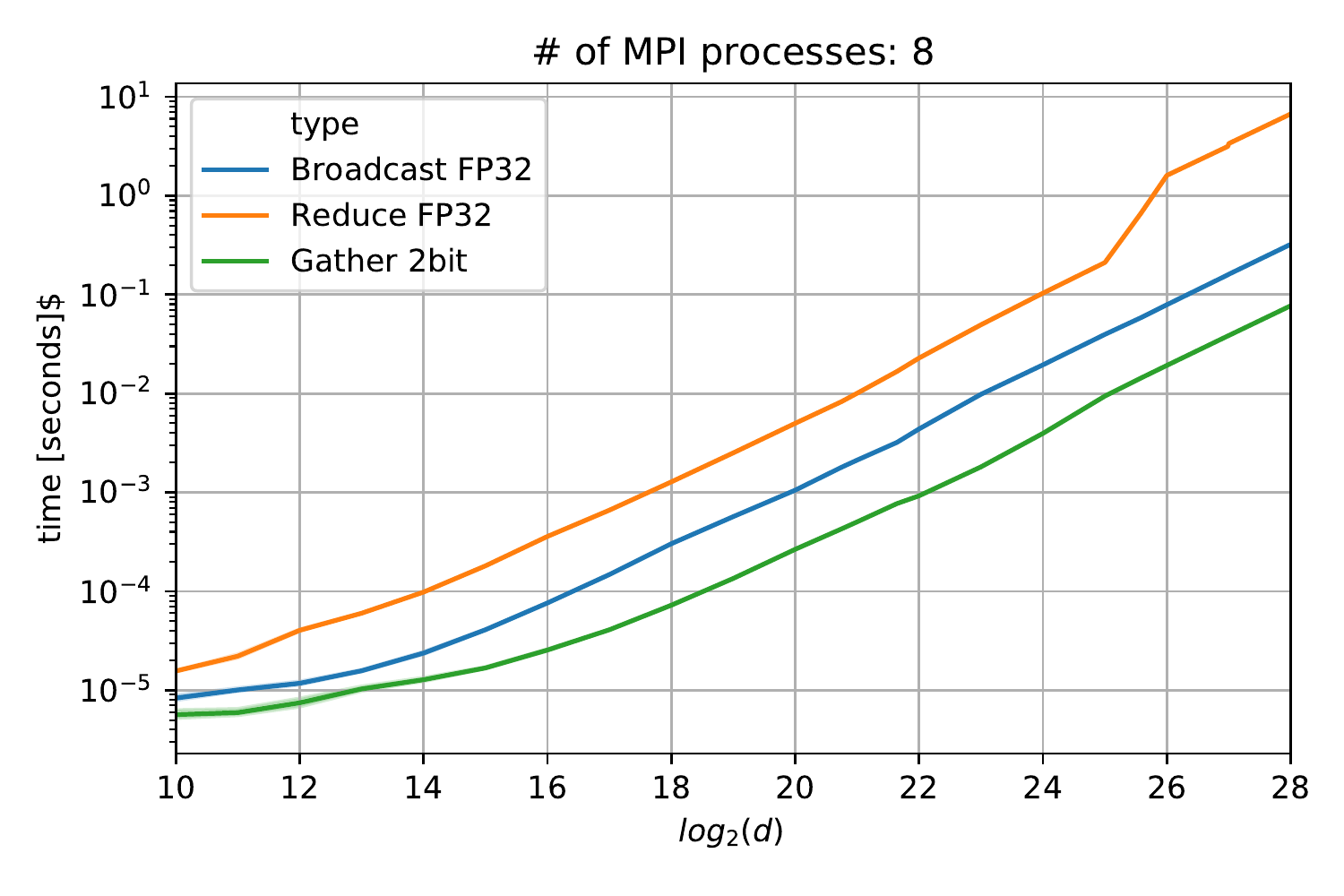}

\includegraphics[scale=0.46]{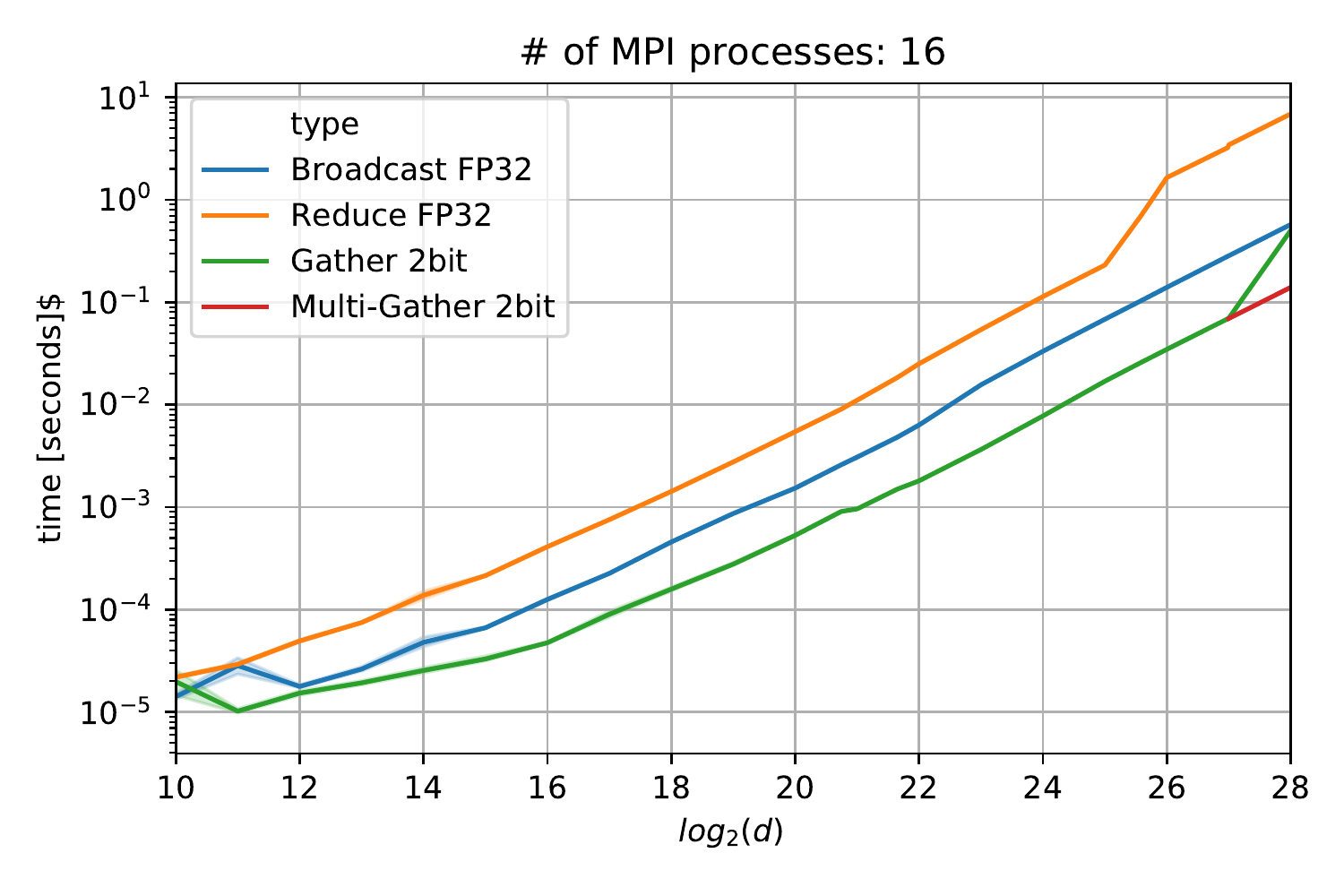}
\includegraphics[scale=0.46]{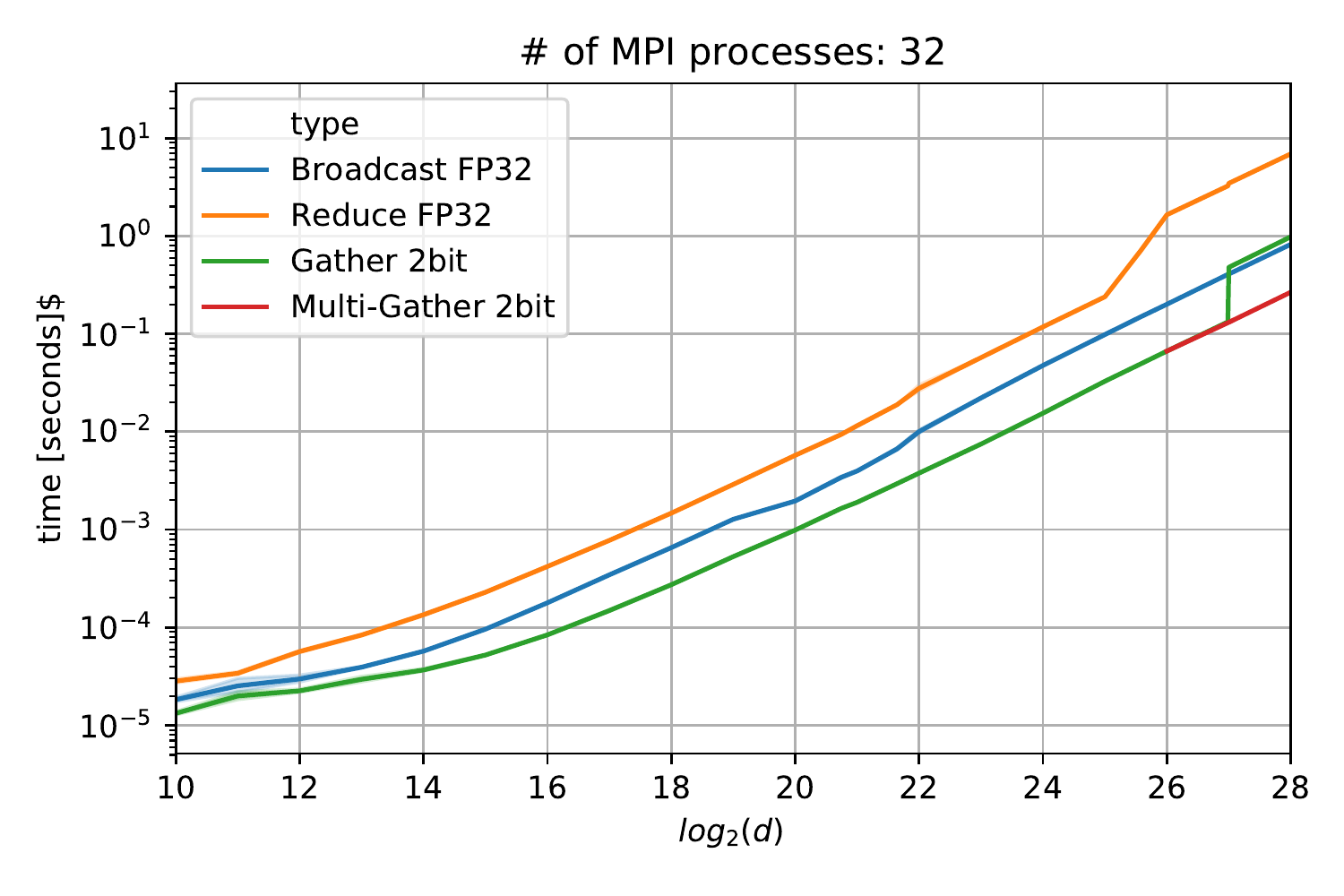}

\includegraphics[scale=0.46]{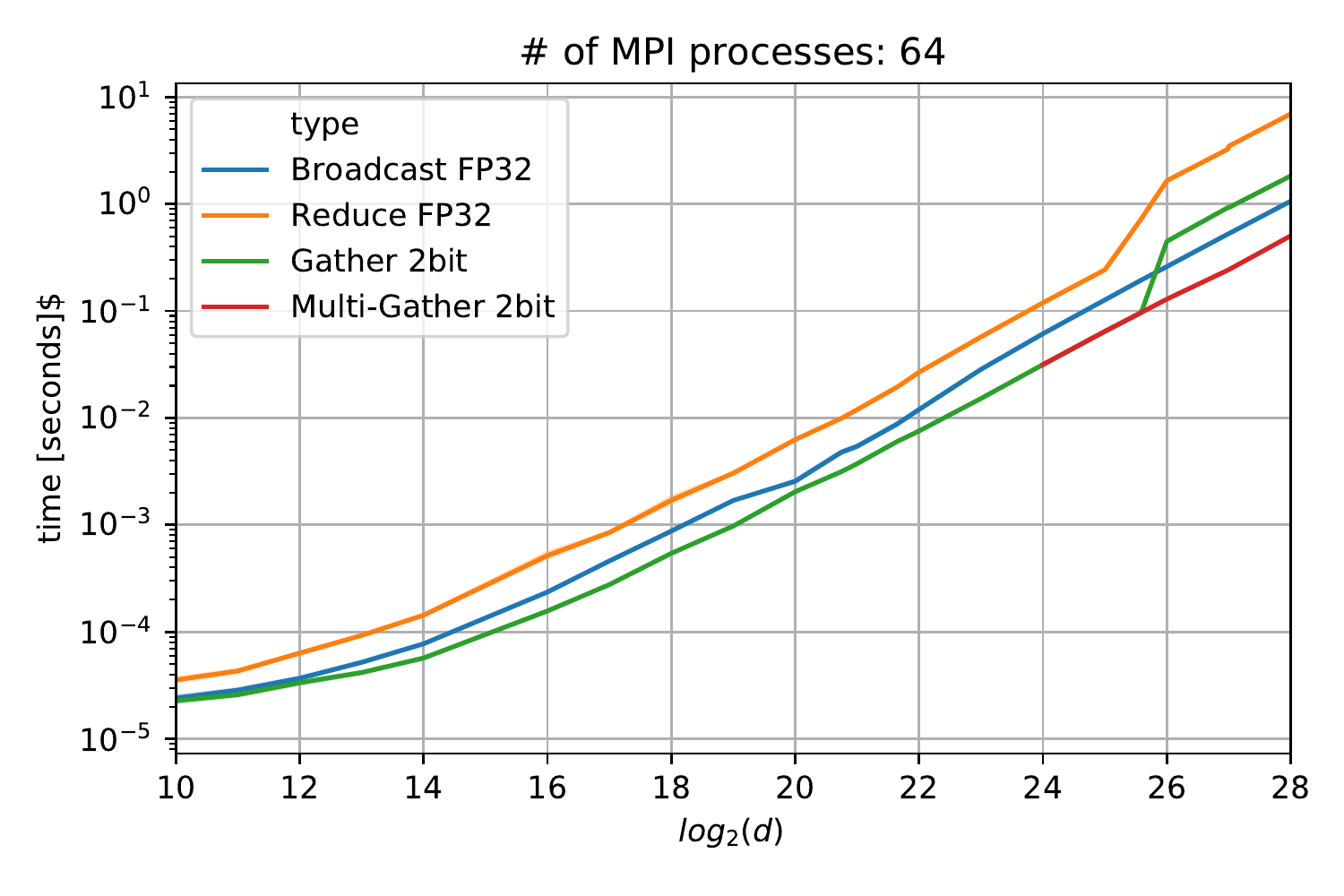}
\includegraphics[scale=0.46]{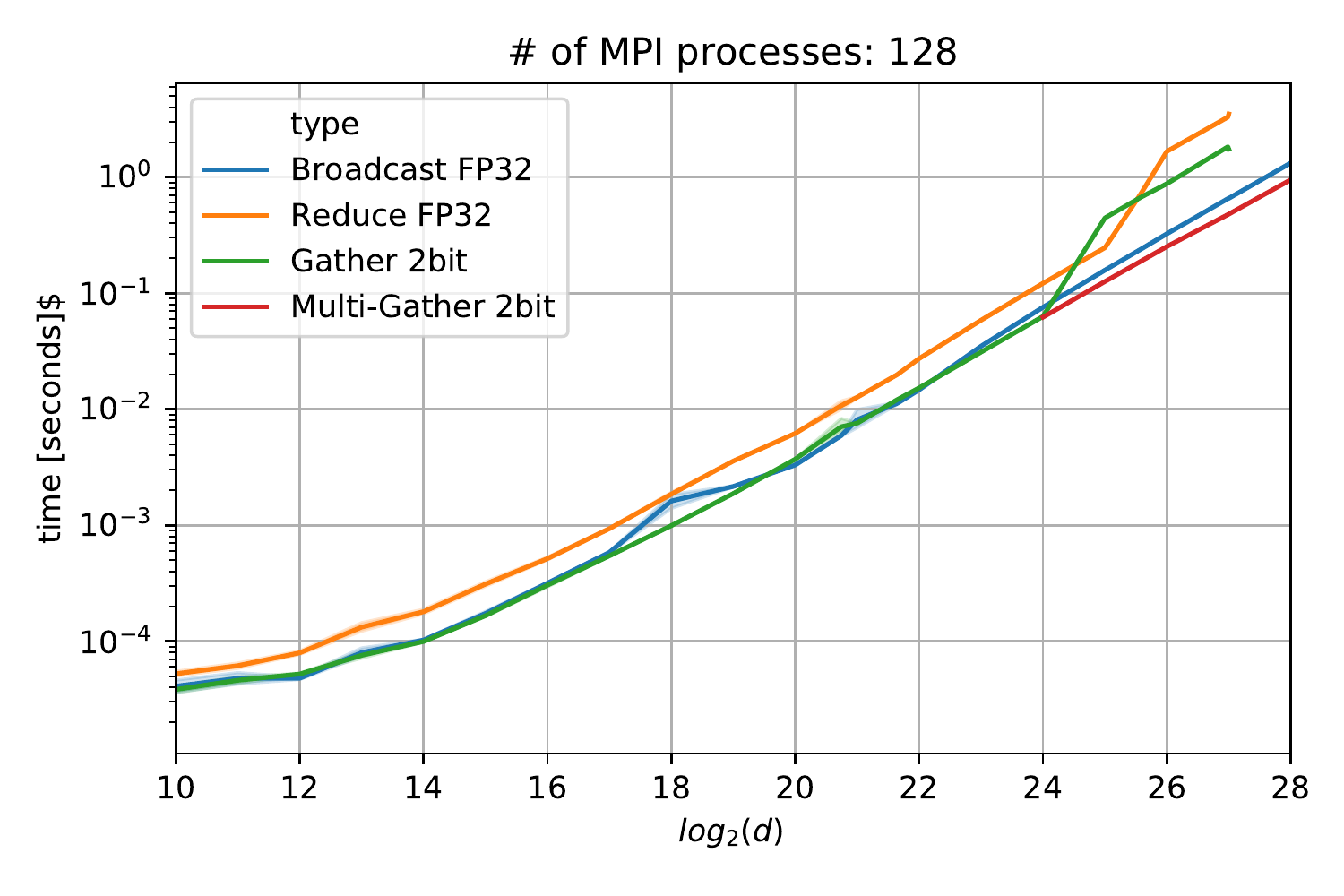}

\caption{The duration of communication for MPI Broadcast, MPI Reduce and MPI Gather. We show how the communication time depends on the size of the vector in $\R^d$ (x-axis) for various \# of MPI processes. In this experiment, we have run 4 MPI processes per computing node. For Broadcast and Reduce we have used a single precision floating point number. For Gather we used 2bits per dimension. For longer vectors and large number of MPI processes, one can observe that Gather has a very weird scaling issue. It turned out to be some weird behaviour of Cray-MPI implementation.}
\label{fig:communication_details}

\end{figure}

\begin{figure}[h!]
\centering

\includegraphics[scale=0.46]{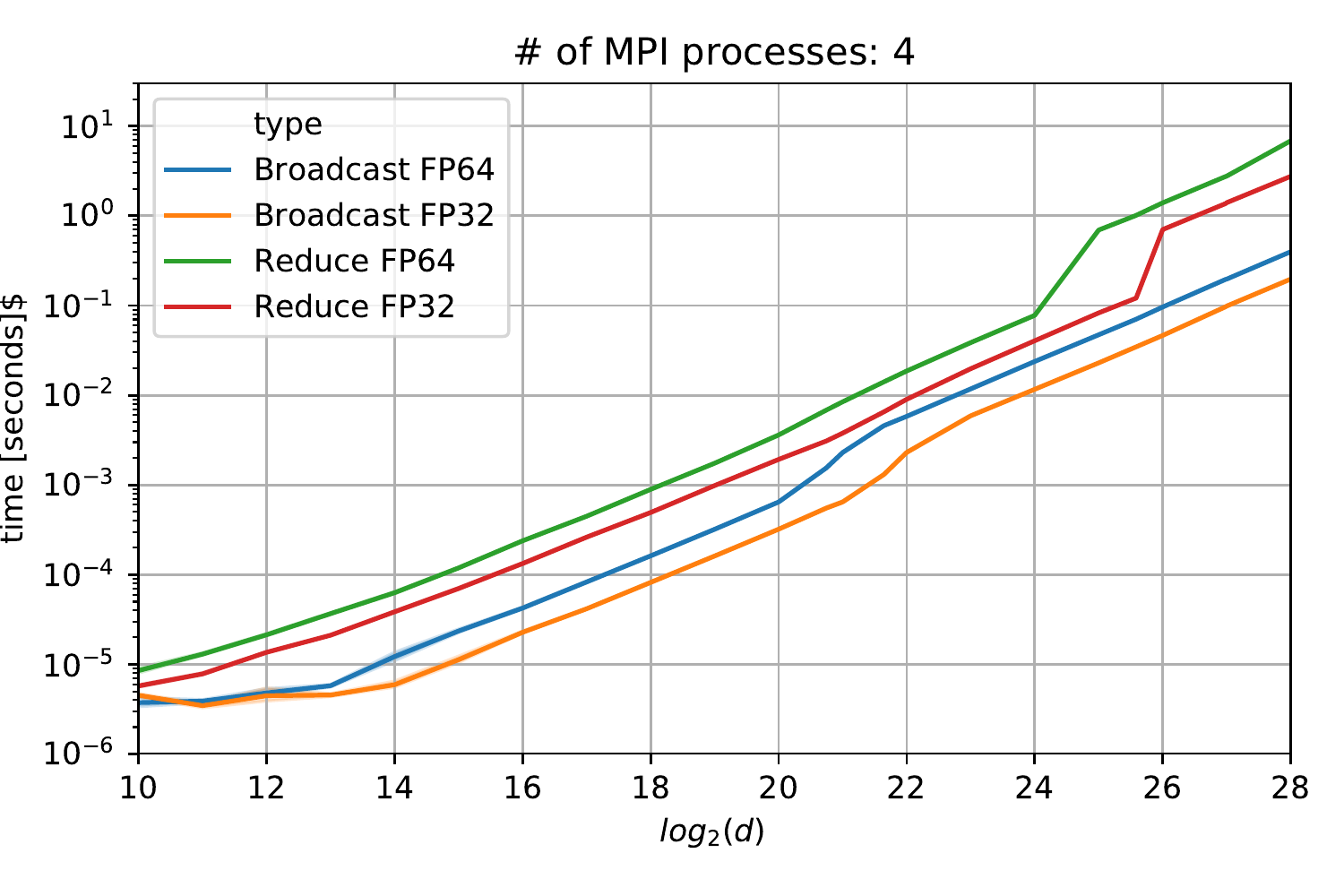}
\includegraphics[scale=0.46]{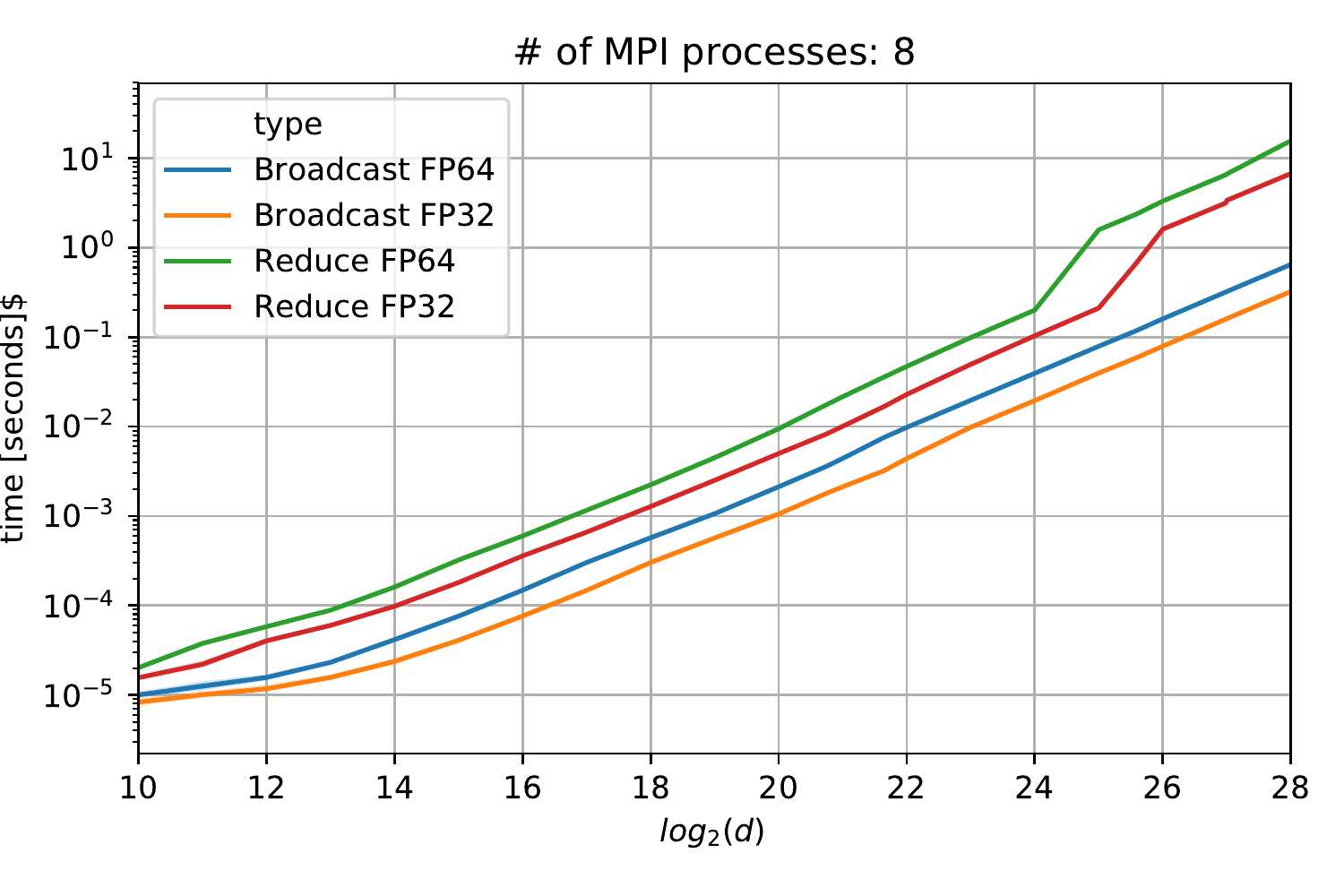}

\includegraphics[scale=0.46]{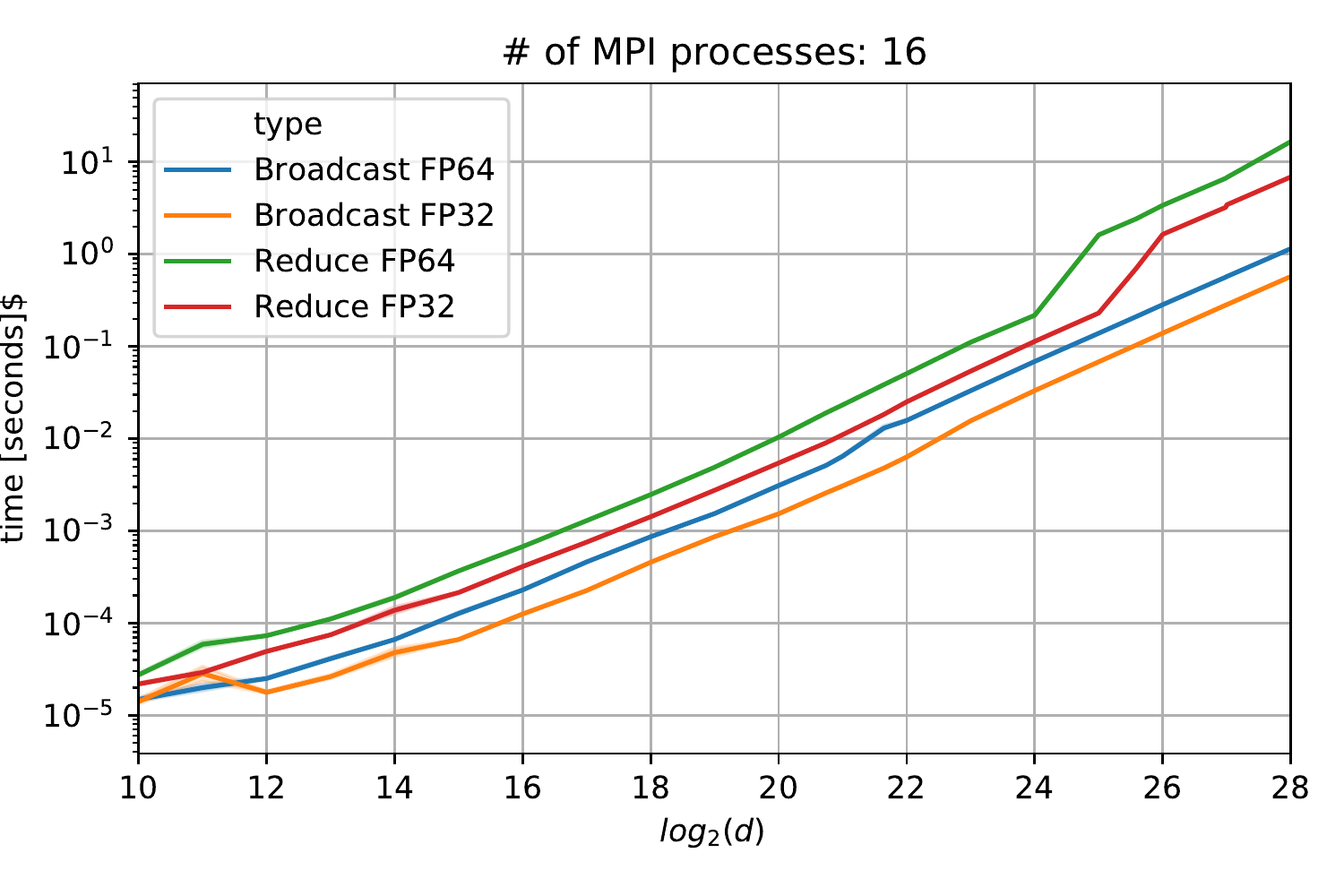}
\includegraphics[scale=0.46]{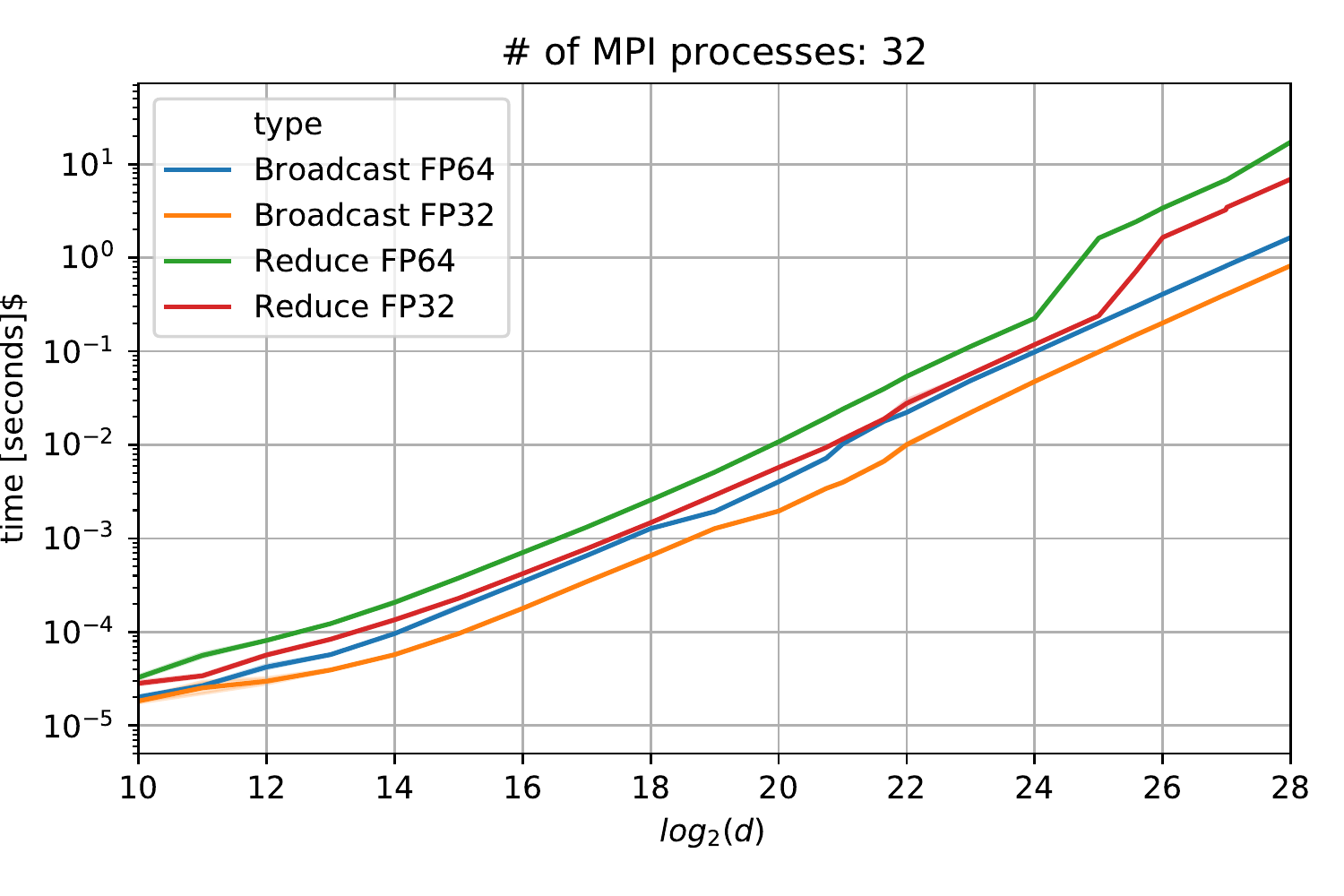}

\includegraphics[scale=0.46]{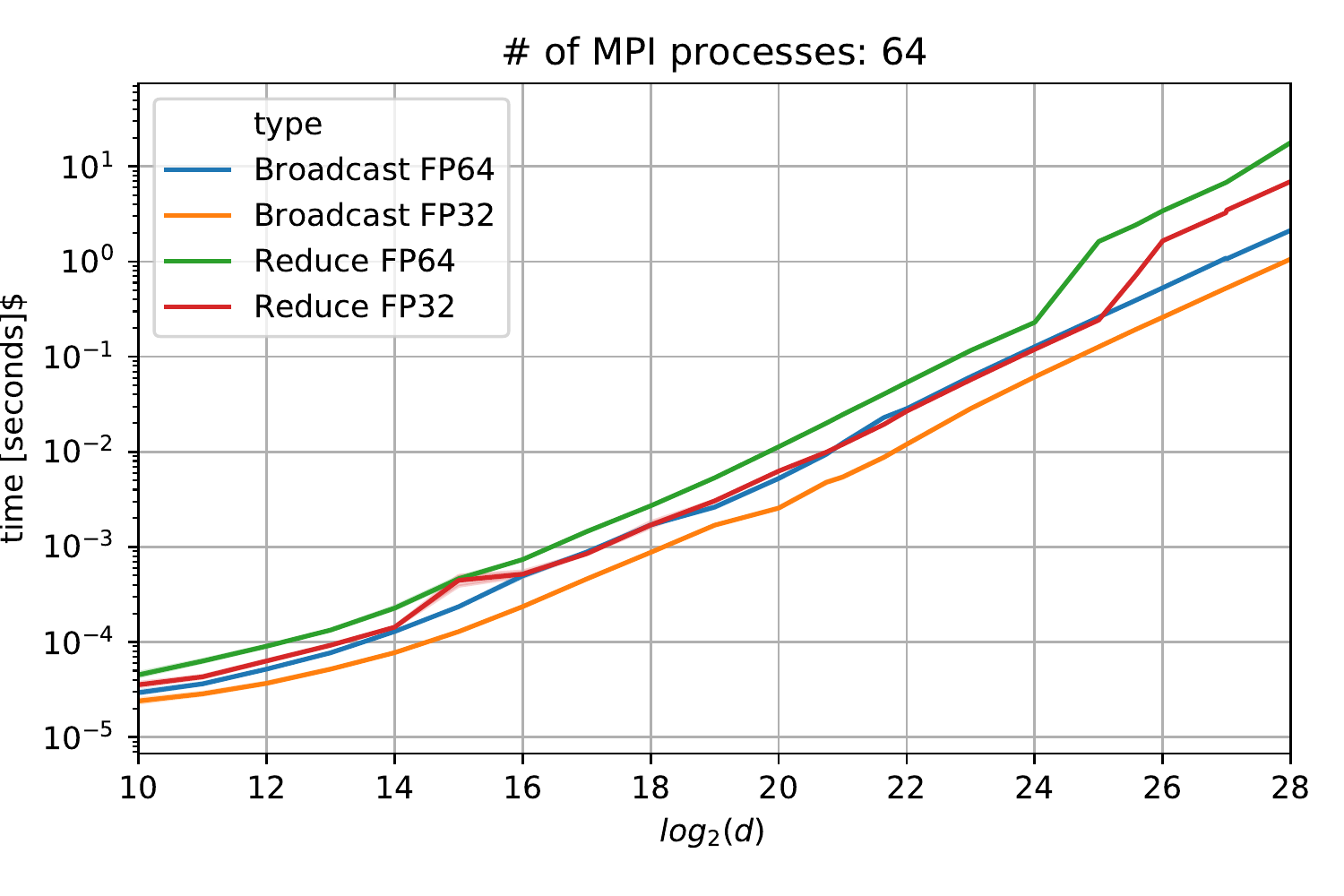}
\includegraphics[scale=0.46]{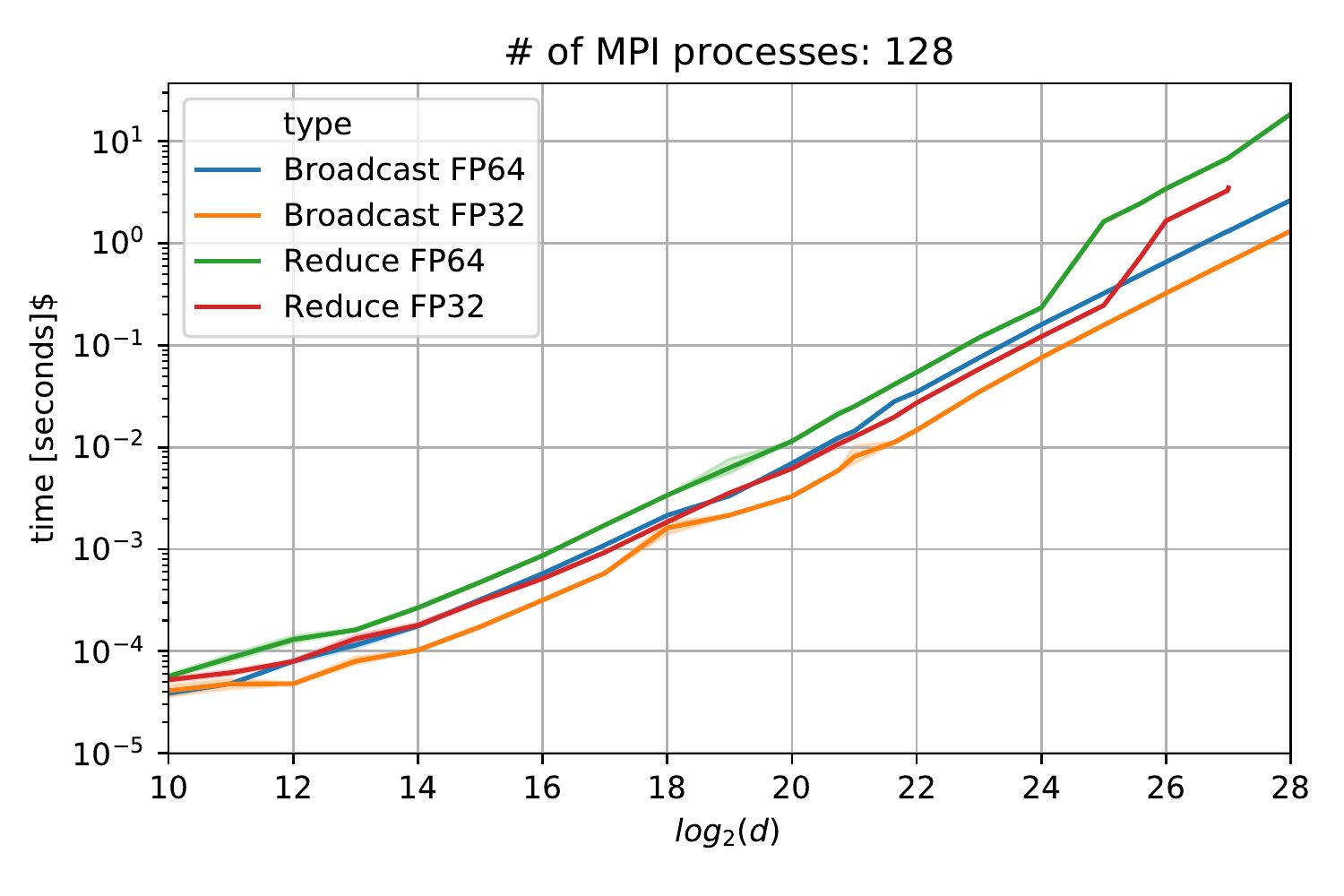}

\caption{The duration of communication for MPI Broadcast, MPI Reduce
for single precision (FP32) and double precision (FP64) floating numbers. We show how the communication time depends on the size of the vector in $\R^d$ (x-axis) for various \# of MPI processes. In this experiment, we have run 4 MPI processes per computing node. We have used Cray implementation of MPI.}
\label{fig:communication_32v64}

\end{figure}

\subsection{Performance of GPU}

In Table~\ref{tbl:networks}, we list the DNN networks that we have experimented with in this chapter.

\begin{table}[h!] 
\caption{Deep Neural Networks used in the experiments.}
\centering
\begin{tabular}{lccc}
\toprule  
Model &  $d$ &  \# classes & Input shape
 \\  \midrule
LeNet & 3.2M & 10 & $28\times28\times3$ 
\\  

CifarNet & 1.7M & 10 & $32\times32\times3$ 
\\

alexnet v2 & 50.3M & 1,000 & $224\times224\times3$ 
\\  
vgg a & 132.8M  & 1,000 & $224\times224\times3$
\\
\bottomrule 
\end{tabular}
\label{tbl:networks}

\end{table}

Figure~\ref{fig:computationCost}
shows the performance of a single P100 GPU
for different batch size, DNN network and operation.
\begin{figure}[h!]
\centering

\includegraphics[scale=0.7]{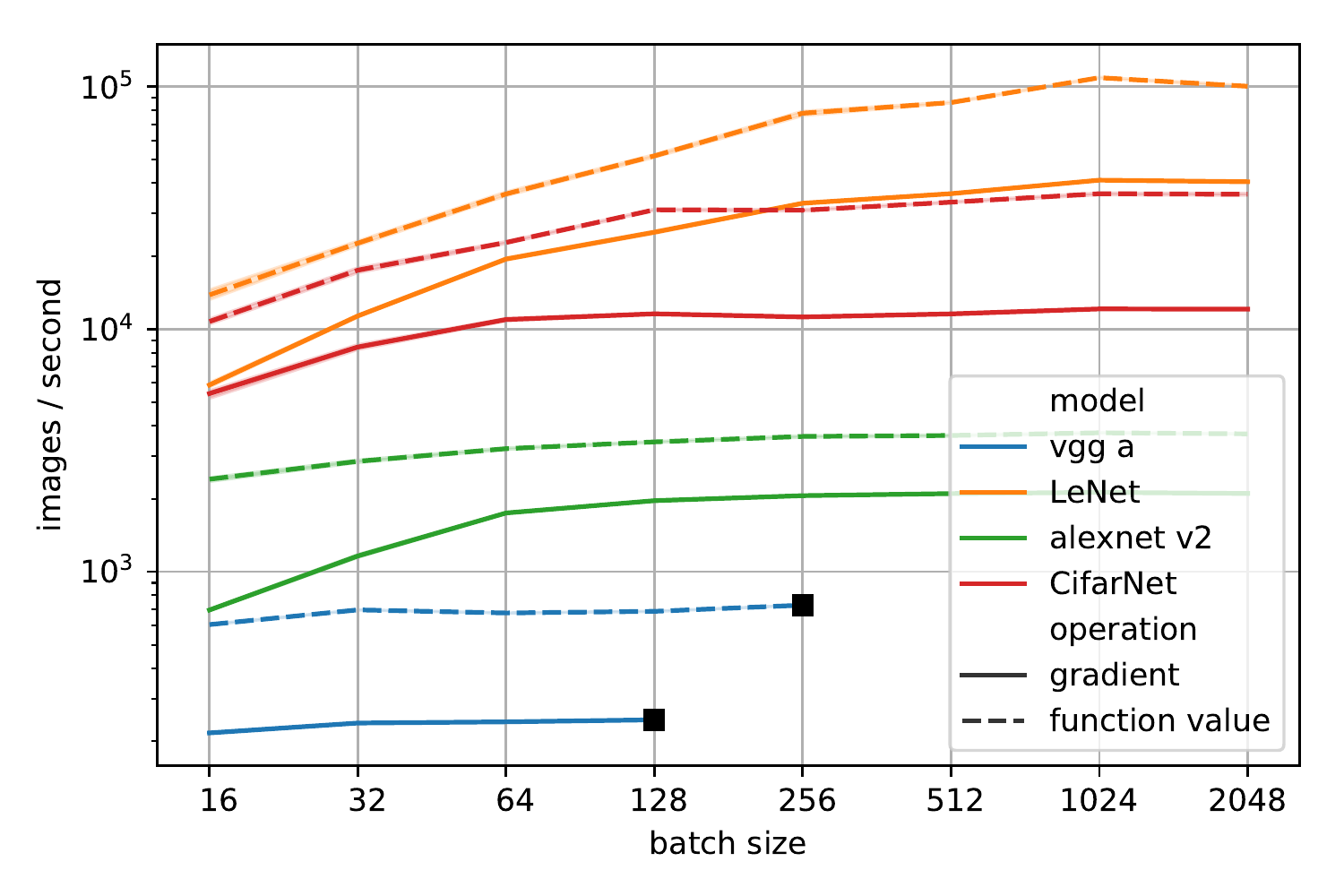}

\caption{The performance (images/second) of NVIDIA Tesla P100 GPU  on 4 different problems as a function of batch size.
We show how different choice of batch size affects the speed of function evaluation and gradient evaluation.
For vgg a, we have run out of memory on GPU for batch size larger than 128 (gradient evaluation) and 256 for function evaluation.
Clearly, this graph suggest that choosing small batch size leads to small utilization of GPU. Note that using larger batch size do not necessary reduce the training process.
}

\label{fig:computationCost}

\end{figure}

\subsection{\algname{DIANA} vs.\ \algname{TernGrad}, \algname{SGD} and \algname{QSGD}}

In Figure~\ref{fig:imagesPerSecond}
we compare the performance of \algname{DIANA} vs.\ doing an MPI reduce operation with 32bit floats. The computing cluster had Cray Aries High Speed Network.
However, for \algname{DIANA} we used 2bit per dimension, we have experienced an weird scaling behaviour, which was documented also in\cite{parker2018performance}.
In our case, this affected speed for alexnet and vgg\_a beyond 64 or 32 MPI processes respectively.
For more detailed experiments, see Section~\ref{sec:A:MPI}.
In order to improve the speed of Gather, we impose a Multi-Gather strategy, when we call Gather multiple-times on shorter vectors. This significantly improved the communication cost of Gather (see 
Figure \ref{fig:communication_details}) and leads to much nicer scaling -- see green bars -- \algname{DIANA}-MultiGather in Figure~\ref{fig:imagesPerSecond}).
\begin{figure}[h!]

\centering 
\includegraphics[scale=0.46]{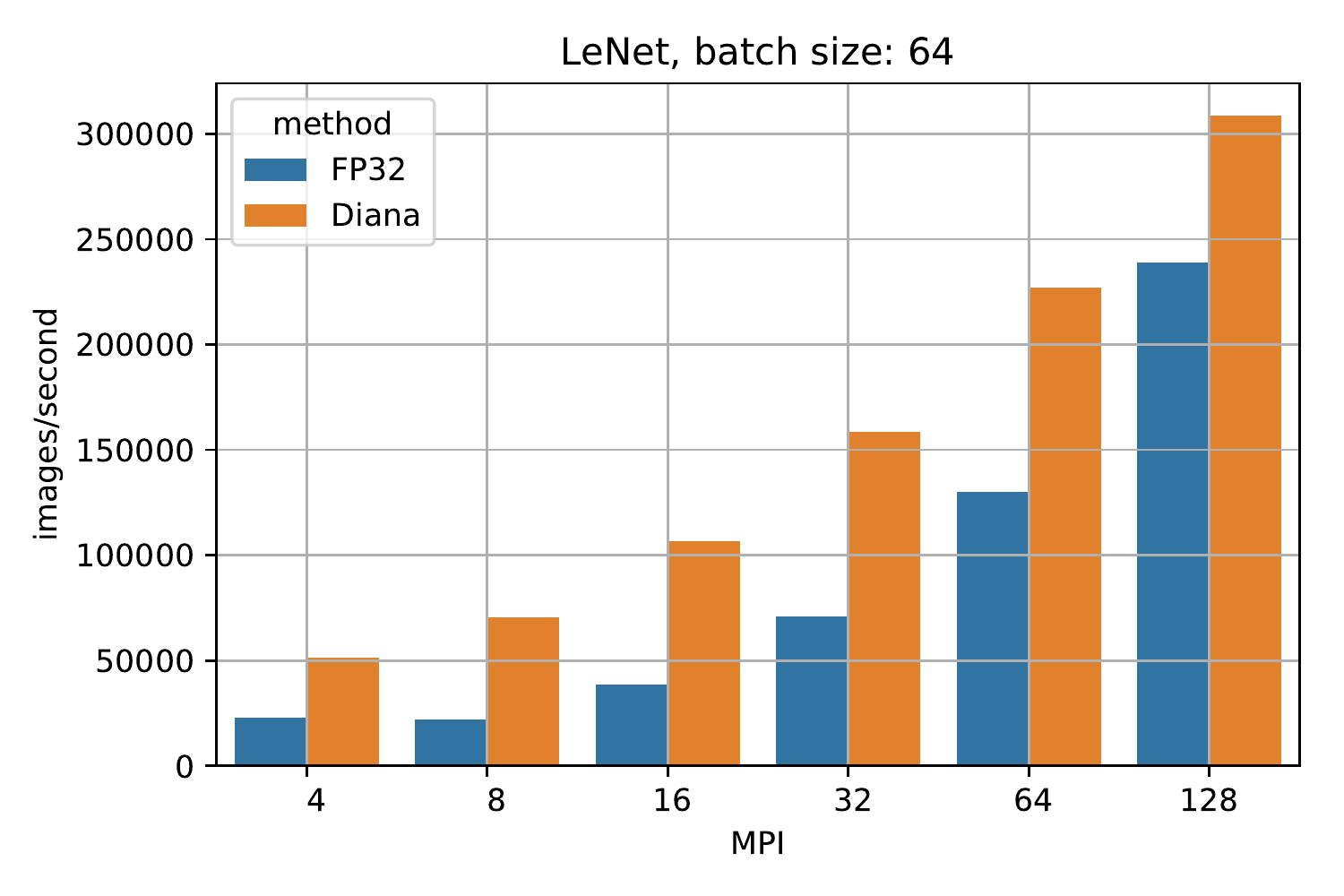}
\includegraphics[scale=0.46]{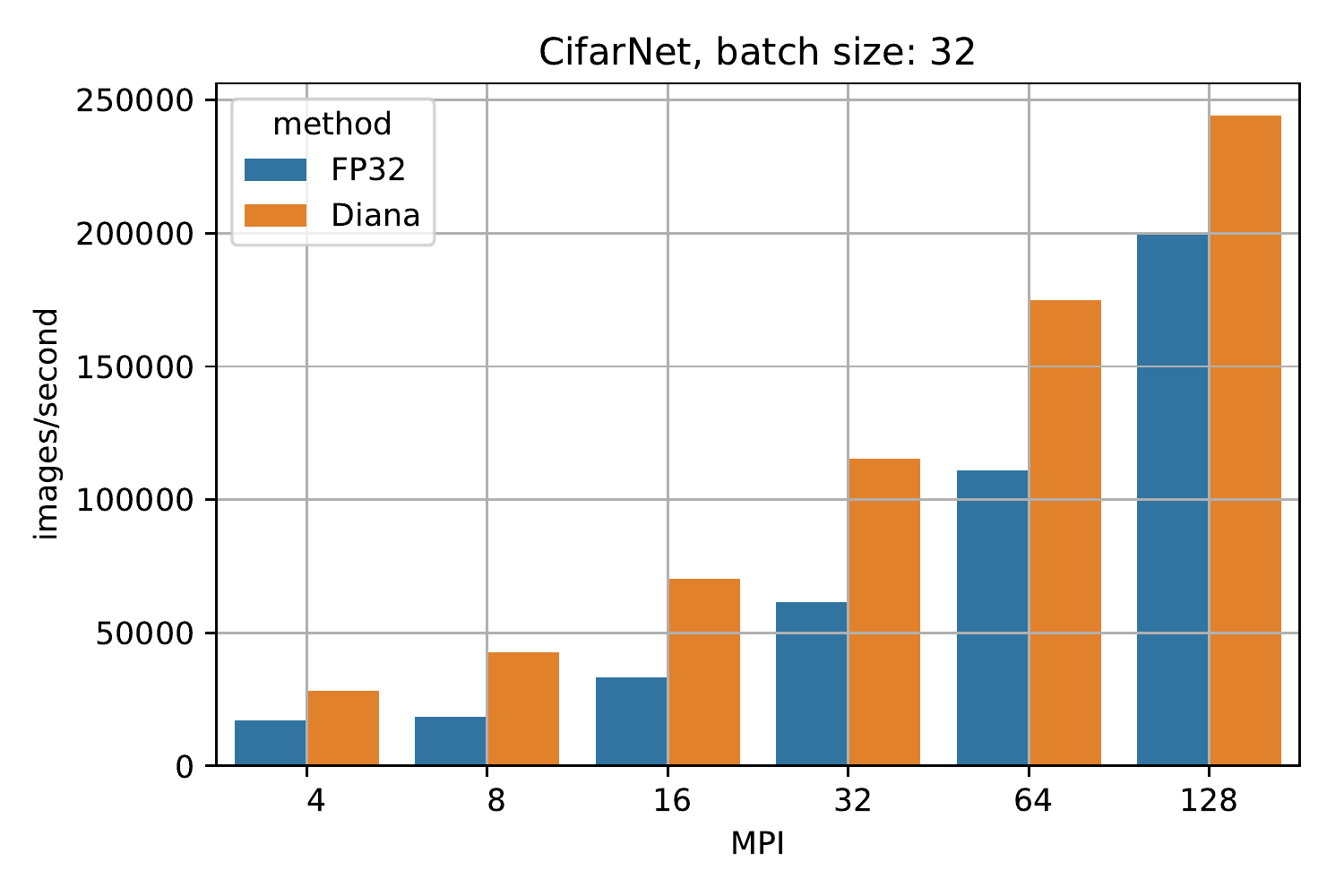}

\includegraphics[scale=0.46]{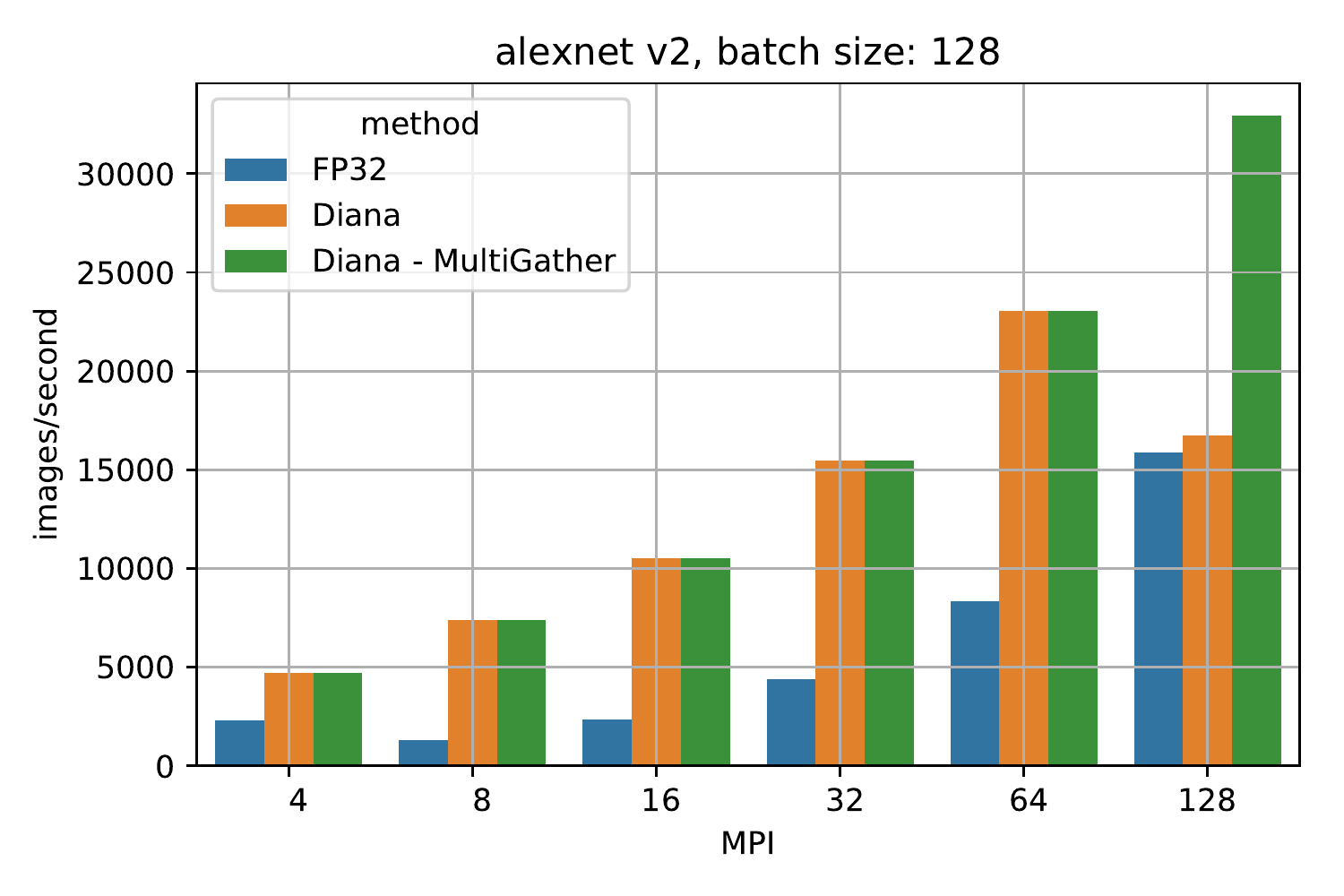}
\includegraphics[scale=0.46]{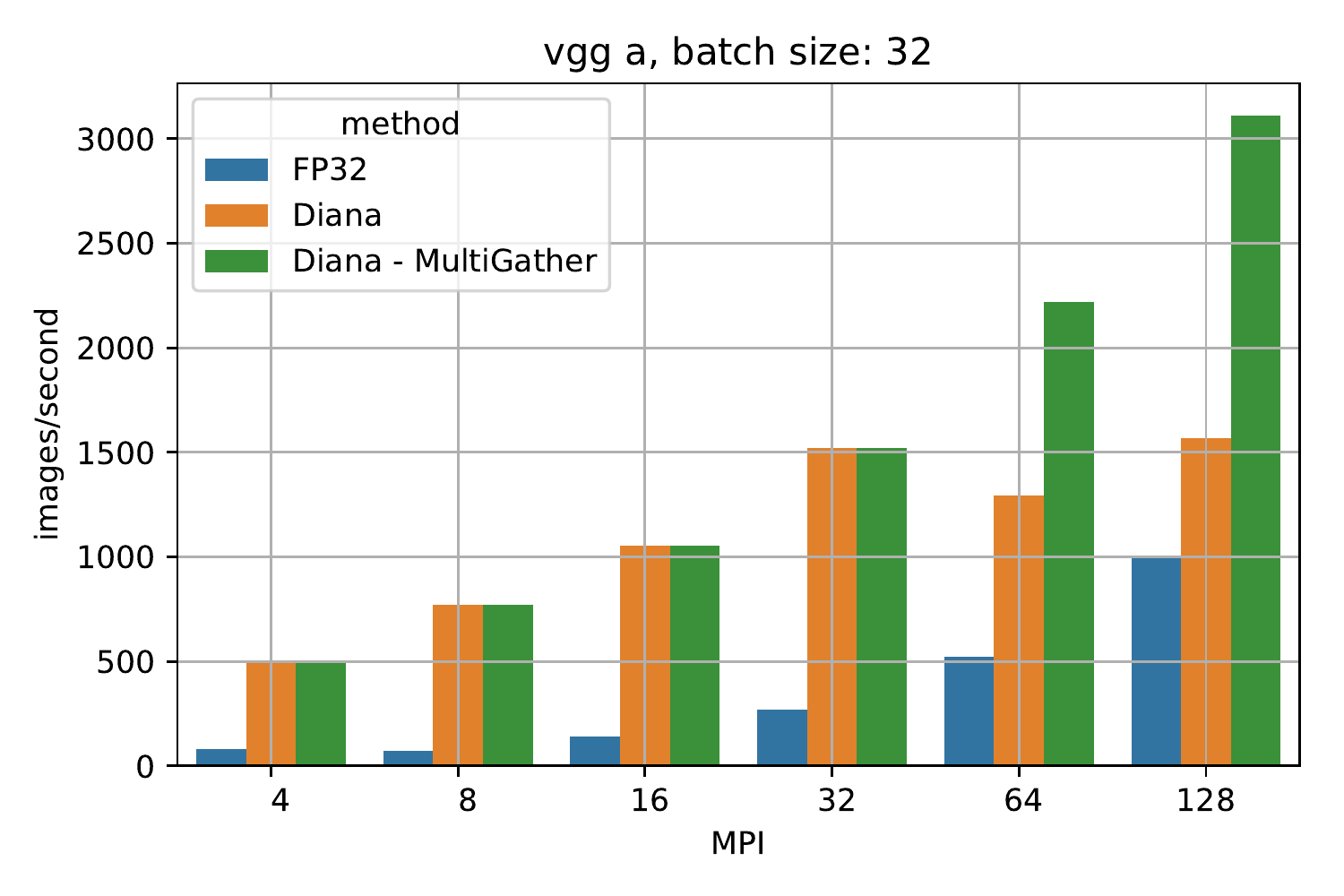}

\caption{Comparison of performance (images/second) for various number of GPUs/MPI processes and sparse communication \algname{DIANA} (2bit) vs.\ Reduce with 32bit float (FP32).
We have run 4 MPI processes on each node. Each MPI process is using single P100 GPU. 
Note that increasing MPI from 4 to 8 will not bring any significant improvement for FP32, because with 8 MPI processes, communication will happen between computing nodes and will be significantly slower when compare to the single node communication with 4MPI processes.
}
\label{fig:imagesPerSecond}

\end{figure}

In the next experiments, we run \algname{QSGD} \cite{alistarh2017qsgd}, \algname{TernGrad} \cite{wen2017terngrad}, \algname{SGD} with momentum and \algname{DIANA} on MNIST dataset and CIFAR-10 dataset for 3 epochs. We have selected 8 workers and run each method with stepsize chosen from $\{0.1, 0.2, 0.05\}$.
For \algname{QSGD}, \algname{DIANA} and \algname{TernGrad}, we also tried various quantization bucket sizes in $\{32, 128, 512\}$.
For \algname{QSGD} we have chosen $2,4,8$ quantization levels.
For \algname{DIANA} we have chosen $\alpha \in 
\{0, 1.0/\sqrt{\mbox{quantization bucket sizes }}\}$
and have selected initial $h = 0$. 
For \algname{DIANA} and \algname{SGD} we also run a momentum version, with a momentum parameter in $\{0, 0.95, 0.99\}$.
For \algname{DIANA} we also run with two choices of norm $\ell_2$ and $\ell_\infty$.
For each experiment we have selected softmax cross entropy loss. MNIST-Convex is a simple DNN with no hidden layer, MNIST-DNN is a convolutional NN described here \url{https://github.com/floydhub/mnist/blob/master/ConvNet.py}
and CIFAR-10-DNN is a convolutional DNN described here
\url{https://github.com/kuangliu/pytorch-cifar/blob/master/models/lenet.py}.
In Figure~\ref{fig:DNN:evolution} we show the best runs over all the parameters for all the methods. 
For MNIST-Convex \algname{SGD} and \algname{DIANA} makes use of the momentum and dominate all other algorithms.
For MNIST-DNN situation is very similar.
For CIFAR-10-DNN  both \algname{DIANA} and \algname{SGD} have significantly outperform other methods.
\begin{figure}[h!]
\centering 
\includegraphics[width=0.46\textwidth]{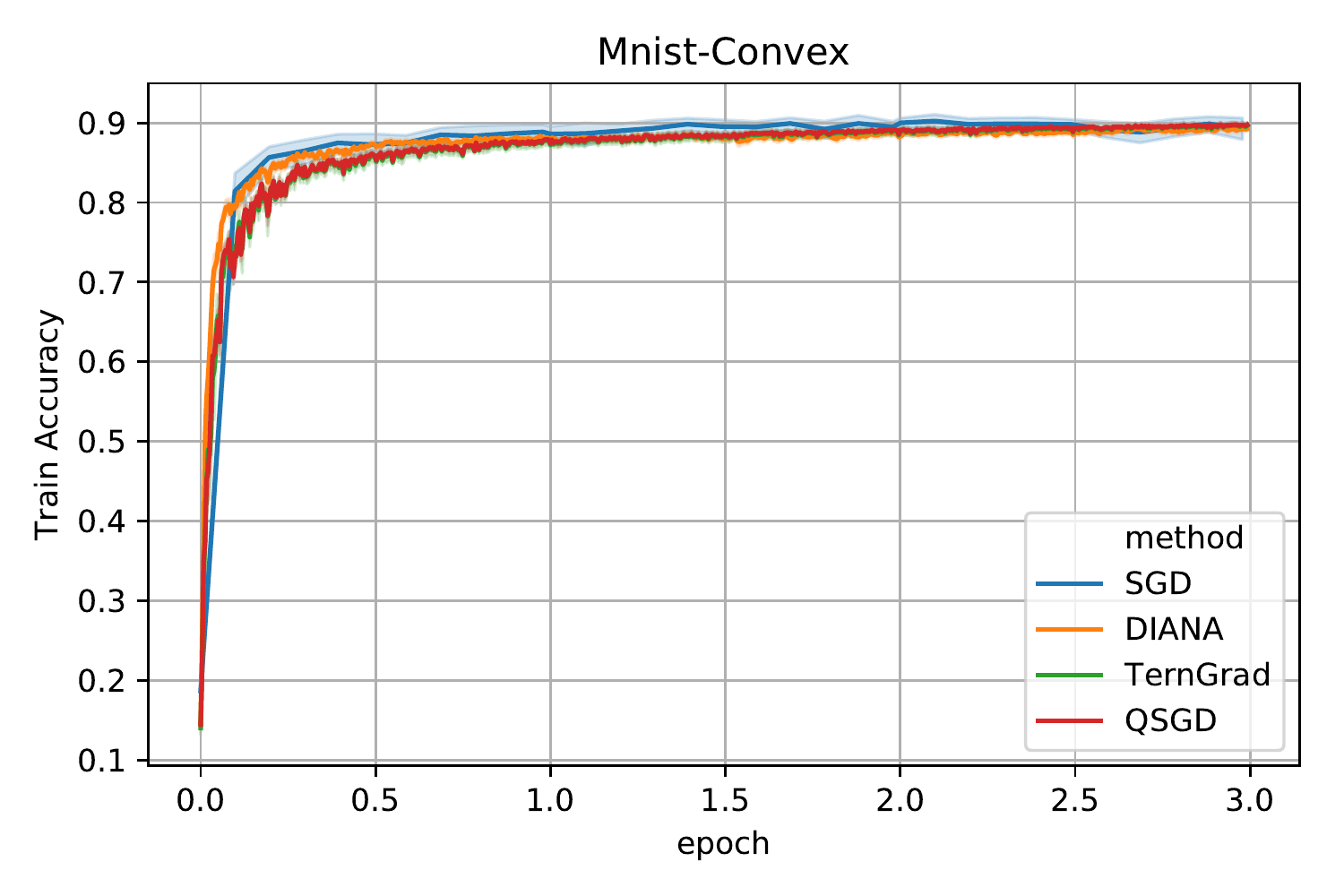}
\includegraphics[width=0.46\textwidth]{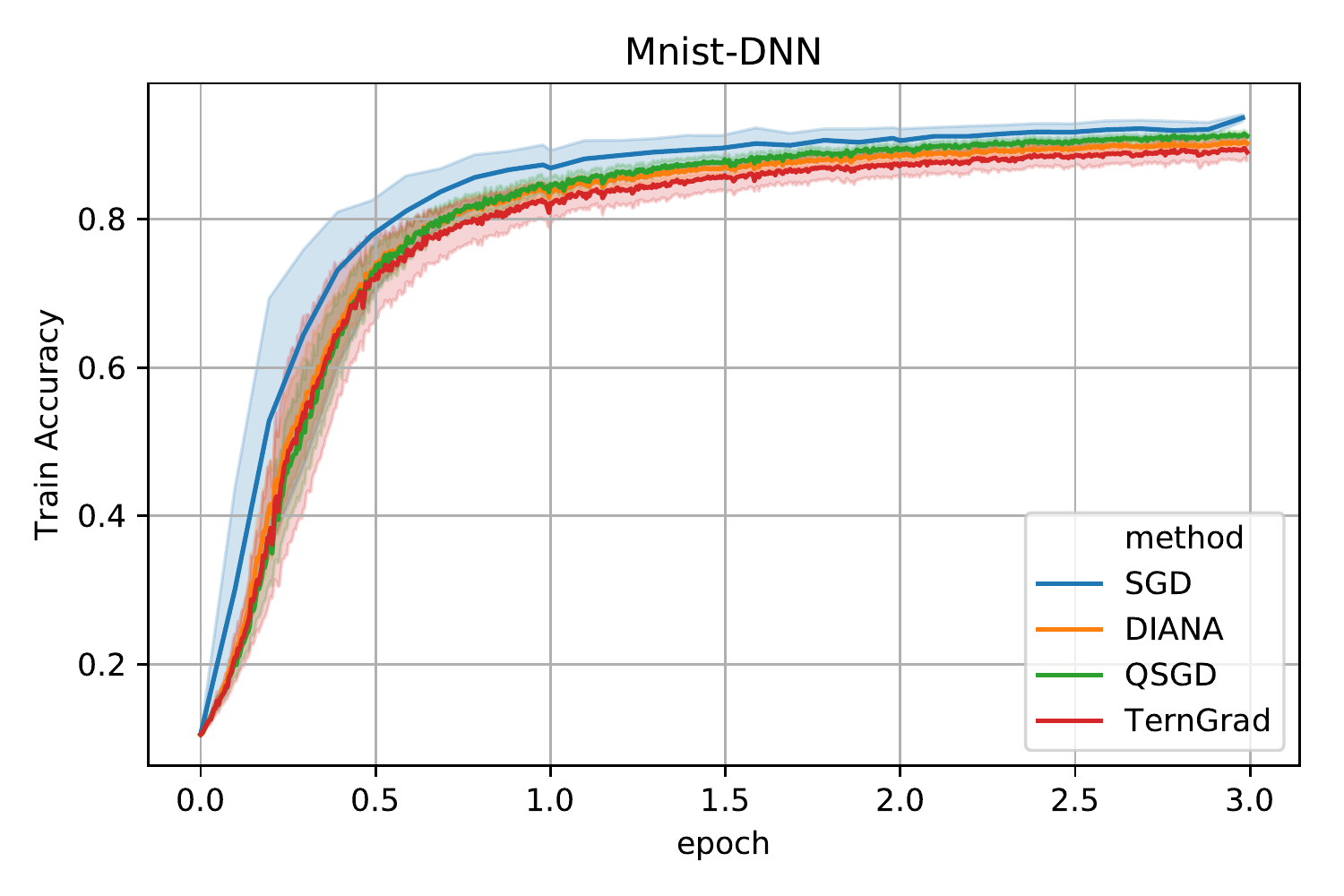}
\includegraphics[width=0.46\textwidth]{dnn_evolution_cifar10_N2_train.pdf}

\includegraphics[width=0.46\textwidth]{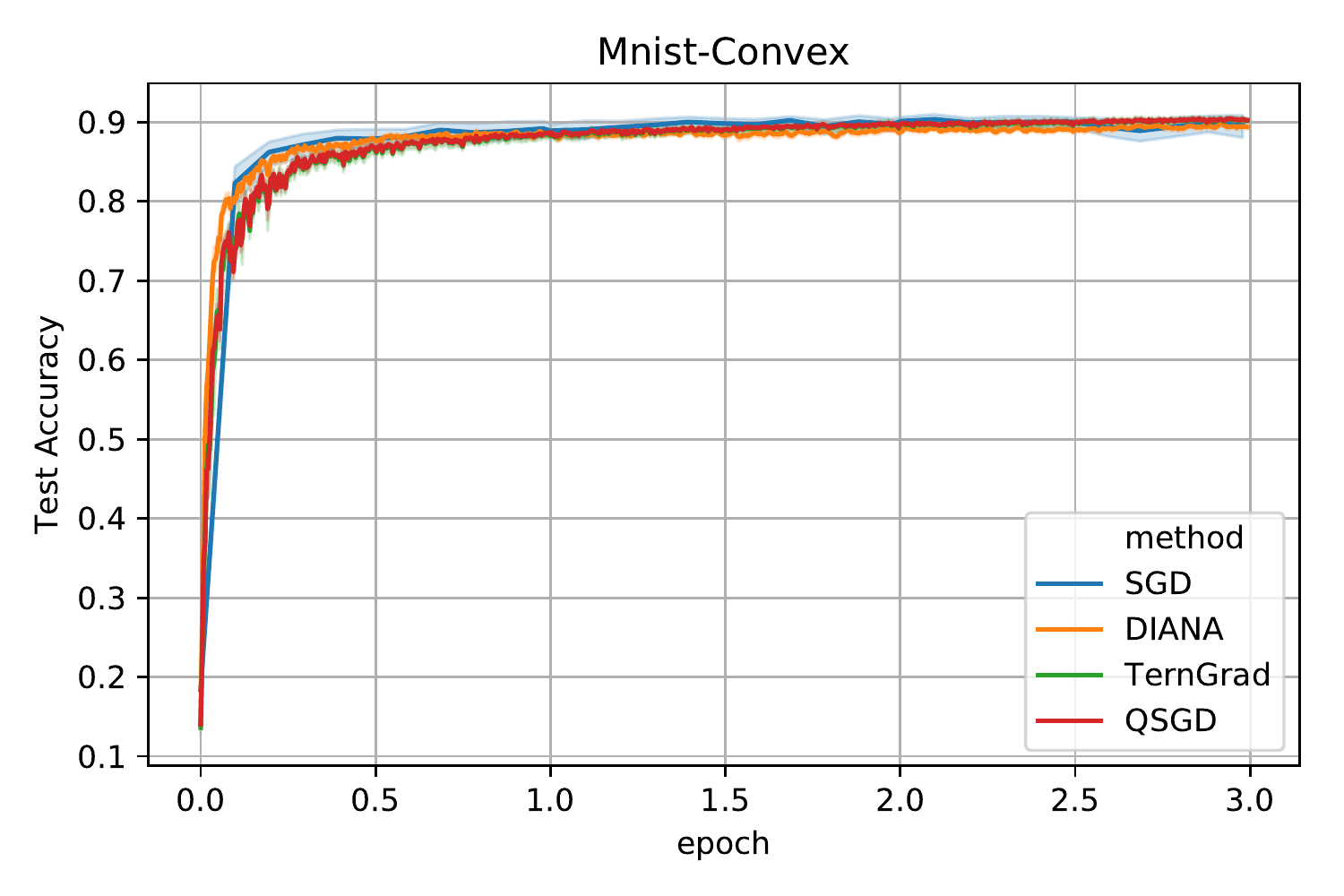}
\includegraphics[width=0.46\textwidth]{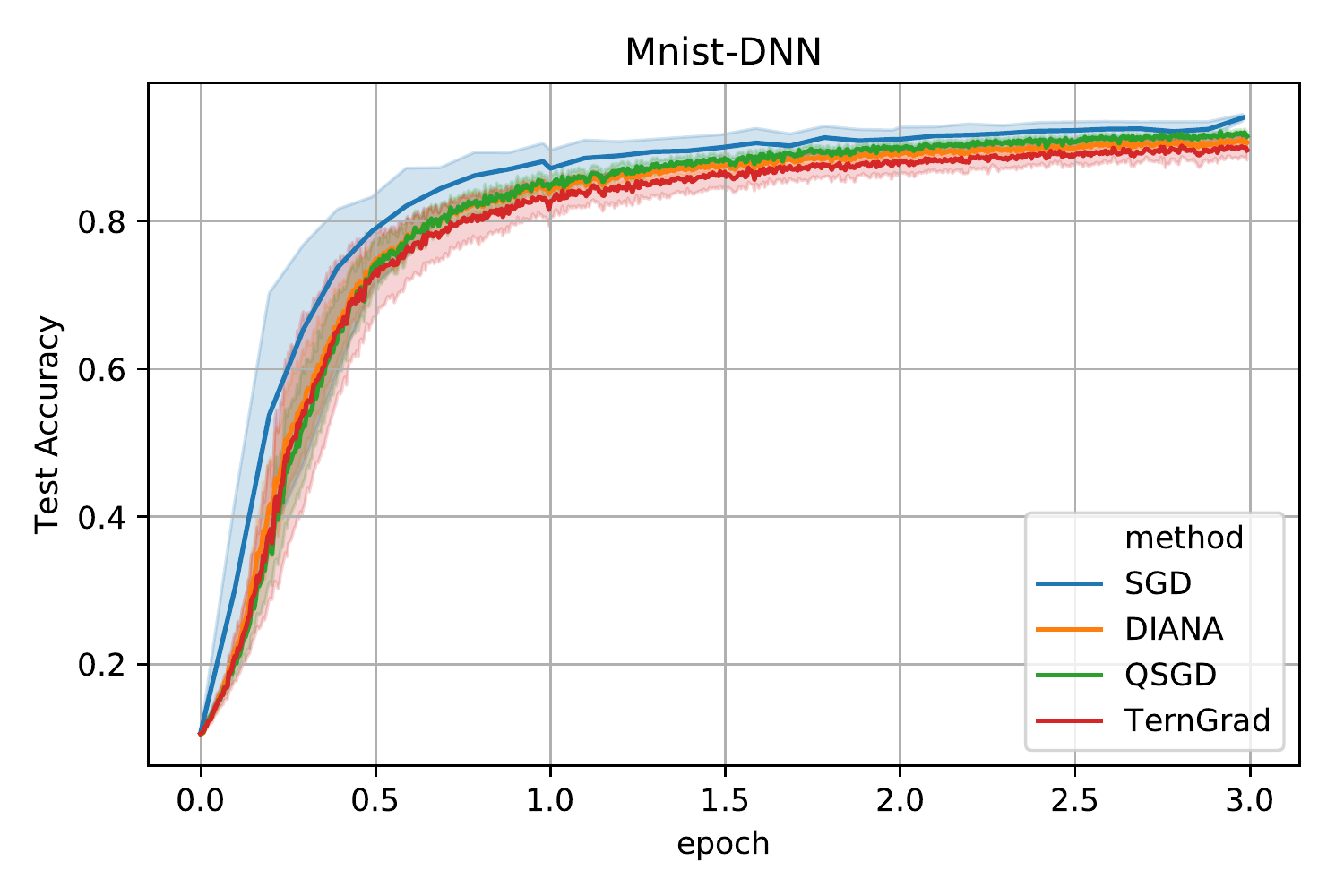}
\includegraphics[width=0.46\textwidth]{dnn_evolution_cifar10_N2_test.pdf}

\caption{Evolution of training and testing accuracy for 3 different problems, using 4 algorithms: \algname{DIANA}, \algname{SGD}, \algname{QSGD} and \algname{TernGrad}. 
We have chosen the best runs over all tested hyper-parameters.}
\label{fig:DNN:evolution}

\end{figure}

In Figure~\ref{fig:DNN:sparsity} show the evolution of sparsity of
the quantized gradient for the 3 problems and \algname{DIANA}, \algname{QSGD} and \algname{TernGrad}. For MNIST-DNN, it seems that the quantized gradients are becoming sparser as the training progresses.

\begin{figure}[h!]

\begin{center}
\includegraphics[width=0.46\textwidth]{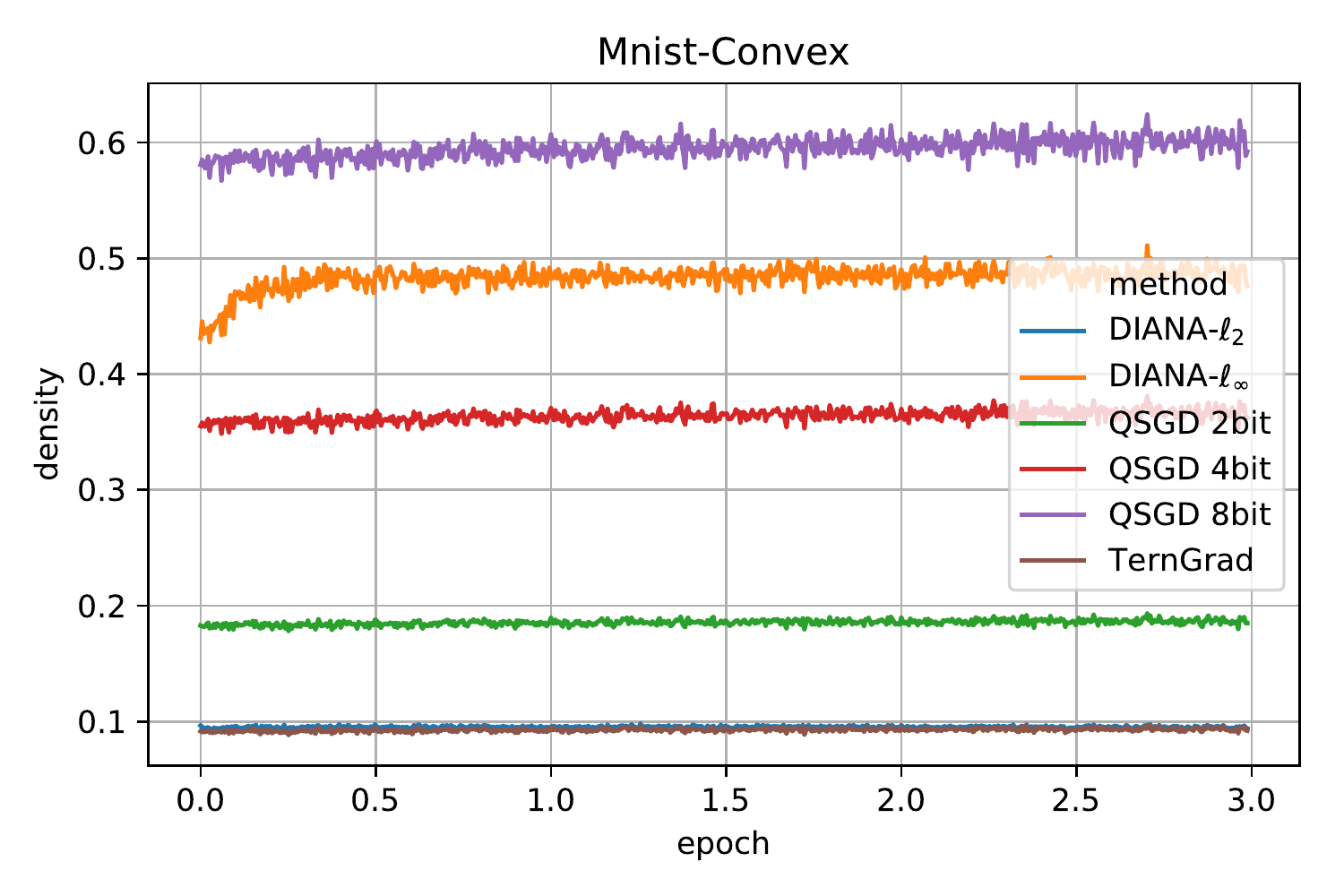}
\includegraphics[width=0.46\textwidth]{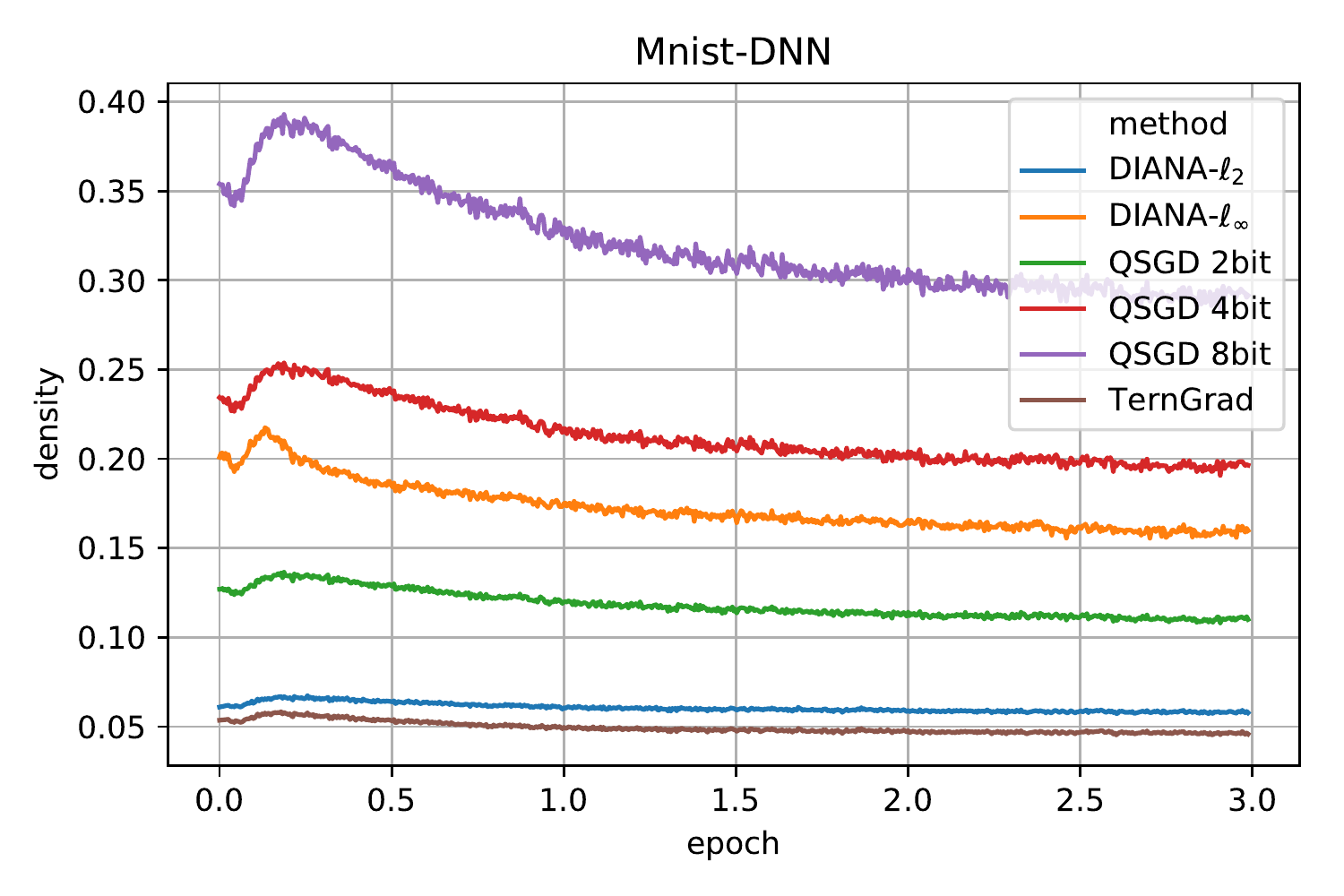}
\includegraphics[width=0.46\textwidth]{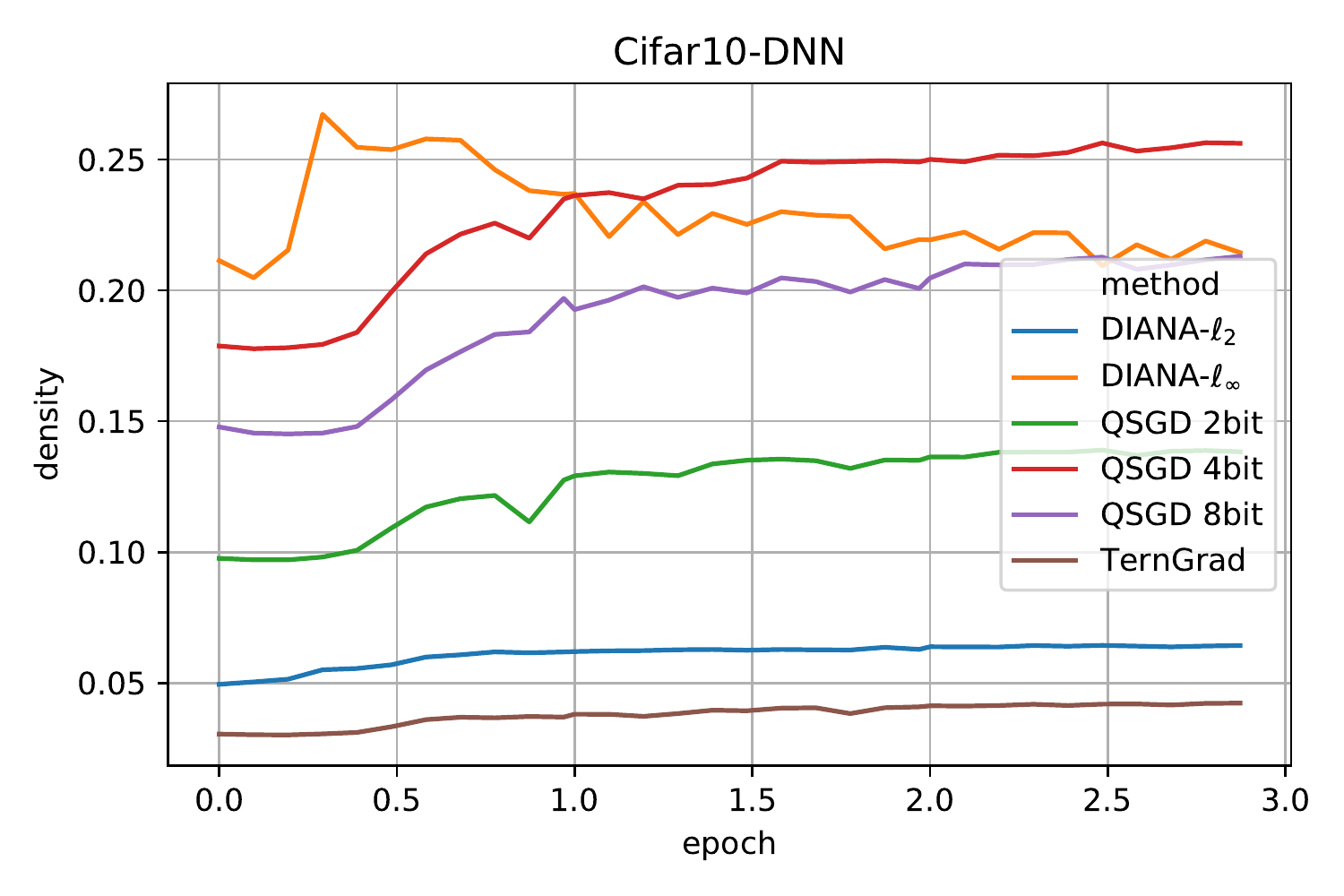}
\end{center}

\caption{Evolution of sparsity of the quantized gradient for 3 different problems and 3 algorithms.}
\label{fig:DNN:sparsity}

\end{figure}

\subsection{Computational Cost}

\begin{figure}[h!]
\centering

\includegraphics[scale=0.5]{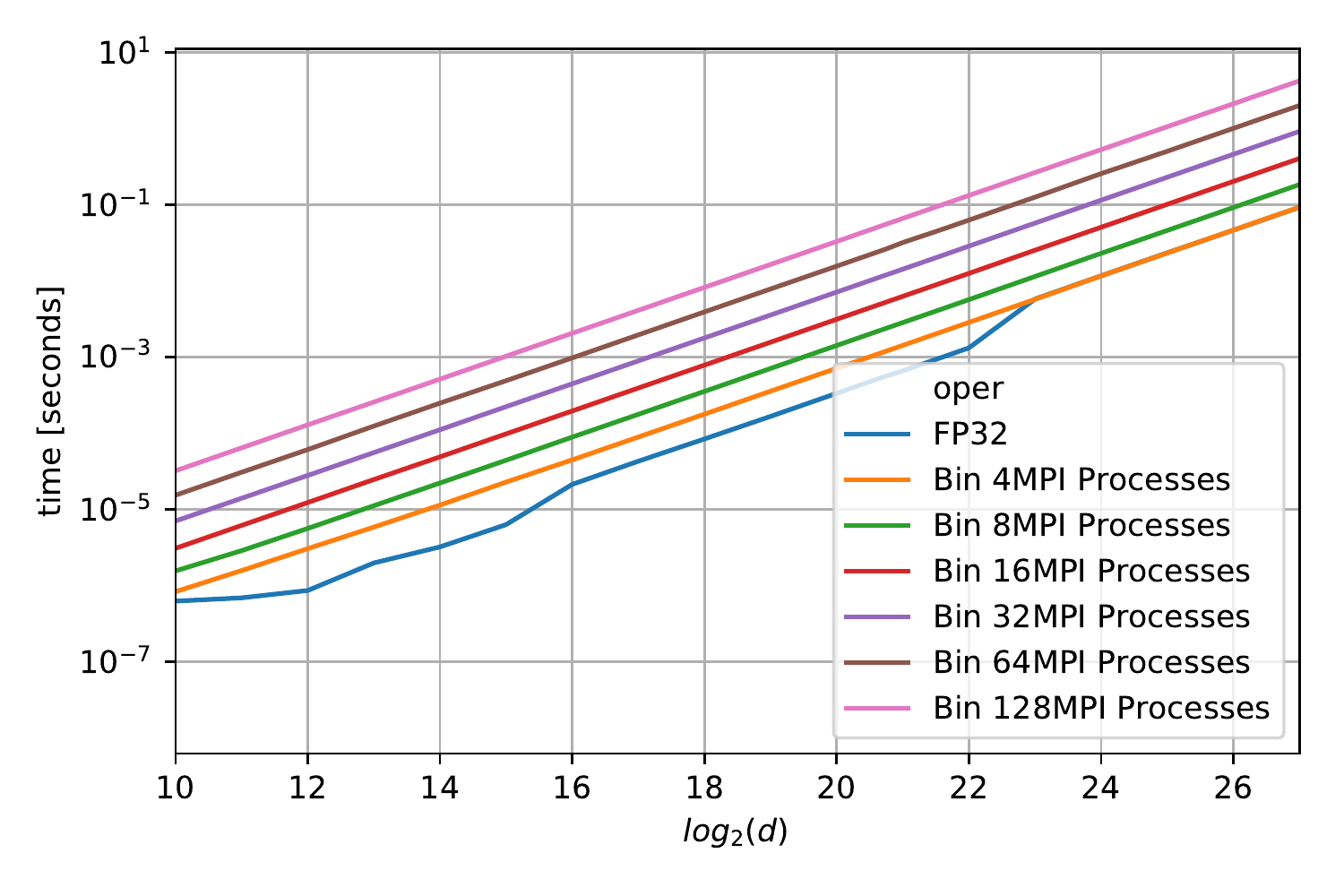}

\caption{Comparison of a time needed to update weights after a reduce vs.\
the time needed to update the weights when using a sparse update from \algname{DIANA} using 4-128 MPI processes and 10\% sparsity.
}

\label{fig:add}

\end{figure}

\chapter{Appendix for Chapter~\ref{chapter:sdm}}
\graphicspath{{sdm/}}

\section{Applications}\label{ap:applications}
In this section  we list a number of selected applications for our method:

\begin{itemize}
    \item Compressed sensing~\cite{candes2006robust}.
    \item Total Generalized Variance (TGV) for image denoising~\cite{bredies2010total}.
    \item Decentralized optimization over networks~\cite{Nedic2009distributed}.
    \item Support-vector machine~\cite{cortes1995support}.
    \item Dantzig selector~\cite{candes2007dantzig}.
    \item Group-Lasso~\cite{yuan2006model}.
	\item Network utility maximization.
    \item Square-root Lasso~\cite{belloni2011square}.
    \item $\ell_1$ trend filtering~\cite{kim2009ell_1}.
    \item Convex relaxation of unsupervised image matching and object discovery~\cite{vo2019unsupervised}.
\end{itemize}

In the rest of this section we formulate some  of them explicitly. A summary of the mapping of these problems to the structure of problem \eqref{eq:pb_general} is provided in Table~\ref{tbl:summary_apps}.

\begin{table}[!h]
\centering
{
\caption{Selected applications of Algorithm~\ref{alg:sdm} for solving problem \eqref{eq:pb_general}.}\label{tbl:summary_apps}
\begin{tabular}{cccc}
\toprule
Special case of problem \eqref{eq:pb_general} & $f(x)$ & $g_j(x)$ & $\psi(x)$ \\
\midrule 
Constrained optimization \eqref{eq:bifg9bd89d} & $f(x)$ & $\ind_{\cC_j}(x)$ & $\psi(x)$  \\
Convex projection & $\frac{1}{2}\|x-x^0\|^2$ & $\ind_{\cC_j}(x)$ & $0$  \\
Convex feasibility & $0$ & $\ind_{\cC_j}(x)$ & $0$  \\
Dantzig selector \eqref{eq:Dantigg9gf98f} & $0$ & $\chi_{\cB_{\lambda}^j}(x)$ & $\|x\|_1$   \\
Decentralized optimization~\eqref{eq:bid79f9hbfss} & $f_i(x_i)$ & $\ind_{\{x\;:\;w_j^\top x = 0\}}(x)$ & 0  \\
Support vector machine \eqref{eq:SVMbuis} & $f(x)=\frac{\lambda}{2} \|x\|^2$, $n=1$ & $\max\{0, 1 - b_j a_j^\top x\}$ & 0   \\
Overlapping group Lasso \eqref{eq:OGL} & $f_i(x)=\frac{1}{2}(a_i^\top x - b_i)^2$ & $\|x\|_{G_j}$ & 0   \\
Fused Lasso  \eqref{eq:FLi0hf88fh0} & $\frac{1}{2}(a_i^\top x - b_i)^2$ & $\ind_{\cC_j^\varepsilon}(x)$ & $\lambda\|x\|_1$  \\
Fused Lasso \eqref{eq:bifub98f0fss} & $\frac{1}{2}(a_i^\top x - b_i)^2$ & $\lambda_2 |\mD_{j:} x|$ & $\lambda_1\|x\|_1$  \\
\bottomrule 
\end{tabular}
}
\end{table}

\subsection{Constrained optimization} 
Let  $\cC_j\subseteq \RR^d$ be closed convex sets with a non-empty intersection and consider the constrained composite optimization problem
\begin{align*}
    \min & \quad f(x) + \psi(x)
    \qquad \text{subject to} \quad \quad x \in \cap_{j=1}^m \cC_j.
\end{align*}
If we let $g_j\equiv \ind_{\cC_j}$ be the characteristic function of $\cC_j$, defined as follows: $\ind_{\cC_j}(x)=0$ for $x\in \cC_j$ and $\ind_{\cC_j}(x) =+\infty$ for $x\notin \cC_j$, this problem can be written in the form 
\begin{equation}\label{eq:bifg9bd89d}
\min_{x\in \RR^d} f(x) +\psi(x) + \frac{1}{m}\sum_{j=1}^m \underbrace{\ind_{\cC_j}(x)}_{g_j(x)}.
\end{equation}


For $f(x)=\frac{1}{2}\|x-x^0\|^2$ and $\psi\equiv 0$, this specialized to the best approximation problem. For $f\equiv 0$  and $\psi\equiv 0$, this problem specializes to the convex feasibility problem.

\subsection{Dantzig selector} Dantzig selector~\cite{candes2007dantzig} solves the problem of estimating sparse parameter $x$ from a linear model. Given an  input matrix $\mA\in\RR^{m\times d}$, output vector $b\in\RR^m$ and threshold parameter $\lambda\geq 0$, define \[\cB_\lambda\eqdef \{x\;:\;\|\mA^\top(b - \mA x)\|_{\infty}\le \lambda \} = \bigcap_{j=1}^m \cB_{\lambda}^j,\] 
where $\cB^j_\lambda \eqdef \bigl\{x\;:\; \bigl| \left(\mA^\top(b - \mA x)\right)_j \bigr| \leq \lambda\bigr\}$. The goal of the Dantzig selector problem is to find the solution to
\begin{align*}
	\min_{x\in \RR^d} \|x\|_1 + \ind_{\cB_\lambda}(x),
\end{align*}
which can equivalently be written in the  finite-sum form
\begin{equation}\label{eq:Dantigg9gf98f}
	\min_{x\in \RR^d}  \underbrace{\|x\|_1}_{\psi(x)} + \avejm \underbrace{\ind_{\cB_\lambda^j}(x)}_{g_j(x)}.
	\end{equation}

\subsection{Decentralized optimization} The problem of minimizing the sum of functions over a network~\cite{Nedic2009distributed} can be reformulated as
\begin{align*}
	\min_{x=(x_1, \dotsc, x_n)} \frac{1}{n}\sumin f_i(x_i) + \ind_{\{x\;:\; \mW x = 0\}}(x),
\end{align*}
where $\mW$ is a matrix such that $\mW x =0$ if and only if $x_1=\dotsb=x_n$. Functions $f_1,\dotsc, f_n$ are stored on different nodes and each node has access only to its own function. Matrix $\mW$ is often derived from a communication graph, which defines how the nodes can communicate with each other. Formally, if $\mW = (w_1^\top, \dotsc, w_m^\top)^\top$, we rewrite the problem above as
\begin{equation}\label{eq:bid79f9hbfss}
	\min_{x=(x_1, \dotsc, x_n)} \frac{1}{n}\sumin \underbrace{f_i(x_i)}_{f_i(x)} + \avejm \underbrace{\ind_{\{x\;:\;w_j^\top x = 0\}}(x)}_{g_j(x)}.
\end{equation}

\subsection{Support-vector machine (SVM)} Support-vector machine~\cite{cortes1995support} is a very popular method for supervised classification. The primal formulation of SVM is given by
\begin{equation}\label{eq:SVMbuis}
	\min_{x\in \RR^d} \underbrace{\frac{\lambda}{2} \|x\|^2}_{f(x) } + \avejm \underbrace{\max\{0, 1 - b_j a_j^\top x\}}_{g_j(x)},
\end{equation}
where $a_1,\dotsc, a_m\in \RR^d$ and $b_1, \dotsc, b_m$ are the features and the outputs. It is easy to verify that for $g_j(x) = \max\{0, 1 - b_j a_j^\top x\}$ the proximal operator is given by
\begin{align*}
	\proxj(x) = x + \Pi_{[0, \eta_j]}\biggl(\frac{1 - b_j a_j^\top x}{\|a_j\|^2} \biggr) b_j a_j.
\end{align*}

The celebrated stochastic subgradient descent method \algname{Pegasos} \cite{Pegasos, Pegasos-MAPR, Pegasos2} for SVMs achieves only slow $\cO(\frac{1}{K})$ rate.

\subsection{Overlapping group Lasso} 
This is a generalization of Lasso proposed in~\cite{yuan2006model} to efficiently select groups of features that are most valuable for the given objective. Let us assume that we are given sets of indices $G_1, \dotsc, G_m \subseteq	\{1,\dotsc, d\}$ and let $\|x\|_{G_j} \eqdef \sqrt{\sum_{i\in G} [x]_i^2}$, where $[x]_i$ is the $i$-th coordinate of vector $x$. Then, assuming that we are given vectors $a_1, \dotsc, a_n\in \RR^d$ and scalars $b_1, \dotsc, b_n$, the objective we want to minimize is 
\begin{equation}\label{eq:OGL}
	\min_{x\in \RR^d}  \frac{1}{n}\sumin \underbrace{\frac{1}{2}(a_i^\top x - b_i)^2}_{f_j(x)} + \avejm \underbrace{\|x\|_{G_j}}_{g_j(x)}.
\end{equation}
It is easy to verify that if $g_j(x) = \|x\|_{G_j}$, then
\begin{align*}
	[\prox_{\eta_j g_j}(x)]_i
	=\begin{cases}
	[x]_i, & \text{if } i\not\in G_j, \\
	\max\left\{0, \left(1 - \frac{\eta_j}{\|x\|_{G_j}}\right) \right\}[x]_i, & \text{if } i\in G_j.
	 \end{cases}
\end{align*}
Vector $y_j^k$ will always have at most $|G_j|$ nonzeros, so one can store in memory only the coordinates of $y_j^k$ from $G_j$.

\subsection{Fused Lasso} The Fused Lasso problem~\cite{tibshirani2005sparsity} is defined as
\begin{equation}\label{eq:FLi0hf88fh0}
	\min_{x\in \RR^d} \frac{1}{n}\sumin \underbrace{\frac{1}{2}(a_i^\top x - b_i)^2}_{f_i(x)} + \underbrace{\lambda\|x\|_1}_{\psi(x)} + \frac{1}{d-1}\sum_{j=1}^{d-1} \underbrace{\ind_{\cC_j^\varepsilon}(x)}_{g_j(x)},
\end{equation}
where $\cC_j^\varepsilon \eqdef \left\{x\;:\; \left| [x]_j - [x]_{j+1} \right| \le \varepsilon \right\},$
$[x]_j$ is the $j$-th entry of vector $x$, $a_1,\dotsc, a_n\in \RR^d$ and $b_1, \dotsc, b_n\in\RR$ are given vectors and scalars, $\varepsilon$ is given thresholding parameter.

Another formulation of the Fused Lasso is done by using penalty functions. Define $\mD$ to be zero everywhere except for $\mD_{i,i}=1$ and $\mD_{i, i+1}=-1$ with $i=1,\dotsc, d-1$. Note that $\|\mD x\|_1 = \sum_{j=1}^m |\mD_{j:} x|$, where $m$ is the number of rows of $\mD$. Then the reformulated objective is
\begin{equation}\label{eq:bifub98f0fss}
	\min_x \frac{1}{n}\sumin \underbrace{\frac{1}{2}(a_i^\top x - b_i)^2}_{f_i(x)} + \underbrace{\lambda_1\|x\|_1}_{\psi(x)} + \avejm \underbrace{\lambda_2 |\mD_{j:} x|}_{g_j(x)}.
\end{equation}
In our notation, this means $\mA = \mD^\top$ and $\mA^\top \mA$ is a tridiagonal matrix given by
\begin{align*}
	\mA^\top \mA
	=
	\begin{pmatrix}
		2   & -1  \\
		-1 &    2 & -1  \\
		     &  -1 & 2   & -1 & \\
		     &       &      & \ddots    & -1 \\
		     &       &      &           -1  & 2
	\end{pmatrix}.
\end{align*} 
Let $\mW$ be a tridiagonal matrix of size $(d-1)\times (d-1)$ with $a$ on its main diagonal and $b$ on the other two diagonals. It can be shown that its eigenvalues  are given by $\lambda_k(\mW) = a+2|b|\cos \left(\frac{k\pi}{d} \right)$, $k=1,\dotsc, d-1$. Thus, $\lambda_{\min}(\mA^\top \mA) = 2 + 2\cos\left(\left(1 - \frac{1}{d} \right)\pi\right) = 2 - 2\cos\left(\frac{\pi}{d}\right) \approx \frac{1}{2d^2}$ and $\min_j \frac{1}{\|\mA_j\|^2} = \frac{1}{6}$. Therefore, if in~\eqref{eq:FLi0hf88fh0} or~\eqref{eq:bifub98f0fss} $\lambda_1=0$, we guarantee linear convergence with the aforementioned constants.

\subsection{Square-root Lasso} The approach gets its name from minimizing the square root of the regular least squares, i.e., $\|\mD w - b\|$ instead of $\|\mD w - b\|^2$. This is then combined with $\ell_1$-penalty for feature selection, which gives the objective
\begin{align*}
		\min_{w\in \RR^d} \|\mD w - b\| + \lambda \|w\|_1.
\end{align*}
Equivalently, by introducing a new variable $z$ we can put constraints $\mD_{j:} x - [z]_j =0$ for $j=1,\dotsc, m$, which can be written as $a_j^\top(w^\top, z^\top)^\top=0$ with $a_j = (\mD_{j:}, e_j^\top)^\top$ and $e_j\eqdef (0, 0, \dotsc, \underbrace{1}_{j}, \dotsc, 0)$. Then, the reformulation is
\begin{align*}
	\min_{x=(w, z)\in \RR^{d+m}} \avejm \underbrace{\ind_{\{x:a_j^\top x =0\}}}_{g_j(x)=g_j(w, z)} + \underbrace{\|z - b\| + \lambda \|w\|_1}_{\psi(x)=\psi(w,z)}.
\end{align*}
The proximal operator of $\psi$ is that of a block-separable function, which is easy to evaluate:
\begin{align*}
	\prox_{\gamma \psi}(x) 
	= \begin{pmatrix}\prox_{\gamma\lambda\|\cdot\|_1}(w) \\ \prox_{\gamma \|\cdot - b\|}(z)\end{pmatrix}.
\end{align*}

\section{Relation to Existing Methods}

\subsection{\algname{SDCA}, \algname{Dykstra's algorithm} and the \algname{Kaczmarz method}}
Here we formulate \algname{SDCA}~\cite{shalev2013stochastic}, \algname{Dykstra's algorithm} and \algname{Kaczmarz method}. \algname{SDCA} is a method for solving
\begin{align*}
	\min_{x\in\RR^d} \avejm g_j(x)+ \frac{1}{2}\|x - x^0\|^2.
\end{align*}
If $j$ is sampled uniformly from $\{1,\dotsc, m\}$, \algname{SDCA} iterates can be defined by the following recursion,
\begin{align*}
	x^{k+1} 
	&= \prox_{\eta g_j}(x^k + \overline y_j^k), \\
	\overline y_j^{k+1}
	&= \overline y_j^k + x^k - x^{k+1},
\end{align*}
If we restrict our attention to characteristic functions, i.e., 
\begin{align*}
	g_j(x) = \ind_{\cC_j}(x) = \begin{cases}0, &\text{if } x\in \cC_j \\ +\infty, & \text{otherwise} \end{cases},
\end{align*}
then the proximal operator step is replaced with projection:
\begin{align*}
	x^{k+1} 
	&= \Pi_{\cC_j}(x^k + \overline y_j^k).
\end{align*}
This is known as \algname{Dykstra's algorithm}. Finally, if $\cC_j= \{x: a_j^\top x=b_j\}$, then it boils down to random projections, i.e.,
\begin{align*}
	x^{k+1} = \Pi_{\{a_j^\top x = b_j\}} (x^k),
\end{align*}
which is the method of \algname{Kaczmarz}.

\begin{theorem}\label{th:kaczmarz}
	Consider the regularized minimization problem of \algname{SDCA}, which is 
	\begin{align*}
		\min_x \avejm g_j(x)	 + \frac{1}{2}\|x - x^0\|^2
	\end{align*}		
	with convex $g_1, \dotsc, g_m$. Then, \algname{SDCA} is a special cases of Algorithm~\ref{alg:sdm} obtained by applying it with $f(x)=\frac{1}{2}\|x - x^0\|^2$, $\psi(x)\equiv 0$, stepsize $\gamma= \frac{1}{m}$ and initialization $y_1^0=\dotsb=y_m^0 =0$. Furthermore, if we consider special case $g_j = \ind_{\cC_j}$, where $\cC_j\neq \emptyset$ is a closed convex set, then we also obtain \algname{Dykstra's algorithm}, and if every $\cC_j$ is a linear subspace, then we recover the \algname{Kaczmarz method}.
\end{theorem}
\begin{proof}
	Consider the iterates of \algname{SDM}. We will show by induction that $y^k = x^0-x^k$ and $x^{k+1} = \proxj(x^k + \eta_j y_j^k)$.
	Indeed, it holds for $y^0$ by initialization, and then by induction assumption we have
	\begin{align*}
		z^k
		= x^k - \gamma(x^k - x^0) - \gamma y^k
		= x^k - \gamma(x^k - x^0) - \gamma (x^0 - x^k)
		= x^k.
	\end{align*}
	Therefore, if we denote $\overline y_j^k \eqdef \eta_j y_j^k$, then
	\begin{align*}
		x^{k+1} 
		=\proxj(x^k + \overline y_j^k),
	\end{align*}
	which is the update rule of $x^{k+1}$ in \algname{SDCA}. Moreover, we have
	\begin{align*}
		\overline y_j^{k+1} 
		= \gamma y_j^{k+1}
		=\gamma y_j^k + x^k - x^{k+1}
		= \overline y_j^k + x^k - x^{k+1}.
	\end{align*}
	Finally, by induction assumption it holds $y^k=x^0 - x^k$, whence
	\begin{align*}
		y^{k+1}
		= y^k + \frac{1}{m}(y_j^{k+1} - y_j^{k})
		= y^k + z^k - x^{k+1}
		= x^0 - x^k + x^k - x^{k+1}
		= x^0 - x^{k+1},
	\end{align*}
	which yields our induction step and the proof itself.
\end{proof}
\subsection{\algname{Accelerated Kaczmarz}}\label{ap:kaczmarz}
\algname{Accelerated Kaczmarz}~\cite{liu2016accelerated} performs the following updates:
\begin{align*}
	z^k
	&= (1 - \alpha_k) x^k - \alpha_k y^k,\\
	x^{k+1}
	&= \Pi_{\{x: a_i^\top x = b_i\}} (z^k), \\
	y^{k+1}
	&= y^k + \nu_k (z^k - x^{k+1}) + (1 - \beta_k) (z^k - y^k)
\end{align*}
with some parameters $\alpha_k, \nu_k, \beta_k$. While the original analysis~\cite{liu2016accelerated} suggests $\beta_k<1$, our method gives the same update when $f(x)=\frac{1}{2}\|x\|^2$, $\psi\equiv 0$, $\alpha_k=\gamma$, $\beta_k=1$, $\nu_k=\frac{1}{\gamma n}$. 
\subsection{\algname{ADMM} and \algname{Douglas--Rachford Splitting}}
\algname{ADMM}, also known as \algname{Douglas--Rachford Splitting}, in its simplest form as presented in~\cite{parikh2014proximal} is a special case of Algorithm~\ref{alg:sdm} when $f\equiv 0$ and $m=1$.
\subsection{\algname{Point--SAGA}, \algname{SAGA}, \algname{SVRG} and \algname{Proximal GD}}
In the trivial case $f\equiv 0$ and $\psi\equiv 0$, we recover \algname{Point--SAGA}. Methods such as \algname{SAGA}, \algname{SVRG} and \algname{Proximal Gradient Descent} are obtained, in contrast, by setting $g\equiv 0$. We would like to mention that introducing $g$ does not change the stepsizes for which those methods work, e.g., \algname{Gradient Descent} works with arbitrary $\gamma < \frac{2}{L}$, which is tight. The similarity suggests that small $\gamma$ should be used when solving this problem and this observation is validated by our experiments.

\subsection{\algname{Stochastic Primal--Dual Hybrid Gradient}}\label{ap:spdhg}
The relation to the \algname{Stochastic Primal--Dual Hybrid Gradient} (\algname{SPDHG}) is complicated. On the one hand, \algname{SPDHG} is a general method with three parameters and it preconditions proximal operators with matrices, so our method cannot be its strict generalization. On the other hand, \algname{SPDHG} does not allow for $f$. Moreover, when $f\equiv 0$ and some parameters are set to specific values in \algname{SPDHG}, the methods coincide, but the guarantees are not the same. In particular, we show below that one of the parameters in \algname{SPDHG}, $\theta$, should be set to 1, in which case linear convergence for smooth $g_1, \dotsc, g_m$ was not known for \algname{SPDHG}. Therefore, the tools developed in this chapter can potentially lead to new discoveries about full version of \algname{SPDHG} as well.

Let us now formulate the method explicitly. After a simple rescaling of the functions, \algname{SPDHG} from~\cite{ehrhardt2017faster} can be formulated as a method to solve the problem
\begin{align}
	\min_{x\in\RR^d} \avejm \phi_j(\mA_j^\top x) + \psi(x).\label{pb:spdhg}
\end{align}
Renaming the variables for our convenience and choosing for simplicity uniform probabilities of sampling $j$ from $\{1,\dotsc, m\}$, the update rules of \algname{SPDHG} can be written as
\begin{align*}
	w^k 
	&= \proxR(w^{k-1} - \gamma \overline y^k), \\
	y_j^{k+1} 
	&= \prox_{\sigma \phi_j^\star}(\sigma\mA_j^\top w^k + y_j^k), \\
	y^{k+1} 
	&= y^k + \frac{1}{m}\mA_j(y_j^{k+1} - y_j^k), \\
	\overline y^{k+1} 
	&= y^k + \theta \mA_j(y_j^{k+1} - y_j^k),
\end{align*}
where $\gamma, \sigma$ and $\theta$ are the method's parameters and $\phi_j^\star$ is the Fenchel conjugate of $\phi_j$. The initialization that we are interested in is with $y^0 = \avejm y_j^0$, $\overline y^0 = y^0$, $w^0 = x^0$.

One can immediately see that one big difference with our approach is that the method puts $\mA_j$ outside of the proximal operator, which also leads to different iteration complexity. In particular, when $\phi_1,\dotsc, \phi_m$ are smooth, the complexity proved in~\cite{chambolle2018stochastic} is 
\begin{align*}
	\cO\left(\left(m + \sumjm \|\mA_j\|\sqrt{\frac{L_\phi}{\mu_\psi}} \right)\log \frac{1}{\varepsilon}\right),
\end{align*}
where $\mu_{\psi}$ is the strong convexity constant of $\psi$ and $L_\phi$ is the smoothness constant of $\phi_1,\dotsc, \phi_m$. Since function $g_j(x) = \phi_j(\mA_j^\top x)$ is at most $L_\phi \|\mA_j\|^2$ smooth, our rate from Corollary~\ref{cor:acc_in_g} with $\mu$-strongly convex and $L$-smooth $f$ is 
\begin{align*}
	\cO\left(\left(n+m + \frac{L}{\mu} +  \sqrt{m\frac{L_\phi}{\mu}}\max_j \|\mA_j\| \right)\log \frac{1}{\varepsilon}\right).
\end{align*}
If, in addition, we use sampling with probabilities proportional to $\|\mA_j\|$, then we can achieve
\begin{align*}
	\cO\left(\left(n+m + \frac{L}{\mu} +  \frac{1}{\sqrt{m}}\sumjm \|\mA_j\|\sqrt{\frac{L_\phi}{\mu}} \right)\log \frac{1}{\varepsilon}\right).
\end{align*}

We do not prove this, but the complexity for our method will be similar if we use strongly convex $\psi$ rather than $f$, so our rates should match or be even be superior to that of \algname{SPDHG}, at the cost of evaluating potentially harder proximal operators.

Now, let us prove that our method is indeed connected to \algname{SPDHG} via choice of $\theta=1$ and $\gamma\sigma = 1$.
\begin{theorem}
	If we apply \algname{SPDHG} with identity matrices $\mA_j = \mI$, i.e., $\phi_j(x) = g_j(x)$, and choose parameters $\theta=1$ and $\gamma\sigma = 1$, then it is algorithmically equivalent to Algorithm~\ref{alg:sdm} with $f\equiv 0$.
\end{theorem}
\begin{proof}
	Since $\phi_j$ and $g_j$ are the same, we will use in the proof $g_j$ only. 
	
	First, mention that it is straightforward to show by induction that $y^k = \avejm y_j^k$, which coincides with our update. Our goal is to show by induction that in \algname{SPDHG} it holds
	\begin{align*}
		w^{k-1} - \gamma \overline y^k
		= x^k - \gamma y^k,
	\end{align*}
	where we define sequence $x^k$ as 
	\begin{align*}
		x^{k+1} 
		\eqdef \prox_{\gamma g_j}(w^k + \gamma y_j^k)
		= \prox_{\frac{1}{\sigma} g_j}(w^k + \gamma y_j^k).
	\end{align*}
	We will see that implicitly $x^{k+1}$ is present in every update of \algname{SPDHG}. To this end, let us first rewrite the update for $y_j^{k+1}$. We have by Moreau's identity
	\begin{align*}
		y_j^{k+1}
		= \prox_{\sigma g_j^\star}(\sigma w^k + y_j^k)
		= \sigma w^k + y_j^k - \sigma\prox_{\frac{1}{\sigma} g_j}\left(\frac{\sigma w^k + y_j^k}{\sigma} \right).
	\end{align*}
	Since we consider $\sigma= \frac{1}{\gamma}$, it transforms into
	\begin{align*}
		y_j^{k+1}
		= y_j^k  + \frac{1}{\gamma}\left( w^k - \prox_{\eta g_j}(w^k + \gamma y_j^k) \right) 
		&= y_j^k  + \frac{1}{\gamma}\left( w^k - x^{k+1} \right)
	\end{align*}
	The only missing thing is rewriting update for $w^k$ in terms of $x^k$ and $y^k$. From the update rule for $y_j^{k+1}$ we derive
	\begin{align*}
		\overline y^{k+1}
		= y^k + \theta (y_j^{k+1} - y_j^k)
		= y^k + \frac{\theta }{\gamma}(w^k - x^{k+1}).
	\end{align*}
	Hence,
	\begin{align*}
		w^{k+1} 
		&= \proxR(w^k - \gamma \overline y^{k+1})
		= \proxR(w^k - \gamma  y^{k+1} - \theta (w^k - x^{k+1}))\\
		&\overset{\theta=1}{=} \proxR(x^{k+1} - \gamma  y^{k+1}).
	\end{align*}
	Thus, updates for $w^k$, $y_j^k$ and $y^k$ completely coincide under this choice of parameters.
\end{proof}
Since our method under $f\equiv 0$ reduces to \algname{Point--SAGA}, we obtain the following result that was unknown.
\begin{corollary}
	\algname{Point--SAGA}~\cite{defazio2016simple} is a special case of \algname{Stochastic Primal--Dual Hybrid Gradient}~\cite{chambolle2018stochastic}.
\end{corollary}

\section{Evaluating Proximal Operators}
For some functions, the proximal operator admits a closed form solution, for instance if $g_j( x) = \ind_{\{x\;:\; a_j^\top x=b_j\}}(x)$, then
\begin{align*}
	\proxj(x)
	= x - \frac{a_j^\top x - b_j}{\|a_j\|^2}a_j.
\end{align*}
If, however, the proximal operator is not given in a closed form, then it is still possible to efficiently evaluate it. If $g_j=\phi_j(\mA_j^\top x)$, $\mA_j\in \RR^{d\times d_j}$, then the proximal operator is the solution of a $d_j$-dimensional strongly convex problem.
\begin{lemma}\label{lem:prox_of_composition}
	Let $\phi_j \colon \RR^{d_j}\to \RR$ be a convex lower semi-continuous function such that $\Range{\mA_j^\top}$ has a point of $\dom(\phi)$. If $g_j(x) = \phi_j(\mA_j^\top x)$, then
	\begin{align*}
		x - \prox_{\eta_j g_j}(x)
		\in \Range{\mA_j}.
	\end{align*}
\end{lemma}
\begin{proof}
	Let us fix $x$. Any vector $z\in\RR^d$ can be decomposed as $z = x + \mA_j\beta + w$, where $\beta\in\RR^{d_j}$ and $\mA_j^\top w = 0$, from which it also follows $g_j(z) = \phi_j(\mA_j^\top x+ \mA_j^\top \mA_j \beta)$. Then
	\begin{align*}
		\prox_{\eta_j g_j}(x)
		&\eqdef \argmin_{z\in\RR^d}\left\{\eta_j \phi_j(\mA_j^\top z) + \frac{1}{2}\|z - x\|^2 \right\}\\
		&= \argmin_{z=x+\mA_j\beta + w}\left\{\eta_j \phi_j(\mA_j^\top x + \mA_j^\top \mA_j \beta) + \frac{1}{2}\|\mA_j\beta + w\|^2 \right\} \\
		&= \argmin_{z=x+\mA_j\beta + w}\left\{\eta_j \phi_j(\mA_j^\top x + \mA_j^\top \mA_j \beta) + \frac{1}{2}\|\mA_j\beta\|^2 + \frac{1}{2}\| w\|^2 \right\} .
	\end{align*}
	Clearly, the last expression achieves its minimum only when $w=0$.
\end{proof}

We can simplify the expression for the proximal operator even further if $\mA_j$ is of full column rank, for instance if it is just a single nonzero row. It is straightforward to verify that for any matrix $\mB\in\RR^{d_1\times d_2}$, constant vector $c\in\RR^{d_2}$ and function $\Phi$ with a unique minimizer and $\dom(\Phi(\mB \beta +c))\neq \emptyset$ it holds
\begin{align}
	\argmin_{\beta=\mB(\alpha + c), \beta\in\RR^{d_2}}\Phi(\beta)
	&= \argmin_{\beta=\mB(\alpha + c), \beta\in\RR^{d_2}} \Phi(\mB(\alpha + c)) \notag \\
	&= \mB\argmin_{u=\alpha + c, \alpha\in\RR^{d_1}} \Phi(\alpha + c) \notag \\
	&= \mB\Bigl(\argmin_{u\in\RR^{d_1}} \Phi( u ) - c\Bigr). \label{eq:argmin_rule}
\end{align}
Since we know by Lemma~\ref{lem:prox_of_composition} that $u\eqdef \proxj(x) = x + \mA_j \beta_j$ for some $\beta_j\in\RR^{d_j}$, we can write the necessary and sufficient optimality condition for $u$ by repeatedly applying~\eqref{eq:argmin_rule}
\begin{align*}
	&\proxj(x) \\
	&= \argmin_{u = x + \mA_j\beta,\; \beta\in\RR^{d_j}} \left\{\phi_j(\mA_j^\top u) + \frac{1}{2\eta_j}\|x - u\|^2 \right\}\\
	&\overset{\eqref{eq:argmin_rule}}{=} x + \mA_j \argmin_{\beta\in\RR^{d_j}} \left\{\phi_j\left(\mA_j^\top (x + \mA_j\beta)\right) + \frac{1}{2\eta_j} \|\mA_j\beta\|^2\right\} \\
	&=x + \mA_j \argmin_{\beta\in\RR^{d_j}} \left\{\phi_j\left(\mA_j^\top x + \mA_j^\top\mA_j\beta\right) + \frac{1}{2\eta_j} \|\mA_j(\mA_j^\top\mA_j)^{-1}\mA_j^\top\mA_j\beta\|^2\right\} \\
	&\overset{\eqref{eq:argmin_rule}}{=}x + \mA_j (\mA_j^\top\mA_j)^{-1}\argmin_{\alpha=\mA_j^\top\mA_j\beta} \left\{\phi_j\left(\mA_j^\top x + \alpha\right) + \frac{1}{2\eta_j} \|\mA_j(\mA_j^\top\mA_j)^{-1}\alpha\|^2\right\} \\
	&\overset{\eqref{eq:argmin_rule}}{=}x + \mA_j (\mA_j^\top\mA_j)^{-1}\Bigl(\argmin_{\theta=\alpha + \mA_j^\top x} \left\{\phi_j\left(\theta\right) + \frac{1}{2\eta_j} \|\mA_j(\mA_j^\top\mA_j)^{-1}(\theta - \mA_j^\top x)\|^2\right\} - \mA_j^\top x\Bigr).
\end{align*}
Note that
\begin{align*}
	\|\mA_j(\mA_j^\top\mA_j)^{-1}(\theta - \mA_j^\top x)\|^2 
	&= (\theta - \mA_j^\top x)^\top (\mA_j^\top\mA_j)^{-1}\mA_j^\top\mA_j(\mA_j^\top\mA_j)^{-1}(\theta - \mA_j^\top x) \\
	&= \|\theta - \mA_j^\top x\|^2_{(\mA_j^\top\mA_j)^{-1}},
\end{align*}
where for any positive semi-definite matrix $\mW$ we denote $\|x\|_{\mW}^2\eqdef x^\top \mW x$. Denoting similarly $\prox^{\mW}_{\eta_j \phi_j}(x)\eqdef \argmin_{\theta}\{\phi_j(\theta) + \frac{1}{2\eta_j}\|\theta - x\|_\mW^2\}$, we obtain
\begin{align*}
	\proxj(x) 
	&=x + \mA_j (\mA_j^\top\mA_j)^{-1}\left( \argmin_{\theta\in\RR^{d_j}} \left\{\phi_j\left(\beta\right) + \frac{1}{2\eta_j} \|\theta - \mA_j^\top x\|_{(\mA_j^\top\mA_j)^{-1}}\right\}  - \mA_j^\top x \right) \\
	&=x + \mA_j (\mA_j^\top\mA_j)^{-1}\left( \prox_{\eta_j\phi_j}^{(\mA_j^\top\mA_j)^{-1}}\left(\mA_j^\top x\right)  - \mA_j^\top x \right).
\end{align*}
Thus, we only need to know how to efficiently evaluate $\prox_{\lambda \phi_j}^{(\mA_j^\top\mA_j)^{-1}}(z)$ for arbitrary $\lambda>0$ and $z\in\RR^{d_j}$, assuming that matrix $(\mA_j^\top\mA_j)^{-1}$ can be precomputed. For example, if $\mA_j=a_j\in\RR^d$, then
\begin{align*}
	\prox_{\eta_j \phi_j}^{(a_j^\top a_j)^{-1}}(x) = \prox_{\eta_j \|a_j\|^2 \phi_j}(x).
\end{align*}

If, in addition, $\phi_j\colon \RR \to \RR$ is given by
\begin{align*}
	\phi_j(z) = \begin{cases} b_jz, & \text{if } z\le 0,\\ c_j z, & \text{otherwise} \end{cases}
\end{align*}
with some $b_j, c_j\in\RR$, $b_j<c_j$, then $\prox_{\lambda \phi_j}(z) = z - \lambda b_j $ for $z\le\lambda b_j$, $\prox_{\lambda \phi_j}(z)=0$ for $z\in(\lambda b_j, \lambda c_j]$ and $\prox_{\lambda \phi_j}(z)=z - \lambda c_j$ for $z> \lambda c_j$. Therefore,
\begin{align*}
	\prox_{\eta_j g_j}(x)
	= \begin{cases} x - \eta_j a_j b_j, &\text{if } a_j^\top x\le \|a_j\|^2 b_j,\\ x - \frac{a_j^\top x}{\|a_j\|^2}a_j, & \text{if } \|a_j\|^2b_j\le a_j^\top x \le \|a_j\|^2c_j , \\ x - \eta_j a_j c_j, &\text{otherwise} \end{cases}.
\end{align*}
Note that if $\|a_j\|^2b_j\le a_j^\top x \le \|a_j\|^2c_j $, then $a_j^\top \prox_{\eta_j g_j}(x) =0 $.

\section{Optimality Conditions} \label{sec:Opt_Cond}

We now comment on the nature of Assumption~\ref{as:optimality}.
In view of the first-order necessary and sufficient condition for the solution of~\eqref{eq:pb_general}, we have
\begin{align*}
    x^\star \in \cX^\star  \quad \Leftrightarrow \quad  0 \in \partial P(x^\star) = \nabla f(x^\star) + \partial (g + \psi)(x^\star).
\end{align*}
By the weak sum rule \cite[Corollary 3.38]{Beck2017}, we have \[\partial P(x)  \supseteq  \nabla f(x) + \frac{1}{m}\sum_{j=1}^m\partial g_j(x) +\partial  \psi(x)\] for all $x\in \dom(P) \supseteq \cX^\star$. Under the regularity condition $\cap_{j=1}^k(\dom(g_j)) \cap_{j=k+1}^m {\rm ri} (\dom (g_j))\cap {\rm ri}( \dom (\psi)) \neq \emptyset$, where $g_1,\dotsc, g_k$ are polyhedral functions, the inclusions becomes an identity \cite[Theorem 23.8]{rockafellar2015convex}, which means that Assumption~\ref{as:optimality} is satisfied.

For functions $g_j$ of the form $g_j(x) = \phi_j(\mA_j^\top x)$, where $\phi_j:\RR^{d_j}\to \RR\cup \{+\infty\}$ are proper closed convex functions and $\mA_j\in \RR^{d\times d_j}$, we shall instead consider the following (slightly stronger) assumption:
\begin{assumption}\label{as:optimality2}   There exists $x^\star\in \cX^\star$ and vectors $y_1^\star\in \mA_1 \partial \phi_1( \mA_1^\top x^\star), \dots, y_m^\star\in  \mA_m \partial \phi_m( \mA_m^\top x^\star)$ and $r^\star\in \partial \psi(x^\star)$ such that $        \nabla f(x^\star) + \avejm y_j^\star + r^\star = 0.$
\end{assumption}

Since $\mA_j \partial \phi_j( \mA_j^\top x) \subseteq \partial g_j(x)$ for all $x\in \dom(g_j)$  \cite[Theorem 3.43]{Beck2017}, Assumption~\ref{as:optimality2}  is indeed stronger than Assumption~\ref{as:optimality}. 
If $\Range{\mA_j^\top}$ contains a point from ${\rm ri}(\dom(g_j))$, or $g_j$ is polyhedral and $\Range{\mA_j^\top}$ contains a point from mere $\dom (g_j)$, then $\partial g_j(x) = \mA_j\partial \phi_j(\mA_j^\top x)$ for any $x$ \cite[Theorem 23.9]{rockafellar2015convex}, and these two assumptions are the same.

Below we provide another stationarity condition that shows why $x^\star$ is a fixed-point of our method.

\begin{lemma}[Optimality conditions]\label{lem:optimality}
    Let $x^\star$ be a solution of~\eqref{eq:pb_general} and let Assumption~\ref{as:optimality} be satisfied. Then for any $\gamma,\eta_j\in \RR$,
    \begin{align*}
        x^\star &= \proxR(x^\star - \gamma\nabla f(x^\star) - \gamma y^\star), \qquad x^\star = \proxj(x^\star + \eta_j y_j^\star).
    \end{align*}    
\end{lemma}
\begin{proof}
    Let 
    \[
    		z = \proxR(x^\star - \gamma\nabla f(x^\star) - \gamma y^\star) = \argmin_u \{\gamma \psi(u) + \frac{1}{2}\|u - (x^\star - \gamma\nabla f(x^\star) - \gamma y^\star)\|^2 \}.
    	\]
    	 Since $\psi$ is convex, the problem inside $\argmin$ is strongly convex, and the necessary and sufficient condition for $z$ to be its solution is
    \begin{align*}
        0\in z - x^\star + \gamma\nabla f(x^\star) + \gamma y^\star + \gamma \partial \psi(z).
    \end{align*}
    By Assumption~\ref{as:optimality} it holds for $z=x^\star$, implying the first equation that we want to prove. The second one follows by exactly the same argument applied to $\argmin_u \{\eta_j g_j(u) + \frac{1}{2}\|u - (x^\star + \eta_j y_j^\star)\|^2 \}$.
\end{proof}


\section{Convergence Proofs}
In this section, we provide the proofs of our convergence results. Each lemma, theorem and corollary is first restated and only then proved to simplify the reading.

\subsection{Proof of Lemma~\ref{lem:gd} (\algname{Gradient Descent})}\label{ap:gd}

Here we prove that \algname{Gradient Descent} update on $f$ satisfies our assumption on the method with the best possible stepsizes.
\begin{lemma}[Same as Lemma~\ref{lem:gd}]
	If $f$ is convex, \algname{Gradient Descent} satisfies Assumption~\ref{as:method}(a) with any $\gamma_{\max} < \frac{2}{L}$, $\omega = 2 - \gamma_{\max} L$ and $\cM^k = 0$. If $f$ is $\mu$-strongly convex, \algname{Gradient Descent} satisfies Assumption~\ref{as:method}(b) with $\gamma_{\max} = \frac{2}{L + \mu}$, $\omega = 1$ and $\cM^k=0$.
\end{lemma}
\begin{proof}
	Since we consider \algname{Gradient Descent}, we have
	\begin{align*}
		w^k
		&= x^k - \gamma \nabla f(x^k).
	\end{align*}
	First, if $f$ is convex and smooth, then for any $\gamma \le \gamma_{\max} < \frac{2}{L}$
	\begin{align*}
		\|w^k - w^\star\|^2
		&= \|x^k - x^\star\|^2 - 2 \gamma\<\nabla f(x^k) - \nabla f(x^\star), x^k - x^\star> + \gamma^2\|\nabla f(x^k) - \nabla f(x^\star)\|^2 \\
		&\le \|x^k - x^\star\|^2 - \gamma(2 - \gamma_{\max} L)\<\nabla f(x^k) - \nabla f(x^\star), x^k - x^\star> \\
		&\quad -  \gamma \gamma_{\max} L\<\nabla f(x^k) - \nabla f(x^\star), x^k - x^\star> + \gamma\gamma_{\max}\|\nabla f(x^k) - \nabla f(x^\star)\|^2 \\
		&\overset{\eqref{eq:grad_dif_scalar_prod}}{\le} \|x^k - x^\star\|^2 - \gamma(2 - \gamma_{\max} L)\<\nabla f(x^k) - \nabla f(x^\star), x^k - x^\star>\\
		 &\overset{\eqref{eq:scal_prod_cvx}}{\le} \|x^k - x^\star\|^2 - \gamma\left(2 - \gamma_{\max}L \right)D_f(x^k, x^\star).	
	\end{align*}
	Now let us consider $\mu$-strongly convex $f$. We have
	\begin{align*}
		\|w^k - w^\star\|^2
		&= \|x^k - x^\star\|^2 - 2 \gamma\<\nabla f(x^k) - \nabla f(x^\star), x^k - x^\star> + \gamma^2\|\nabla f(x^k) - \nabla f(x^\star)\|^2 \\
		&\overset{\eqref{eq:scal_prod_tight_str_cvx}}{\le} \left(1 - \frac{2 \gamma \mu L}{L + \mu}\right)\|x - y\|^2 - \gamma\left(\frac{2}{L + \mu} - \gamma\right) \|\nabla f(x^k) - \nabla f(x^\star)\|^2 \\
		&\overset{\eqref{eq:scal_prod_str_cvx}}{\le} \left(1 - \frac{2 \gamma \mu L}{L + \mu}\right)\|x - y\|^2 - \gamma\left(\frac{2}{L + \mu} - \gamma\right) \mu^2\|x^k - x^\star\|^2  \\
		&= (1 - \gamma\mu)^2\|x^k - x^\star\|^2 \\
		&\le (1 - \gamma\mu)\|x^k - x^\star\|^2.
	\end{align*}
	The last step simply uses $1 - \gamma\mu\le 1$, which, of course, makes our guarantees slightly weaker, but, on the other hand, puts \algname{Gradient Descent} under the umbrella of Assumption~\ref{as:method}.
\end{proof}

\subsection{Key lemma}\label{ap:key_lemma}
The result below is  the most important lemma of this chapter as it lies at the core of our analysis. It provides a very generic statement about the step with stochastic proximal operators.  At the same time, it is a mere corollary of firm nonexpansiveness of the proximal operator.

\begin{lemma} \label{lem:key}
    Let $z^k = \proxR(w^k - \gamma y^k)$ and $x^{k+1} = \proxj(z^{k} + \eta_j  y_j^k)$, where $j$ is sampled from $\{1,\dotsc, m\}$ with probabilities $\{p_1,\dotsc, p_m\}$, $\eta_j = \frac{\gamma}{m p_j}$ and $\gamma$ is a positive number. If $y_j^{k+1} = y_j^k + \frac{1}{\eta_j}(z^k - x^{k+1})$ and $y_l^{k+1} = y_l^k$ for all $l\neq j$, it holds
    \begin{align*}
        &\ec{\|x^{k+1} - x^\star\|^2 + \cY^{k+1}} \\
        &\le \ec{\|w^k - w^\star\|^2 + \left(1 - \frac{\nu}{m(1+\nu)}\right)\cY^k - \|z^k - w^k - (x^\star - w^\star)\|^2},
    \end{align*}
    where $\nu \eqdef \min_{j=1,\dotsc,m} \frac{1}{\eta_j L_j}$ and $L_j\in\RR\cup \{+\infty\}$ is the smoothness constant of $g_j$.
\end{lemma}
\begin{proof}
    Mention that $x^\star = \proxj(x^\star + \eta_j  y_j^\star)$ by optimality condition. In addition, it holds by definition $y_j^{k+1} = \frac{1}{\eta_j}(z^k + \eta_j y_j^k - x^{k+1})$, so property~\eqref{eq:nonexp} yields
    \begin{align*}
        \|x^{k+1} - x^\star\|^2 + \left(1 + \frac{1}{\eta_j L_j}\right)\eta_j^2 \|y_j^{k+1} - y_j^\star\|^2
        &\le \|z^k + \eta_j y_j^k - (x^\star + \eta_j y_j^\star)\|^2
    \end{align*}
    and we can replace $1+\frac{1}{\eta_j L_j}$ with $1+\nu$ since $\nu\le \frac{1}{\eta_j L_j}$.
    
    Let $\mathbb{E}_j$ be the expectation with respect to sampling of $j$. Then, we observe
    \begin{align}
        &\mathbb{E}_j \left[\|z^k + \eta_j y_j^k - (x^\star + \eta_j y_j^\star)\|^2 \right]\notag \\
        &= \|z^k - x^\star\|^2 + \mathbb{E}_j\left[\frac{\gamma^2}{m^2p_j^2}\|y_j^k - y_j^\star\|^2\right] + 2\<z^k - x^\star, \gamma\mathbb{E}_j \left[\frac{1}{mp_j}(y_j^k - y_j^\star) \right]> \nonumber\\
        &= \|z^k - x^\star\|^2 + \frac{\gamma^2}{m^2}\sumlm \frac{1}{p_l}\|y_l^k - y_l^\star\|^2 + 2\gamma\<z^k - x^\star, y^k - y^\star>. \label{eq:firm_nonexp_for_stoch_prox}
    \end{align}
    Denote $w^\star\eqdef x^\star - \gamma \nabla f(x^\star)$. Another optimality condition from Lemma~\ref{lem:optimality} is $x^\star=\proxR(w^\star - \gamma y^\star)$, so let us use~\eqref{eq:nonexp} one more time to obtain
    \begin{align*}
        \|z^k - x^\star\|^2
        &\le \|w^k - \gamma y^k - (w^\star - \gamma y^\star)\|^2 - \|w^k - \gamma y^k - z^k - (w^\star - \gamma y^\star - x^\star)\|^2 \\
        &= \gamma^2\|y^k - y^\star\|^2 - 2\gamma\<w^k - w^\star , y^k - y^\star> + \|w^k - w^\star\|^2\\
        &\qquad - \|w^k - \gamma y^k - z^k - (w^\star - \gamma y^\star - x^\star)\|^2.
    \end{align*}
    Furthermore,
    \begin{align*}
        \|w^k - \gamma y^k - z^k - (w^\star - \gamma y^\star - x^\star)\|^2
        &= \|w^k - z^k - (w^\star - x^\star)\|^2 + \gamma^2\|y^k - y^\star\|^2\\
        &\quad - 2\gamma\<w^k - z^k - (w^\star - x^\star), y^k - y^\star>,
    \end{align*}
    so
    \begin{align*}
        \|z^k - x^\star\|^2 
        &\le - 2\gamma\<w^k - w^\star , y^k - y^\star>  + 2\gamma\<w^k - z^k - (w^\star - x^\star), y^k - y^\star> \\
        &\quad + \|w^k - w^\star\|^2 - \|w^k - z^k - (w^\star - x^\star)\|^2 \\
        &=  \|w^k - w^\star\|^2 - 2\gamma\<z^k - x^\star, y^k - y^\star> - \|w^k - z^k - (w^\star - x^\star)\|^2.
    \end{align*}
    Together with the previously obtained bounds it adds up to
    \begin{align*}
        \mathbb{E}_j \left[\|z^k + \eta_j y_j^k - (x^\star + \eta_j y_j^\star)\|^2 \right]
        &\le \|w^k - w^\star\|^2 + \frac{\gamma^2}{m^2}\sumlm \frac{1}{p_l}\|y_l^k - y_l^\star\|^2\\
        &\qquad - \|z^k - w^k - (x^\star - w^\star)\|^2.
    \end{align*}
    To get the expression in the left-hand side of this lemma's statement, let us add the missing sum and evaluate its expectation:
    \begin{align*}
        \mathbb{E} \left[\sumlm \eta_l^2\|y_l^{k+1} - y_l^\star\|^2 \right]
        = \ec{ \|y_j^{k+1} - y_j^\star\|^2 }+ \mathbb{E} \left[\sum_{l\neq j} \eta_l^2\|y_l^{k+1} - y_l^\star\|^2 \right].
    \end{align*}
    Clearly, all summands in the last sum were not changed at iteration $k$, so
    \begin{align*}
        \mathbb{E}_j \left[\sum_{l\neq j} \eta_l^2\|y_l^{k+1} - y_l^\star\|^2 \right]
        &= \mathbb{E}_j \left[\sum_{l\neq j} \eta_l^2\|y_l^{k} - y_l^\star\|^2 \right]\\
        &= \sumlm (1 - p_l)\eta_l^2\|y_l^{k} - y_l^\star\|^2 \\ 
        &= \sumlm \eta_l^2\|y_l^{k} - y_l^\star\|^2 - \frac{\gamma^2}{m^2}\sumlm \frac{1}{p_l}\|y_l^{k} - y_l^\star\|^2.
    \end{align*}
    The negative sum will cancel out with the same in equation~\eqref{eq:firm_nonexp_for_stoch_prox} and we conclude the proof.
\end{proof}

\subsection{Convergence of Bregman divergence to 0 almost surely}\label{ap:almost_sure}
Here we formulate a result that we only briefly mentioned in the main text. It states that for convex problems, Bregman divergence $D_f(x^k, x^\star)$ almost surely converges to 0. To show it, let us borrow the classical result on supermartingale  convergence.
\begin{proposition}[\cite{bertsekas2015convex}, Proposition~A.4.5]\label{pr:supmart}
	Let $(X^k)_k$, $(Y^k)_k$, $(Z^k)_k$ be three sequences of nonnegative random variables and let $(\mathcal{F}^k)_k$ be a sequence of $\sigma$-algebras such that $\mathcal{F}^k\subset \mathcal{F}^{k+1}$ for all $k.$ Assume that:
		\begin{itemize}
			\item The random variables $X^k, Y^k, Z^k$ are nonnegative and $\mathcal{F}^k$-measurable.
			\item For each $k$, we have $\EE[X^{k+1}\mid \mathcal{F}^k] \le X^k - Y^k + Z^k$.
			\item There holds, with probability 1,
			\begin{align*}
				\sum_{k=0}^{\infty} Z^k < \infty.
			\end{align*}
		\end{itemize}
	Then $X^k$ converges to a nonnegative random variable $X$ and we have $\sum_{k=0}^{\infty} Y^k < \infty$ with probability 1.
\end{proposition}

\begin{theorem}
	Take a method that satisfies Assumption~\ref{as:method}(a), a stepsize $\gamma\le \gamma_{\max}$ and an optimum $x^\star$ satisfying Assumption~\ref{as:optimality}. Then, with probability 1 it holds $D_f(x^k, x^\star)\to 0$.
\end{theorem}
\begin{proof}
	Fix any solution $x^\star$, $y_1^\star, \dotsc, y_m^\star$. Let $\mathcal{F}^k=\sigma(x^0, y_1^0, \dotsc, y_m^0, \dotsc, x^k, y_1^k, \dotsc, y_m^k)$ be the $\sigma$-algebra generated by all random variables prior to moment $k$, and let $\overline \cM^k$ be $\cM^k$ conditioned on $\mathcal{F}^k$, i.e., $\overline \cM^k \eqdef \cM^k | \mathcal{F}^k$, from which it follows $\cM^k = \ec{ \overline \cM^k}$. Then, the assumptions of Proposition~\ref{pr:supmart} are satisfied for sequences 
	\begin{align*}
		X^k 
		&= \|x^k - x^\star\|^2 + \overline \cM^k + (1 + \nu)\sumlm \eta_l^2\|y_l^k - y_l^\star\|^2, \\
		Y^k
		&= \omega \gamma D_f(x^k, x^\star), \\
		Z^k
		&=0.
	\end{align*}
	Therefore, we have that $\sum_{k=0}^{\infty} Y^k < \infty$ and $Y^k \to 0$ almost surely, from which it follows $D_f(x^k, x^\star)\to 0$.
\end{proof}
The almost sure guarantee is not applicable to \algname{SGD} which has $\cM^0$ proportional to the number of iterations. We leave its analysis as well as analysis of convergence of $x^k$ to an optimum for future work.
\subsection{Proof of Theorem~\ref{th:1_over_t_rate} ($\cO\left(\frac{1}{K}\right)$ rate)}\label{ap:1_t_rate}
Below we provide the proof of $\cO\left(\frac{1}{K}\right)$ rate for general convex functions.
\begin{theorem}[Same as Theorem~\ref{th:1_over_t_rate}]
    Assume $f$ is $L$-smooth and $\mu$-strongly convex, $g_1, \dotsc, g_m, \psi$ are convex, closed and lower semi-continuous. Take a method satisfying Assumption~\ref{as:method} and $\gamma\le \gamma_{\max}$, then
    \begin{align*}
        \ec{ D_f(\overline x^K, x^\star)}
        \le \frac{1}{\omega \gamma K}\cL^0,
    \end{align*}
    where $\cL^0\eqdef \|x^0 - x^\star\|^2 + \cM^0 + \sumlm  \eta_l^2 \|y_l^0 - y_l^\star\|^2$ and $\overline x^K \eqdef \frac{1}{K}\sum_{k=0}^{K-1} x^k$.
\end{theorem}
\begin{proof}
    Recall that
    \begin{align*}
        \cL^k
        \eqdef  \ec{\|x^k - x^\star\|^2} + \cM^k + \cY^k,
    \end{align*}
    and by Assumption~\ref{as:method} combined with Lemma~\ref{lem:key}
    \begin{align*}
    		\cL^{k+1}
    		\le \cL^k - \omega \gamma \ec{ D_f(x^k, x^\star)}.
    \end{align*}
	Telescoping this inequality from $0$ to $K-1$, we obtain
    \begin{align*}
        \mathbb{E} \left[\sum_{l=0}^{K-1} D_f(x^l, x^\star)\right]
        \le \frac{1}{\omega\gamma}(\cL^0 - \cL^{K})
        \le \frac{1}{\omega\gamma}\cL^0.
    \end{align*}    
    By convexity of $f$, the left-hand side is lower bounded by $K\ec{ D_f(\overline x^K, x^\star)}$, so dividing both sides by $K$ finishes the proof.
\end{proof}

\subsection{Proof of Theorem~\ref{th:1_t2_rate} ($\cO(\frac{1}{K^2})$ rate)}\label{ap:1_t2_rate}

In this subsection, we show the $\cO\left(\frac{1}{K^2}\right)$ rate.
\begin{theorem}[Same as Theorem~\ref{th:1_t2_rate}]
	Consider updates with time-varying stepsizes, $\gamma_{k-1} = \frac{2}{\omega\mu(a + k)}$ and $\eta_{k,j} = \frac{\gamma_k}{m p_j}$ for $j=1,\dotsc, m$, where $a\ge 2\max\left\{\frac{1}{\omega\mu\gamma_{\max}}, \frac{1}{\rho} \right\}$. Then, it holds
	\begin{align*}
		\ec{ \|x^K - x^\star\|^2}
		\le \frac{a^2}{(K+a-1)^2}\cL^0,
	\end{align*}
	where $\cL^0 = \|x^0 - x^\star\|^2 + \cM^0 + \sumlm \eta_{0, l}^2 \|y_l^0 - y_l^\star\|^2$.
\end{theorem}
\begin{proof}
	For this proof, we redefine the sequence $\cY^k$ to have time-varying stepsizes:
	\begin{align*}
		\cY^k
		\eqdef \sumlm \eta_{k, l}^2\ec{\|y_l^{k} - y_l^\star\|^2}.
	\end{align*}
	Before writing a new recurrence, let us note that 
	\begin{align*}
		(1 - \omega\gamma_k\mu)\left(\frac{\gamma_{k-1}}{\gamma_{k}}\right)^2
		= \frac{\left(1 - \frac{2}{a + k}\right)(a+k)^2}{(a+k-1)^2}
		= \frac{(a+k-2)(a+k)}{(a+k-1)^2}
		< 1,
	\end{align*}
	so $1 - \omega\gamma_k\mu \le \left(\frac{\gamma_k}{\gamma_{k-1}}\right)^2$. Then, Lemma~\ref{lem:key} gives a similar recurrence to what we have seen in other proofs, but the stepsizes in the right-hand side are now time-dependent:
	\begin{align*}
		\cL^{k+1}
		&=\mathbb{E} \left[\|x^{k+1} - x^\star\|^2 \right] + \cM^{k+1} + \cY^{k+1}\\
		&\le (1 - \omega\gamma_k \mu) \ec{\|x^k - x^\star\|^2} + (1 - \rho)\cM^k + \sumlm \eta_{k,l}^2\ec{\|y_l^{k} - y_l^\star\|^2} \\
		&\le (1 - \omega\gamma_k\mu)\mathbb{E}\left[\|x^k - x^\star\|^2 + \cM^k\right] + \left(\frac{\gamma_k}{\gamma_{k-1}}\right)\sumlm\eta_{k,l}^2\mathbb{E}\|y_l^k - y_l^\star\|^2 \\
		&\le \left(\frac{\gamma_k}{\gamma_{k-1}}\right)^2\mathbb{E}\left[\|x^k - x^\star\|^2 + \cM^k\right] + \left(\frac{\gamma_k}{\gamma_{k-1}}\right)^2\sumlm\eta_{k,l}^2\mathbb{E}\|y_l^k - y_l^\star\|^2  \\
		&= \left(\frac{\gamma_k}{\gamma_{k-1}}\right)^2\cL^k.
	\end{align*}
	Recursing this inequality yields
	\begin{align*}
		\cL^{k+1} 
		\le \cL^0\prod_{k=1}^k\left(\frac{\gamma_k}{\gamma_{k-1}}\right)^2
		= \left(\frac{\gamma_k}{\gamma_{0}}\right)^2\cL^0
		= \left(\frac{a}{a+k}\right)^2\cL^0.
	\end{align*}
	It remains to plug-in $k=K-1$.
\end{proof}
\subsection{Proof of Theorem~\ref{th:sgd_str_cvx} ($\cO(\frac{1}{K})$ rate of \algname{SGD})}\label{ap:sgd_str_cvx}
Here we consider the case where $f(x)$ is given as expectation parameterized by a random variable~$\xi$,
\begin{align*}
	f(x)
	= \EE_\xi \left[f(x; \xi) \right].
\end{align*}
While it is often assumed in the literature that $\mathbb{E} \|\nabla f(x; \xi) - \nabla f(x)\|^2\le \sigma^2$ uniformly over $x$, we do not need this assumption and bound the variance using the following lemma.
\begin{lemma}\label{lem:sgd_variance_sdm}
	Let $w^k = x^k - \gamma \nabla f(x^k; \xi^k)$, where random function $f(x; \xi)$ is almost surely convex and $L$-smooth. Then,
	\begin{align}
		\ec{ \|\nabla f(x^k; \xi^k) - \nabla f(x^\star)\|^2}
		\le 4L \ec{ D_f(x^k, x^\star)} + 2\sigma_{\star}^2,\label{eq:sgd_variance_sdm}
	\end{align}
	where $\sigma_{\star}^2\eqdef \ec{ \|\nabla f(x^\star; \xi) - \nabla f(x^\star)\|^2}$, i.e., $\sigma_{\star}^2$ is the variance at an optimum. If more than  one $x^\star$ exists, take the one that minimizes $\sigma_{\star}^2$.
\end{lemma}
\begin{proof}
	This proof is based on existing results for \algname{SGD} and goes in a very standard way. By Young's inequality
	\begin{align*}
		\mathbb{E} \left[\|\nabla f(x^k; \xi^k) - \nabla f(x^\star)\|^2\right]
		&\le 2 \mathbb{E} \left[\|\nabla f(x^k; \xi^k) - \nabla f(x^\star; \xi^k)\|^2\right]\\
		&\qquad + 2 \mathbb{E} \left[\|\nabla f(x^\star; \xi^k)  - \nabla f(x^\star)\|^2\right] \\
		&\overset{\eqref{eq:grad_dif_bregman}}{\le} 4L \mathbb{E} \left[D_{f(\cdot; \xi^k)}(x^k, x^\star)\right] + 2\sigma_{\star}^2 \\
		&= 4L \mathbb{E} \left[D_{f}(x^k, x^\star)\right] + 2\sigma_{\star}^2.
	\end{align*}
\end{proof}
In the proof of Theorem~\ref{th:sgd_str_cvx} we will again need time-varying stepsize and $\cY^k$ should be defined as
\begin{align*}
		\cY^k
		\eqdef \sum_{j=1}^m \eta_{k,j}^2\ec{\|y_j^{k} - y_j^\star\|^2}.
	\end{align*}
But before let us prove a simple statement about sequences with contraction and additive error.
\begin{lemma}
\label{lem:recurrence_with_sigma}
	Assume that sequence $(\cL^k)_k$ satisfies inequality $\cL^{k+1}\le \left(\frac{\gamma_k}{\gamma_{k-1}}\right)^2 \cL^k + 2\gamma_k^2 \sigma_{\star}^2$ with some constant $\sigma_{\star}\ge 0$. Then, it holds
	\begin{align*}
		\cL^{k} 
		\le \left(\frac{\gamma_{k-1}}{\gamma_0}\right)^2 \cL^0 + 2k \gamma_{k-1}^2 \sigma_{\star}^2.
	\end{align*}
\end{lemma}
\begin{proof}
	We will prove the inequality by induction. For $k=0$ it is straightforward. The induction step follows from
	\begin{align*}
		\cL^{k+1}
		&\le \left(\frac{\gamma_k}{\gamma_{k-1}}\right)^2 \cL^k + 2\gamma_k^2 \sigma_{\star}^2 \\
		&\le \left(\frac{\gamma_k}{\gamma_{k-1}}\right)^2\left(\frac{\gamma_{k-1}}{\gamma_0}\right)^2 \cL^0 + 2\left(\frac{\gamma_k}{\gamma_{k-1}}\right)^2(\gamma_{k-1})^2k\sigma_{\star}^2+ 2k \gamma_{k-1}^2 \sigma_{\star}^2 \\
		&= \left(\frac{\gamma_{k}}{\gamma_0}\right)^2 \cL^0 + 2(k+1) \gamma_{k-1}^2 \sigma_{\star}^2.
	\end{align*}
\end{proof}
Now we are ready to prove the theorem.
\begin{theorem}[Same as Theorem~\ref{th:sgd_str_cvx}]
	Assume $f$ is $\mu$-strongly convex, $f(\cdot; \xi)$ is almost surely convex and $L$-smooth. Let the update be produced by \algname{SGD}, i.e., $v^k = \nabla f(x^k; \xi^k)$, and let us use time-varying stepsizes $\gamma_{k-1} = \frac{2}{a + \mu k}$ with $a\ge 4L$. Then, it holds
	\begin{align*}
		\ec{ \|x^K - x^\star\|^2}
		\le \frac{8\sigma_{\star}^2}{\mu(a + \mu K)} + \frac{a^2}{(a + \mu K)^2}\cL^0.
	\end{align*}
\end{theorem}
\begin{proof}
	It holds by Lemma~\ref{lem:sgd_variance_sdm}
	\begin{align*}
		\ec{ \|\nabla f(x^k; \xi^k)  - \nabla f(x^\star)\|^2}
		\le 4L \ec{ D_f(x^k, x^\star)} + 2\sigma_{\star}^2.
	\end{align*}
	Therefore, for $w^k \eqdef x^k - \gamma_k v^k = x^k - \gamma_k \nabla f(x^k; \xi^k) $ and $w^\star\eqdef x^\star - \gamma_k \nabla f(x^\star)$ we have
	\begin{align*}
		&\ec{ \|w^k - w^\star\|^2}\\
		&= \mathbb{E} \left[\|x^k - x^\star\|^2 - 2\gamma_k\<\nabla f(x^k) - \nabla f(x^\star), x^k - x^\star> + \gamma_k^2\|\nabla f(x^k; \xi^k) - \nabla f(x^\star)\|^2 \right] \\
		&\overset{\eqref{eq:sgd_variance_sdm}}{\le} \mathbb{E} \left[\|x^k - x^\star\|^2 - 2\gamma_k\<\nabla f(x^k) - \nabla f(x^\star), x^k - x^\star> + 4\gamma_k^2LD_f(x^k, x^\star) + 2\gamma_k^2\sigma_{\star}^2 \right]\\
		&\overset{\eqref{eq:scal_prod_cvx}}{\le} \mathbb{E}\Bigl[(1 - \gamma_k\mu)\|x^k - x^\star\|^2 - 2\gamma_k(\underbrace{1 - 2\gamma_k L}_{\ge 0})D_f(x^k, x^\star) + 2\gamma_k^2\sigma_{\star}^2 \Bigr] \\
		&\le (1 - \gamma_k\mu)\ec{\|x^k - x^\star\|^2} + 2\gamma_k^2\sigma_{\star}^2.
	\end{align*}
	Using the same argument as in the proof of Theorem~\ref{th:1_t2_rate}, we can show that $1 - \gamma_k\mu\le \left(\frac{\gamma_k}{\gamma_{k-1}}\right)^2$. Combining these results with Lemma~\ref{lem:key}, we obtain for $\cL^{k+1}\eqdef \ec{\|x^{k+1} - x^\star\|^2} + \cY^{k+1}$ the following bound:
	\begin{align*}
		\cL^{k+1}
		&\le (1 - \gamma_k\mu)\ec{\|x^k - x^\star\|^2} + \sumlm \eta_{k,l}^2\ec{\|y_l^{k} - y_l^\star\|^2} + 2\gamma_k^2\sigma_{\star}^2 \\
		&\le \left(\frac{\gamma_k}{\gamma_{k-1}}\right)^2\ec{\|x^k - x^\star\|^2} + \left(\frac{\gamma_k}{\gamma_{k-1}}\right)^2\cY^k + 2\gamma_k^2\sigma_{\star}^2 \\
		&= \left(\frac{\gamma_k}{\gamma_{k-1}}\right)^2\cL^k + 2\gamma_k^2\sigma_{\star}^2.
	\end{align*}
	By Lemma~\ref{lem:recurrence_with_sigma}
	\begin{align*}
		\ec{\|x^K - x^\star\|^2}
		\le \cL^K
		\le \left(\frac{\gamma_{K-1}}{\gamma_0}\right)^2 \cL^0 + 2K \gamma_{k-1}^2 \sigma_{\star}^2
		\le \frac{a^2}{(a + \mu K)^2}\cL^0 + \frac{8k}{(a + \mu K)\mu K} \sigma_{\star}^2.
	\end{align*}
\end{proof}

\subsection{Proof of Theorem~\ref{th:lin_conv_lin_model} (linear rate for $g_j=\phi_j(\mA_j^\top x)$)}\label{ap:lin_conv_lin_model}

Let us now show linear convergence of our method when the consider problem has linear structure, i.e., $g_j(x) = \phi_j(\mA_j^\top x)$.

We first need a lemma on the nature of $y_1^k,\dotsc, y_m^k$ in the considered case.

\begin{lemma}\label{lem:linear}
    Let the proximal sum be of the form $\frac{1}{m}\sumjm \phi_j(\mA_j^\top x)$ with some matrices $\mA_j\in \RR^{d\times d_j}$, and $y_j^0 = \mA_j\beta_j^0$ for $j=1,\dotsc, m$. Then, if Assumption~\ref{as:optimality2} is satisfied, for any $k$ and $j$ we have
    \begin{align*}
        y_j^k = \mA_j \beta_j^k, \quad  y^k = \frac{1}{m} \sumjm y_j^k= \frac{1}{m}\mA \beta^k, \quad    y_j^\star =  \mA_j\beta_j^\star, \quad  y^\star = \frac{1}{m}\sumjm y_j^\star= \frac{1}{m}\mA \beta^\star
    \end{align*}
    with some vectors $\beta_i^k, \beta_i^\star\in\RR^{d_i}$ with $i=1,\dotsc, m$, $\beta^k\eqdef ((\beta_1^k)^\top, \dotsc, (\beta_m^k)^\top)^\top$, $\beta^\star \eqdef ((\beta_1^\star)^\top, \dotsc, (\beta_m^\star)^\top)^\top$ and $\mA\eqdef [\mA_1, \dotsc, \mA_m]$.
\end{lemma}
\begin{proof}
    By definition $y_j^{k+1} = y_j^k + \frac{1}{\eta_j}(z^k - x^{k+1}) = \frac{1}{\eta_j}(z^k + \eta_j y_j^k - x^{k+1})$. In addition, by Lemma~\ref{lem:prox_of_composition} there exists $\beta_j^{k+1} \in \partial \phi_j (\mA_j^\top x^{k+1})$ such that $x^{k+1} = \proxj(z^k + \eta_j y_j^k) \in z^k + \eta_j y_j^k - \eta_j\mA_j \beta_j^{k+1}$
     and, thus, $y_j^{k+1} = \mA_j\beta_j^{k+1} $. Therefore, we also have $y^k = \frac{1}{m}\mA\beta^k$.
    
    The claims about $y_1^\star,\dotsc, y_m^\star$  and $y^\star$ follow from Assumption~\ref{as:optimality2}.
\end{proof}

Now it is time to prove Theorem~\ref{th:lin_conv_lin_model}.

\begin{theorem}[Same as Theorem~\ref{th:lin_conv_lin_model}]
    Assume that $f$ is $\mu$-strongly convex, $\psi\equiv 0$, $g_j(x) = \phi_j(\mA_j^\top x)$ for $j=1,\dotsc, m$ and take a method satisfying Assumption~\ref{as:method} with $\rho>0$. Then, if $\gamma\le \gamma_{\max}$,
    \begin{align*}
        \ec{ \|x^K - x^\star\|^2}
        \le \left(1 - \min\{\rho, \omega\gamma\mu, \rho_{\mA}\} \right)^K \cL^0,
    \end{align*}
    where $\rho_{\mA} \eqdef \lambda_{\min}(\mA^\top\mA) \min_j \left(\frac{p_j}{\|\mA_j\|}\right)^2$, and $\cL^0\eqdef \|x^0 - x^\star\|^2 + \cM^0 + \sumlm  \eta_l^2 \|y_l^0 - y_l^\star\|^2$.
\end{theorem}

\begin{proof}
    Lemma~\ref{lem:key} with Assumption~\ref{as:method} yields
    \begin{align*}
        \cL^{k+1}
        &\le (1 - \min\{\rho, \omega\gamma\mu\}) \left(\mathbb{E} \left[\|x^k - x^\star\|^2\right] + \cM^k\right)  + \cY^k - \mathbb{E} \left[\|z^k - w^k - (x^\star - w^\star)\|^2 \right].
    \end{align*}
    By Lemma~\ref{lem:linear}
    \begin{align*}
        \cY^k 
        = \sumlm  \eta_l^2 \mathbb{E}\left[\|y_l^k - y_l^\star\|^2 \right]
        = \sumlm \eta_l^2 \mathbb{E}\left[\|\mA_l(\beta_l^k - \beta_l^\star)\|^2 \right].
    \end{align*}
    Since we assume $\psi\equiv \mathrm{const}$, we have $z^k - w^k = x^k - \gamma v^k - \gamma y^k - (x^k - \gamma v^k) = -\gamma y^k$ and
    \begin{align*}
        \|z^k - w^k - (x^\star - w^\star)\|^2
        &= \gamma^2 \|y^k - y^\star\|^2\\
        &= \frac{\gamma^2}{m^2} \|\mA (\beta^k - \beta^\star)\|^2\\
        &\ge \lambda_{\min} (\mA^\top\mA) \frac{\gamma^2}{m^2}\|\beta^k - \beta^\star\|^2\\
        &= \lambda_{\min} (\mA^\top\mA) \sumlm \frac{p_l^2}{\|\mA_l\|^2}\eta_l^2\|\mA_l\|^2\|\beta_l^k - \beta_l^\star\|^2\\
        &\ge \lambda_{\min} (\mA^\top\mA) \min_l \left(\frac{p_l^2}{\|\mA_l\|^2}\right) \sumlm \eta_l^2\|\mA_l\|^2\|\beta_l^k - \beta_l^\star\|^2\\
        &\ge \lambda_{\min} (\mA^\top\mA) \min_l \left(\frac{p_l^2}{\|\mA_l\|^2}\right) \sumlm \eta_l^2\|\mA_l(\beta_l^k - \beta_l^\star)\|^2\\
        &=\rho_\mA \sumlm  \eta_l^2 \|y_l^k - y_l^\star\|^2.
    \end{align*}
    Therefore,
    \begin{align*}
        \cY^{k} - \ec{ \|z^k - w^k - (x^\star - w^\star)\|^2}
        &\le (1 - \rho_\mA)\sumlm  \eta_l^2 \ec{ \left[\|y_l^k - y_l^\star\|^2 \right]} \\
        &\le (1 - \rho_\mA)\cY^k.
    \end{align*}
    Putting the pieces together, we obtain
    \begin{align*}
        \cL^{k+1}
        \le (1 - \min\{\rho, \omega\gamma\mu, \rho_\mA\})\cL^k,
    \end{align*}
    from which it follows that $\cL^k$ converges to 0 linearly. Finally, note that 
    \[
    		\ec{ \|x^K - x^\star\|^2}\le \cL^K\le (1 - \min\{\rho, \omega\gamma\mu, \rho_\mA\})^k\cL^0.
    	\]
\end{proof}

\subsection{Proof of Theorem~\ref{th:lin_constr} (linear constraints)}\label{ap:lin_constr}
Here we discuss the problem of linearly constrained minimization
\begin{align*}
	\min_x \{f(x): \mA^\top x = b \}.
\end{align*}
We split matrix $\mA =  [\mA_1, \dots, \mA_m]$ and vector $b=(b_1^\top, \dotsc, b_m^\top)^\top$ and define projection operator $\Pi_j(\cdot)\eqdef \Pi_{\{x:\mA_j^\top x= b_j\}}(\cdot)$ . Since $y_j^k\in \Range{\mA_j}$, it is orthogonal to the hyperplane $\{x: \mA_j^\top x = b_j\}$ for any $x$ it holds
\begin{align*}
	\Pi_j(x + y_j^k) 
	=\Pi_j(x).
\end{align*}
This allows us to write a memory-efficient version of Algorithm~\ref{alg:sdm} as given in Algorithm~\ref{alg:sdm_lin_system}. If only a subset of functions $g_1, \dotsc, g_m$ is linear equality constraints, then similarly the corresponding vectors $y_j^k$ are not needed in the update, although they are still useful for the analysis.
\begin{algorithm}[t]
   \caption{\algname{Stochastic Decoupling Method} for linearly constrained problem}
   \label{alg:sdm_lin_system}
\begin{algorithmic}[1]
   \Require Stepsize $\gamma$, initial vectors $x^0, y^0\in \RR^d$, probabilities $p_1,\dotsc, p_m$, oracle that gives gradient estimates, number of steps $K$
   \For{$k=0,1,\dotsc, K-1$}
	   \State Produce an estimate $v^k$ of $\nabla f(x^t)$
	   \State $z^{k} = \proxR(x^k - \gamma v^k - \gamma y^k)$
	   \State Sample $j$ from $\{1,\dotsc, m\}$ with probabilities $\{p_1, \dotsc, p_m\}$
	   \State $x^{k+1} = \Pi_{j}(z^k)$
	   \State $y^{k+1} = y^k + \frac{p_j}{\gamma}(z^k - x^{k+1})$
   \EndFor
\end{algorithmic}
\end{algorithm}

Here we show that if $f$ is strongly convex and the non-smooth part is constructed of linear constraints, then we can guarantee linear rate of convergence. Moreover, the rate will depend only on the smallest nonzero eigenvalue of $\mA^\top \mA$, implying that even if $\mA^\top \mA$ is degenerate, convergence will be linear.
\begin{theorem}[Same as Theorem~\ref{th:lin_constr}]
	Under the same assumptions as in Theorem~\ref{th:lin_conv_lin_model} and assuming, in addition, that $g_j(x) = \ind_{\{x\;:\;\mA_j^\top x = b_j\}}$ it holds 
	\[
		\ec{ \|x^K - x^\star\|^2} \le (1 - \min\{\rho, \omega\gamma\mu, \rho_\mA\})^K\cL^0
	\]
	with $\rho_\mA=\lambda_{\min}^+(\mA^\top\mA)\min_{j}\left(\frac{p_j}{\|\mA_j\|}\right)^2$, i.e., $\rho_\mA$ depends only on  the smallest positive eigenvalue of $\mA^\top\mA$.
\end{theorem}
\begin{proof}
	The main reason we get an improved guarantee for linear constraints is that one can write a closed form expression for the proximal operators:
	\begin{align*}
		\proxj(x)
		= x - \mA_j(\mA_j^\top \mA_j)^\dagger (\mA_j^\top x - b_j).
	\end{align*}
	Assume that $j$ was sampled at iteration $k$, then
	\begin{align*}
		y_j^{k+1}
		&= \frac{1}{\eta_j}\left(z^k + \eta_j y_j^k - \proxj(z^k + \eta_j y_j^k) \right)\\
		&= \mA_j(\mA_j^\top \mA_j)^\dagger (\mA_j^\top (z^k + \eta_i y_j^k) - b_j).
	\end{align*}
	Therefore, for any $j$ and $k$ there exists a vector $x_j^k\in\RR^d$ such that 
	\begin{align*}
		y_j^{k} 
		&= \mA_j(\mA_j^\top \mA_j)^\dagger (\mA_j^\top x_j^k - b_j) \\
		&= \mA_j(\mA_j^\top \mA_j)^\dagger \mA_j^\top (x_j^k - x^\star),
	\end{align*} 
	where the second step is by the fact that $x^\star$ is from the set $\{x:\mA_j^\top x = b_j\}$.
	One can show using SVD that $\Range{(\mA_j^\top \mA_j)^\dagger\mA_j^\top}=\Range{\mA_j^\top}$. Then, $y_j^k - y_j^\star = \frac{1}{m}\mA_j(\beta_j^k - \beta_j^\star)$ with $\beta_j^k - \beta_j^\star\in \Range{\mA_j^\top}$. This, in  turn, implies $\beta^k - \beta^\star \in\Range{\mA^\top}$, so
    \begin{align*}
        \|z^k - w^k - (x^\star - w^\star)\|^2
        &= \gamma^2 \|y^k - y^\star\|^2\\
        &= \frac{\gamma^2}{m^2} \|\mA (\beta^k - \beta^\star)\|^2\\
        &\ge \lambda_{\min}^+ (\mA^\top\mA) \frac{\gamma^2}{m^2}\|\beta^k - \beta^\star\|^2.
    \end{align*}
    The rest of the proof goes the same way as that of Theorem~\ref{th:lin_conv_lin_model} in Appendix~\ref{ap:lin_conv_lin_model}.
\end{proof}

\subsection{Proof of Theorem~\ref{th:lin_conv_smooth} (smooth $g_j$)}\label{ap:lin_conv_smooth}
This is the only proof where Lemma~\ref{lem:key} is used with finite smoothness constants, i.e., $\max_{j=1,\dotsc, m}L_j < +\infty$. On the other hand, we will not use the negative square term from Lemma~\ref{lem:key}, which is rather needed in the case $g_j(x) = \phi_j(\mA_j^\top x)$.
\begin{theorem}[Same as Theorem~\ref{th:lin_conv_smooth}]
    Assume that $f$ is $L$-smooth and $\mu$-strongly convex, $g_j$ is $L_j$-smooth for $j=1,\dotsc, m$ and Assumption~\ref{as:method}(b) is satisfied. Then, Algorithm~\ref{alg:sdm} converges as
	\begin{align*}
		\ec{\|x^K - x^\star\|^2}
		\le \left(1 - \min\left\{\omega\gamma\mu, \rho, \frac{\nu}{m(1+\nu)}\right\}\right)^K\cL^0,
	\end{align*}	    
    where $\nu \eqdef \min_{j=1,\dotsc,m} \frac{1}{\eta_j L_j}$.
\end{theorem}
\begin{proof}
	Following the same lines as in the proof of Theorem~\ref{th:lin_conv_lin_model}, we get a contraction in $\cY^k$. Now we obtain it from the fact that functions $g_1,\dotsc, g_m$ are smooth, so the recursion is
	\begin{align*}
		\cL^{k+1}
		&\le (1 - \omega\gamma\mu)\mathbb{E}\|x^k-x^\star\|^2 +  (1 - \rho)\cM^k + \left(1 - \frac{\nu}{m(1+\nu)}\right)\cY^k \\
		&\le \left(1 - \min\left\{\omega\gamma\mu, \rho, \frac{\nu}{m(1+\nu)}\right\}\right)\cL^k.
	\end{align*}
	This is sufficient to show the claimed result.
\end{proof}

\subsection{Proof of Corollary~\ref{cor:acc_in_g} (optimal stepsize)}\label{ap:acc_in_g}
Corollary~\ref{cor:acc_in_g} is a statement about the optimal stepsizes for the case where $g_1, \dotsc, g_m$ are smooth functions. Its proof is a mere check that the choice of stepsizes gives the claimed complexity. 
\begin{corollary}[Same as Corollary~\ref{cor:acc_in_g}]
	Choose as solver for $f$ \algname{SVRG} or \algname{SAGA} without mini-batching, which satisfy Assumption~\ref{as:method} with $\gamma_{\max}=\frac{1}{5L}$ and $\rho = \frac{1}{3n}$, and consider for simplicity situation where $L_1= \dotsb = L_m \eqdef L_g$ and $p_1=\dotsb=p_m$. Define $\gamma_{\mathrm{best}} \eqdef \frac{1}{\sqrt{\omega \mu m L_g}}$,
	and set the stepsize to be $\gamma=\min\{\gamma_{\max}, \gamma_{\mathrm{best}}\}$. Then the complexity to get $\ec{\|x^K - x^\star\|^2}\le \varepsilon$ is
\[
		K=\cO\left(\left(n + m + \frac{L}{\mu} + \sqrt{m\frac{L_g}{\mu}}\right) \log\frac{1}{\varepsilon}\right).
\]
\end{corollary}
\begin{proof}
	According to Theorem~\ref{th:lin_conv_smooth}, in general, for any $\gamma\le \gamma_{\max}$ the complexity to get $\ec{\|x^K - x^\star\|^2}\le \varepsilon$ is 
	\begin{align*}
		\cO\left(\left(\frac{1}{\rho} + m + \frac{1}{\omega\gamma\mu} + \frac{1}{\nu}m\right) \log\frac{1}{\varepsilon}\right),
	\end{align*}
	where $\frac{1}{\nu}$ simplifies to $\gamma L_g$ when $L_1=\dotsb=L_m=L_g$ and $p_1=\dotsb=p_m=\frac{1}{m}$. In addition, for \algname{SVRG} and \algname{SAGA}, $\omega$ is a constant close to 1, so we can ignore it. Since $m$ and $\frac{1}{\rho}=3n$ do not depend on $\gamma$, we only need to simplify the other two terms. One of them decreases with $\gamma$ and the other increases, so the best complexity is achieved when the two quantities are equal to each other. The corresponding equation is
		$\omega \gamma^2 \mu m L_g = 1$,
	whose solution is
	\begin{align*}
		\gamma
		=\gamma_{\mathrm{best}}
		= \frac{1}{\sqrt{\omega \mu m L_g}}.
	\end{align*}
	Thus, we see that $\gamma_{\mathrm{best}}$ is optimal. Moreover, if $\gamma_{\mathrm{best}}\le \gamma_{\max}$ and $\gamma=\gamma_{\mathrm{best}}$, the two terms in the complexity both become equal to
	\begin{align*}
		\frac{1}{\omega \gamma_{\mathrm{best}}\mu}
		= m \gamma_{\mathrm{best}} L_g
		= \sqrt{\frac{m L_g}{\omega \mu}}.
	\end{align*}
	However, if $\gamma_{\mathrm{best}} > \gamma_{\max}$, then $\gamma=\min\{\gamma_{\max}, \gamma_{\mathrm{best}}\}=\gamma_{\max}$ is relatively small and the dominating term in the complexity is $\frac{1}{\omega\gamma\mu}$ rather than $\gamma L_g m$. Therefore, the complexity is
	\begin{align*}
		\cO\left(n + m + \frac{1}{\gamma_{\max}\mu} \right)
		= \cO\left(n + m + \frac{L}{\mu} \right).
	\end{align*}
	Combining the two complexities into one, we get the result.
\end{proof}
\subsection{Proof of Lemma~\ref{lem:svrg_saga} (\algname{SVRG} and \algname{SAGA})}\label{ap:svrg_saga}
\begin{algorithm}[t]
   \caption{\algname{Stochastic Decoupling Method} with \algname{SVRG}}
   \label{alg:svrg}
\begin{algorithmic}[1]
   \Require Stepsize $\gamma$, initial vectors $x^0$, $u^0$, $\nabla f(u^0)$, $y_1^0, \dotsc, y_m^0$, $y^0=\avejm y_j^0$, mini-batch size $\tau$, number of steps $K$
   \For{$k=0,1,\dotsc, K-1$}
	   \State Sample subset $S$ from $\{1,\dotsc, n\}$ of size $\tau$
	   \State $v^k= \frac{1}{\tau}\sum_{i\in S}\left(\nabla f_i(x^k) - \nabla f_i(u^k) + \nabla f(u^k)\right)$
	   \State $z^{k} = \proxR(x^k - \gamma v^k - \gamma y^k)$
	   \State $u^{k+1}= \begin{cases} x^k, &\text{with probability } \frac{\tau}{n},\\ u^k, & \text{with probability } 1 - \frac{\tau}{n} \end{cases}$
	   \State Sample $j$ from $\{1,\dotsc, m\}$ with probabilities $\{p_1, \dotsc, p_m\}$ and set $\eta_j = \frac{\gamma}{mp_j}$
	   \State $x^{k+1} = \proxj\left(z^k + \eta_j y_j^k\right)$
	   \State $y_j^{k+1} = y_j^k + \frac{1}{\eta_j}(z^{k} - x^{k+1})$
	   \State $ y^{k+1} =  y^k + \frac{1}{m}(y_j^{k+1} - y_j^k)$
   \EndFor
\end{algorithmic}
\end{algorithm}
Here we consider the update rule of \algname{SVRG} and \algname{SAGA} with mini-batch of size $\tau$. Following~\cite{hofmann2015variance}, we analyze \algname{SVRG} and \algname{SAGA} together by treating them both as memorization methods. More precisely, \algname{SAGA} stores each gradient estimate, $\nabla f_i(u_i^k)$ individually, and \algname{SVRG} stores only the reference point, $u^k$, itself and every iteration reevaluates $\nabla f_i(u^k)$ for all sampled $i$ to compute $v^k$. To avoid any confusion, we provide the explicit formulation of our method with the \algname{SVRG} solver in Algorithm~\ref{alg:svrg}.

First of all, let us show that the estimate that we use, $v^k$, is unbiased.
\begin{lemma}
	Let us sample a set of indices $S$ of size $\tau$ from $\{1, \dotsc, n\}$. Then, it holds for
    \begin{align*}
	    v^k
	    \eqdef \frac{1}{\tau}\sum_{i\in S} \left(\nabla f_i(x^k) - \nabla f_i(u_i^k) + \alpha^k\right)
    \end{align*} 
    that it is unbiased
    \begin{align}
	    \ec{v^k}
	    = \ec{\nabla f(x^k)}. \label{eq:svrg_saga_unbiased}
    \end{align} 
\end{lemma}
\begin{proof}
	Clearly, since $i$ is sampled with probability $\frac{\tau}{n}$, it holds
	\begin{align*}
		\frac{1}{\tau}\mathbb{E} \left[\nabla f_i(x^k) - \nabla f_i(u_i^k)\right]
		&= \mathbb{E}\left[\frac{1}{n}\sumin (\nabla f_i(x^k) - \nabla f_i(u_i^k)) \right]\\
		&= \mathbb{E}\left[\nabla f(x^k) - \avein \nabla f_i(u_i^k) \right].
	\end{align*}
	Therefore, $\mathbb{E} [v^k]  = \mathbb{E} [\nabla f(x^k)]$. 
\end{proof}

We continue our analysis with the following lemma.
\begin{lemma}
	Consider \algname{SVRG} and \algname{SAGA} solver for $f$. Assume that every $f_i$ is convex and $L$-smooth and define
	\begin{align}
		\cM^k
		\eqdef \frac{3\gamma^2}{\tau}\sumin \ec{\|\nabla f_i(u_i^k) - \nabla f_i(x^\star)\|^2},\label{eq:mt_svrg_saga}
	\end{align}
	where for \algname{SVRG} $u_1^k=\dotsb=u_n^k$ is the reference point at moment $k$ and for \algname{SAGA} it is the point whose gradient is saved in memory for function $f_i$. Then,
	\begin{align*}
		\cM^{k+1}
		\le \left(1 - \frac{\tau }{n}\right) \cM^k + 6\gamma^2L \mathbb{E} \left[D_f(x^k, x^\star)\right].
	\end{align*}
\end{lemma}
\begin{proof}
	We have for \algname{SVRG} that $\cM^{k+1}$ changes with probability $\frac{\tau}{n}$ and with probability $1 - \frac{\tau}{n}$ it remains the same. Then,
	\begin{align*}
		\mathbb{E} \left[\sumin \|\nabla f_i(u_i^{k+1}) - \nabla f_i(x^\star)\|^2 \right]
		&= \frac{\tau}{n} \sumin \mathbb{E}\left[\|\nabla f_i(x^{k}) - \nabla f_i(x^\star)\|^2\right] + \left(1  - \frac{\tau}{n}\right)\cM^k.
	\end{align*}
	Similarly, for \algname{SAGA} we update exactly $\tau$ out of $n$ gradient in the memory, which leads to the following identity:
    \begin{align*}
        &\ec{ \sumin \|\nabla f_i(u_i^{k+1}) - \nabla f_i(x^\star)\|^2 }\\
        &\qquad = \ec{ \sum_{i\in S}\|\nabla f_i(u_i^{k+1}) - \nabla f_i(x^\star)\|^2} + \ec{ \sum_{i\not \in S} \|\nabla f_i(u_i^{k+1}) - \nabla f_i(x^\star)\|^2}\\
        &\qquad = \frac{\tau}{n}\sumin\ec{\|\nabla f_i(x^k) - \nabla f_i(x^\star)\|^2} + \left(1  - \frac{\tau}{n}\right)\cM^k.
    \end{align*}
    In both cases, we obtained the same recursion. Now let us bound the gradient difference in the identity above:
    \begin{align*}
    		\frac{1}{n}\sumin\|\nabla f_i(x^k) - \nabla f_i(x^\star)\|^2
    		\overset{\eqref{eq:grad_dif_bregman}}{\le} \frac{1}{n}\sumin 2L D_{f_i}(x^k, x^\star) 
    		= 2L D_f(x^k, x^\star).
    \end{align*}
    This gives us the claimed inequality.
\end{proof}
Now let us show how the recursion looks like when $\cM^{k+1}$ is combined with $\|w^{k} - w^\star\|^2$.
\begin{lemma}\label{lem:saga}
    Consider the iterates of Algorithm~\ref{alg:sdm} with \algname{SVRG} or \algname{SAGA} estimate $v^k$. Let $f_1, \dotsc, f_n$ be convex and $L$-smooth, $f$ be $\mu$-strongly convex, $S$ be a subset of $\{1,\dotsc, n\}$ of size $\tau$ sampled with equal probabilities, $\alpha^k = \avein \nabla f_i(u_i^k)$ and $w^k = x^k - \gamma v^k$ with 
    \begin{align*}
	    v^k
	    =\frac{1}{\tau} \sum_{i\in S} \left(\nabla f_i(x^k) - \nabla f_i(u_i^k) + \alpha^k\right)
    \end{align*} 
    then we have
    \begin{align*}
        \mathbb{E} \left[\|w^{k} - w^\star\|^2 + \cM^{k+1} \right]
        \le (1 - \rho)\mathbb{E}\left[\|x^k - x^\star\|^2 + \cM^k\right],
    \end{align*}
    where $w^\star \eqdef x^\star - \gamma \nabla f(x^\star)$ and  $\rho\eqdef \min \left\{\gamma \mu, \frac{\tau}{3n}\right\}$.
\end{lemma}
\begin{proof}
    It holds
    \begin{align*}
        &\mathbb{E} \left[\|w^{k} - w^\star\|^2 \right] \\
        &= \mathbb{E} \left[\|x^k - x^\star - \gamma(v^k - \nabla f(x^\star))\|^2 \right] \\
        &= \mathbb{E}\left[\|x^k - x^\star\|^2 - 2\gamma\<x^k - x^\star, \mathbb{E} [v^k\mid x^k] - \nabla f(x^\star)> + \gamma^2 \|v^k - \nabla f(x^\star)\|^2 \right]\\
        &\overset{\eqref{eq:svrg_saga_unbiased}}{=} \mathbb{E}\left[\|x^k - x^\star\|^2 - 2\gamma\<x^k - x^\star, \nabla f(x^k) - \nabla f(x^\star)> + \gamma^2 \|v^k - \nabla f(x^\star)\|^2 \right].\\
        &\overset{\eqref{eq:scal_prod_cvx}}{\le}
         \mathbb{E}\left[(1 - \gamma \mu)\|x^k - x^\star\|^2 - 2\gamma D_f(x^k, x^\star) + \gamma^2 \|v^k - \nabla f(x^\star)\|^2 \right].
    \end{align*}
    On the other hand, by Jensen's and Young's inequalities
    \begin{align*}
        &\ec{ \|v^k - \nabla f(x^\star)\|^2 }\\
        &= \ec{\biggl\|\frac{1}{\tau}\sum_{i\in S}(\nabla f_i(x^k) - \nabla f_i(u_i^k)) + \alpha^k - \nabla f(x^\star)) \biggr\|^2 }\\
        &\le \frac{1}{\tau}\ec{\sum_{i\in S}\left\|\nabla f_i(x^k) - \nabla f_i(x^\star) + \nabla f_i(x^\star)- \nabla f_i(u_i^k) + \alpha^k - \nabla f(x^\star) \right\|^2}\\
        &\le \frac{2}{\tau}\ec{\sum_{i\in S}\left\|\nabla f_i(x^k) - \nabla f_i(x^\star)\right\|^2}\\
        &\qquad + \frac{2}{\tau}\ec{\sum_{i\in S}\|\nabla f_i(x^\star) - \alpha_i^k + \alpha^k - \nabla f(x^\star)\|^2 }\\
        &\overset{\eqref{eq:grad_dif_bregman}}{\le} 4L\ec{ D_f(x^k, x^\star)} + \frac{2}{n}\sumin\ec{\left\|\nabla f_i(x^\star) - \alpha_i^k + \alpha^k - \nabla f(x^\star) \right\|^2}.
    \end{align*}
    Using inequality $\ec{\|X - \ec{ X }\|^2} \le \ec{ \|X\|^2}$ that holds for any random variable $X$, the second term can be simplified to
    \begin{align*}
        &\frac{2}{n}\sumin\ec{\left\|\nabla f_i(x^\star) - \alpha_i^k + \alpha^k - \nabla f(x^\star) \right\|^2 }\\
        &\le \frac{2}{n}\sumin\ec{\|\nabla f_i(u_i^k) - \nabla f_i(x^\star)\|^2 }\\
        &= \frac{2\tau }{3n}\cM^k.
    \end{align*}
    Thus,
    \begin{align}
        \mathbb{E} \left[\|w^{k} - w^\star\|^2 + \cM^{k+1} \right]
        &\le (1 - \gamma \mu)\mathbb{E}\left[\|x^k - x^\star\|^2\right] + \left(\left(1 - \frac{\tau}{n}\right) + \frac{2\tau }{3n}\right) \cM^k  \notag \\
        &\quad - 2\gamma\left(1 - 2\gamma L - 3\gamma  L\right)\ec{ D_f(x^k, x^\star) }. \label{eq:saga_precise}
    \end{align}
    The second term in the right-hand side can be dropped as $1 - 2\gamma L - \frac{c L}{n\gamma}= 1 - 2\gamma L - 3\gamma L\le 0$. In addition, $\rho\le \gamma\mu$ and $\rho \le \frac{\tau}{3n}$, so the claim follows.
\end{proof}
Now we are ready to prove Lemma~\ref{lem:svrg_saga}.
\begin{lemma}[Same as Lemma~\ref{lem:svrg_saga}]
	In \algname{SVRG} and \algname{SAGA}, if $f_i$ is $L$-smooth and convex for all $i$, Assumption~\ref{as:method} is satisfied with $\gamma_{\max} = \frac{1}{6L}$, $\omega = \frac{1}{3}$, $\rho = \frac{1}{3n}$ and 
	\begin{align*}
		\cM^k 
		= \frac{3\gamma^2}{n}\sumin \mathbb{E}\left[\|\nabla f_i(u_i^k) - \nabla f_i(x^\star)\|^2 \right],
	\end{align*}
	where in \algname{SVRG} $u_i^k=u^k$ is the reference point of the current loop, and in \algname{SAGA} $u_i^k$ is the point whose gradient is stored in memory for function $f_i$. If $f$ is also strongly convex, then Assumption~\ref{as:method} holds with $\gamma_{\max} = \frac{1}{5L}$, $\omega = 1$, $\rho = \frac{1}{3n}$ and the same $\cM^k$.
\end{lemma}
\begin{proof}
	Equation~\eqref{eq:saga_precise} gives immediately the second part of the claim. 
	
	Similarly, if $\gamma\le \frac{1}{6L}$, from Equation~\ref{eq:saga_precise} we obtain by mentioning $1 - 2\gamma L - \frac{c L}{n\gamma} = 1 - 5\gamma L \le \frac{1}{6}$ that
	\begin{align*}
        \mathbb{E} \left[\|w^{k} - w^\star\|^2 + \frac{c}{n}\sumin \|\alpha_i^{k+1} - \alpha_i^\star\|^2 \right]
        &\le \mathbb{E}\left[\|x^k - x^\star\|^2 + \frac{c}{n}\sumin \|\alpha_i^k - \alpha_i^\star\|^2\right] \\
        &\qquad -\frac{\gamma}{3}\ec{ D_f(x^k, x^\star)}.
    \end{align*}
\end{proof}
\subsection{Proof of Lemma~\ref{lem:sgd} (\algname{SGD})}\label{ap:sgd}
\begin{lemma}[Same as Lemma~\ref{lem:sgd}]
	Assume that at an optimum $x^\star$ the variance of stochastic gradients is finite, i.e., 
	\[
		\sigma_{\star}^2\eqdef \ec{ \|\nabla f(x^\star; \xi^k)  - \nabla f(x^\star)\|^2} 
		< +\infty.
	\]
	Then, \algname{SGD} that terminates after at most $K$ iterations satisfies Assumption~\ref{as:method}(a) with $\gamma_{\max}=\frac{1}{4L}$, $\omega=1$ and $\rho=0$. If $f$ is strongly convex, it satisfies Assumption~\ref{as:method}(b) with $\gamma_{\max} = \frac{1}{2L}$, $\omega=1$ and $\rho=0$. In both cases, sequence $(\cM^k)_{k=0}^{K}$ is given by
	\begin{align*}
		\cM^k = 2\gamma^2(K - k)\sigma_{\star}^2.
	\end{align*}
\end{lemma}
\begin{proof}
	Clearly, we have
	\begin{align*}
		\mathbb{E} \left[\|w^k - w^\star\|^2 \right]
		&= \mathbb{E}\left[\|x^k - x^\star\|^2 - 2\gamma\<\nabla f(x^k) - \nabla f(x^\star), x^k - x^\star>\right]\\
		&\qquad + \gamma^2 \ec{\|\nabla f(x^k; \xi^k)  - \nabla f(x^\star)\|^2} \\
		&\overset{\eqref{eq:sgd_variance_sdm}}{\le} \mathbb{E}\left[\|x^k - x^\star\|^2 - 2\gamma\<\nabla f(x^k) - \nabla f(x^\star), x^k - x^\star> \right]\\
		&\qquad + \gamma^2\ec{4L D_f(x^k, x^\star) + 2\sigma_{\star}^2} \\
		&\overset{\eqref{eq:scal_prod_cvx}}{\le}\mathbb{E}\left[(1 - \gamma\mu)\|x^k - x^\star\|^2 - 2\gamma(1 - 2\gamma L)D_f(x^k, x^\star)\right] +2\gamma^2 \sigma_{\star}^2 .
	\end{align*}
	If $f$ is not strongly convex, then $\mu=0$ and by assuming $\gamma \le \gamma_{\max} = \frac{1}{4L}$ we get $1 - 2\gamma L\ge \frac{1}{2}$ and
	\begin{align*}
		\ec{\|w^k - w^\star\|^2}
		\le \mathbb{E}\left[\|x^k - x^\star\|^2 - \gamma D_f(x^k, x^\star)\right] + 2\gamma^2 \sigma_{\star}^2.
	\end{align*}
	In case $\mu=0$, by defining $\{\cM^{k}\}_{k=0}^{K}$ with recursion 
	\begin{align*}
		\cM^{k+1}
		= \cM^k - 2\gamma^2\sigma_{\star}^2,
	\end{align*}
	we can verify Assumption~\ref{as:method}(a) as long as $\cM^{K} = \cM^0 - 2K \gamma^2\sigma_{\star}^2\ge 0$. This is the reason we choose $\cM^0=2K \gamma^2\sigma_{\star}^2$.
	
	On the other hand, when $\mu>0$ and $\sigma_{\star}=0$, it follows from $\gamma\le \gamma_{\max}=\frac{1}{2L}$ that 
	\begin{align*}
		\mathbb{E} \left[\|w^k - w^\star\|^2 \right]
		\le (1 - \gamma\mu)\mathbb{E} \left[\|x^k - x^\star\|^2 \right].
	\end{align*}
\end{proof}

\section{Additional Experiments}\label{sec:add_experiments}
Here we want to see how changing $m$ and $n$ affects the comparison between \algname{SVRG} with exact projection and decoupled \algname{SVRG} with one stochastic projection. The problem that we consider is again $\ell_2$-regularized constrained linear regression. We took Gisette dataset from LIBSVM, whose dimension is $d=5000$, and used its first 1000 observations to construct $f$ and $g$. In particular, we split these observations into soft loss $f_i(x)=\frac{1}{2}\|a_i^\top x - b_i\|^2$ and hard constraints $g_j(x) = \ind_{\{x:a_j^\top x = b_j\}}$ with $n+m=1000$ and we considered three choices of $n$: 250, 500 and 750. To make sure that the constraints can be satisfied, we generated a random vector $x_0$ from normal distribution $\cN(0, \frac{1}{\sqrt{d}})$ and set $b=\mA x_0$. In all cases, first part of data was used in $f$ and the rest in $g$. To better see the effect of changing $n$, we used fixed $\ell_2$ penalty of order $\frac{1}{(n+m)}$ for all choices of $n$.

Computing the projection of a point onto the intersection of all constraints as at least as expensive as $m$ individual projections and we count it as such for \algname{SVRG}. In practice it might be by orders of magnitude slower than this estimate for big matrices, but the advantage of our method can be seen even without taking it into account. On the other hand, to make the comparison fair in terms of computation trade-off, we use \algname{SVRG} with mini-batch 20 and our method with mini-batch 1. The stepsize for both methods is $\frac{1}{(2L)}$.
\begin{figure}[t]
     \centering
     \begin{subfigure}[t]{0.32\textwidth}
         \centering
         \includegraphics[width=1\textwidth]{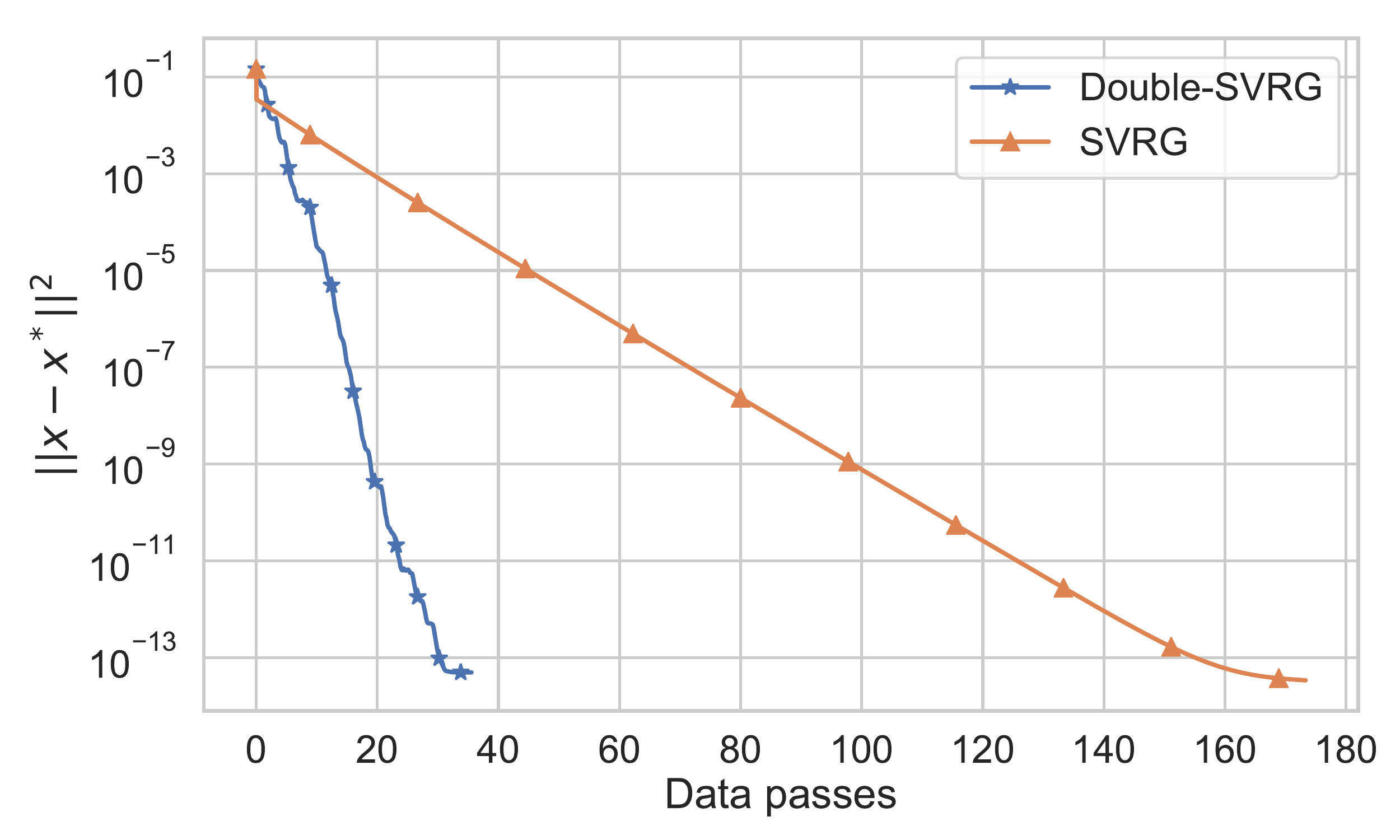}
         \caption{$m=100$, $n=900$}
     \end{subfigure}
     \begin{subfigure}[t]{0.32\textwidth}
         \centering
         \includegraphics[width=1\textwidth]{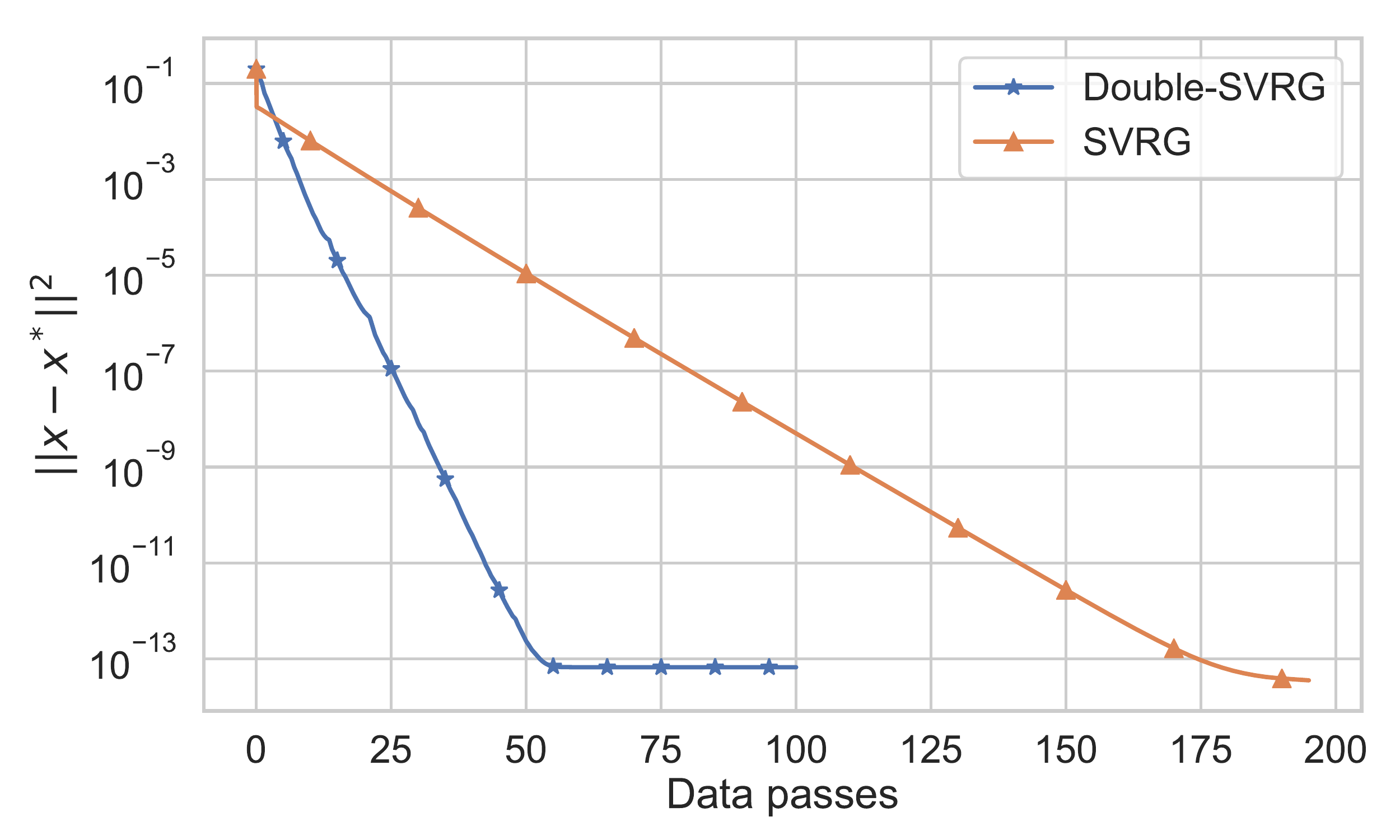}
         \caption{$m=200$, $n=800$}
     \end{subfigure}
     \begin{subfigure}[t]{0.32\textwidth}
         \centering
         \includegraphics[width=1\textwidth]{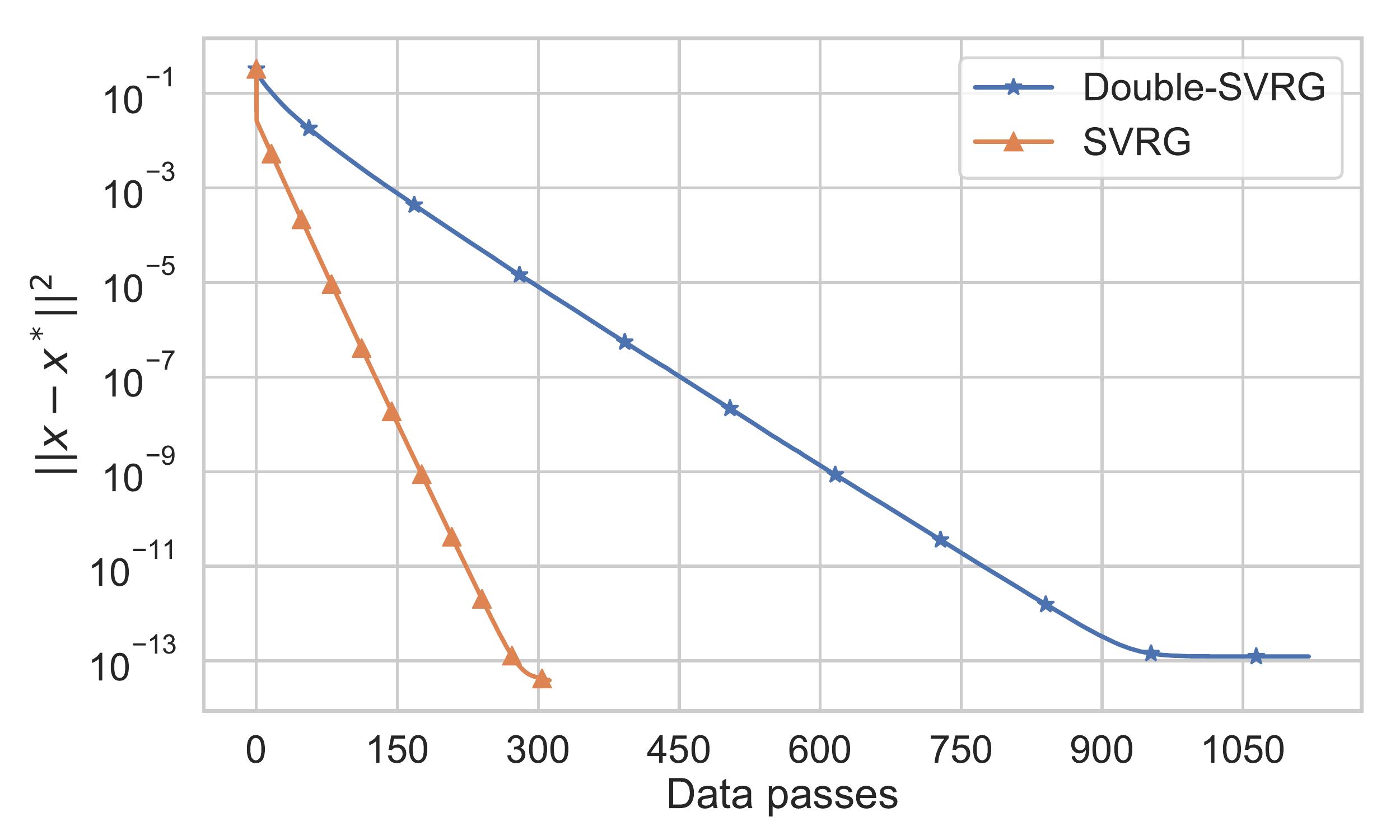}
         \caption{$m=500$, $n=500$}
     \end{subfigure}
     \\
     \begin{subfigure}[t]{0.32\textwidth}
         \centering
         \includegraphics[width=1\textwidth]{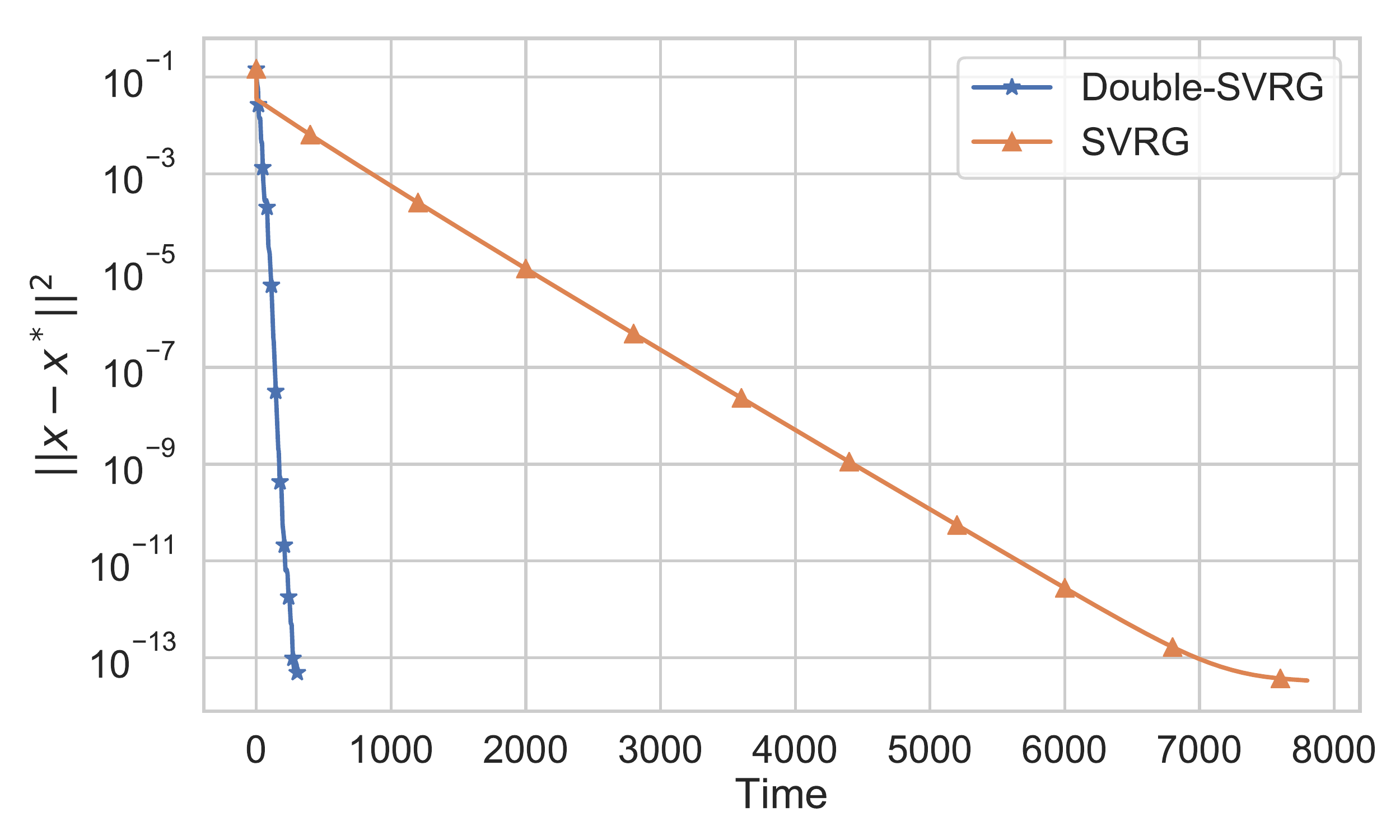}
         \caption{$m=100$, $n=900$}
     \end{subfigure}
     \begin{subfigure}[t]{0.32\textwidth}
         \centering
         \includegraphics[width=1\textwidth]{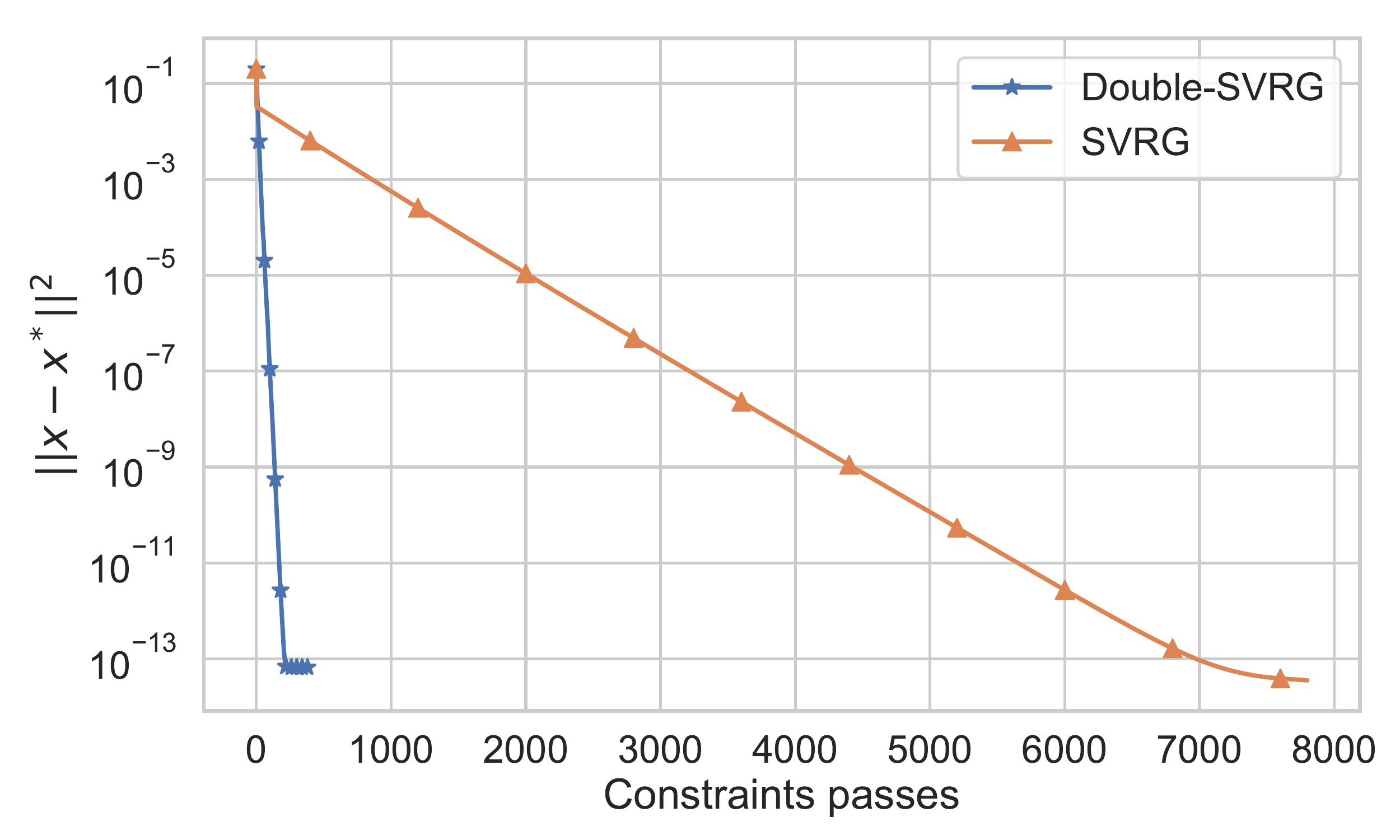}
         \caption{$m=200$, $n=800$}
     \end{subfigure}
     \begin{subfigure}[t]{0.32\textwidth}
         \centering
         \includegraphics[width=1\textwidth]{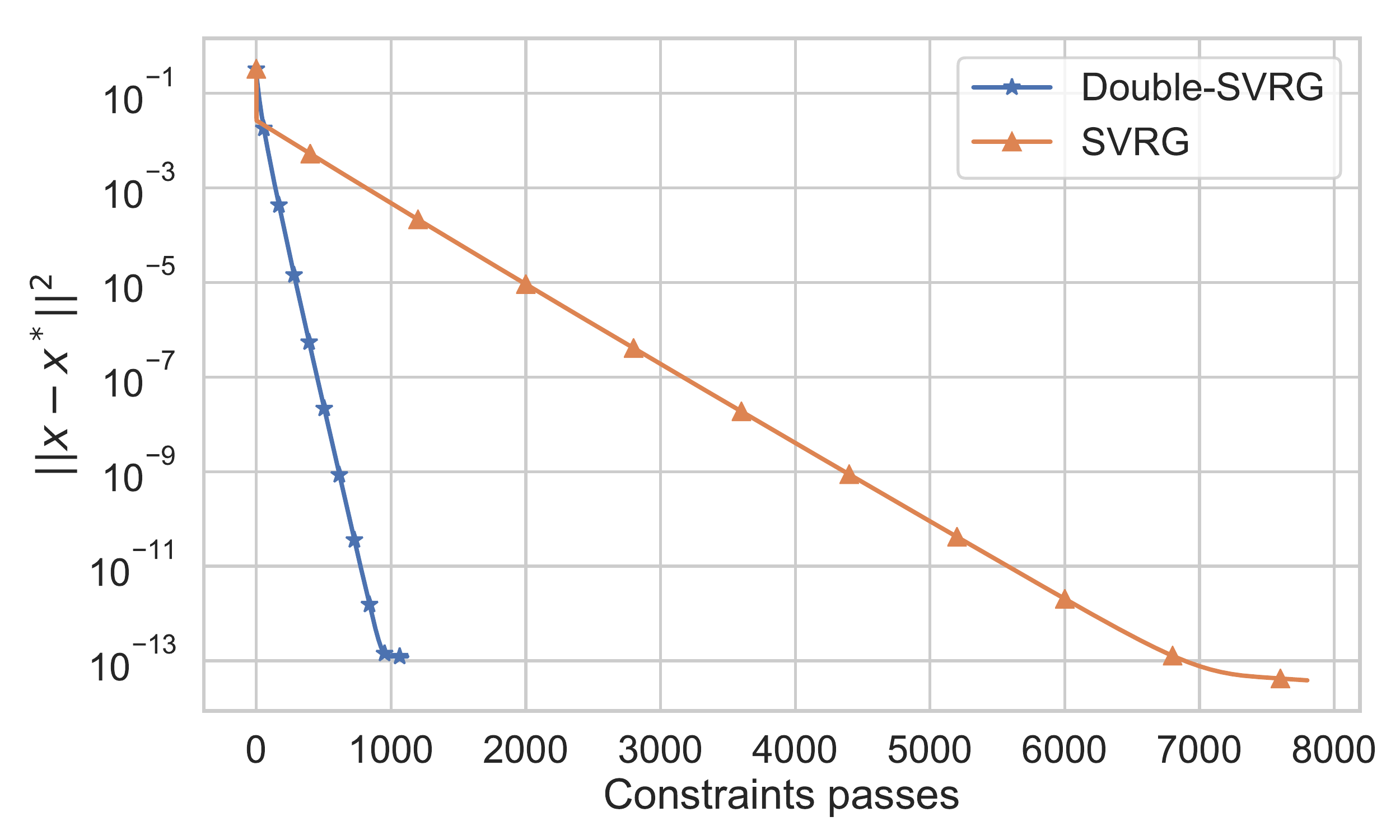}
         \caption{$m=500$, $n=500$}
     \end{subfigure}
     \caption{Comparison of \algname{SVRG} with precise projection onto all constraints (labeled as 'SVRG') to our stochastic version of \algname{SVRG} (labeled as 'Double-SVRG').}
     \label{fig:dif_m_and_n}
\end{figure}

As we can from Figure~\ref{fig:dif_m_and_n}, the trade-off between projections and gradients changes dramatically when $m$ increases. When $m=100$, which implies that the term corresponding to $\mA$ in the complexity is small, the difference is tremendous, partially because mini-batching for \algname{SVRG} improves only part of its complexity~\cite{gower2018stochastic}. In the setting $m=n=500$, we see that the number of data passes  taken by our method to solve the problem is a few times bigger than than that taken by \algname{SVRG}. Clearly, this happens because the term related to $\mA$ becomes dominating in the complexity and \algname{SVRG} uses $m=500$ times more constraints at each iteration than our method.

\graphicspath{{pddy/}}
\chapter{Appendix for Chapter~\ref{chapter:pddy}}
\section{Experiments}
\label{sec:exp}
   
\begin{figure*}[t]
	\centering
	\includegraphics[scale=0.24]{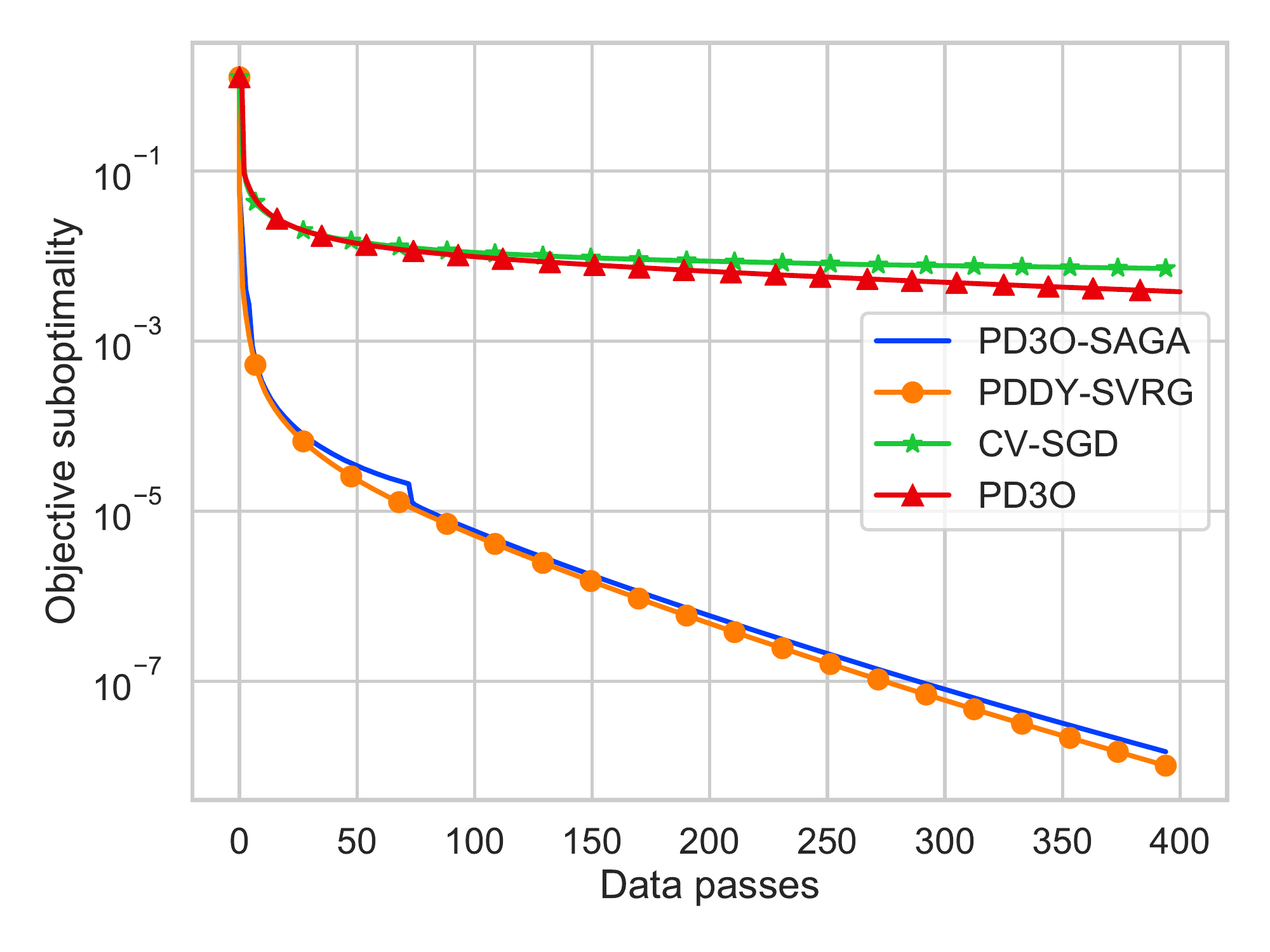}
	\includegraphics[scale=0.24]{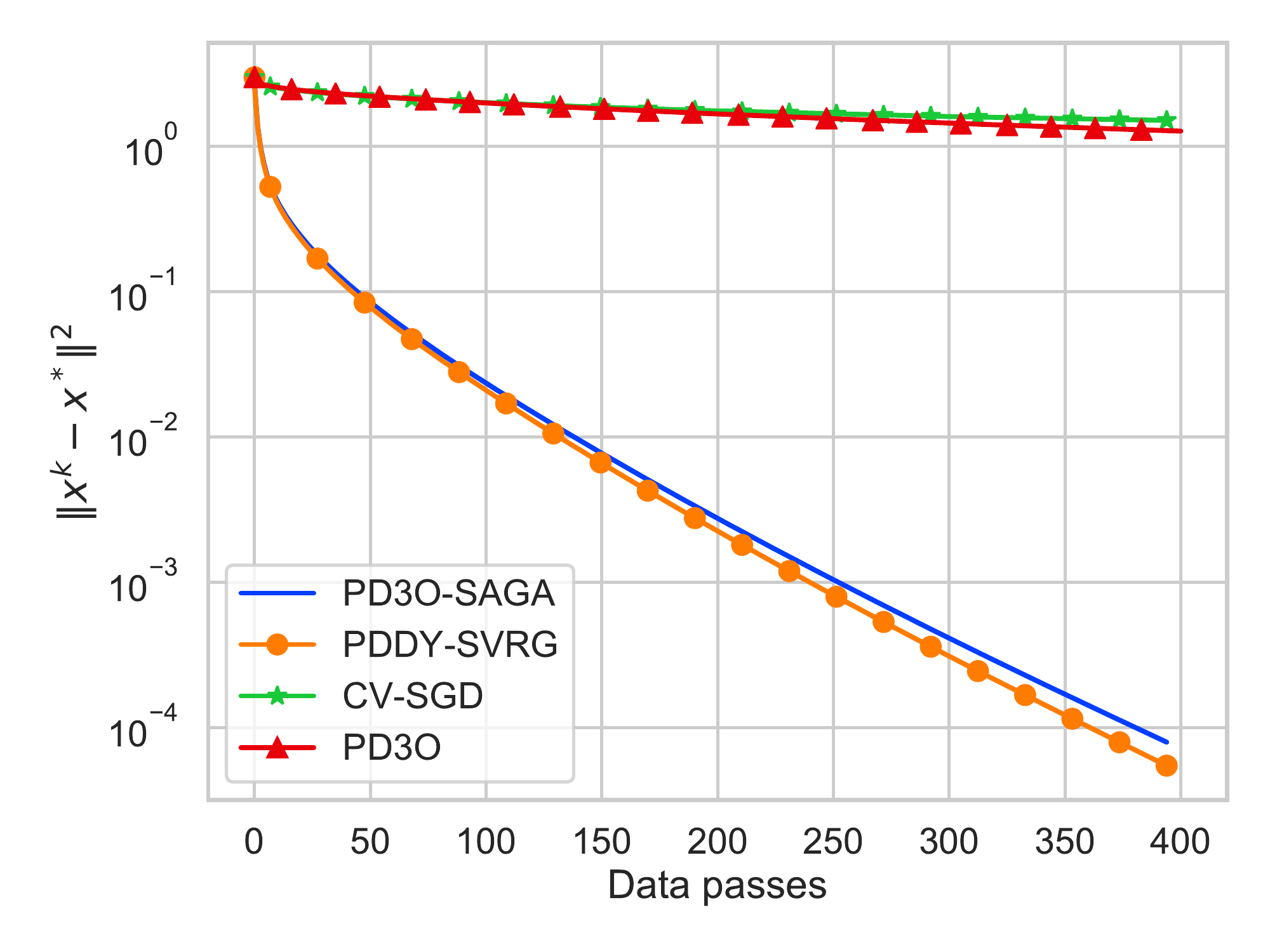}
	\includegraphics[scale=0.24]{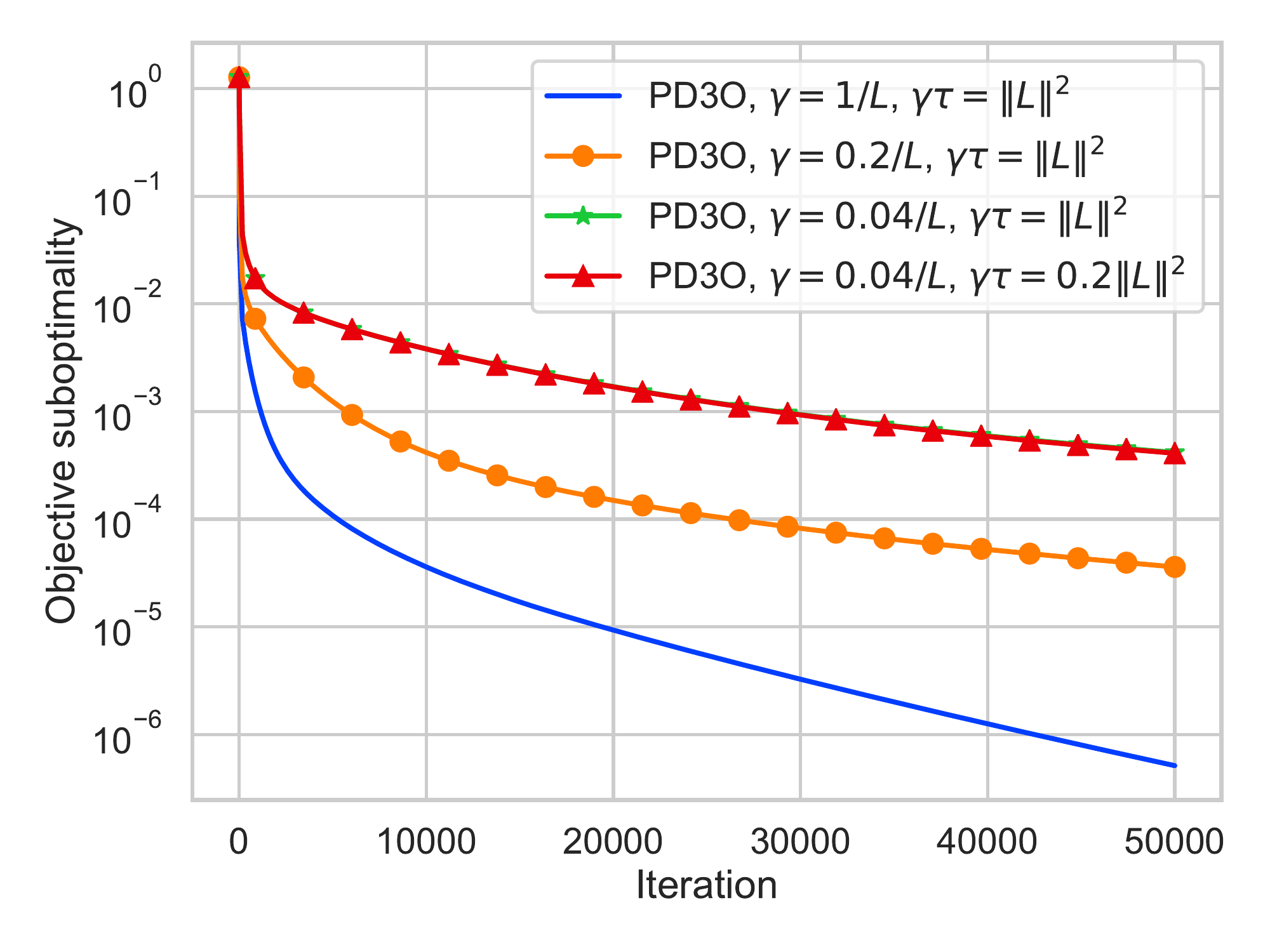}
	\caption{Results for the PCA-Lasso experiment. Left: convergence in the objective, middle: convergence in norm, right: the effect of the stepsizes.}
\end{figure*}

In this section, we present numerical experiments for the \algname{PDDY}, \algname{PD3O} and \algname{Condat--V\~u} (\algname{CV})~\cite[Algorithm 3.1]{condat2013primal} algorithms. \algname{SGD}
was always used with a small $\gamma$, 
such as $\frac{0.01}{\nu}$. where $\nu$ is the smoothness constant of $f$. For stochastic methods, we used a batch size of 16 for better parallelism, while the sampling type is specified in the figures. The stepsizes were tuned with log-grid-search for all methods. We used closed-form expressions to compute $\nu$ for all problems and tuned the stepsizes for all methods by running logarithmic grid search with factor 1.5 over multiples of $\frac{1}{\nu}$.

We observed that the performances of these algorithms are nearly identical, when the same stepsizes are used, so we do not provide their direct comparison in the plots. Instead, we 1) compare different stochastic oracles, 2) illustrate how convergence differs in functional suboptimality and distances, and 3) show how the stepsizes affect the performance.\medskip

\begin{figure*}[t]
	\centering
	\includegraphics[scale=0.24]{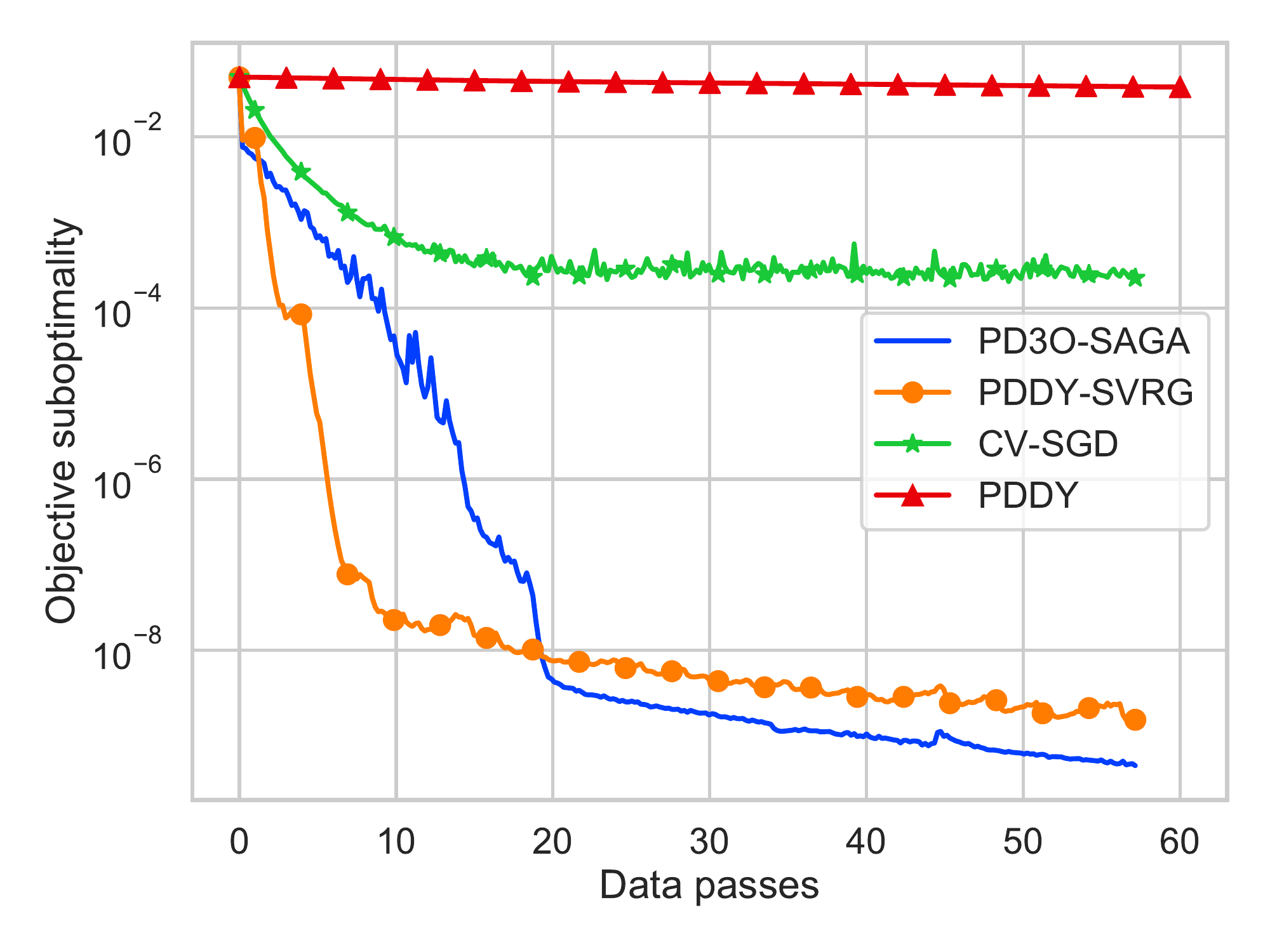}
	\includegraphics[scale=0.24]{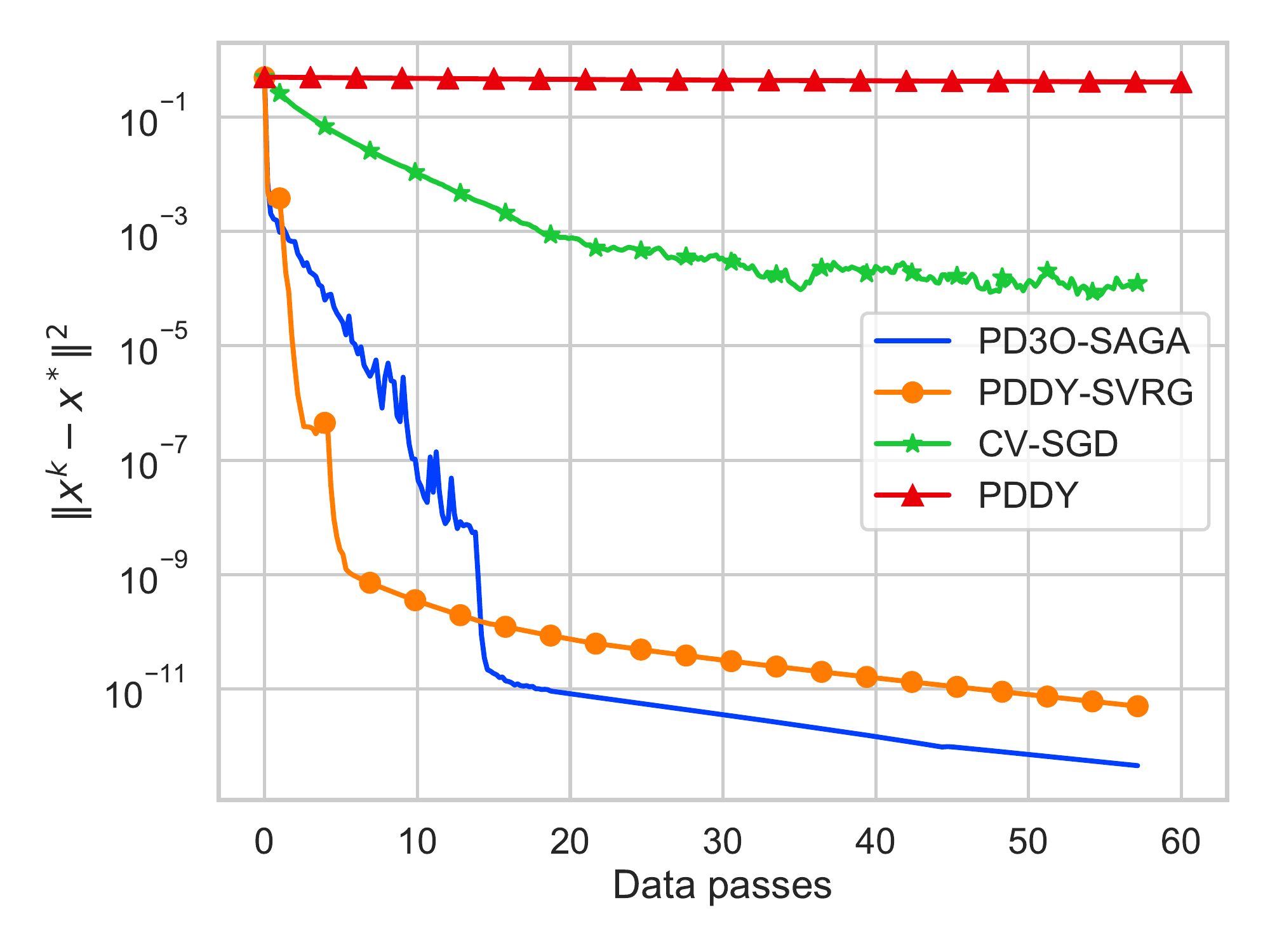}
	\includegraphics[scale=0.24]{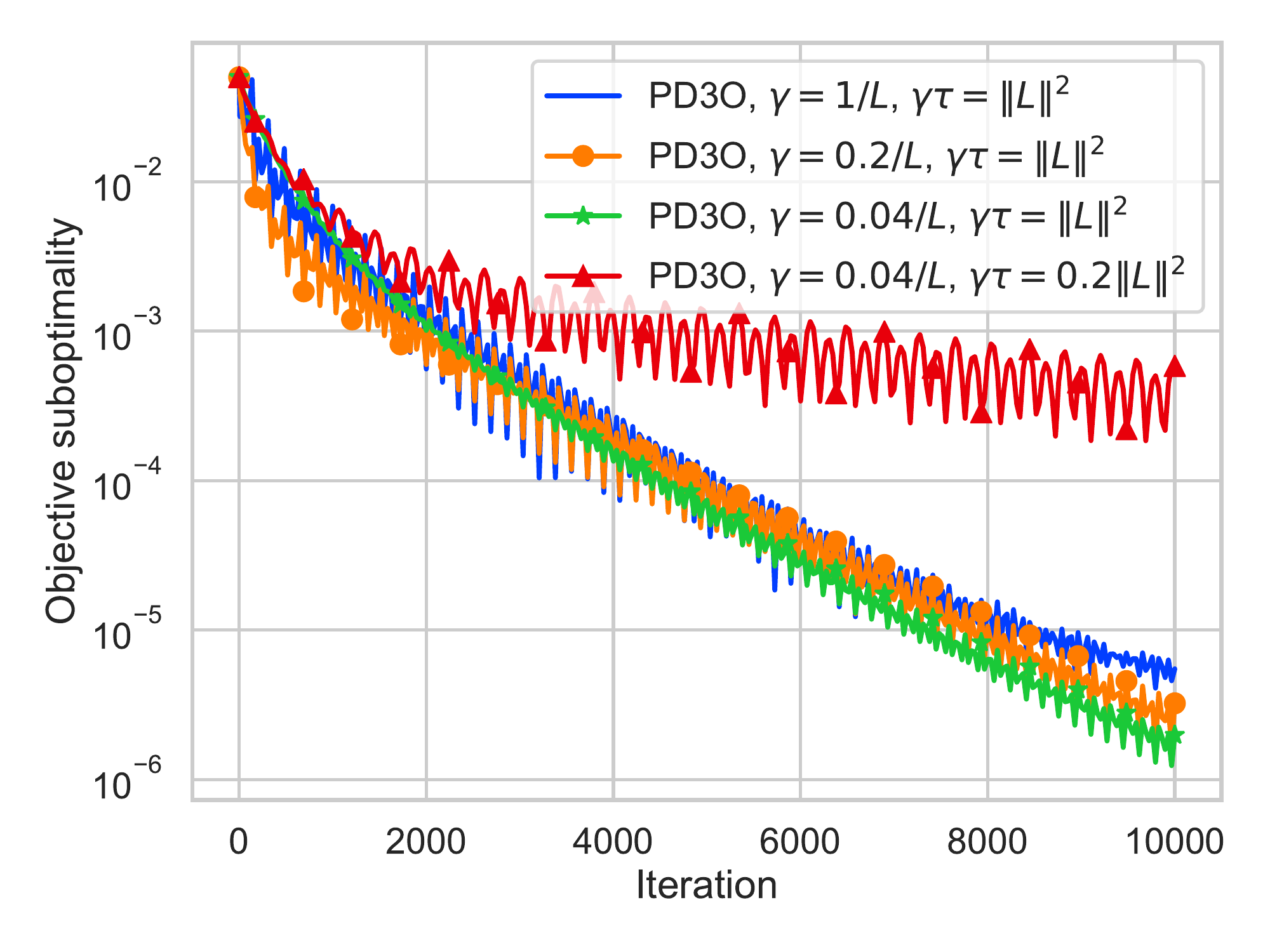}
	\caption{Results for the MNIST experiment. Left: convergence in the objective, middle: convergence in norm, right: the effect of the stepsizes.}
\end{figure*}

\textbf{PCA-Lasso}\ \ \ In a recent work~\cite[Equation~(12)]{tay2018principal} the following difficult PCA-based Lasso problem was introduced: 
$\min_x \frac{1}{2}\|\mW x - a\|^2 + \lambda\|x\|_1 + \lambda_1\sum_{i=1}^m \|\mL _i x\|$, 
where $\mW\in\mathbb{R}^{n\times d}$, $a\in\mathbb{R}^n$, $\lambda, \lambda_1>0$ are given.
 We generate 10 matrices $\mL_i$  randomly with standard normal i.i.d.\ entries, each with 20 rows. $\mW$ and $y$ are taken from the 'mushrooms' dataset from the LIBSVM package~\cite{chang2011libsvm}. We chose $\lambda=\frac{\nu}{10n}$ and $\lambda_1=\frac{2\nu}{nm}$, where $\nu$, the smoothness of $f$, is needed to compensate for the fact that we do not normalize the objective.\medskip

\textbf{MNIST with Overlapping Group Lasso}\ \ \ 
Now we consider the problem where $f$ is the $\ell_2$-regularized logistic loss and a group Lasso penalty. Given the data matrix $\mW\in \mathbb{R}^{n\times d}$ and vector of labels $a\in\{0,1\}^n$, $f(x)=\frac{1}{n}\sum_{i=1}^n f_i(x) + \frac{\lambda}{2}\|x\|^2$ is a finite sum,
	$f_i(x)=-\big(a_i \log \big(h(w_i^\top x)\big) + (1-a_i)\log\big(1-h(w_i^\top x)\big)\big)$,
where, $\lambda=\frac{2\nu}{n}$, $w_i\in\mathbb{R}^d$ is the $i$-th row of $\mW$ and $h:t\to1/(1+e^{-t})$ is the sigmoid function. The non-smooth regularizer, in turn, is given by
	$\lambda_1\sum_{j=1}^m \|x\|_{G_j},$
where $\lambda_1=\frac{\nu}{5n}$, $G_j\subset \{1,\dotsc, p\}$ is a given subset of coordinates and $\|x\|_{G_j}$ is the $\ell_2$-norm of the corresponding block of $x$. To apply splitting methods, we use $\mL=(\mI_{G_1}^\top, \dotsc, \mI_{G_m}^\top)^\top$, where $\mI_{G_j}$ is the operator that takes $x\in\mathbb{R}^d$ and returns only the entries from block $G_j$. Then, we can use $H(y)=\lambda_1\sum_{j=1}^m\|y\|_{G_j}$, which is separable in $y$ and, thus, proximable.
We use the MNIST datasetw~\cite{lecun2010mnist} of 70000 black and white $28\times 28$ images. For each pixel, we add a group of pixels $G_j$ adjacent to it, including the pixel itself. Since there are some border pixels, groups consist of 3, 4 or 5 coordinates, and there are 784 penalty terms in total.\medskip

\begin{figure*}[t]
	\centering
	\includegraphics[scale=0.24]{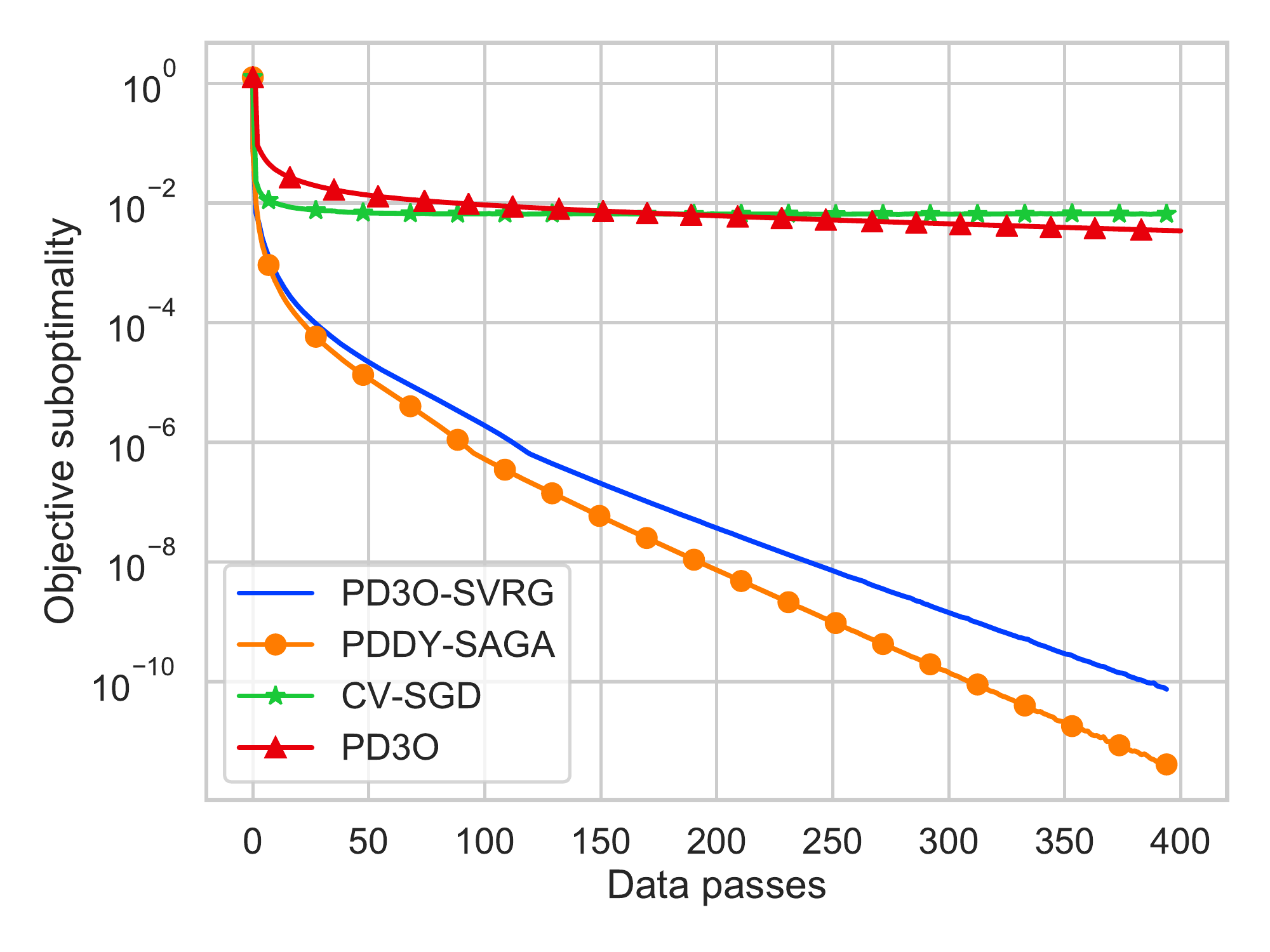}
	\includegraphics[scale=0.24]{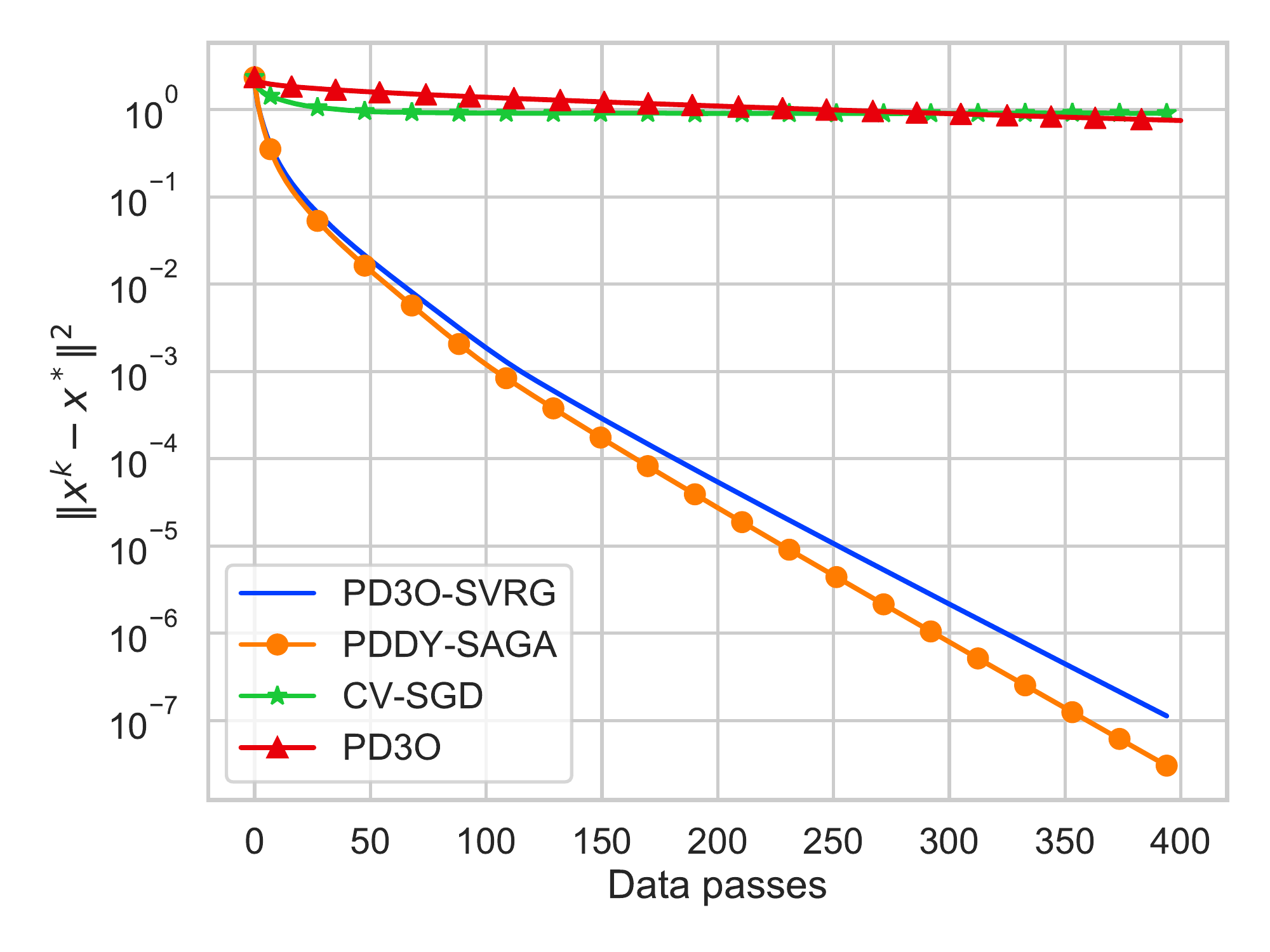}
	\includegraphics[scale=0.24]{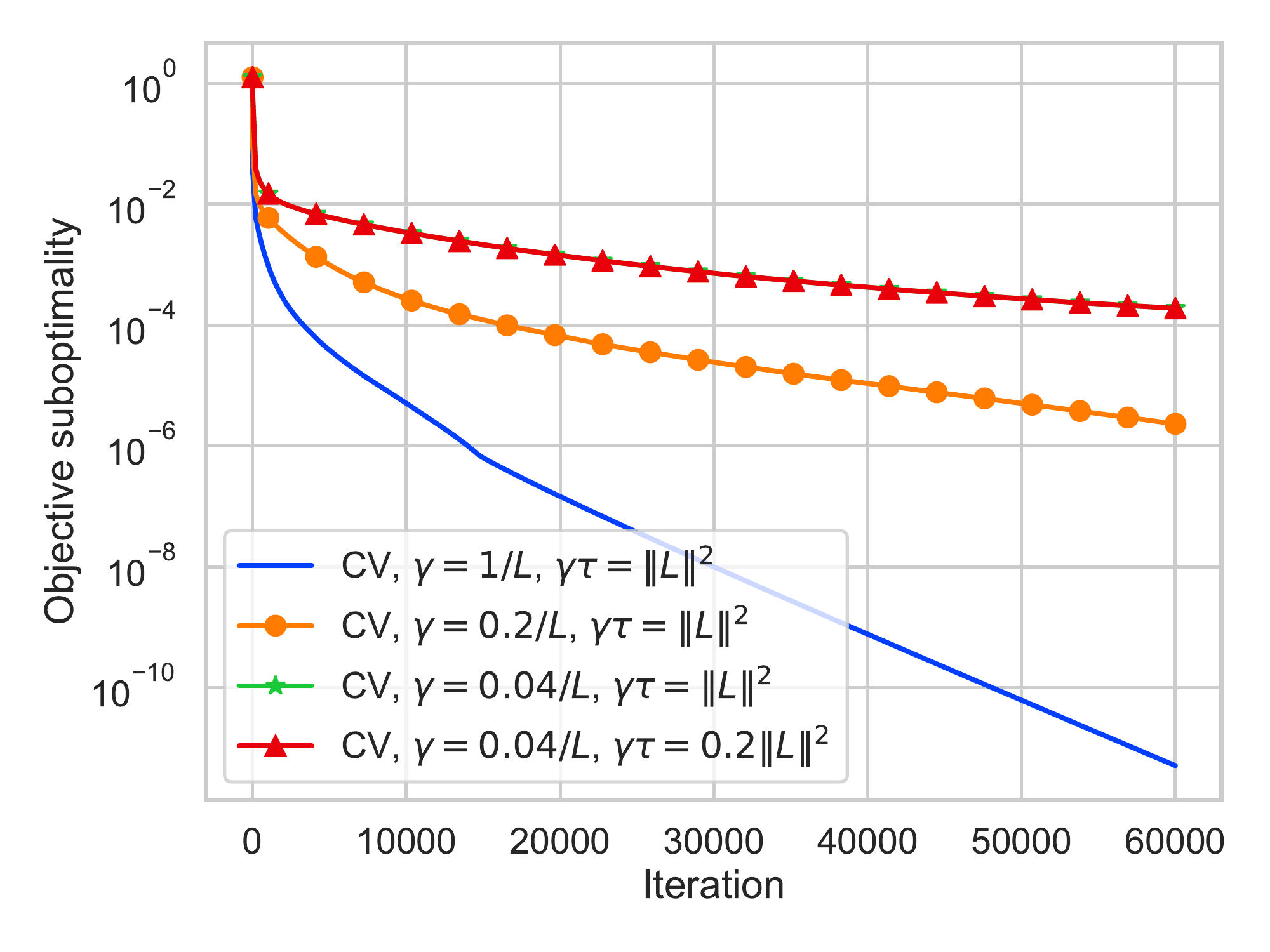}
	\caption{Results for the Fused Lasso experiment. Left: convergence with respect to the objective function, middle: convergence in norm, right: illustration of the effect of the stepsizes.}
	\label{fig:fused}
\end{figure*}

\textbf{Fused Lasso Experiment}\ \ \ In the Fused Lasso problem, we are given a feature matrix $\mW\in\mathbb{R}^{n\times d}$ and an output vector $a$, which define the least-squares smooth objective $f(x)=\frac{1}{2}\|\mW x-a\|^2$. This function is regularized with $\frac{\lambda}{2}\|x\|^2$ and $\lambda_1\|\mD x\|_1$, where $\lambda=\frac{\nu}{n}$, $\lambda_1=\frac{\nu}{10n}$ and $\mD\in\mathbb{R}^{(d-1)\times d}$ has entries $\mD_{i,i}=1$, $\mD_{i, i+1}=-1$, for $i=1,\dotsc, p-1$, and $\mD_{ij}=0$ otherwise. We use the 'mushrooms' dataset from the LIBSVM package. Our numerical findings for this problem are very similar to the ones for PCA-Lasso. In particular, larger values of $\gamma$ seem to perform significantly better and the value of the objective function does not oscillate, unlike in the MNIST experiment. The results are shown in Figure~\ref{fig:fused}. The proposed \algname{Stochastic PDDY} algorithm with the \algname{SAGA} estimator performs best in this setting.\medskip

\textbf{Summary of results}\ \ \ 
We can see from the plots that stochastic updates make the convergence extremely faster, sometimes even without variance reduction. The stepsize plots suggest that it is best to keep $\gamma\tau$ close to $\frac{1}{\|\mL\|^2}$, 
while the optimal value of $\gamma$ might sometimes be smaller than $\frac{1}{\nu}$. This is especially clearly seen from the fact that \algname{SGD} works sufficiently fast even despite using $\gamma$ inversely proportional to the number of iterations.

\section{Proofs Related to Primal--Dual Optimality}
\subsection{Optimality conditions}
Let $x^\star$ be a minimizer of Problem~\eqref{eq:original-pb}. Assuming a standard qualification condition, for instance that $0$ belongs to the relative interior of $\dom(H) - \mL\dom(\psi )$, then for every $x \in \cX$,
\begin{equation*}
\partial(f+\psi +H\circ \mL)(x) = \nabla f(x) + \partial \psi(x) + \mL^* \partial H(\mL  x),\end{equation*}
see for instance Theorem 16.47 of~\cite{bauschke2017convex}. Then,
\begin{align*}
    \label{eq:existsy}
        &x^\star \in \argmin_{x \in \mathcal{X}} \{ f(x)+\psi(x)+H(\mL x) \}\\
        \Leftrightarrow{}& 0 \in \nabla f(x^\star) + \partial \psi(x^\star) + \mL^* \partial H (\mL  x^\star)\\
        \Leftrightarrow{}& \exists y^\star \in \partial H (\mL  x^\star) \text{ such that } 0 \in \nabla f(x^\star) + \partial \psi(x^\star) + \mL^* y^\star\\
        \Leftrightarrow{}& \exists y^\star \in \mathcal{Y} \text{ such that } 0 \in \nabla f(x^\star) + \partial \psi(x^\star) + \mL^* y^\star \text{ and } 0 \in -\mL  x^\star + \partial H^*(y^\star),
    \end{align*}
where we used $\partial H^* = (\partial H)^{-1}$.

\subsection{Proof of Lemma~\ref{lem:duality-gap}}

Using the optimality conditions~\eqref{eq:saddle}, we have
\begin{align*}
    D_f(x,x^\star)+D_\psi(x,x^\star) &= (f+\psi )(x) - (f+\psi )(x^\star) -\ps{\nabla f(x^\star) + r^\star,x-x^\star}\notag\\
    &= (f+\psi )(x) - (f+\psi )(x^\star) +\ps{\mL^* y^\star,x-x^\star}\\
    &= (f+\psi )(x) - (f+\psi )(x^\star) +\ps{y^\star,\mL x} - \ps{y^\star, \mL x^\star}.\notag
\end{align*}
We also have
\begin{align*}
    D_{H^*}(y,y^\star) &= H^*(y) - H^*(y^\star) -\ps{h^\star,y-y^\star}\notag\\
    &= H^*(y) - H^*(y^\star) -\ps{\mL x^\star,y-y^\star}\\
    &= H^*(y) - H^*(y^\star) -\ps{\mL x^\star,y} +\ps{y^\star,\mL x^\star}.\notag
\end{align*}
Summing the two last equations, we have
\begin{align*}
&D_f(x,x^\star)+D_\psi(x,x^\star)+D_{H^*}(y,y^\star)\notag\\
{}={}& (f+\psi )(x) - (f+\psi )(x^\star) + H^*(y) - H^*(y^\star) -\ps{\mL x^\star,y} +\ps{y^\star,\mL x}\\
{}={}& \cL(x,y^\star) - \cL(x^\star,y).\notag
\end{align*}

\section{Proof of Lemma~\ref{lem:funda-DYS}}
\label{sec:funda-DYS}
Since $z^k = J_{\gamma \tilde{B}}(v^k)$, $z^k \in v^k - \gamma \tilde{B}(z^k)$ by definition of the resolvent operator. Therefore, there exists $b^k \in \tilde{B}(z^k)$ such that $z^k = v^k - \gamma b^k$. Similarly, 
\begin{equation*}
u^{k+1} \in 2z^k - v^k - \gamma \tilde{C}(z^{k}) - \gamma \tilde{A}(u^{k+1}) = v^k - 2\gamma b^k - \gamma \tilde{C}(z^{k}) - \gamma \tilde{A}(u^{k+1}).
\end{equation*}
 Therefore, there exists $a^{k+1} \in \tilde{A}(u^{k+1})$ such that
\begin{equation}
    \label{eq:DYS2}
     \left\{\begin{array}{l}
        z^{k} = v^k - \gamma b^k \\
        u^{k+1} = v^k - 2 \gamma b^k - \gamma \tilde{C}(z^{k}) - \gamma a^{k+1} \\
        v^{k+1} = v^k + u^{k+1} - z^{k}.
        \end{array}\right.
\end{equation}
Moreover,
\begin{equation}
    \label{eq:DYS3}
    v^{k+1} = v^{k} - \gamma b^k - \gamma \tilde{C}(z^{k}) - \gamma a^{k+1}.
\end{equation}
Similarly, there exist $a^\star \in \tilde{A}(u^\star), b^\star \in \tilde{B}(z^\star)$ such that
\begin{equation}
    \label{eq:DYS2-fix}
    \left\{\begin{array}{l}
        z^\star = v^\star - \gamma b^\star \\
        u^\star = v^\star - 2 \gamma b^\star - \gamma \tilde{C}(z^\star) - \gamma a^\star \\
        v^\star = v^\star + u^\star - z^\star,
        \end{array}\right.
\end{equation}
and
\begin{equation}
    \label{eq:DYS3-fix}
    v^\star = v^\star - \gamma b^\star - \gamma \tilde{C}(z^\star) - \gamma a^\star.
\end{equation}
Therefore, using~\eqref{eq:DYS3} and~\eqref{eq:DYS3-fix},
\begin{align*}
    \|v^{k+1} - v^\star\|^2 ={}& \|v^k - v^\star\|^2 -2\gamma\ps{a^{k+1}+b^k + \tilde{C}(z^{k}) - \left(a^{\star}+b^\star + \tilde{C}(z^{\star})\right),v^{k} - v^\star}\\
    &+\gamma^2\|a^{k+1}+b^k + \tilde{C}(z^{k}) - \left(a^{\star}+b^\star + \tilde{C}(z^{\star})\right)\|^2.\notag
\end{align*}
By expanding the last square at the right-hand side, and by using~\eqref{eq:DYS2} and~\eqref{eq:DYS2-fix} in the inner product, we get
\begin{align*}
    \|v^{k+1} - v^\star\|^2 ={}& \|v^k - v^\star\|^2\notag\\
    &-2\gamma\ps{b^k + \tilde{C}(z^{k})- \left(b^\star + \tilde{C}(z^\star)\right),z^{k} - z^\star}\notag\\
    &-2\gamma\ps{a^{k+1}- a^\star,u^{k+1} - u^\star}\notag\\
    &-2\gamma\ps{b^k + \tilde{C}(z^{k})- \left(b^\star + \tilde{C}(z^\star)\right),\gamma b^k - \gamma b^\star}\notag\\
    &-2\gamma\ps{a^{k+1}- a^\star,2 \gamma b^k + \gamma \tilde{C}(z^{k}) + \gamma a^{k+1} - \left(2 \gamma b^\star + \gamma \tilde{C}(z^{\star}) + \gamma a^\star\right)}\\
    &+\gamma^2\|a^{k+1}+b^k- \left(a^\star+b^\star\right)\|^2\notag\\
    &+\gamma^2\|\tilde{C}(z^{k})- \tilde{C}(z^\star)\|^2\notag\\
    &+2\gamma^2\ps{a^{k+1}+b^k - \left(a^\star+b^\star\right),\tilde{C}(z^{k})- \tilde{C}(z^\star)}.\notag
\end{align*}
Then, the last five terms at the right-hand side simplify to
\begin{equation*}
\gamma^2\|\tilde{C}(z^{k}) - \tilde{C}(z^\star)\|^2-\gamma^2\|a^{k+1}+b^k - \left(a^\star+b^\star\right)\|^2,
\end{equation*}
and we get the result.

\section{Proofs Related to the \algname{PDDY} and \algname{PD3O} Algorithms}

\subsection{Resolvent calculus}

For the sake of completeness, we reproduce a resolvent computation that can be found in~\cite{con19}, showing that the \algname{PD3O} algorithm is an instance of \algname{DYS}.

In the notations of Section \ref{sec:pdalgos}, let us state the lemma:
\begin{lemma}
    \label{lem:res-pd}
    $J_{\gamma \mP^{-1}A}$ maps $(x,y)$ to $(x',y')$, such that
    \begin{align}
 &\left\lfloor
 \begin{array}{l}
 y'
= \mathrm{prox}_{\tau H^*}\big(y+\tau \mL (x - \gamma \mL^* y)\big)\\
x'=x -\gamma \mL^* y'.
\end{array}\right.
\end{align}
\end{lemma}
\begin{proof}
Let $(x,y)$ and $(x',y') \in \cZ$, such that
\begin{equation*}
\mP \begin{bmatrix} x' - x\\ y' - y\end{bmatrix} \in - \gamma \begin{bmatrix} \!\!\!\!\!\!&{}  + \mL^* y' \\ -\mL  x' &{} + \partial H^*(y')\end{bmatrix},
\end{equation*}
where
\begin{equation*}
\mP = \begin{bmatrix} \mI &  0 \\ 0 & \frac{\gamma}{\tau}\mI - \gamma^2 \mL \mL^*   \end{bmatrix}.
\end{equation*}

We shall express $(x',y')$ as a function of $(x,y)$.
First,
$$x' = x - \gamma \mL^* y'.$$
Moreover, $y'$ is given by
\begin{align*}
    \left(\frac{\gamma}{\tau}\mI - \gamma^2  \mL\mL^*\right)(y') &\in \left(\frac{\gamma}{\tau}\mI - \gamma^2 \mL \mL^*\right)(y) +\gamma \mL x' -\gamma \partial H^*(y')\\
    &\in \left(\frac{\gamma}{\tau}\mI - \gamma^2 \mL \mL^*\right)(y) +\gamma \mL\left(x - \gamma \mL^* y'\right) -\gamma \partial H^*(y').
\end{align*}
Therefore, the term $\gamma^2 \mL\mL^* y'$ disappears from both sides and
\begin{align*}
    y' &\in y -\gamma\tau \mL\mL^* y - \tau \partial H^*(y') +\tau \mL x.
\end{align*}
Finally,
\begin{equation*}
y -\gamma\tau \mL \mL^* y +\tau \mL x \in y' + \tau \partial H^*(y'),
\end{equation*}
and
\begin{equation*}
y' = \prox_{\tau H^*}(y -\gamma\tau \mL \mL^* y +\tau \mL x).
\end{equation*}
\end{proof}

\section{Proofs Related to the \algname{Condat--V\~u} Algorithm}

\subsection{Resolvent calculus}
The results of this section rely on the following resolvent computation, which is new to our knowledge.

In the notations of Section \ref{seccv}, let us state the lemma:
\begin{lemma}
    \label{lem:res-cv}
    $J_{\gamma \mQ^{-1}A}$ maps $(x,y)$ to $(x',y')$, such that
\begin{align}
 &\left\lfloor
    \begin{array}{l}
    x'
= \mathrm{prox}_{\tau \psi }\big((\mI-\tau\gamma \mL^* \mL) x-\tau \mL^* y \big),\\
y'=y +\gamma \mL x'.
\end{array}\right.
\end{align}
\end{lemma}
\begin{proof}
Let $(x,y)$ and $(x',y') \in \cZ$ be such that
\begin{equation*}
\mQ \begin{bmatrix} x' - x\\ y' - y\end{bmatrix} \in - \gamma \begin{bmatrix} \partial \psi(x')\!\!\!\!\!\!&{}  + \mL^* y' \\ -\mL x' &\end{bmatrix},
\end{equation*}
where
\begin{equation*}
\mQ = \begin{bmatrix} \frac{\gamma}{\tau}\mI - \gamma^2 \mL^* \mL&  0 \\ 0 & \mI  \end{bmatrix}.
\end{equation*}
We shall express $(x',y')$ as a function of $(x,y)$.
First,
$$y' = y + \gamma \mL x'.$$
Moreover, $x'$ is given by
\begin{align*}
    \left(\frac{\gamma}{\tau}\mI - \gamma^2 \mL^* \mL\right)(x') &\in \left(\frac{\gamma}{\tau}\mI - \gamma^2 \mL^* \mL\right)(x) - \gamma \partial \psi(x') -\gamma \mL^* y'\\
    &\in \left(\frac{\gamma}{\tau}\mI - \gamma^2 \mL^* \mL\right)(x) - \gamma \partial \psi(x') -\gamma \mL^* y - \gamma^2 \mL^* \mL x'.
\end{align*}
Therefore, the term $\gamma^2 \mL^* \mL x'$ disappears from both sides and
\begin{align*}
    x' &\in x -\gamma\tau \mL^* \mL x - \tau \partial \psi(x') -\tau \mL^* y.
\end{align*}
Finally,
\begin{equation*}
x -\gamma\tau \mL^* \mL x -\tau \mL^* y \in x' + \tau \partial \psi(x'),
\end{equation*}
and
\begin{equation*}
x' = \prox_{\tau \psi }(x -\gamma\tau \mL^* \mL x -\tau \mL^* y).
\end{equation*}
\end{proof}

\subsection{Algorithm 3.2 of~\cite{condat2013primal} as an instance of Davis--Yin Splitting}
We apply $\algname{DYS}(\mQ^{-1}B,\mQ^{-1}A,\mQ^{-1}C)$:
\begin{equation*}
\left\lfloor\begin{array}{l}
    x^{k} = p^k\\
    y^k = \prox_{\gamma H^*}(q^k)\\
    s^{k+1} = \prox_{\tau \psi }\Big((\mI - \gamma\tau \mL^* \mL) \big(2x^k - p^k - \gamma (\frac{\gamma}{\tau}\mI - \gamma^2 \mL^* \mL)^{-1} \nabla f(x^k)\big) - \tau \mL^* \left(2y^k - q^k\right)\Big) \\
    h^{k+1} = 2y^k - q^k + \gamma \mL s^{k+1}\\
    p^{k+1} = p^k + s^{k+1} - x^k\\
    q^{k+1} = q^k + h^{k+1} - y^k.
        \end{array}\right.
\end{equation*}
Note that $p^{k} = s^{k} = x^k$ and $q^{k+1} = y^k + \gamma \mL x^{k+1}$. So,
\begin{align*}
    x^{k+1} &= \prox_{\tau \psi }\big(x^k - \tau\nabla f(x^k) - \tau \mL^* (2y^k - q^k + \gamma \mL x^k)\big),\\
   &= \prox_{\tau \psi }\big(x^k - \tau\nabla f(x^k) - \tau \mL^* (2y^k - y^{k-1})\big),
\end{align*}
and
\begin{equation*}
    y^{k} = \prox_{\gamma H^*}(y^{k-1} + \gamma \mL x^{k}).
\end{equation*}
The sequence $(y^{k},x^{k+1})$ follows the updates of Algorithm 3.2 in~\cite{condat2013primal}.

\subsection{Algorithm 3.1 of~\cite{condat2013primal} as an instance of Davis--Yin Splitting}\label{ap:cv1}
We apply $\algname{DYS}(\mQ^{-1}A,\mQ^{-1}B,\mQ^{-1}C)$, which yields:
\begin{equation*}
\left\lfloor\begin{array}{l}
    x^{k} = \prox_{\tau \psi }\big((\mI - \gamma\tau \mL^* \mL) p^k - \tau \mL^* q^k\big)\\
    y^k = q^k + \gamma \mL x^k\\
    s^{k+1} = 2x^k - p^k - \gamma (\frac{\gamma}{\tau}\mI - \gamma^2 \mL^* \mL)^{-1} \nabla f(x^k)\\
    h^{k+1} = \prox_{\gamma H^*}(2y^k - q^k)\\
    p^{k+1} = p^k + s^{k+1} - x^k\\
    q^{k+1} = q^k + h^{k+1} - y^k.
    \end{array}\right.
\end{equation*}
Note that $p^{k+1} = x^k - (\frac{1}{\tau}\mI - \gamma \mL^* \mL)^{-1} \nabla f(x^k)$ and $q^{k+1} = h^{k+1} - \gamma \mL x^k$. Thus,
\begin{align*}
    x^{k} &= \prox_{\tau \psi }\big(x^{k-1} - \gamma\tau \mL^* \mL x^{k-1} - \tau\nabla f(x^{k-1}) - \tau \mL^* (h^{k} - \gamma \mL x^{k-1})\big)\notag\\
    &= \prox_{\tau \psi }\big(x^{k-1} - \tau \nabla f(x^{k-1}) - \tau \mL^* h^{k}\big),
\end{align*}
and
\begin{align*}
    h^{k+1} &= \prox_{\gamma H^*}(2y^k - q^k)\notag\\
    &= \prox_{\gamma H^*}(q^k + 2\gamma \mL x^k)\\
    &= \prox_{\gamma H^*}(h^{k} + \gamma \mL (2x^k - x^{k-1})).\notag
\end{align*}
The sequence $(x^k,h^{k+1})$ follows the updates of Algorithm 3.1 of~\cite{condat2013primal}.

\subsection{Cocoercivity parameter of $\mQ^{-1} C$}
Since $\mK$ is positive definite, $\|\mK^{-1/2}\|^2 = \|\mK^{-1}\|$. Let $z = (x,y), z'=(x',y') \in \cZ$. Then,
\begin{align*}
    \|\mQ^{-1}C(z) - \mQ^{-1}C(z')\|_\mQ^2 &= \|\mK^{-1}\nabla f(x) - \mK^{-1} \nabla f(x')\|_\mK^2\\
    &= \|\mK^{-1/2}\nabla f(x) - \mK^{-1/2} \nabla f(x')\|^2\\
    &\leq  \|\mK^{-1/2}\|^2 \|\nabla f(x) - \nabla f(x')\|^2\\
    &=  \|\mK^{-1}\| \|\nabla f(x) - \nabla f(x')\|^2\\
    &\leq \|\mK^{-1}\| \nu \ps{\nabla f(x) - \nabla f(x'),x-x'}\\
    &= \|\mK^{-1}\| \nu \ps{\mK^{-1}\nabla f(x) - \mK^{-1}\nabla f(x'),x-x'}_\mK\\
    &= \|\mK^{-1}\| \nu \ps{\mQ^{-1}C(z) - \mQ^{-1}C(z'),z-z'}_\mQ.
\end{align*}
Since $\|\mK^{-1}\|$ is the inverse of the smallest eigenvalue of $\mK$, $$\|\mK^{-1}\|= \frac{1}{\frac{\gamma}{\tau} - \gamma^2\|\mL\|^2}.$$
Therefore, $\mQ^{-1}C$ is $\xi = \frac{\frac{\gamma}{\tau} - \gamma^2\|\mL\|^2}{\nu}$-cocoercive. Moreover, the condition $\gamma < 2\xi$ of Lemma~\ref{lem:DYS-cv} is equivalent to $\nu/2 < \frac{1}{\tau} - \gamma\|\mL\|^2$, which is exactly the condition (i) of Theorem 3.1 in~\cite{condat2013primal}.

\section{Linear Convergence Results}
In this section, we provide linear convergence results for the \algname{Stochastic PD3O} and the \algname{Stochastic PDDY} algorithm, in addition to Theorem \ref{th:LV0}.

\paragraph{About the linear convergence of \algname{DYS}$(\tilde{A},\tilde{B},\tilde{C})$.} In general, operator splitting methods like \algname{DYS}$(\tilde{A},\tilde{B},\tilde{C})$ require $\tilde{A}+\tilde{B}+\tilde{C}$ to be strongly monotone to converge linearly. Besides, to converge linearly \textit{in general}, \algname{DYS}$(\tilde{A},\tilde{B},\tilde{C})$ requires the stronger assumption that $\tilde{A}$ or $\tilde{B}$ or $\tilde{C}$ is strongly monotone \footnote{This assumption is stronger than assuming $\tilde{A}+\tilde{B}+\tilde{C}$ strongly monotone.} and that $\tilde{A}$ or $\tilde{B}$ is cocoercive\footnote{For linear convergence, this assumption is proved to be necessary in general in~\cite{davis2017three}.}, see~\cite{davis2017three}.

Recall that the \algname{PDDY} algorithm is equivalent to 
\algname{DYS}$(\mP^{-1}B,\mP^{-1}A,\mP^{-1}C)$ and the \algname{PD3O} algorithm is equivalent to \algname{DYS}$(\mP^{-1}A,\mP^{-1}B,\mP^{-1}C)$, see Section~\ref{sec:pdalgos}. However, $\mP^{-1}A$, $\mP^{-1}B$ and $\mP^{-1}C$ are not strongly monotone. In spite of this, we shall prove the linear convergence of the (stochastic) \algname{PDDY} algorithm and the (stochastic) \algname{PD3O} algorithm. 

\paragraph{Our assumptions.} For both algorithms, we shall make the weaker assumption that $\mP^{-1}A+\mP^{-1}B+\mP^{-1}C$ is strongly monotone (which turns out to be equivalent to assuming $M = A+B+C$ strongly monotone, i.e., $f+\psi $ strongly convex and $H$ smooth). Indeed, the algorithms need to be contractive in both the primal and the dual space. For instance, \algname{Chambolle--Pock}~\cite{chambolle2011first,chambolle2016ergodic} algorithm, which is a particular case of the \algname{PD3O} and the \algname{PDDY} algorithms, requires $\psi $ strongly convex and $H$ smooth to converge linearly in general. Therefore, we shall assume $M$ strongly monotone for both algorithms, which is weaker than assuming $\mP^{-1}A$, $\mP^{-1}B$ or $\mP^{-1}C$ strongly monotone.

Moreover, for the \algname{PD3O} algorithm we shall add a cocoercivity assumption, as suggested by the general linear convergence theory of \algname{DYS}. More precisely, we shall assume $\psi $ smooth (i.e., $\mP^{-1}B$ cocoercive). Our first result is therefore an extension of~\cite[Theorem 3]{yan2018new} to the stochastic setting.

For the \algname{PDDY} algorithm, we shall not make any cocoercivity assumption on $\mP^{-1}A$ and $\mP^{-1}B$, but we shall make an assumption of the stepsize instead. Our second result is new, even in the deterministic case $g^{k+1} = \nabla f(x^k)$.

We denote by $\|\cdot\|_{\gamma,\tau}$ the norm induced by $\frac{\gamma}{\tau}\mI - \gamma^2 \mL \mL^*$ on $\cY$.

\subsection{The \algname{Stochastic PD3O} algorithm}

\begin{theorem}[$M$ strongly monotone and $\psi $ smooth]
    \label{th:Hsmooth:PD3O}
    Suppose that Assumption~\ref{as:sto-grad} holds. Also, suppose that $H$ is $1/\mu_{H^*}$-smooth, $f$ is $\mu_f$-strongly convex, and $\psi $ is $\mu_\psi $-strongly convex and $\lambda$-smooth, where $\mu \eqdef \mu_f + 2\mu_\psi  >0$ and $\mu_{H^*} >0$. For every $\kappa > \beta/\rho$ and every $\gamma, \tau >0$ such that $\gamma \leq 1/(\alpha+\kappa\delta)$ and $\gamma\tau\|\mL\|^2 < 1$, define
    \begin{equation}
        \label{eq:lyapunov-smooth}
        V^k \eqdef \|p^{k} - p^\star\|^2+ \left(1+2\tau\mu_{H^*}\right)\|y^{k} - y^\star\|_{\gamma,\tau}^2 + \kappa\gamma^2 \sigma_{k}^2,
    \end{equation}
    and
    \begin{equation}
        \label{eq:rate-smooth}
        r \eqdef \max\left(1-\frac{\gamma\mu}{(1+\gamma \lambda)^2},\left(1-\rho+\frac{\beta}{\kappa}\right),\frac{1}{1+2\tau\mu_{H^*}}\right).
    \end{equation}
    Then,
    \begin{equation}
        \ec{ V^{k}}
        \leq r^k V^0.
    \end{equation}
\end{theorem}
Under smoothness and strong convexity assumptions, Theorem~\ref{th:Hsmooth:PD3O} implies  linear convergence of the dual variable $y^k$ to $y^\star$, with convergence rate given by $r$.
Since $\|x^k - x^\star\| \leq \|p^k - p^\star\|$, Theorem~\ref{th:Hsmooth:PD3O} also implies linear convergence of the primal variable $x^k$ to $x^\star$, with same convergence rate. 
For primal--dual algorithms to converge linearly on Problem \eqref{eq:original-pb} with any $\mL$, it seems unavoidable that the primal term $f+\psi$ is strongly convex and that the dual term $H^*$ is strongly convex too; this means that $H$ must be smooth. We can notice that if $H$ is smooth, it is tempting to use its gradient instead of its proximity operator. We can then use the proximal gradient algorithm with $\nabla (f+H\circ \mL)(x)=\nabla f(x) + \mL^*\nabla H (\mL x)$. However, in practice, it is often faster to use the proximity operator instead of the gradient, see a recent analysis of this topic \cite{combettes2019proximal}.

\begin{remark}[Particular case] In the case where $\psi = H =0$ and $\mL =0$, then the \algname{Stochastic PD3O} algorithm boils down to \algname{Stochastic Gradient Descent} (\algname{SGD}), where the stochastic gradient oracle satisfies Assumption~\ref{as:sto-grad}. Moreover, the value of $r$ boils down to $r = \max\left(1-\gamma\mu,\left(1-\rho+\frac{\beta}{\kappa}\right)\right)$. Consider the applications of \algname{SGD} covered by Assumption~\ref{as:sto-grad}, and mentioned
  in Sect.~\ref{sec:grad-estimators}. Then, as proved in~\cite{gorbunov2020unified}, the value $r = \max\left(1-\gamma\mu,\left(1-\rho+\frac{\beta}{\kappa}\right)\right)$ matches the best known convergence rates for these applications, with an exception for some coordinate descent algorithms. However, if $H = 0$ and $\mL =0$ but $\psi  \neq 0$, then the \algname{Stochastic PD3O} algorithm  boils down to \algname{Proximal SGD}, and $r$ boils down to $r = \max\left(1-\frac{\gamma\mu}{(1+\gamma\lambda)^2},\left(1-\rho+\frac{\beta}{\kappa}\right)\right)$, whereas the best known rates for \algname{Proximal SGD} under Assumption~\ref{as:sto-grad} is $\max\left(1-\gamma\mu,\left(1-\rho+\frac{\beta}{\kappa}\right)\right)$. Finally, if $g^{k+1} = \nabla f(x^k)$, the \algname{Stochastic PD3O} algorithm boils down to the \algname{PD3O} algorithm and Theorem~\ref{th:Hsmooth:PD3O} provides a convergence rate similar to Theorem 3 in~\cite{yan2018new}. In this case, by taking $\kappa = 1$, we obtain
\begin{equation*}
  r = \max\left(1-\gamma\frac{\mu_f + 2\mu_\psi }{(1+\gamma \lambda)^2},\frac{1}{1+2\tau\mu_{H^*}}\right),
\end{equation*}
whereas Theorem 3 in~\cite{yan2018new} provides the rate\footnote{The reader might not recognize the rate given in Theorem 3 of~\cite{yan2018new} because of some typos in its Equation (39).} 
\[
	\max\left(1-\gamma\frac{2(\mu_f+\mu_\psi ) - \gamma\alpha\mu_f}{(1+\gamma \lambda)^2},\frac{1}{1+2\tau\mu_{H^*}}\right).
\]
\end{remark}

\subsection{The \algname{Stochastic PDDY} algorithm}

\begin{theorem}[$M$ strongly monotone]
    \label{th:Hsmooth:PDDY}
    Suppose that Assumption~\ref{as:sto-grad} holds. Also, suppose that
     $H$ is $1/\mu_{H^*}$-smooth, $f$ is $\mu_f$-strongly convex and $\psi $ is $\mu_\psi $-strongly convex, where $\mu_\psi  >0$ and $\mu_{H^*} >0$.
    For every $\kappa > \beta/\rho$ and every $\gamma, \tau >0$ such that $\gamma \leq 1/(\alpha+\kappa\delta)$, $\gamma\tau\|\mL\|^2 < 1$ and $\gamma^2 \leq \frac{\mu_{H^*}}{\|\mL\|^2 \mu_\psi }$, define $\eta \eqdef 2\left(\mu_{H^*} -\gamma^2\|\mL\|^2\mu_\psi \right) \geq 0$,
    \begin{equation}
        V^k \eqdef (1+\gamma\mu_\psi )\|p^{k} - p^\star\|^2 + (1+\tau\eta)\|y^{k} - y^\star\|_{\gamma,\tau}^2 + \kappa \gamma^2\sigma_{k}^2,
    \end{equation}
     and
     \begin{equation}
        r \eqdef \max\left(\frac{1}{1+\gamma\mu_\psi },1-\rho+\frac{\beta}{\kappa},\frac{1}{1+\tau\eta}\right)
    \end{equation}
    Then,
    \begin{equation}
     \bE \left[ V^{k} \right]
     \leq r^k V^0.
    \end{equation}
\end{theorem}
Theorem~\ref{th:Hsmooth:PDDY} provides a linear convergence result for the \algname{Stochastic PDDY} algorithm, without assuming $\psi $ smooth.

\section{Proofs Related to the \algname{Stochastic PDDY} Algorithm}

Recall that the \algname{PDDY} algorithm is equivalent to \algname{DYS}$(\mP^{-1}B,\mP^{-1}A,\mP^{-1}C)$. We denote by $v^k = (p^k,q^k)$, $z^k = (x^k,y^k)$, $u^k = (s^k,h^k)$ the iterates of \algname{DYS}$(\mP^{-1}B,\mP^{-1}A,\mP^{-1}C)$, where $p^k, x^k, s^k \in \cX$ and $q^k, y^k, h^k \in \cY$.

Using~\eqref{algorithm_reslv}, the step
\begin{equation*}
    z^k = J_{\gamma \mP^{-1}A}(v^k),
\end{equation*}
is equivalent to
\begin{equation*}
\left\lfloor\begin{array}{l}
    x^k = p^k -\gamma \mL^* y^{k}\notag\\
    y^{k} = \prox_{\tau  H^*}\big((\mI-\tau\gamma \mL\mL^*)q^{k} +\tau \mL p^k \big).
    \end{array}\right.
\end{equation*}
Then, the step
\begin{equation*}
    u^{k+1} = J_{\gamma \mP^{-1}B}\big(2z^k - v^k - \gamma \mP^{-1}C(z^k)\big)
\end{equation*}
is equivalent to
\begin{equation*}\left\lfloor\begin{array}{l}
    s^{k+1} =\prox_{\gamma \psi }\big(2x^k-p^k-\gamma \nabla f(x^k)\big)\\
    h^{k+1} = 2y^k - q^k.
\end{array}\right.\end{equation*}
Finally, the step
\begin{equation*}
    v^{k+1} = v^k + u^{k+1} - z^k
\end{equation*}
is equivalent to
\begin{equation*}\left\lfloor\begin{array}{l}
    p^{k+1} = p^k + s^{k+1} - x^k\\
    q^{k+1} = q^k + h^{k+1} - y^k.
\end{array}\right.\end{equation*}
Similarly, the fixed points $v^\star = (p^\star,q^\star), z^\star = (x^\star,y^\star), u^\star = (s^\star,h^\star)$ of \algname{DYS}$(\mP^{-1}B,\mP^{-1}A,\mP^{-1}C)$ satisfy
\begin{equation*}\left\{\begin{array}{l}
    x^\star = p^\star -\gamma \mL^* y^{\star}\\
    y^{\star} = \prox_{\tau  H^*}\big((\mI-\tau\gamma \mL\mL^*)q^{\star} +\tau \mL p^\star \big)\\
    s^{\star} =\prox_{\gamma \psi }\big(2x^\star-p^\star-\gamma \nabla f(x^\star)\big)\\
    h^{\star} = 2y^\star - q^\star\\
    p^\star = p^\star + s^\star - x^\star\\
    q^\star = q^\star + h^\star - y^\star,
\end{array}\right.\end{equation*}
and the iterates of the \algname{Stochastic PDDY} algorithm satisfy
\begin{equation*}\left\lfloor\begin{array}{l}
    x^k = p^k -\gamma \mL^* y^{k}\\
    y^{k} = \prox_{\tau  H^*}\big((\mI-\tau\gamma \mL\mL^*)q^{k} +\tau \mL p^k \big)\\
    s^{k+1} =\prox_{\gamma \psi }\big(2x^k-p^k-\gamma g^{k+1}\big)\\
    h^{k+1} = 2y^k - q^k\\
    p^{k+1} = p^k + s^{k+1} - x^k\\
    q^{k+1} = q^k + h^{k+1} - y^k.
\end{array}\right.\end{equation*}

\begin{lemma}
\label{lem:PDDY}
Suppose that $(g^k)_k$ satisfies Assumption~\ref{as:sto-grad}. Then, the iterates of the \algname{Stochastic PDDY} algorithm satisfy
\begin{align*}
    \bE_k \left[\|v^{k+1} - v^\star\|_\mP^2 \right] + \kappa \gamma^2\bE_k\left[\sigma_{k+1}^2\right] 
    \leq{}& \|v^k - v^\star\|_\mP^2 + \kappa \gamma^2\left(1-\rho+\frac{\beta}{\kappa }\right)\sigma_k^2\notag\\
    &-2\gamma(1-\gamma(\alpha  +\kappa \delta)) D_f(x^{k},x^\star)\notag\\
    &-2\gamma\ps{\partial H^*(y^{k}) - \partial H^*(y^\star),y^{k} - y^\star}\\
    &-2\gamma\bE_k\left[\ps{\partial \psi(s^{k+1}) - \partial \psi(s^\star),s^{k+1} - s^\star} \right].\notag
\end{align*}
\end{lemma}

\begin{proof}

Applying Lemma~\ref{lem:funda-DYS} for $\algname{DYS}(\mP^{-1}B,\mP^{-1}A,\mP^{-1}C)$ using the norm induced by $\mP$, we have
\begin{align*}
    \|v^{k+1} - v^\star\|_\mP^2 ={}& \|v^k - v^\star\|_\mP^2 \\
    &-2\gamma\ps{\mP^{-1}A(z^{k}) - \mP^{-1}A(z^\star),z^{k} - z^\star}_\mP\\
    &-2\gamma\ps{\mP^{-1}C(z^{k}) - \mP^{-1}C(z^\star),z^{k} - z^\star}_\mP\\
    &-2\gamma\ps{\mP^{-1}B(u^{k+1}) - \mP^{-1}B(u^\star),u^{k+1} - u^\star}_\mP\\
    &+\gamma^2\|\mP^{-1}C(z^{k}) - \mP^{-1}C(z^\star)\|_\mP^2\\
    &-\gamma^2\|\mP^{-1}B(u^{k+1})+\mP^{-1}A(z^{k}) - \left(\mP^{-1}B(u^{\star})+\mP^{-1}A(z^{\star})\right)\|_\mP^2\\
    ={}& \|v^k - v^\star\|_\mP^2 \\
    &-2\gamma\ps{A(z^{k}) - A(z^\star),z^{k} - z^\star}\\
    &-2\gamma\ps{C(z^{k}) - C(z^\star),z^{k} - z^\star}\\
    &-2\gamma\ps{B(u^{k+1}) - B(u^\star),u^{k+1} - u^\star}\\
    &+\gamma^2\|\mP^{-1}C(z^{k}) - \mP^{-1}C(z^\star)\|_\mP^2\\
    &-\gamma^2\|\mP^{-1}B(u^{k+1})+\mP^{-1}A(z^{k}) - \left(\mP^{-1}B(u^{\star})+\mP^{-1}A(z^{\star})\right)\|_\mP^2.
\end{align*}
Using
\begin{align*}
    A(z^{k}) &=  \begin{bmatrix} &   \mL^* y^k \\ -\mL  x^k\!\!\!\!\!\!& {}+\partial H^*(y^k)\end{bmatrix},&
    B(u^{k+1}) &=  \begin{bmatrix} \partial \psi(s^{k+1}) \\ 0 \end{bmatrix}, &
    C(z^k) &= \begin{bmatrix} g^{k+1}\\ 0\end{bmatrix},
\end{align*}
and
\begin{align*}
    A(z^{\star}) &=  \begin{bmatrix} &   \mL^* y^\star \\ -\mL  x^\star\!\!\!\!\!\!& {}+\partial H^*(y^\star)\end{bmatrix}, &
    B(u^{\star}) &=  \begin{bmatrix} \partial \psi(s^{\star}) \\ 0 \end{bmatrix}, &
    C(z^\star) &= \begin{bmatrix} \nabla f(x^\star)\\ 0\end{bmatrix},\notag
\end{align*}
we have,
\begin{align*}
    \|v^{k+1} - v^\star\|_\mP^2 \leq{}& \|v^k - v^\star\|_\mP^2\notag\\
    &-2\gamma\ps{\partial H^*(y^{k}) - \partial H^*(y^\star),y^{k} - y^\star}\notag\\
    &-2\gamma\ps{g^{k+1} - \nabla f(x^\star),x^{k} - x^\star}\\
    &-2\gamma\ps{\partial \psi(s^{k+1}) - \partial \psi(s^\star),s^{k+1} - s^\star}\notag\\
    &+\gamma^2\|g^{k+1} - \nabla f(x^\star)\|^2.\notag
\end{align*}
Applying the conditional expectation w.r.t.\ $\cF_k$ and using Assumption~\ref{as:sto-grad},
\begin{align*}
    \bE_k \left[\|v^{k+1} - v^\star\|_\mP^2\right]
    \leq{}& \|v^k - v^\star\|_\mP^2\notag\\
    &-2\gamma\ps{\partial H^*(y^{k}) - \partial H^*(y^\star),y^{k} - y^\star}\notag\\
    &-2\gamma\ps{\nabla f(x^{k}) - \nabla f(x^\star),x^{k} - x^\star}\\
    &-2\gamma\bE_k\left[ \ps{\partial \psi(s^{k+1}) - \partial \psi(s^\star),s^{k+1} - s^\star} \right] \notag\\
    &+\gamma^2 \left(2\alpha   D_f(x^{k},x^\star) + \beta\sigma_k^2\right).\notag
\end{align*}

Using the convexity of $f$,
\begin{align*}
    \bE_k\left[ \|v^{k+1} - v^\star\|_\mP^2 \right]
    \leq{}& \|v^k - v^\star\|_\mP^2\notag\\
    &-2\gamma\ps{\partial H^*(y^{k}) - \partial H^*(y^\star),y^{k} - y^\star}\notag\\
    &-2\gamma\bE_k \left[ \ps{\partial \psi(s^{k+1}) - \partial \psi(s^\star),s^{k+1} - s^\star} \right] \\
    &-2\gamma D_f(x^{k},x^\star)\notag\\
    &+\gamma^2 \left(2\alpha   D_f(x^{k},x^\star) + \beta\sigma_k^2\right).\notag
\end{align*}
Using Assumption~\ref{as:sto-grad},
\begin{align*}
    \bE_k \left[\|v^{k+1} - v^\star\|_\mP^2  \right] + \kappa \gamma^2\bE_k\left[\sigma_{k+1}^2 \right]
    \leq{}& \|v^k - v^\star\|_\mP^2 + \kappa \gamma^2\left(1-\rho+\frac{\beta}{\kappa }\right)\sigma_k^2\notag\\
    &-2\gamma(1-\gamma(\alpha  +\kappa \delta)) D_f(x^{k},x^\star)\notag\\
    &-2\gamma\ps{\partial H^*(y^{k}) - \partial H^*(y^\star),y^{k} - y^\star}\\
    &-2\gamma\bE_k \left[ \ps{\partial \psi(s^{k+1}) - \partial \psi(s^\star),s^{k+1} - s^\star} \right] .\notag
\end{align*}
\end{proof}

\subsection{Proof of Theorem~\ref{th:cvx:PDDY}}
Using Lemma~\ref{lem:PDDY} and the convexity of $f,\psi ,H^*$,
\begin{align*}
    &\bE_k \left[\|v^{k+1} - v^\star\|_\mP^2 \right]+ \kappa \gamma^2\bE_k \left[\sigma_{k+1}^2 \right] \\
    &\leq \|v^k - v^\star\|_\mP^2 + \kappa \gamma^2\left(1-\rho+\frac{\beta}{\kappa }\right)\sigma_k^2\\
    &-2\gamma\big(1-\gamma(\alpha  +\kappa \delta)\big) \left(D_f(x^{k},x^\star) + D_{H^*}(y^{k},y^\star)+\bE_k \left[D_\psi(s^{k+1},s^\star) \right] \right).
\end{align*}

Since $1 - \rho + \beta/\kappa  = 1$, $\gamma \leq 1/2(\alpha  +\kappa \delta)$. Set 
\begin{equation*}
V^k = \|v^{k} - v^\star\|_\mP^2 + \kappa \gamma^2\sigma_{k}^2.
\end{equation*}
Then
\begin{equation*}
    \bE_k \left[ V^{k+1} \right]
    \leq V^k -\gamma \bE_k\left[D_f(x^{k},x^\star) + D_{H^*}(y^{k},y^\star)+D_\psi(s^{k+1},s^\star)\right].
\end{equation*}
Taking the expectation,
\begin{equation*}
    \gamma\bE \left[D_f(x^{k},x^\star) + D_{H^*}(y^{k},y^\star)+D_\psi(s^{k+1},s^\star)\right] \leq \bE \left[ V^k \right] - \bE \left[ V^{k+1} \right].
\end{equation*}
Iterating and using the nonnegativity of $V^k$,
\begin{equation}
    \gamma \sum_{j = 0}^{k-1} \bE \left[D_f(x^{k},x^\star) + D_{H^*}(y^{k},y^\star)+D_\psi(s^{k+1},s^\star)\right] 
    \leq \bE \left[ V^0 \right].
\end{equation}
We conclude using the convexity of the Bregman divergence in its first variable.

\subsection{Proof of Theorem~\ref{th:Hsmooth:PDDY}}

We first use Lemma~\ref{lem:PDDY} along with the strong convexity of $\psi ,H^*$. Note that $y^{k} = q^{k+1}$. We have
\begin{align*}
    \bE_k \left[ \|v^{k+1} - v^\star\|_\mP^2 \right] + \kappa \gamma^2\bE_k\left[ \sigma_{k+1}^2 \right] 
    \leq{}& \|v^k - v^\star\|_\mP^2 + \kappa \gamma^2\left(1-\rho+\frac{\beta}{\kappa }\right)\sigma_k^2\notag\\
    &-2\gamma\mu_{H^*}\bE_k \left[ \|q^{k+1} - q^\star\|^2 \right] -2\gamma\mu_{\psi }\bE_k\left[ \|s^{k+1} - s^\star\|^2 \right].
\end{align*}
Note that $s^{k+1} = p^{k+1} - \gamma \mL^* y^k$. Therefore, $s^{k+1} - s^\star = (p^{k+1} - p^\star) - \gamma \mL^*  (y^k - y^\star)$. Using Young's inequality $-\|a+b\|^2 \leq -\frac{1}{2}\|a\|^2 + \|b\|^2$, we have
\begin{equation*}
	-\bE_k \left[ \|s^{k+1} - s^\star\|^2 \right] 
	\leq -\frac{1}{2}\bE_k \left[ \|p^{k+1} - p^\star\|^2 \right] + \gamma^2\|\mL\|^2 \bE_k \left[ \|q^{k+1} - q^\star\|^2 \right].
\end{equation*}
Hence,
\begin{align*}
    \bE_k\left[ \|v^{k+1} - v^\star\|_\mP^2 \right] + \kappa \gamma^2\bE_k \left[ \sigma_{k+1}^2 \right]
    \leq{}& \|v^k - v^\star\|_\mP^2 + \kappa \gamma^2\left(1-\rho+\frac{\beta}{\kappa }\right)\sigma_k^2\\
    &-2\gamma\left(\mu_{H^*} -\gamma^2\|\mL\|^2\mu_\psi \right)\bE_k \left[ \|q^{k+1} - q^\star\|^2\right]\\
    & - \gamma\mu_{\psi }\bE_k \left[\|p^{k+1} - p^\star\|^2 \right] . \notag
\end{align*}
Set $\eta = 2\left(\mu_{H^*} -\gamma^2\|\mL\|^2\mu_\psi \right) \geq 0$. Then
\begin{align*}
    &(1+\gamma\mu_\psi )\bE_k \left[ \|p^{k+1} - p^\star\|^2 \right] + (1+\tau\eta)\bE_k \left[\|q^{k+1} - q^\star\|_{\gamma,\tau}^2 \right] + \kappa \gamma^2\bE_k \left[\sigma_{k+1}^2 \right] \notag\\
     \leq{}& \|v^k - v^\star\|_\mP^2 + \kappa \gamma^2\left(1-\rho+\frac{\beta}{\kappa }\right)\sigma_k^2.
\end{align*}
Set
\begin{equation*}
    V^k = (1+\gamma\mu_\psi )\|p^{k} - p^\star\|^2 + (1+\tau\eta)\|q^{k} - q^\star\|_{\gamma,\tau}^2 + \kappa \gamma^2\sigma_{k}^2
\end{equation*}
and
\begin{equation*}
    r = \max\left(\frac{1}{1+\gamma\mu_\psi },1-\rho+\frac{\beta}{\kappa},\frac{1}{1+\tau\eta}\right).
\end{equation*}
Then
\begin{equation*}
    \bE_k \left[ V^{k+1}\right] 
     \leq r V^k.
\end{equation*}

\section{Proofs Related to the \algname{Stochastic PD3O} algorithm}

Recall that the \algname{PD3O} algorithm is equivalent to \algname{DYS}$(\mP^{-1}A,\mP^{-1}B,\mP^{-1}C)$. We denote by $v^k = (p^k,q^k)$, $z^k = (x^k,y^k)$, $u^k = (s^k,h^k)$ the variables in \algname{DYS}$(\mP^{-1}A,\mP^{-1}B,\mP^{-1}C)$, with $p^k, x^k, s^k \in \cX$ and $q^k, y^k, h^k \in \cY$.

Then, the step
\begin{equation*}
    z^k = J_{\gamma \mP^{-1}B}(v^k),
\end{equation*}
is equivalent to
\begin{equation*}\left\lfloor\begin{array}{l}
    x^k = \prox_{\gamma \psi }(p^k)\\
    y^k = q^k.
    \end{array}\right.
\end{equation*}
Using~\eqref{algorithm_reslv}, the step
\begin{equation*}
    u^{k+1} = J_{\gamma \mP^{-1}A}(2z^k - v^k - \gamma \mP^{-1}C(z^k)),
\end{equation*}
is equivalent to
\begin{equation*}\left\lfloor\begin{array}{l}
    s^{k+1} = (2 x^k - p^k - \gamma \nabla f(x^k)) - \gamma \mL^* h^{k+1}\\
    h^{k+1} = \prox_{\tau H^*}\left((\mI-\gamma\tau \mL \mL^*)(2 y^k - q^k) + \tau \mL (2 x^k - p^k - \nabla f(x^k))\right).
    \end{array}\right.
\end{equation*}
Finally, the step
\begin{equation*}
    v^{k+1} = v^k + u^{k+1} - z^k,
\end{equation*}
is equivalent to
\begin{equation*}\left\lfloor\begin{array}{l}
    p^{k+1} = p^k + s^{k+1} - x^k\\
    q^{k+1} = q^k + h^{k+1} - y^k.
    \end{array}\right.
\end{equation*}
Similarly, the fixed points $v^\star = (p^\star,q^\star), z^\star = (x^\star,y^\star), u^\star = (s^\star,h^\star)$ of \algname{DYS}$(\mP^{-1}A,\mP^{-1}B,\mP^{-1}C)$ satisfy
\begin{equation*}\left\lfloor\begin{array}{l}
    x^\star = \prox_{\gamma \psi }(p^\star)\\
    y^\star = q^\star\\
    s^\star = (2 x^\star - p^\star - \gamma \nabla f(x^\star)) - \gamma \mL^* h^\star\\
    h^\star = \prox_{\tau H^*}\left((\mI-\gamma\tau \mL \mL^*)(2 y^\star - q^\star) + \tau \mL (2 x^\star - p^\star - \nabla f(x^\star))\right)\\
    p^\star = p^\star + s^\star - x^\star\\
    q^\star = q^\star + h^\star - y^\star.
    \end{array}\right.
\end{equation*}
and the iterates of the \algname{Stochastic PD3O} algorithm satisfy
\begin{equation*}\left\lfloor\begin{array}{l}
    x^k = \prox_{\gamma \psi }(p^k)\\
    y^k = q^k\\
    s^{k+1} = (2 x^k - p^k - \gamma g^{k+1}) - \gamma \mL^* h^{k+1}\\
    h^{k+1} = \prox_{\tau H^*}\left((\mI-\gamma\tau \mL \mL^*)(2 y^k - q^k) + \tau \mL (2 x^k - p^k - g^{k+1})\right)\\
    p^{k+1} = p^k + s^{k+1} - x^k\\
    q^{k+1} = q^k + h^{k+1} - y^k.
    \end{array}\right.
\end{equation*}

\begin{lemma}
\label{lem:ineq:PD3O}
Assume that $f$ is $\mu_f$-strongly convex, for some $\mu_f \geq 0$, and that $(g^k)_k$ satisfies Assumption~\ref{as:sto-grad}. Then, the iterates of the \algname{Stochastic PD3O} algorithm satisfy
\begin{align}
    \bE_k \left[ \|v^{k+1} - v^\star\|_\mP^2  \right] + \kappa \gamma^2\bE_k\left[ \sigma_{k+1}^2 \right] 
    \leq{}& \|v^k - v^\star\|_\mP^2 + \kappa \gamma^2\left(1-\rho+\frac{\beta}{\kappa }\right)\sigma_k^2\notag \\
    &-2\gamma(1-\gamma(\alpha  +\kappa \delta)) D_f(x^{k},x^\star)-\gamma\mu_f\|x^{k}-x^\star\|^2\notag\\
    &-2\gamma\ps{\partial \psi(x^{k}) - \partial \psi(x^\star),x^{k} - x^\star}\\
    &-2\gamma\bE_k \left[\ps{\partial H^*(h^{k+1}) - \partial H^*(h^\star),h^{k+1} - h^\star} \right] \notag\\
    &-\gamma^2\bE_k \Bigl[\big\|\mP^{-1}A(u^{k+1})+\mP^{-1}B(z^{k})\notag\\
    & \quad\quad\quad\quad- \left(\mP^{-1}A(u^{\star})+\mP^{-1}B(z^{\star})\right)\big\|_\mP^2\Bigr].\notag
\end{align}
\end{lemma}

\begin{proof}
Applying Lemma~\ref{lem:funda-DYS} for $\algname{DYS}(\mP^{-1}A,\mP^{-1}B,\mP^{-1}C)$ using the norm induced by $\mP$ we have
\begin{align*}
    \|v^{k+1} - v^\star\|_\mP^2 ={}& \|v^k - v^\star\|_\mP^2 \\
    &-2\gamma\ps{\mP^{-1}B(z^{k}) - \mP^{-1}B(z^\star),z^{k} - z^\star}_\mP\\
    &-2\gamma\ps{\mP^{-1}C(z^{k}) - \mP^{-1}C(z^\star),z^{k} - z^\star}_\mP\\
    &-2\gamma\ps{\mP^{-1}A(u^{k+1}) - \mP^{-1}A(u^\star),u^{k+1} - u^\star}_\mP\\
    &+\gamma^2\|\mP^{-1}C(z^{k}) - \mP^{-1}C(z^\star)\|_\mP^2\\
    &-\gamma^2\|\mP^{-1}A(u^{k+1})+\mP^{-1}B(z^{k}) - \left(\mP^{-1}A(u^{\star})+\mP^{-1}B(z^{\star})\right)\|_\mP^2\\
    {}={}& \|v^k - v^\star\|_\mP^2 \\
    &-2\gamma\ps{B(z^{k}) - B(z^\star),z^{k} - z^\star}\\
    &-2\gamma\ps{C(z^{k}) - C(z^\star),z^{k} - z^\star}\\
    &-2\gamma\ps{A(u^{k+1}) - A(u^\star),u^{k+1} - u^\star}\\
    &+\gamma^2\|\mP^{-1}C(z^{k}) - \mP^{-1}C(z^\star)\|_\mP^2\\
    &-\gamma^2\|\mP^{-1}A(u^{k+1})+\mP^{-1}B(z^{k}) - \left(\mP^{-1}A(u^{\star})+\mP^{-1}B(z^{\star})\right)\|_\mP^2.
\end{align*}
Using
\begin{align*}
    A(u^{k+1}) &=  \begin{bmatrix} &   \mL^* h^{k+1} \\ -\mL s^{k+1}\!\!\!\!\!\!& {}+\partial H^*(h^{k+1})\end{bmatrix}, &
    B(z^k) &=  \begin{bmatrix} \partial \psi(x^{k}) \\ 0 \end{bmatrix}, &
    C(z^k) &= \begin{bmatrix} g^{k+1}\\ 0\end{bmatrix},
\end{align*}
and
\begin{align*}
	A(u^{\star}) &=  \begin{bmatrix} &   \mL^* h^{\star} \\ -\mL  s^{\star}\!\!\!\!\!\!& {}+\partial H^*(h^{\star})\end{bmatrix},&
	B(z^\star) &=  \begin{bmatrix} \partial \psi(x^{\star}) \\ 0 \end{bmatrix},&
	C(z^\star) &= \begin{bmatrix} \nabla f(x^\star)\\ 0\end{bmatrix},
\end{align*}
we derive
\begin{align*}
    \|v^{k+1} - v^\star\|_\mP^2 ={}& \|v^k - v^\star\|_\mP^2 \notag\\
    &-2\gamma\ps{\partial \psi(x^{k}) - \partial \psi(x^\star),x^{k} - x^\star}\notag\\
    &-2\gamma\ps{g^{k+1} - \nabla f(x^\star),x^{k} - x^\star}\notag\\
    &-2\gamma\ps{\partial H^*(h^{k+1}) - \partial H^*(h^\star),h^{k+1} - h^\star}\\
    &+\gamma^2\|g^{k+1} - \nabla f(x^\star)\|^2\notag\\
    &-\gamma^2\|\mP^{-1}A(u^{k+1})+\mP^{-1}B(z^{k}) - \left(\mP^{-1}A(u^{\star})+\mP^{-1}B(z^{\star})\right)\|_\mP^2.\notag
\end{align*}
Taking conditional expectation w.r.t.\ $\cF_k$ and using Assumption~\ref{as:sto-grad},
\begin{align*}
    \bE_k \left[ \|v^{k+1} - v^\star\|_\mP^2 \right]
    \leq{}& \|v^k - v^\star\|_\mP^2 \notag\\
    &-2\gamma\ps{\partial \psi(x^{k}) - \partial \psi(x^\star),x^{k} - x^\star}\notag\\
    &-2\gamma\ps{\nabla f(x^{k}) - \nabla f(x^\star),x^{k} - x^\star}\notag\\
    &-2\gamma\bE_k \left[\ps{\partial H^*(h^{k+1}) - \partial H^*(h^\star),h^{k+1} - h^\star} \right]\\
    &+\gamma^2 \left(2\alpha   D_f(x^{k},x^\star) + \beta\sigma_k^2\right)\notag\\
    &-\gamma^2\bE_k\left[ \|\mP^{-1}A(u^{k+1})+\mP^{-1}B(z^{k}) - \left(\mP^{-1}A(u^{\star})+\mP^{-1}B(z^{\star})\right)\|_\mP^2 \right] .\notag
\end{align*}
Using strong convexity of $f$,
\begin{align*}
    \bE_k \left[\|v^{k+1} - v^\star\|_\mP^2 \right]
    \leq{}& \|v^k - v^\star\|_\mP^2 \notag\\
    &-\gamma\mu_f\|x^{k}-x^\star\|^2\notag\\
    &-2\gamma D_f(x^{k},x^\star)\notag\\
    &+\gamma^2 \left(2\alpha  D_f(x^{k},x^\star) + \beta\sigma_k^2\right)\\
    &-2\gamma\ps{\partial \psi(x^{k}) - \partial \psi(x^\star),x^{k} - x^\star}\notag\\
    &-2\gamma\bE_k \left[ \ps{\partial H^*(h^{k+1}) - \partial H^*(h^\star),h^{k+1} - h^\star} \right] \notag\\
    &-\gamma^2\bE_k \left[\|\mP^{-1}A(u^{k+1})+\mP^{-1}B(z^{k}) - \left(\mP^{-1}A(u^{\star})+\mP^{-1}B(z^{\star})\right)\|_\mP^2 \right].\notag
\end{align*}
Using Assumption~\ref{as:sto-grad},
\begin{align*}
    \bE_k \left[\|v^{k+1} - v^\star\|_\mP^2 \right] + \kappa \gamma^2\bE_k\left[ \sigma_{k+1}^2 \right] 
    \leq{}& \|v^k - v^\star\|_\mP^2 + \kappa \gamma^2\left(1-\rho+\frac{\beta}{\kappa }\right)\sigma_k^2 \notag\\
    &-\gamma\mu_f\|x^{k}-x^\star\|^2\notag\\
    &-2\gamma(1-\gamma(\alpha  +\kappa \delta)) D_f(x^{k},x^\star)\notag\\
    &-2\gamma\ps{\partial \psi(x^{k}) - \partial \psi(x^\star),x^{k} - x^\star}\\
    &-2\gamma\bE_k \left[ \ps{\partial H^*(h^{k+1}) - \partial H^*(h^\star),h^{k+1} - h^\star} \right] \notag\\
    &-\gamma^2\bE_k \Bigl[\big\|\mP^{-1}A(u^{k+1})+\mP^{-1}B(z^{k})\\
    &\quad\quad\quad\quad - \left(\mP^{-1}A(u^{\star})+\mP^{-1}B(z^{\star})\right)\big\|_\mP^2 \Bigr].\notag
\end{align*}
\end{proof}

\subsection{Proof of Theorem~\ref{th:cvx:PD3O}}

Using Lemma~\ref{lem:ineq:PD3O}, convexity of $f,\psi ,H^*$, and Lemma~\ref{lem:duality-gap},
\begin{align*}
    \bE_k \left[\|v^{k+1} - v^\star\|_\mP^2\right] + \kappa \gamma^2\bE_k\left[\sigma_{k+1}^2\right] 
    &\leq \|v^k - v^\star\|_\mP^2 + \kappa \gamma^2\left(1-\rho+\frac{\beta}{\kappa }\right)\sigma_k^2\notag\\
    &\qquad -2\gamma(1-\gamma(\alpha  +\kappa \delta))\bE_k\left[\cL(x^k,h^{\star}) - \cL(x^\star,h^{k+1})\right].
\end{align*}

Recall that $1 - \rho + \beta/\kappa  = 1$, $\gamma \leq 1/2(\alpha  +\kappa \delta)$. Set 
\begin{equation*}
V^k = \|v^{k} - v^\star\|_\mP^2 + \kappa \gamma^2\sigma_{k}^2.
\end{equation*}
Then,
\begin{equation*}
    \bE_k \left[V^{k+1}\right]
    \leq \left[V^k\right] -\gamma \bE_k\left[\cL(x^{k},h^\star) - \cL(x^\star,h^{k+1})\right].
\end{equation*}
Taking the expectation,
\begin{equation*}
    \gamma\bE\left[\cL(x^{k},h^\star) - \cL(x^\star,h^{k+1})\right] 
    \leq \bE \left[V^k\right] - \bE \left[V^{k+1} \right].
\end{equation*}
Iterating and using the nonnegativity of $V^k$,
\begin{equation*}
    \gamma \sum_{j = 0}^{k-1} \bE\left[\cL(x^{j},h^\star) - \cL(x^\star,h^{j+1})\right] \leq \bE \left[V^0 \right].
\end{equation*}
We conclude using the convex-concavity of $\mL$.

\subsection{Proof of Theorem~\ref{th:Hsmooth:PD3O}}

We first use Lemma~\ref{lem:ineq:PD3O} along with the strong convexity of $\psi ,H^*$. Note that $y^{k} = q^k$ and therefore $q^{k+1} = q^k + h^{k+1} - q^{k} = h^{k+1}$. We have
\begin{align*}
    &\bE_k \left[ \|p^{k+1} - p^\star\|^2 \right] + \bE_k \left[\|q^{k+1} - q^\star\|_{\gamma,\tau}^2 \right] + 2\gamma\mu_{H^*}\bE_k \left[ \|q^{k+1} - q^\star\|^2 \right] + \kappa \gamma^2\bE_k \left[\sigma_{k+1}^2 \right] \\
    \leq{}& \|p^{k} - p^\star\|^2 + \|q^{k} - q^\star\|_{\gamma,\tau}^2 -\gamma\mu\|x^{k}-x^\star\|^2 \\
    &+ \kappa \gamma^2\left(1-\rho+\frac{\beta}{\kappa }\right)\sigma_k^2-2\gamma(1-\gamma(\alpha  +\kappa \delta)) D_f(x^{k},x^\star)
\end{align*}
Noting that for every $q \in \mathcal{Y}$, $\|q\|_{\gamma,\tau}^2 = \frac{\gamma}{\tau}\|q\|^2 - \gamma^2\|\mL ^* q\|^2 \leq \frac{\gamma}{\tau}\|q\|^2$, and taking $\gamma \leq 1/(\alpha  +\kappa \delta)$,
\begin{align}
    &\bE_k \left[ \|p^{k+1} - p^\star\|^2 \right] + \left(1+2\tau\mu_{H^*}\right)\bE_k \left[\|q^{k+1} - q^\star\|_{\gamma,\tau}^2 \right] + \kappa \gamma^2\bE_k \left[\sigma_{k+1}^2 \right] \notag\\
    \leq{}& \|p^{k} - p^\star\|^2 + \|q^{k} - q^\star\|_{\gamma,\tau}^2 -\gamma\mu\|x^{k}-x^\star\|^2 + \kappa \gamma^2\left(1-\rho+\frac{\beta}{\kappa }\right)\sigma_k^2.\notag
\end{align}
Finally, since $\psi $ is $\lambda$-smooth, $\|p^k - p^\star\|^2 \leq (1+2 \gamma \lambda + \gamma^2 \lambda^2)\|x^{k}-x^\star\|^2$. Indeed, in this case, applying Lemma~\ref{lem:funda-DYS} with $\tilde{A} =0$, $\tilde{C} = 0$ and $\tilde{B} = \nabla \psi $, we obtain that if $x^k = \prox_{\gamma \psi }(p^k)$ and $x^\star = \prox_{\gamma \psi }(p^\star)$, then
\begin{align*}
\|x^k - x^\star\|^2 ={}& \|p^k - p^\star\|^2-2\gamma\ps{\nabla \psi(x^k) - \nabla \psi(x^\star),x^k - x^\star} - \gamma^2\|\nabla \psi(x^k) - \nabla \psi(x^\star)\|^2\\
\geq{}& \|p^k - p^\star\|^2-2\gamma\lambda \|x^k - x^\star\|^2 - \gamma^2\lambda^2\|x^{k} - x^\star\|^2.\notag
\end{align*}
Hence,
\begin{align*}
    &\bE_k \left[ \|p^{k+1} - p^\star\|^2 \right] + \left(1+2\tau\mu_{H^*}\right)\bE_k \left[ \|q^{k+1} - q^\star\|_{\gamma,\tau}^2 \right] + \kappa \gamma^2\bE_k \left[ \sigma_{k+1}^2 \right] \\
    \leq{}& \|p^{k} - p^\star\|^2 + \|q^{k} - q^\star\|_{\gamma,\tau}^2 -\frac{\gamma\mu}{(1+\gamma \lambda)^2}\|p^{k}-p^\star\|^2+ \kappa \gamma^2\left(1-\rho+\frac{\beta}{\kappa }\right)\sigma_k^2.\notag
\end{align*}
Thus, set
\begin{equation*}
    V^k = \|p^{k} - p^\star\|^2+ \left(1+2\tau\mu_{H^*}\right)\|q^{k} - q^\star\|_{\gamma,\tau}^2 + \kappa \gamma^2 \sigma_{k}^2,
\end{equation*}
and
\begin{equation*}
    r = \max\left(1-\frac{\gamma\mu}{(1+\gamma \lambda)^2},\left(1-\rho+\frac{\beta}{\kappa }\right),\frac{1}{1+2\tau\mu_{H^*}}\right).
\end{equation*}
Then,
\begin{equation*}
    \bE_k \left[  V^{k+1} \right]    
     \leq r V^k.
\end{equation*}

\section{Convergence Results for \algname{LiCoSGD}}
\begin{algorithm}[t]
\caption{\algname{LiCoSGD} $\big($deterministic version: $g^{k+1}=\nabla f(x^k)\big)$}
 \begin{algorithmic}[1]
    \State \textbf{Input:} $x^0 \in \cX, y^0 \in \cY $, $\gamma>0$, $\tau>0$
 	\For{$k = 0,1,2,\dots$}
 	\State $w^k = x^k -\gamma g^{k+1}$
	\State $y^{k+1}=y^k + \tau \mL (w^k-\gamma \mL^* y^k)-\tau b$%
	\State $x^{k+1} =w^k - \gamma \mL^* y^{k+1}$
 	\EndFor
 \end{algorithmic}
\end{algorithm}
 
 We consider the problem
\begin{equation}
\min_{x\in\mathcal{X}} \,f(x)\quad \mbox{s.t.}\quad \mL x=b,\label{pbconstr}
\end{equation}
where $\mL\colon\mathcal{X}\rightarrow \mathcal{Y}$ is a linear operator,
 $\mathcal{X}$ and  $\mathcal{Y}$ are real Hilbert spaces, $f$ is a $\nu$-smooth convex function, for some $\nu>0$, 
 and  $b \in \Range{\mL}$. This is a particular case of Problem~\eqref{eq:original-pb}  with $\psi =0$ and $H:y\mapsto (0$ if $y=b$, $+\infty$ else$)$. The \algname{Stochastic PD3O} and \algname{PDDY} algorithms both revert to the same algorithm, shown above, which   
    we call \algname{Linearly Constrained Stochastic Gradient Descent} (\algname{LiCoSGD}).

    Theorem \ref{th:pddy-cv} becomes:
    
    \begin{theorem}[Convergence of \algname{LiCoSGD}, deterministic case]
    \label{th:lico-cv}
    Suppose that $\gamma\in (0,2/\nu)$ and that $\tau\gamma\|\mL\|^2\leq 1$. Then in \algname{LiCoSGD},    $(x^k)_{k}$ converges to some
    solution $x^\star$ to the problem \eqref{pbconstr} and $(y^k)_{k}$ converges to some dual solution $y^\star \in \argmin_y f^* (-\mL ^*y) + \langle y,b\rangle$.
\end{theorem}

Note that the case $\tau\gamma\|\mL\|^2= 1$ is not covered by Theorem \ref{th:pddy-cv} but follows from convergence of the \algname{PD3O} algorithm in that case, as proved in \cite{o2020equivalence}.

Theorem \ref{th:cvx:PD3O} becomes:
\begin{theorem}[Convergence of \algname{LiCoSGD}, stochastic case]
\label{th:cvx:lico}
Suppose that Assumption~\ref{as:sto-grad} holds.
Let $\kappa \eqdef \beta/\rho$, $\gamma, \tau >0$ be such that $\gamma \leq 1/{2(\alpha+\kappa\delta)}$ and $\gamma\tau\|\mL\|^2 < 1$.
Set $V^0 \eqdef \|v^{0} - v^\star\|_\mP^2 + \gamma^2 \kappa \sigma_{0}^2,$ where $v^0 = (w^0,y^0)$.
Then,
\begin{equation}
     \ec{f(\bar{x}^{k})-f(x^\star)+\langle \mL\bar{x}^{k}-b,y^\star\rangle }
     \leq \frac{V^0}{k \gamma},
\end{equation}
where $\bar{x}^{k} = \frac{1}{k} \sum_{j = 0}^{k-1} x^j$, 
$x^\star$ and $y^\star$ are some primal and dual solutions.
\end{theorem}
Note that the convex function $x\mapsto f(x)-f(x^\star)+\langle \mL x-b,y^\star\rangle$ is nonnegative and its minimum is 0, attained at $x^\star$; under mild conditions, this function takes value zero only if $f(x)=f(x^\star)$ and $\mL x=b$, so that $x$ is a solution.

Replacing the variable $y^k$ by the variable $a^k=\mL^*y^k$ in \algname{LiCoSGD} yields \algname{PriLiCoSGD}.  In the conditions of Theorem \ref{th:lico-cv}, 
$(a^k)_{k}$ converges to $a^\star=-\nabla f(x^\star)$.

\subsection{Proof of Theorem~\ref{th:LV0}}

We first derive the following lemma:

\begin{lemma}
    \label{lem:linalg}
    Let $x \in \Range{\mL^*}$, the range space of $\mL ^*$. There exists an unique $y \in \Range{\mL}$ such that $\mL ^*y = x$. Moreover, for every $y \in \Range{\mL}$, \begin{equation}\lambda_{\min}^+(\mL)\|y\|^2 \leq \|\mL ^* y\|^2,\end{equation}
        where $\lambda_{\min}^+(\mL)$ is the smallest positive eigenvalue of $\mL\mL^*$ (or $\mL^*\mL$).
\end{lemma}
\begin{proof}
    Using basic linear algebra, $\mL \mL^* x = 0$ implies $\mL ^*x \in \Range{\mL^*} \cap \ker(\mL)$ therefore $\mL ^* x = 0$. Hence, $\ker(\mL\mL^*) \subset \ker(\mL ^*)$ and therefore $\Range{\mL} \subset \Range{\mL\mL^*}$. Since $\mL\mL^*$ is real symmetric, for every $y \in \Range{\mL\mL^*}$, $\ps{y, \mL\mL^* y} \geq \lambda_{\min}^+(\mL)\|y\|^2$, where $\lambda_{\min}^+(\mL)$ is the smallest positive eigenvalue of $\mL\mL^*$. Therefore, for every $y \in \Range{\mL},$ $\|\mL ^* y\|^2 \geq \lambda_{\min}^+(\mL)\|y\|^2$. Moreover, $\mL ^* y = 0$ implies $y = 0$ on $\Range{\mL}$, therefore there is at most one solution $y$ in $\Range{\mL}$ to the equation $\mL ^* y = x$. The existence of a solution follows from $x \in \Range{\mL^*}$.
\end{proof}

Now, we prove Theorem~\ref{th:LV0}. 
First, we define $y^\star$. In the case $\psi  =0$ and $H = \ind_b$, Equation~\eqref{eq:saddle0} states that $\nabla f(x^\star) \in \Range{\mL^*}$. Using Lemma~\ref{lem:linalg}, there exists an unique $y^\star \in \Range{\mL}$ such that $\nabla f(x^\star) + \mL^* y^\star = 0$. Noting that $y^\star = h^\star = q^\star$ and applying Lemma~\ref{lem:ineq:PD3O} with $\gamma \leq (\alpha  +\kappa \delta)$,
\begin{align*}
    &\bE_k  \left[  \|p^{k+1} - p^\star\|^2 \right] + \bE_k \left[ \|q^{k+1} - q^\star\|_{\gamma,\tau}^2 \right] +  \kappa \gamma^2\bE_k \left[ \sigma_{k+1}^2 \right] \\
    &\leq{} \|p^{k} - p^\star\|^2 + \|q^{k} - q^\star\|_{\gamma,\tau}^2 -\gamma\mu_f\|x^{k}-x^\star\|^2 + \kappa \gamma^2\left(1-\rho+\frac{\beta}{\kappa }\right)\sigma_k^2\\
    &\qquad -\gamma^2\|\mP^{-1}A(u^{k+1}) - \mP^{-1}A(u^{\star})\|_\mP^2.\notag
\end{align*}
Since the component of $\mP^{-1}A(u^{k+1}) - \mP^{-1}A(u^{\star})$ in $\cX$ is $\mL^* h^{k+1} - \mL^* h^{\star}$, we have
\begin{align*}
    \bE_k  \left[ \|p^{k+1} - p^\star\|^2 \right] + \bE_k  \left[ \|q^{k+1} - q^\star\|_{\gamma,\tau}^2 \right] +  \kappa \gamma^2\bE_k \left[ \sigma_{k+1}^2 \right] 
    \leq{}& \|x^{k} - x^\star\|^2 + \|q^{k} - q^\star\|_{\gamma,\tau}^2 \\
    &-\gamma\mu_f\|p^{k}-p^\star\|^2\notag\\
    &+ \kappa \gamma^2\left(1-\rho+\frac{\beta}{\kappa }\right)\sigma_k^2\\
    &-\gamma^2\|\mL^* h^{k+1} - \mL^* h^{\star}\|^2.\notag
\end{align*}
Inspecting the iterations of the algorithm, one can see that $h^{0} \in \Range{\mL}$ implies $h^{k+1} \in \Range{\mL}$. Since $h^\star \in \Range{\mL}$, $h^{k+1} - h^\star \in \Range{\mL}$. Therefore, using Lemma~\ref{lem:linalg},
$\lambda_{\min}^+(\mL)\|h^{k+1} - h^{\star}\|^2 \leq \|\mL ^* h^{k+1} - \mL^* h^{\star}\|^2$. Since $q^{k+1} = h^{k+1} = y^{k+1}$ and $x^k = p^k$,
\begin{align*}
    &\bE_k  \left[ \|x^{k+1} - x^\star\|^2 \right] + (1+\gamma\tau\lambda_{\min}^+(\mL))\bE_k \left[ \|y^{k+1} - y^\star\|_{\gamma,\tau}^2 \right]  + \kappa \gamma^2\bE_k\left[\sigma_{k+1}^2\right]\notag\\
   \leq{}& (1-\gamma\mu_f)\|x^{k} - x^\star\|^2 + \|y^{k} - y^\star\|_{\gamma,\tau}^2+ \kappa \gamma^2\left(1-\rho+\frac{\beta}{\kappa}\right)\sigma_k^2.
    \end{align*}
Setting
\begin{equation*}
    V^k =\|x^{k} - x^\star\|^2+ \big(1+\tau\gamma\lambda_{\min}^+(\mL)\big)\|y^{k} - y^\star\|_{\gamma,\tau}^2 + \kappa \gamma^2 \sigma_{k}^2,
\end{equation*}
and
\begin{equation*}
    r = \max\left(1-\gamma\mu,1-\rho+\frac{\beta}{\kappa },\frac{1}{1+\tau\gamma\lambda_{\min}^+(\mL)}\right),
\end{equation*}
we have
\begin{equation*}
    \bE_k \left[V^{k+1} \right]
    \leq r V^k.
\end{equation*}

\section{Application of \algname{PriLiCoSGD} to Stochastic Decentralized Optimization}

Consider a connected undirected graph $G = (V,E)$, where $V = \{1,\ldots,N\}$ is the set of nodes and $E$ the set of edges. Consider a family $(f_i)_{i \in V}$ of $\mu$-strongly convex and $\nu$-smooth functions $f_i$, for some $\mu\geq 0$ and $\nu>0$. In this section, we consider solving the minimization problem
\begin{equation}
    \label{eq:decentralized}
    \min_{x \in \cX} \,\sum_{i \in V} f_i(x).
\end{equation}
Consider a gossip matrix of the graph $G$, i.e., a $N \times N$ symmetric positive semidefinite matrix $\widehat{\mW} = (\widehat{\mW}_{i,j})_{i,j \in V}$, such that $\ker(\widehat{\mW}) = \Span([1\ \cdots\ 1]^\mathrm{T})$ and $\widehat{\mW}_{i,j} \neq 0$ if and only if $i=j$ or $\{i,j\} \in E$ is an edge of the graph. $\widehat{\mW}$ can be the Laplacian matrix of the graph $G$, for instance.
Set $\mW = \widehat{\mW} \otimes \mI$, where $\otimes$ is the Kronecker product and $\mI$ the identity of $\cX$; $\mW$ is a positive linear operator over $\cX^V$ and $\mW(x_1,\ldots,x_N) = 0$ if and only if $x_1 = \ldots = x_N$. Therefore, Problem~\eqref{eq:decentralized} is equivalent to the lifted problem
\begin{equation}
    \label{eq:dec2}
    \min_{\tilde{x} \in \cX^V} f(\tilde{x}) \quad \text{such that} \quad \mW^{1/2} \tilde{x} = 0, 
\end{equation}
where for every $\tilde{x}=(x_1,\ldots,x_N) \in \cX^V$, $f(\tilde{x}) = \sum_{i=1}^N f_i(x_i)$. \algname{PriLiCoSGD} can be applied to Problem~\eqref{eq:dec2}. It involves
only one multiplication by $\mW$ per iteration and generates the sequence  $(\tilde{x}^k)_{k}$, where $\tilde{x}^k = (x_1^k,\ldots,x_N^k) \in \cX^V$. The update of each $x_i^k$ consists in local computations involving $f_i$ and communication steps involving $x_j^k$, where $j$ is a neighbor of $i$. We called the instance of \algname{PriLiCoSGD} applied to this setting the Decentralized Stochastic Optimization Algorithm (\algname{DESTROY}), shown above.
In details, at each iteration of \algname{DESTROY}, an estimate $\tilde{g}^{k+1}=(g_1^{k+1},\ldots,g_N^{k+1})$ of $\nabla f(\tilde{x}^k)=\big(\nabla f_1(x_1^k),\ldots,\nabla f_N(x_N^k)\big)$ is computed. Assuming that each sequence $(g_i^{k+1})_{k}$ satisfies Assumption~\ref{as:sto-grad} as an estimator of $\nabla f_i$, $(\tilde{g}^{k+1})_{k}$ satisfies Assumption~\ref{as:sto-grad} as an estimator of $\nabla f$. Moreover, the computation of $\tilde{g}^{k+1}$ boils down to the `local' computation of $g_i^{k+1}$ at each node $i \in V$, independently on each other. 
After this, 
decentralized communication in the network $G$ is performed, modeled by an application of $\mW$. 

For instance, the variance reduced estimator $g_i^k$ can be the \algname{L-SVRG} estimator when $f_i$ is a finite sum, or a compressed version of $\nabla f_i$. Such compressed gradients are suitable when communication (i.e., applications of $\mW$) is expensive, see Section~\ref{sec:grad-estimators}. 

\begin{algorithm}[t]
    \caption{
    \algname{DESTROY}
    $\big($deterministic version: $g_i^{k+1}=\nabla f_i(x^k)\big)$}
 \begin{algorithmic}[1]
    \State \textbf{Input:} $x_i^0 \in \mathcal{X}$ and $a_i^0 \in \mathcal{X}$,  for every $i \in V$, such that $\sum_{i\in V}a_i^0=0$, $\gamma>0$, $\tau>0$
 	\For{$k = 0,1,2,\dots$}
	\For{all $i \in V$ in parallel}
      \State $t_i^{k+1}=x_i^k - \gamma g_i^{k+1}$
	   \State $a_i^{k+1} = (1-\tau\gamma\widehat{\mW}_{i,i}) a_i^k +\tau\widehat{\mW}_{i,i} t_i^{k+1} + \tau \sum_{j\neq i:\{i,j\}\in V}\widehat{\mW}_{i,j} (t_j^{k+1}-\gamma a_j^k)$
	  \State $x_i^{k+1} = t_i^{k+1}-\gamma a_i^{k+1}$.
	  \EndFor
  		\EndFor
 \end{algorithmic}
 \end{algorithm}

As an application of the convergence results for \algname{LiCoSGD}, we obtain the following results for \algname{DESTROY} from Theorem \ref{th:lico-cv}.
    
    \begin{theorem}[Convergence of \algname{DESTROY}, deterministic case]
    \label{th:destroy-cv}
    Suppose that $\gamma\in (0,2/\nu)$ and that $\tau\gamma\|\widehat{\mW}\|\leq 1$. Then in \algname{DESTROY},    each $(x_i^k)_{k}$ converges to the same
    solution $x^\star$ to the problem \eqref{eq:decentralized} and each $(a_i^k)_{k}$ converges to $a_i^\star=-\nabla f_i(x^\star)$.
\end{theorem}

Theorem \ref{th:cvx:lico} can be applied to the stochastic case, with $\cO(1/k)$ convergence of the Lagrangian gap, where $\mathcal{Y}=\mathcal{X}$ and $\mL =\mL^* = \mW^{1/2}$.

Similarly, Theorem \ref{th:LV0} yields linear convergence of \algname{DESTROY} in the strongly convex case $\mu>0$, with $\mL^*\mL$  replaced by $\mW$ and $\|\mL\|^2$ replaced by $\|\mW\|=\|\widehat{\mW}\|$. In particular, in the deterministic case $g_i^{k+1}=\nabla f_i(x^k)$, with $\gamma=1/\nu$ and $\tau\gamma =\theta/\|\mW\|$ for some fixed $\theta\in(0,1)$, 
$\varepsilon$-accuracy is reached after
   \[
   	\cO\left(\left(\frac{\nu}{\mu}+\frac{ \|\mW \|}{\lambda_{\min}^+(\mW)}\right)\log\frac{1}{\varepsilon}\right)
   \]
   iterations.


\chapter{Published Papers}

\begin{quote}
	\cite{hanzely2018sega} \bibentry{hanzely2018sega}.
\end{quote}

\begin{quote}
	\cite{mishchenko201999} \bibentry{mishchenko201999}.
\end{quote} 

\begin{quote}
	\cite{khaled2020tighter} \bibentry{khaled2020tighter}.
\end{quote} 

\begin{quote}
	\cite{malitsky2019adaptive} \bibentry{malitsky2019adaptive}.
\end{quote} 

\begin{quote}
	\cite{mishchenko2018delay} \bibentry{mishchenko2018delay}.
\end{quote} 

\begin{quote}
	\cite{mishchenko2019revisiting} \bibentry{mishchenko2019revisiting}.
\end{quote} 

\begin{quote}
	\cite{mishchenko2020distributed} \bibentry{mishchenko2020distributed}.
\end{quote} 

\begin{quote}
	\cite{MKR2020rr} \bibentry{MKR2020rr}.
\end{quote} 

\begin{quote}
	\cite{soori2020dave} \bibentry{soori2020dave}.
\end{quote} 

\chapter{Preprints}

\begin{quote}
	\cite{kovalev2019stochastic} \bibentry{kovalev2019stochastic}.
\end{quote} 

\begin{quote}
	\cite{mishchenko2018stochastic} \bibentry{mishchenko2018stochastic}.
\end{quote} 

\begin{quote}
	\cite{mishchenko2019distributed} \bibentry{mishchenko2019distributed}.
\end{quote} 

\begin{quote}
	\cite{mishchenko2019sinkhorn} \bibentry{mishchenko2019sinkhorn}.
\end{quote} 

\begin{quote}
	\cite{mishchenko2019self} \bibentry{mishchenko2019self}.
\end{quote} 

\begin{quote}
	\cite{mishchenko2019stochastic} \bibentry{mishchenko2019stochastic}.
\end{quote} 

\begin{quote}
	\cite{mishchenko2021proximal} \bibentry{mishchenko2021proximal}.
\end{quote} 

\begin{quote}
	\cite{mishchenko2021intsgd} \bibentry{mishchenko2021intsgd}.
\end{quote} 

\begin{quote}
	\cite{qian2019miso} \bibentry{qian2019miso}.
\end{quote} 

\begin{quote}
	\cite{salim2020dualize} \bibentry{salim2020dualize}.
\end{quote} 

\begin{quote}
	\cite{horvath2019stochastic} \bibentry{horvath2019stochastic}.
\end{quote} 

\begin{quote}
	\cite{khaled2019analysis} \bibentry{khaled2019analysis}.
\end{quote}

\end{document}